\documentclass[english, 12pt]{amsart}

\usepackage{amsmath}
\usepackage{amssymb} 
\usepackage{amscd}
\usepackage{tabularx}
\usepackage[all,cmtip,line]{xy}

\newtheorem{thm}{Theorem}[section]
\newtheorem{prop}[thm]{Proposition}
\newtheorem{lem}[thm]{Lemma}
\newtheorem{cor}[thm]{Corollary}

\newtheorem{ithm}{Theorem}

\numberwithin{equation}{subsection}

\newcommand{\onto}{\rightarrow \hspace{-.86em} \rightarrow}

\newcommand{\ra}{\rightarrow}
\newcommand{\lra}{\longrightarrow}
\newcommand{\la}{\leftarrow}
\newcommand{\lla}{\longleftarrow}
\newcommand{\into}{\hookrightarrow}

\newcommand{\ua}{\uparrow}
\newcommand{\da}{\downarrow}
\newcommand{\iso}{\stackrel{\sim}{\ra}}
\newcommand{\liso}{\stackrel{\sim}{\lra}}

\newcommand{\pfbegin}{{{\em Proof:}\;}}
\newcommand{\pfend}{$\Box$ \medskip}
\newcommand{\mat}[4]{\left( \begin{array}{cc} {#1} & {#2} \\ {#3} & {#4}
\end{array} \right)}

\newlength{\ownl}

\newcommand{\norm}{{\mbox{\bf N}}}
\newcommand{\ndiv}{{\mbox{$\not| $}}}

\newcommand{\Aff}{{\operatorname{Aff}}}

\newcommand{\Aut}{{\operatorname{Aut}\,}}
\newcommand{\BC}{{\operatorname{BC}\,}}

\newcommand{\Coord}{{\operatorname{Coord}}}

\newcommand{\End}{{\operatorname{End}\,}}

\newcommand{\Frob}{{\operatorname{Frob}}}
\newcommand{\Gal}{{\operatorname{Gal}\,}}

\newcommand{\Herm}{{\operatorname{Herm}}}
\newcommand{\Hom}{{\operatorname{Hom}\,}}

\renewcommand{\Im}{{\operatorname{Im}\,}}
\newcommand{\Ind}{{\operatorname{Ind}\,}}

\newcommand{\Lie}{{\operatorname{Lie}\,}}
\newcommand{\coker}{{\operatorname{coker}\,}}

\newcommand{\Pol}{{\operatorname{Pol}}}

\newcommand{\RS}{{\operatorname{RS}}}

\newcommand{\WD}{{\operatorname{WD}}}

\newcommand{\Std}{{\operatorname{Std}}}

\newcommand{\Spec}{{\operatorname{Spec}\,}}

\newcommand{\Spf}{{\operatorname{Spf}\,}}
\newcommand{\Sp}{{\operatorname{Sp}\,}}
\newcommand{\sw}{{\operatorname{sw}}}

\newcommand{\ad}{{\operatorname{ad}\,}}
\newcommand{\dlog}{{\operatorname{dlog}\,}}

\newcommand{\gr}{{\operatorname{gr}\,}}

\newcommand{\rec}{{\operatorname{rec}}}
\newcommand{\tr}{{\operatorname{tr}\,}}

\newcommand{\nind}{{\operatorname{n-Ind}\,}}

\newcommand{\diag}{{\operatorname{diag}}}

\newcommand{\std}{{\operatorname{Std}}}
\newcommand{\hasse}{{\operatorname{Hasse}}}

\newcommand{\an}{{\operatorname{an}}}

\newcommand{\can}{{\operatorname{can}}}

\newcommand{\isog}{{\operatorname{isog}}}

\newcommand{\KS}{{\operatorname{KS}}}
\newcommand{\herm}{{\operatorname{herm}}}

\newcommand{\lin}{{\operatorname{lin}}}
\newcommand{\Int}{{\operatorname{Int}}}

\newcommand{\mini}{{\operatorname{min}}}

\newcommand{\nc}{{\operatorname{nc}}}

\newcommand{\ord}{{\operatorname{ord}}}
\newcommand{\nord}{{\operatorname{n-ord}}}

\newcommand{\rig}{{\operatorname{rig}}}
\newcommand{\semis}{{\operatorname{ss}}}
\newcommand{\sub}{{\operatorname{sub}}}
\newcommand{\Fsemis}{{\operatorname{F-ss}}}

\newcommand{\TF}{{\operatorname{TF}}}
\newcommand{\tor}{{\operatorname{tor}}}

\newcommand{\univ}{{\operatorname{univ}}}

\newcommand{\A}{{\mathbb{A}}}

\newcommand{\C}{{\mathbb{C}}}

\newcommand{\F}{{\mathbb{F}}}
\newcommand{\G}{{\mathbb{G}}}
\newcommand{\HH}{{\mathbb{H}}}

\newcommand{\PP}{{\mathbb{P}}}
\newcommand{\Q}{{\mathbb{Q}}}
\newcommand{\R}{{\mathbb{R}}}

\newcommand{\T}{{\mathbb{T}}}
\newcommand{\U}{{\mathbb{U}}}

\newcommand{\Z}{{\mathbb{Z}}}

\newcommand{\cA}{{\mathcal{A}}}

\newcommand{\cD}{{\mathcal{D}}}
\newcommand{\cE}{{\mathcal{E}}}
\newcommand{\cF}{{\mathcal{F}}}
\newcommand{\cG}{{\mathcal{G}}}
\newcommand{\cH}{{\mathcal{H}}}
\newcommand{\cI}{{\mathcal{I}}}
\newcommand{\cJ}{{\mathcal{J}}}

\newcommand{\cL}{{\mathcal{L}}}
\newcommand{\cM}{{\mathcal{M}}}
\newcommand{\cN}{{\mathcal{N}}}
\newcommand{\cO}{{\mathcal{O}}}
\newcommand{\cP}{{\mathcal{P}}}

\newcommand{\cS}{{\mathcal{S}}}
\newcommand{\cT}{{\mathcal{T}}}

\newcommand{\cX}{{\mathcal{X}}}
\newcommand{\cY}{{\mathcal{Y}}}
\newcommand{\cZ}{{\mathcal{Z}}}

\newcommand{\co}{{{o}}}
\newcommand{\cm}{{{m}}}

\newcommand{\gA}{{\mathfrak{A}}}

\newcommand{\gC}{{\mathfrak{C}}}

\newcommand{\gE}{{\mathfrak{E}}}
\newcommand{\gF}{{\mathfrak{F}}}

\newcommand{\gH}{{\mathfrak{H}}}

\newcommand{\gQ}{{\mathfrak{Q}}}

\newcommand{\gT}{{\mathfrak{T}}}
\newcommand{\gU}{{\mathfrak{U}}}

\newcommand{\gX}{{\mathfrak{X}}}
\newcommand{\gY}{{\mathfrak{Y}}}

\newcommand{\gm}{{\mathfrak{m}}}

\newcommand{\gp}{{\mathfrak{p}}}
\newcommand{\gq}{{\mathfrak{q}}}

\newcommand{\gt}{{\mathfrak{t}}}

\newcommand{\barA}{\overline{{A}}}

\newcommand{\barK}{\overline{{K}}}

\newcommand{\barT}{\overline{{T}}}

\newcommand{\barX}{\overline{{X}}}
\newcommand{\barY}{\overline{{Y}}}
\newcommand{\barZ}{\overline{{Z}}}

\newcommand{\barFF}{\overline{{\F}}}

\newcommand{\barQQ}{\overline{{\Q}}}

\newcommand{\bare}{\overline{{e}}}

\newcommand{\bark}{\overline{{k}}}

\newcommand{\barv}{\overline{{v}}}

\newcommand{\bary}{{\overline{{y}}}}

\newcommand{\tC}{\widetilde{{C}}}

\newcommand{\tG}{\widetilde{{G}}}

\newcommand{\tL}{\widetilde{{L}}}

\newcommand{\tN}{\widetilde{{N}}}

\newcommand{\tP}{\widetilde{{P}}}

\newcommand{\tS}{\widetilde{{S}}}
\newcommand{\tT}{\widetilde{{T}}}
\newcommand{\tU}{\widetilde{{U}}}

\newcommand{\ta}{\widetilde{{a}}}

\newcommand{\te}{\widetilde{{e}}}

\newcommand{\tg}{\widetilde{{g}}}

\newcommand{\ty}{\widetilde{{y}}}

\newcommand{\bargamma   }{\overline{\gamma}}

\newcommand{\barSigma   }{\overline{\Sigma}}

 \newcommand{\tpi    }{{\widetilde{\pi}} }

\newcommand{\tDelta         }{\widetilde{\Delta}}

\newcommand{\tXi   }{\widetilde{\Xi}}   
\newcommand{\tPi   }{\widetilde{\Pi}}   
\newcommand{\tSigma   }{\widetilde{\Sigma}}

\newcommand{\tOmega   }{\widetilde{\Omega}}   

\newcommand{\hatotimes}{{\widehat{\otimes}}}

\newcommand{\hata}{{\widehat{a}}}
\newcommand{\hatZ}{{\widehat{\Z}}}
\newcommand{\Zhat}{{\widehat{\Z}}}
\newcommand{\tcE}{{\widetilde{\cE}}}
\newcommand{\tcI}{{\widetilde{\cI}}}

\newcommand{\tcT}{{\widetilde{\cT}}}
\newcommand{\tcG}{{\widetilde{\cG}}}

\newcommand{\barcI}{\overline{{\cI}}}

\newcommand{\tcS}{\widetilde{{\cS}}}
\newcommand{\hS}{\widehat{{S}}}
\newcommand{\hcS}{\widehat{{\cS}}}

\begin{document}

\title{On the Rigid Cohomology of Certain Shimura Varieties.}

 \date{\today}

\author{Michael Harris}\email{harris@math.jussieu.fr}\address{Institut de Math\'{e}matiques de Jussieu,
Paris, France.}
\author{Kai-Wen Lan} \email{kwlan@math.umn.edu}\address{School of Mathematics, University of Minnesota, MN, USA.}
\author{Richard Taylor}\email{rtaylor@ias.edu}\address{School of Mathematics, IAS, Princeton, NJ, USA.} \thanks{We would all like to thank the Institute for Advanced Study for its support and hospitality. This project was begun, and the key steps completed, while we were all attending the special IAS special year on `Galois representations and automorphic forms'. M.H.'s research received funding from the European Research Council under the European Community's Seventh Framework Programme (FP7/2007-2013) / ERC Grant agreement number 290766 (AAMOT). K.-W. L.'s research was partially supported by NSF Grants DMS-1069154 and DMS-1258962.   R.T.'s research was partially supported by NSF
Grants DMS-0600716, DMS-1062759 and DMS-1252158, and by the IAS Oswald Veblen and Simonyi Funds.  During some of the period when this research was being written up J.T. served as a Clay Research Fellow.  }
\author{Jack Thorne}\email{thorne@math.harvard.edu}\address{Dept. of Mathematics, Harvard University, Cambridge, MA, USA.}

\begin{abstract}  
We construct the compatible system of $l$-adic representations associated to a regular algebraic cuspidal automorphic representation of $GL_n$ over a CM (or totally real) field and check local-global compatibility for the $l$-adic representation  away from $l$ and finite number of rational primes above which the CM field or the automorphic representation ramify. The main innovation is that we impose no self-duality hypothesis on the automorphic representation.
\end{abstract}

\maketitle

\section*{Introduction}

Our main theorem is as follows (see corollary \ref{mt3}).

\begin{ithm}\label{thma} Let $p$ denote a rational prime and let $\imath: \barQQ_p \iso \C$. Suppose that $E$ is a CM (or totally real) field and that $\pi$ is a cuspidal automorphic representation of $GL_n(\A_E)$ such that $\pi_{\infty}$ has the same infinitesimal character as an irreducible algebraic representation $\rho_\pi$ of $\RS^E_\Q GL_n$. Then there is a unique continuous semi-simple representation 
\[ r_{p,\imath}(\pi):G_E \lra GL_{n}(\barQQ_p) \]
such that, if $q \neq p$ is a prime above which $\pi$ and $E$ are unramified and if $v|q$ is a prime of $E$, then $r_{p,\imath}(\pi)$ is unramified at $v$ and
\[ r_{p,\imath}(\pi)|_{W_{E_v}}^\semis = \imath^{-1} \rec_{E_v}(\pi_{v}|\det|_v^{(1-n)/2}). \]
\end{ithm}

Here $\rec_{E_v}$ denotes the local Langlands correspondence for $E_v$. It may be possible to extend the local-global compatibility to other primes $v$. Ila Varma is considering this question.

The key point is that we make no self-duality assumption on $\pi$. In the presence of such a self-duality assumption (`polarizability', see \cite{blggt}) the existence of $r_{p,\imath}(\pi)$ has been known for some years (see \cite{shin} and \cite{ch}). In almost all polarizable cases $r_{p,\imath}(\pi)$ is realized in the cohomology of a Shimura variety, and in all polarizable cases $r_{p,\imath}(\pi)^{\otimes 2}$ is realized in the cohomology of a Shimura variety (see \cite{ana}).  In contrast, according to unpublished computations of one of us (M.H.) and of Laurent Clozel, in the non-polarizable case the representation $r_{p,\imath}(\pi)$ will never occur in the cohomology of a Shimura variety. Rather we construct it as a $p$-adic limit of representations which do occur in the cohomology of Shimura varieties. 

We sketch our argument. We may easily reduce to the case of an imaginary CM field $F$ which contains an imaginary quadratic field in which $p$ splits. For all sufficiently large integers $N$, we construct a $2n$-dimensional representation $R_p(\imath^{-1}(\pi ||\det||^N)^\infty)$ such that for good primes $v$ we have
\[ \begin{array}{l} R_p(\imath^{-1}(\pi ||\det||^N)^\infty)|_{W_{F_v}}^\semis \cong \\ \imath^{-1} \rec_{F_v}(\pi_v |\det|_v^{N+(1-n)/2}) \oplus \imath^{-1} \rec_{F_v}(\pi_v |\det|_v^{N+(1-n)/2})^{\vee,c} \epsilon_p^{1-2n}, \end{array} \]
as a $p$-adic limit of (presumably irreducible) $p$-adic representations associated to polarizable, regular algebraic cuspidal automorphic representations of $GL_{2n}(\A_F)$. It is then elementary algebra to reconstruct $r_{p,\imath}(\pi)$.

We work on the quasi-split unitary similitude group $G_n$ associated to $F^{2n}$. Note that $G_n$ has a maximal parabolic subgroup $P_{n,(n)}^+$ with Levi component 
\[ L_{n,(n)}\cong GL_1 \times RS^F_\Q GL_n. \]
(We will give all these groups integral structures.)
We set
\[ \Pi(N) = \Ind_{P_{n,(n)}^+(\A^{p,\infty})}^{G_n(\A^{p,\infty})} (1 \times \imath^{-1}(\pi ||\det||^N)^{p,\infty}). \]
Then our strategy is to realize $\Pi(N)$, for sufficiently large $N$, in a space of overconvergent $p$-adic {\em cusp} forms for $G_n$ of  finite slope. It is a space of forms of a weight for which we expect no classical forms. Once we have done this, we can use an argument of Katz (see \cite{katz}) to find  congruences modulo arbitrarily high powers of $p$ to classical (holomorphic) cusp forms on $G_n$ (of other weights). (Alternatively it is presumably possible to construct an eigenvariety in this setting, but we have not carried this out.) One can attach Galois representations to these classical cusp forms by using the trace formula to lift them to polarizable, regular algebraic, discrete automorphic representations of $GL_{2n}(\A_F)$ (see e.g. \cite{sws}) and then applying the results of \cite{shin} and \cite{ch}.

We learnt the idea that one might try to realize $\Pi(N)$ in a space of overconvergent $p$-adic {\em cusp} forms for $G_n$ (of finite slope) from Chris Skinner.
The key problem was how to achieve such a realization. To sketch our approach we must first establish some more notation.

To a neat open compact subgroup $U$ of $G_n$ we can associate a Shimura variety $X_{n,U}/\Spec \Q$. It is a moduli space for abelian $n[F:\Q]$-folds with an isogeny action of $F$ and certain additional structures. It is not proper. It has a canonical normal compactification $X_{n,U}^\mini$ and, to certain auxiliary data $\Delta$, one can attach a smooth compactification $X_{n,U,\Delta}$ which naturally lies over $X_{n,U}^\mini$ and whose boundary is a simple normal crossings divisor. To a representation $\rho$ of $L_{n,(n)}$ (over $\Q$) one can attach a locally free sheaf $\cE_{U,\rho}/X_{n,U}$ together with a canonical (locally free) extension $\cE_{U,\Delta,\rho}$ to $X_{n,U,\Delta}$, whose global sections are holomorphic automorphic forms on $G_n$ `of weight $\rho$ and level $U$'. (The space of global sections does not depend on $\Delta$.) The product of $\cE_{U,\Delta,\rho}$ with the ideal sheaf of the boundary of $X_{n,U,\Delta}$, which we denote $\cE_{U,\Delta,\rho}^\sub$, is again locally free and its global sections are holomorphic cusp forms on $G_n$ `of weight $\rho$ and level $U$' (and again the space of global sections does not depend on $\Delta$).

To the schemes $X_{n,U}$, $X_{n,U}^\mini$ and $X_{n,U,\Delta}$ one can naturally attach dagger spaces $X_{n,U}^\dag$, $X_{n,U}^{\mini,\dag}$ and $X_{n,U,\Delta}^\dag$ in the sense of \cite{gkcrelle}. These are like rigid analytic spaces except that one consistently works with overconvergent sections. If $U$ is the product of a neat open compact subgroup of $G_n(\A^{\infty,p})$ and a suitable open compact subgroup of $G_n(\Q_p)$, then one can define admissible open sub-dagger spaces (`the ordinary loci')
\[ X_{n,U}^{\ord,\dag} \subset X_{n,U}^\dag \]
and
\[ X_{n,U}^{\mini,\ord,\dag} \subset X_{n,U}^{\mini,\dag} \]
and
\[ X_{n,U,\Delta}^{\ord,\dag} \subset X_{n,U,\Delta}^\dag. \]
By an overconvergent cusp form of weight $\rho$ and level $U$ one means a section of $\cE_{U,\rho}^\sub$ over $X_{n,U,\Delta}^{\ord,\dag}$. (Again this definition does not depend on the choice of $\Delta$.) 

We write $G_n^{(m)}$ for the semi-direct product of $G_n$ with the additive group with $\Q$-points $\Hom_F(F^m,F^{2n})$, and $P_{n,(n)}^{(m),+}$ for the pre-image of $P_{n,(n)}^+$ in $G_n^{(m)}$. We also write $L_{n,(n)}^{(m)}$ for the semi-direct product of $L_{n,(n)}$ with the additive group with $\Q$-points $\Hom_F(F^m,F^n)$, which is naturally a quotient of $P_{n,(n)}^{(m),+}$. (Again we will give these groups integral structures.) 
To a neat open compact subgroup $U \subset G_n^{(m)}(\A^\infty)$ with projection $U'$ in $G_n(\A^\infty)$ one can attach a (relatively smooth, projective) Kuga-Sato variety $A_{n,U}^{(m)}/X_{n,U'}$. For a cofinal set of $U$ it is an abelian scheme isogenous to the $m$-fold self product of the universal abelian variety over $X_{n,U'}$. To certain auxiliary data $\Sigma$ one can attach a smooth compactification $A_{n,U,\Sigma}^{(m)}$ of $A_{n,U}^{(m)}$ whose boundary is a simple normal crossings divisor; which lies over $X_{n,U}^\mini$; and which, for suitable $\Sigma$ depending on $\Delta$, lies over $X_{n,U',\Delta}$. Thus
\[ \begin{array}{ccc} A_{n,U} & \into & A_{n,U,\Sigma} \\ \da && \da \\ X_{n,U'} & \into & X_{n,U',\Delta} \\ || & & \da \\ X_{n,U'} & \into & X_{n,U'}^\mini .\end{array} \]
We define $A_{n,U}^{(m),\ord,\dag}$ and $A_{n,U,\Sigma}^{(m),\ord,\dag}$ to be the pre-image of $X_{n,U'}^{\mini,\ord,\dag}$ in the dagger spaces associated to $A_{n,U}^{(m)}$ and $A_{n,U,\Sigma}^{(m)}$. 

We will define 
\[ H^i_{c-\partial}(\barA^{(m),\ord}_{n,U}, \barQQ_p) \]
to be the hypercohomology of the complex on $A^{(m),\ord,\dag}_{n,U,\Sigma}$ which is the tensor product of the de Rham complex with log poles towards the boundary, $A^{(m),\ord,\dag}_{n,U,\Sigma}-A^{(m),\ord,\dag}_{n,U}$, and the ideal sheaf defining the boundary. We believe it is a sort of rigid cohomology of the ordinary locus $\barA^{(m),\ord}_{n,U}$ in the special fibre of an integral model of $A_{n,U}^{(m)}$. More specifically cohomology with compact support towards the toroidal boundary, but not towards the non-ordinary locus. Hence our notation. However we have not bothered to verify that this group only depends on ordinary locus in the special fibre. The theory of Shimura varieties provides us with sufficiently canonical lifts that this will not matter to us. Our proof that for $N$ sufficiently large $\Pi(N)$ occurs in the space of overconvergent $p$-adic  cusp forms for $G_n$ proceeds by evaluating $H^i_{c-\partial}(\barA^{(m),\ord}_{n,U}, \barQQ_p)$ in two ways. 

Firstly we use the usual Hodge spectral sequence. The higher direct images from 
$A_{n,U,\Sigma}^{(m)}$ to $X_{n,U',\Delta}$ of the tensor product of the ideal sheaf of the boundary and the sheaf of differentials of any degree with log poles along the boundary, is canonically filtered with graded pieces sheaves of the form $\cE_{U',\Delta,\rho}^\sub$. Thus $H^i_{c-\partial}(\barA^{(m),\ord}_{n,U}, \barQQ_p)$ can be computed in terms of the groups
\[ H^j(X_{n,U,\Delta}^{\ord,\dag}, \cE_{U,\Delta,\rho}^\sub) \]
A crucial observation for us is that for $j>0$ this group vanishes
(see theorem \ref{ttbb} and proposition \ref{admis}).
This observation seems to have been made independently, at about the same time, by Andreatta, Iovita and Pilloni (see \cite{iovita1} and \cite{iovita2}). It seems quite surprising to us. It is false if one replaces $\cE_{U,\Delta,\rho}^\sub$ with $\cE_{U,\Delta,\rho}^\can$. Its proof depends on a number of apparently unrelated facts, including:
\begin{itemize}
\item $X_{n,U}^{\mini,\ord,\dag}$ is affinoid.
\item The stabilizer in $GL_n(\cO_F)$ of a positive definite hermitian $n \times n$ matrix over $F$ is finite.
\item Certain line bundles on self products $A$ of the universal abelian variety over $X_{n',U'}$ (for $n'<n$) are relatively ample for $A/X_{n',U'}$.
\end{itemize}
This observation implies that $H_{c-\partial}^i(\barA^{(m),\ord}_{n,U}, \barQQ_p)$ can be computed by a complex whose terms are spaces of overconvergent cusp forms. Hence it suffices for us to show that, for $N$ sufficiently large, $\Pi(N)$ occurs in
\[ H^i_{c-\partial}(\barA^{(m),\ord}_{n}, \barQQ_p) = \lim_{\ra U, \Sigma} H^i_{c-\partial}(\barA^{(m),\ord}_{n,U}, \barQQ_p) \]
for some $m$ and $i$ (depending on $N$). 

To achieve this we use a second spectral sequence which computes the cohomology group $H^i_{c-\partial}(\barA^{(m),\ord}_{n,U}, \barQQ_p)$ in terms of the rigid cohomology of $\barA^{(m),\ord}_{n,U,\Sigma}$ and its various boundary strata. See section \ref{secrc}. This is an analogue of the spectral sequence 
\[ E_1^{i,j}=H^j(Y^{(i)},\C) \Rightarrow H^{i+j}_c(Y-\partial Y,\C), \]
where $Y$ is a proper smooth variety over $\C$, where $\partial Y$ is a simple normal crossing divisor on $Y$, and where $Y^{(i)}$ is the disjoint union of the $i$-fold intersections of irreducible components of $\partial Y$. (So $Y^{(0)}=Y$.)
Some of the terms in this spectral sequence seem a priori to be hard to control, e.g. $H^i_\rig( 
\barA^{(m),\ord}_{n,U,\Sigma})$. However employing theorems about rigid cohomology due to Berthelot and Chiarellotto, we see that the eigenvalues of Frobenius on $H^i_{c-\partial}(\barA^{(m),\ord}_{n,U,\Sigma}, \barQQ_p)$ are all Weil $p^j$-numbers for $j\geq 0$. Moreover the weight $0$ part, $W_0 H^i_{c-\partial}(\barA^{(m),\ord}_{n,U}, \barQQ_p)$, equals the cohomology of a complex only involving the rigid cohomology in degree $0$ of $\barA^{(m),\ord}_{n,U}$ and its various boundary strata. (See proposition \ref{toptorig}.) This should have a purely combinatorial description. More precisely we define a simplicial complex $\cS(\partial \barA^{(m),\ord}_{n,U,\Sigma})$ whose vertices correspond to boundary components of $\barA^{(m),\ord}_{n,U,\Sigma}$ and whose $j$-faces correspond to $j$-boundary components with non-trivial intersection. For $i>0$ we obtain an isomorphism
\[ H^i(|\cS(\partial \barA^{(m),\ord}_{n,U,\Sigma})|,\barQQ_p) \cong W_0 H_{c-\partial}^{i+1}(\barA^{(m),\ord}_{n,U}, \barQQ_p). \]
Thus it suffices to show that for $N$ sufficiently large $\Pi(N)$ occurs in
\[ H^i(|\cS(\partial \barA^{(m),\ord}_{n})|,\barQQ_p)=\lim_{\ra U, \Sigma} H^i(|\cS(\partial \barA^{(m),\ord}_{n,U,\Sigma})|,\barQQ_p) \]
for some $m$ and some $i>0$ (possibly depending on $N$).

The boundary of $\barA^{(m),\ord}_{n,U,\Sigma}$ comes in pieces indexed by the conjugacy classes of maximal parabolic subgroups of $G_n$. We shall be interested in the union of the irreducible components which are associated to $P_{n,(n)}^+$. These correspond to an open subset $|\cS(\partial \barA^{(m),\ord}_{n,U,\Sigma})|_{=n}$ of $|\cS(\partial \barA^{(m),\ord}_{n,U,\Sigma})|$. As $|\cS(\partial \barA^{(m),\ord}_{n,U,\Sigma})|$ is compact, the interior cohomology 
\[ H^i_\Int(|\cS(\partial \barA^{(m),\ord}_{n})|_{=n},\barQQ_p) = \lim_{\ra U, \Sigma} H^i_\Int(|\cS(\partial \barA^{(m),\ord}_{n,U,\Sigma})|_{=n},\barQQ_p) \]
is naturally a sub-quotient of $H^i(|\cS(\partial \barA^{(m),\ord}_{n})|,\barQQ_p)$. (By interior cohomology we mean the image of the cohomology with compact support in the cohomology. The interior cohomology of an open subset of an ambient compact Hausdorff space is naturally a sub-quotient of the cohomology of that ambient space.) Thus it even suffices to show that for $N$ sufficiently large $\Pi(N)$ occurs in
\[ H^i_\Int(|\cS(\partial \barA^{(m),\ord}_{n})|_{=n},\barQQ_p) \]
for some $m$ and some $i>0$ (possibly depending on $N$).

However the data $\Sigma$ is a $G_n^{(m)}(\Q)$-invariant (glued) collection of polyhedral cone decompositions and $\cS(\partial \barA^{(m),\ord}_{n})$ is obtained from $\Sigma$ by replacing $1$-cones by vertices, $2$-cones by edges etc. The cones corresponding to $|\cS(\partial \barA^{(m),\ord}_{n})|_{=n}$ are a disjoint union of polyhedral cones in the space of positive definite hermitian forms on $F^n$. From this one obtains an equality
\[ |\cS(\partial \barA^{(m),\ord}_{n})|_{=n}= \coprod_{h \in P_{n,(n)}^{+}(\A^{p,\infty}) \backslash G_n(\A^{p,\infty})/U^p} \gT^{(m)}_{(n), hUh^{-1} \cap P_{n,(n)}^{(m),+}(\A^\infty)}, \]
where
\[ \gT^{(m)}_{(n), U'} = 
L_{n,(n)}^{(m)}(\Q) \backslash L_{n,(n)}^{(m)}(\A)/ U' (\R^\times_{>0} \times (U(n)^{[F^+:\Q]}\R^\times_{>0})), \]
with $U(n)$ denoting the usual $n \times n$ compact unitary group. We deduce that
\[ H^i_\Int(|\cS(\partial \barA^{(m),\ord}_{n})|_{=n},\barQQ_p) = \Ind_{P_{n,(n)}^{+}(\A^{\infty,p})}^{G_n(\A^{p,\infty})} H^i_\Int (\gT_{(n)}^{(m)},\barQQ_p)^{\Z_p^\times}, \]
where 
\[ H^i_\Int (\gT_{(n)}^{(m)},\barQQ_p) = \lim_{\ra U'} H^i_\Int ( \gT^{(m)}_{(n), U'}, \barQQ_p) \]
as $U'$ runs over neat open compact subgroups of $L_{n,(n)}^{(m)}(\A^\infty)$. (The $\Z_p^\times$-invariants results from a restriction on the open compact subgroups of $G_n(\A^\infty)$ that we are considering.) Thus it suffices to show that for all sufficiently large $N$, the representation $1 \times (\pi ||\det||^N)^{p, \infty}$ occurs in $H^i_\Int(T_{(n)}^{(m)}, \C)$ for some $i>0$ and some $m$ (possibly depending on $N$). 

We will write simply $\gT_{(n),U'}$ for $\gT^{(0)}_{(n),U'}$, a locally symmetric space associated to $L_{n,(n)}\cong GL_1 \times \RS_\Q^F GL_n$. If $\rho$ is a representation of $L_{n,(n)}$ over $\C$, then it gives rise to a locally constant sheaf $\cL_{\rho,U'}$ over $\gT_{(n),U'}$. We set
\[ H^i_\Int(\gT_{(n)},\cL_\rho) = \lim_{\ra U'} H^i_\Int(\gT_{(n),U'},\cL_{\rho,U'}), \]
a smooth $L_{n,(n)}(\A^\infty)$-module. 
The space $\gT^{(m)}_{(n), U'}$ is an $(S^1)^{nm[F:\Q]}$-bundle over the locally symmetric space $\gT^{(0)}_{(n),U'}$ and if $\pi^{(m)}$ denotes the fibre map then
\[ R^j\pi^{(m)}_* \C \cong \cL_{\wedge^j\Hom_F(F^m,F^n)^\vee \otimes_\Q \C,U'}, \]
where $L_{n,(n)}$ acts on $\Hom_F(F^m,F^n)$ via projection to $\RS^F_\Q GL_n$. Moreover the Leray spectral sequence
\[ E_2^{i,j}=H^i_\Int(\gT_{(n)},\cL_{\wedge^j\Hom_F(F^m,F^n)^\vee \otimes_\Q \C}) \Rightarrow H^{i+j}_\Int(\gT^{(m)}_{(n)} ,\C) \]
degenerates at the second page. (This can be seen by considering the action of the centre of $L_{n,(n)}(\A^\infty)$.) Thus it suffices to show that for all sufficiently large $N$, we can find non-negative integers $j$ and $m$ and an irreducible constituent $\rho$ of $\wedge^j\Hom_F(F^m,F^n)^\vee \otimes_\Q \C$ such that 
the representation $1 \times (\pi ||\det||^N)^{p, \infty}$ occurs in $H^i_\Int(T_{(n)}, \cL_\rho)$ for some $i \in \Z_{>0}$.  Clozel \cite{clozelaa} checked that (for $n>1$) this will be the case as long as $1 \times (\pi||\det||^N)_\infty$ has the same infinitesimal character as some irreducible constituent of $\wedge^j\Hom_F(F^m,F^n) \otimes_\Q \C$, i.e. if $\rho_\pi \otimes (\norm_{F/\Q} \circ \det)^N$ occurs in $\wedge^j\Hom_F(F^m,F^n) \otimes_\Q \C$. 
From Weyl's construction of the irreducible representations of $GL_n$, for large enough $N$ this will indeed be the case for some $m$ and $j$.

We remark it is essential to work with $N$ sufficiently large. It is not an artifact of the fact that we are working with Kuga-Sato varieties rather than local systems on the Shimura variety. We can twist a local system on the Shimura variety by a power of the multiplier character of $G_n$. However the restriction of the multiplier factor of $G_n$ to $L_{n,(n)} \cong GL_1 \times \RS^F_\Q GL_n$ factors through the $GL_1$-factor and does not involve the $\RS^F_\Q GL_n$ factor.

We learnt from the series of papers \cite{hz1}, \cite{hz2}, \cite{hz3}, the key observation that $|\cS(\partial A_{n,U,\Sigma}^{(m)})|$ has a nice geometric interpretation involving the locally symmetric space for $L_{n,(n)}$ and that this could be used to calculate cohomology. 

Although the central argument we have sketched above is not long, this paper has unfortunately become very long. If we had only wanted to construct $r_{p,\imath}(\pi)$ for all but finitely many primes $p$, then the argument would have been significantly shorter as we could have worked only with Shimura varieties $X_{n,U}$ which have good integral models at $p$. The fact that we want to construct $r_{p,\imath}(\pi)$ for all $p$ adds considerable technical complications and also requires appeal to the recent work \cite{kw2}. (Otherwise we would only need to appeal to \cite{kw1} and \cite{kw1.5}.) 

Another reason this paper has grown in length is the desire to use a language to describe toroidal compactifications of mixed Shimura varieties that is different from the language used in \cite{kw1}, \cite{kw1.5} and \cite{kw2}. We do this because at least one of us (R.T.) finds this language clearer. In any case it would be necessary to establish a substantial amount of notation regarding toroidal compactifications of Shimura varieties, which would require significant space. We hope that the length of the paper, and the technicalities with which we have to deal, won't obscure the main line of the argument.

After we announced these results, but while we were writing up this paper, Scholze found another proof of theorem \ref{thma}, relying on his theory of perfectoid spaces. His arguments seem to be in many ways more robust. For instance he can handle torsion in the cohomology of the locally symmetric varieties associated  to $GL_n$ over a CM field. Scholze's methods have some similarities with ours. Both methods first realize the Hecke eigenvalues of interest in the cohomology with compact support of the open Shimura variety by an analysis of the boundary and then show that they also occur in some space of $p$-adic cusp forms. We work with the ordinary locus of the Shimura variety, which for the minimal compactification is affinoid. Scholze works with the whole Shimura variety, but at infinite level. He (very surprisingly) shows that at infinite level, as a perfectoid space, the Shimura variety has a Hecke invariant affinoid cover.
\newpage

\subsection*{Notation}

{\em If $G \onto H$ is a  surjective group homomorphism and if $U$ is a subgroup of $G$ we will sometimes use $U$ to also denote the image of $U$ in $H$.

If $f:X \ra Y$ and $f':Y \ra Z$ then we will denote by $f' \circ f:X \ra Z$ the composite map $f$ followed by $f'$.
In this paper we will use both left and right actions. Suppose that $G$ is a group acting on a set $X$ and that $g,h \in G$. If $G$ acts on $X$ on the left we will write $gh$ for $g \circ h$. If $G$ acts on $X$ on the right we will write $hg$ for $g \circ h$. 
}

If $G$ is a group (or group scheme) then $Z(G)$ will denote its centre. 

We will write $S_n$ for the symmetric group on $n$ letters. 
We will write $U(n)$ for the group of $n \times n$ complex matrices $h$ with ${}^th{}^ch=1_n$. 

If $G$ is an abelian group we will write $G[\infty]$ for the torsion subgroup of $G$,  $G[\infty^p]$ for the subgroup of elements of order prime to $p$, and $G^\TF = G/G[\infty]$. We will write
$TG= \lim_{\la N} G[N]$ and $T^pG=\lim_{\la p \ndiv N}G[N]$. We will also write $VG=TG \otimes_\Z \Q$ and $V^pG=T^pG \otimes_\Z \Q$.

If $A$ is a ring, if $B$ is a locally free, finite $A$-algebra, and if $X/\Spec B$ is a quasi-projective scheme; then we will let $\RS^B_A X$ denote the restriction of scalars (or Weil restriction) of $X$ from $B$ to $A$. (See for instance section 7.6 of \cite{blr}.)

By a $p$-adic formal scheme we mean a formal scheme such that $p$ generates an ideal of definition.

If $X$ is an $A$-module and $B$ is an $A$-algebra, we will sometimes write $X_B$ for $X \otimes_A B$. We will also use $X$ to denote the abelian group scheme over $A$ defined by 
\[ X(B)=X \otimes_A B = X_B \]
for all $A$-algebras $B$. 

If $Y$ is a scheme and if $G_1,G_2/Y$ are group schemes then we will let 
\[ \cH\co\cm(G_1,G_2)\]
 denote the Zariski sheaf on $Y$ whose sections over an open $W$ are 
\[ \Hom(G_1|_W,G_2|_W).\]
 If in addition $R$ is a ring then we will let 
$\cH\co\cm(G_1,G_2)_R$ denote the tensor product $\cH\co\cm(G_1,G_2)\otimes_\Z R$ and we will let $\Hom(G_1,G_2)_R$ denote the $R$-module of global sections of $\cH\co\cm(G_1,G_2)_R$. If $Y$ is noetherian this is the same as $\Hom(G_1,G_2) \otimes_\Z R$, but for a general base $Y$ it may differ.

If $\cS$ is a simplicial complex we will write $|\cS|$ for the corresponding topological space. To a scheme or formal scheme $Z$, we will associate to it a simplicial complex $\cS(Z)$ as follows: Let $\{ Z_j\}_{j \in J}$ denote the set of irreducible components of $Z$. Then $J$ will be the set of vertices of $\cS(Z)$ and a subset $K \subset J$ will span a simplex if and only if $\bigcap_{j \in K} Z_j \neq \emptyset$. 

If $F$ is a field then $G_F$ will denote its absolute Galois group. If $F$ is a number field and $F_0 \subset F$ is a subfield and $S$ is a finite set of primes of $F_0$, then we will denote by $G_F^S$ the maximal continuous quotient of $G_F$ in which all primes of $F$ not lying above an element of $S$ are unramified.

Suppose that $F$ is a number field and that $v$ is a place of $F$. If $v$ is finite we will write $\varpi_v$ for a uniformizer in $F_v$ and $k(v)$ for the residue field of $v$. We will write $| \,\,\,|_v$ for the absolute value on $F$ associated to $v$ and normalized as follows: 
\begin{itemize}
\item if $v$ is finite then $|\varpi_v|_v = (\# k(v))^{-1}$;
\item if $v$ is real then $|x|_v=\pm x$;
\item if $v$ is complex then $|x|_v={}^cxx$.
\end{itemize}
We write 
\[ ||\,\,\,||_F = \prod_v |\,\,\,|_v: \A_F^\times \lra \R_{>0}^\times . \]
We will write $\cD_F^{-1}$ for the inverse different of $\cO_F$.

If $w \in \Z_{>0}$ and $p$ is a prime number then by a Weil $p^w$-number we mean an element $\alpha \in \barQQ$ which is an integer away from $p$ and such that for each infinite place $v$ of $\barQQ$ we have $|\alpha|_v=p^w$. 

Suppose that $v$ is finite and that
\[ r: G_{F_v} \lra GL_n(\barQQ_l) \]
is a continuous representation, which in the case $v|l$ we assume to be de Rham. Then we will write $\WD(r)$ for the corresponding Weil-Deligne representation of of the Weil group $W_{F_v}$ of $F_v$ (see for instance
section 1 of \cite{ty}.) If $\pi$ is an irreducible smooth representation of $GL_n(F_v)$ over $\C$ we will write $\rec_{F_v}(\pi)$ for the Weil-Deligne representation of $W_{F_v}$ corresponding to $\pi$ by the local Langlands conjecture (see for instance the introduction to \cite{ht}). If $\pi_i$ is an irreducible smooth representation of $GL_{n_i}(F_v)$ over $\C$ for $i=1,2$ then there is an irreducible smooth representation $\pi_1 \boxplus \pi_2$ of $GL_{n_1+n_2}(F_v)$ over $\C$ satisfying
\[ \rec_{F_v}(\pi_1 \boxplus \pi_2)=\rec_{F_v}(\pi_1) \oplus \rec_{F_v}(\pi_2). \]

Suppose that $G$ is a reductive group over $F_v$ and that $P$ is a parabolic subgroup of $G$ with unipotent radical $N$ and Levi component $L$. Suppose also that $\pi$ is a smooth representation of $L(F_v)$ on a vector space $W_\pi$ over a field $\Omega$ of characteristic $0$. We will define 
\[ \Ind_{P(F_v)}^{G(F_v)} \pi \]
to be the representation of $G(F_v)$ by right translation on the set of locally constant functions
\[ \varphi: G(F_v) \lra W_\pi \]
such that
\[ \varphi(hg)=\pi(h) \varphi(g) \]
for all $h \in P(F_v)$ and $g \in G(F_v)$. In the case $\Omega=\C$ we also define
\[ \nind_{P(F_v)}^{G(F_v)} \pi = \Ind_{P(F_v)}^{G(F_v)} \pi \otimes \delta_P^{1/2} \]
where 
\[ \delta_P(h)^{1/2} = | \det( \ad(h)|_{\Lie N})|^{1/2}_v. \]

If $G$ is a linear algebraic group over $F$ then the concept of a {\em neat} open compact subgroup of $G(\A_F^\infty)$ is defined for instance in section 0.6 of \cite{pink}.
\newpage

\section{Some algebraic groups and automorphic forms.}\label{gaf}

{\em For the rest of this paper 
fix the following notation. Let $F^+$ be a totally real field and $F_0$ an imaginary quadratic field, and set $F=F_0F^+$.
Write $c$ for
the non-trivial element of
$\Gal(F/F^+)$. 
Also choose a rational prime $p$ which splits in $F_0$ 
and choose an element $\delta_F \in \cO_{F,(p)}$ with $\tr_{F/F^+} \delta_F=1$ (which is possible as $p$ is unramified in $F/F^+$).

Fix an isomorphism $\imath:\barQQ_p \iso \C$. Fix a choice of $\sqrt{p} \in \barQQ_p$ by $\imath \sqrt{p} >0$. If $v$ is a prime of $F$ and $\pi$ an irreducible admissible representation of $GL_m(F_v)$ over $\barQQ_l$ define
\[ \rec_{F_v}(\pi)=\imath^{-1} \rec_{F_v}(\imath \pi) \]
a Weil-Deligne representation of $W_{F_v}$ over $\barQQ_l$.
 
 Let $n$ be a non-negative integer. 
  We will often attach $n$ as a subscript to other notation, when we need to record the particular choice of $n$ we are working with, but, at other times when the choice of $n$ is clear, we may drop it from the notation.}
 
\subsection{Three algebraic groups.}\label{s1.1}

Write $\Psi_n$ for the $n \times n$-matrix with $1$'s on the anti-diagonal and
$0$'s elsewhere, and set
\[ J_{n}= \mat{0}{\Psi_n}{-\Psi_n}{0} \in GL_{2n}(\Z). \]
Let 
\[  \Lambda_n=(\cD^{-1}_F)^n \oplus \cO_{F}^{n}, \]
and define a perfect pairing 
\[ \langle \,\,\, ,\,\,\, \rangle_n: \Lambda_n \times \Lambda_n \lra \Z \]
by
\[ \langle x,y \rangle_n = \tr_{F/\Q} {}^tx J_{n} {}^cy. \]
We will write $V_n$ for $\Lambda_n \otimes \Q$.
Let $G_n$ denote denote the group scheme over $\Z$
defined by
\[ G_n(R) =\{ (g,\mu) \in \Aut((\Lambda_n \otimes_{\Z} R)/(\cO_F \otimes_\Z R)) \times R^\times:\,\,
{}^tgJ_n {}^cg = \mu J_n\}, \]
for any ring $R$, and let $\nu:G_n \ra GL_1$ denote the multiplier character which sends $(g,\mu)$
to $\mu$. Then $G_n$ is a quasi-split connected reductive group scheme over
$\Z[1/D_{F/\Q}]$ (where $D_{F/\Q}$ denotes the discriminant of $F/\Q$) and
splits over $\cO_{F^\nc}[1/D_{F/\Q}]$ (where $F^{\nc}$ denotes the normal closure of $F/\Q$).  In particular $G_0$ will denote $GL_1$
and $\nu:G_0 \ra GL_1$ is the identity map. 

If $n>0$ set
\[ C_n = \G_m \times \ker(N_{F/F^+}: \RS_{\Z}^{\cO_F} \G_m \lra \RS^{\cO_{F^+}}_{\Z} \G_m). \]
Then there is a natural map
\[ \begin{array}{rcl} G_n & \lra & C_n \\
(g,\mu) & \longmapsto & (\mu, \mu^{-n} \det g). \end{array} \]
If $n=0$ we set $C_0 = \G_m$ and let $G_0 \lra C_0$ denote the map $\nu$. In either case this map identifies $C_n$ with $G_n/[G_n,G_n]$. 

We will write $\Lambda_{n,(i)}$ for the submodule of $\Lambda_n$ consisting of elements whose last $2n-i$ entries are $0$, and $V_{n,(i)}$ for $\Lambda_{n,(i)} \otimes \Q$. If $W$ is a submodule of $\Lambda_n$ we will write $W^\perp$ for its orthogonal complement with respect to $\langle \,\,\,,\,\,\,\rangle_n$. Thus $\Lambda_{n,(i)}^\perp$ is the submodule of $\Lambda_n$ consisting of vectors whose last $i$ entries are $0$.
Also write 
\[ \Lambda_n^{(m)}=\Hom(\cO_F^m,\Z) \oplus \Lambda_n, \]
and set $V_n^{(m)}=\Lambda_n^{(m)} \otimes_\Z \Q$.

Define an additive group scheme $\Hom_n^{(m)}$ over $\Z$ by
\[ \Hom^{(m)}_n(R)=\Hom_{\cO_F}(\cO_F^m,\Lambda_n) \otimes_\Z R. \]
Then $\Hom_n^{(m)}$ has an action of $G_n \times \RS^{\cO_F}_\Z GL_m$ given by
\[ (g,h)f=g \circ f \circ h^{-1}. \]
Also define a perfect pairing
\[ \langle \,\,\,,\,\,\, \rangle_n^{(m)}: \Hom_n^{(m)}(R) \times \Hom_n^{(m)}(R) \lra R \]
by
\[ \langle f,f'\rangle_n^{(m)} = \sum_{i=1}^m \langle fe_i, f'e_i \rangle_n, \]
where $e_1,...,e_m$ denotes the standard basis of $\cO_F^m$. We have
\[ \langle (g,h)f,f'\rangle_n^{(m)}= \nu(g) \langle f, (g^{-1},{}^{c,t}h)f'\rangle_n^{(m)}. \]
Moreover $G_n(R)$ is identified with the set of pairs
\[ (g,\mu) \in GL(\Hom_{\cO_F}(\cO_F^m, \Lambda_n))(R) \times R^\times \]
such that $g$ commutes with the action of $GL_m(\cO_F \otimes_{\Z} R)$ and such that
\[ \langle gf,gf'\rangle_n^{(m)} = \mu \langle f,f'\rangle_n^{(m)} \]
for all $f,f' \in \Hom_{\cO_F}(\cO_F^m, \Lambda_n))(R)$. We set
\[ G_n^{(m)}=G_n \ltimes \Hom_n^{(m)}. \]
Then $G_n^{(m)}$ has an action of $\RS^{\cO_F}_\Z GL_m$ by
\[ h(g,f)=(g,(1,h)f). \]

Moreover $G_n^{(m)}$ acts on $\Lambda_n^{(m)}$, by letting $f \in \Hom_n^{(m)}$ act by
\[ f: (h,x) \longmapsto (h+\langle x, f \rangle_n, x) \]
and $g \in G_n$ act by
\[ g: (h,x) \longmapsto (h,gx). \]
Moreover $\RS_\Z^{\cO_F} GL_n$ acts on $\Lambda_n^{(m)}$ by
\[ \gamma: (h ,x) \longmapsto (h \circ \gamma^{-1}, x). \]
We have $\gamma \circ g = \gamma(g) \circ \gamma$.

If $m_1 \geq m_2$ we embed $\cO_F^{m_2} \into \cO_F^{m_1}$ via
\[ i_{m_2,m_1}:(x_1,...,x_{m_2}) \longmapsto (x_1,...,x_{m_2},0,...,0). \]
This gives rise to maps 
\[ i_{m_2,m_1}^*: \Hom_n^{(m_1)} \lra \Hom_n^{(m_2)} \]
and
\[ i_{m_2,m_1}^*:G_n^{(m_1)} \lra G_n^{(m_2)}. \]
It also gives rise to
\[ i_{m_2,m_1}^*: \Lambda_n^{(m_1)} \onto \Lambda_n^{(m_2)}. \]

Suppose that $R$ is a ring and that $X$ is an $\cO_F \otimes_\Z R$-module. We will write $\Herm_X$ for the $R$-module of $R$-bilinear pairings
\[ (\,\,\, ,\,\,\,): X \times X \lra R \]
which satisfy
\begin{enumerate}
\item\label{cond1} $(ax,y)=(x,{}^cay)$ for all $a \in \cO_F$ and $x,y \in X$;
\item \label{cond2} $(x,y)=(y,x)$ for all $x,y \in X$.
\end{enumerate}
If $z \in \Herm_X$ we will sometimes denote the corresponding pairing $(\,\,\,,\,\,\,)_z$. 
If $S$ is an $R$-algebra we have a natural map
\[ \Herm_X \otimes_R S \lra \Herm_{X \otimes_R S}. \]
If $X=\cO_F^m \otimes_\Z R$ then we will write 
\[ \Herm^{(m)}(R)=\Herm_{\cO_F^m} \otimes_\Z R \liso \Herm_{\cO_F^m \otimes_\Z R}. \]
If $X \ra Y$ then there is a natural map $\Herm_Y \ra \Herm_X$. In particular if $m_1 \geq m_2$, then there is a natural map 
\[ \Herm^{(m_1)} \lra \Herm^{(m_2)} \]
induced by the map $\cO_F^{(m_2)} \into \cO_F^{(m_1)}$ described in the last paragraph.
The group $GL(X/\cO_F)$ acts on the left on $\Herm_X$ by
\[ (x,y)_{h z}=(h^{-1}x,h^{-1}y)_z. \]
There is a natural isomorphism
\[ \Herm_{X \oplus Y} \cong \Herm_X  \oplus \Hom_R(X \otimes_{\cO_F \otimes R, c \otimes 1} Y, R)\oplus \Herm_Y, \]
under which an element $(z,f,w)$ of the right hand side corresponds to
\[ ((x,y),(x',y'))_{(z,f,w)}= (x,x')_z + f(x \otimes y')+f(x' \otimes y)+(y,y')_w. \]

Set $N_n^{(m)}(\Z)$ to be the set of pairs
\[ (f,z) \in \Hom_{\cO_F}(\cO_F^m,\Lambda_n) \oplus (\frac{1}{2}\Herm^{(m)}(\Z))  \]
such that 
\[ (x,y)_z - \frac{1}{2} \langle f x ,fy \rangle_n\in \Z \]
for all $x,y \in \cO_F^m$.
We define a group scheme $N_{n}^{(m)}/\Spec \Z$ by setting $N_{n}^{(m)}(R)$ to be the set of pairs
\[ (f,z) \in N_n^{(m)}(\Z) \otimes_\Z R \]
with group law given by
\[ (f,z) (f',z') = (f+f', z+z'+\frac{1}{2}( \langle f\,\,\,,f'\,\,\,\rangle_n - \langle f'\,\,\,,f\,\,\,\rangle_n )), \]
where by $\langle f\,\,\,,f'\,\,\,\rangle_n - \langle f'\,\,\,,f\,\,\,\rangle_n$ we mean the hermitian form
\[ (x,y) \longmapsto \langle f(x),f'(y)\rangle_n - \langle f'(x),f(y)\rangle_n. \]
Note that $(f,z)^{-1}=(-f,-z)$.
Thus there is an exact sequence
\[ (0) \lra \Herm^{(m)} \lra N_{n}^{(m)} \lra \Hom^{(m)}_n \lra (0). \]
In fact $Z(N_{n}^{(m)})=\Herm^{(m)}$. The commutator in $N_n^{(m)}$ induces an alternating 
map
\[  \Hom_n^{(m)}(R) \times \Hom_n^{(m)}(R) \lra \Herm^{(m)}(R) \]
under which $(f,f')$ maps to the pairing
\[ (x,y) \longmapsto  \langle f(x),f'(y)\rangle_n-\langle f'(x),f(y)\rangle_n.  \]
If $m_1 \geq m_2$ there is a natural map
\[ N_n^{(m_1)} \lra N_n^{(m_2)} \]
compatible with the previously described maps 
\[ \Hom_n^{(m_1)} \ra \Hom_n^{(m_2)} \] 
and 
\[ \Herm^{(m_1)} \ra \Herm^{(m_2)}. \]

Note that $G_n \times \RS^{\cO_F}_\Z GL_m$ acts on $N_{n}^{(m)}$ from the left by 
\[ (g,h)(f,z)=(g \circ f \circ h^{-1}, \nu(g) hz). \]
If $2$ is invertible in $R$ we see that
\[ \Herm^{(m)}(R)=\{ g\in N_n^{(m)}(R): \,\, (-1_m)(g)=g \} \]
and
\[ \Hom_n^{(m)}(R)=\{ g\in N_n^{(m)}(R): \,\, (-1_m)(g)=g^{-1} \}. \]
Set 
\[ \tG^{(m)}_n = G_n \ltimes N_n^{(m)},\]
which has an $ \RS^{\cO_F}_\Z GL_m$-action via
 \[ h(g,u)=(g,h(u)). \]
 If $m_1 \geq m_2$ then we get a natural map $\tG^{(m_1)}_n \ra \tG_n^{(m_2)}$.
 Note that
\[ G^{(m)}_n\cong \tG^{(m)}_n/\Herm^{(m)}.\]

Let $B_n$ denote the subgroup of $G_n$ consisting of elements which preserve the chain $\Lambda_{n,(n)} \supset \Lambda_{n,(n-1)} \supset ... \supset \Lambda_{n,(1)} \supset \Lambda_{n,(0)}$ and let $N_n$ denote the normal subgroup of $B_n$ consisting of elements with $\nu=1$, which also act trivially on $\Lambda_{n,(i)}/\Lambda_{n,(i-1)}$ for all $i=1,...,n$. Let $T_n$ denote the group consisting of the diagonal elements of $G_n$ and let $A_n$ denote the image of $\G_m$ in $G_n$ via the embedding that sends $t$ onto $t1_{2n}$. 
Over $\Q$ we see that $T_n$ is a maximal torus in a Borel subgroup $B_n$ of $G_n$, and that $N_n$ is the unipotent radical of $B_n$. Moreover $A_n$ is a maximal split torus in the centre of $G_n$.

If $\Omega$ is an algebraically closed field of characteristic $0$ then set 
\[ X^*(T_{n,/\Omega})=\Hom (T_n \times \Spec \Omega, \G_m \times \Spec \Omega).\] 
Also let $\Phi_n \subset X^*(T_{n,/\Omega})$ denote the set of roots of $T_n$ on $\Lie G_n$; let
$\Phi^+_n
\subset \Phi_n$ denote the set of positive roots with respect to $B_n$ and let
$\Delta_n \subset \Phi_n^+$ denote the set of simple positive roots. We will write $\varrho_n^+$ for half the sum of the elements of $\Phi_n^+$.
If $R \subset \R$ is a subring then $X^*(T_{n,/\Omega})^+_R$ will denote the subset of $X^*(T_{n,/\Omega})_R$ consisting of elements which pair non-negatively with the coroot $\check{\alpha} \in X_*(T_{n,/\Omega})$ corresponding to each $\alpha \in \Delta_n$. We will write simply $X^*(T_{n,/\Omega})^+$ for $X^*(T_{n,/\Omega})_\Z^+$. If $\lambda \in X^*(T_{n,/\Omega})^+$  we will let $\rho_{n,\lambda}$ (or simply $\rho_\lambda$) denote the irreducible
representation of $G_n$ with highest weight $\lambda$.
{\em When $\rho_\lambda$ is used as a subscript we
will sometimes replace it by just $\lambda$.} 

There is a natural identification 
\[ G_n \times \Spec \Omega \cong \left\{ (\mu,g_\tau) \in \G_m \times GL_{2n}^{\Hom(F,\Omega)}: \,\, g_{\tau c} =\mu J_n {}^tg_\tau^{-1} J_n \,\, \forall \tau\right\} . \]
This gives rise to the further identification
\[ T_n \times \Spec \Omega \cong \left\{ (t_0,(t_{\tau,i})) \in \G_m \times (\G_m^{2n})^{\Hom(F,\Omega)}: \,\, t_{\tau,i}t_{\tau c, 2n+1-i}=t_0 \,\, \forall \tau, i\right\}. \]
We will use this to identify $X^*(T_{n,/\Omega})$ with a quotient of 
\[ X^*(\G_m \times (\G_m^{2n})^{\Hom(F,\Omega)}) \cong \Z \oplus (\Z^{2n})^{\Hom(F,\Omega)}. \]
Under this identification $X^*(T_{n,/\Omega})^+$ is identified to the image of the set of
\[ (a_0,(a_{\tau,i})) \in \Z \oplus (\Z^{2n})^{\Hom(F,\Omega)} \]
with
\[ a_{\tau,1} \geq a_{\tau,2} \geq ... \geq a_{\tau,2n} \]
for all $\tau$.

If $R$ is a subring of $\R$ and $H$ an algebraic subgroup of $\tG_n^{(m)}$ we will write $H(R)^+$ for the subgroup of $H(R)$ consisting of elements with positive multiplier. Thus $G_n(\R)^+$ (resp. $G_n^{(m)}(\R)^+$, resp. $\tG_n^{(m)}(\R)^+$) is the connected component of the identity in $G_n(\R)$ (resp. $G_n^{(m)}(\R)$, resp. $\tG_n^{(m)}(\R)$).

Let 
\[ U_{n,\infty} = (U(n)^2)^{\Hom(F^+,\R)} \rtimes \{ 1,j \} \]
with $j^2=1$ and $j(A_\tau,B_\tau)j=(B_\tau,A_\tau)$. Embed $U_{n,\infty}$ in $G_n(\R)$ by sending 
$(A_\tau,B_\tau) \in (U(n)^2)^{\Hom(F^+,\R)}$ to
\[ \begin{array}{r} \left( 1,\left( \mat{(A_\tau+B_\tau)/2}{(A_\tau-B_\tau)\Psi_n/2i}{(\Psi_n(B_\tau-A_\tau)/2i}{\Psi_n(A_\tau+B_\tau)\Psi_n/2} \right)_\tau \right) \\ \\ \in G_n(\R) \subset \R^\times \times \prod_{\tau \in \Hom(F^+,\R)} GL_{2n}(F \otimes_{F^+,\tau} \R), \end{array} \]
and sending $j$ to 
\[ \left( -1,\left(\mat{-1_n}{0}{0}{1_n}\right)_\tau \right). \]
(This map depends on identifications $F \otimes_{F^+,\tau} \R \cong \C$, but the image of the map does not, and this image is all that will concern us.)
Then $U_{n,\infty}$ is a maximal compact subgroup of $G_n(\R)$ (and even of $\tG_n^{(m)}(\R)$). If $L \supset T_n \times \Spec \R$ is a Levi component of a parabolic subgroup $P \supset B_n \times \Spec \R$ then $U_{n,\infty} \cap L(\R)$ is a maximal compact subgroup of $L(\R)$. The connected component of the identity of $U_{n,\infty}$ is $U_{n,\infty}^0=U_{n,\infty} \cap G_n(\R)^+$.

We will write $\gp_n$ for the set of elements of $\Lie G_n(\R)$ of the form
\[ \left(0, \left( \mat{A_\tau}{B_\tau \Psi_n}{\Psi_nB_\tau}{\Psi_nA_\tau\Psi_n} \right)_{\tau \in \Hom(F^+,\R)} \right), \]
where ${}^tA_\tau^c=A_\tau$ and ${}^tB_\tau^c=B_\tau$ for all $\tau$. Then
\[ \Lie G_n(\R) = \gp_n \oplus \Lie (U_{n,\infty} A_n(\R)). \]
We give the real vector space $\gp_n$ a complex structure by letting $i$ act by 
\[ i_0: (A_\tau,B_\tau)_{\tau \in \Hom(F^+,\R)} \longmapsto (B_\tau,-A_\tau)_{\tau \in \Hom(F^+,\R)} . \]
We decompose
\[ \gp_n \otimes_\R \C = \gp_n^+ \oplus \gp_n^- \]
by setting 
\[ \gp_n^{\pm} = (\gp_n \otimes_\R\C)^{i_0\otimes 1=\pm 1 \otimes i}. \]
We also set
\[ \gq_n= \gp_n^- \oplus \Lie (U_{n,\infty} A_n(\R)) \otimes_\R \C. \]
It is a parabolic sub-algebra of $(\Lie G_n(\R)) \otimes_\R \C$ with unipotent radical $\gp_n^-$ and Levi component $\Lie (U_{n,\infty} A_n(\R)) \otimes_\R \C$. We will write $\gQ_n$ for the parabolic subgroup of $G_n \times_\Q \C$ with Lie algebra $\gq_n$. Note that
\[ \gQ_n(\C) \cap G_n(\R) = U_{n,\infty}^0 A_n(\R). \]

Let $\gH_n^+$ (resp. $\gH_n^{\pm}$) denote the set of $I$ in $G_n(\R)$ with multiplier $1$ such that
$I^2=-1_{2n}$ and
such that the symmetric bilinear form $\langle I \,\,\, ,\,\,\, \rangle_n$ on $\Lambda_n
\otimes_{\Z} \R$ is positive definite (respectively positive or negative definite). Then $G_n(\R)$ (resp. $G_n(\R)^+$) acts transitively on
$\gH_n^\pm$ (resp. $\gH_n^+$) by conjugation. Moreover $J_n \in \gH_n^+$ has stabilizer
$U_{n,\infty}^0A_n(\R)^0$ and so we get an identification of $\gH_n^{\pm}$ (resp. $\gH_n^+$) with
$G_n(\R)/U_{n,\infty}^0A_n(\R)^0$ (resp. $G_n(\R)^+/U_{n,\infty}^0A_n(\R)^0$). The natural map
\[ \gH_n^\pm = G_n(\R)/U_{n,\infty}^0A_n(\R)^0 \into G_n(\C)/\gQ_n(\C) \]
is an open embedding and gives $\gH_n^\pm$ the structure of a complex manifold. The action of $G_n(\R)$ is holomorphic and the complex structure induced on the tangent space $T_{J} \gH_n^\pm \cong \gp_n$ is the complex structure described in the previous paragraph.

If $\rho$ is an algebraic representation of $\gQ_n$ on a $\C$-vector space $W_\rho$, then there is a holomorphic vector bundle $\gE_\rho/\gH_n^\pm$ together with a holomorphic action of $G_n(\R)$, defined as the pull back to $\gH^\pm$ of 
$(G_n(\C) \times W_\rho)/\gQ_n(\C)$, where
\begin{itemize}
\item $h \in \gQ_n(\C)$ sends $(g,w)$ to $(gh,h^{-1}w)$,
\item and where $h \in G_n(\R)$ sends $[(g,w)]$ to $[(hg,w)]$.
\end{itemize}

If $N_2 \geq N_1 \geq 0$ are integers we will write $U_p(N_1,N_2)_n$ for the subgroup of $G_n(\Z_p)$ consisting of elements whose reduction modulo $p^{N_2}$ preserves
\[ \Lambda_{n,(n)} \otimes_\Z (\Z/p^{N_2}\Z) \subset \Lambda_n \otimes  _\Z (\Z/p^{N_2}\Z) \]
and acts trivially on $\Lambda_n/(\Lambda_{n,(n)} + p^{N_1} \Lambda_n)$. If $N_1 \geq N_1' \geq 0$ then
$U_p(N_1,N_2)_n$ is a normal subgroup of $U_p(N_1',N_2)_n$ and
\[ U_p(N_1',N_2)_n/U_p(N_1,N_2)_n \cong \ker (GL_n(\cO_F/p^{N_1}) \ra  GL_n(\cO_F/p^{N_1'}) ). \]
We will also set
\[ \begin{array}{rcl} U_p(N_1,N_2)^{(m)}_n&=&U_p(N_1,N_2)_n \ltimes 
\Hom_{\cO_{F,p}}(\cO_{F,p}^m, \Lambda_{n,(n)} +p^{N_1}\Lambda_n)\\ & \subset &G_n^{(m)}(\Z_p) \end{array} \]
and set $\tU_p(N_1,N_2)^{(m)}_n$ to be the pre-image of $U_p(N_1,N_2)^{(m)}_n$ in $\tG_n^{(m)}(\Z_p)$.
If $U^p$ is an open compact subgroup of $G_n(\A^{p,\infty})$ (resp. $G_n^{(m)}(\A^{p,\infty})$, resp. $\tG_n^{(m)}(\A^{p,\infty})$) we will set $U^p(N_1,N_2)$ to be $U^p \times U_p(N_1,N_2)_n$ (resp. $U^p \times U_p(N_1,N_2)_n^{(m)}$, resp. $U^p \times \tU_p(N_1,N_2)_n^{(m)}$), a compact open subgroup of $G_n(\A^{\infty})$ (resp. $G_n^{(m)}(\A^{\infty})$, resp. $\tG_n^{(m)}(\A^{\infty})$).

\newpage \subsection{Maximal parabolic subgroups.}\label{mps}

We will write $P_{n,(i)}^+$ (resp. $P^{(m),+}_{n,(i)}$, resp. $\tP^{(m),+}_{n,(i)}$) for the subgroup of $G_n$ (resp. $G_n^{(m)}$, resp. $\tG_n^{(m)}$) consisting of elements which (after projection to $G_n$) take $\Lambda_{n,(i)}$ to itself. We will also write $N_{n,(i)}^+$ (resp. $N^{(m),+}_{n,(i)}$, resp. $\tN^{(m),+}_{n,(i)}$) for the subgroups of $P_{n,(i)}^+$ (resp. $P^{(m),+}_{n,(i)}$, resp. $\tP^{(m),+}_{n,(i)}$) consisting of elements which act trivially on $\Lambda_{n,(i)}$ and $\Lambda_{n,(i)}^\perp/\Lambda_{n,(i)}$ and $\Lambda_n/\Lambda_{n,(i)}^\perp$. Over $\Q$ the groups  $P_{n,(i)}^+$ (resp. $P^{(m),+}_{n,(i)}$, resp. $\tP^{(m),+}_{n,(i)}$) are maximal parabolic subgroups of $G_n$ (resp. $G_n^{(m)}$, resp. $\tG_n^{(m)}$)
containing the pre-image of $B_n$. The groups $N_{n,(i)}^+$ (resp. $N^{(m),+}_{n,(i)}$, resp. $\tN^{(m),+}_{n,(i)}$) are their unipotent radicals.

In some instances it will be useful to replace these groups by their `Hermitian part'. We will write $P_{n,(i)}$ for the normal subgroup of $P_{n,(i)}^+$  consisting of elements which act trivially on $\Lambda_n/\Lambda_{n,(i)}^\perp$. We will also write $P^{(m)}_{n,(i)}$ for the normal subgroup
\[ P_{n,(i)} \ltimes \Hom_{\cO_F}(\cO_F^m,\Lambda_{n,(i)}^\perp) \]
of $P^{(m),+}_{n,(i)}$, and $\tP^{(m)}_{n,(i)}$ for the pre-image of $P^{(m)}_{n,(i)}$ in $\tP^{(m),+}_{n,(i)}$. We will let
\[ N_{n,(i)}=N_{n,(i)}^+ \]
and
\[ N^{(m)}_{n,(i)}=N^{(m),+}_{n,(i)} \cap P^{(m)}_{n,(i)} \]
and
\[ \tN^{(m)}_{n,(i)}=\tN^{(m),+}_{n,(i)} \cap \tP^{(m)}_{n,(i)}. \]
Over $\Q$ these are the unipotent radicals of $P_{n,(i)}$ (resp. $P^{(m)}_{n,(i)}$, resp. $\tP^{(m)}_{n,(i)}$). 

We have an isomorphism
\[ P_{n,(i)} \cong \tG_{n-i}^{(i)}. \]
To describe it let $\Lambda_{n,(i)}'$ denote the subspace of $\Lambda_n$ consisting of vectors with their first $2n-i$ entries $0$, so that 
\[ \Lambda_{n,(i)}' \cong \cO_F^i \]
and 
\[ \Lambda_{n-i} \cong \Lambda_{n,(i)}^\perp \cap (\Lambda_{n,(i)}')^\perp \liso \Lambda_{n,(i)}^\perp/\Lambda_{n,(i)}. \]
We define 
\[ G_{n-i} \into P_{n,(i)} \]
by letting $g \in G_{n-i}$ act as $\nu(g)$ on $\Lambda_{n,(i)}$, as $g$ on $\Lambda_{n-i} \cong \Lambda_{n,(i)}^\perp \cap (\Lambda_{n,(i)}')^\perp$ and as $1$ on $\Lambda_{n,(i)}'$, i.e.
\[ g \longmapsto  \left( \begin{array}{ccc} \nu(g) 1_i & 0&0 \\ 0 & g & 0 \\ 0 & 0 & 1_i \end{array} \right) \in P_{n,(i)}. \]
We define 
\[ N_{n,(i)} \lra \Hom_{n-i}^{(i)} \]
by sending $h$ to the map
\[ \cO_F^i \cong \Lambda_{n,(i)}' \stackrel{h-1_{2n}}{\lra} \Lambda_{n,(i)}^\perp \onto \Lambda_{n-i}. \]
We also define 
\[ Z(N_{n,(i)}) \liso \Herm^{(i)} \]
by sending $z$ to the pairing 
\[ (x,y)_z=\langle (z-1_{2n})x,y\rangle_n \]
on $\Lambda_{n,(i)}'$. In the other direction $(f,z) \in N_{n-i}^{(i)}$ is mapped to
\[ \left( \begin{array}{ccc} 1_i & \Psi_{i} {}^tf^cJ_{n-i}& \Psi_i{}^t(z-\frac{1}{2} {}^tfJ_{n-i}f^c)\\ 0 & 1_{2(n-i)} & f \\ 0 & 0 & 1_i \end{array} \right) \in N_{n,(i)}, \]
where we think of $f \in M_{2(n-i) \times i}(F)$ with first $n-i$ rows in $(\cD_F^{-1})^i$ and second $(n-i)$ rows in $\cO_F^i$, and we think of $z \in M_{i \times i}(F)^{t=c}$.

We also have isomorphisms
\[ P_{n,(i)}^{(m)} \cong \tG_{n-i}^{(i+m)}/\Herm^{(m)} \]
and
\[ \tP_{n,(i)}^{(m)} \cong \tG_{n-i}^{(i+m)}. \]
We will describe the second of these isomorphisms. Suppose $f \in \Hom_{n-i}^{(i)}$ and $g \in \Hom_{n-i}^{(m)}$. Also suppose that $z \in \frac{1}{2}\Herm^{(i)}$ and $w \in  \frac{1}{2}\Herm^{(m)}$ and $u \in  \frac{1}{2}\Hom(\cO_F^i \otimes_{\cO_F,c} \cO_F^m,\Z)$, so that
\[ ((f,g),(z,u,w) ) \in N_{n-i}^{(i+m)}. \]
Let $h(f,z)$ denote the element of $P_{n,(i)}$ corresponding to $(f,z) \in N_{n,(i)}$. Think of $g$ as a map
\[ g: \cO_F^m \lra \Lambda_{n-i} \cong \Lambda_{n,(i)}^\perp \cap (\Lambda_{n,(i)}')^\perp \subset \Lambda_n. \]
Define $j(f,g,u) \in \Hom(\cO_F^m, \Lambda_{n,(i)})$ by
\[ \langle y, j(f,g,u)(x)\rangle_n = 1/2 \langle f(y),g(x)\rangle_{n-i} -u(y \otimes x) \]
for all $x \in \cO_F^m$ and $y \in \Lambda_{n,(i)}' \cong \cO_F^i$. Then
\[ ((f,g),(z,u,w)) \longmapsto h(f,z) (g+j(f,g,u),w) \in N_{n,(i)} \ltimes \tN_n^{(m)}. \]

Note that 
\[ Z(\tN^{(m)}_{n,(i)}) \cong \Herm_{i+m} \]
and that
\[ Z(N^{(m)}_{n,(i)}) \cong \Herm_{i+m}/\Herm_m. \]

Write $L_{n,(i),\lin}$ for the subgroup of $P_{n,(i)}^+$ consisting of elements with $\nu=1$ which preserve $\Lambda_{n,(i)}'\subset \Lambda_n$ and act trivially on $\Lambda_{n,(i)}^\perp/\Lambda_{n,(i)}$. We set $N(L_{n,(i),\lin}^{(m)})$ to be the additive group scheme over $\Z$ associated to
\[ \Hom_{\cO_F}(\cO_F^m, \Lambda_{n,(i)}'), \]
and write $L_{n,(i),\lin}^{(m)}$ for
\[ L_{n,(i),\lin} \ltimes N(L_{n,(i),\lin}^{(m)}) \subset P_{n,(i)}^{(m),+} \]
and $\tL_{n,(i),\lin}^{(m)}$ for
\[ L_{n,(i),\lin} \ltimes N(L_{n,(i),\lin}^{(m)}) \subset \tP_{n,(i)}^{(m),+}. \]
Note that
\[ P_{n,(i)}^+ = L_{n,(i),\lin} \ltimes P_{n,(i)} \]
and
\[ P_{n,(i)}^{(m),+} = L_{n,(i),\lin}^{(m)} \ltimes P_{n,(i)}^{(m)} \]
and
\[ \tP_{n,(i)}^{(m),+} = \tL_{n,(i),\lin}^{(m)} \ltimes \tP_{n,(i)}^{(m)}. \]
Also note that
\[ L_{n,(i),\lin} \cong \RS^{\cO_F}_\Z GL_i \]
via its action on $\Lambda_{n,(i)}'\cong \cO_F^i$,
and that
\[ \tL_{n,(i),\lin}^{(m)} \liso L_{n,(i),\lin}^{(m)} \cong (\RS^{\cO_F}_\Z GL_i) \ltimes \Hom_{\cO_F}(\cO_F^m,\cO_F^i). \]

We let $L_{n,(i),\herm}$ denote the subgroup of $P_{n,(i)}$ consisting of elements which preserve $\Lambda_{n,(i)}'$. Thus
\[ L_{n,(i),\herm} \cong G_{n-i}. \]
In particular
\[ \nu: L_{n,(n),\herm} \liso \G_m.\]
Over $\Q$ it is a Levi component for $P_{n,(i)}$ and $P^{(m)}_{n,(i)}$ and $\tP^{(m)}_{n,(i)}$, so in particular 
\[ P_{n,(i)} = L_{n,(i),\herm} \ltimes N_{n,(i)} \]
and
\[ P_{n,(i)}^{(m)} = L_{n,(i),\herm} \ltimes N_{n,(i)}^{(m)} \]
and
\[ \tP_{n,(i)}^{(m)} = L_{n,(i),\herm} \ltimes \tN_{n,(i)}^{(m)}. \]
We also set
\[ L_{n,(i)} = L_{n,(i),\herm} \times L_{n,(i),\lin} \]
and
\[ L_{n,(i)}^{(m)} = L_{n,(i),\herm} \times L_{n,(i),\lin}^{(m)}. \]
and
\[ \tL_{n,(i)}^{(m)} = L_{n,(i),\herm} \times \tL_{n,(i),\lin}^{(m)}. \]
Over $\Q$ we see that $L_{n,(i)}$ is a Levi component for each of $P_{n,(i)}^+$ and $P^{(m),+}_{n,(i)}$ and $\tP^{(m),+}_{n,(i)}$. Moreover
\[ P_{n,(i)}^+ = L_{n,(i)} \ltimes N_{n,(i)} \]
and
\[ P_{n,(i)}^{(m),+} = L_{n,(i)}^{(m)} \ltimes N_{n,(i)}^{(m)} \]
and
\[ \tP_{n,(i)}^{(m),+} = \tL_{n,(i)}^{(m)} \ltimes \tN_{n,(i)}^{(m)}. \]

We will occasionally write $P_{n,(i)}^{(m),-}$ (resp. $L_{n,(i),\herm}^-$) for the kernel of the map $P_{n,(i)}^{(m)}\ra C_{n-i}$ (resp. $L_{n,(i),\herm} \ra C_{n-i}$).  

We will write $R_{n,(n),(i)}$ for the subgroup of $L_{n,(n)}$ mapping $\Lambda_{n,(i)}'$ to itself. We will write $N(R_{n,(n),(i)})$ for the subgroup of $R_{n,(n),(i)}$ which acts trivially on $\Lambda_{n,(i)}'$ and $(\Lambda_{n,(i)}')^\perp/\Lambda_{n,(i)}'$ and $\Lambda_n/(\Lambda_{n,(i)}')^\perp$.

We will also write $R_{n,(n)}^{(m)}$ for the semi-direct product
\[ L_{n,(n)} \ltimes \Hom_{\cO_F}(\cO_F^{(m)},\Lambda_{n,(n)}). \]

If $m'\leq m$ we will fix $\Z^m \onto \Z^{m'}$ to be projection onto the last $m'$-coordinates and define $Q_{m,(m')}$ for the subgroup of $GL_m$ consisting of elements preserving the kernel of this map. We also define $Q_{m,(m')}'$ to be the subgroup of $Q_{m,(m')}$ consisting of elements which induce $1_{\Z^{m'}}$ on $\Z^{m'}$.
Thus there is an exact sequence
\[ (0) \lra \Hom(\Z^{m'},\Z^{m-m'}) \lra Q_{m,(m')}' \lra GL_{m-m'} \lra \{ 1\}. \]
Moroever
\[ \tL_{n,(i),\lin}^{(m)} \cong L_{n,(i),\lin}^{(m)} \cong \RS^{\cO_F}_\Z Q'_{m+i,(m)}. \]

We will also write $A_{n,(i),\lin}$ (resp. $A_{n,(i),\herm}$) for the image of the map from $\G_m$ to $L_{n,(i),\lin}$ (resp. $L_{n,(i),\herm}$) sending $t$ to $t 1_i$ (resp. $(t^2,t 1_{2(n-i)})$). Moreover write $A_{n,(i)}$ for $A_{n,(i),\lin} \times A_{n,(i),\herm}$. The group $A_{n,(i)}$ (resp. $A_{n,(i),\lin}$, resp. $A_{n,(i),\herm}$) is the maximal split torus in the centre of $L_{n,(i)}$  (resp. $L_{n,(i),\lin}$, resp. $L_{n,(i),\herm}$).

Again suppose that $\Omega$ is an algebraically closed field of characteristic $0$. 
Let $\Phi_{(n)} \subset \Phi_n$ denote the set of roots of $T_n$ on $\Lie L_{n,(n)}$, and set $\Phi_{(n)}^+=\Phi_n^+ \cap \Phi_{(n)}$ and $\Delta_{(n)}=\Delta_n \cap \Phi_{(n)}$. We will write $\varrho_{n,(n)}$ for half the sum of the elements of $\Phi_{(n)}^+$. If $R \subset \R$ then $X^*(T_{n,/\Omega})^+_{(n),R}$ will denote the subset of $X^*(T_{n,/\Omega})_R$ consisting of elements which pair non-negatively with the coroot $\check{\alpha} \in X_*(T_{n,/\Omega})$ corresponding to each $\alpha \in \Delta_{(n)}$. We write $X^*(T_{n,/\Omega})_{(n)}^+$ for $X^*(T_{n,/\Omega})_{(n),\Z}^+$. If $\lambda \in X^*(T_{n,/\Omega})_{(n)}^+$  we will let $\rho_{(n),\lambda}$ denote the irreducible
representation of $L_{n,(n)}$ with highest weight $\lambda$.
{\em When $\rho_{(n),\lambda}$ is used as a subscript we
will sometimes replace it by just $(n),\lambda$.} 

Note that $\Lie P_{n,(n)}(\C)$ and $\gq_n$ are conjugate under $G_n(\C)$ and hence we obtain an identification (`Cayley transform') of $(\Lie U_{n,\infty} A_n(\R)) \otimes_\R \C$ and $\Lie L_{n,(n)}(\C)$, which is well defined up to conjugation by $L_{n,(n)}(\C)$. Similarly $\gQ_n$ and $P_{n,(n)}(\C)$ are conjugate in $G_n \times_\Q \C$. Thus $L_{n,(n)}(\C)$ can be identified with $\gQ_n$ modulo its unipotent radical, canonically up to $L_{n,(n)}(\C)$-conjugation. Thus if $\rho$ is an algebraic representation of $L_{n,(n)}$ over $\C$, we can associate to it a representation of $\gQ_n$ and of $\gq_n$, and hence a holomorphic vector bundle $\gE_\rho/\gH_n^\pm$ with $G_n(\R)$-action.

The isomorphism $L_{n,(n)} \cong GL_1 \times \RS^{\cO_F}_\Z GL_n$ gives rise to a natural identification
\[ L_{n,(n)} \times \Spec \Omega \cong GL_1 \times GL_n^{\Hom(F,\Omega)}, \]
and hence to identifications
\[ T_n \times \Spec \Omega \cong GL_1 \times (GL_1^{n})^{\Hom(F,\Omega)} \]
and
\[ X^*(T_{n,/\Omega}) \cong \Z \oplus (\Z^{n})^{\Hom(F,\Omega)}. \]
Under this identification 
$X^*(T_{n,/\Omega})_{(n)}^+$ is identified to the set of
\[ (b_0,(b_{\tau,i})) \in \Z \oplus (\Z^{n})^{\Hom(F,\Omega)} \]
with
\[ b_{\tau,1} \geq b_{\tau,2} \geq ... \geq b_{\tau,n} \]
for all $\tau$. 

To compare this parametrization of $X^*(T_{n,/\Omega})$ with the one introduced in section \ref{s1.1} note that the map
\[ GL_1 \times GL_n^{\Hom(F,\Omega)} \into  \left\{ (\mu,g_\tau) \in \G_m \times GL_{2n}^{\Hom(F,\Omega)}: \,\, g_{\tau c} =\mu J_n {}^tg_\tau^{-1} J_n \,\, \forall \tau\right\} \]
coming from $L_{n,(n)} \into G_n$ sends
\[ (\mu,(g_\tau)_{\tau \in \Hom(F,\Omega)}) \longmapsto  \left( \mu, \left( \mat{\mu \Psi_n {}^t g_{\tau c}^{-1} \Psi_n}{0}{0}{g_\tau} \right)_{\tau \in \Hom(F,\Omega)} \right). \]
Thus the map
\[ \Z \oplus (\Z^{2n})^{\Hom(F,\Omega)} \onto X^*(T_{n,/\Omega}) \cong \Z \oplus (\Z^{n})^{\Hom(F,\Omega)} \]
sends
\[ (a_0, (a_{\tau,i})_{\tau \in \Hom(F,\Omega);\;\;\; i=1,...,2n}) \longmapsto \left(a_0+\sum_\tau\sum_{j=1}^n a_{\tau,j}, \,\, \left(a_{\tau,n+i}-a_{\tau c, n+1-i}\right)_{\tau,i}\right). \]
A section is provided by the map 
\[ (b_0,(b_{\tau,i})) \longmapsto (b_0,(0,...,0,b_{\tau,1},...,b_{\tau,n})_\tau). \]
In particular we see that $X^*(T_{n,/\Omega})^+\subset X^*(T_{n,/\Omega})_{(n)}^+$ is identified with the set of 
\[ (b_0,(b_{\tau,i})) \in \Z \oplus (\Z^{n})^{\Hom(F,\Omega)} \]
with
\[ b_{\tau,1} \geq b_{\tau,2} \geq ... \geq b_{\tau,n} \]
and
\[ b_{\tau,1}+b_{\tau c,1} \leq 0 \]
for all $\tau$.

Note that
\[ 2(\varrho_n-\varrho_{n,(n)})= (n^2[F^+:\Q],(-n)_{\tau,i}). \]

 We write $\Std$ for the representation of $L_{n,(n)}$ on $\Lambda_n/\Lambda_{n,(n)}$ over $\Z$, and if $\tau: F \into \barQQ$ we write $\Std_\tau$ for the representation of $L_{n,(n)}$ on $(\Lambda_n/\Lambda_{n,(n)}) \otimes_{\cO_F,\tau}\cO_{\barQQ}$. If $\Omega$ is an algebraically closed field of characteristic $0$ and if $\tau: F \into \Omega$ we will sometimes write $\Std_\tau$ for the representation of $L_{n,(n)}$ on $(\Lambda_n/\Lambda_{n,(n)}) \otimes_{\cO_F,\tau}\Omega$. We hope that context will make clear the distinction between these two slightly different meaning of $\Std_\tau$.
 We also let
$\KS$ denote the unique representation of $L_{n,(n)}$ over $\Z$ such that
\[ W_\KS = (W_{\nu} \otimes_\Q \bigoplus_{\tau\in \Hom(F,\barQQ)/\{1,c\}} W_{\Std_\tau}^\vee \otimes W_{\Std_{\tau \circ c}}^\vee)^{G_\Q}. \]
Note that over $\barQQ$ the representation $\Std_\tau^\vee$ is irreducible and in our normalizations has highest weight $(0,b_{\tau'})$ where 
 \[ b_\tau = (0,...,0,-1) \]
 but $b_{\tau'}=0$ for $\tau' \neq \tau$. Similarly the representation $\wedge^{n[F:\Q]} \Std^\vee$ is irreducible with highest weight 
 \[ (0,(-1,...,-1)_\tau). \]
 Finally $\KS$ is the direct sum of the $[F^+:\Q]$ irreducible representations indexed by $\tau\in \Hom(F^+, \barQQ)$ with highest weights $(1,b_{\tau'})$, where 
 \[ b_{\tau'} = (0,...,0,-1) \]
 if $\tau'$ extends $\tau$, and $b_{\tau'}=0$ otherwise. 

We will let $\varsigma_p \in L_{n,(n),\herm}(\Q_p)\cong \Q_p^\times$ denote the unique element with multiplier $p^{-1}$. 

Set 
\[ U_{p}(N)_{n,(i)}=\ker(L_{n,(i),\lin}(\Z_p) \ra L_{n,(i),\lin}(\Z/p^N\Z))\]
 and 
 \[ U_{p}(N)_{n,(i)}^{(m)}=\ker(L_{n,(i),\lin}^{(m)}(\Z_p) \ra L_{n,(i),\lin}^{(m)}(\Z/p^N\Z)).\]
 Also set
\[ U_p(N_1,N_2)_{n,(i)}^{(m)}=U_p(N_1,N_2)_{n-i} \times U_p(N_1)_{n,(i)}^{(m)} \subset L_{n,(i)}^{(m)}(\Z_p) \]
and
\[ \tU_p(N_1,N_2)_{n,(i)}^{(m),+}=
U_p(N_1)^{(m)}_{n,(i)} \ltimes \tU_p(N_1,N_2)_{n-i}^{(m+i)} \subset \tP^{(m),+}_{n,(i)}(\Z_p) \]
and
\[ U_p(N_1,N_2)_{n,(i)}^{(m),+}=\tU_p(N_1,N_2)_{n,(i)}^{(m),+} /\Herm_{\cO_{F,p}^m} \subset P^{(m),+}_{n,(i)}(\Z_p). \]
If $m=0$ we will drop it from the notation. 
If $U^p$ is an open compact subgroup of $L_{n,(i)}(\A^{p,\infty})$ (resp. $(P_{n,(i)}^{(m),+}/Z(N_{n,(i)}^{(m)}))(\A^{p,\infty})$, resp. $P_{n,(i)}^{(m),+}(\A^{p,\infty})$, resp. $\tP_{n,(i)}^{(m),+}(\A^{p,\infty})$) then set 
\[ U^p(N_1,N_2) = U^p \times U_p(N_1,N_2)_{(i)}^{(m)} \subset L_{n,(i)}(\A^\infty) \]
(resp.
\[ U^p(N_1,N_2) = U^p \times (U_p(N_1,N_2)_{(i)}^{(m),+}/Z(N_{n,(i)}^{(m)})(\Z_p)) \subset P_{n,(i)}^{(m),+}/Z(N^{(m)}_{n,(i)})(\A^\infty),\]
resp. 
\[ U^p(N_1,N_2) =U^p \times U_p(N_1,N_2)_{(i)}^{(m),+} \subset P_{n,(i)}^{(m),+}(\A^\infty),\]
resp.
\[ U^p(N_1,N_2) =U^p \times \tU_p(N_1,N_2)^{(m),+}_{n,(i)} \subset \tP_{n,(i)}^{(m),+}(\A^\infty)). \]
In the case $i=n$ these groups do not depend on $N_2$, so we will write simply $U^p(N_1)$.

For the study of the ordinary locus we will need a variant of $G_n(\A^\infty)$ and $G_n^{(m)}(\A^\infty)$ and $\tG_n^{(m)}(\A^\infty)$. More specifically define a semigroup
\[ \tG_n^{(m)}(\A^\infty)^\ord = \tG^{(m)}_n(\A^{p,\infty}) \times (\varsigma_p^{\Z_{\geq 0}}\tP^{(m),+}_{n,(n)}(\Z_p)). \]
Its maximal sub-semigroup that is also a group is 
\[ \tG_n^{(m)}(\A^\infty)^{\ord,\times} = \tG^{(m)}_n(\A^{p,\infty}) \times \tP^{(m),+}_{n,(n)}(\Z_p). \]
If $H$ is an algebraic subgroup of $\tG_n^{(m)}$ (over $\Spec \Q$) we set
 \[ H(\A^\infty)^\ord = H(\A^{\infty}) \cap \tG_n^{(m)}(\A^\infty)^\ord. \]
 Its maximal sub-semigroup that is also a group is
 \[ H(\A^\infty)^{\ord,\times} = H(\A^{\infty}) \cap (\tG^{(m)}_n(\A^{p,\infty})\times \tP^{(m),+}_{n,(n)}(\Z_p)). \]
 Thus
 \[ G_n(\A^\infty)^{\ord,\times} = G_n(\A^{p,\infty}) \times P_{n,(n)}^+(\Z_p) \]
and
\[ G_n^{(m)}(\A^\infty)^{\ord,\times} = G^{(m)}_n(\A^{p,\infty}) \times P^{(m),+}_{n,(n)}(\Z_p). \]
If $U^p$ is an open compact subgroup of $H(\A^{p,\infty})$, we set  
 \[ U^p(N)=H(\A^\infty)^{\ord,\times} \cap (U^p \times \tU_p(N,N')_n^{(m)}) \]
 for any $N' \geq N$. The group does not depend on the choice of $N'$.  
\newpage

\subsection{Base change.}\label{secbc}

We will write $B_{GL_m}$ for the subgroup of upper triangular elements of $GL_m$ and $T_{GL_m}$ for the subgroup of diagonal elements of $B_{GL_m}$. 

We will also let $G_n^1$ denote the group scheme over $\cO_{F^+}$ defined by
\[ G_n^1(R) = \{ g \in \Aut((\Lambda_n \otimes_{\cO_{F^+}} R)/(\cO_F \otimes_{\cO_{F^+}} R)) :\,\,
{}^tgJ_n {}^cg = J_n\}. \]
Thus
\[ \ker \nu \cong \RS^{\cO_{F^+}}_\Z G_n^1. \]
We will write $B_n^1$ for the subgroup of $G_n^1$ consisting of upper triangular matrices and $T_n^1$ for the subgroup of $B_n^1$ consisting of diagonal matrices. There is a natural projection $B_n^1 \onto T_n^1$ obtained by setting the off diagonal entries of an element of $B_n^1$ to $0$. 

Suppose that $q$ is a rational prime. Let $u_1,...,u_r$ denote the primes of $F^+$ above $\Q$ which split $u_i=w_i{}^cw_i$ in $F$ and let $v_1,...,v_s$ denote the primes of $F^+$ above $q$ which do not split in $F$. Then
\[ G_n(\Q_q) \cong \prod_{i=1}^r GL_{2n}(F_{w_i}) \times H \]
where 
\[ H= \left\{ (\mu,g_i) \in \Q_q^\times \times \prod_{i=1}^s GL_{2n}(F_{v_i}): \,\, {}^tg_i J_n{}^cg_i=\mu J_n \,\, \forall i\right\} \supset \prod_{i=1}^s G_n^1(F^+_{v_i}). \] 
Suppose that $\Pi$ is an irreducible smooth representation of $G_n(\Q_q)$ then
\[ \Pi = \left( \bigotimes_{i=1}^r \Pi_{w_i} \right) \otimes \Pi_H. \]
We define $\BC(\Pi)_{w_i}=\Pi_{w_i}$ and $\BC(\Pi)_{c w_i}= \Pi_{w_i}^{c,\vee}$. Note that this does not depend on the choice of primes $w_i|u_i$. 
We will say that $\Pi$ is {\em unramified} at $v_i$ if $v_i$ is unramified in $F$ and 
\[ \Pi^{G_n^1(\cO_{F^+,v_i})} \neq (0). \]
If $\Pi$ is unramified at $v_i$ then there is a character $\chi$ of $T_n^1(F^+_{v_i})/T_n^1(\cO_{F^+,v_i})$ such that $\Pi|_{G_n^1(F^+_{v_i})}$ and $\nind_{B_n^1(F^+_{v_i})}^{G_n^1(F^+_{v_i})} \chi$ share an irreducible sub-quotient with a $G_n^1(\cO_{F^+,v_i})$-fixed vector. Moreover this character $\chi$ is unique modulo the action of the normalizer $N_{G_n^1(F_{v_i}^+)}(T_n^1(F_{v_i}^+))/T_n^1(F_{v_i}^+)$. (If $\pi$ and $\pi'$ are two irreducible subquotients of $\Pi_H|_{G_n^1(F_{v_i}^+)}$ then we must have $\pi' \cong \pi^{\varsigma_{v_i}^{-1}}$ where 
\[ \varsigma_{v_i} = \mat{\varpi_{v_i}^{-1} 1_n}{0}{0}{1_n} \in GL_{2n}(F_{v_i}). \]
However 
\[ \left. \left( \nind_{B_n^1(F^+_{v_i})}^{G_n^1(F^+_{v_i})} \chi \right)^{\varsigma_v^{-1}} \cong \nind_{B_n^1(F^+_{v_i})}^{G_n^1(F^+_{v_i})} \chi .\right) \]
Let 
\[ \begin{array}{rcl} \norm: T_{GL_{2n}}(F_{v_i}) &\lra &T_n^1(F_{v_i}^+) \\ \diag(t_1,...,t_{2n})& \longmapsto & \diag(t_1/{}^ct_{2n},...,t_{2n}/{}^ct_1). \end{array} \]
Then we define $\BC(\Pi)_{v_i}$ to be the unique subquotient of
\[ \nind_{B_{GL_{2n}}(F_{v_i})}^{GL_{2n}(F_{v_i})} \chi \circ \norm \]
with a $GL_{2n}(\cO_{F,v_i})$-fixed vector.
The next lemma is easy to prove.

\begin{lem}\label{indbc} Suppose that $\psi \otimes \pi$ is an irreducible smooth representation of 
\[ L_{(n)}(\Q_q) \cong L_{(n),\herm}(\Q_q) \times L_{(n),\lin}(\Q_q) = \Q_q^\times \times GL_n(F_q). \]
\begin{enumerate} 
\item If $v$ is unramified over $F^+$ and $\pi_v$ is unramified then $\nind_{P_{(n)}(\Q_q)}^{G_n(\Q_q)} (\psi \otimes \pi)$ has a subquotient $\Pi$ which is unramified at $v$. Moreover $\BC(\Pi)_v$ is the unramified irreducible subquotient of $\nind_{Q_{2n,(n)}(F_v)}^{GL_{2n}(F_v)} (\pi^{c,\vee}_v \otimes \pi_v)$.
\item If $v$ is split over $F^+$ and $\Pi$ is an irreducible sub-quotient of the normalized induction $\nind_{P_{(n)}(\Q_q)}^{G_n(\Q_q)} (\psi \otimes \pi)$, then $\BC(\Pi)_v$ is an irreducible subquotient of $\nind_{Q_{2n,(n)}(F_v)}^{GL_{2n}(F_v)} ((\pi_{{}^cv})^{c,\vee} \otimes \pi_v)$.
\end{enumerate} 
Note that in both cases $\BC(\Pi_v)$ does not depend on $\psi$. \end{lem} 

In this paragraph let $K$ be a number field, $m\in \Z_{>0}$, and write $U_{K,\infty}$ for a maximal compact subgroup of $GL_m(K_\infty)$. We shall (slightly abusively) refer to an admissible 
\[ G_n(\A^\infty) \times ((\Lie G_n(\R))_\C, U_{n,\infty})\]
 (resp.
\[ L_{n,(i)}(\A^\infty) \times ((\Lie L_{n,(i)}(\R))_\C,U_{n,\infty}\cap L_{n,(i)}(\R)),\]
 resp. 
 \[ GL_m(\A^\infty_K) \times ((\Lie GL_m(K_\infty))_\C,
U_{K,\infty} ))\]
 module as an admissible $G_n(\A)$-module (resp.
$L_{n,(i)}(\A)$-module, resp. $GL_m(\A_K)$-module). By a square-integrable automorphic representation of 
$G_n(\A)$ (resp. $L_{n,(i)}(\A)$, resp. $GL_m(\A_K)$) we shall mean the twist by a character of an irreducible admissible $G_n(\A)$-module (resp. $L_{n,(i)}(\A)$-module, resp. $GL_m(\A_K)$-module) that occurs discretely in the space of square integrable automorphic forms on the double coset space
$G_n(\Q) \backslash G_n(\A)/A_n(\R)^0$ (resp. $L_{n,(i)}(\Q)\backslash L_{n,(i)}(\A)/A_{n,(i)}(\R)^0$, resp. $GL_m(K)\backslash GL_m(\A_K)/\R^\times_{>0}$). By a cuspidal automorphic representation of 
$G_n(\A)$ (resp. $L_{n,(i)}(\A)$, resp. $GL_m(\A_K)$) we shall mean an irreducible admissible $G_n(\A)$-sub-module (resp. $L_{n,(i)}(\A)$-sub-module, resp. $GL_m(\A_K)$-sub-module) of the space of cuspidal automorphic forms on $G_n(\A)$ (resp. $L_{n,(i)}(\A)$, resp. $GL_m(\A_K)$).

\begin{prop} Suppose that $\Pi$ is a square integrable automorphic representation of $G_n(\A)$ and that $\Pi_\infty$ is cohomological. Then there is an expression
\[ 2n=m_1n_1+...+m_rn_r \]
with $m_i,n_i \in \Z_{>0}$ and cuspidal automorphic representations $\tPi_i$ of $GL_{m_i}(\A_F)$ such that
\begin{itemize}
\item $\tPi_i^\vee \cong \tPi_i^c$;
\item $\tPi_i ||\det||^{(m_i+n_i-1)/2}$ is cohomological;
\item if $v$ is a prime of $F$ above a rational prime $q$ such that
\begin{itemize}
\item either $q$ splits in $F_0$,
\item or $F$ and $\Pi$ are unramified above $q$,
\end{itemize}
 then  
\[ \BC(\Pi_q)_v = \boxplus_{i=1}^r \boxplus_{j=0}^{n_i-1} \tPi_{i,v} |\det|_v^{(n_i-1)/2-j}. \]
\end{itemize} \end{prop}

\pfbegin This follows from the main theorem of \cite{sws} and the classification of square integrable automorphic representations of $GL_m(\A_F)$ in \cite{mw}. \pfend

\begin{cor}\label{galois1} Keep the assumptions of the proposition. Then there is a continuous, semi-simple, algebraic (i.e. unramified almost everywhere and de Rham above $p$) representation 
\[ r_{p,\imath}(\Pi): G_F \lra GL_{2n}(\barQQ_p) \]
with the following property: If $v$ is a prime of $F$ above a rational prime $q$ such that
\begin{itemize}
\item either $q$ splits in $F_0$,
\item or $F$ and $\Pi$ are unramified above $q$,
\end{itemize}
then 
\[ \imath \WD(r_{p,\imath}(\Pi)|_{G_{F_v}})^\semis \cong \rec_{F_v}(\BC(\Pi_q)_v|\det|_v^{(1-2n)/2}). \]
\end{cor}

\pfbegin Combine the proposition with for instance theorem 1.2 of \cite{blght} and theorem A of \cite{blggtcomp2}. (These results are due to many people and we simply choose these particular references for convenience.)
\pfend

\newpage \subsection{Spaces of Hermitian forms.}\label{herm}

There is a natural pairing
\[ \begin{array}{rcl} (X \otimes_{\cO_F\otimes R,c \otimes 1} X) \times \Herm_{X} & \lra & R \\ (x \otimes y, z)& \longmapsto & (x,y)_z. \end{array} \]
We further define
\[ \begin{array}{rcl}  \sw: (X \otimes_{\cO_F \otimes R,c \otimes 1} X)  & \lra & (X \otimes_{\cO_F \otimes R,c \otimes 1} X) \\ x \otimes y & \longmapsto & y \otimes x, \end{array} \]
and
\[ S(X)=(X \otimes_{\cO_F\otimes R,c \otimes 1} X)/(\sw -1). \]
There is a natural map in the other direction
\[ \begin{array}{rcl} S(X) & \lra & X \otimes_{\cO_F \otimes R,c \otimes 1} X \\ w & \longmapsto & w+\sw(w), \end{array} \]
such that the composite $S(X) \ra X \otimes_{\cO_F \otimes R,c \otimes 1} X \ra S(X)$ is multiplication by $2$. 
Note that if $F/F^+$ is ramified above $2$ then $S(\cO_F^m)$ can have $2$-torsion, but that $S(\cO_{F,(p)}^m)$ is torsion free. (Either $p>2$ or by assumption $F/F^+$ is not ramified above $2$.) There is a perfect duality
\[ S(\cO_F^m)^\TF \times \Herm^{(m)}(\Z) \lra \Z. \]
We will write
\[ e = \sum_{i=1}^m e_i \otimes e_i \in \cO_F^m \otimes_{\cO_F,c} \cO_F^m, \]
where $e_1,...,e_m$ denotes the standard basis of $\cO_F^m$.

If $R \subset \R$ then we will denote by $\Herm_X^{>0}$ (resp. $\Herm_X^{\geq 0}$) the set of pairings $(\,\,\,,\,\,\,)$ in $\Herm_X$ such that 
\[ (x,x) > 0 \]
(resp. $\geq 0$) for all $x \in X-\{0\}$. We will denote by $S(F^m)^{>0}$ the set of elements $a \in S(F^m)$ such that for each $\tau:F \into \C$ the image of $a$ under the map
\[ \begin{array}{rcl} S(F^m) & \lra & M_m(F)^{t=c} \\ x \otimes y & \longmapsto & x{}^ty^c+y{}^tx^c \end{array} \]
is positive definite, i.e. all the roots of its characteristic polynomial are strictly positive real numbers. Then $S(F^m)^{>0}$ is the set of elements of $S(F^m)$ whose pairing with every element of $\Herm_{F^m}^{>0}$ is strictly positive; and $\Herm_{F^m}^{>0}$ is the set of elements of  $\Herm_{F^m}$ whose pairing with every element of $S(F^m)^{>0}$ is strictly positive.

Suppose that $W \subset V_n$  is an isotropic $F$-direct summand. We set 
\[ \gC^{(m)}(W)= (\Herm_{V_n/W^\perp} \oplus \Hom_F(F^m,W)) \otimes_\Q \R. \]
If $m=0$ we will drop it from the notation. There is a natural map 
\[ \gC^{(m)}(W) \lra \gC(W). \]
Note that if $f \in \Hom_F(F^m,W)$ we can define $f' \in \Hom(F^m \otimes_{F,c} (V_n/W^\perp),\Q)$ by
\[ f'(x \otimes y) = \langle f(x),y\rangle_n. \]
This establishes an isomorphism
\[ \Hom_F(F^m,W) \liso \Hom(F^m \otimes_{F,c} (V_n/W^\perp),\Q) \]
and hence an isomorphism
\[ \gC^{(m)}(W) \liso (\Herm_{V_n/W^\perp \oplus F^m}/\Herm_{F^m}) \otimes_\Q\R. \] 
Thus 
\[ \gC^{(m)}(V_{n,(i)}) \cong Z(N_{n,(i)}^{(m)})(\R). \]
If $g \in G_n(\Q)$ we define
\[ \begin{array}{rcl} g: \gC^{(m)}(W) &\lra &\gC^{(m)}(gW) \\ (z,f) & \longmapsto & (gz,g \circ f), \end{array} \]
where
\[ (x,y)_{gz}=|\nu(g)| (g^{-1}x,g^{-1}y)_z. \]
We extend this to an action of $G_n^{(m)}(\Q)$ as follows:
 If $g \in \Hom_F(F^m,V_n)$ then we set
\[ g(z,f)=(z,f - \theta_z \circ g) \]
where $\theta_z: V_n \ra W$ satisfies
\[ (x \bmod W^\perp , y \bmod W^\perp)_z = \langle \theta_z(x),y \rangle_n \]
for all $x,y \in V_n$. If $W' \subset W$ there is a natural embedding
\[ \gC^{(m)}(W') \into \gC^{(m)}(W). \]

We will write $\gC^{>0}(W)=\Herm_{V_n/W^\perp \otimes_\Q \R}^{>0}$ and $\gC^{\geq 0}(W)=\Herm_{V_n/W^\perp\otimes_\Q \R}^{\geq 0}$. We will also write $\gC^{(m),>0}(W)$ (resp. $\gC^{(m),\geq 0}(W)$) for the pre-image of $\gC^{>0}(W)$ (resp. $\gC^{\geq 0}(W)$) in $\gC^{(m)}(W)$. Moreover we will set 
\[ \gC^{(m),\succ 0}(W) = \bigcup_{W' \subset W} \gC^{(m),>0}(W'). \]
Thus
\[ \gC^{(m),>0}(W) \subset \gC^{(m),\succ 0}(W) \subset \gC^{(m),\geq 0}(W). \]
Note that the natural map $\gC^{(m)}(W) \onto \gC(W)$ gives rise to a surjection
\[ \gC^{(m),\succ 0}(W) \onto \gC^{\succ 0}(W) \]
and that the pre image of $(0)$ is $(0)$. Also note that if $W' \subset W$ then there is a closed embedding
\[ \gC^{(m),\succ 0}(W') \into \gC^{(m),\succ 0}(W). \]
Finally note that the action of $G_n^{(m)}(\Q)$ takes $\gC^{(m),\succ 0}(W)$ (resp. $\gC^{(m),> 0}(W)$, resp. $\gC^{(m),\geq 0}(W)$) to $\gC^{(m),\succ 0}(gW)$ (resp. $\gC^{(m),> 0}(gW)$, resp. $\gC^{(m),\geq 0}(gW)$).

Note that $L_{n,(i)}^{(m)}(\R)$ acts on 
\[ \pi_0(L_{n,(i),\herm}(\R)) \times \gC^{(m)}(V_{n,(i)})\]
and preserves 
\[ \pi_0(L_{n,(i),\herm}(\R)) \times \gC^{(m),>0}(V_{n,(i)}).\]
 Moreover $L_{n,(i)}^{(m)}(\Q)$ preserves 
 \[ \pi_0(L_{n,(i),\herm}(\R)) \times \gC^{(m),\succ 0}(V_{n,(i)}). \]
   In fact
$L_{n,(i)}^{(m)}(\R)$ acts transitively on $\pi_0(L_{n,(i),\herm}(\R)) \times \gC^{(m),>0}(V_{n,(i)})$. For this paragraph let $(\,\,\, ,\,\,\,)_0\in \gC^{>0}(V_{n,(i)})$ denote the pairing on  $(V_n/V_{n,(i)}^\perp \otimes_\Q \R)^2$ induced by $\langle J_n\,\,\, ,\,\,\,\rangle_n$. Then the stabilizer of $1 \times ((\,\,\, ,\,\,\,)_0,0)$ in $L_{n,(i)}^{(m)}(\R)$ is 
\[ L_{n,(i),\herm}(\R)^{\nu =1} (U_{n,\infty} \cap L_{n,(i), \lin}^{(m)}(\R)) A_n(\R)^0.\] Thus we get an $L_{n,(i)}^{(m)}(\R)$-equivariant identification
\[ \begin{array}{l} \pi_0(L_{n,(i),\herm}(\R)) \times \gC^{(m),>0}(V_{n,(i)})/\R^\times_{>0} \cong \\ L_{n,(i)}^{(m)}(\R) / L_{n,(i),\herm}(\R)^+ (U_{n,\infty} \cap L_{n,(i),\lin}^{(m)}(\R))^0 A_{n,(i)}(\R)^0. \end{array} \]

We define $\gC^{(m)}$ to be the topological space
\[ \left(\bigcup_W \gC^{(m),\succ 0}(W) \right)  / \sim, \]
where $\sim$ is the equivalence relation generated by the identifications of $\gC^{(m),\succ 0}(W')$ with its image in $\gC^{(m),\succ 0}(W)$ whenever $W' \subset W$. (This is sometimes referred to as the `conical complex'.) Thus as a set
\[ \gC^{(m)} = \coprod_W \gC^{(m),>0}(W). \]
We will let $\gC^{(m)}_{=i}$ denote
\[ \coprod_{\dim_F W=i} \gC^{(m),>0}(W). \]
Note that $\gC^{(m)}_{=n}$ is a dense open subset of $\gC^{(m)}$. 

The space $\gC^{(m)}$ has a natural, continuous, left action of $G^{(m)}_n(\Q) \times \R^\times_{>0}$. (The second factor acts on each $\gC^{(m),\succ 0}(W)$ by scalar multiplication.)

We have homeomorphisms
\[  \begin{array}{ll} & G_n^{(m)}(\Q) \backslash (G_n^{(m)}(\A^\infty)/U \times \pi_0(G_n(\R)) \times \gC^{(m)}_{=i}/\R_{>0}^\times) \\ \cong & P_{n,(i)}^{(m),+}(\Q) \backslash \left( G_n^{(m)}(\A^\infty)/U \times \pi_0(G_n(\R)) \times (\gC^{(m),>0}(V_{n,(i)})/\R^\times_{>0}) \right) \\ \cong & \coprod_{h \in P_{n,(i)}^{(m),+}(\A^\infty) \backslash G_n^{(m)}(\A^\infty)/U} L_{n,(i)}^{(m)}(\Q) \backslash L_{n,(i)}^{(m)}(\A) / \\ & (hUh^{-1} \cap P_{n,(i)}^{(m),+}(\A^\infty))L_{n,(i),\herm}(\R)^+ (L_{n,(i),\lin}^{(m)}(\R) \cap U_{n,\infty}^0)A_{n,(i)}(\R)^0 . \end{array}\]
(Use the fact, strong approximation for unipotent groups, that 
\[ N_{n,(i)}^{(m),+}(\A^\infty)=V+N_{n,(i)}^{(m),+}(\Q) \]
for any open compact subgroup $V$ of $N_{n,(i)}^{(m),+}(\A^\infty)$.)
If $g \in G_n^{(m)}(\A^\infty)$ and if $g^{-1}Ug \subset U'$ then the right translation map
\[ \begin{array}{l} g: G_n^{(m)}(\Q) \backslash (G_n^{(m)}(\A^\infty)/U \times \pi_0(G_n(\R)) \times \gC^{(m)}_{=i}/\R_{>0}^\times) \lra \\ G_n^{(m)}(\Q) \backslash (G_n^{(m)}(\A^\infty)/U' \times \pi_0(G_n(\R)) \times \gC^{(m)}_{=i}/\R_{>0}^\times) \end{array} \]
corresponds to the coproduct of the right translation maps
\[ \begin{array}{l} g': L_{n,(i)}^{(m)}(\Q) \backslash L_{n,(i)}^{(m)}(\A) / \\ (hUh^{-1} \cap P_{n,(i)}^{(m),+}(\A^\infty))L_{n,(i),\herm}(\R)^+ (L_{n,(i),\lin}^{(m)}(\R) \cap U_{n,\infty}^0)A_{n,(i)}(\R)^0 \\ \lra  \\
L_{n,(i)}^{(m)}(\Q) \backslash L_{n,(i)}^{(m)}(\A) / \\(h'U'(h')^{-1} \cap P_{n,(i)}^{(m),+}(\A^\infty))L_{n,(i),\herm}(\R)^+ (L_{n,(i),\lin}^{(m)}(\R) \cap U_{n,\infty}^0)A_{n,(i)}(\R)^0 \end{array} \]
where $hg=g'h'u'$ with $g' \in P_{n,(i)}^{(m),+}(\A^\infty)$ and $u' \in U'$.

We set
\[ (G_n^{(m)}(\A^\infty) \times \pi_0(G_n(\R)) \times \gC^{(m)})^\ord \]
to be the subset of $G_n^{(m)}(\A^\infty) \times \pi_0(G_n(\R)) \times \gC^{(m)}$ consisting of elements $(g,\delta,x)$ such that for some $W$ we have
\[ x \in \gC^{(m),\succ 0}(W) \]
and
\[ W \otimes_\Q \Q_p = g_p (V_{n,(n)} \otimes_\Q \Q_p). \]
It has a left action of $G_n^{(m)}(\Q)$ and a right action of $G_n^{(m)}(\A^\infty)^\ord \times \R^\times_{>0}$. 
We define
\[ (G_n^{(m)}(\A^\infty) \times \pi_0(G_n(\R)) \times \gC^{(m)}_{=i})^\ord \]
similarly. We also set 
\[ G_n^{(m)}(\Q) \backslash (G_n^{(m)}(\A^\infty)/U^p(N_1,N_2) \times \pi_0(G_n(\R)) \times \gC^{(m)})^\ord \]
(resp.
\[ G_n^{(m)}(\Q) \backslash (G_n^{(m)}(\A^\infty)/U^p(N_1,N_2) \times \pi_0(G_n(\R)) \times \gC^{(m)}_{=i})^\ord) \]
to be the image of $(G_n^{(m)}(\A^\infty) \times \pi_0(G_n(\R)) \times \gC^{(m)})^\ord$ in 
\[ G_n^{(m)}(\Q) \backslash (G_n^{(m)}(\A^\infty)/U^p(N_1,N_2) \times \pi_0(G_n(\R)) \times \gC^{(m)})\]
 (resp. 
 \[ G_n^{(m)}(\Q) \backslash (G_n^{(m)}(\A^\infty)/U^p(N_1,N_2) \times \pi_0(G_n(\R)) \times \gC^{(m)}_{=i})).\]
  Then as a set
\[ \begin{array}{l} G_n^{(m)}(\Q) \backslash (G_n^{(m)}(\A^\infty)/U^p(N_1,N_2) \times \pi_0(G_n(\R)) \times \gC^{(m)})^\ord = \\ \coprod_i G_n^{(m)}(\Q) \backslash (G_n^{(m)}(\A^\infty)/U^p(N_1,N_2) \times \pi_0(G_n(\R)) \times \gC^{(m)}_{=i})^\ord. \end{array} \]

\begin{lem} 
\[ \begin{array}{l} G_n^{(m)}(\Q) \backslash (G_n^{(m)}(\A^\infty) \times \pi_0(G_n(\R)) \times \gC^{(m)}_{=n})^\ord /U^p(N_1) \liso \\ G_n^{(m)}(\Q) \backslash (G_n^{(m)}(\A^\infty)/U^p(N_1,N_2) \times \pi_0(G_n(\R)) \times \gC^{(m)}_{=n})^\ord .\end{array} \]
\end{lem}

\pfbegin There is a natural surjection. We must check that it is also injective. The right hand side equals
\[ \begin{array}{l} 
P_{n,(n)}^{(m),+}(\Q) \backslash (G_n^{(m)}(\A^{p,\infty})/U^p \times (P_{n,(n)}^{(m),+}(\Q_p) U_p(N_1,N_2)_n^{(m)})/U_p(N_1,N_2)_n^{(m)} \times \\ \pi_0(G_n(\R)) \times \gC^{(m),>0}(V_{n,(n)})) \cong \\
P_{n,(n)}^{(m),+}(\Q) \backslash (G_n^{(m)}(\A^{\infty})^\ord/U^p(N_1)  \times \pi_0(G_n(\R)) \times \gC^{(m),>0}(V_{n,(n)})), \end{array} \]
which is clearly isomorphic to the left hand side.
\pfend

We set 
\[ \gT_{U^p(N_1),=n}^{(m),\ord} := G_n^{(m)}(\Q) \backslash (G_n^{(m)}(\A^\infty)/U^p(N_1,N_2) \times \pi_0(G_n(\R)) \times (\gC^{(m)}_{=n}/\R_{>0}^\times))^\ord. \]

There does not seem to be such a simple description of 
\[ G_n^{(m)}(\Q) \backslash (G_n^{(m)}(\A^\infty)/U^p(N_1,N_2) \times \pi_0(G_n(\R)) \times \gC^{(m)}_{=i})^\ord \]
 for $i \neq n$. However we do have the lemma below.

\begin{lem}\label{orddeco} There is a natural homeomorphism
\[ \begin{array}{ll} & G_n^{(m)}(\Q) \backslash (G_n^{(m)}(\A^\infty)/U^p(N_1,N_2) \times \pi_0(G_n(\R)) \times \gC^{(m)}_{=i})^\ord \\ \cong &

 \coprod_{h \in P_{n,(i)}^{(m),+}(\A^\infty)^{\ord,\times} \backslash G_n^{(m)}(\A^\infty)^{\ord,\times}/U^p(N_1)} L_{n,(i)}^{(m)}(\Q) \backslash L^{(m)}_{n,(i)}(\A)/  \\ 
 \multicolumn{2}{l}{
 \!\! ( (hU^p(N_1)h^{-1}\! \cap\! P^{(m),+}_{n,(i)}(\A^{\infty})^{\ord,\times}) L_{n,(i),\herm}^-(\Z_p)
 L_{n,(i),\herm}(\R)^0 (L_{n,(i),\lin}^{(m)}(\R) \! \cap \! U_{n,\infty}^0) )}

\end{array} \]
where $U^p(N_1) \subset G_n^{(m)}(\A^\infty)^{\ord,\times}$.

In particular 
\[ G_n^{(m)}(\Q) \backslash (G_n^{(m)}(\A^\infty)/U^p(N_1,N_2) \times \pi_0(G_n(\R)) \times \gC^{(m)}_{=i})^\ord \]
and
\[ G_n^{(m)}(\Q) \backslash (G_n^{(m)}(\A^\infty)/U^p(N_1,N_2) \times \pi_0(G_n(\R)) \times \gC^{(m)})^\ord  \]
 are independent of $N_2 \geq N_1$. 
\end{lem}

\pfbegin
Firstly we have that 
\[ \begin{array}{ll} & G_n^{(m)}(\Q) \backslash (G_n^{(m)}(\A^\infty)/U^p(N_1,N_2) \times \pi_0(G_n(\R)) \times \gC^{(m)}_{=i})^\ord \\ \cong &

P_{n,(i)}^{(m),+}(\Q) \backslash (G_n^{(m)}(\A^{p,\infty})/U^p \times \\ & (P_{n,(i)}^{(m),+}(\Q) P_{n,(n)}^{(m),+}(\Q_p) U_p(N_1,N_2)_n^{(m)})/U_p(N_1,N_2)_n^{(m)} \times \\ &\pi_0(G_n(\R)) \times \gC^{(m),>0}(V_{n,(i)})). \end{array} \]
We can replace the second $P_{n,(i)}^{(m),+}(\Q)$ by $P_{n,(i)}^{(m),+}(\Q_p)$, and then, using in particular the Iwasawa decomposition for $L_{n,(n)}(\Q_p)$, replace $P_{n,(n)}^{(m),+}(\Q_p)$ by $P_{n,(n)}^{(m),+}(\Z_p)$. Next we can replace $P_{n,(i)}^{(m),+}(\Q_p)$ by $P_{n,(i)}^{(m),+}(\Z_p)$ as long as we also replace $P_{n,(i)}^{(m),+}(\Q)$ by $P_{n,(i)}^{(m),+}(\Z_{(p)})$. This gives
\[ \begin{array}{ll} & G_n^{(m)}(\Q) \backslash (G_n^{(m)}(\A^\infty)/U^p(N_1,N_2) \times \pi_0(G_n(\R)) \times \gC^{(m)}_{=i})^\ord \\ \cong &

P_{n,(i)}^{(m),+}(\Z_{(p)}) \backslash (G_n^{(m)}(\A^{p,\infty})/U^p \times \\ & (P_{n,(i)}^{(m),+}(\Z_p) P_{n,(n)}^{(m),+}(\Z_p) U_p(N_1,N_2)_n^{(m)})/U_p(N_1,N_2)_n^{(m)} \times \\ & \pi_0(G_n(\R)) \times \gC^{(m),>0}(V_{n,(i)})). \end{array} \]

Note that
\[ P_{n-i,(n-i)}^{+}(\Z_p) \onto C_{n-i}(\Z_p). \]
[This follows from the fact that primes above $p$ on $F^+$ are unramified in $F$, which implies that
\[ \ker (\norm_{F/F^+}: \cO_{F,p}^\times \ra \cO_{F^+,p}^\times)=\{ x^{c-1}: \,\, x \in \cO_{F,p}^\times \}. ] \]
Thus
\[ L_{n,(i),\herm}^-(\Z_p) P_{n-i,(n-i)}^+(\Z_p) = L_{n,(i),\herm}(\Z_p) \]
and
\[ P_{n,(i)}^{(m),+}(\Z_p) P_{n,(n)}^{(m),+}(\Z_p) = P_{n,(i)}^{(m),-}(\Z_p) P_{n,(n)}^{(m),+}(\Z_p). \]
Moreover, by strong approximation, $P_{n,(i)}^{(m),-}(\Z_{(p)})$ (resp. $L_{n,(i),\herm}^-(\Z_{(p)})$) is dense in $P_{n,(i)}^{(m),-}(\A^{p,\infty}\times \Z_{p})$ (resp. $L_{n,(i),\herm}^-(\A^{p,\infty}\times \Z_{p})$). Thus 
\[ \begin{array}{ll} & G_n^{(m)}(\Q) \backslash (G_n^{(m)}(\A^\infty)/U^p(N_1,N_2) \times \pi_0(G_n(\R)) \times \gC^{(m)}_{=i})^\ord \\ \cong &

P_{n,(i)}^{(m),+}(\Z_{(p)}) \backslash (G_n^{(m)}(\A^{p,\infty})/U^p \times \\ & (P_{n,(i)}^{(m),-}(\Z_p) P_{n,(n)}^{(m),+}(\Z_p) U_p(N_1,N_2)_n^{(m)})/U_p(N_1,N_2)_n^{(m)} \times \\ & \pi_0(G_n(\R)) \times \gC^{(m),>0}(V_{n,(i)})) \\ \cong & 

L_{n,(i)}^{(m)}(\Z_{(p)}) \backslash ((P_{n,(i)}^{(m), -}(\A^{p,\infty})\backslash G_n^{(m)}(\A^{p,\infty})/U^p) \times \\ & (P_{n,(i)}^{(m),-}(\Z_{p}) \backslash (P_{n,(i)}^{(m),-}(\Z_p) P_{n,(n)}^{(m),+}(\Z_p) U_p(N_1,N_2)_n^{(m)})/U_p(N_1,N_2)_n^{(m)}) \times \\ & \pi_0(G_n(\R)) \times \gC^{(m),>0}(V_{n,(i)})). \end{array} \]

Next we claim that the natural map
\[ \begin{array}{rl} & (P_{n,(i)}^{(m),-} \cap P_{n,(n)}^{(m),+}) (\Z_{p}) \backslash P_{n,(n)}^{(m),+}(\Z_p) /(U_p(N_1)_{n,(n)}^{(m)}\Z_p^\times N_{n,(n)}^{(m)}(\Z_p)) \\
\lra & P_{n,(i)}^{(m),-}(\Z_{p}) \backslash (P_{n,(i)}^{(m),-}(\Z_p) P_{n,(n)}^{(m),+}(\Z_p) U_p(N_1,N_2)_n^{(m)})/U_p(N_1,N_2)_n^{(m)} \end{array} \]
is an isomorphism. It suffices to check this modulo $p^{N_2}$, where the map becomes
\[ \begin{array}{l}  (P_{n,(i)}^{(m),-} \cap P_{n,(n)}^{(m),+}) (\Z/p^{N_2}\Z) \backslash P_{n,(n)}^{(m),+}(\Z/p^{N_2}\Z) /(U_p(N_1)_{n,(n)}^{(m)}\Z_p^\times N_{n,(n)}^{(m)}(\Z_p)) \lra \\
 P_{n,(i)}^{(m),-}(\Z/p^{N_2}\Z) \backslash (P_{n,(i)}^{(m),-}(\Z/p^{N_2}\Z) P_{n,(n)}^{(m),+}(\Z/p^{N_2}\Z) /(U_p(N_1)_{n,(n)}^{(m)}\Z_p^\times N_{n,(n)}^{(m)}(\Z_p)), \end{array} \]
which is clearly an isomorphism. Thus we have
\[ \begin{array}{ll} & G_n^{(m)}(\Q) \backslash (G_n^{(m)}(\A^\infty)/U^p(N_1,N_2) \times \pi_0(G_n(\R)) \times \gC^{(m)}_{=i})^\ord \\ \cong &

L_{n,(i)}^{(m)}(\Z_{(p)}) \backslash ((P_{n,(i)}^{(m), -}(\A^{p,\infty})\backslash G_n^{(m)}(\A^{p,\infty})/U^p) \times \\ & 
((P_{n,(i)}^{(m),-} \cap P_{n,(n)}^{(m),+}) (\Z_{p}) \backslash P_{n,(n)}^{(m),+}(\Z_p) /(U_p(N_1)_{n,(n)}^{(m)}\Z_p^\times N_{n,(n)}^{(m)}(\Z_p)))
\times \\ & \pi_0(G_n(\R)) \times \gC^{(m),>0}(V_{n,(i)})),  

\end{array} \]
where $\gamma \in L_{n,(i)}^{(m)}(\Z_{(p)})$ acts on $(P_{n,(i)}^{(m),-} \cap P_{n,(n)}^{(m),+}) (\Z_{p}) \backslash P_{n,(n)}^{(m),+}(\Z_p)$ via an element of $P_{n-i,(n-i)}^+(\Z_p) \times L_{n,(i),\lin}(\Z_p)$ with the same image in $C_{n-i}(\Z_p) \times L_{n,(i),\lin}(\Z_p)$. 

Note that
\[ \begin{array}{rcl} P_{n,(i)}^{(m), -}(\A^{p,\infty})\backslash G_n^{(m)}(\A^{p,\infty})/U^p &=& \coprod_{h \in P_{n,(i)}^{(m),+}(\A^{p,\infty})\backslash G_n^{(m)}(\A^{p,\infty})/U^p}\\ \multicolumn{3}{r}{ L_{n,(i),\herm}^{(m), -}(\A^{p,\infty})\backslash L_{n,(i)}^{(m)}(\A^{p,\infty})/(hU^ph^{-1} \cap P_{n,(i)}^{(m),+}(\A^{p,\infty})).} \end{array} \]
Also note that, if we set $U_p=(U_p(N_1)_{n,(n)}^{(m)}\Z_p^\times N_{n,(n)}^{(m)}(\Z_p))$, then 
\[ \begin{array}{rcl} (P_{n,(i)}^{(m),-} \cap P_{n,(n)}^{(m),+}) (\Z_{p}) \backslash P_{n,(n)}^{(m),+}(\Z_p) /U_p &=& \coprod_{h \in (P_{n,(i)}^{(m),+} \cap P_{n,(n)}^{(m),+}) (\Z_{p}) \backslash P_{n,(n)}^{(m),+}(\Z_p) /U_p} \\ \multicolumn{3}{r}{( L_{n,(i),\lin}^{(m)}(\Z_p) \times \Im(P_{n-i,(n-i)}(\Z_p) \ra C_{n-i}(\Z_p)))/(hU_ph^{-1} \cap P_{n,(i)}^{(m),+}(\Z_p)).} \end{array}\]
However as the primes above $p$ split in $F^+$ split in $F$ we see that
\[ \Im(P_{n-i,(n-i)}(\Z_p) \ra C_{n-i}(\Z_p)) = L_{n,(i),\herm}(\Z_p)/L_{n,(i),\herm}^-(\Z_p), \]
and so
\[ \begin{array}{rcl} (P_{n,(i)}^{(m),-} \cap P_{n,(n)}^{(m),+}) (\Z_{p}) \backslash P_{n,(n)}^{(m),+}(\Z_p) /U_p& =& \coprod_{h \in (P_{n,(i)}^{(m),+} \cap P_{n,(n)}^{(m),+}) (\Z_{p}) \backslash P_{n,(n)}^{(m),+}(\Z_p) /U_p} \\ \multicolumn{3}{r}{L_{n,(i)}^{(m)}(\Z_p) /L_{n,(i),\herm}^-(\Z_p)(hU_ph^{-1} \cap P_{n,(i)}^{(m),+}(\Z_p)).} \end{array} \]

Thus we see that
\[ \begin{array}{ll} & G_n^{(m)}(\Q) \backslash (G_n^{(m)}(\A^\infty)/U^p(N_1,N_2) \times \pi_0(G_n(\R)) \times \gC^{(m)}_{=i})^\ord \\ \cong &

\coprod_{h \in P_{n,(i)}^{(m),+}(\A^\infty)^{\ord,\times} \backslash G_n^{(m)}(\A^\infty)^{\ord,\times}/U^p(N_1)} L_{n,(i)}^{(m)}(\Z_{(p)}) \backslash \\ &
\left( L^{(m)}_{n,(i)}(\A^{p,\infty} \times \Z_p)/ L_{n,(i),\herm}^-(\A^{p,\infty} \times \Z_p)
 (hU^p(N_1)h^{-1} \cap P^{(m),+}_{n,(i)}(\A^{\infty})^{\ord,\times})  \right. \\ & 
\left.  \times \pi_0(G_n^{(m)}(\R)) \times \gC^{(m),>0}(V_{n,(i)})\right). \end{array} \]
As $L_{n,(i),\herm}^{-}(\Z_{(p)})$ acts trivially on 
\[ (L^{(m)}_{n,(i)}(\Z_p)/ L_{n,(i),\herm}^-(\Z_p)) \times \pi_0(G_n^{(m)}(\R)) \times \gC^{(m),>0}(V_{n,(i)}) \]
and is dense in $L_{n,(i),\herm}^{-}(\A^{p,\infty})$, we further see that
\[ \begin{array}{ll} & G_n^{(m)}(\Q) \backslash (G_n^{(m)}(\A^\infty)/U^p(N_1,N_2) \times \pi_0(G_n(\R)) \times \gC^{(m)}_{=i})^\ord \\ \cong &

\coprod_{h \in P_{n,(i)}^{(m),+}(\A^\infty)^{\ord,\times} \backslash G_n^{(m)}(\A^\infty)^{\ord,\times}/U^p(N_1)} L_{n,(i)}^{(m)}(\Z_{(p)}) \backslash \\ &
\left( L^{(m)}_{n,(i)}(\A^{p,\infty} \times \Z_p)/ L_{n,(i),\herm}^-(\Z_p)
 (hU^p(N_1)h^{-1} \cap P^{(m),+}_{n,(i)}(\A^{\infty})^{\ord,\times})  \right. \\ & 
\left.  \times \pi_0(G_n^{(m)}(\R)) \times \gC^{(m),>0}(V_{n,(i)})\right) \\ \cong &
 
\coprod_{h \in P_{n,(i)}^{(m),+}(\A^\infty)^{\ord,\times} \backslash G_n^{(m)}(\A^\infty)^{\ord,\times}/U^p(N_1)} L_{n,(i)}^{(m)}(\Z_{(p)}) \backslash L^{(m)}_{n,(i)}(\A^{p} \times \Z_p)/ \\ &
\left(  
 (hU^p(N_1)h^{-1} \cap P^{(m),+}_{n,(i)}(\A^{\infty})^{\ord,\times}) L_{n,(i),\herm}^-(\Z_p) L_{n,(i),\herm}(\R)^0 \right. \\ \multicolumn{2}{r}{ \left. (L_{n,(i),\lin}^{(m)}(\R)  \cap U_{n,\infty}^0) \right), } \\ \cong &
 
\coprod_{h \in P_{n,(i)}^{(m),+}(\A^\infty)^{\ord,\times} \backslash G_n^{(m)}(\A^\infty)^{\ord,\times}/U^p(N_1)} L_{n,(i)}^{(m)}(\Q) \backslash L^{(m)}_{n,(i)}(\A)/  \\ 
\multicolumn{2}{l}{ \!\! ( (hU^p(N_1)h^{-1} \!\cap\! P^{(m),+}_{n,(i)}(\A^{\infty})^{\ord,\times}\!) L_{n,(i),\herm}^-(\Z_p)
 L_{n,(i),\herm}(\R)^0 (L_{n,(i),\lin}^{(m)}(\R) \!\cap \! U_{n,\infty}^0) ), }

 \end{array} \]
as desired.
\pfend

Abusing notation slightly, we will write
\[ G_n^{(m)}(\Q) \backslash (G_n^{(m)}(\A^\infty)/U^p(N_1) \times \pi_0(G_n(\R)) \times \gC^{(m)}_{=i})^\ord \]
for 
\[ G_n^{(m)}(\Q) \backslash (G_n^{(m)}(\A^\infty)/U^p(N_1,N_2) \times \pi_0(G_n(\R)) \times \gC^{(m)}_{=i})^\ord, \]
and
\[ G_n^{(m)}(\Q) \backslash (G_n^{(m)}(\A^\infty)/U^p(N_1) \times \pi_0(G_n(\R)) \times \gC^{(m)})^\ord \]
for
\[ G_n^{(m)}(\Q) \backslash (G_n^{(m)}(\A^\infty)/U^p(N_1,N_2) \times \pi_0(G_n(\R)) \times \gC^{(m)})^\ord. \]

If $U$ is a neat open compact subgroup of $L_{n,(i)}^{(m)}(\A^\infty)$, set 
\[ \gT^{(m)}_{(i),U} = L_{n,(i)}^{(m)}(\Q) \backslash L_{n,(i)}^{(m)}(\A) / U L_{n,(i),\herm}(\R)^0 (L_{n,(i),\lin}^{(m)}(\R) \cap U_{n,\infty}^0)A_{n,(n)}(\R)^0. \]

\begin{cor} 
\[  \gT^{(m),\ord}_{U^p(N_1),=n}  \cong 
\coprod_{h \in P_{n,(n)}^{(m),+}(\A^\infty)^\ord \backslash G_n^{(m)}(\A^\infty)^\ord/U^p(N_1)} 
\gT_{(n),(hU^ph^{-1} \cap P_{n,(n)}^{(m),+}(\A^{p,\infty})) \Z_p^\times U_p(N_1)^{(m)}_{n,(n)}}^{(m)}. \]
\end{cor}

If $Y$ is a locally compact, Hausdorff topological space then we write $H^i_\Int(Y,\C)$ for the image of 
\[ H^i_c(Y,\C) \lra H^i(Y,\C). \]

We define 
\[  H^i_\Int(\gT_{=n}^{(m),\ord},\barQQ_p) = \lim_{\ra U^p,N}H^i_\Int(\gT^{(m),\ord}_{U^p(N),=n},\barQQ_p) \]
a smooth $G^{(m)}_n(\A^\infty)^\ord$-module, and 
\[ H^i_\Int(\gT_{(n)}^{(m)},\barQQ_p) = \lim_{\ra U}H^i_\Int(\gT^{(m)}_{(n),U},\barQQ_p) \]
a smooth $L_{n,(n)}^{(m)}(\A^\infty)$-module. Note that
\[ H^i_\Int(\gT_{(n)}^{(m)},\barQQ_p)^{\Z_p^\times} = \lim_{\ra U^p,N}H^i_\Int(\gT^{(m)}_{(n),U^pU_p(N)^{(m)}_{n,(n)} \Z_p^\times},\barQQ_p) \]
as $N$ runs over positive integers and $U^p$ runs over neat open compact subgroups of $L_{n,(n)}^{(m)}(\A^{p,\infty})$. With these definitions we have the following corollary.

\begin{cor}\label{indu} There is a $G_n^{(m)}(\A^\infty)^\ord$-equivariant isomorphism
\[ \Ind^{G_n^{(m)}(\A^{p,\infty})}_{P_{n,(n)}^{(m),+}(\A^{p,\infty})} H^i_\Int(\gT_{(n)}^{(m)},\barQQ_p)^{\Z_p^\times} \cong H^i_\Int(\gT^{(m),\ord}_{=n},\barQQ_p). \] \end{cor}

Interior cohomology has the following property which will be key for us.

\begin{lem}\label{intcoh} Suppose that $G$ is a locally compact, totally disconnected topologyical group. Suppose that for any sufficiently small open compact subgroup $U \subset G$ we are given a compact Hausdorff space $Z_U$ and an open subset $Y_U \subset Z_U$. Suppose moreover that whenever $U$, $U'$ are sufficiently small open compact subgroups of $G$ and $g \in G$ with $g^{-1} U g \subset U'$, then there is a proper continuous map
\[ g: Z_U \lra Z_{U'} \]
with $g Y_U \subset Y_{U'}$. Also suppose that $g \circ h = hg$ whenever these maps are all defined and that, if $g\in U$ then the map $g: Z_U \ra Z_U$ is the identity. 

If $\Omega$ is a field, set
\[ H^i(Z,\Omega) = \lim_{\ra U} H^i(Z_U,\Omega) \]
and
\[ H^i_\Int(Y,\Omega) = \lim_{\ra U} H^i_\Int(Y_U,\Omega). \]
These are both smooth $G$-modules. Moreover $H^i_\Int(Y,\Omega)$ is a sub-quotient of $H^i(Z,\Omega)$ as $G$-modules.
\end{lem}

\pfbegin
Note that the diagram
\[ \begin{array}{rcl} H^i_c(Y_U,\Omega) & \lra & H^i(Y_U,\Omega) \\ \da && \ua \\
H^i_c(Z_U,\Omega) & = & H^i(Z_U,\Omega) \end{array} \]
is commutative. 
Set
\[ A = \lim_{\ra U} \Im\left( H^i_c(Y_U,\Omega) \lra H^i_c(Z_U,\Omega)=H^i(Z_U,\Omega)\right) \]
and
\[ B= \lim_{\ra U} \Im\left( \ker\left( H^i_c(Y_U,\Omega) \lra H^i(Y_U,\Omega)\right) \lra H^i(Z_U,\Omega)\right). \]
Then 
\[ B \subset A \subset H^i(Z,\Omega) \]
are $G$-invariant subspaces with
\[ A/B \liso H^i_\Int(Y,\Omega). \]
\pfend

\newpage \subsection{Locally symmetric spaces.}

In this section we will calculate $H^i_\Int(\gT_{(n)}^{(m)},\barQQ_p)$ in terms of automorphic forms on $L_{n,(n)}(\A)$. 

If $m=0$ we will write $\gT_{(n)}$ for $\gT^{(0)}_{(n)}$. Let $\Omega$ denote an algebraically closed field of characteristic $0$. 
If $\rho$ is a finite dimensional algebraic representation of $L_{n,(n)}$ on a $\Omega$-vector space $W_\rho$ then we define a locally constant sheaf $\cL_{\rho,U}/\gT_{(n),U}$ as
\[ \begin{array}{c} L_{n,(n)}(\Q) \backslash \left(W_\rho \times L_{n,(n)}(\A) /U (L_{n,(n)}(\R) \cap U_{n,\infty}^0)A_{n,(n)}(\R)^0\right) \\ \da \\
 L_{n,(n)}(\Q) \backslash L_{n,(n)}(\A)/U (L_{n,(n)}(\R) \cap U_{n,\infty}^0)A_{n,(n)}(\R)^0. \end{array} \]
The system of sheaves $\cL_{\rho,U}$ has a right action of $L_{n,(n)}(\A^\infty)$. We define
\[ H^i_\Int(\gT_{(n)},\cL_\rho)= \lim_{\ra U} H^i_\Int(\gT_{(n),U},\cL_{\rho,U}), \]
a smooth $L_{n,(n)}(\A^\infty)$-module. Note that if $\rho$ has a central character $\chi_\rho$ then,
\[ \alpha \in Z(L_{n,(n)})(\Q)^+ \subset L_{n,(n)}(\A^\infty) \]
acts on $H^i_\Int(\gT_{(n)},\cL_\rho)$ via $\chi_\rho(\alpha)^{-1}$. (Use the fact that $Z(L_{n,(n)})(\Q)^+ \subset 
(L_{n,(n)}(\R) \cap U_{n,\infty}^0)A_{n,(n)}(\R)^0$.)

The natural map $L_{n,(n)}^{(m)} \ra L_{n,(n)}$ gives rise to continuous maps
\[ \pi^{(m)}:\gT^{(m)}_{(n),U} \lra \gT_{(n), U} \]
compatible with the action of $L_{n,(n)}^{(m)}(\A^\infty)$.

\begin{lem}\label{pf} \begin{enumerate}
\item The maps $\pi^{(m)}$ are real-torus bundles (i.e. $(S^1)^r$-bundles for some $r$), and in particular are proper maps. 

\item There are $L^{(m)}_{n,(n)}(\A^\infty)$-equivariant identifications
\[ R^i\pi_{*}^{(m)} \Omega \cong \cL_{\wedge^i \left(\bigoplus_{\tau:F \into \Omega} \Std_\tau^{\oplus m}\right)^\vee}. \]
In particular the action of  $L^{(m)}_{n,(n)}(\A^\infty)$ on the relative cohomology sheaf $R^i\pi_{*}^{(m)} \Omega$ factors through $L_{n,(n)}(\A^\infty)$.
\end{enumerate}\end{lem}

\pfbegin
Recall that
\[ N(L_{n,(n),\lin}^{(m)})= \ker(L_{n,(n)}^{(m)} \ra L_{n,(n)}). \]
Suppose that $U$ is a neat open compact subgroup of $L^{(m)}_{n,(n)}(\A^\infty)$ with image $U'$ in $L_{n,(n)}(\A^\infty)$. Then $L_{n,(n)}(\Q) \times U'$ acts freely on
\[ L_{n,(n)}(\A)/ (U_{n,\infty}^0 \cap L_{n,(n)}(\R)) A_{n,(n)}(\R)^0. \]
Thus it suffices to prove that the map $\tpi^{(m)}$
\[ \begin{array}{c} N(L_{n,(n),\lin}^{(m)})(\Q) \backslash L^{(m)}_{n,(n)}(\A)/(U \cap N(L_{n,(n),\lin}^{(m)})(\A^\infty)) (U_{n,\infty}^0 \cap L_{n,(n)}(\R)) A_{n,(n)}(\R)^0 \\ \da \\ L_{n,(n)}(\A)/(U_{n,\infty}^0 \cap L_{n,(n)}(\R)) A_{n,(n)}(\R)^0 \end{array} \]
is a real torus bundle and that there are $L_{n,(n)}(\Q) \times L^{(m)}_{n,(n)}(\A^\infty) $-equivariant isomorphisms
\[ R^i\tpi_{*}^{(m)} \Omega \cong \cL_{\wedge^i \left( \bigoplus_\tau \Std_\tau^{\oplus m}\right)^\vee}. \]

Using the identification of spaces (but not of groups) that comes from the group product
\[ L_{n,(n)}^{(m)}(\A)=N(L_{n,(n),\lin}^{(m)})(\A) \times L_{n,(n)}(\A), \]
we see that $\tpi^{(m)}$ can be identified with the map
\[ \begin{array}{c} \left( N(L_{n,(n),\lin}^{(m)})(\Q) \backslash N(L_{n,(n),\lin}^{(m)})(\A) / (U \cap N(L_{n,(n),\lin}^{(m)})(\A^\infty))    \right) \\ \times \left( L_{n,(n)}(\A)/ (U_{n,\infty}^0 \cap L_{n,(n)}(\R)) A_{n,(n)}(\R)^0\right)  \\ \da \\ L_{n,(n)}(\A)/(U_{n,\infty}^0 \cap L_{n,(n)}(\R)) A_{n,(n)}(\R)^0, \end{array} \]
or, using the equality $N(L_{n,(n),\lin}^{(m)})(\A^\infty)=N(L_{n,(n),\lin}^{(m)})(\Q) (U \cap N(L_{n,(n),\lin}^{(m)})(\A^\infty))$, even with
\[ \begin{array}{c} \!\!\! \left( \!(N(L_{n,(n),\lin}^{(m)})(\Q)\! \cap\! U) \backslash N(L_{n,(n),\lin}^{(m)})(\R) \!   \right)\!\!  \times \!\!  \left( L_{n,(n)}(\A)/ (U_{n,\infty}^0\! \cap L_{n,(n)}(\R)) A_{n,(n)}(\R)^0\right)  \\ \da \\ L_{n,(n)}(\A)/(U_{n,\infty}^0  \cap L_{n,(n)}(\R)) A_{n,(n)}(\R)^0, \end{array} \]
The right $L_{n,(n)}^{(m)}(\A^\infty)$-action is by right translation on the second factor. The left action of $L_{n,(n)}(\Q)$ is via conjugation on the first factor and left translation on the second. 

The first part of the lemma follows, and we see that
\[ R^i\tpi^{(m)}_* \Omega  \]
is $L_{n,(n)}(\Q) \times L^{(m)}_{n,(n)}(\A^\infty)$ equivariantly identified with the locally constant sheaf
\[ \begin{array}{c} \left( \wedge^i N(L_{n,(n),\lin}^{(m)})(\Omega)^\vee    \right) \times \left( L_{n,(n)}(\A)/ (U_{n,\infty}^0 \cap \! L_{n,(n)}(\R)) A_{n,(n)}(\R)^0\right)  \\ \da \\ L_{n,(n)}(\A)/(U_{n,\infty}^0 \cap L_{n,(n)}(\R)) A_{n,(n)}(\R)^0. \end{array} \]
The lemma follows.
\pfend

\begin{lem} There is an $L_{n,(n)}^{(m)}(\A^\infty)$-equivariant isomorphism
\[ H^{k}_\Int(\gT^{(m)}_{(n)},\Omega) \cong \bigoplus_{i+j=k} H^i_\Int(\gT_{(n)}, \cL_{\wedge^j(\bigoplus_\tau \Std_\tau^{\oplus m})^\vee}). \]
\end{lem}

\pfbegin
There is an $L_{n,(n)}^{(m)}(\A^\infty)$-equivariant spectral sequence
\[ E_{2}^{i,j}=H^i(\gT_{(n)}, \cL_{\wedge^j(\bigoplus_\tau \Std_\tau^{\oplus m})^\vee}) \Rightarrow H^{i+j}(\gT^{(m)}_{(n)},\Omega). \]
If $\alpha \in \Q^\times_{>0} \subset Z(L_{n,(n),\lin})(\A^\infty)$, then $\alpha$ acts on $E_{2}^{i,j}$ via $\alpha^j$. 
We deduce that all the differentials (on the second and any later page) vanish, i.e. the spectral sequence degenerates on the second page. Moreover the $\alpha \mapsto \alpha^j$ eigenspace in $H^{i+j}(\gT^{(m)}_{(n)},\Omega)$ is naturally identified with $H^i(\gT_{(n)}, \cL_{\wedge^j(\bigoplus_\tau \Std_\tau^{\oplus m})^\vee})$. (This standard argument is sometimes referred to as `Lieberman's trick'.)

As the maps $\pi^{(m)}$ are proper, there is also a $L_{n,(n)}^{(m)}(\A^\infty)$-equivariant spectral sequence 
\[ E_{c,2}^{i,j}=H^i_c(\gT_{(n)}, \cL_{\wedge^j(\bigoplus_\tau \Std_\tau^{\oplus m})^\vee}) \Rightarrow H^{i+j}_c(\gT^{(m)}_{(n)},\Omega) \]
and $\alpha \in \Q^\times_{>0} \subset Z(L_{n,(n),\lin})(\A^\infty)$ acts on $E_{c,2}^{i,j}$ via $\alpha^j$. 
Again we see that the spectral sequence degenerates on the second page, and that the $\alpha \mapsto \alpha^j$ eigenspace in $H^{i+j}_c(\gT^{(m)}_{(n)},\Omega)$ is naturally identified with $H^i_c(\gT_{(n)}, \cL_{\wedge^j(\bigoplus_\tau \Std_\tau^{\oplus m})^\vee})$.

The lemma follows.
\pfend

\begin{cor}\label{cor196} Suppose that $\rho$ is an irreducible representation of $L_{n,(n),\lin}$ over $\Omega$, which we extend to a representation of $L_{n,(n)}$ by letting it be trivial on $L_{n,(n),\herm}$. Let $d=\norm_{F/\Q}\circ \det: L_{n,(n),\lin} \ra \G_m$. Then for all $N$ sufficiently large there are $j(N),\, m(N) \in \Z_{\geq 0}$ such that, for all $i$, 
\[ H^i_\Int( \gT_{(n)},\cL_{\rho \otimes d^{-N}}) \]
is an $L_{n,(n)}(\A^\infty)$-direct summand of 
\[ H^{i+j(N)}_\Int(\gT_{(n)}^{(m(N))},\Omega). \] 
\end{cor} 

\pfbegin It follows from Weyl's construction of the irreducible representations of $GL_n$ that, for $N$ sufficiently large, $\rho \otimes d^{-N}$ is a direct summand of 
\[ \bigotimes_\tau (\Std_\tau^{\vee})^{\otimes m_\tau(N)} \]
 for certain non-negative integers $m_\tau(N)$. Hence for $N$ sufficiently large and $m(N)=\max \{m_\tau(N)\}$ the representation $\rho \otimes d^{-N}$ is also a direct summand of 
\[ \wedge^{\sum_\tau m_\tau(N)}(\bigoplus_\tau \Std_\tau^{\oplus m(N)})^\vee. \]
\pfend

\begin{lem}\label{cuspcoh} Suppose that $\rho$ is an irreducible algebraic representation of $L_{n,(n)}$ on a finite dimensional $\C$-vector space.\begin{enumerate}
\item Then 
\[ \bigoplus_\Pi \Pi^\infty \otimes H^i(\Lie L_{n,(n)}, (U_{n,\infty}^0 \cap L_{n,(n)}(\R))  A_{n,(n)}(\R)^0, \Pi_\infty \otimes \rho)  \into  H^i_\Int(\gT_{(n)}, \cL_\rho), \]
where $\Pi$ runs over cuspidal automorphic representations of $L_{n,(n)}(\A)$.

\item If $n>1$ and if $\Pi$ is a cuspidal automorphic representation of $L_{n,(n)}(\A)$ such that $\Pi_\infty$ has the same infinitesimal character as $\rho^\vee$, then
\[ H^i(\Lie L_{n,(n)}, (U_{n,\infty}^0 \cap L_{n,(n)}(\R))  A_{n,(n)}(\R)^0, \Pi_\infty \otimes \rho) \neq (0) \]
for some $i>0$.
\end{enumerate}\end{lem}

\pfbegin The first part results from \cite{borel}, more precisely from combining theorem 5.2, the discussion in section 5.4 and corollary 5.5 of that paper. The second part results from \cite{clozelaa}, see the proof of theorem 3.13, and in particular lemma 3.14, of that paper.
\pfend

Combining this lemma and corollary \ref{cor196} we obtain the following consequence.

\begin{cor}\label{toptocone} Suppose that $n>1$ and that $\rho$ is an irreducible algebraic representation of $L_{n,(n),\lin}$ on a finite dimensional $\C$-vector space.
Suppose also that
$\pi$ is a cuspidal automorphic representation of $L_{n,(n),\lin}(\A)$ so that $\pi_\infty$ has the same infinitesimal character as $\rho^\vee$ and that $\psi$ is a continuous character of $\Q^\times \backslash \A^\times/ \R^\times_{>0}$. Then for all sufficiently large integers $N$ there are integers $m(N) \in \Z_{\geq 0}$ and $i(N) \in \Z_{>0}$, and a $L_{n,(n)}(\A^\infty)$-equivariant embedding
\[ (\pi ||\det||^N)\times \psi \into H^{i(N)}_\Int(\gT_{(n)}^{(m(N))},\C). \] 
\end{cor}

\newpage

\section{Tori, torsors and torus embeddings.}

Throughout this section let $R_0$ denote an irreducible noetherian ring (i.e. a noetherian ring with a unique minimal prime ideal). In the applications of this section elsewhere in this paper it will be either $\Q$ or $\Z_{(p)}$ or $\Z/p^r\Z$ for some $r$. We will consider $R_0$ endowed with the discrete topology so that $\Spf R_0 \cong \Spec R_0$.

\subsection{Tori and torsors.}

If $S/Y$ is a torus (i.e. a group scheme etale locally on $Y$ isomorphic to $\G_m^N$ for some $N$) then we can define its sheaf of characters $X^*(S)=\Hom(S,\G_m)$ and its sheaf of cocharacters $X_*(S)=\Hom(\G_m,S)$. These are locally constant sheaves of free $\Z$-modules in the etale topology on $Y$. They are naturally $\Z$-dual to each other. More generally if $S_1/Y$ and $S_2/Y$ are two tori then $\Hom(S_1,S_2)$ is a locally constant sheaf of free $\Z$-modules in the etale topology on $Y$. In fact
\[ \Hom(S_1,S_2)=\Hom(X_*(S_1),X_*(S_2))= \Hom(X^*(S_2),X^*(S_1)). \]
By a {\em quasi-isogeny} (resp. {\em isogeny}) from $S_1$ to $S_2$ we shall mean a global section of the sheaf $\Hom(S_1,S_2)_\Q$ (resp. $\Hom(S_1,S_2)$) with an inverse in $\Hom(S_2,S_1)_\Q$. We will write $[S]_\isog$ for the category whose objects are tori over $Y$ quasi-isogenous to $S$ and whose morphisms are isogenies. The sheaves $X_*(S)_\Q$ and $X^*(S)_\Q$ only depend on the quasi-isogeny class of $S$ so we will write $X_*([S]_\isog)_\Q$ and $X^*([S]_\isog)_\Q$. 

If $\bary$ is a geometric point of $Y$ then we define 
\[ TS_\bary = \lim_{\substack{\lla \\ N}} S[N](k(\bary)) \]
and
\[ T^pS_\bary= \lim_{\substack{\lla \\ p \ndiv N}} S[N](k(\bary)) \]
with the transition map from $MN$ to $N$ being multiplication by $M$. (The Tate modules of $S$.) Also define
\[ VS_\bary=TS_\bary \otimes_\Z \Q \]
and
\[ V^pS_\bary=TS_\bary \otimes_\Z \Q. \]
If $Y$ is a scheme over $\Spec \Q$ then 
\[ TS_\bary\cong X_*(S)_\bary \otimes_\Z \hatZ(1). \]
If $Y$ is a scheme over $\Spec \Z_{(p)}$ then
\[ T^pS_\bary\cong X_*(S)_\bary \otimes_\Z \hatZ^p(1). \]

Now suppose that $S$ is split, i.e. isomorphic to $\G_m^N$ for some $N$. By an {\em $S$-torsor} $T/Y$ we mean a scheme $T/Y$ with an action of $S$, which locally in the Zariski topology on $Y$ is isomorphic to $S$. By a {\em rigidification} of $T$ along $e: Y' \ra Y$ we mean an isomorphism of $S$-torsors $e^* T \cong S$ over $Y'$. If $U$ is a connected open subset of $Y$ then
\[ T|_U = \underline{\Spec} \bigoplus_{\chi \in X^*(S)(U)} \cL_T(\chi), \]
where $\cL_T(\chi)$ is a line bundle on $U$.  If $Z$ is any open subset of $Y$ and if $\chi \in X^*(S)(Z)$ then there is a unique line bundle $\cL_T(\chi)$ on $Z$ whose restriction to any connected open subset $U \subset Z$ is $\cL_T(\chi|_U)$. Multiplication gives isomorphisms
\[\cL_{T}(\chi_1) \otimes \cL_{T}(\chi_2) \liso \cL_{T}(\chi_1+\chi_2). \]
 Note that if $U$ has infinitely many connected components then it may not be the case that $T|_U = \underline{\Spec} \bigoplus_{\chi \in X^*(S)(U)} \cL_T(\chi)$. 
The map
\[ T \longmapsto \cL_{T,1}^\vee \]
gives a bijection between isomorphism classes of $\G_m$-torsors and isomorphism classes of line bundles on $Y$. The inverse map sends $\cL$ to
\[ \underline{\Spec} \bigoplus_{N \in \Z} \cL^{\vee, \otimes N}. \]

If $\alpha: S \ra S'$ is a morphism of split tori and if $T/Y$ is an $S$-torsor we can form a pushout $\alpha_*T$, an $S'$ torsor on $Y$ defined as the quotient
\[ (S' \times_Y T)/S \]
where $S$ acts by
\[ s:(s',t) \longmapsto (s's,s^{-1}t). \]
There is a natural map $T \ra \alpha_* T$ compatible with $\alpha: S \ra S'$. If $\alpha$ is an isogeny then
\[ \alpha_* T = (\ker \alpha) \backslash T. \]

If $T_1$ and $T_2$ are $S$-torsors over $Y$ we define 
\[ (T_1 \otimes_S T_2)/Y \]
to be the $S$-torsor
\[ (T_1 \times_Y T_2) /S \]
where $S$ acts by 
\[ s:(t_1,t_2) \longmapsto (st_1,s^{-1}t_2). \]
If $T$ is an $S$-torsor on $Y$ we define an $S$-torsor $T^\vee/Y$ by taking $T^\vee=T$ as schemes but defining an $S$ action $.$ on $T^\vee$ by
\[ s.t=s^{-1}t, \]
i.e. $T^\vee = [-1]_{S,*} T$.
Then
\[ T^\vee \otimes_S T \cong S \]
via the map that sends $(t_1,t_2)$ to the unique section $s$ of $S$ with $st_1=t_2$.

\newpage \subsection{Log structures.}\label{logstruct}

 We will call a formal scheme 
\[ \gX \lra \Spf R_0 \]
{\em suitable} if it has a cover by affine opens $\gU_i = \Spf (A_i)^\wedge_{I_i}$, where $A_i$ is a finitely generated $R_0$-algebra and $I_i$ is an ideal of $A_i$ whose inverse image in $R_0$ is $(0)$. 

By a {\em log structure} on a scheme $X$ (resp. formal scheme $\gX$) we mean a sheaf of monoids $\cM$ on $X$ (resp. $\gX$) together with a morphism
\[ \alpha: \cM \lra (\cO_X,\times) \]
(resp.
\[ \alpha: \cM \lra (\cO_\gX,\times)) \]
such that the induced map
\[ \alpha^{-1} \cO_X^\times \lra \cO_X^\times \]
(resp.
\[ \alpha^{-1} \cO_\gX^\times \lra \cO_\gX^\times) \]
is an isomorphism. We will refer to a scheme (resp. formal scheme) endowed with a log structure as a {\em log scheme} (resp. {\em log formal scheme}). By a {\em morphism of log schemes} (resp. {\em morphism of log formal schemes})
\[ (\phi,\psi):(X,\cM,\alpha) \lra (Y,\cN,\beta)\]
(resp. 
\[ (\phi,\psi):(\gX,\cM,\alpha) \lra (\gY,\cN,\beta)\,) \]
we shall mean a morphism $\phi:X \ra Y$ (resp. $\phi:\gX \ra \gY$) and a map
\[ \psi: \phi^{-1} \cN \lra \cM \]
such that $\phi^* \circ \phi^{-1}(\beta) = \alpha \circ \psi$. We will consider $R_0$ endowed with the trivial log structure $(\cO_{\Spec R_0}^\times, 1)$ (resp. $(\cO_{\Spf R_0}^\times, 1)$). We will call a log formal scheme $(\gX,\cM,\alpha)/\Spf R_0$ {\em suitable} if $\gX/\Spf R_0$ is suitable and if, locally in the Zariski topology, $\cM/\alpha^* \cO_\gX^\times$ is finitely generated. (In the case of schemes these definitions are well known. We have not attempted to optimize the definition in the case of formal schemes. We are simply making a definition which works for the limited purposes of this article.) 

If $X/\Spec R_0$ is a scheme of finite type and if $Z \subset X$ is a closed sub-scheme which is flat over $\Spec R_0$, then the formal completion $X_Z^\wedge$ is a suitable formal scheme. Let $i^\wedge$ denote the map of ringed spaces $X_Z^\wedge \ra X$. If $(\cM,\alpha)$ is a log structure on $X$, then we get a map
\[ (i^\wedge)^{-1}(\alpha): (i^\wedge)^{-1} \cM \lra \cO_{X_Z^\wedge}. \]
It induces a log structure $(\cM^\wedge,\alpha^\wedge)$ on $X_Z^\wedge$, where $\cM^\wedge$ denotes the push out  
\[ \begin{array}{rcl} ((i^\wedge)^{-1}(\alpha))^{-1} \cO_{X_Z^\wedge}^\times & \into & (i^\wedge)^{-1} \cM \\ \da &&\da \\ \cO_{X_Z^\wedge}^\times & \lra & \cM^\wedge. \end{array} \]

If
\[ (\phi,\psi):(X,\cM,\alpha) \lra (Y,\cN,\beta)\]
is a morphism of schemes with log structures over $\Spec R_0$ then there is a right exact sequence
\[ \phi^* \Omega^1_Y(\log \cN) \lra \Omega^1_X(\log \cM) \lra \Omega^1_{X/Y}(\log \cM/\cN) \lra (0) \]
of sheaves of log differentials.
If the map $(\phi,\psi)$ is log smooth then this sequence is also left exact and the sheaf $\Omega^1_{X/Y}(\log \cM/\cN)$ is locally free. As usual, we write $\Omega^i_X(\log \cM)=\wedge^i \Omega^1_X(\log \cM)$ and $\Omega^i_{X/Y}(\log \cM/\cN)=\wedge^i\Omega^1_{X/Y}(\log \cM/\cN)$.

By a {\em coherent sheaf of differentials} on a formal scheme $\gX/\Spf R_0$ we will mean a coherent sheaf $\Omega/\gX$ together with a differential $d:\cO_\gX \ra \Omega$ which vanishes on $R_0$. By a {\em coherent sheaf of log differentials} on a log formal scheme $(\gX,\cM,\alpha)/\Spf R_0$ we shall mean a coherent sheaf $\Omega/\gX$ together with a differential, which vanishes on $R_0$,
\[ d: \cO_\gX \lra \Omega, \]
and a homomorphism
\[ \dlog : \cM \lra \Omega \]
such that
\[ \alpha (m) \dlog m = d( \alpha (m)) . \]
By a {\em universal coherent sheaf of differentials} (resp. {\em universal coherent sheaf of log differentials}) we shall mean a coherent sheaf of differentials $(\Omega,d)$ (resp. a coherent sheaf of log differentials $(\Omega,d,\dlog)$) such that for any other coherent sheaf of differentials $(\Omega',d')$ (resp. a coherent sheaf of log differentials $(\Omega',d',\dlog')$) there is a unique map $f:\Omega \ra \Omega'$ such that $f \circ d=d'$ (resp. $f \circ d = d'$ and $f \circ \dlog = \dlog'$). 

Note that if a universal coherent sheaf of differentials (resp. universal coherent sheaf of log differentials) exists, it is unique up to unique isomorphism.

\begin{lem} Suppose that $R_0$ is a discrete, noetherian topological ring.
\begin{enumerate}
\item A universal sheaf of coherent differentials $\Omega^1_{\gX/\Spf R_0}$ exists for any suitable formal scheme  $\gX/\Spf R_0$. 
\item If $X/\Spec R_0$ is a scheme of finite type and if $Z \subset X$ is flat over $R_0$ then 
\[ \Omega_{X_Z^\wedge/\Spf R_0}^1 \cong (\Omega^1_{X/\Spec R_0})^\wedge. \]
\item A universal sheaf of coherent log differentials $\Omega^1_{\gX/\Spf R_0}(\log\cM)$ exists for any suitable log formal scheme  $(\gX,\cM,\alpha)/\Spf R_0$. 
\item Suppose that $X/\Spec R_0$ is a scheme of finite type, that $Z \subset X$ is flat over $R_0$ and that $(\cM,\alpha)$ is a log structure on $X$ such that Zariski locally $\cM/\alpha^{-1} \cO_X^\times$ is finitely generated. Then 
\[ \Omega_{X_Z^\wedge/\Spf R_0}^1(\log \cM^\wedge) \cong (\Omega^1_{X/\Spec R_0}(\log \cM))^\wedge. \]
\end{enumerate} \end{lem} 

\pfbegin
Consider the first part.
Suppose that $\gU = \Spf A_I^\wedge$ is an affine open in $\gX$, where $A$ is a finitely generated $R_0$-algebra and $I$ is an ideal of $A$ with inverse image $(0)$ in $R_0$. Then there exists a universal finite module of differentials $\Omega^1_{\gU}$ for $\gU$, namely the coherent sheaf of $\cO_\gU$-modules associated to $(\Omega^1_{A/R_0})_I^\wedge$. (See sections 11.5 and 12.6 of \cite{kunz}.) We must show that if $\gU' \subset \gU$ is open then $\Omega^1_{\gU}|_{\gU'}$ is a universal finite module of differentials for $\gU'$. For then uniqueness will allow us to glue the coherent sheaves $\Omega^1_{\gU}$ to form $\Omega^1_{\gX}$.

So suppose that $(\Omega',d')$ is a finite module of differentials for $\gU'$. We must show that there is a unique map of $\cO_{\gU'}$-modules
\[ f: \Omega^1_{\gU}|_{\gU'} \lra \Omega' \]
such that $d'=f \circ d$. We may cover $\gU'$ by affine opens of the form $\Spf (A_g)^\wedge_I$ and it will suffice to find, for each $g$, a unique 
\[  f_g: \Omega^1_{\gU}|_{\Spf (A_g)^\wedge_I} \lra \Omega'|\Spf (A_g)^\wedge_I \]
with $d'=f_g \circ d$. Thus we may assume that $\gU'=\Spf (A_g)^\wedge_I$. But in this case we know $\Omega^1_{\gU'}$ exists, and is the coherent sheaf associated to
\[ (\Omega^1_{A_g/R_0})^\wedge_I \cong (\Omega^1_{A/R_0} \otimes_A A_g)^\wedge_I. \]
On the other hand $\Omega^1_{\gU}|_{\gU'}$ is the coherent sheaf associated to
\[ (\Omega^1_{A/R_0})^\wedge_I \otimes_{A_I^\wedge} (A_g)_I^\wedge. \]
Thus
\[ \Omega^1_{\gU'} \liso \Omega^1_{\gU}|_{\gU'} \]
and the first part follows. 
The second part also follows from the proof of the first part.

For the third part, because of uniqueness, it suffices to work locally. Thus we may assume that there are finitely many sections $m_1,...,m_r \in \cM(\gX)$, which together with $\alpha^{-1} \cO_\gX^\times$ generate $\cM$. Then we define $\Omega^1_{(\gX,\cM,\alpha)}$ to be the cokernel of the map
\[ \begin{array}{rcl} \cO_\gX^{\oplus r} & \lra & \Omega^1_\gX \oplus \cO_{\gX}^{\oplus r} \\
(f_i)_i & \longmapsto & (-\sum_i f_i d \alpha(m_i),(f_i \alpha(m_i))_i). \end{array} \]
It is elementary to check that this has the desired universal property. The fourth part is also elementary to check.
\pfend

If
\[ (\phi,\psi): (\gX,\cM,\alpha) \lra (\gY,\cN,\beta) \]
is a map of suitable log formal schemes over $\Spf R_0$ then we set
\[ \Omega^1_{\gX/\gY}(\log \cM/\cN) = \Omega^1_{\gX/\Spf R_0}(\log \cM)/\phi^* \Omega^1_{\gY/\Spf R_0}(\log \cN). \]
We also set 
\[ \Omega^i_{\gX/\Spf R_0}=\wedge^i \Omega^1_{\gX/\Spf R_0} \]
and
\[ \Omega^i_{\gX/\Spf R_0}(\log \cM) = \wedge^i \Omega^1_{\gX/\Spf R_0}(\log \cM) \]
and
\[ \Omega^i_{\gX/\gY}(\log \cM/\cN)= \wedge^i \Omega^1_{\gX/\gY}(\log \cM/\cN). \]

\begin{cor} Suppose that $R_0$ is a discrete, noetherian topological ring; that $(X,\cM,\alpha) \ra (Y,\cN,\beta)$ is a map of log schemes over $\Spec R_0$; and that $Z \subset X$ and $W \subset Y$ are closed sub-schemes flat over $\Spec R_0$ which map to each other under $X \ra Y$. Suppose moreover that $X$ and $Y$ have finite type over $\Spec R_0$ and that $\cM/\alpha^{-1} \cO_\gX^\times$ and $\cN/\beta^{-1} \cO_\gY^\times$ are locally (in the Zariski topology) finitely generated. Then
\[ \Omega^1_{(X_Z^\wedge,\cM^\wedge,\alpha^\wedge)/(Y_W^\wedge,\cN^\wedge,\beta^\wedge)} \cong ( \Omega^1_{X/Y}(\log \cM/\cN))_X^\wedge. \]
\end{cor}

\pfbegin This follows from the lemma and from the exactness of completion. \pfend

If $Y$ is a scheme we will let 
\[ \Aff^n_Y=\underline{\Spec} \cO_Y[T_1,...,T_n] \]
 denote affine $n$-space over $Y$ and
 \[ \Coord^n_Y=\underline{\Spec} \cO_Y[T_1,...,T_n]/(T_1...T_n) \subset \Aff^n_Y \]
 denote the union of the coordinate hyperplanes in $\Aff^n_Y$.
Now suppose that $X \ra Y$ is a smooth map of schemes of relative dimension $n$. By a {\em simple normal crossing divisor in $X$ relative to $Y$} we shall mean a closed subscheme $D \subset X$ such that 
$X$ has an affine Zariski-open cover $\{ U_i\}$ such that each $U_i$ admits an etale map $f_i:U_i \ra \Aff^n_Y$ so that $D|_{U_i}$ is the (scheme-theoretic) preimage of $\Coord^n_Y$.
In the case that $Y$ is just the spectrum of a field we will refer simply to a {\em simple normal crossing divisor in $X$}.

Suppose that $Y$ is locally noetherian and separated, and that the connected components of $Y$ are irreducible. 
If $S$ is a finite set of irreducible components of $D$ we will set
\[ D_S=\bigcap_{E \in S} E. \]
It is smooth over $Y$.
We will also set 
\[ D^{(s)} = \coprod_{\# S = s} D_S. \]
If $E$ is an irreducible component of $D^{(s)}$ then the set $S(E)$ of irreducible components of $D$ containing $E$ has cardinality $s$. 
If $\geq$ is a total order on the set of irreducible components of $D$, we can define a delta set $\cS(D, \geq)$, or simply $\cS(D)$, as follows. (For the definition of delta set, see for instance \cite{rs}. We can, if we prefer to be more abstract, replace $\cS(D,\geq)$ by the associated simplicial set.)  The $n$ cells consist of all irreducible components of $D^{(n+1)}$. If $E$ is such an irreducible component and if $i\in \{0,...,n\}$ then the image of $E$ under the face map $d_i$ is the unique irreducible component of 
\[ \bigcap_{F \in S(E)_i} F \]
which contains $E$. Here $S(E)_i$ equals $S(E)$ with its $(i+1)^{th}$ smallest element removed.
The topological realization $|\cS(D,\geq)|$ does not depend on the total order $\geq$, so we will often write $|\cS(D)|$.

If $D$ is a simple normal crossing divisor in $X$ relative to $Y$ we define a log structure $\cM(D)$ on $X$ by setting 
\[ \cM(D)(U) = \cO_X(U) \cap \cO_X(U-D)^\times. \]

We record a general observation about log de Rham complexes and divisors with simple normal crossings, which is probably well known. We include a proof because it is of crucial importance for our argument.

\begin{lem}\label{log} Suppose that $Y$ is a smooth scheme of finite type over a field $k$ and that $Z \subset Y$ is a divisor with simple normal crossings. Let $Z_1,...,Z_m$ denote the distinct irreducible components of $Z$ and set
\[ Z_S = \bigcap_{j \in S} Z_j \subset Y\]
(in particular $Z_\emptyset =Y$), and
\[ Z^{(s)} = \coprod_{\# S=s} Z_S. \]
Let $i_S$ (resp. $i^{(s)}$) denote the natural maps $Z_S \ra Y$ (resp. $Z^{(s)} \ra Y$). Also let $\cI_Z$ denote the ideal of definition of $Z$.

There is a double complex 
\[ i^{(s)}_* \Omega^r_{Z^{(s)}} \]
with maps
\[ d: i^{(s)}_* \Omega^r_{Z^{(s)}} \lra i^{(s)}_* \Omega^{r+1}_{Z^{(s)}} \]
and
\[ i^{(s)}_*  \Omega^r_{Z^{(s)}} \lra i^{(s+1)}_* \Omega^{r}_{Z^{(s+1)}} \]
being the sum of the maps
\[ i_{S,*}\Omega^r_{Z_S} \lra i_{S',*} \Omega^{r}_{Z_{S'}}, \]
which are 
\begin{itemize}
\item $0$ if $ S \not\subset S'$,
\item and $(-1)^{\# \{ i \in S: \,\,\, i<j\} }$ times the natural pull-back if $S \cup \{ j\} =S'$.
\end{itemize}

The natural inclusions 
\[ \Omega^r_Y(\log \cM(Z)) \otimes \cI_Z \lra \Omega^r_Y \]
give rise to a map of complexes
\[ \Omega^\bullet_Y(\log \cM(Z)) \otimes \cI_Z \lra \Omega^r_Y = i^{(0)}_* \Omega^r_{Z^{(0)}}. \]

For fixed $r$ the simple complexes 
\[ (0) \lra \Omega^r_Y(\log \cM(Z)) \otimes \cI_Z \lra  i^{(0)}_* \Omega^r_{Z^{(0)}} \lra  i^{(1)}_* \Omega^r_{Z^{(1)}} \lra ... \]
are exact.
\end{lem}

\pfbegin
Only the last assertion is not immediate. So consider the last assertion. We can work Zariski locally, so we may assume that the complex is pulled back from the corresponding complex for the case $Y=\Spec k[X_1,...,X_d]$ and $Z$ is given by $X_1X_2...X_m=0$. In this case we take $Z_j$ to be the scheme $X_j=0$, for $j=1,...,m$. In this case
\[ \Omega^r_Y(\log \cM(Z)) \otimes \cI_Z = \bigoplus_{T} k[X_1,...,X_d]  \left( \prod_{j=1, \,\, j \not\in T}^m X_j \right) \bigwedge_{j \in T} dX_j \]
where $T$ runs over r element subsets of $\{1,...,d\}$. On the other hand 
\[ i_{S,*} \Omega^r_{Z_S} = \bigoplus_T k[X_1,...,X_d]/(X_j)_{j \in S} \bigwedge_{j \in T} dX_j \]
where $T$ runs over $r$ element subsets of $\{ 1,...,d\}-S$. Thus it suffices to show that, for each subset $T \subset \{ 1,...,d\}$ the sequence
\[ \begin{array}{r} (0) \lra \left( \prod_{j=1, \,\, j \not\in T}^m X_j \right) k[X_1,...,X_d] \lra k[X_1,...,X_d] \lra ... \\ ...\lra \bigoplus_{\# S=s, \,\, S\cap T = \emptyset} k[X_1,...,X_d]/(X_j)_{j \in S} \lra ... \end{array} \]
is exact, where $S \subset \{ 1,...,m\}$. The sequence for $T \subset \{1,...,d\}$ is obtained from the sequence for $\emptyset \subset \{1,...,m\}-T$ by tensoring over $k$ with $k[X_j]_{j \in T \cup \{ m+1,...,d\} }$, and so we only need treat  the case $m=d$ and $T=\emptyset$. 

If $\mu$ is a monomial in the variables $X_1,...,X_m$, let $R(\mu)$ denote the subset of $\{ 1,...,m\}$ consisting of the indices $j$ for which $X_j$ does not occur in $\mu$. Then our complex is the direct sum over $\mu$ of the complexes
\[ (0) \lra A_\mu \lra k \lra ... \lra \bigoplus_{S \subset R(\mu),\,\, \# S = s} k \lra ... \]
where $A_\mu=k$ if $R(\mu)=\emptyset$ and $=(0)$ otherwise. So it suffices to prove this latter complex exact for all $\mu$. If $R(\mu)=\emptyset$ then it becomes
\[ (0) \lra k \lra k \lra (0) \lra (0) \lra ..., \]
which is clearly exact. If $R(\mu) \neq \emptyset$, our complex becomes 
\[ (0) \lra k \lra \bigoplus_{S \subset R(\mu),\,\, \# S = 1} k  \lra ... \lra \bigoplus_{S \subset R(\mu),\,\, \# S = s} k \lra ... . \]
If we suppress the first $k$, this is the complex that computes the simplicial cohomology with $k$-coefficients of the simplex with $\# R(\mu)$ vertices. Thus it is exact everywhere except $\bigoplus_{S \subset R(\mu),\,\, \# S = 1} k$ and the kernel of 
\[  \bigoplus_{S \subset R(\mu),\,\, \# S = 1} k  \lra \bigoplus_{S \subset R(\mu),\,\, \# S = 2} k \]
is just $k$. The desired exactness follows.
\pfend

\newpage \subsection{Torus embeddings.}

We will now discuss relative torus embeddings. We will suppose that $Y/\Spec R_0$ is flat and locally of finite type. To simplify the notation, for now we will restrict to the case of a split torus $S/Y$ with $Y$ connected. We will record the (trivial) generalization to the case of a disconnected base below. Thus we can think of $X^*(S)$ and $X_*(S)$ as abelian groups, rather than as locally constant sheaves on $Y$, i.e. we replace the sheaf by its global sections over $Y$. We will let $T/Y$ denote an $S$-torsor. 

By a {\em rational polyhedral cone} $\sigma \subset X_*(S)_\R$ we mean a non-empty subset consisting of all $\R_{\geq 0}$-linear combinations of a finite set of elements of $X_*(S)$, but which contains no complete line through $0$. (We include the case $\sigma=\{0\}$. The notion we define here is sometimes called a `non-degenerate rational polyhedral cone'.) By the {\em interior} $\sigma^0$ of $\sigma$ we shall mean the complement in $\sigma$ of all its proper faces. (We consider $\sigma$ as a face of $\sigma$, but not a proper face.) We call $\sigma$ {\em smooth} if it consists of all $\R_{\geq 0}$-linear combinations of a subset of a $\Z$-basis of $X_*(S)$. Note that any face of a smooth cone is smooth. Then we define $\sigma^\vee$ to be the set of elements of $X^*(S)_\R$ which have non-negative pairing with every element of $\sigma$ and $\sigma^{\vee,0}$ to be the set of elements of $X^*(S)_\R$ which have strictly positive pairing with every element of $\sigma-\{0\}$. Moreover we set
\[ T_\sigma = \underline{\Spec} \bigoplus_{\chi \in X^*(S) \cap \sigma^\vee} \cL_T(\chi). \]
Then $T_\sigma$ is a scheme over $Y$ with an action of $S$ and there is a natural $S$-equivariant dense open embedding $T \into T_\sigma$. If $\sigma' \subset \sigma$ there is a natural map $T_{\sigma'} \ra T_\sigma$ compatible with the embeddings of $T$.  If $f:Y' \ra Y$ then $T_\sigma/Y$ pulls back under $f$ to $(f^*T)_\sigma/Y'$ compatibly with the maps $T_{\sigma'}\into T_\sigma$ for $\sigma'\subset \sigma$. 

Suppose that $\Sigma_0$ is a non-empty set of faces of $\sigma$ such that 
\begin{itemize}
\item $\{0\} \not\in \Sigma_0$,
\item and, if $\tau' \supset \tau \in \Sigma_0$, then $\tau' \in \Sigma_0$. 
\end{itemize}
In this case define
\[ |\Sigma_0|^0=\sigma-\bigcup_{\tau \not\in \Sigma_0} \tau \]
and
\[ |\Sigma_0|^{\vee,0} = \{ \chi \in X^*(S)_\R:\,\, \chi>0\,\, {\rm on}\,\, |\Sigma_0|^0\}.\]
Then we define $\partial_{\Sigma_0} T_\sigma \subset T_\sigma$ to be the closed sub-scheme defined by the sheaf of ideals
\[ \bigoplus_{\chi \in X^*(S) \cap |\Sigma_0|^{\vee,0}} \cL_T(\chi) \subset \bigoplus_{\chi \in X^*(S) \cap \sigma^\vee} \cL_T(\chi). \]
If $\Sigma_0$ contains all the faces of $\sigma$ other than $\{0\}$ we will write $\partial T_\sigma$ for $\partial_{\Sigma_0} T_\sigma$.
If $\sigma'$ is a face of $\sigma$ then under the open embedding
\[ T_{\sigma'} \into T_\sigma \]
$\partial_{\Sigma_0}T_\sigma$ pulls back to $\partial_{\{ \tau \in \Sigma_0: \,\, \tau \subset \sigma'\}} T_{\sigma'}$. 

By a {\em fan} in $X_*(S)_\R$ we shall mean a non-empty collection $\Sigma$ of rational polyhedral cones $\sigma \subset X_*(S)_\R$ which satisfy
\begin{itemize}
\item if $\sigma \in \Sigma$, so is each face of $\sigma$,
\item if $\sigma, \sigma' \in \Sigma$ then $\sigma \cap \sigma'$ is a face of $\sigma$ and of $\sigma'$.
\end{itemize}
We call $\Sigma$ {\em smooth} if each $\sigma \in \Sigma$ is smooth. 
We will call $\Sigma$ {\em full} if every element of $\Sigma$ is contained in an element of $\Sigma$ with the same dimension as $X_*(S)_\R$. 
Define
\[ |\Sigma|= \bigcup_{\sigma \in \Sigma} \sigma. \]
Also define 
\[ |\Sigma|^{\vee}  = \{ \chi \in X^*(S)_\R:\,\, \chi(|\Sigma|) \subset \R_{\geq 0}\} = \bigcap_{\sigma \in \Sigma} \sigma^{\vee} \]
and
\[ |\Sigma|^{\vee,0}  = \{ \chi \in X^*(S)_\R:\,\, \chi(|\Sigma|-\{0\}) \subset \R_{>0}\} = \bigcap_{\sigma \in \Sigma} \sigma^{\vee,0}. \]
We call $\Sigma'$ a {\em refinement} of $\Sigma$ if each $\sigma' \in \Sigma'$ is a subset of some element of $\Sigma$ and each element $\sigma \in \Sigma$ is a finite union of elements of $\Sigma'$. 

\begin{lem} \begin{enumerate}
\item If $\Sigma$ is a fan and $\Sigma' \subset \Sigma$ is a finite cardinality sub-fan then there is a refinement $\tSigma$ of $\Sigma$ with the following properties:
\begin{itemize}
\item any element of $\Sigma$ which is smooth also lies in $\tSigma$;
\item any element of $\tSigma$ contained in an element of $\Sigma'$ is smooth;
\item and if $\sigma' \in \Sigma -\tSigma$ then $\sigma'$ has a non-smooth face lying in $\Sigma'$.
\end{itemize}
\item Any fan $\Sigma$ has a smooth refinement $\Sigma'$ such that every smooth cone $\sigma \in \Sigma$ also lies in $\Sigma'$.  \end{enumerate} \end{lem}

\pfbegin The first part is proved just as for finite fans by making a finite series of elementary subdivisions by $1$ cones that lie in some element $\sigma' \in \Sigma'$ but not in any of its smooth faces. See for instance section 2.6 of \cite{fulton}. 

For the second part, consider the $\cS$ the set of pairs $(\tSigma,\Delta)$ where $\tSigma$ is a refinement of $\Sigma$ and $\Delta$ is a sub-fan of $\Sigma$ such that 
\begin{itemize}
\item every smooth element of $\Sigma$ lies in $\tSigma$;
\item and if $\sigma \in \tSigma$ is contained in an element of $\Delta$ then $\sigma$ is smooth. 
\end{itemize}
It suffices to show that $\cS$ contains an element $(\tSigma,\Delta)$ with $\Delta =\Sigma$. 

If $(\tSigma,\Delta) \in \cS$ and $\sigma \in \Sigma$ we define $\tSigma(\sigma)$ to be the set of elements of $\tSigma$ contained in $\sigma$.
We define a partial order on $\cS$ by decreeing that $(\tSigma,\Delta) \geq (\tSigma',\Delta')$ if and only if the following conditions are satisfied:
\begin{itemize}
\item $\tSigma$ refines $\tSigma'$;
\item $\Delta \supset \Delta'$;
\item $\tSigma'(\sigma) =\tSigma(\sigma)$ unless $\sigma$ has a face that is contained in an element of $\Delta$ but in no element of $\Delta'$.
\end{itemize}

Suppose that $\cS' \subset \cS$ is totally ordered.  Set
\[ \Delta = \bigcup_{(\tSigma',\Delta') \in \cS'} \Delta', \]
and let $\tSigma$ denote the set of cones $\sigma'$ which lie in  $\tSigma'$ for all sufficiently large elements of $(\tSigma',\Delta') \in \cS'$. If $\sigma \in \Sigma$ then we can choose $(\tSigma',\Delta') \in \cS'$ so that the number of faces of $\sigma$ in $\Delta'$ is maximal. If $(\tSigma',\Delta') \leq (\tSigma'',\Delta'') \in \cS'$ then $\tSigma'(\sigma)=\tSigma''(\sigma)$. Thus $\tSigma(\sigma)=\tSigma'(\sigma)$. We conclude that $\tSigma$ is a refinement of $\Sigma$. Thus $(\tSigma,\Delta) \in \cS$ and it is an upper bound for $\cS'$. 

By Zorn's lemma $\cS$ has a maximal element $(\tSigma,\Delta)$. We will show that $\Delta = \Sigma$, which will complete the proof of the lemma. Suppose not. Choose $\sigma \in \Sigma-\Delta$. Set $\Delta'$ to be the union of $\Delta$ and the faces of $\sigma$. Let $\tSigma'$ be a refinement of $\tSigma$ such that
\begin{itemize}
\item any element of $\tSigma$ which is smooth also lies in $\tSigma'$;
\item any element of $\tSigma'$ contained in $\sigma$ is smooth;
\item and if $\sigma' \in \tSigma -\tSigma'$ then $\sigma'$ has a non-smooth face contained in $\sigma$.
\end{itemize}
Then $(\tSigma',\Delta') \in \cS$ and $(\tSigma',\Delta') > (\tSigma,\Delta)$, a contradiction.
\pfend

To a fan $\Sigma$ one can attach a connected scheme $T_\Sigma$ that is separated, locally (on $T_\Sigma$) of finite type and flat over $Y$ of relative dimension $\dim_\R X_*(S)_\R$, together with an action of $S$ and an $S$-equivariant dense open embedding $T \into T_\Sigma$ over $Y$. The scheme $T_\Sigma$ has an open cover by the $T_\sigma$ for $\sigma \in \Sigma$ such that $T_{\sigma'} \subset T_{\sigma}$ if and only if $\sigma' \subset \sigma$. We write $\cO_{T_\Sigma}$ for the structure sheaf of $T_\Sigma$.  If $\Sigma$ is smooth then $T_\Sigma/Y$ is smooth. If $\Sigma$ is finite and $|\Sigma|=X_*(S)_\R$, then $T_\Sigma/Y$ is proper.  If $\Sigma' \subset \Sigma$ then $T_{\Sigma'}$ can be identified with an open sub-scheme of $T_\Sigma$. If $\Sigma'$ refines $\Sigma$ then there is an $S$-equivariant proper map 
\[ T_{\Sigma'} \ra T_\Sigma \]
which restricts to the identity on $T$: its restriction to $T_{\sigma'}$ equals the map 
\[ T_{\sigma'} \lra T_\sigma \into T_\Sigma \]
where $\sigma' \subset \sigma \in \Sigma$. 

By {\em boundary data for $\Sigma$} we shall mean a proper subset $\Sigma_0 \subset \Sigma$ such that 
$\Sigma - \Sigma_0$ is a fan. (Note that $\Sigma_0$ may not be closed under taking faces.)
If $\Sigma_0$ is boundary data we define $\partial_{\Sigma_0}T_\Sigma$ to be the closed subscheme of $T_\Sigma$ with
\[ (\partial_{\Sigma_0} T_\Sigma) \cap T_\sigma = \partial_{\{ \tau \in \Sigma_0:\,\, \tau \subset \sigma\}} T_\sigma. \]
Note that
\[ \partial_{\Sigma_0} T_\Sigma \subset \bigcup_{\sigma \in \Sigma_0} T_\sigma. \]
Thus $\partial_{\Sigma_0} T_\Sigma$ has an open cover by the sets 
\[ (\partial_{\Sigma_0}T_\Sigma)_\sigma = T_\sigma \cap \partial_{\Sigma_0}T_\Sigma \]
as $\sigma$ runs over $\Sigma_0$.
We write $\cI_{\partial_{\Sigma_0} T_\Sigma}$ for the ideal sheaf in $\cO_{T_\Sigma}$ defining $\partial_{\Sigma_0} T_\Sigma$. If $\Sigma_0 \subset \Sigma' \subset \Sigma$ then 
\[ \partial_{\Sigma_0} T_{\Sigma'} \liso \partial_{\Sigma_0} T_\Sigma. \]

Note that $\cI_{\partial_{\Sigma_0}T_\Sigma}|_{T_\sigma}$ corresponds to the ideal
\[ \bigoplus_{\chi \in \gX_{\Sigma_0,\sigma,1}} \cL_T(\chi) \]
of
\[ \bigoplus_{\chi \in X^*(S) \cap \sigma^\vee} \cL_T(\chi), \]
where
\[ \gX_{\Sigma_0,\sigma,1}= X^*(S) \cap \sigma^\vee - \bigcup_{\tau \in \Sigma_0, \tau \subset \sigma} \tau^\perp \]
and $\tau^\perp$ denotes the annihilator of $\tau$ in $X^*(S)_\R$. If we let $\gX_{\Sigma_0,\sigma,m}$ denote the set of sums of $m$ elements of $\gX_{\Sigma_0,\sigma,1}$, then $\cI_{\partial_{\Sigma_0}T_\Sigma}^m|_{T_\sigma}$ corresponds to the ideal
\[ \bigoplus_{\chi \in \gX_{\Sigma_0,\sigma,m}} \cL_T(\chi). \]
If $\sigma \not\in \Sigma_0$ then 
\[ \gX_{\Sigma_0,\sigma,m}= X^*(S) \cap \sigma^\vee \]
for all $m$. If on the other hand $\sigma \in \Sigma_0$ then
\[ \bigcap_m \gX_{\Sigma_0,\sigma,m} = \emptyset. \]
(For if $\chi \in \sigma^0$ then $\chi \geq m$ on $\gX_{\Sigma_0,\sigma,m}$.)

In the special case $\Sigma_0=\Sigma -\{ \{0\}\}$ we will write $\partial T_\Sigma$ for $\partial_{\Sigma_0} T_\Sigma$ and $\cI_{\partial T_\Sigma}$ for $\cI_{\partial_{\Sigma_0} T_\Sigma}$. Then
\[ T = T_\Sigma - \partial T_\Sigma. \]
We will write $\cM_\Sigma \ra \cO_{T_\Sigma}$ for the log structure corresponding to the closed embedding $\partial T_\Sigma \into T_\Sigma$. We will write $\Omega^1_{T_\Sigma/\Spec R_0}(\log \infty)$ for the log differentials $\Omega^1_{T_\Sigma/\Spec R_0}(\log \cM_\Sigma)$.

If $\Sigma$ is smooth then $\partial T_\Sigma$ is a simple normal crossings divisor on $T_\Sigma$ relative to $Y$.  

If $\Sigma_0$ is boundary data for $\Sigma$ we will set
\[ |\Sigma_0| = \bigcup_{\sigma \in \Sigma_0} \sigma .  \]
and
\[ |\Sigma_0|^0 = |\Sigma_0|-|\Sigma-\Sigma_0|.  \]
We will call $\Sigma_0$ 
\begin{itemize}
\item {\em open} if $|\Sigma_0|^0$ is open in $X_*(S)_\R$;
\item {\em finite} if it has finite cardinality;
\item {\em locally finite} if for every rational polyhedral cone $\tau \subset |\Sigma_0|$ (not necessarily in $\Sigma_0$) the intersection $\tau \cap |\Sigma_0|^0$ meets only finitely many elements of $\Sigma_0$. (We remark that although this condition may be intuitive in the case $|\Sigma_0|^0=|\Sigma_0|$, in other cases it may be less so.)
\end{itemize}

Let $\Sigma$ continue to denote a fan and $\Sigma_0$ boundary data for $\Sigma$. If $\sigma \in \Sigma$ we write 
\[ \Sigma(\sigma) = \{ \tau \in \Sigma: \,\,\tau \supset \sigma \} .\]
This is an example of boundary data for $\Sigma$.
If $\sigma \in \Sigma_0$ then
\[ \Sigma(\sigma) = \{ \tau \in \Sigma_0: \,\,\tau \supset \sigma \} \]
and we will sometimes denote it $\Sigma_0(\sigma)$. If $\Sigma_0$ is locally finite then $\Sigma_0(\sigma)$ is finite for all $\sigma \in \Sigma_0$.
If $\{0\} \neq \sigma \in \Sigma$ we write
\[ \partial_\sigma T_\Sigma = \partial_{\Sigma(\sigma)} T_\Sigma \]
and
\[ \partial^0_\sigma T_\Sigma = \partial_\sigma T_\Sigma - \bigcup_{\sigma' \supsetneq \sigma} \partial_{\sigma'} T_\Sigma \]
Sometimes we also write
\[ \partial_{\{0\}}^0T_\Sigma=T. \]
If $\Sigma_0$ is locally finite then the $\partial_\sigma T_\Sigma$ for $\sigma \in \Sigma_0$ form a locally finite closed cover of $\partial_{\Sigma_0} T_\Sigma$.
Set theoretically we have
\[ \partial_\sigma T_\Sigma = \coprod_{\sigma' \in \Sigma(\sigma)} \partial_{\sigma'}^0T_\Sigma \]
and
\[ (\partial_{\Sigma_0}T_\Sigma)_\sigma = \coprod_{\substack{ \sigma' \in \Sigma_0 \\ \sigma' \subset \sigma}} \partial_{\sigma'}^0 T_\Sigma \]
and
\[ T_\sigma = \coprod_{\sigma' \subset \sigma} \partial_{\sigma'}^0 T_\Sigma \]
and
\[ \partial_{\Sigma_0} T_\Sigma = \coprod_{\sigma' \in \Sigma_0} \partial_{\sigma'}^0T_\Sigma. \]
If $\dim \sigma=1$ then $\partial_\sigma^0 T_\Sigma= \partial T_\sigma$.

Keep the notation of the previous paragraph.
We define $S(\sigma)$ to be the split torus with co-character group $X_*(S)$ divided by the subgroup generated by $\sigma \cap X_*(S)$, and $T(\sigma)$ to be the push-out of $T$ to $S(\sigma)$. We also define $\barSigma(\sigma)$  to be the set of images in $X_*(S(\sigma))_\R$ of elements of $\Sigma(\sigma)$. It is a fan for $X_*(S)_\R/\langle \sigma \rangle_\R$. [The main point to check is that if $\tau, \tau' \in \Sigma(\sigma)$ then $(\tau \cap \tau') +\langle \sigma \rangle_\R= (\tau+\langle \sigma \rangle_\R) \cap (\tau' + \langle \sigma \rangle_\R)$. To verify this suppose that $x \in \tau$ and $y \in \tau'$ with $x-y \in \langle \sigma \rangle_\R$. Then $x-y=z-w$ with $z ,w \in \sigma$. Thus $x+w = y+z \in \tau \cap \tau'$ and $x+\langle \sigma \rangle_\R = (x+w) +\langle \sigma \rangle_\R$.] If $\sigma \in \Sigma_0$ we will sometimes write $\barSigma_0(\sigma)$ for $\barSigma(\sigma)$, as it depends only on $\Sigma_0$ and not on $\Sigma$. Then
\[ \partial^0_\sigma T_\Sigma \cong T(\sigma) \subset T(\sigma)_{\barSigma(\sigma)} \cong \partial_\sigma T_\Sigma. \]
Thus $\partial_\sigma T_\Sigma$ is separated, locally (on the source) of finite type and flat over $Y$. The closed subscheme $\partial_\sigma T_\Sigma$ has codimension in $T_\Sigma$ equal to the dimension of $\sigma$. If $\Sigma(\sigma)$ is smooth then $\partial_\sigma T_\Sigma$ is smooth over $Y$. 

If $\Sigma(\sigma)$ is open then $\partial_\sigma T_\Sigma$ satisfies the valuative criterion of properness over $Y$. 
If in addition $\Sigma(\sigma)$ is finite then $\partial_\sigma T_\Sigma$ is proper over $Y$. 
If $\Sigma_0$ is open, then $\partial_{\Sigma_0} T_\Sigma$ satisfies the valuative criterion of properness over $Y$. If in addition $\Sigma_0$ is finite then $\partial_{\Sigma_0}T_\Sigma$ is proper over $Y$.

The schemes $\partial_{\sigma_1} T_\Sigma,...,\partial_{\sigma_s} T_\Sigma$ intersect if and only if $\sigma_1, ..., \sigma_s$ are all contained in some $\sigma \in \Sigma$. In this case the intersection equals $\partial_\sigma T_\Sigma$ for the smallest such $\sigma$. We set
\[ \partial_i T_\Sigma = \coprod_{\dim \sigma =i} \partial_\sigma T_\Sigma. \]

If $Y$ is irreducible then $T_\Sigma$ and each $\partial_\sigma T_\Sigma$ is irreducible. Moreover the irreducible components of $\partial T_\Sigma$ are the $\partial_\sigma T_\Sigma$ as $\sigma$ runs over one dimensional elements of $\Sigma$. 
If $\Sigma$ is smooth then we see that $\cS(\partial T_\Sigma)$ is the delta complex with cells in bijection with the elements of $\Sigma-\{ \{0\}\}$ and with the same `face relations'. In particular it is in fact a simplicial complex and
\[ |\cS(\partial T_\Sigma)|=(|\Sigma|-\{0\})/\R^{\times}_{>0}. \]

We say that $(\Sigma',\Sigma_0')$ {\em refines} $(\Sigma,\Sigma_0)$ if $\Sigma'$ refines $\Sigma$ and $\Sigma'-\Sigma_0'$ is the set of elements of $\Sigma'$ contained in some element of $\Sigma-\Sigma_0$. In this case $\partial_{\Sigma_0'}T_{\Sigma'}$ maps to $\partial_{\Sigma_0} T_{\Sigma}$, and in fact set theoretically $\partial_{\Sigma_0'}T_{\Sigma'}$ is the pre-image of $\partial_{\Sigma_0} T_{\Sigma}$ in $T_{\Sigma'}$. 

If $\Sigma$ is a fan, then by {\em line bundle data for $\Sigma$} we mean a  continuous function $\psi:|\Sigma| \ra \R$, such that for each cone $\sigma \in \Sigma$, the restriction $\psi|_{\sigma}$ equals some $\psi_{\sigma} \in X^*(S)$.
To $\psi$ we can attach a line bundle $\cL_\psi$ on $T_\Sigma$: On $T_{\sigma}$ (with $\sigma \in \Sigma$) it corresponds to the  $\bigoplus_{\chi \in \sigma^\vee \cap X^*(S)} \cL_T(\chi)$-module
\[ \bigoplus_{\substack{ \chi \in X^*(S) \\ \chi-\psi \geq 0 \,\,{\rm on}\,\, \sigma}} \cL_T(\chi). \]
Note that there are natural isomorphisms
\[ \cL_\psi \otimes \cL_{\psi'}\cong \cL_{\psi+\psi'}, \]
and that
\[ \cL_\psi^{\otimes -1} \cong \cL_{-\psi}. \]

We have the following examples of line bundle data.
\begin{enumerate}
\item $\cO_{T_\Sigma}$ is the line bundle associated to $\psi\equiv 0$. 

\item If $\Sigma$ is smooth then $\cI_{\partial T_\Sigma}$ is the line bundle associated to the unique such function $\psi_\Sigma$ which for every one dimensional cone $\sigma \in \Sigma$ satisfies
\[ \psi_{\Sigma}(X_*(S) \cap \sigma) = \Z_{\geq 0}. \]

\end{enumerate}

Suppose that $\alpha:S \onto S'$ is a surjective map of split tori over $Y$. Then $X^*(\alpha): X^*(S') \into X^*(S)$ and $X_*(\alpha):X_*(S) \ra X_*(S')$, the latter with finite cokernel. We call fans $\Sigma$ for $X_*(S)$ and $\Sigma'$ for $X_*(S')$ {\em compatible} if for all $\sigma \in \Sigma$ the image $X_*(\alpha) \sigma$ is contained in some element of $\Sigma'$. In this case the map $\alpha:T \ra \alpha_* T$ extends to an $S$-equivariant map
\[ \alpha:T_\Sigma \lra (\alpha_*T)_{\Sigma'}. \] 
We will write
\[ \Omega_{T_\Sigma/(\alpha_* T)_{\Sigma'}}^1(\log \infty) = \Omega_{T_\Sigma/(\alpha_* T)_{\Sigma'}}^1(\log \cM_\Sigma/\cM_{\Sigma'}). \]
If for all $\sigma' \in \Sigma'$ the pre-image $X_*(\alpha)^{-1}(\sigma')$ is a finite union of elements of $\Sigma$, then $\alpha: T_\Sigma \ra (\alpha_*T)_{\Sigma'}$ is proper.

If $\alpha$ is an isogeny, if $\Sigma$ and $\Sigma'$ are compatible, and if every element of $\Sigma'$ is a finite union of elements of $\Sigma$, then we call $\Sigma$ a {\em quasi-refinement} of $\Sigma'$. In that case the map $T_{\Sigma} \ra T_{\Sigma'}$ is proper.

\begin{lem}\label{logs} If $\alpha$ is surjective and $\# \coker X_*(\alpha)$ is invertible on $Y$ then $\alpha:(T_\Sigma,\cM_\Sigma) \ra ((\alpha_*T)_{\Sigma'},\cM_{\Sigma'})$ is log smooth, and there is a natural isomorphism
\[ (X^*(S)/X^*(\alpha)X^*(S')) \otimes_\Z \cO_{T_\Sigma} \liso \Omega^1_{T_\Sigma/(\alpha_*T)_{\Sigma'}}(\log \infty). \]
\end{lem}

\pfbegin
We can work Zariski locally on $T_\Sigma$. Thus we map replace $T_\Sigma$ by $T_\sigma$ and $(\alpha_* T)_{\Sigma'}$ by $(\alpha_*T)_{\sigma'}$ for cones $\sigma$ and $\sigma'$ with $X_*(\alpha)\sigma \subset \sigma'$. We may also replace $Y$ by an affine open subset $U$ so that $T|_U$ is trivial, i.e. each $\cL_T(\chi) \cong \cO_Y$ compatibly with $ \cL_T(\chi) \otimes \cL_T(\chi') \iso \cL_T(\chi+\chi')$. Then the log structure on $T_\sigma$ has a chart $\Z[\sigma^\vee \cap X^*(S)] \ra \cO_{T_{\sigma}}$ sending $\chi$ to 
\[ 1 \in \cO_Y(Y) \cong \cL_T(\chi). \]
Similarly the log structure on $(\alpha_*T)_{\sigma'}$ has a chart $\Z[(\sigma')^\vee \cap X^*(S')] \ra \cO_{(\alpha_*T)_{\sigma'}}$ sending $\chi$ to 
\[ 1 \in \cO_Y(Y) \cong \cL_{\alpha_*T}(\chi). \]
The lemma follows because
\[ X^*(\alpha):X^*(S') \lra X^*(S) \]
is injective and the torsion subgroup of the kernel is finite with order invertible on $Y$.
\pfend

We will call pairs $(\Sigma,\Sigma_0)$ and $(\Sigma',\Sigma_0')$ of fans and boundary data for $S$ and $S'$, respectively, {\em compatible} if $\Sigma$ and $\Sigma'$ are compatible and if no cone of $\Sigma_0$ maps into any cone of $\Sigma'-\Sigma_0'$. In this case 
\[ \partial_{\Sigma_0} T_\Sigma \lra \partial_{\Sigma_0'} (\alpha_*T)_{\Sigma'}. \]
We will call them {\em strictly compatible} if they are compatible and
$\Sigma-\Sigma_0$ is the set of cones in $\Sigma$ mapping into some element of $\Sigma'-\Sigma_0'$. 

\begin{lem}\label{compat} Suppose that $\alpha:S \onto S'$ is a surjective map of split tori, that $T/Y$ is an $S$-torsor, and that $(\Sigma,\Sigma_0)$ and $(\Sigma',\Sigma_0')$ are strictly compatible fans with boundary data for $S$ and $S'$ respectively. Then 
locally on $T_\Sigma$ there is a strictly positive integer $m$ such that
\[ \alpha^* \cI_{\partial_{\Sigma_0'}(\alpha_* T)_{\Sigma'}} \supset \cI_{\partial_{\Sigma_0} T_\Sigma}^m \]
and
\[ \cI_{\partial_{\Sigma_0} T_\Sigma} \supset (\alpha^* \cI_{\partial_{\Sigma_0'}(\alpha_* T)_{\Sigma'}})^m. \]
\end{lem}

\pfbegin We may work locally on $Y$ and so we may suppose that $Y=\Spec A$ is affine and that each $\cL_T(\chi)$ is trivial.
It also suffices to check the lemma locally on $T_{\Sigma}$. Thus we may suppose that $\Sigma$ consists of a cone $\sigma$ and all its faces. Let $\sigma'$ denote the smallest element of $\Sigma'$ containing the image of $\sigma$. Then we may further suppose that $\Sigma'$ consists of $\sigma'$ and all its faces. We may further suppose that $\sigma \in \Sigma_0$ and $\sigma' \in \Sigma_0'$, else there is nothing to prove.

Then 
\[ T_\Sigma = \Spec \bigoplus_{\chi \in X^*(S) \cap \sigma^\vee} \cL_T(\chi) \]
and $\partial_{\Sigma_0} T_\Sigma$ is defined by
\[ \bigoplus_{\chi \in X^*(S) \cap|\Sigma_0|^{\vee,0}} \cL_T(\chi). \]
Moreover $T_\Sigma \times_{(\alpha_* T)_{\Sigma'}} \partial_{\Sigma_0'} (\alpha_*T)_{\Sigma'}$ is defined by
\[ \bigoplus_{\substack{\chi_1 \in X^*(S') \cap|\Sigma_0'|^{\vee,0} \\ \chi_2 \in X^*(S) \cap \sigma^\vee}} \cL_T(X^*(\alpha)\chi_1 + \chi_2). \]
Thus it suffices to show that for some positive integer $m$ we have
\[ \begin{array}{rcl} (1/m)(X^*(S) \cap|\Sigma_0|^{\vee,0}) &\supset &X^*(\alpha)(X^*(S') \cap |\Sigma_0'|^{\vee,0}) + (X^*(S) \cap \sigma^\vee)\\ & \supset & m (X^*(S) \cap|\Sigma_0|^{\vee,0}). \end{array} \]
This is equivalent to
\[ |\Sigma_0|^{\vee,0}= X^*(\alpha)|\Sigma_0'|^{\vee,0} + \sigma^\vee. \]

Suppose that $\chi_1 \in |\Sigma_0'|^{\vee,0}$ and $\chi_2 \in \sigma^\vee$. Then 
\[ X^*(\alpha)(\chi_1) (\sigma - |\Sigma-\Sigma_0|) = \chi_1(X_*(\alpha)(\sigma - |\Sigma-\Sigma_0|))\subset \chi_1(\sigma'-|\Sigma'-\Sigma_0'|) \subset \R_{>0} \]
and so 
\[ (X^*(\alpha)(\chi_1)+\chi_2) (\sigma - |\Sigma-\Sigma_0|) \subset \R_{>0}. \]
Thus
\[ |\Sigma_0|^{\vee,0}\supset X^*(\alpha)|\Sigma_0'|^{\vee,0} + \sigma^\vee. \]

Conversely suppose that $\chi \in |\Sigma_0|^{\vee,0}$. Let $\tau$ denote the face of $\sigma$, where $\chi=0$. Then $\tau \in \Sigma-\Sigma_0$. Let $\tau'$ denote the smallest face of $\sigma'$ containing $X_*(\alpha) \tau$. Then $\tau' \in \Sigma'-\Sigma_0'$. We can find $\chi_1 \in |\Sigma_0'|^{\vee,0}$ with $\chi_1(\tau')=\{0\}$. 
Note that if $a \in \sigma$ and $\chi(a)=0$ then $(X^*(\alpha)(\chi_1))(a)=0$. Thus 
we can find $\epsilon >0$ such that
\[ \chi - X^*(\alpha)(\epsilon \chi_1) \in \sigma^\vee. \]
It follows that 
\[ |\Sigma_0|^{\vee,0}\subset X^*(\alpha)|\Sigma_0'|^{\vee,0} + \sigma^\vee. \]
The lemma follows.
\pfend

Suppose that $(\Sigma,\Sigma_0)$ and $(\Sigma',\Sigma_0')$ are strictly compatible. We will say that
\begin{itemize}
\item $\Sigma_0$ is {\em open over} $\Sigma_0'$ if $|\Sigma_0|^0$ is open in $X_*(\alpha)^{-1} |\Sigma_0'|^0$;
\item and that $\Sigma_0$ is {\em finite over} $\Sigma_0'$ if only finitely many elements of $\Sigma_0$ map into any element of $\Sigma_0'$.
\end{itemize}
If $\alpha$ is an isogeny, if $\Sigma$ is a quasi-refinement of $\Sigma'$, and if $(\Sigma,\Sigma_0)$ and $(\Sigma',\Sigma_0')$ are strictly compatible, then we call $(\Sigma,\Sigma_0)$ a {\em quasi-refinement} of $(\Sigma',\Sigma_0')$. In this case $\Sigma_0$ is open and finite over $\Sigma_0'$. 

\begin{lem}\label{prop} Suppose that $\alpha:S \onto S'$ is a surjective map of split tori, that $T/Y$ is an $S$-torsor, and that $(\Sigma,\Sigma_0)$ and $(\Sigma',\Sigma_0')$ are strictly compatible fans with boundary data for $S$ and $S'$ respectively. If $\Sigma_0$ is locally finite and $\Sigma_0$ is open over $\Sigma_0'$ then
\[ \partial_{\Sigma_0} T_\Sigma \lra \partial_{\Sigma_0'} T_\Sigma \]
satisfies the valuative criterion of properness. If in addition $\Sigma_0$ is finite over $\Sigma_0'$ then this morphism is proper.
\end{lem}

\pfbegin
It suffices to show that if $\sigma \in \Sigma_0$ and if $\sigma'$ is the smallest element of $\Sigma'_0$ containing $X_*(\alpha) \sigma$, then 
\[ \partial_{\sigma} T_\Sigma \lra \partial_{\sigma'} (\alpha_*T)_{\Sigma'} \]
satisfies the valuative criterion of properness. However this is the map of toric varieties
\[ T(\sigma)_{\barSigma_0(\sigma)} \lra (\alpha_*T)(\sigma')_{\barSigma_0'(\sigma')}. \]
As $\barSigma_0(\sigma)$ is finite, it suffices to check that
\[ \bigcup_{\substack{ \tau' \supset \sigma' \\ \tau' \in \Sigma_0'}} X_*(\alpha)^{-1}((\tau')^0+ \langle \sigma' \rangle_\R) = \bigcup_{\substack{\tau \supset \sigma \\ \tau \in \Sigma_0}} (\tau^0+\langle \sigma \rangle_\R). \]

Choose a point $P \in \sigma^0$ such that 
\[ X_*(\alpha)P \in (X_*(\alpha)\sigma)^0 \subset (\sigma')^0. \]
Then
\[ \langle \sigma' \rangle_R=\sigma' + \R X_*(\alpha)(P). \]
[To see this choose non-zero vectors $v_i$ in each one dimensional face of $\sigma'$. Then we can write $X_*(\alpha)(P)=\sum_i \mu_i v_i$ with each $\mu_i > 0$. If $\lambda_i \in \R$, then for $\lambda$ sufficiently large $\lambda_i+\lambda \mu_i \in \R_{>0}$ for all $i$, and so
\[ \sum_i \lambda_i v_i =  \sum_i (\lambda_i + \lambda \mu_i)v_i -\lambda X_*(\alpha)(P) \in  \sigma' + \R X_*(\alpha)(P). ]\]
Thus
\[ \langle \sigma'\rangle_\R = \sigma' + X_*(\alpha)\langle \sigma \rangle_\R.  \]
Hence for all $\tau' \in \Sigma_0'$ with $\tau' \supset \sigma'$, we have
\[ (\tau')^0+ \langle \sigma'\rangle_\R = (\tau')^0 + X_*(\alpha)\langle \sigma\rangle_\R \]
and so
\[ X_*(\alpha)^{-1}((\tau')^0+\langle \sigma' \rangle_\R)= \langle \sigma\rangle_\R+ X_*(\alpha)^{-1} (\tau')^0. \]
We deduce that it suffices to check that
\[ \langle \sigma\rangle_\R + \bigcup_{\substack{ \tau' \supset \sigma' \\ \tau' \in \Sigma_0'}} X_*(\alpha)^{-1}(\tau')^0 = \langle \sigma \rangle_\R +\bigcup_{\substack{\tau \supset \sigma \\ \tau \in \Sigma_0}} \tau^0. \]
The left hand side certainly contains the right hand side, so it suffices to prove that for all $\tau' \in \Sigma_0'$ with $\tau' \supset X_*(\alpha)\sigma$ we have
\[ \langle \sigma\rangle_\R + X_*(\alpha)^{-1} \tau' \subset \langle \sigma \rangle_\R + 
\bigcup_{\substack{\tau \supset \sigma \\ \tau \in \Sigma_0}} \tau^0. \]

Let $\pi$ denote the map
\[ \pi: X_*(S)_\R \onto X_*(S)_\R/\langle \sigma \rangle_\R. \]
Because $X_*(\alpha)^{-1} \tau'$ and $\bigcup_{\substack{\tau \supset \sigma \\ \tau \in \Sigma_0}} \tau^0$ are invariant under the action of $\R^\times_{>0}$ it suffices to find an open set $U \subset X_*(S)_\R$ containing $P$ such that
\[ (\pi U) \cap \pi X_*(\alpha)^{-1} \tau' \subset \pi \bigcup_{\substack{\tau \supset \sigma \\ \tau \in \Sigma_0}} \tau^0, \]
or equivalently such that
\[ U \cap (\langle \sigma\rangle_\R+ X_*(\alpha)^{-1} \tau') \subset \langle \sigma \rangle_\R + \bigcup_{\substack{\tau \supset \sigma \\ \tau \in \Sigma_0}} \tau^0. \]

Thus it suffices to find an open set $U \subset X_*(S)_\R$ containing $P$ such that
\begin{enumerate}
\item $U \cap X_*(\alpha)^{-1}|\Sigma_0'|^0 \subset \bigcup_{\substack{\tau \supset \sigma \\ \tau \in \Sigma_0}} \tau^0$;
\item $U \cap X_*(\alpha)^{-1} \tau' \subset X_*(\alpha)^{-1} |\Sigma_0'|^0$;
\item and for all open $U' \subset U$ containing $P$ we have $U' \cap (\langle \sigma \rangle_\R + X_*(\alpha)^{-1} \tau') =  U' \cap X_*(\alpha)^{-1} \tau'$.
\end{enumerate}
Moreover in order to find such a $U \ni P$ it suffices to find one satisfying each property independently and take their intersection.

One can find an open set $U \ni P$ satisfying the first property because 
\[ \bigcup_{\substack{\tau \supset \sigma \\ \tau \in \Sigma_0}} \tau^0 \subset |\Sigma_0|^0 \subset X^*(\alpha)^{-1}|\Sigma_0'|^0 \]
are both open inclusions. 

To find $U \ni P$ satisfying the second condition we just need to avoid the faces of $X_*(\alpha)^{-1}\tau'$ which do not contain $P$. 

It remains to check that we can find an open $U \ni P$ satisfying the last condition. 
Suppose that $X_*(\alpha)^{-1}\tau'$ is defined by inequalities $\chi_i \geq 0$ for $i=1,...,r$ with $\chi_i \in X^*(S)_\R$. Suppose that $\chi_i=0$ on $\sigma$ for $i=1,...,s$, but that $\chi_i(P)>0$ for $i=s+1,...,r$. It suffices to choose $U$ so that $\chi_i>0$ on $U$ for $i=s+1,...,r$. For then if $x \in X_*(\alpha)^{-1} \tau'$ and $y \in \langle \sigma \rangle_\R$ with $x+y \in U$ we see that 
\[ \chi_i(x+y)=\chi_i(x) \geq 0 \]
for $i=1,...,s$, while $\chi_i(x+y)>0$ for $i=s+1,...,r$. Thus for $U' \subset U$ we have
\[ U' \cap (\langle \sigma \rangle_\R + X_*(\alpha)^{-1}\tau' ) = U' \cap X_*(\alpha)^{-1} \tau', \]
as desired.
\pfend

By a {\em partial fan} we will mean a collection $\Sigma_0$ of rational polyhedral cones satisfying
\begin{itemize}
\item $(0) \not\in \Sigma_0$;
\item if $\sigma_1, \sigma_2 \in \Sigma_0$, then $\sigma_1 \cap \sigma_2$ is a face of $\sigma_1$ and of $\sigma_2$;
\item if $\sigma_1, \sigma_2 \in \Sigma_0$, and if $\sigma \supset \sigma_2$ is a face of $\sigma_1$, then $\sigma \in \Sigma_0$.
\end{itemize}
(Again note that $\Sigma_0$ may not be closed under taking faces.)
In this case we will let $\tSigma_0$ denote the set of faces of elements of $\Sigma_0$. Then $\tSigma_0$ and $\tSigma_0-\Sigma_0$ are fans, and $\Sigma_0$ is boundary data for $\tSigma_0$. [For suppose that $\tau_i$ is a face of $\sigma_i \in \Sigma_0$ for $i=1,2$. Then $\sigma_1 \cap \sigma_2$ is a face of $\sigma_1$ and so $\tau_1 \cap \sigma_2 =\tau_1 \cap (\sigma_1 \cap \sigma_2)$ is a face of $\sigma_1 \cap \sigma_2$ and hence of $\sigma_2$. Thus $\tau_1 \cap \tau_2 = \tau_2 \cap (\tau_1 \cap \sigma_2)$ is a face of $\tau_2$.] If $\Sigma$ is a fan and $\Sigma_0$ is boundary data for $\Sigma$, then $\Sigma_0$ is a partial fan, and $\Sigma \supset \tSigma_0$. Thus
\[ \partial_{\Sigma_0} T_\Sigma \cong \partial_{\Sigma_0} T_{\tSigma_0}. \]

If $\Sigma_0$ and $\Sigma_0'$ are partial fans we will say that $\Sigma_0$ {\em refines} $\Sigma_0'$ if every element of $\Sigma_0$ is contained in an element of $\Sigma_0'$ and if every element of $\Sigma_0'$ is a finite union of elements of $\Sigma_0$. In this case $\tSigma_0$ also refines $\tSigma_0'$.

If $\Sigma_0$ is a partial fan we will set
\[ |\Sigma_0| = \bigcup_{\sigma \in \Sigma_0} \sigma =|\tSigma_0|.  \]
and
\[ |\Sigma_0|^0 = |\Sigma_0|-|\tSigma_0-\Sigma_0|.  \]
We will call $\Sigma_0$ 
\begin{itemize}
\item {\em smooth} if each $\sigma \in \Sigma_0$ is smooth;
\item {\em full} if every element of $\Sigma_0$ which is not a face of any other element of $\Sigma_0$, has the same dimension as $S$;
\item {\em open} if $|\Sigma_0|^0$ is open in $X_*(S)_\R$;
\item {\em finite} if it has finite cardinality;
\item {\em locally finite} if for every rational polyhedral cone $\tau \subset |\Sigma_0|$ (not necessarily in $\Sigma_0$) the intersection $\tau \cap |\Sigma_0|^0$ meets only finitely many elements of $\Sigma_0$. 
\end{itemize}
If $\Sigma_0$ is smooth, so is $\tSigma_0$.

Suppose that $\Sigma_0$ is a partial fan. If $\Sigma \supset \tSigma_0$ is a fan then the natural maps
\[ \partial_{\Sigma_0} T_{\tSigma_0} \lra \partial_{\Sigma_0} T_\Sigma \]
and
\[ (T_{\tSigma_0})^\wedge_{\partial_{\Sigma_0}T} \lra (T_{\Sigma})^\wedge_{\partial_{\Sigma_0}T} \]
are isomorphisms, and we will denote these schemes/formal schemes $\partial_{\Sigma_0}T$ and $T_{\Sigma_0}^\wedge$ respectively.
Moreover the log structures induced on $T_{\Sigma_0}^\wedge$ by $\cM_{\tSigma_0}$ and by $\cM_\Sigma$ are the same and we will denote them $\cM^\wedge_{\Sigma_0}$. If $\Sigma_0' \subset \Sigma_0$ is also a partial fan, then $T_{\Sigma_0'}^\wedge$ can be identified with the completion of $T_{\Sigma_0}^\wedge$ along $\partial_{\Sigma_0'}T$, and $\cM_{\Sigma_0}^\wedge$ induces $\cM_{\Sigma_0'}^\wedge$. 
 If $\sigma \in \tSigma_0$ then we will let 
\[ (T_{\Sigma_0}^\wedge)_\sigma \]
denote the restriction of $T_{\Sigma_0}^\wedge$ to the topological space $(\partial_{\Sigma_0} T_{\tSigma_0})_\sigma$. Thus the $(T_{\Sigma_0}^\wedge)_\sigma$ for $\sigma \in \Sigma_0$ form an affine open cover of $T_{\Sigma_0}^\wedge$. We have
\[ (T_{\Sigma_0}^\wedge)_{\{ 0\}}= \emptyset \]
and
\[ (T_{\Sigma_0}^\wedge)_{\sigma_1} \cap (T_{\Sigma_0}^\wedge)_{\sigma_2} = (T_{\Sigma_0}^\wedge)_{\sigma_1 \cap \sigma_2}. \]
If $\Sigma_0'$ refines $\Sigma_0$ then there is an induced map
\[ T^\wedge_{\Sigma_0'} \lra T^\wedge_{\Sigma_0}. \] 

We will call $\Sigma_1 \subset \Sigma_0$ {\em boundary data} if, whenever $\sigma \in \Sigma_0$ contains $\sigma' \in \Sigma_1$, then $\sigma \in \Sigma_1$. In this case $\Sigma_1$ is a partial fan and $T_{\Sigma_1}^\wedge$ is canonically identified with the completion of $T_{\Sigma_0}^\wedge$ along $\partial_{\Sigma_1} T_{\tSigma_0}$. 

We will also use the following notation.
\begin{itemize}
\item $\cO_{T_{\Sigma_0}^\wedge}$ will denote the structure sheaf of $T_{\Sigma_0}^\wedge$.
\item $\cI_{T_{\Sigma_0}^\wedge}$ will denote the completion of $\cI_{\partial_{\Sigma_0} T_{\tSigma_0}}$, an ideal of definition for $T_{\Sigma_0}^\wedge$.
\item $\cI_{\partial, \Sigma_0}^\wedge$ will denote the completion of $\cI_{\partial T_{\tSigma_0}}$. Thus $\cI_{\partial_{\Sigma_0} T_{\tSigma_0}} \supset \cI_{\partial ,\Sigma_0}^\wedge$. 
\item  $\Omega^1_{T_{\Sigma_0}^\wedge/\Spf R_0}(\log \infty)$ will denote $\Omega^1_{T_{\Sigma_0}^\wedge/\Spf R_0}(\log \cM_\Sigma^\wedge)$, which is isomorphic to the completion of $\Omega^1_{T_{\tSigma_0}/\Spec R_0}(\log \infty)$. 
\end{itemize}
For $\sigma \in \tSigma_0$ recall that $\cI_{\partial_{\Sigma_0}T_{\tSigma_0}}^m|_{T_\sigma}$ corresponds to the ideal
\[ \bigoplus_{\chi \in \gX_{\Sigma_0,\sigma,m}} \cL_T(\chi) \]
of
\[ \bigoplus_{\chi \in \sigma^\vee \cap X^*(S)} \cL_T(\chi). \]
Also recall that if $\sigma \not\in \Sigma_0$ then 
\[ \gX_{\Sigma_0,\sigma,m} = \sigma^\vee \cap X^*(S) \]
for all $m$, while if $\sigma \in \Sigma_0$ then
\[ \bigcap_m \gX_{\Sigma_0,\sigma,m} = \emptyset. \]

By {\em line bundle data for $\Sigma_0$} we mean a  continuous functions $\psi:|\Sigma_0| \ra \R$, such that for each cone $\sigma \in \Sigma_0$, the restriction $\psi|_{\sigma}$ equals some $\psi_{\sigma} \in X^*(S)$. This is the same as line bundle data for the fan $\tSigma_0$, and we will write $\cL_\psi^\wedge$ for the line bundle on $T_{\Sigma_0}^\wedge$, which is the completion of $\cL_\psi/T_{\tSigma_0}$.
Note that
\[ \cL_\psi^\wedge \otimes \cL_{\psi'}^\wedge= \cL_{\psi+\psi'}^\wedge, \]
and that
\[ (\cL_\psi^\wedge)^{\otimes -1} = \cL_{-\psi}^\wedge. \]

We have the following examples of line bundle data.
\begin{enumerate}
\item $\cO_{T_{\Sigma_0}^\wedge}$ is the line bundle associated to $\psi\equiv 0$. 

\item If $\Sigma_0$ is smooth then $\cI_{\partial, \Sigma_0}^\wedge$ is the line bundle associated to the unique such function $\psi_{\tSigma_0}$ which for every one dimensional cone $\sigma \in \tSigma_0$ satisfies
\[ \psi_{\tSigma_0}(X_*(S) \cap \sigma) = \Z_{\geq 0}. \]

\end{enumerate}

Suppose that $\alpha:S \onto S'$ is a surjective map of tori, and that $\Sigma_0$ (resp. $\Sigma_0'$) is a partial fan for $S$ (resp. $S'$). We call $\Sigma_0$ and $\Sigma_0'$ {\em compatible} if for every $\sigma \in \Sigma_0$ the image $X_*(\alpha) \sigma$ is contained in some element of $\Sigma_0'$ but in no element of $\tSigma_0'-\Sigma_0'$. In this case $(\tSigma_0,\Sigma_0)$ and $(\tSigma_0',\Sigma_0')$ are compatible, and there is a natural morphism
\[ \alpha: (T^\wedge_{\Sigma_0},\cM_{\Sigma_0}^\wedge) \lra ((\alpha_*T)_{\Sigma_0'}^\wedge,\cM_{\Sigma_0'}^\wedge). \]
We will write
\[ \Omega_{T_{\Sigma_0}^\wedge/(\alpha_* T)_{\Sigma_0'}^\wedge}^1(\log \infty) = \Omega_{T_{\Sigma_0}^\wedge/(\alpha_* T)^\wedge_{\Sigma_0'}}^1(\log \cM_{\Sigma_0}^\wedge/\cM^\wedge_{\Sigma_0'}). \]
The following lemma follows immediately from lemma \ref{logs}.

\begin{lem}\label{clogs} If $\alpha$ is surjective and $\# \coker X_*(\alpha)$ is invertible on $Y$ then there is a natural isomorphism
\[ (X^*(S)/X^*(\alpha)X^*(S')) \otimes_\Z \cO_{T_{\Sigma_0}^\wedge} \liso \Omega_{T_{\Sigma_0}^\wedge/(\alpha_* T)_{\Sigma_0'}^\wedge}^1(\log \infty) . \]
\end{lem}

We will call $\Sigma_0$ and $\Sigma_0'$ {\em strictly compatible} if they are compatible and if an element of $\tSigma_0$ lies in $\Sigma_0$ if and only if it maps to no element of $\tSigma_0'-\Sigma_0'$. In this case $(\tSigma_0,\Sigma_0)$ and $(\tSigma_0',\Sigma_0')$ are strictly compatible. We will say that
\begin{itemize}
\item $\Sigma_0$ is {\em open over} $\Sigma_0'$ if $|\Sigma_0|^0$ is open in $X_*(\alpha)^{-1} |\Sigma_0'|^0$;
\item and that $\Sigma_0$ is {\em finite over} $\Sigma_0'$ if only finitely many elements of $\Sigma_0$ map into any element of $\Sigma_0'$.
\end{itemize}
If $\alpha$ is an isogeny, if $\Sigma_0$ and $\Sigma_0'$ are strictly compatible, and if every element of $\Sigma_0'$ is a finite union of elements of $\Sigma_0$, then we call $\Sigma_0$ a {\em quasi-refinement} of $\Sigma_0'$. In this case $\Sigma_0$ is open and finite over $\Sigma_0'$. 
The next lemma follows immediately from lemma \ref{compat} and \ref{prop}. 

\begin{lem}\label{proper} Suppose that $\Sigma_0'$ and $\Sigma_0$ are strictly compatible.
\begin{enumerate}
\item $T_{\Sigma_0}^\wedge$ is the formal completion of $T_{\tSigma_0}$ along $\partial_{\Sigma_0'}(\alpha_* T)$, and $T_{\Sigma_0}^\wedge$ is locally (on the source) topologically of finite type over $(\alpha_*T)^\wedge_{\Sigma_0'}$.
\item If $\Sigma_0$ is locally finite and if it is open and finite over $\Sigma_0'$ then $T_{\Sigma_0}^\wedge$ is proper over $(\alpha_* T)_{\Sigma_0'}$.
\end{enumerate} \end{lem}

\begin{cor} If $\alpha$ is an isogeny, if $\Sigma_0$ is locally finite, and if $\Sigma_0$ is a quasi-refinement of $\Sigma_0'$ then $T_{\Sigma_0}^\wedge$ is proper over $(\alpha_* T)_{\Sigma_0'}$.\end{cor}

If $\Sigma_0$ and $\Sigma_0'$ are compatible partial fans and if $\Sigma_1' \subset \Sigma_0'$ is boundary data then $\Sigma_0(\Sigma_1')$ will denote the set of elements $\sigma \in \Sigma_0$ such that $X_*(\alpha)\sigma$ is contained in no element of $\Sigma_0'-\Sigma_1'$. It is boundary data for $\Sigma_0$. Moreover the formal completion of $T_{\Sigma_0}^\wedge$ along the reduced sub-scheme of $(\alpha_*T)_{\Sigma_1'}^\wedge$ is canonically identified with $T^\wedge_{\Sigma_0(\Sigma_1')}$. If $\Sigma_1'=\{ \sigma'\}$ is a singleton we will write $\Sigma_0(\sigma')$ for $\Sigma_0(\{ \sigma'\})$.

\newpage \subsection{Cohomology of line bundles.}

In this section we will compute the cohomology of line bundles on formal completions of torus embeddings. {\em We will work throughout over a base scheme $Y$ which is connected, separated, and flat and locally of finite type over $\Spec R_0$.}

We start with some definitions. We continue to assume that $S/Y$ is a split torus, that $T/Y$ is an $S$-torsor, that $\Sigma_0$ is a partial fan and that $\psi$ is line bundle data for $\Sigma_0$. If $\sigma \in \tSigma_0$ then we set
\[ \gX_{\Sigma_0,\psi,\sigma,0} = \{ \chi \in X^*(S) \cap \sigma^\vee: \,\, \chi \geq \psi \,\, {\rm on} \,\, \sigma \}. \]
For $m>0$ we define $\gX_{\Sigma_0,\psi,\sigma,m}$ to be the set of sums of an element of $\gX_{\Sigma_0,\psi,\sigma,0}$ and an element of $\gX_{\Sigma_0,\sigma,m}$. If $\sigma \not\in \Sigma_0$ then
\[ \gX_{\Sigma_0,\psi,\sigma,m} = \gX_{\Sigma_0,\psi,\sigma,0} \]
for all $m$, while if $\sigma \in \Sigma_0$ 
\[ \bigcap_m \gX_{\Sigma_0,\psi,\sigma,m} = \emptyset. \]
Further suppose that $\chi \in X^*(S)$.
\begin{itemize}
\item Set $Y_\psi(\chi)=\{ x \in X^*(S)_\R
:\,\, (\psi-\chi)(x)>0\} $.
\item If $U \subset Y$ is open let $H^j_{\Sigma_0,\psi,m}(\chi)(U)$ denote the $j^{th}$ cohomology of the Cech complex
with $i^{th}$ term
\[ \prod_{\substack{(\sigma_0,...,\sigma_i) \in \Sigma_0^{i+1} \\ 
\chi \in \gX_{\Sigma_0,\psi,\sigma_0 \cap...\cap \sigma_i,0} \\ \chi \not\in \gX_{\Sigma_0,\psi,\sigma_0\cap...\cap \sigma_i,m}}}
\cL_T(\chi)(U). \]
\end{itemize}
Note the examples:
\begin{enumerate}
\item $Y_0(\chi)\cap |\Sigma_0|^0 = \emptyset$ if and only if $\chi \in |\Sigma_0|^\vee$.

\item $Y_{\psi_{\tSigma_0}}(\chi) \cap |\Sigma_0|^0=\emptyset$ if and only if $\chi \in |\Sigma_0|^{\vee,0}$.
\end{enumerate}
Also note that if $\Sigma_0$ is finite then, for $m$ large enough, $H^j_{\Sigma_0,\psi,m}(\chi)(U)$ does not depend on $m$. We will denote it simply $H^j_{\Sigma_0,\psi}(\chi)(U)$. It equals the cohomology of the Cech complex
\[ \prod_{\substack{(\sigma_0,...,\sigma_i) \in \Sigma_0^{i+1} \\ 
\sigma_0 \cap ... \cap \sigma_i \in \Sigma_0}}
\cL_T(\chi)(U). \]

\begin{lem}\label{calc} If $U$ is connected 
then
\[ H_{\Sigma_0,\psi}^i(\chi)(U) = H^i_{|\Sigma_0|^0-Y_\psi(\chi)}(|\Sigma_0|^0,\cL_T(\chi)(U)). \]
\end{lem}

\pfbegin
Write $M$ for $\cL_T(\chi)(U)$.
We follow the argument of section 3.5 of \cite{fulton}. As $\sigma_0 \cap ... \cap \sigma_i \cap |\Sigma_0|^0$ and $\sigma_0 \cap ... \cap \sigma_i \cap |\Sigma_0|^0 \cap Y_\psi(\chi)$ are convex, we see that
\[ \begin{array}{rl} &H^j_{(\sigma_0 \cap ... \cap \sigma_i \cap |\Sigma_0|^0)- Y_\psi(\chi)} (\sigma_0 \cap ... \cap \sigma_i  \cap |\Sigma_0|^0, M) \\ 
=& \left\{ \begin{array}{ll} M & {\rm if}\,\,  j=0 \,\, {\rm and} \,\, (\sigma_0 \cap ... \cap \sigma_i \cap |\Sigma_0|^0)\cap Y_\psi(\chi)= \emptyset \\ (0) & {\rm otherwise}.\end{array} \right. \end{array}\]
(See the first paragraph of section 3.5 of \cite{fulton}.) Thus the $i^{th}$ term of our Cech complex becomes
\[  \prod_{(\sigma_0,...,\sigma_i)\in \Sigma_0^{i+1}} H^0_{(\sigma_0 \cap ... \cap \sigma_i \cap |\Sigma_0|^0)- Y_\psi(\chi)} (\sigma_0 \cap ... \cap \sigma_i \cap |\Sigma_0|^0, M) . \]
Thus it suffices to show that the Cech complex with $i^{th}$ term
\[ \prod_{(\sigma_0,...,\sigma_i )\in \Sigma_0^{i+1}} H^0_{(\sigma_0 \cap ... \cap \sigma_i \cap |\Sigma_0|^0)- Y_\psi(\chi)} (\sigma_0 \cap ... \cap \sigma_i \cap |\Sigma_0|^0, M) \]
computes
\[ H^i_{|\Sigma_0|^0- Y_\psi(\chi)} (|\Sigma_0|^0, M).\]

To this end choose an injective resolution
\[ M \lra \cI^0 \lra \cI^1 \lra ... \]
as sheaves of abelian groups on $|\Sigma_0|^0$, and consider the double complex
\[ \prod_{(\sigma_0,...,\sigma_i) \in \Sigma_0^{i+1}} H^0_{(\sigma_0 \cap ... \cap \sigma_i \cap |\Sigma_0|^0)- Y_\psi(\chi)} (\sigma_0 \cap ... \cap \sigma_i \cap |\Sigma_0|^0, \cI^j). \]
We compute the cohomology of the corresponding total complex in two ways. Firstly 
the $j^{th}$ cohomology of the complex 
\[ \begin{array}{c}  H^0_{(\sigma_0 \cap ... \cap \sigma_i \cap |\Sigma_0|^0)- Y_\psi(\chi)} (\sigma_0 \cap ... \cap \sigma_i \cap |\Sigma_0|^0, \cI^0)\\ \da \\ H^0_{(\sigma_0 \cap ... \cap \sigma_i \cap |\Sigma_0|^0)- Y_\psi(\chi)} (\sigma_0 \cap ... \cap \sigma_i \cap |\Sigma_0|^0, \cI^1)\\ \da \\ \vdots \end{array} \]
equals
\[ H^j_{(\sigma_0 \cap ... \cap \sigma_i \cap |\Sigma_0|^0)- Y_\psi(\chi)} (\sigma_0 \cap ... \cap \sigma_i \cap |\Sigma_0|^0,M). \]
(See theorem 4.1, proposition 5.3 and theorem 5.5 of chapter II of \cite{bredon}.) This vanishes for $j>0$, and so the cohomology of our total complex is the same as the cohomology of the Cech complex with $i^{th}$ term 
\[ \prod_{(\sigma_0,...,\sigma_i) \in \Sigma_0^{i+1}} H^0_{(\sigma_0 \cap ... \cap \sigma_i \cap |\Sigma_0|^0)- Y_\psi(\chi)} (\sigma_0 \cap ... \cap \sigma_i \cap |\Sigma_0|^0, M). \]
Thus it suffices to identify the cohomology of our double complex with 
\[ H^i_{|\Sigma_0|^0- Y_\psi(\chi)} (|\Sigma_0|^0, M).\] 

For this it suffices to show that
\[ \begin{array}{l} (0) \lra H^0_{|\Sigma_0|^0-Y_\psi(\chi)}(|\Sigma_0|^0,\cI^j) \lra \prod_{\sigma_0 \in \Sigma_0} H^0_{\sigma_0 \cap |\Sigma_0|^0- Y_\psi(\chi)} (\sigma_0 \cap |\Sigma_0|^0, \cI^j) \lra \\ \lra \prod_{(\sigma_0,\sigma_1) \in \Sigma_0^2} H^0_{(\sigma_0 \cap \sigma_1 \cap |\Sigma_0|^0)- Y_\psi(\chi)} (\sigma_0 \cap \sigma_1 \cap |\Sigma_0|^0, \cI^j) \lra ... \end{array} \]
is exact for all $j$. Let $\tcI^j$ denote the sheaf of discontinuous sections of $\cI^j$, i.e. $\tcI^j(V)$ denotes the set of functions which assign to each point of $x \in V$ an element of the stalk $\cI^j_x$ of $\cI^j$ at $x$. Then $\cI^j$ is a direct summand of $\tcI^j$ so it suffices to show that
\[ \begin{array}{l} (0) \lra H^0_{|\Sigma_0|^0-Y_\psi(\chi)}(|\Sigma_0|^0,\tcI^j) \lra \prod_{\sigma_0 \in \Sigma_0} H^0_{\sigma_0 \cap |\Sigma_0|^0- Y_\psi(\chi)} (\sigma_0 \cap |\Sigma_0|^0, \tcI^j) \lra \\ \lra \prod_{(\sigma_0,\sigma_1) \in \Sigma_0^2} H^0_{(\sigma_0 \cap \sigma_1 \cap |\Sigma_0|^0)- Y_\psi(\chi)} (\sigma_0 \cap \sigma_1 \cap |\Sigma_0|^0, \tcI^j) \lra ... \end{array} \]
is exact for all $j$. However $\tcI^j$ is the direct product over $x$ in $|\Sigma_0|^0$ of the sky-scraper $\barcI^j_x$ sheaf at $x$ with stalk $\cI^j_x$. Thus it suffices to show
that
\[ \begin{array}{l} (0) \lra H^0_{|\Sigma_0|^0-Y_\psi(\chi)}(|\Sigma_0|^0,\barcI_x^j) \lra \prod_{\sigma_0 \in \Sigma_0} H^0_{\sigma_0 \cap |\Sigma_0|^0- Y_\psi(\chi)} (\sigma_0 \cap |\Sigma_0|^0, \barcI_x^j) \lra \\ \lra \prod_{(\sigma_0,\sigma_1) \in \Sigma_0^2} H^0_{(\sigma_0 \cap \sigma_1 \cap |\Sigma_0|^0)- Y_\psi(\chi)} (\sigma_0 \cap \sigma_1 \cap |\Sigma_0|^0, \barcI_x^j) \lra ... \end{array} \]
is exact for all $x \in |\Sigma_0|^0$ and for all $j$. If $x \in Y_\psi(\chi)\cap |\Sigma_0|^0$ all the terms in this sequence are $0$, so the sequence is certainly exact. If $x \in |\Sigma_0|^0-Y_\psi(\chi)$, this sequence equals
\[ (0) \lra \cI_x^j \lra \prod_{\substack{\sigma_0 \in \Sigma_0 \\ x \in \sigma}} \cI_x^j \lra \prod_{\substack{(\sigma_0,\sigma_1) \in \Sigma_0^2 \\ x \in (\sigma_0 \cap \sigma_1)}} \cI_x^j \lra ... \]

A standard argument shows that this is indeed exact: Choose $\sigma \in \Sigma_0$ with $x \in \sigma$. Suppose
\[ (a(\sigma_0,...,\sigma_i)) \in \ker \left( \prod_{\substack{(\sigma_0,...,\sigma_i) \in \Sigma_0^{i+1} \\ x \in \sigma_0 \cap ... \cap \sigma_i}} \cI_x^j \lra \prod_{\substack{(\sigma_0,...,\sigma_{i+1}) \in \Sigma_0^{i+2} \\ x \in \sigma_0 \cap ... \cap \sigma_{i+1}}} \cI_x^j \right) .\]
Define 
\[ (a'(\sigma_0,...,\sigma_{i-1})) \in \prod_{\substack{(\sigma_0,...,\sigma_{i-1}) \in \Sigma_0^{i} \\ x \in \sigma_0 \cap ... \cap \sigma_{i-1}}} \cI_x^j \]
by
\[ a'(\sigma_0,...,\sigma_{i-1})=a(\sigma_0,...,\sigma_{i-1},\sigma). \]
If $\partial a'$ denotes the image of $a'$ in 
\[ \prod_{\substack{(\sigma_0,...,\sigma_i) \in \Sigma_0^{i+1} \\ x \in \sigma_0 \cap ... \cap \sigma_i}} \cI_x^j \]
then
\[ (\partial a')(\sigma_0,...,\sigma_i) = \sum_{k=0}^i (-1)^k a(\sigma_0,...,\widehat{\sigma_k},...,\sigma_i,\sigma) = (-1)^i a(\sigma_0,...,\sigma_i), \]
i.e. $a=(-1)^i \partial a'$.
\pfend

In general we will let $H^i_{|\Sigma_0|^0-Y_\psi(\chi)}(|\Sigma_0|^0,\cL_T(\chi))$ denote the sheaf of $\cO_Y$-modules on $Y$ associated to the pre-sheaf
\[ U \longmapsto 
H^i_{|\Sigma_0|^0-Y_\psi(\chi)}(|\Sigma_0|^0,\cL_T(\chi)(U)). \] 

\begin{lem} \label{ctorem0} Let $Y$ be a connected, locally noetherian, separated scheme, let $S/Y$ be a split torus, let $T/Y$ be an $S$-torsor, let $\Sigma_0$ be a partial fan for $S$, let $\psi$ be line bundle data for $\Sigma_0$, and let $\pi_{\Sigma_0}^\wedge$ denote the map $T^\wedge_{\Sigma_0} \ra Y$. Suppose that $\Sigma_0$ is finite and open.
Then
\[ R^i\pi_{\Sigma_0,*}^\wedge \cL_\psi^\wedge = \prod_{\chi \in X^*(S)} H^i_{|\Sigma_0|^0-Y_\psi(\chi)}(|\Sigma_0|^0,\cL_T(\chi)). \]

(Note that $R^i\pi_{\Sigma_0,*}^\wedge \cL_\psi^\wedge$ may not be quasi-coherent on $Y$. Infinite products of quasi-coherent sheaves may not be quasi-coherent.)
\end{lem}

\pfbegin
The left hand side is the sheaf associated to the pre-sheaf
\[ U \longmapsto H^i(T_{\Sigma_0}^\wedge|_U,\cL_\psi^\wedge) \]
and the right hand side is the sheaf associated to the pre sheaf
\[ U \longmapsto \prod_{\chi \in X^*(S)(U)}
H^i_{|\Sigma_0|^0-Y_\psi(\chi)}(|\Sigma_0|^0, \cL_T(\chi)(U)). \]
Thus it suffices to establish isomorphisms
\[ H^i(T_{\Sigma_0}^\wedge|_U,\cL_\psi^\wedge) \cong \prod_{\chi \in X^*(S)(U)}
H^i_{|\Sigma_0|^0-Y_\psi(\chi)}(|\Sigma_0|^0,\cL_T(\chi)(U)), \]
compatibly with restriction, for $U=\Spec A$, with $A$ noetherian and $\Spec A$
 connected.

Write $\partial_{\Sigma_0,m} T_{\tSigma_0}$ for the closed subscheme of $T_{\tSigma_0}$ defined by $\cI_{\partial_{\Sigma_0} T_{\tSigma_0}}^m$. It has the same underlying topological space as $\partial_{\Sigma_0}T_{\tSigma_0}$.
We will first compute 
\[ H^i(\partial_{\Sigma_0,m} T_{\tSigma_0}|_U,\cL_\psi/\cI_{\partial_{\Sigma_0} T_{\tSigma_0}}^m \cL_\psi), \]
using the affine cover of $\partial_{\Sigma_0,m} T_{\tSigma_0}$ by the open sets $T_\sigma$ for $\sigma \in \Sigma_0$. This gives rise to a Cech complex with terms
\[ \prod_{(\sigma_0,...,\sigma_i) \in \Sigma_0^{i+1}} \bigoplus_{\substack{ \chi \in X^*(S) \\ \chi  \in \gX_{\Sigma_0,\psi,\sigma_0 \cap...\cap \sigma_i,0} \\ \chi \not\in\gX_{\Sigma_0,\psi,\sigma_0 \cap...\cap \sigma_i,m}}}
\cL_T(\chi)(U). \]
As $\Sigma_0$ is finite, we see that
\[ H^i(\partial_{\Sigma_0,m} T_{\tSigma_0}|_U,\cL_\psi/\cI_{\partial_{\Sigma_0} T_{\tSigma_0}}^m \cL_\psi) = \bigoplus_{\chi \in X^*(S)} H_{\Sigma_0,\psi,m}^i(\chi)(U). \]

Because $A$ is noetherian, because $\partial_{\Sigma_0,m} T_{\tSigma_0}$ is proper over
$\Spec A$ and because $\cL_\psi/\cI_{\partial_{\Sigma_0} T_{\tSigma_0}}^m \cL_\psi$ is a 
coherent sheaf on $\partial_{\Sigma_0,m} T_{\tSigma_0}$, we see
 that the cohomology group $H^i(\partial_{\Sigma_0,m} T_{\tSigma_0}|_U,\cL_\psi/\cI_{\partial_{\Sigma_0} T_{\tSigma_0}}^m \cL_\psi)$ 
is a finitely generated $A$-module, and hence, for fixed $m$ and $i$, we see that the groups
 $H_{\Sigma_0,\psi,m}^i(\chi)(U)=(0)$ for all but finitely many $\chi$. In particular
\[ H^i(\partial_{\Sigma_0,m} T_{\tSigma_0}|_U,\cL_\psi/\cI_{\partial_{\Sigma_0} T_{\tSigma_0}}^m \cL_\psi) = \prod_{\chi \in X^*(S)} H_{\Sigma_0,\psi,m}^i(\chi)(U). \]
Moreover, combining this observation with the fact that $\{ H^i_{\Sigma_0,\psi,m}(\chi)(U)\}$ satisfies the Mittag-Leffler condition, we see that  the system 
\[ \{ H^i(\partial_{\Sigma_0,m} T_{\tSigma_0}|_U,\cL_\psi/\cI_{\partial_{\Sigma_0} T_{\tSigma_0}}^m \cL_\psi)\}\]
 satisfies the Mittag-Leffler condition. Hence from proposition 0.13.3.1 of \cite{ega3} we see that
\[ \begin{array}{rcl} H^i(T_{\Sigma_0}^\wedge|_U,\cL_\psi^\wedge) &\cong &\lim_{\la m} H^i(\partial_{\Sigma_0,m} T_{\tSigma_0}|_U,\cL_\psi/\cI_{\partial_{\Sigma_0} T_{\tSigma_0}}^m \cL_\psi) \\ & \cong & \prod_{\chi \in X^*(S)} \lim_{\la m} H^i_{\Sigma_0,\psi,m}(\chi)(U),
\end{array} \]
and the present lemma follows from lemma \ref{calc}.
\pfend

\begin{lem} \label{ctorem1} Let $Y$ be a connected, locally noetherian, separated  scheme, $S/Y$ be a split torus, let $T/Y$ be an $S$-torsor, let $\Sigma_\infty$ be a partial fan for $S$, let 
\[ \Sigma_1 \subset \Sigma_2 \subset ... \]
be a nested sequence of partial fans with $\Sigma_\infty = \bigcup_i \Sigma_i$ and let $\psi$ be line bundle data for $\Sigma_\infty$. For $i=1,2,3,...,\infty$ let $\pi_{\Sigma_i}^\wedge$ denote the map $T_{\Sigma_i}^\wedge \ra Y$. 

Suppose that for $i\in \Z_{>0}$ the partial fan $\Sigma_i$ is finite and open. Suppose also that for all $i \in \Z_{\geq 0}$ and all connected, noetherian, affine open sets $U \subset Y$, the inverse system 
\[ \{ H^{i}_{|\Sigma_j|^0-Y_\psi(\chi)}(|\Sigma_j|^0,\cO_Y(U))\} \]
 satisfies the Mittag-Leffler condition. Then
\[ R^i\pi_{\Sigma_\infty,*}^\wedge \cL_\psi^\wedge \cong \prod_{\chi \in X^*(S)} \lim_{\la j} H^i_{|\Sigma_j|^0-Y_\psi(\chi)}(|\Sigma_j|^0,\cL_T(\chi)). \]
\end{lem}

\pfbegin
The left hand side is the sheaf associated to the pre-sheaf
\[ U \longmapsto H^i(T_{\Sigma_\infty}^\wedge|_U, \cL_\psi^\wedge) \]
and the right hand side is the sheaf associated to the pre-sheaf
\[ U \longmapsto \prod_{\chi \in X^*(S)} \lim_{\la j} H^i_{|\Sigma_j|^0-Y_\psi(\chi)}(|\Sigma_j|^0,\cO_Y(U)) \otimes \cL_T(\chi)(U). \]
Thus it suffices to establish isomorphisms 
\[ H^i(T_{\Sigma_\infty}^\wedge|_U, \cL_\psi^\wedge) \cong \prod_{\chi \in X^*(S)(Y)} \lim_{\la j} H^i_{|\Sigma_j|^0-Y_\psi(\chi)}(|\Sigma_j|^0,\cO_Y(U)) \otimes \cL_T(\chi)(U), \]
compatibly with restriction, for $U=\Spec A$, with $A$ noetherian and $\Spec A$ connected. 

We can compute $H^i(T_{\Sigma_\infty}^\wedge|_U, \cL_\psi^\wedge)$ as the cohomology of the Cech complex with $i^{th}$ term
\[ \prod_{(\sigma_0,...,\sigma_i) \in \Sigma_\infty^{i+1}} \cL_\psi^\wedge((T_{\Sigma_\infty}^\wedge)_{(\sigma_0 \cap ... \cap \sigma_i)}|_U), \]
and we can compute $H^i(T_{\Sigma_j}^\wedge|_U, \cL_\psi^\wedge)$ as the cohomology of the Cech complex with $i^{th}$ term
\[ \prod_{(\sigma_0,...,\sigma_i) \in \Sigma_j^{i+1}} \cL_\psi^\wedge((T_{\Sigma_j}^\wedge)_{(\sigma_0 \cap ... \cap \sigma_i)}|_U). \]
Note that as soon as the faces of $\sigma$ in $\Sigma_j$ equals the faces of $\sigma$ in $\Sigma_\infty$ then $(T_{\Sigma_\infty}^\wedge)_\sigma=(T_{\Sigma_j}^\wedge)_\sigma$. Thus 
\[ \lim_{\la j} \prod_{(\sigma_0,...,\sigma_i) \in \Sigma_j^{i+1}} \cL_\psi^\wedge((T_{\Sigma_j}^\wedge)_{(\sigma_0 \cap ... \cap \sigma_i)}|_U) \cong \prod_{(\sigma_0,...,\sigma_i) \in \Sigma_\infty^{i+1}} \cL_\psi^\wedge((T_{\Sigma_\infty}^\wedge)_{(\sigma_0 \cap ... \cap \sigma_i)}|_U), \]
and
\[ \left\{ \prod_{(\sigma_0,...,\sigma_i) \in \Sigma_j^{i+1}} \cL_\psi^\wedge((T_{\Sigma_j}^\wedge)_{(\sigma_0 \cap ... \cap \sigma_i)(U)}) \right\} \]
satisfies the Mittag-Leffler condition (with $j$ varying but $i$ fixed). 

From theorem 3.5.8 of \cite{weibel} we see that there is a short exact sequence
\[ \begin{array}{l} (0) \lra {\lim_{\la j}}^1 H^{i-1}(T_{\Sigma_j}^\wedge|_U, \cL_\psi^\wedge)\lra
H^i(T_{\Sigma_\infty}^\wedge|_U, \cL_\psi^\wedge)
\lra \\ \lra \lim_{\la j} H^{i}(T_{\Sigma_j}^\wedge|_U, \cL_\psi^\wedge) \lra (0).\end{array} \]
Applying lemma \ref{ctorem0} and the fact that $\lim_{\la}$ and ${\lim_{\la}}^1$ in the category of abelian groups commute with arbitrary products, the present lemma follows. (It follows easily from 
definition 3.5.1 of \cite{weibel} and the exactness of infinite products in the category of abelian groups, that $\lim_{\la}$ and ${\lim_{\la}}^1$ commute with arbitrary products
in the category of abelian groups.)
\pfend

We now turn to two specific line bundles: $\cO_{T_{\Sigma_0}^\wedge}$ and, in the case that $\Sigma_0$ is smooth, $\cI_{\partial,\Sigma_0}^\wedge$.

\begin{lem}\label{ctorem2} Let $Y$ be a connected, locally noetherian, separated scheme, let $S/Y$ be a split torus, let $T/Y$ be an $S$-torsor, let $\Sigma_0$ be a partial fan for $S$, and let $\pi_{\Sigma_0}^\wedge$ denote the map $T^\wedge_{\Sigma_0} \ra Y$. Suppose that $\Sigma_0$ is finite and open and that $|\Sigma_0|^0$ is convex. 
\begin{enumerate}
\item Then 
\[ R^i\pi^\wedge_{\Sigma_0,*} \cO_{T_{\Sigma_0}^\wedge} = \left\{ \begin{array}{ll} \prod_{\chi \in |\Sigma_0|^\vee}  \cL(\chi)& {\rm if}\,\, i=0
\\ (0) & {\rm otherwise.}\end{array} \right. \]

\item If in addition $\Sigma_0$ is smooth then
\[ R^i\pi_{\Sigma_0,*}^\wedge \cI_{\partial, \Sigma_0}^\wedge = \left\{ \begin{array}{ll} \prod_{\chi \in |\Sigma_0|^{\vee,0}} \cL(\chi) & {\rm if}\,\, i=0
\\ (0) & {\rm otherwise.}\end{array} \right. \]
\end{enumerate} 
\end{lem}

\pfbegin 
The first part follows from lemma \ref{ctorem0} because $Y_0(\chi)\cap |\Sigma_0|^0$ is empty if $\chi \in |\Sigma_0|^\vee$ and otherwise, being the intersection of two convex sets, it is convex.

For the second part we have that $Y_{\psi_{\tSigma_0}}(\chi) \cap |\Sigma_0|^0 =\emptyset$ if and only if $\chi \in |\Sigma_0|^{\vee,0}$. Thus it suffices to show that each $Y_{\psi_{\tSigma_0}}(\chi)\cap |\Sigma_0|^0$ is empty or contractible.

To this end, consider the sets
\[ Y'(\chi)= \bigcup_{\substack{\sigma \in \Sigma_0 \\ \chi \leq 0 \,\,{\rm on}\,\,\sigma}} \sigma  \]
and
\[ Y''(\chi)= \bigcup_{\substack{\sigma \in \Sigma_0 \\ \chi \not> 0 \,\,{\rm on}\,\,\sigma-\{0\}}} \sigma . \]
If $\chi>0$ on $\sigma-\{0\}$ then $\chi \geq \psi_{\tSigma_0}$ on $\sigma$ so that $\sigma \cap Y_{\psi_{\tSigma_0}}(\chi) = \emptyset$.
Thus 
\[ Y''(\chi) \supset Y_{\psi_{\tSigma_0}}(\chi)\cap |\Sigma_0|^0 \supset Y'(\chi) \cap |\sigma_0|^0 \]
and
\[ Y''(\chi) \supset \{ x \in |\Sigma_0|^0:\,\,\chi(x)\leq 0\} \supset Y'(\chi) \cap |\sigma_0|^0. \]
We will describe a deformation retraction
\[ H: Y''(\chi) \times [0,1] \lra Y''(\chi) \]
from $Y''(\chi)$ to $Y'(\chi)$, which restricts to deformation retractions
\[ (Y_{\psi_{\tSigma_0}}(\chi) \cap |\Sigma_0|^0)\times [0,1] \lra Y_{\psi_{\tSigma_0}} \cap |\sigma_0|^0(\chi) \]
from $Y_{\psi_{\tSigma_0}}(\chi)\cap |\Sigma_0|^0$ to $Y'(\chi)\cap |\Sigma_0|^0$, and
\[ \{ x \in |\Sigma_0|^0:\,\,\chi(x)\leq 0\} \times [0,1] \lra \{ x \in |\Sigma_0|^0:\,\,\chi(x)\leq 0\} \]
from $\{ x \in |\Sigma_0|^0:\,\,\chi(x)\leq 0\}$ to $Y'(\chi)\cap |\Sigma_0|^0$. (Recall that in particular $H_{Y'(\chi) \times [0,1]}$ is just projection to the first factor.) As $\{ x \in |\Sigma_0|^0:\,\,\chi(x)\leq 0\}$ is empty or convex, it would follow that $Y_{\psi_{\tSigma_0}}(\chi)\cap |\Sigma_0|^0$ is empty or contractible and the second part of the corollary would follow.

To define $H$ it suffices to define, for each $\sigma \in \tSigma_0$ with $\sigma \subset Y''(\chi)$, a deformation retraction
\[ H_\sigma: \sigma \times [0,1] \lra \sigma \]
from $\sigma$ to $\sigma \cap Y'(\chi)$ with the following properties:
\begin{itemize}
\item If $\sigma' \subset \sigma$ then $H_\sigma|_{\sigma' \times [0,1]}=H_{\sigma'}$.
\item $H_\sigma|_{(\sigma \cap Y_{\psi_{\tSigma_0}}(\chi)\cap |\Sigma_0|^0) \times [0,1]}$ is a deformation retraction from $\sigma \cap Y_{\psi_{\tSigma_0}}(\chi)\cap |\Sigma_0|^0$ to $Y'(\chi)\cap |\Sigma_0|^0$.
\item $H_\sigma|_{(\sigma \cap \{ x \in |\Sigma_0|^0:\,\,\chi(x)\leq 0\}) \times [0,1]}$ is a deformation retraction from $\sigma \cap \{ x \in |\Sigma_0|^0:\,\,\chi(x)\leq 0\}$ to $Y'(\chi)\cap |\Sigma_0|^0$.
\end{itemize}
To define $H_\sigma$, let $v_1,...,v_r,w_1,...,w_s$ denote the shortest non-zero vectors in $X_*(S)$ on each of the one dimensional faces of $\sigma$ (so that $r+s=\dim \sigma$), with the notation chosen such that $\chi(v_i)\leq 0$ for all $i$ and $\chi(w_j)>0$ for all $j$. Note that $1-\chi(v_i) >0$ for all $i$ and $1-\chi(w_j) \leq 0$ for all $j$. We set
\[ H_\sigma(\sum_i a_i v_i + \sum_j b_i w_j,t)=\sum_i a_i v_i +(1-t) \sum_j b_i w_j .\]
Because 
\[ \begin{array}{rl} Y_{\psi_{\tSigma_0}}(\chi)\cap \sigma \cap |\Sigma_0|^0= & \{ \sum_i a_i v_i+\sum_j b_jw_j: \,\, a_i,b_j \in \R_{\geq 0} \, {\rm and}\\ & \sum_ia_i(1-\chi(v_i))+\sum_jb_j(1-\chi(w_j)) >  0\} \cap |\Sigma_0|^0 \end{array} \]
and
\[ \begin{array}{l} \{ x \in |\Sigma_0|^0:\,\,\chi(x)\leq 0\} \cap \sigma = \\ \{ \sum_i a_i v_i+\sum_j b_jw_j:\,\, a_i,b_j \in \R_{\geq 0} \,\, {\rm and} \sum_ia_i\chi(v_i)+\sum_jb_j\chi(w_j) \leq 0\}\cap |\Sigma_0|^0 \end{array} \]
are convex sets, and because
\[ Y'(\chi) \cap \sigma = \{ \sum_i a_i v_i:\,\, a_i \in \R_{\geq 0}\}, \]
it is easy to check that it has the desired properties and the proof of the lemma is complete.
\pfend

\begin{lem}\label{refine} Let $Y$ be a connected scheme, let $\alpha:S \ra S'$ be an isogeny of split tori over $Y$, and let $\Sigma_0'$ (resp. $\Sigma_0$) be a locally finite partial fan for $S'$ (resp. $S$). Suppose that $Y$ is separated and locally noetherian and that $\Sigma_0'$ is full. Also suppose that $\Sigma_0$ is a quasi-refinement of $\Sigma_0'$, and let $\pi^\wedge$ denote the map $T_{\Sigma_0}^\wedge \ra (\alpha_*T)_{\Sigma_0'}^\wedge$.

Then for $i>0$ we have
\[ R^i\pi^\wedge_{*} \cO_{T_{\Sigma_0}^\wedge} = (0), \]
while
\[ \cO_{(\alpha_*T)_{\Sigma_0'}^\wedge} \liso (\pi^\wedge_{*} \cO_{T_{\Sigma_0}^\wedge})^{\ker \alpha}. \]
If moreover $\Sigma_0$ and $\Sigma_0'$ are smooth then, for $i>0$ we have 
\[ R^i\pi^\wedge_* \cI_{\partial,\Sigma_0}^\wedge = (0). \]
while
\[ \cI_{\partial,\Sigma_0'}^\wedge \liso (\pi^\wedge_* \cI_{\partial, \Sigma_0}^\wedge)^{\ker \alpha}. \]
\end{lem}

\pfbegin 
We may reduce to the case that $Y=\Spec A$ is affine. The map $\pi^\wedge$ is proper and hence
\[ R^i\pi^\wedge_{*} \cO_{T_{\Sigma_0}^\wedge}  \]
and
\[ R^i\pi^\wedge_* \cI_{\partial,\Sigma_0}^\wedge \]
and
\[ \coker( \cO_{(\alpha_*T)_{\Sigma_0'}^\wedge} \lra (\pi^\wedge_{*} \cO_{T_{\Sigma_0}^\wedge})^{\ker \alpha}) \]
and
\[ \coker( \cI_{\partial,\Sigma_0'}^\wedge \lra (\pi^\wedge_* \cI_{\partial,\Sigma_0}^\wedge)^{\ker \alpha}) \]
are coherent sheaves. Thus they have closed support. Their support is also $S$-invariant. Thus it suffices to show that for each maximal element $\sigma' \in \Sigma_0'$ the space $\partial_{\sigma'}(\alpha_*T)_{\tSigma_0'}$ does not lie in the support of these sheaves.  Let $\Sigma_0(\sigma')$ denote the subset of elements $\sigma \in \Sigma_0$ which lie in $\sigma'$, but in no face of $\sigma'$. Then $\Sigma_0(\sigma')$ is a partial fan and $T^\wedge_{\Sigma_0(\sigma')}$ equals the formal completion of $T^\wedge_{\Sigma_0}$ along $\partial_{\sigma'}T_{\tSigma_0'}$. Thus the formal completion of the above four sheaves along $\partial_{\sigma'}T_{\tSigma_0'}$ equal the corresponding sheaf for the pair $\Sigma_0(\sigma')$ and $\{ \sigma'\}$, so that we are reduced to the case that $\Sigma_0'=\{ \sigma'\}$ is a singleton.

In the case that $\Sigma_0'=\{\sigma'\}$ then $(\alpha_*T)_{\Sigma_0'}$ and $Y$ have the same underlying topological space. Let $\pi_1^\wedge$ denote the map
of ringed spaces $T^\wedge_{\Sigma_0} \ra Y$. Then it suffices to show that for $i>0$ we have
\[ R^i\pi^\wedge_{1,*} \cO_{T_{\Sigma_0}^\wedge} = (0) \]
and
\[ R^i\pi^\wedge_{1,*} \cI_{\partial,\Sigma_0}^\wedge = (0); \]
and that we have
\[ \cO_{(\alpha_*T)_{\Sigma_0'}^\wedge} \liso (\pi^\wedge_{1,*} \cO_{T_{\Sigma_0}^\wedge})^{\ker \alpha} \]
and
\[ \cI_{\partial,\Sigma_0'}^\wedge \liso (\pi^\wedge_{1,*} \cI_{\partial,\Sigma_0}^\wedge)^{\ker \alpha}. \]
This follows from lemma \ref{ctorem2}. (Note that
\[ \prod_{\chi \in |\Sigma_0|^\vee \cap X^*(S)} \cL(\chi) = \bigoplus_{\xi \in (\ker \alpha)^\vee} \prod_{\substack{\chi \in |\Sigma_0|^\vee \cap X^*(S) \\ \chi|_{\ker \alpha = \xi}}} \cL(\chi), \]
where $\ker \alpha$ acts on the $\xi$ term via $\xi$; and that
\[ \{ \chi \in |\Sigma_0|^\vee \cap X^*(S):\,\, \chi|_{\ker \alpha}=1\}= |\Sigma_0|^\vee \cap X^*(S') = |\{ \sigma'\}|^{\vee} \cap X^*(S'). \]
These assertions remain true with $|\Sigma_0|^{\vee,0}$ replacing $|\Sigma_0|^\vee$ and $|\{ \sigma'\}|^{\vee,0}$ replacing $|\{\sigma'\}|^{\vee}$.)
\pfend

\begin{lem}\label{ctorem3} Let $Y$ be a connected scheme, let $S/Y$ be a split torus, let $T/Y$ be an $S$-torsor, let $\Sigma_0$ be a partial fan for $S$, and let $\pi_{\Sigma_0}^\wedge$ denote the map of ringed spaces $T^\wedge_{\Sigma_0} \ra Y$. Suppose that $Y$ is  separated and locally noetherian, that $\Sigma_0$ is full, locally finite and open, and that $|\Sigma_0|^0$ is convex. 
\begin{enumerate}
\item Then
\[ R^i\pi^\wedge_{\Sigma_0,*} \cO_{T_{\Sigma_0}^\wedge} = \left\{ \begin{array}{ll} \prod_{\chi \in |\Sigma_0|^\vee}  \cL(\chi)& {\rm if}\,\, i=0
\\ (0) & {\rm otherwise.}\end{array} \right. \]

\item If in addition $\Sigma_0$ is smooth then
\[ R^i\pi_{\Sigma_0,*}^\wedge \cI_{\partial,\Sigma_0}^\wedge = \left\{ \begin{array}{ll} \prod_{\chi \in |\Sigma_0|^{\vee,0}} \cL(\chi) & {\rm if}\,\, i=0
\\ (0) & {\rm otherwise.}\end{array} \right. \]
\end{enumerate} 
\end{lem}

\pfbegin Let $\sigma_1,\sigma_2,...$ be an enumeration of the $1$ cones in $\tSigma_0$. Let $\Delta^{(i)} \subset |\Sigma|$ denote the convex hull of $\sigma_1,...,\sigma_i$. It is a rational polyhedral cone contained in $|\Sigma_0|$, and there exists $i_0$ such that for $i \geq i_0$ the cone $\Delta^{(i)}$ will have the same dimension as $X_*(S)_\R$. Let $\partial \Delta^{(i)}$ denote the union of the proper faces of $\Delta^{(i)}$; and let $\Delta^{(i),c}$ denote the closure of $|\Sigma_0|-\Delta^{(i)}$ in $|\Sigma_0|$. 

Define recursively fans $\Sigma^{(i)}$ and boundary data $\Sigma^{(i)}_0$ as follows. We set $\Sigma^{(i_0-1)}=\tSigma_0$ and $\Sigma^{(i_0-1)}_0=\Sigma_0$. For $i\geq i_0$ set
\[ \Sigma^{(i)} = \{ \sigma \cap \Delta^{(i)},\,\, \sigma \cap \partial \Delta^{(i)}, \,\, \sigma  \cap \Delta^{(i),c}:\,\,\, \sigma \in \Sigma^{(i-1)} \}. \]
Then $\Sigma^{(i)}$ refines $\Sigma^{(i-1)}$ and we choose $\Sigma^{(i)}_0$ to be the unique subset of $\Sigma^{(i)}$ such that $(\Sigma^{(i)},\Sigma_0^{(i)})$ refines 
$(\Sigma^{(i-1)},\Sigma_0^{(i-1)})$. Then $\widetilde{\Sigma_0^{(i)}}=\Sigma^{(i)}$. We also check by induction on $i$ that
\begin{itemize}
\item $\Sigma^{(i)} \cup \Sigma^{(i-1)}-(\Sigma^{(i)}\cap \Sigma^{(i-1)})$ is finite;
\item and $\Sigma^{(i)}_0$ is locally finite.
\end{itemize}
(The point being that the local finiteness of $\Sigma^{(i-1)}_0$ implies that only finitely many elements of $\Sigma^{(i-1)}_0$, and hence of $\Sigma^{(i-1)}$, meet both $\Delta^{(i)}-\partial \Delta^{(i)}$ and $X_*(S)_\R-\Delta^{(i)}$.)

Now define $\Sigma^{(\infty)}$ (resp. $\Sigma^{(\infty)}_0$) to be the set of cones that occur in $\Sigma^{(i)}$ (resp. $\Sigma_0^{(i)}$) for infinitely many $i$. Alternatively
\[ \Sigma^{(\infty)}=\bigcup_i \{ \sigma \in \Sigma^{(i)}: \,\, \sigma \subset \Delta^{(i)} \}. \]
Then $\Sigma^{(\infty)}$ is a fan, $\Sigma_0^{(\infty)}$ provides locally finite boundary data for $\Sigma^{(\infty)}$, we have $\widetilde{\Sigma_0^{(\infty)}}=\Sigma^{(\infty)}$, and $(\Sigma^{(\infty)},\Sigma^{(\infty)}_0)$ refines $(\Sigma,\Sigma_0)$. Moreover $\Sigma^{(\infty)}_0$ is open. We also define
$\Sigma^{(\infty)}_i$ to be the set of $\sigma \in \Sigma^{(i)}_0$ such that $\sigma \subset \Delta^{(i)}$ but $\sigma \not\subset \partial \Delta^{(i)}$. Note that:
\begin{itemize}
\item $\Sigma^{(\infty)}_i$ is finite and open;
\item $\Sigma^{(\infty)}_i \supset \Sigma^{(\infty)}_{i-1}$; 
\item $|\Sigma^{(\infty)}_i|^0=\Delta^{(i)} - \partial \Delta^{(i)}$ is convex;
\item and $\Sigma^{(\infty)}_0 = \bigcup_{i>0} \Sigma^{(\infty)}_i$.
\end{itemize}
(For the last of these properties use the fact that 
$\Sigma^{(\infty)}_0$ is open.) 

By lemma \ref{refine} it suffices to prove this lemma after replacing the pair $\Sigma_0$ by $\Sigma^{(\infty)}_0$. This lemma then follows from lemmas \ref{ctorem1} and \ref{ctorem2}.
\pfend

\newpage \subsection{The case of a disconnected base.}

{\em Throughout this section we will continue to assume that $Y$ is a separated scheme, flat and locally of finite type over $\Spec R_0$.}

Let $S$ be a split torus over $Y$ and let $T/Y$ be an $S$-torsor. By a rational polyhedral cone $\sigma$ in $X_*(S)_\R$ we shall mean a locally constant sheaf of subsets $\sigma \subset X_*(S)_\R$, such that 
\begin{itemize}
\item for each connected open $U \subset Y$ the set $\sigma(U) \subset X_*(S)_\R(U)$ is either empty or a rational polyhedral cone,
\item and the locus where $\sigma \neq \emptyset$ is non-empty and connected.
\end{itemize}
We call this locus the {\em support} of $\sigma$.
We call $\sigma'$ a {\em face} of $\sigma$ if for each open connected $U$ either $\sigma(U)=\sigma'(U)=\emptyset$ or the cone $\sigma'(U)$ is a face of $\sigma(U)$. We call $\sigma$ {\em smooth} if each $\sigma(U)$ is smooth.
By a {\em fan} $\Sigma$ in $X_*(S)_\R$ we mean a set of  rational polyhedral cones in $X_*(S)_\R$, such that
\begin{itemize}
\item if $\sigma \in \Sigma$ then so is any face $\sigma'$ of $\sigma$;
\item if $\sigma,\sigma' \in \Sigma$ then $\sigma \cap \sigma'$ is either empty or a face of $\sigma$ and $\sigma'$.
\end{itemize}
Thus to give a fan in $X_*(S)_\R$ is the same as giving a fan in $X_*(S)_\R(Z)$ for each connected component $Z$ of $Y$. If $U$ is a connected open in $Y$ then we set
\[ \Sigma(U)=\{ \sigma(U): \,\, \sigma \in \Sigma\} -\{ \emptyset\} . \]
It is a fan for $X_*(S)_\R(U)$.

We call $\Sigma$ {\em smooth} (resp. {\em full}, resp. {\em finite}) if each $\Sigma(U)$ is. 
We define a locally constant sheaf $|\Sigma|$ 
of subsets of $X_*(S)_\R$ by setting
\[ |\Sigma|(U)= \bigcup_{\sigma \in \Sigma} \sigma(U) \]
(resp.
\[ |\Sigma|^*(U)= \bigcup_{\sigma \in \Sigma} (\sigma(U) -\{0\})) \]
for $U$ any connected open subset of $Y$. 
We will call $|\Sigma|$ (resp. $|\Sigma|^*$) convex if $|\Sigma|(U)$ (resp. $|\Sigma|^*(U)$) is for each connected open $U\subset Y$.
We also define locally constant sheaves of subsets $|\Sigma|^{\vee}$ and $|\Sigma|^{\vee,0}$ of $X^*(S)_\R$ by setting
\[ |\Sigma|^{\vee}(U)  = \{ \chi \in X^*(S)_\R(U):\,\, \chi(|\Sigma|(U)) \subset \R_{\geq 0}\} = \bigcap_{\sigma \in \Sigma} \sigma(U)^{\vee} \]
and
\[ |\Sigma|^{\vee,0}(U)  = \{ \chi \in X^*(S)_\R(U):\,\, \chi(|\Sigma|^*(U)) \subset \R_{>0}\} = \bigcap_{\sigma \in \Sigma} \sigma(U)^{\vee,0}. \]
We call $\Sigma'$ a {\em refinement} of $\Sigma$ if each $\Sigma'(U)$ is a refinement of $\Sigma(U)$ for each open, connected $U$.
Any fan $\Sigma$ has a smooth refinement $\Sigma'$ such that every smooth cone $\sigma \in \Sigma$ also lies in $\Sigma'$. 

To a fan $\Sigma$ one can attach a scheme $T_\Sigma$ flat and separated over $Y$ and locally (on $T_\Sigma$) of finite type over $Y$, together with an action of $S$ and an $S$-equivariant embedding $T \into T_\Sigma$.  It has an open cover $\{ T_\sigma\}_{\sigma \in \Sigma}$, with each $T_\sigma$ relatively affine over $Y$. 
Over a connected open $U \subset Y$ this restricts to the corresponding picture for ${\Sigma(U)}$.
We write $\cO_{T_\Sigma}$ for the structure sheaf of $T_\Sigma$.
If $\Sigma$ is smooth then $T_\Sigma/Y$ is smooth. 
If $\Sigma$ is finite and $|\Sigma|=X_*(S)_\R$, then $T_\Sigma$ is proper over $Y$. If $\Sigma'$ refines $\Sigma$ then there is an $S$-equivariant proper map 
\[ T_{\Sigma'} \ra T_\Sigma \]
which restricts to the identity on $T$.

By {\em boundary data} we shall mean a proper subset $\Sigma_0 \subset \Sigma$ such that 
$\Sigma - \Sigma_0$ is a fan. If $U\subset Y$ is a connected open we set
\[ \Sigma_0(U)=\{ \sigma(U):\,\, \sigma \in \Sigma_0\} -\{ \emptyset\}. \]
If $\Sigma_0$ is boundary data, then we can associate to it a closed sub-sheme 
$\partial_{\Sigma_0} T_\Sigma \subset  T_\Sigma$, which over a connected open $U\subset Y$ restricts to $\partial_{\Sigma_0(U)} (T|_U)_{\Sigma(U)} \subset  (T|_U)_{\Sigma(U)}$.

In the case that $\Sigma_0$ is the set of elements of $\Sigma$ of dimension bigger than $0$ we shall simply write $\partial T_\Sigma$ for $\partial_{\Sigma_0}T_\Sigma$. Thus $T=T_\Sigma-\partial T_\Sigma$. We will write $\cI_{\partial T_\Sigma}$ for the ideal of definition of $\partial T_\Sigma$ in $\cO_{T_\Sigma}$. We will also write $\cM_\Sigma \ra \cO_{T_\Sigma}$ for the associated log structure and $\Omega^1_{T_\Sigma/\Spec R_0}(\log \infty)$ for the log differentials $\Omega^1_{T_\Sigma/\Spec R_0}(\log \cM_\Sigma)$. 

If $\Sigma$ is smooth then $\partial T_\Sigma$ is a simple normal crossing divisor on $T_\Sigma$ relative to $Y$.

If $\sigma \in \Sigma$ has positive dimension and if $\Sigma_0$ denotes the set of elements of $\Sigma$ which have $\sigma$ for a face, then we will write $\partial_\sigma T_\Sigma$ for $\partial_{\Sigma_0} T_\Sigma$.
It is connected and flat over $Y$ of codimension in $T_\Sigma$ equal to the dimension of $\sigma$. If $\Sigma$ is smooth then each $\partial_\sigma T_\Sigma$ is smooth over $Y$. 
The schemes $\partial_{\sigma_1} T_\Sigma,...,\partial_{\sigma_s} T_\Sigma$ intersect if and only if $\sigma_1, ..., \sigma_s$ are all contained in some $\sigma \in \Sigma$. In this case the intersection equals $\partial_\sigma T_\Sigma$ for the smallest such $\sigma$. We set
\[ \partial_i T_\Sigma = \coprod_{\dim \sigma =i} \partial_\sigma T_\Sigma. \]

If the connected components of $Y$ are irreducible then each $\partial_\sigma T_\Sigma$ is irreducible. Moreover the irreducible components of $\partial T_\Sigma$ are the $\partial_\sigma T_\Sigma$ as $\sigma$ runs over one dimensional elements of $\Sigma$. 
If $\Sigma$ is smooth then we see that $\cS(\partial T_\Sigma)$ is the delta complex with cells in bijection with the elements of $\Sigma$ with dimension bigger than $0$ and with the same `face relations'. In particular it is in fact a simplicial complex and
\[ |\cS(\partial T_\Sigma)|= \coprod_{Z \in \pi_0(Y)} |\Sigma|^*(Z)/\R^{\times}_{>0}. \]

We will call $\Sigma_0$ {\em open} (resp. {\em finite}, resp. {\em locally finite}) if $\Sigma_0(U)$ is for each open connected $U \subset Y$. If $\Sigma_0$ is finite and open, then $\partial_{\Sigma_0} T_\Sigma$ is proper over $Y$. 

By a {\em partial fan} in $X_*(S)$ we mean a collection $\Sigma_0$ of rational polyhedral cones in $X_*(S)$ such that 
\begin{itemize}
\item $\Sigma_0$ does not contain $(0) \subset X_*(S)(U)_\R$ for any open connected $U$;
\item if $\sigma_1, \sigma_2 \in \Sigma_0$ then $\sigma_1 \cap \sigma_2$ is either empty or a face of $\sigma_1$ and of $\sigma_2$;
\item if $\sigma_1, \sigma_2 \in \Sigma_0$ and if $\sigma \supset \sigma_2$ is a face of $\sigma_1$, then $\sigma \in \Sigma_0$.
\end{itemize}
In this case we will let $\tSigma_0$ denote the set of faces of elements of $\Sigma_0$. It is a fan, and $\Sigma_0$ is boundary data for $\tSigma_0$. By {\em boundary data} $\Sigma_1$ for $\Sigma_0$ we shall mean a subset $\Sigma_1 \subset \Sigma_0$ such that if $\sigma \in \Sigma_0$ contains $\sigma_1 \in \Sigma_1$, then $\sigma \in \Sigma_1$. In this case $\Sigma_1$ is again a partial fan and boundary data for $\tSigma_0$. We say that a partial fan $\Sigma_0$ for $X_*(S)$ {\em refines} a partial fan $\Sigma_0'$ for $X_*(S)$ if every element of $\Sigma_0$ lies in an element of $\Sigma_0'$ and if every element of $\Sigma_0'$ is a finite union of elements of $\Sigma_0$. 

If $\Sigma_0$ is a partial fan we define locally constant sheaves of subsets $|\Sigma_0|$, $|\Sigma_0|^*$, $|\Sigma_0|^\vee$ and $|\Sigma_0|^{\vee,0}$ of $X_*(S)_\R$ or $X^*(S)_\R$ to be $|\tSigma_0|$, $|\tSigma_0|^*$, $|\tSigma_0|^\vee$ and $|\tSigma_0|^{\vee,0}$ respectively. We also define a sheaf of subsets $
|\Sigma_0|^0$ by
\[ |\Sigma_0|^0(U)=|\tSigma_0|(U)-|\tSigma_0-\Sigma_0|(U) \]
for any connected open set $U \subset Y$. We will call $|\Sigma_0|$ (resp. $|\Sigma_0|^0$) {\em convex} if $|\Sigma_0|(U)$ (resp. $|\Sigma_0|^0(U)$) is convex for all open connected subsets $U \subset Y$. 

We will call $\Sigma_0$ {\em smooth} (resp. {\em full}, resp. {\em open}, resp. {\em finite}, resp. {\em locally finite}) if $\Sigma_0(U)$ is for each $U \subset Y$ open and connected. 

If $\Sigma_0$ is a partial fan we will write 
\[ \partial_{\Sigma_0} T \]
for $\partial_{\Sigma_0} T_{\tSigma_0}$;
\[ T^\wedge_{\Sigma_0} \]
for the completion of $T_{\tSigma_0}$ along $\partial_{\Sigma_0} T_{\tSigma_0}$; and 
\[ \cM^\wedge_{\Sigma_0} \lra \cO_{T_{\Sigma_0}^\wedge} \]
for the log structure induced by $\cM_{\tSigma_0}$.  We make the following definitions. 
\begin{itemize}
\item $\cI_{T_{\Sigma_0}^\wedge}$ will denote the completion of $\cI_{\partial_{\Sigma_0} T_{\tSigma_0}}$, the sheaf of ideals defining $\partial_{\Sigma_0} T_{\tSigma_0}$. It is an ideal of definition for $T_{\Sigma_0}^\wedge$.
\item $\cI_{\partial, \Sigma_0}^\wedge$ will denote the completion of $\cI_{\partial T_{\tSigma_0}}$, the sheaf of ideals defining $\partial_{\Sigma_0} T_{\tSigma_0}$. Thus $\cI_{\partial, \Sigma_0}^\wedge \subset \cI_{T_{\Sigma_0}^\wedge}$.
\item  $\Omega^1_{T_{\Sigma_0}^\wedge/\Spf R_0}(\log \infty)$ will denote $\Omega^1_{T_{\Sigma_0}^\wedge/\Spf R_0}(\log \cM_\Sigma^\wedge)$. 
\end{itemize}

We will write 
\[ \prod_{\chi \in |\Sigma_0|^\vee} \cL_T(\chi) \]
(resp.
\[ \prod_{\chi \in |\Sigma_0|^{\vee,0}} \cL_T(\chi)) \]
for the sheaf (of abelian groups) on $Y$ such that for any connected open subset $U \subset Y$ we have
\[ \left. \left(\prod_{\chi \in |\Sigma_0|^\vee} \cL_T(\chi)\right)\right|_U = \prod_{\chi \in |\Sigma_0|^\vee(U) \cap X^*(S)(U)} \cL_T(\chi) \]
(resp.
\[ \left. \left(\prod_{\chi \in |\Sigma_0|^{\vee,0}} \cL_T(\chi)\right)\right|_U = \prod_{\chi \in |\Sigma_0|^{\vee,0}(U) \cap X^*(S)(U)} \cL_T(\chi)). \]

Suppose that $\alpha:S \ra S'$ is a surjective map of split tori over $Y$. Then $X^*(\alpha): X^*(S') \into X^*(S)$ and $X_*(\alpha):X_*(S)_\R \onto X_*(S')_\R$. We call fans $\Sigma$ for $X_*(S)$ and $\Sigma'$ for $X_*(S')$ {\em compatible} if for all $\sigma \in \Sigma$ the image $X_*(\alpha) \sigma$ is contained in some element of $\Sigma'$. In this case the map $\alpha:T \ra \alpha_* T$ extends to an $S$-equivariant map
\[ \alpha:T_\Sigma \lra (\alpha_*T)_{\Sigma'}. \] 
We will write
\[ \Omega_{T_\Sigma/(\alpha_* T)_{\Sigma'}}^1(\log \infty) = \Omega_{T_\Sigma/(\alpha_* T)_{\Sigma'}}^1(\log \cM_\Sigma/\cM_{\Sigma'}). \]
The following lemma is an immediate consequence of lemma \ref{logs}.

\begin{lem}\label{dlogs} If $\alpha$ is surjective and $\# \coker X_*(\alpha)$ is invertible on $Y$ then $\alpha:(T_\Sigma,\cM_\Sigma) \ra ((\alpha_*T)_{\Sigma'},\cM_{\Sigma'})$ is log smooth, and there is a natural isomorphism
\[ (X^*(S)/X^*(\alpha)X^*(S')) \otimes_\Z \cO_{T_\Sigma} \liso \Omega^1_{T_\Sigma/(\alpha_*T)_{\Sigma'}}(\log \infty). \]
\end{lem}

If $\alpha$ is an isogeny, if $\Sigma$ and $\Sigma'$ are compatible, and if every element of $\Sigma'$ is a finite union of elements of $\Sigma$, then we call $\Sigma$ a {\em quasi-refinement} of $\Sigma'$. In that case the map $T_{\Sigma} \ra T_{\Sigma'}$ is proper. 

Suppose that $\alpha:S \onto S'$ is a surjective map of tori, and that $\Sigma_0$ (resp. $\Sigma_0'$) is a partial fan for $S$ (resp. $S'$). We call $\Sigma_0$ and $\Sigma_0'$ {\em compatible} if for every $\sigma \in \Sigma_0$ the image $X_*(\alpha) \sigma$ is contained in some element of $\Sigma_0'$ but in no element of $\tSigma_0'-\Sigma_0'$. In this case there is a natural morphism
\[ \alpha: (T^\wedge_{\Sigma_0},\cM_{\Sigma_0}^\wedge) \lra ((\alpha_*T)_{\Sigma_0'}^\wedge,\cM_{\Sigma_0'}^\wedge). \]
We will write
\[ \Omega_{T_{\Sigma_0}^\wedge/(\alpha_* T)_{\Sigma_0'}^\wedge}^1(\log \infty) = \Omega_{T_{\Sigma_0}^\wedge/(\alpha_* T)^\wedge_{\Sigma_0'}}^1(\log \cM_{\Sigma_0}^\wedge/\cM^\wedge_{\Sigma_0'}). \]
The following lemma follows immediately from lemma \ref{clogs}.

\begin{lem}\label{dclogs} If $\alpha$ is surjective and $\# \coker X_*(\alpha)$ is invertible on $Y$ then there is a natural isomorphism
\[ (X^*(S)/X^*(\alpha)X^*(S')) \otimes_\Z \cO_{T_{\Sigma_0}^\wedge} \liso \Omega_{T_{\Sigma_0}^\wedge/(\alpha_* T)_{\Sigma_0'}^\wedge}^1(\log \infty) . \]
\end{lem}

We will call $\Sigma_0$ and $\Sigma_0'$ {\em strictly compatible} if they are compatible and if an element of $\tSigma_0$ lies in $\Sigma_0$ if and only if it maps to no element of $\tSigma_0'-\Sigma_0'$. We will say that
\begin{itemize}
\item $\Sigma_0$ is {\em open over} $\Sigma_0'$ if $|\Sigma_0|^0(U)$ is open in $X_*(\alpha)^{-1} |\Sigma_0'|^0(U)$ for all connected opens $U \subset Y$;
\item and that $\Sigma_0$ is {\em finite over} $\Sigma_0'$ if only finitely many elements of $\Sigma_0$ map into any element of $\Sigma_0'$.
\end{itemize}
If $\alpha$ is an isogeny, if $\Sigma_0$ and $\Sigma_0'$ are strictly compatible, and if every element of $\Sigma_0'$ is a finite union of elements of $\Sigma_0$, then we call $\Sigma_0$ a {\em quasi-refinement} of $\Sigma_0'$. In this case $\Sigma_0$ is open and finite over $\Sigma_0'$. 
The next lemma follows immediately from lemma \ref{proper}. 

\begin{lem}\label{dproper} Suppose that $\Sigma_0'$ and $\Sigma_0$ are strictly compatible.
\begin{enumerate}
\item $T_{\Sigma_0}^\wedge$ is the formal completion of $T_{\tSigma_0}$ along $\partial_{\Sigma_0'}(\alpha_* T)$, and $T_{\Sigma_0}^\wedge$ is locally (on the source) topologically of finite type over $(\alpha_*T)^\wedge_{\Sigma_0'}$.
\item If $\Sigma_0$ is locally finite and if it is open and finite over $\Sigma_0'$ then $T_{\Sigma_0}^\wedge$ is proper over $(\alpha_* T)_{\Sigma_0'}^\wedge$.
\end{enumerate} \end{lem}

\begin{cor} If $\alpha$ is an isogeny, if $\Sigma_0$ is locally finite, and if $\Sigma_0$ is a quasi-refinement of $\Sigma_0'$ then $T_{\Sigma_0}^\wedge$ is proper over $(\alpha_* T)_{\Sigma_0'}^\wedge$.\end{cor}

If $\Sigma_0$ and $\Sigma_0'$ are compatible partial fans and if $\Sigma_1' \subset \Sigma_0'$ is boundary data then $\Sigma_0(\Sigma_1')$ will denote the set of elements $\sigma \in \Sigma_0$ such that $X_*(\alpha)\sigma$ is contained in no element of $\Sigma_0'-\Sigma_1'$. It is boundary data for $\Sigma_0$. Moreover the formal completion of $T_{\Sigma_0}^\wedge$ along the reduced sub-scheme of $(\alpha_*T)_{\Sigma_1'}^\wedge$ is canonically identified with $T^\wedge_{\Sigma_0(\Sigma_1')}$. If $\Sigma_1'=\{ \sigma'\}$ is a singleton we will write $\Sigma_0(\sigma')$ for $\Sigma_0(\{ \sigma'\})$.

The next two lemmas follow immediately from lemmas \ref{refine} and \ref{ctorem3} respectively.

\begin{lem}\label{dcrefine} Let $Y$ be a scheme, let $\alpha:S \ra S'$ be an isogeny of split tori over $Y$, let $\Sigma_0'$ (resp. $\Sigma_0$) be a locally finite partial fan for $S'$ (resp. $S$). Suppose that $Y$ is separated and locally noetherian, that $\Sigma_0'$ is full and that $\Sigma_0$ is locally finite. Also suppose that $\Sigma_0$ is a quasi-refinement of $\Sigma_0'$. Let $\pi^\wedge$ denote the map $T_{\sigma_0}^\wedge \ra (\alpha_*T)^\wedge_{\sigma_0'}$. 

Then for $i>0$ we have
\[ R^i\pi^\wedge_{*} \cO_{T_{\Sigma_0}^\wedge} = (0), \]
while
\[ \cO_{(\alpha_*T)_{\Sigma_0'}^\wedge} \liso (\pi^\wedge_{*} \cO_{T_{\Sigma_0}^\wedge})^{\ker \alpha}. \]
If moreover $\Sigma_0$ and $\Sigma_0'$ are smooth then, for $i>0$ we have 
\[ R^i\pi^\wedge_* \cI_{\partial,\Sigma_0}^\wedge = (0). \]
while
\[ \cI_{\partial,\Sigma_0'}^\wedge \liso (\pi^\wedge_* \cI_{\partial, \Sigma_0}^\wedge)^{\ker \alpha}. \]
\end{lem}

\begin{lem}\label{dctorem3} Let $Y$ be a scheme, let $S/Y$ be a split torus, let $T/Y$ is an $S$-torsor, let $\Sigma_0$ be a partial fan for $S$, and let $\pi_{\Sigma_0}^\wedge$ denote the map $T^\wedge_{\Sigma_0} \ra Y$. Suppose that $Y$ is  separated and locally noetherian, that $\Sigma_0$ is full, locally finite and open, and that $|\Sigma_0|^0$ is convex. 
\begin{enumerate}
\item Then
\[ R^i\pi^\wedge_{\Sigma_0,*} \cO_{T_{\Sigma_0}^\wedge} = \left\{ \begin{array}{ll} \prod_{\chi \in |\Sigma_0|^\vee}  \cL(\chi)& {\rm if}\,\, i=0
\\ (0) & {\rm otherwise.}\end{array} \right. \]

\item If in addition $\Sigma_0$ is smooth then
\[ R^i\pi_{\Sigma_0,*}^\wedge \cI_{\partial , \Sigma_0}^\wedge = \left\{ \begin{array}{ll} \prod_{\chi \in |\Sigma_0|^{\vee,0}} \cL(\chi) & {\rm if}\,\, i=0
\\ (0) & {\rm otherwise.}\end{array} \right. \]
\end{enumerate} 
\end{lem}

\newpage

\section{Shimura Varieties.}

In this section we will describe the Shimura varieties associated to $G_n$ and the mixed Shimura varieties associated to $G_n^{(m)}$ and 
$\tG_n^{(m)}$. We assume that all schemes discussed in this section are locally noetherian.

\subsection{Some Shimura varieties.}\label{introshim}

By a {\em $G_n$-abelian scheme} over a 
scheme $Y/\Q$ we shall mean an abelian
scheme $A/Y$ of relative dimension $n[F:\Q]$ together with an embedding
\[ i:F \into \End(A/Y)_\Q \]
such that $\Lie A$ is a locally free right $\cO_Y 
\otimes_\Q F$-module of rank $n$.  By a {\em morphism (resp. quasi-isogeny) of $G_n$-abelian schemes} we mean a morphism (resp. quasi-isogeny) of abelian schemes which commutes with the $F$-action. If $(A,i)$ is a $G_n$-abelian scheme then we give $A^\vee$ the structure $(A^\vee,i^\vee)$ of a $G_n$-abelian scheme by setting $i^\vee(a)=i(a^c)^\vee$. By a {\em quasi-polarization of a $G_n$-abelian scheme}
$(A,i)/Y$ we shall mean a quasi-isogeny $\lambda: A \ra A^\vee$ of $G_n$-abelian schemes, some $\Q^\times$-multiple of which is a polarization. If $Y=\Spec k$ with $k$ a field, we will let $\langle \,\,\,,\,\,\,\rangle_\lambda$ denote the Weil pairing induced on the adelic Tate module $VA$ (see section 23 of \cite{mumav}).

\begin{lem}\label{freet} If $k$ is a field of characteristic $0$ and if $(A,i,\lambda)/k$ is a $G_n$-abelian scheme, 
then $V_p(A \times \bark)$ is a free $F_p$-module of rank $2n$.
\end{lem}

\pfbegin
We may suppose that $k$ is a finitely generated field extension of $\Q$, which we may embed into $\C$.
Then 
\[ (V_p(A \times \bark) \otimes_{\Q_p,\imath} \C) \cong (\Lie A_{\bary} \otimes_{k} \C) \oplus  (\Lie A_{\bary} \otimes_{k,c} \C), \]
so that $V_p(A \times \bark) \otimes_{\Q_p,\imath} \C$ is a free $F \otimes_\Q \C$-module. As $F \otimes_\Q \C = F_p \otimes_{\Q_p, \imath} \C$ we deduce that $V_p (A \times \bark)$ is a free $F_p$-module, as desired.
\pfend

By an {\em ordinary $G_n$-abelian scheme} over a $\Z_{(p)}$-
scheme $Y$ we shall mean an abelian
scheme $A/Y$ of relative dimension $n[F:\Q]$, such that for each geometric point $\bary$ of $Y$ we have $\# A[p](k(\bary)) \geq p^{n[F:\Q]}$, together with an embedding
\[ i:\cO_{F,(p)} \into \End(A/Y)_{\Z_{(p)}} \]
such that $\Lie A$ is a locally free right $\cO_Y 
\otimes_{\Z_{(p)}} \cO_{F,(p)}$-module of rank $n$.  By a {\em morphism of ordinary $G_n$-abelian schemes} we mean a morphism of abelian schemes which commutes with the $\cO_{F,(p)}$-action. If $(A,i)$ is an ordinary  $G_n$-abelian scheme then we give $A^\vee$ the structure, $(A^\vee,i^\vee)$, of a $G_n$-abelian scheme by setting $i^\vee(a)=i(a^c)^\vee$. By a {\em prime-to-$p$ quasi-polarization of an ordinary $G_n$-abelian scheme}
$(A,i)/Y$ we shall mean a prime-to-$p$ quasi-isogeny $\lambda: A \ra A^\vee$ of ordinary $G_n$-abelian schemes, some $\Z^\times_{(p)}$-multiple of which is a prime-to-$p$ polarization.

If $U$ is an open
compact subgroup of $G_n(\A^{\infty})$ then by a {\em $U$-level structure} on a
quasi-polarized $G_n$-abelian variety $(A,i,\lambda)$ over a connected scheme $Y/\Spec
\Q$ with a geometric point $\bary$, we mean a
$\pi_1(Y,\bary)$-invariant $U$-orbit $[\eta]$ of pairs
$(\eta_0,\eta_1)$ of $\A^\infty$-linear isomorphisms
\[ \eta_0: \A^{\infty}_\bary \liso \A^{\infty}(1)_\bary=V\G_{m,\bary} \]
and
\[ \eta_1:V_n \otimes_{\Q} \A^{\infty} \liso VA_\bary \]
such that 
\[ \eta_1(ax) = i(a) \eta_1(x) \]
for all $a \in F$ and $x \in V_n \otimes_{\Q} \A^{\infty}$,
and such that
\[ \langle \eta_1 x , \eta_1 y \rangle_\lambda = \eta_0 \langle x,y \rangle_n \]
for all $x,y \in V_n \otimes_{\Q} \A^{\infty}$. This definition is
independent of the choice of geometric point $\bary$ of $Y$. By a $U$-level structure on a quasi-polarized
$G_n$-abelian scheme $(A,i,\lambda)$ over a general (locally noetherian) scheme $Y/\Spec
\Q$, we mean the collection of a $U$-level structure over each connected component of $Y$.  If $[(\eta_0,\eta_1)]$ is a level structure we define $||\eta_0|| \in \Q^\times_{>0}$ by
\[ ||\eta_0|| \eta_0 \Zhat = \Zhat(1). \]

Now suppose that $U^p$ is an open compact subgroup of $G_n(\A^{p,\infty})$ and that $N_1 \leq N_2$ are non-negative integers. 
By a $U^p(N_1,N_2)$-level structure on an ordinary, prime-to-$p$
quasi-polarized, $G_n$-abelian scheme $(A,i,\lambda)$ over a connected scheme $Y/\Spec
\Z_{(p)}$ with a geometric point $\bary$, we mean a
$\pi_1(Y,\bary)$-invariant $U^p$-orbit $[\eta]$ of four-tuples
$(\eta_0^p,\eta_1^p,C,\eta_p)$ consisting of
\begin{itemize}
\item an $\A^{p,\infty}$-linear isomorphism $\eta_0^p: \A^{p,\infty}_\bary \liso \A^{p,\infty}(1)_\bary=V^p\G_{m,\bary}$;
\item an $\A_F^{p,\infty}$-linear isomorphism
\[ \eta_1^p: V_n \otimes_{\Q} \A^{p,\infty} \liso V^pA_\bary \]
such that
\[ \langle \eta_1^p x , \eta_1^p y \rangle_\lambda = \eta_0 \langle x,y \rangle_n \]
for all $x,y \in V_n \otimes_{\Q} \A^{p,\infty}$;
\item a locally free sub-$\cO_{F,(p)}$-module scheme $C \subset A[p^{N_2}]$,
such that for every geometric point $\ty$ of $Y$ there is an $\cO_{F,(p)}$-invariant sub-Barsotti-Tate group $\tC_{\ty} \subset A_{\ty}[p^\infty]$ with the following properties
\begin{itemize}
\item $C_{\ty}=\tC_{\ty}[p^{N_2}]$,
\item for all $N$ the sub-group scheme $\tC_{\ty}[p^N]$ is isotropic in $A[p^N]_{\ty}$ for the $\lambda$-Weil pairing,
\item $A_{\ty}[p^\infty]/\tC_{\ty}$ is ind-etale,
\item the Tate module $T(A_{\ty}[p^\infty]/\tC_{\ty})$ is free over $\cO_{F,p}$ of rank $n$;
\end{itemize}
\item and an isomorphism 
\[  \eta_p: p^{-N_1}\Lambda_n/(p^{-N_1}\Lambda_{n,(n)}+ \Lambda_n) \liso A[p^{N_1}]/(A[p^{N_1}] \cap C) \]
such that 
\[ \eta_p(ax) = i(a) \eta_p(x) \]
for all $a \in \cO_{F,(p)}$ and $x \in p^{-N_1}\Lambda_n/(p^{-N_1}\Lambda_{n,(n)}+ \Lambda_n)$.
\end{itemize}
This definition is
independent of the choice of geometric point $\bary$ of $Y$. By a $U^p(N_1,N_2)$-level structure on an ordinary, prime-to-$p$ quasi-polarized,
$G_n$-abelian scheme $(A,i,\lambda)$ over a general (locally noetherian) scheme $Y/\Spec
\Z_{(p)}$, we mean the collection of a $U^p(N_1,N_2)$-level structure over each connected component of $Y$. 
If $[(\eta_0^p,\eta_1^p,C,\eta_p)]$ is a level structure we define $||\eta_0^p|| \in \Z^\times_{(p),>0}$ by
\[ ||\eta_0^p|| \eta_0^p \Zhat^p = \Zhat^p(1). \]

If $(A,i,\lambda,[\eta])/Y$ is a quasi-polarized, $G_n$-abelian scheme with $U$-level structure and if $g \in G_n(\A^{\infty})$ 
with $U' \supset g^{-1} U g$, then we can
define a quasi-polarized, $G_n$-abelian scheme with
$U'$-level structure $(A,i,\lambda,[\eta])g /Y$ by
\[ (A,i,\lambda,[(\eta_0,\eta_1)])g=(A,i,\lambda,[(\nu(g) \eta_0, \eta_1 \circ
g]). \]

If $(A,i,\lambda,[\eta])/S$ is an ordinary, prime-to-$p$ quasi-polarized, $G_n$-abelian scheme with $U^p(N_1,N_2)$-level structure and if $g\in G_n(\A^\infty)^{\ord,\times}$ with $(U')^p(N_1',N_2') \supset g^{-1} U^p(N_1,N_2) g$ (so that in particular $N_i \geq N_i'$ for $i=1,2$),
then we can
define an ordinary, prime-to-$p$ quasi-polarized, $G_n$-abelian scheme with
$(U')^p(N_1',N_2')$-level structure $(A,i,\lambda,[\eta])g /Y$ by
\[ (A,i,\lambda,[(\eta_0^p,\eta_1^p,C,\eta_p)])g=(A,i,\lambda,[(\nu(g^p) \eta_0^p, \eta_1^p \circ
g^p,C[p^{N_2'}],\eta_p \circ g_p)]). \]
If $(U')^p(N_1',N_2') \supset \varsigma_p^{-1} U^p(N_1,N_2) \varsigma_p$ (so that in particular $N_1 \geq N_1'$ and $N_2> N_2'$), then we can
define an ordinary, prime-to-$p$ quasi-polarized, $G_n$-abelian scheme with
$(U')^p(N_1',N_2')$-level structure $(A,i,\lambda,[\eta])\varsigma_p/Y$ by
\[ \begin{array}{l} (A,i,\lambda,[(\eta_0^p,\eta_1^p,C,\eta_p)])\varsigma_p=\\ (A/C[p],i,F(\lambda) ,[(p\eta_0^p, F(\eta_1^p), C[p^{1+N_2'}]/C[p], F(\eta_p))]); \end{array}\] 
where 
\[ F(\lambda): A/C[p] \stackrel{\lambda}{\lra} A^\vee/\lambda C[p] = A^\vee/C[p]^\perp \liso (A/C[p])^\vee  \]
with the latter isomorphism being induced by the dual of the map $A/C[p] \ra A$ induced by multiplication by $p$ on $A$; 
where $F(\eta_1^p)$ is the composition of $\eta_1^p$ with the natural map $V^pA \iso V^p(A/C[p])$; and where $F(\eta_p)$ is the composition of $\eta_p$ with the natural identification  
\[ A[p^{N_1'}]/(C \cap A[p^{N_1'}]) = (A/C[p])[p^{N_1'}]/(C[p^{1+N_2'}]/C[p] \cap (A/C[p])[p^{N_1'}]). \]
Together these two definitions give an action of $G_n(\A^\infty)^\ord$.

By a quasi-isogeny (resp. prime-to-$p$ quasi-isogeny) between quasi-polarized,
$G_n$-abelian schemes with $U$-level structures (resp. ordinary, prime-to-$p$ quasi-polarized,
$G_n$-abelian schemes with $U^p(N_1,N_2)$-level structures)
\[ (\beta,\delta): (A, i, \lambda, [\eta]) \lra (A',i',\lambda',[\eta']) \]
we mean a quasi-isogeny (resp. prime-to-$p$ quasi-isogeny)  of abelian schemes $\beta \in \Hom(A,A')_\Q$ (resp. $\beta \in \Hom(A,A')_{\Z_{(p)}}$) and $\delta \in \Q^\times$ (resp. $\delta \in \Z_{(p)}^\times$) such that
\begin{itemize}
\item $\beta \circ i(a)=i'(a) \circ \beta$ for all $a \in F$ (resp. $\cO_{F,(p)}$);
\item $\delta \lambda=\beta^\vee \circ \lambda' \circ \beta$;
\item $[(\delta \eta_0,(V \beta) \circ \eta_1)]=[\eta']$ (resp. $[(\delta \eta_0^p,(V^p \beta) \circ \eta_1^p, \beta C, \beta \circ \eta_p)]=[\eta']$).
\end{itemize}
The action of $G_n^{(m)}(\A^\infty)$ (resp. $G_n^{(m)}(\A^{\infty})^\ord$) takes one quasi-isogeny (resp. prime-to-$p$ quasi-isogeny) class to another. 

\begin{lem}\label{linalg1} Suppose that $T$ is an $\cO_{F,p}$-module, which is free over $\cO_{F,p}$ of rank $2n$, with a perfect alternating pairing
\[ \langle \,\,\,,\,\,\, \rangle: T \times T \lra \Z_p \]
such that
\[ \langle ax,y \rangle = \langle x,a^cy \rangle \]
for all $x,y \in T$ and $a \in \cO_{F,p}$. Also suppose that $\tT \subset T$ is a sub-$\cO_{F,p}$-module which is isotropic for $\langle \,\,\,,\,\,\,\rangle$ and such that $T/\tT$ is free of rank $n$ over $\cO_{F,p}$. Finally suppose that
\[  \eta_p: p^{-N_1}\Lambda_n/(p^{-N_1}\Lambda_{n,(n)}+ \Lambda_n)  \liso p^{-N_1}T/(p^{-N_1}\tT+T) \]
is an $\cO_{F,p}$-module isomorphism.

Consider the set $[\eta]$ of isomorphisms
\[ \eta: \Lambda_n \otimes \Z_p \liso T \]
such that
\begin{itemize}
\item $\eta(a x)=a \eta(x)$ for all $a \in \cO_{F,(p)}$;
\item there exists $\delta \in \Z_p^\times$ such that
\[ \langle \eta x, \eta y \rangle = \delta \langle x,y\rangle_n \]
for all $x,y \in \Lambda_n \otimes \Z_p$;
\item $\eta ((p^{-N_2} \Lambda_{n,(n)})\otimes \Z_p + \Lambda_{n}\otimes \Z_p) = p^{-N_2} \tT+T$;
\item the map
\[  p^{-N_1}\Lambda_n/(p^{-N_1}\Lambda_{n,(n)}+ \Lambda_n)  \liso p^{-N_1}T/(p^{-N_1}\tT+T)\]
induced by $\eta$ equals $\eta_p$.
\end{itemize}

Then $[\eta]$ is non-empty and a single $U_p(N_1,N_2)$-orbit.\end{lem}

\pfbegin
Let $e_1,...,e_n$ denote a $\cO_{F,p}$-basis of $T/\tT$.
Note that $\langle\,\,\,,\,\,\,\rangle$ induces a perfect pairing between $\tT$ and $T/\tT$.
We recursively lift the $e_i$ to elements $\te_i \in T$ with $\te_i$ orthogonal to the $\cO_{F,p}$ span of the $\te_j$ for $j\leq i$. Suppose that $\te_1,...,\te_{i-1}$ have already been chosen. Choose some lift $e_i'$ of $e_i$. Then choose $t \in \tT$ such that
\begin{itemize}
\item $\langle t , x \rangle = \langle e_i',x \rangle$ for all $x \in \bigoplus_{j=1}^{i-1} \cO_{F,p} \te_j$,
\item and $\langle t, \alpha e_i' \rangle = \langle e_i'/2,\alpha e_i'\rangle$ for all $\alpha \in \cO_{F,p}^{c=-1}$.
\end{itemize}
(If $p=2$ some explanation is required as to why we can do this. The map
\[ \begin{array}{rcl}
\cO_{F,p} & \lra & \Z_p \\ \alpha & \longmapsto & \langle e_i',\alpha e_i'\rangle \end{array} \]
is of the form
\[ \alpha \longmapsto \tr_{F/\Q} \beta \alpha \]
for some $\beta \in (\cD_{F,p}^{-1})^{c=-1}$. Because $p=2$ is unramified in $F/F^+$, we can write $\beta = \gamma -\gamma^c$ for some $\gamma \in \cD_{F,p}^{-1}$. Thus the second condition can be replaced by the condition
\[ \langle t, \alpha e_i'\rangle = \tr_{F/\Q} \gamma \alpha \]
for all $\alpha \in \cO_{F,p}$. Now it is clear that the required element $t$ exists.)
Then take $\te_i=e_i'-t$. Then $\te_i$ is orthogonal to $\bigoplus_{j=1}^{i-1} \cO_{F,p} \te_j$. Moreover for $\alpha \in \cO_{F,p}$ we have
\[ \begin{array}{rcl} \langle \te_i,\alpha \te_i \rangle &=& \langle e_i',\alpha e_i'\rangle -\langle t, (\alpha-\alpha^c)e_i' \rangle \\ &=& \langle e_i',\alpha e_i'\rangle -\langle e_i'/2, (\alpha-\alpha^c)e_i' \rangle \\ &=& (\langle e_i', \alpha e_i' \rangle + \langle e_i', \alpha^c e_i' \rangle)/2 \\ &=& 0. \end{array} \]

Thus we can write 
\[ T = \tT \oplus \tT' \]
with $\tT'$ an isotropic $\cO_{F,p}$-subspace of $T$, which is free over $\cO_{F,p}$ of rank $n$. We see that
\[ \tT' \cong \Hom_{\Z_p}(\tT, \Z_p). \]
The lemma now follows without difficulty.
\pfend

\begin{cor}\label{linalg2} If $Y$ is a $\Q$-scheme with geometric point $\bary$, if $(A,i,\lambda)/Y$ is an ordinary $G_n$-abelian scheme, and if $[(\eta_0^p,\eta_1^p,C, \eta_p)]$ is a $U^p(N_1,N_2)$-level structure on $(A,i,\lambda)$, then there is a unique $U_p(N_1,N_2)$-orbit of pairs of isomorphisms
\[ \eta_{0,p}: \Z_{p,\bary} \liso \Z_p(1)_\bary \]
and
\[ \eta_{1,p}:\Lambda_n \otimes \Z_p \liso T_pA_\bary \]
such that
\begin{itemize}
\item $\eta_{1,p}(a x)=a \eta_{1,p}(x)$ for all $a \in \cO_{F,(p)}$,
\item $\langle \eta_{1,p}x,\eta_{1,p}y\rangle_\lambda = \eta_{0,p}\langle x,y \rangle_n$ for all $x,y \in  \Lambda_n \otimes \Z_p$,
\item $\eta_{1,p} p^{-N_2} \Lambda_{n,(n)}/\Lambda_{n,(n)} = C$,
\item $\eta_{p,1}$ induces $\eta_p$.
\end{itemize} \end{cor}

\pfbegin
This follows on combining the lemmas \ref{freet} and \ref{linalg1}.
\pfend

\begin{cor} Suppose that $Y$ is a scheme over $\Spec \Q$. Then there is a natural bijection between prime-to-$p$ isogeny classes of ordinary, prime-to-$p$ quasi-polarized $G_n$-abelian schemes with $U^p(N_1,N_2)$-level structure and
isogeny classes of quasi-polarized $G_n$-abelian schemes with $U^p(N_1,N_2)$-level structure.
This bijection is $G_n(\A^{\infty})^\ord$-equivariant.
\end{cor}

\pfbegin
We may assume that $Y$ is connected with geometric point $\bary$. We will show both sets are in natural bijection with the set of prime-to-p isogeny classes of 4-tuples $(A,i,\lambda,[\eta])$, where $(A,i)$ is a $G_n$-abelian variety, $\lambda$ is a prime-to-$p$ quasi-polarization of $(A,i)$, and $[\eta]$ is a $\pi_1(Y,\bary)$-invariant $U^p(N_1,N_2)$-orbit of pairs $(\eta_0,\eta_{1})$, where
\begin{itemize}
\item $\eta_0: \A^{p,\infty}\times \Z_p \iso \A^{p,\infty}(1) \times \Z_p(1)$,
\item and $\eta_1: \Lambda_n \otimes (\A^{\infty,p} \times \Z_p) \iso V^pA_\bary \times T_pA_\bary$ satisfies
\[ \eta_0 \langle x,y \rangle_n = \langle \eta_1 x, \eta_1 y\rangle_\lambda .\]
\end{itemize}
There is a natural map from this set to the set of isogeny classes of quasi-polarized $G_n$-abelian schemes with $U^p(N_1,N_2)$-level structure, which is easily checked to be a bijection.
The bijection between this set and the set of prime-to-$p$ isogeny classes of ordinary, prime-to-$p$ quasi-polarized $G_n$-abelian schemes with $U^p(N_1,N_2)$-level structure, follows by the usual arguments (see for instance section III.1 of \cite{ht}) from corollary \ref{linalg2}. 
\pfend

If $U$ is a neat open compact subgroup of $G_n(\A^\infty)$ then the functor that
sends a (locally noetherian) scheme $S/\Q$ to the set of quasi-isogeny classes of polarized $G_n$-abelian schemes  with
$U$-level structures is represented by a quasi-projective
scheme $X_{n,U}$ which is smooth of relative dimension
$n^2[F^+:\Q]$ over $\Q$. 
Let
\[ [(A^\univ,i^\univ,\lambda^\univ,[\eta^\univ])]/ X_{n,U} \]
denote the universal equivalence class of polarized $G_n$-abelian varieties with
$U$-level structure. If $U' \supset g^{-1}Ug$ then there is a map $g:X_{n,U} \ra X_{n,U'}$  arising from $(A^\univ,i^\univ,\lambda^\univ,[\eta^\univ])g/X_{n,U}$ and the universal
property of $X_{n,U'}$. This makes $\{ X_{n,U}\}$ an inverse system of schemes with right $G_n(\A^{\infty})$-action.
The maps $g$ are finite etale. If $U_1 \subset U_2$ is a normal subgroup then $X_{n,U_1}/X_{n,U_2}$ is Galois with group $U_2/U_1$.

There are identifications of topological spaces:
\[ X_{n,U}(\C) \cong G_n(\Q)^+\backslash (G_n(\A^\infty)/U \times \gH_n^+) \cong  G_n(\Q)\backslash (G_n(\A^\infty)/U \times \gH_n^\pm)\] 
compatible with the right action of $G_n(\A^{\infty})$. (Note that
$\ker^1(\Q,G_n)=(0)$.) More precisely we associate to $(g,I) \in
G_n(\A^{\infty})/U \times \gH_n^+$ the torus $\Lambda  \otimes_\Z \R/\Lambda$ with complex
structure coming from $I$; with polarization
corresponding to the Riemann form given by $\langle \,\,\, ,\,\,\,\rangle$; 
and with level structure coming from
\[ \eta_1:\Lambda_n \otimes \A^{\infty} \stackrel{g}{\lra} \Lambda_n \otimes \A^{\infty} = V(\Lambda_n  \otimes_\Z \R/\Lambda_n) \]
and
\[ \begin{array}{rcl} \eta_0: \A^{\infty} &\liso &\A^{\infty}(1) \\ x &
\longmapsto & -x \zeta, \end{array} \]
where $\zeta = \lim_{\la N} e^{2 \pi i/N} \in \widehat{\Z}(1)$. We deduce that
\[ \begin{array}{rcl} \pi_0(X_{n,U} \times \Spec \barQQ) & \cong & G_n(\Q) \backslash G_n(\A)/(U G_n(\R)^+) \\
& \cong & G_n(\Q) \backslash (G_n(\A^\infty)/U \times \pi_0(G_n(\R))) \\ & \cong & 
C_n(\Q) \backslash C_n(\A) / U C_n(\R)^0 . \end{array} \]

If $U^p$ is neat then the functor that
sends a scheme $Y/\Z_{(p)}$ to the set of prime-to-$p$ quasi-isogeny classes of ordinary, prime-to-$p$ quasi-polarized, $G_n$-abelian schemes with
$U^p(N_1,N_2)$-level structure is represented by a scheme $\cX_{n,U^p(N_1,N_2)}^\ord$ quasi-projective over $\Z_{(p)}$. Let
\[ [(\cA^\univ,i^\univ,\lambda^\univ,[\eta^\univ])]/ \cX_{n,U^p(N_1,N_2)}^\ord \]
denote the universal equivalence class. If $g \in G_n(\A^{\infty})^\ord$ and $(U^p)'(N_1',N_2') \supset g^{-1}U^p(N_1,N_2)g$, then there is a quasi-finite map 
\[ g:\cX_{n,U^p(N_1,N_2)}^\ord \lra \cX_{(U^p)'(N_1',N_2')}^\ord \]
arising from $(\cA^\univ,i^\univ,\lambda^\univ,[\eta^\univ])g/\cX_{n,U^p(N_1,N_2)}$ and the universal
property of $\cX_{n,(U^p)'(N_1',N_2')}^\ord$. 
If $g \in G_n(\A^\infty)^{\ord,\times}$ then the map $g$ is etale, and, if further $N_2=N_2'$, then it is finite etale. 
If $U^p(N_1,N_2)$ is a normal subgroup of $(U^p)'(N_1',N_2)$ then $\cX_{n,U^p(N_1,N_2)}^\ord/\cX_{n,(U^p)'(N_1',N_2)}^\ord$ is Galois with group $(U^p)'(N_1')/U^p(N_1)$.
There are $G_n(\A^{\infty})^\ord$-equivariant identifications
\[ \cX_{n,U^p(N_1,N_2)}^\ord \times \Spec \Q \cong X_{n,U^p(N_1,N_2)}. \]
The scheme $\cX_{n,U^p(N_1,N_2)}^\ord$ is smooth over $\Z_{(p)}$ of relative dimension $n^2[F^+:\Q]$. (By the Serre-Tate theorem (see \cite{katzst}) the formal completion of $\cX_{n,U^p(N_1,N_2)}^\ord$ at a point $x$ in the special fibre is isomorphic to 
\[ \Hom_{\Z_p}(S(T_p\cA^\univ_x),\widehat{\G}_m). \]
This is formally smooth as long as $S(T_p\cA^\univ_x) \cong S(\cO_{F,p}^n)$ is torsion free. This module is torsion free because in the case $p=2$ we are assuming that $p=2$ is unramified in $F$.)
Suppose that $g \in G_n(\A^{\infty})^\ord$ and $(U^p)'(N_1',N_2') \supset g^{-1}U^p(N_1,N_2)g$, then the quasi-finite map 
\[ g:\cX_{n,U^p(N_1,N_2)}^\ord \lra \cX_{(U^p)'(N_1',N_2')}^\ord \]
is in fact flat, because the it is a quasi-finite map between locally noetherian regular schemes which are equi-dimensional of the same dimension. (See pages 507 and 508 of \cite{km}.)

On $\F_p$-fibres the map  
\[ \varsigma_p\cX^\ord_{n,U^p(N_1,N_2+1)} \times \Spec \F_p \lra \cX^\ord_{n,U^p(N_1,N_2)} \times \Spec \F_p \]
is the absolute Frobenius map composed with the forgetful map $1: \cX^\ord_{n,U^p(N_1,N_2)} \ra \cX^\ord_{n,U^p(N_1,N_2-1)}$ (for any $N_2 \geq N_1 \geq 0$). Thus if $N_2>0$, then the quasi-finite, flat map
\[ \varsigma_p:\cX^\ord_{n,U^p(N_1,N_2+1)}  \lra \cX^\ord_{n,U^p(N_1,N_2)}  \]
has all its fibres of degree $p^{n^2[F^+:\Q]}$ and hence is finite flat of this degree.  (A flat, quasi-finite morphism $f:X \ra Y$ between noetherian schemes with constant fibre degree is proper and hence, by theorem 8.11.1 of \cite{ega4}, finite. We give the argument for properness. By the valuative criterion we may reduce to the case $Y=\Spec B$ for a DVR $B$ with fraction field $L$.  By theorem 8.12.6 of \cite{ega4} $X$ is a dense open subscheme of $\Spec A$, for $A$ a finite $B$ algebra. Let  $I$ denote the ideal of $A$ consisting of all $\gm_B$-torsion elements. If $f \in A$ and $\Spec A_f \subset X$, then by flatness the map $A \ra A_f$ factors through $A/I$. Thus $X \subset \Spec A/I$ and in fact $I=(0)$, so that $A/B$ is finite flat. Because an open subscheme is determined by its points, we conclude that we must have $X=\Spec A'$ for some $A \subset A' \subset A \otimes_B L$. By the constancy of the fibre degree we conclude that $A'$ is  
 finite over $B$.) 
We deduce that for any $g \in G_n(\A^\infty)^{\ord}$, if $N_2'>0$ and $p^{N_2-N_2'}\nu(g_p) \in \Z_p^\times$, then the map 
\[ g:\cX_{n,U^p(N_1,N_2)}^\ord \lra \cX_{(U^p)'(N_1',N_2')}^\ord  \]
 is finite.

\begin{lem} Write $\cX_{n,U^p(N_1,N_2)}^{\ord,\wedge}$ for the completion of $\cX_{U^p(N_1,N_2)}^{\ord}$ along its $\F_p$-fibre. If $N_2'>N_2 \geq N_1$ then the map
\[ 1: \cX_{n,U^p(N_1,N_2')}^{\ord,\wedge} \lra \cX_{n,U^p(N_1,N_2)}^{\ord,\wedge} \]
is an isomorphism.
\end{lem}

\pfbegin The map has an inverse which sends $[(\cA^\univ,i^\univ,\lambda^\univ, [\eta^\univ])]$ over $\cX_{n,U^p(N_1,N_2)}^{\ord,\wedge}$ to 
\[ [(\cA^\univ,i^\univ,\lambda^\univ, [(\eta_0^{\univ, p}, \eta_1^{\univ,p},\cA^\univ[p^{N_2'}]^0 , \eta_p^\univ )])] \]
over $\cX_{n,U^p(N_1,N_2')}^{\ord,\wedge}$.
\pfend

Thus we will denote $\cX_{n,U^p(N_1,N_2)}^{\ord,\wedge}$ simply
\[ \gX^\ord_{n,U^p(N_1)}. \]
Then $\{ \gX^\ord_{n,U^p(N)} \}$ is a system of $p$-adic formal schemes with right $G_n(\A^{\infty})^\ord $-action. 
We will write $\barX^\ord_{n,U^p(N)}$ for the reduced sub-scheme of  $\gX_{n,U^p(N)}^\ord$. 

\newpage \subsection{Some Kuga-Sato varieties.}\label{s32}
 
 Recall that a semi-abelian scheme is a smooth separated commutative group scheme such that each geometric fibre is the extension of an abelian variety by a torus. To an semi-abelian scheme $G/Y$ one can associate an etale constructible sheaf of abelian groups $X^*(G)$, the `character group of the toric part of $G$'. See theorem I.2.10 of \cite{cf}. If $X^*(G)$ is constant then $G$ is an extension
 \[ (0) \lra S_G \lra G \lra A_G \lra (0) \]
 of a uniquely determined abelian scheme $A_G$ by a uniquely determined split torus $S_G$. By an {\em isogeny} (resp. {\em prime-to-$p$ isogeny}) of semi-abelian schemes we mean a morphism which is quasi-finite and surjective (resp. quasi-finite and surjective and whose geometric fibres have orders relatively prime to $p$). If $Y$ is locally noetherian, then by a {\em quasi-isogeny} (resp. {\em prime-to-$p$ quasi-isogeny}) $\alpha:G \ra G'$ we mean an element of $\Hom(G,G')_{\Q}$ (resp. $\Hom(G,G')_{\Z_{(p)}}$) with an inverse in $\Hom(G',G)_{\Q}$ (resp. $\Hom(G',G)_{\Z_{(p)}}$).
 
 Suppose that $Y/\Spec \Q$ is a locally noetherian scheme. By a {\em $G_n^{(m)}$-semi-abelian scheme} $G$ over  $Y$ we mean a triple $(G,i,j)$ where 
 \begin{itemize}
 \item $G/Y$ is a semi-abelian scheme, 
 \item $i: F \into \End(G)_\Q$ such that $\Lie A_G$ is a free $\cO_Y \otimes_\Q F$ module of rank $n[F:\Q]$,
 \item and $j:F^m \iso X^*(G)_\Q$ is an $F$-linear isomorphism.
 \end{itemize}
  Then $A_G$ is a $G_n$-abelian scheme. By a {\em quasi-isogeny of $G_n^{(m)}$-semi-abelian schemes} we mean a quasi-isogeny of semi-abelian schemes
 \[ \beta:G \ra G' \] such that
 \[ i'(a) \circ \beta = \beta \circ i(a) \]
 for all $a \in F$, and
 \[ j = X^*(\beta) \circ j'. \]
 Note that, if $\bary$ is a geometric point of $Y$, then $j$ induces a map
 \[ j^*: VS_{G,\bary} \liso \Hom_\Q(F^m,V\G_{m,\bary}). \]
 By a {\em quasi-polarization} of $(G,i,j)$ we mean a quasi-polarization of $A_G$. 
 
 If $Y$ is connected and $\bary$ is a geometric point of $Y$ and if $U \subset G_n^{(m)}(\A^\infty)$ is a neat open compact subgroup then by a {\em $U$ level structure} on a quasi-polarized $G_n^{(m)}$-semi-abelian scheme $(G,i,j,\lambda)$ we mean a $\pi_1(Y,\bary)$-invariant $U$-orbit of pairs $(\eta_0,\eta_1)$ where
 \[ \eta_0: \A^\infty \liso V\G_{m,\bary} \]
 is an $\A^\infty$-linear map,  and where
 \[ \eta_1: \Lambda_n^{(m)} \otimes_\Z \A^\infty \lra VG_\bary \]
 is an $\A_F^\infty$-linear map such that
 \[ \eta_1|_{\Hom_\Z(\cO_F^m,\A^\infty)} = (j^*)^{-1}\circ \Hom(1_{F^m},\eta_0)  \]
 and
 \[ [(\eta_0, \eta_1 \bmod VS_{G,\bary})] \]
 is a $U$-level structure on $A_G$. This is canonically independent of $\bary$. We define a $U$ level structure on a $G_n^{(m)}$-semi-abelian scheme over a general locally noetherian scheme $Y$ to be such a level structure over each connected component of $Y$. By a quasi-isogeny between two quasi-polarized, $G_n^{(m)}$-semi-abelian schemes with $U$-level structure
 \[ (\beta,\delta):  (G,i,j,\lambda,[(\eta_0,\eta_1)]) \lra (G',i',j',\lambda',[(\eta_0',\eta_1')]) \]
 we mean a quasi-isogeny
 \[ \beta: (G,i,j) \lra (G',i',j') \]
 and an element $\delta \in \Q^\times$ such that
 \[ \delta \lambda = \beta^\vee \circ \lambda' \circ \beta \]
 and
 \[ [(\eta_0',\eta_1')] = [( \delta \eta_0, V(\beta) \circ \eta_1)]. \]
 
  If $(G,i,j,\lambda,[(\eta_0,\eta_1)])$ is a quasi-polarized, $G_n^{(m)}$-semi-abelian scheme with $U$-level structure, if $g \in G_n^{(m)}(\A^\infty)$ and if $U'\supset g^{-1}Ug$ then we define a quasi-polarized, $G_n^{(m)}$-semi-abelian scheme with $U'$-level structure
 \[ (G,i,j,\lambda,[(\eta_0,\eta_1)])g=(G,i,j,\lambda,[(\nu(g)\eta_0,\eta_1\circ g)]). \]
 The quasi-isogeny class of $(G,i,j,\lambda,[(\eta_0,\eta_1)])g$ only depends on the quasi-isogeny class of 
 $(G,i,j,\lambda,[(\eta_0,\eta_1)])$. If $(G,i,j,\lambda,[(\eta_0,\eta_1)])$ is a quasi-polarized, $G_n^{(m)}$-semi-abelian scheme with $U$-level structure, if $\gamma \in GL_m(F)$ and $U' \supset \gamma U$ then we define a  
quasi-polarized, $G_n^{(m)}$-semi-abelian scheme with $U'$-level structure
 \[ \gamma (G,i,j,\lambda,[(\eta_0,\eta_1)])=(G,i,j \circ \gamma^{-1},\lambda,[(\eta_0,\eta_1 \circ \gamma)]). \]
 The quasi-isogeny class of $\gamma(G,i,j,\lambda,[(\eta_0,\eta_1)])$ only depends on the quasi-isogeny class of 
 $(G,i,j,\lambda,[(\eta_0,\eta_1)])$. We have $\gamma \circ g = \gamma(g) \circ \gamma$.
 If $(G,i,j,\lambda,[(\eta_0,\eta_1)])$ is a quasi-polarized, $G_n^{(m)}$-semi-abelian scheme with $U$-level structure, if $m'\leq m$ and if $U'\supset i_{m',m}^*U$, then we define a  
quasi-polarized, $G_n^{(m')}$-semi-abelian scheme with $U'$-level structure
 \[ \pi_{m,m'} (G,i,j,\lambda,[(\eta_0,\eta_1)])=(G/S,i,j \circ i_{m',m},\lambda,[(\eta_0,\eta_1')]), \]
 where $S \subset S_G$ is the subtorus with 
 \[ X^*(S) = X^*(S_G) / (X^*(S_G) \cap j \circ i_{m',m} F^{m'}) \]
 and where
 \[ \eta_1' \circ i_{m',m}^* = \eta_1 \bmod VS. \]
 The quasi-isogeny class of $\pi_{m,m'}(G,i,j,\lambda,[(\eta_0,\eta_1)])$ only depends on the quasi-isogeny class of 
 $(G,i,j,\lambda,[(\eta_0,\eta_1)])$. If $\gamma \in Q_{m,(m-m')}(F)$ then $\pi_{m,m'} \circ \gamma = \bargamma \circ \pi_{m,m'}$, where $\bargamma$ denotes the image of $\gamma$ in $GL_{m'}(F)$. If $g \in G_n^{(m)}(\A^\infty)$ then $\pi_{m,m'} \circ g = i_{m',m}^*(g)\circ \pi_{m,m'}$.
 
 If $U$ is a neat open compact subgroup of $G_n^{(m)}(\A^\infty)$ then the functor which sends a locally noetherian scheme $Y/\Q$ to the set of quasi-isogeny classes of quasi-polarized $G_n^{(m)}$-semi-abelian schemes with $U$-level structure is represented by a quasi-projective scheme $A_{n,U}^{(m)}$, which is smooth of dimension $n(n+2m)[F^+:\Q]$. (See proposition 1.3.2.14 of \cite{kw2}.) Let 
 \[ [(G^\univ,i^\univ,j^\univ,\lambda^\univ,[\eta^\univ])]/A_{n,U}^{(m)} \]
 denote the universal quasi-isogeny class of quasi-polarized $G_n^{(m)}$-semi-abelian \linebreak schemes with $U$-level structure. If $g \in G_n^{(m)}(\A^\infty)$ and $U_1,U_2$ are neat open compact subgroups of $G_n^{(m)}(\A^\infty)$ with $U_2 \supset g^{-1}U_1g$ then
 there is a map
 \[ g: A^{(m)}_{n,U_1} \lra A^{(m)}_{n,U_2} \]
 arising from $(G^\univ,i^\univ,j^\univ,\lambda^\univ,[\eta^\univ])g/A^{(m)}_{n,U_1}$ and the universal property of $A_{n,U_2}^{(m)}$. Similarly if $\gamma \in GL_m(F)$ and $U_1,U_2$ are neat open compact subgroups of $G_n^{(m)}(\A^\infty)$ with $U_2 \supset \gamma U_1$ then there is a map
 \[ \gamma: A^{(m)}_{n,U_1} \lra A^{(m)}_{n,U_2} \]
 arising from $\gamma(G^\univ,i^\univ,j^\univ,\lambda^\univ,[\eta^\univ])/A^{(m)}_{n,U_1}$ and the universal property of $A_{n,U_2}^{(m)}$. Moreover if $m' \leq m$, if $U \subset G_n^{(m)}(\A^\infty)$ and if $U'$ denotes the image of $U$ in $G_n^{(m')}(\A^\infty)$, then there is a smooth projective map 
 \[ \pi_{A^{(m)}_n/A^{(m')}_n}: A^{(m)}_{n,U} \lra A_{n,U'}^{(m')} \]
 arising from $\pi_{m,m'} (G^\univ,i^\univ,j^\univ,\lambda^\univ,[\eta^\univ])/A^{(m)}_{n,U}$ and the universal property of $A_{n,U'}^{(m')}$. We see that these actions have the following properties.
 \begin{itemize}
  \item $A^{(0)}_{n,U}=X_{n,U}$. (We will sometimes write $\pi_{A^{(m)}_n/X_n}$ for $\pi_{A^{(m)}_n/A^{(0)}_n}$.) This identification is $G_n(\A^\infty)$-equivariant.
  \item $g_1 \circ g_2=g_2g_1$ (i.e. this is a right action) and $\gamma_1 \circ \gamma_2 = \gamma_1 \gamma_2$ (i.e. this is a left action) and $\gamma \circ g = \gamma(g) \circ \gamma$.
 \item If $\gamma \in Q_{m,(m-m')}(F)$ then $\pi_{A^{(m)}_n/A^{(m')}_n} \circ \gamma = \bargamma \circ \pi_{A^{(m)}_n/A^{(m')}_n}$, where $\bargamma$ denotes the image of $\gamma$ in $GL_{m'}(F)$.
 \item $\pi_{A^{(m)}_n/A^{(m')}_n} \circ g = g'\circ \pi_{A^{(m)}_n/A^{(m')}_n}$, where $g'$ denotes the image of $g$ in $G_n^{(m')}(\A^\infty)$.
 \end{itemize}
 
 Moreover we have the following properties.
 \begin{itemize}
 \item The maps $g$ and $\gamma$ are finite etale. The maps $\pi_{m.m'}$ are smooth and projective.

 \item If $U_1 \subset U_2$ is an open normal subgroup of a neat open compact then $A_{n,U_1}^{(m)}/A_{n,U_2}^{(m)}$ is Galois with group $U_2/U_1$.
\item If $U=U' \ltimes M$ with $U' \subset G_n(\A^\infty)$ and $M \subset \Hom_n^{(m)}(\A^\infty)$ then $A^{(m)}_{n,U}/X_{n,U'}$ is an abelian scheme of relative dimension $mn[F:\Q]$. 
 \item In general $A^{(m)}_{n,U}$ is a principal homogenous space for $A^{(m)}_{n,U' \ltimes (U \cap \Hom^{(m)}_n(\A^\infty))}$ over $X_{n,U'}$, where $U'$ denotes the image of $U$ in $G_n(\A^\infty)$.
 \item There are $G^{(m)}_{n}(\A^\infty)$ and $GL_m(F)$ equivariant homeomorphisms
\[ A^{(m)}_{n,U}(\C) \cong G^{(m)}_{n}(\Q) \backslash G^{(m)}_{n}(\A) /(U \times U_{n,\infty}^0A_n(\R)^0) . \] 
 \end{itemize}
 Moreover in the case $U=U' \ltimes M$, if $G^\univ/A_{n,U}^{(m)}$ and $A^\univ/X_{n,U'}$ are chosen so that $\pi_{A_n^{(m)}/X_{n}}^* A^\univ \cong A_{G^\univ}$, then there is a $\Q$-linear map 
 \[ i^{(m)}_{A^\univ}: F^m \lra \Hom(A_{n,U}^{(m)}, (A^\univ)^\vee)_\Q \]
 with the following properties.
 \begin{itemize}
 \item If $a \in F$ then
 \[ i^{(m)}_{A^\univ}(ax)=i^{\univ,\vee}({}^ca) \circ i^{(m)}_{A^\univ}(x). \]
 
 \item If $(\beta,\delta)$ is a quasi-isogeny
\[ (G^\univ,i^\univ,j^\univ,\lambda^\univ,[\eta^\univ]) \lra (G^{\univ,\prime},i^{\univ,\prime},j^{\univ,\prime},\lambda^{\univ,\prime},[\eta^{\univ,\prime}]), \]
then
\[ \beta^\vee \circ i^{(m)}_{(A^\univ)'}(x) = i^{(m)}_{A^\univ}(x). \]
In particular $i^{(m)}_{A^\univ}$ depends only on $A^\univ$ and not on $G^\univ$.

\item If $g \in G_n^{(m)}(\A^\infty)$ and $\gamma \in GL_m(F)$ then 
\[ i^{(m)}_{A^\univ}(x) \circ g = i^{(m)}_{g^* A^\univ}(x) \]
and
\[ i^{(m)}_{A^\univ}(x) \circ \gamma = i^{(m)}_{\gamma^* A^\univ}(\gamma^{-1} x). \]

\item If $e_1,...,e_m$ denotes the standard basis of $F^m$ then
\[ i_{A^\univ}=||\eta_0^\univ||^{-1}( (\lambda^\univ)^{-1}\circ i^{(m)}(e_1),..., (\lambda^\univ)^{-1} \circ  i^{(m)}(e_m)): A_{n,U}^{(m)} \lra (A^\univ)^m \]
is a quasi-isogeny. If $(\beta,\delta)$ is a quasi-isogeny
\[ (G^\univ,i^\univ,j^\univ,\lambda^\univ,[\eta^\univ]) \lra (G^{\univ,\prime},i^{\univ,\prime},j^{\univ,\prime},\lambda^{\univ,\prime},[\eta^{\univ,\prime}]), \]
then
\[ \beta^{\oplus m} \circ i_{A^\univ}=i_{(A^\univ)'}. \]

\item The map
\[ \begin{array}{ccccc} \!\! \eta_{n,U}^{(m)}:\Hom_F(F^m,V_n) \otimes_\Q \A^\infty &\!\!\! \iso& \!\!\! V(A^\univ)^m &\!\!\!\iso& \!\!\! VA_{n,U}^{(m)} \\
f & \!\!\!\mapsto & \!\!\!(\eta_1^\univ(f(e_1)),...,\eta_1^\univ(f(e_m))) && \\
&& x & \!\!\! \mapsto & \!\!\! V(i_{A^\univ})^{-1} x \end{array} \]
is an isomorphism, which does not depend on the choice of $G^\univ$. It satisfies
\[ \eta_{n,U}^{(m)} M = TA_{n,U}^{(m)}. \]
 \end{itemize}
 (See lemmas 1.3.2.7 and 1.3.2.50, proposition 1.3.2.55, theorem 1.3.3.15, and remark 1.3.3.33 of \cite{kw2}; and section 3.5 of \cite{kwan}.)
 
 Note that 
 \[ i_{A^\univ} \circ g = i_{g^* A^\univ} \]
 and
 \[ i_{A^\univ} \circ \gamma = {}^t\gamma^{-1} \circ i_{\gamma^* A^\univ}. \]
 
Define
\[ i_\lambda^{(m)} : F^m \otimes_{F,c}F^m \lra \Hom(A_{n,U}^{(m)}, A_{n,U}^{(m),\vee})_\Q \]
by 
\[ i_\lambda^{(m)}(x \otimes y) = ||\eta_0^\univ||^{-1} i^{(m)}_{A^\univ}(x)^\vee \circ \lambda^{\univ,-1} \circ i^{(m)}_{A^\univ}(y). \]
This does not depend on the choice of $A^\univ$. We have
\[ i_\lambda^{(m)}(x \otimes y)^\vee = i_\lambda^{(m)}(y \otimes x). \]
Moreover
\[ (i_{A^\univ}^{-1})^\vee \circ i_\lambda^{(m)}(x \otimes y) \circ i_{A^\univ}^{-1} = (\lambda^\univ)^{\oplus m} \circ i^\univ ({}^{c,t}xy).  \]
If $a \in (F^m \otimes_{F,c} F^m)^{\sw =1}$ has image in $S(F^m)$ lying in $S(F^m)^{>0}$ then 
\[ (i_{A^\univ}^{-1})^\vee \circ i_\lambda^{(m)}(a) \circ i_{A^\univ}^{-1} = (\lambda^\univ)^{\oplus m} \circ i^\univ (a') \]
for some matrix $a' \in M_{m \times m}(F)^{t=c}$ all whose eigenvalues are positive real numbers. (See section \ref{herm} for the definition of $\sw$.) Thus $i_\lambda^{(m)}(a)$ is a quasi-polarization. (See the end of section 21 of \cite{mumav}.)

 Now suppose that $Y/\Spec \Z_{(p)}$ is a locally noetherian scheme. By an {\em ordinary $G_n^{(m)}$-semi-abelian scheme} $G$ over  $Y$ we mean a triple $(G,i,j)$ where 
 \begin{itemize}
 \item $G/Y$ is a semi-abelian scheme such that $\# G[p](k(\bary)) \geq p^{n[F:\Q]}$ for each geometric point $\bary$ of $Y$,
 \item $i: \cO_{F,(p)} \into \End(G)_{\Z_{(p)}}$ such that $\Lie A_G$ is a free $\cO_Y \otimes_{\Z_{(p)}} \cO_{F,(p)}$ module of rank $n[F:\Q]$,
 \item and $j:\cO_{F,(p)}^m \iso X^*(G)_{\Z_{(p)}}$ is a $\cO_{F,(p)}$-linear isomorphism.
 \end{itemize}
  Then $A_G$ is an ordinary $G_n$-abelian scheme. By a {\em prime-to-$p$ quasi-isogeny of ordinary $G_n^{(m)}$-semi-abelian schemes} we mean a prime-to-$p$ quasi-isogeny of semi-abelian schemes
 \[ \beta:G \ra G' \] such that
 \[ i'(a) \circ \beta = \beta \circ i(a) \]
 for all $a \in \cO_{F,(p)}$, and
 \[ j = X^*(\beta) \circ j'. \]
 Note that, if $\bary$ is a geometric point of $Y$, then $j$ induces a map
 \[ j^*: V^pS_{G,\bary} \liso \Hom_{\Z_{(p)}}(\cO_{F,(p)}^m,V^p\G_{m,\bary}). \]
 By a {\em prime-to-$p$ quasi-polarization} of $(G,i,j)$ we shall mean a prime-to-$p$ quasi-polarization of $A_G$. 
 
 If $Y$ is connected and $\bary$ is a geometric point of $Y$, if $U^p \subset G_n^{(m)}(\A^{p,\infty})$ is a neat open compact subgroup, and if $N_2 \geq N_1 \geq 0$ then by a {\em $U^p(N_1,N_2)$ level structure} on a prime-to-$p$ quasi-polarized ordinary $G_n^{(m)}$-semi-abelian scheme $(G,i,j,\lambda)$ we mean a $\pi_1(Y,\bary)$-invariant $U^p$-orbit $[\eta]$ of five-tuples
 $(\eta_0^p,\eta_1^p,C,D,\eta_p)$ consisting of
 \begin{itemize}
 
 \item an $\A^{p,\infty}$-linear isomorphism $\eta_0^p: \A^{p,\infty} \liso \A^{p,\infty}(1)_\bary=V^p\G_{m,\bary}$;

\item an $\A_F^{p,\infty}$-linear isomorphism
\[ \eta_1^p: \Lambda_n^{(m)} \otimes_{\Z} \A^{p,\infty} \liso V^pG_\bary \]
such that $\eta_1^p|_{\Hom_\Z(\cO_F^m,\A^{p,\infty})} = (j^*)^{-1} \circ \Hom(1_{\cO_F^m},\eta_0^p)$;

\item a locally free sub-$\cO_{F,(p)}$-module scheme $C \subset G[p^{N_2}]$,
such that for every geometric point $\ty$ of $Y$ there is an $\cO_{F,(p)}$-invariant sub-Barsotti-Tate group $\tC_{\ty} \subset G_{\ty}[p^\infty]$ with the following properties
\begin{itemize}
\item $C_{\ty}=\tC_{\ty}[p^{N_2}]$,
\item $\tC_{\ty} \supset S_{G,\ty}[p^\infty]$,
\item for all $N$ the sub-group scheme $\tC_{\ty}[p^N]/S_{G,\ty}[p^N]$ is isotropic in $A_G[p^N]_{\ty}$ for the $\lambda$-Weil pairing,
\item $G_{\ty}[p^\infty]/\tC_{\ty}$ is ind-etale,
\item the Tate module $T(G_{\ty}[p^\infty]/\tC_{\ty})$ is free over $\cO_{F,p}$ of rank $n$;
\end{itemize}

\item a locally free sub-$\cO_{F,(p)}$-module scheme $D \subset C[p^{N_1}]$ such that $D \iso C[p^{N_1}]/S_G[p^{N_1}]$;

\item and an isomorphism 
\[  \eta_p: p^{-N_1}\Lambda_n/(p^{-N_1}\Lambda_{n,(n)}+ \Lambda_n) \liso G[p^{N_1}]/C[p^{N_1}] \]
such that 
\[ \eta_p(ax) = i(a) \eta_p(x) \]
for all $a \in \cO_{F,(p)}$ and $x \in p^{-N_1}\Lambda_n/(p^{-N_1}\Lambda_{n,(n)}+ \Lambda_n)$;
\end{itemize}
such that
\[ [(\eta^p_0, \eta^p_1 \bmod V^pS_G, C/S_G[p^{N_2}], \eta_p)] \]
is a $U^p(N_1,N_2)$-level structure for $(A_G,i,\lambda)$. 
This definition is
independent of the choice of geometric point $\bary$ of $Y$. By a $U^p(N_1,N_2)$-level structure on an ordinary, prime-to-$p$ quasi-polarized,
$G_n^{(m)}$-semi-abelian scheme $(G,i,j,\lambda)$ over a general (locally noetherian) scheme $Y/\Spec
\Z_{(p)}$, we mean the collection of a $U^p(N_1,N_2)$-level structure over each connected component of $Y$. 

By a prime-to-$p$ quasi-isogeny between two quasi-polarized, ordinary $G_n^{(m)}$-semi-abelian schemes with $U^p(N_1,N_2)$-level structure
 \[ (\beta,\delta):  (G,i,j,\lambda,[(\eta_0,\eta_1)]) \lra (G',i',j',\lambda',[(\eta_0',\eta_1')]) \]
 we mean a prime-to-$p$ quasi-isogeny
 \[ \beta: (G,i,j) \lra (G',i',j') \]
 and an element $\delta \in \Z_{(p)}^\times$ such that
 \[ \delta \lambda = \beta^\vee \circ \lambda' \circ \beta \]
 and
 \[ [((\eta_0^p)',(\eta_1^p)',C',D',\eta_p')] = [( \delta \eta^p_0, V^p(\beta) \circ \eta^p_1, \beta C, \beta D, \beta \circ \eta_p)]. \]
 
  If $(G,i,j,\lambda,[(\eta_0^p,\eta_1^p,C,D,\eta_p)])$ is a prime-to-$p$ quasi-polarized, ordinary $G_n^{(m)}$-semi-abelian scheme with $U^p(N_1,N_2)$-level structure, if $g \in G_n^{(m)}(\A^\infty)^{\ord,\times}$ and if $(U^p)'(N_1',N_2')\supset g^{-1}U^p(N_1,N_2)g$ then we define a prime-to-$p$ quasi-polarized, ordinary $G_n^{(m)}$-semi-abelian scheme with $(U^p)'(N_1',N_2')$-level structure
 \[ (G,i,j,\lambda,[(\eta^p_0,\eta^p_1,C,D,\eta_p)])g=(G,i,j,\lambda,[(\nu(g)\eta_0^p,\eta^p_1\circ g^p,C,D,\eta_p \circ g_p)]). \]
 The prime-to-$p$ quasi-isogeny class of $(G,i,j,\lambda,[(\eta_0^p,\eta_1^p,C,D,\eta_p)])g$ only depends on the prime-to-$p$ quasi-isogeny class of 
 $(G,i,j,\lambda,[(\eta_0^p,\eta_1^p,C,D,\eta_p)])$.  Similarly, if $(G,i,j,\lambda,[(\eta_0^p,\eta_1^p,C,D,\eta_p)])$ is a prime-to-$p$ quasi-polarized, ordinary $G_n^{(m)}$-semi-abelian scheme with $U^p(N_1,N_2)$-level structure and if 
 \[ (U^p)'(N_1',N_2') \supset \varsigma_p^{-1} U^p(N_1,N_2) \varsigma_p,\]
  then we define a  prime-to-$p$
quasi-polarized, ordinary $G_n^{(m)}$-semi-abelian scheme with $(U^p)'(N_1',N_2')$-level structure
 \[ \begin{array}{l} (G,i,j,\lambda,[(\eta_0^p,\eta_1^p,C,D,\eta_p)])\varsigma_p= \\ (G/C[p],i,pj ,F(\lambda),[(p\eta_0^p,F(\eta^p_1),C[p^{1+N_2'}]/C[p],D'[p^{N_1'}],F(\eta_p))]); \end{array} \]
where 
\[ F(\lambda): A_G/C[p] \stackrel{\lambda}{\lra} A_G^\vee/\lambda C[p] = A_G^\vee/C[p]^\perp \liso (A_G/C[p])^\vee  \]
with the latter isomorphism being induced by the dual of the map $A_G/C[p] \ra A_G$ induced by multiplication by $p$ on $A_G$; 
where $F(\eta_1^p)$ is the composition of $\eta_1^p$ with the natural map $V^pG \iso V^p(G/C[p])$; 
where $D'$ denotes the pre-image of $D$ under the multiplication by $p$ map $C \ra C$ modulo $C[p]$;
and where $F(\eta_p)$ is the composition of $\eta_p$ with the natural identification  
\[ G[p^{N_1'}]/(C \cap G[p^{N_1'}]) = (G/C[p])[p^{N_1'}]/(C[p^{1+N_2'}]/C[p] \cap (G/C[p])[p^{N_1'}]). \]
Together these two definitions give an action of $G_n(\A^\infty)^\ord$.
 
 If $(G,i,j,\lambda,[(\eta^p_0,\eta^p_1,C,D,\eta_p)])$ is a prime-to-$p$ quasi-polarized, ordinary $G_n^{(m)}$-semi-abelian scheme with $U^p(N_1,N_2)$-level structure, if $\gamma \in GL_m(\cO_{F,(p)})$ and $(U^p)'(N_1',N_2') \supset \gamma U^p(N_1,N_2)$ then we define a  prime-to-$p$
quasi-polarized, ordinary $G_n^{(m)}$-semi-abelian scheme with $(U^p)'(N_1',N_2')$-level structure
 \[ \gamma (G,i,j,\lambda,[(\eta_0^p,\eta_1^p,C,D,\eta_p)])=(G,i,j \circ \gamma^{-1},\lambda,[(\eta_0^p,\eta^p_1 \circ \gamma,C,D,\eta_p)]). \]
 The prime-to-$p$ quasi-isogeny class of $\gamma(G,i,j,\lambda,[(\eta^p_0,\eta^p_1,C,D,\eta_p)])$ only depends on the quasi-isogeny class of 
 $(G,i,j,\lambda,[(\eta^p_0,\eta^p_1,C,D,\eta_p)])$. We have $\gamma \circ g = \gamma(g) \circ \gamma$.
 If $(G,i,j,\lambda,[(\eta_0^p,\eta_1^p,C,D,\eta_p)])$ is a prime-to-$p$ quasi-polarized, ordinary $G_n^{(m)}$-semi-abelian scheme with $U^p(N_1,N_2)$-level structure, if $m'\leq m$ and if $(U^p)'(N_1',N_2')\supset i_{m',m}^*U^p(N_1,N_2)$, then we define a  
quasi-polarized, ordinary $G_n^{(m')}$-semi-abelian scheme with $(U^p)'(N_1',N_2')$-level structure
 \[ \pi_{m,m'} (G,i,j,\lambda,[(\eta^p_0,\eta^p_1,C,D,\eta_p)])=(G/S,i,j \circ i_{m',m},\lambda,[(\eta^p_0,(\eta^p_1)',C',D',\eta_p)]), \]
 where $S \subset S_G$ is the subtorus with 
 \[ X^*(S) = X^*(S_G) / (X^*(S_G) \cap j \circ i_{m',m} \cO_{F,(p)}^{m'}) \]
 and where
 \[ (\eta^p_1)' \circ i_{m',m}^* = \eta^p_1 \bmod V^pS \]
 and $C'$ (resp. $D'$) denotes the image of $C$ (resp. $D$).
 The prime-to-$p$ quasi-isogeny class of $\pi_{m,m'}(G,i,j,\lambda,[(\eta^p_0,\eta^p_1,C,D,\eta_p)])$ only depends on the quasi-isogeny class of 
 $(G,i,j,\lambda,[(\eta^p_0,\eta^p_1,C,D,\eta_p)])$. If $\gamma \in Q_{m,(m-m')}(\cO_{F,(p)})$ then $\pi_{m,m'} \circ \gamma = \bargamma \circ \pi_{m,m'}$, where $\bargamma$ denotes the image of $\gamma$ in $GL_{m'}(\cO_{F,(p)})$. If $g \in G_n^{(m)}(\A^\infty)$ then $\pi_{m,m'} \circ g = i_{m',m}^*(g)\circ \pi_{m,m'}$.

For each $m \geq 0$ there is a system of $\Z_{(p)}$-schemes $\{ \cA^{(m),\ord}_{n,U^p(N_1,N_2)} \}$ as $U^p$ runs over neat open compact subgroups of $G_n^{(m)}(\A^{p,\infty})$ and $N_1,N_2$ run over integers with $N_2 \geq N_1 \geq 0$, together with the following extra structures:
 \begin{itemize}
 \item If $g \in G_n^{(m)}(\A^{\infty})^\ord$ and $U_2^p(N_{21},N_{22}) \supset g^{-1}U_1^p(N_{11},N_{12})g$ then there is a quasi-finite, flat map
 \[ g: \cA^{(m),\ord}_{n,U_1^p(N_{11},N_{12})} \lra \cA^{(m),\ord}_{n,U_2^p(N_{21},N_{22})}. \]
 \item If $m'\leq m$ and if $(U^p)'$ denotes the image of $U^p$ in $G_n^{(m')}(\A^{p,\infty})$, then there is a smooth projective map with geometrically connected fibres
 \[ \pi_{\cA^{(m),\ord}_n/\cA^{(m'),\ord}_n}: \cA^{(m),\ord}_{n,U^p(N_1,N_2)} \lra \cA^{(m'),\ord}_{n,(U^p)'(N_1,N_2)}. \]
 \item If $\gamma \in GL_m(\cO_{F,(p)})$ and $U^p_2 \supset \gamma U_1^p$ then there is a finite etale map
 \[ \gamma: \cA^{(m),\ord}_{n,U_1^p(N_1,N_2)} \lra \cA^{(m),\ord}_{n,U_2^p(N_1,N_2)}. \]
 \end{itemize}
 Moreover there is a canonical prime-to-$p$ quasi-isogeny class of ordinary $G_n^{(m)}$-semi-abelian schemes with $U^p(N_1,N_2)$ level structure
 \[ (\cG^\univ, i^\univ, j^\univ, \lambda^\univ, [\eta^\univ])/ \cA^{(m),\ord}_{n,U^p(N_1,N_2)} \]
 These enjoy the following properties: 
   \begin{itemize}
  \item $\cA^{(0),\ord}_{n,U^p(N_1,N_2)}=\cX^\ord_{n,U^p(N_1,N_2)}$. (We will sometimes write $\pi_{\cA^{(m),\ord}_n/\cX^\ord_n}$ for $\pi_{\cA^{(m),\ord}_n/A^{(0),\ord}_n}$.) This identification is $G_n(\A^\infty)^\ord$ equivariant.
  
  \item $g_1 \circ g_2=g_2g_1$ (i.e. this is a right action) and $\gamma_1 \circ \gamma_2 = \gamma_1 \gamma_2$ (i.e. this is a left action) and $\gamma \circ g = \gamma(g) \circ \gamma$.
  
 \item If $\gamma \in Q_{m,(m-m')}(\cO_{F,(p)})$ then $\pi_{\cA^{(m),\ord}_n/\cA^{(m'),\ord}_n} \circ \gamma = \bargamma \circ \pi_{\cA^{(m),\ord}_n/\cA^{(m'),\ord}_n}$, where $\bargamma$ denotes the image of $\gamma$ in $GL_{m'}(\cO_{F,(p)})$.
 
 \item $\pi_{\cA^{(m),\ord}_n/\cA_n^{(m'),\ord}} \circ g = g'\circ \pi_{\cA_n^{(m),\ord}/\cA_n^{(m'),\ord}}$, where $g'$ denotes the image of $g$ in $G_n^{(m')}(\A^{\infty})^\ord$.
 
 \item If $g \in G_n^{(m)}(\A^\infty)^\ord$, then the induced map
 \[ g: \cA^{(m),\ord}_{n,U_1^p(N_{11},N_{12})} \lra g^* \cA^{(m),\ord}_{n,U_2^p(N_{21},N_{22})} \]
 over $\cX^\ord_{n,U_1^p(N_{11},N_{12})}$ is finite flat of degree $p^{nm[F:\Q]}$. If $g \in G_n^{(m)}(\A^\infty)^{\ord,\times}$, then this map is also etale.

\item If $U_1^p \subset U_2^p$ is an open normal subgroup of a neat open compact of $G_n^{(m)}(\A^{p,\infty})$ and if $N_{11}\geq N_{21}$, then $\cA_{n,U_1^p(N_{11},N_2)}^{(m),\ord}/\cA_{n,U_2^p(N_{21},N_2)}^{(m),\ord}$ is Galois with group $U_2^p(N_{21})/U^p_1(N_{11})$.

 \item On $\F_p$-fibres the map 
 \[ \varsigma_p:  \cA^{(m),\ord}_{n,U^p(N_1,N_2+1)} \times \Spec \F_p \lra \cA^{(m),\ord}_{n,U^p(N_1,N_2)} \times \Spec \F_p\]
 equals the composition of the absolute Frobenius map with the forgetful map (for any $N_2 \geq N_1 \geq 0$).
 
  \item If $g \in G_n^{(m)}(\A^\infty)^\ord$ and $U_2^p(N_{21},N_{22}) \supset g^{-1}U^p_1(N_{11},N_{12})g$ then the pull back
 $g^* (\cG^\univ_2,i_2^\univ, j_2^\univ,\lambda_2^\univ, [\eta_2^\univ])$ is prime-to-$p$ quasi-isogenous to \linebreak $(\cG^\univ_1,i_1^\univ, j_1^\univ,\lambda_1^\univ, [\eta_1^\univ])g$.
 
\item If $\gamma \in GL_m(\cO_{F,(p)})$ and $U_2^p(N_{21},N_{22}) \supset \gamma U^p_1(N_{11},N_{12})$ then the pull-back $\gamma^* (\cG^\univ_2,i_2^\univ, j_2^\univ,\lambda_2^\univ, [\eta_2^\univ])$ is prime-to-$p$ quasi-isogenous to \linebreak $\gamma (\cG^\univ_1,i_1^\univ, j_1^\univ,\lambda_1^\univ, [\eta_1^\univ])$.

\item If $m' \leq m$ and if $U_2^p(N_{21},N_{22}) \supset i_{m',m}^* U^p_1(N_{11},N_{12})$ then the pull-back $\pi_{\cA_n^{(m)}/\cA_n^{(m')}}^* (\cG^\univ_2,i_2^\univ, j_2^\univ,\lambda_2^\univ, [\eta_2^\univ])$ is prime-to-$p$ quasi-isogenous to $\pi_{m,m'} (\cG^\univ_1,i_1^\univ, j_1^\univ,\lambda_1^\univ, [\eta_1^\univ])$.

\item If $U^p=(U^p)' \ltimes M^p$ with $(U^p)' \subset G_n(\A^{p,\infty})$ and $M^p \subset \Hom_n^{(m)}(\A^{p,\infty})$ then $\cA^{(m),\ord}_{n,U^p(N_1,N_2)}/\cX^\ord_{n,(U^p)'(N_1,N_2)}$ is an abelian scheme of relative dimension $mn[F:\Q]$.

 \item In general $\cA^{(m),\ord}_{n,U^p(N_1,N_2)}$ is a principal homogenous space for the abelian scheme $\cA^{(m),\ord}_{n,((U^p)' \ltimes M^p)(N_1,N_2)}$ over $\cX^\ord_{n,(U^p)'(N_1,N_2)}$, where $(U^p)'$ denotes the image of $U^p$ in $G_n(\A^{p,\infty})$ and $M^p=U^p \cap \Hom^{(m)}_n(\A^{p,\infty})$. 
 
  \item There are natural identifications 
\[ \cA^{(m),\ord}_{n,U^p(N_1,N_2)} \times \Spec \Q \cong A^{(m)}_{n,U^p(N_1,N_2)}. \]
These identifications are compatible with the identifications 
\[ \cX_{n,(U^p)'(N_1,N_2)}^{\ord} \times \Spec \Q \cong X_{n,(U^p)'(N_1,N_2)}\]
 and the maps $\pi_{\cA_n^{(m),\ord}/\cA_n^{(m'),\ord}}$ and $\pi_{A^{(m)}_n/A_n^{(m')}}$.
They are also equivariant for the actions of the semi-group $G_n^{(m)}(\A^{\infty})^\ord$ and the group $GL_m(\cO_{F,(p)})$. 
\end{itemize}
Moreover in the case $U^p=(U^p)' \ltimes M^p$, if $\cG^\univ/\cA^{(m),\ord}_{n,U^p(N_1,N_2)}$ and $\cA^\univ/\cX^\ord_{n,U^p(N_1,N_2)}$ are chosen so that 
$\pi^*_{\cA_n^{(m),\ord}/\cX_n^\ord} \cA^\univ \cong A_{\cG^\univ}$, then there is a $\Z_{(p)}$-linear map
\[ i^{(m)}_{\cA^\univ}: \cO_{F,(p)}^m \lra \Hom(\cA_{n,U^p(N_1,N_2)}^{(m),\ord}, (\cA^\univ/C^\univ[p^{N_1}])^\vee)_{\Z_{(p)}} \]
with the following properties. 
\begin{itemize}
\item If $a \in \cO_{F,(p)}$ then
\[ i^{(m)}_{\cA^\univ}(ax)=i^{\univ,\vee}({}^ca) \circ i^{(m)}_{\cA^\univ}(x). \]

\item If $(\beta,\delta)$ is a prime-to-$p$ quasi-isogeny
\[ (\cG^\univ,i^\univ,j^\univ,\lambda^\univ,[\eta^\univ]) \lra (\cG^{\univ,\prime},i^{\univ,\prime},j^{\univ,\prime},\lambda^{\univ,\prime},[\eta^{\univ,\prime}]), \]
then
\[ \beta^\vee \circ i^{(m)}_{(\cA^\univ)'}(x) = i^{(m)}_{\cA^\univ}(x). \]
In particular $i^{(m)}_{\cA^\univ}$ depends only on $\cA^\univ$ and not on $\cG^\univ$.

\item If $g \in G_n^{(m)}(\A^\infty)^\ord$ and $\gamma \in GL_m(\cO_{F,(p)})$ then 
\[ i^{(m)}_{A^\univ}(x) \circ g = i^{(m)}_{g^* A^\univ}(x) \]
and
\[ i^{(m)}_{A^\univ}(x) \circ \gamma = i^{(m)}_{\gamma^* A^\univ}(\gamma^{-1} x). \]

\item If $e_1,...,e_m$ denotes the standard basis of $\cO_{F,(p)}^m$ then
\[ i_{\cA^\univ}=||\eta^{p,\univ}_0||^{-1}( (\lambda(N_1)^\univ)^{-1}\circ i^{(m)}(e_1),..., (\lambda(N_1)^\univ)^{-1} \circ  i^{(m)}(e_m)) \]
is a prime-to-$p$ quasi-isogeny
\[  \cA_{n,U^p(N_1,N_2)}^{(m),\ord} \lra (\cA^\univ/C^\univ[p^{N_1}])^m. \]
Here $\lambda(N_1)^\univ$ refers to the prime-to-$p$ quasi-polarization $\cA^\univ/C[p^{N_1}] \ra (\cA^\univ/C[p^{N_1}])^\vee$ for which the composite
\[ \cA^\univ \lra \cA^\univ/C[p^{N_1}] \stackrel{\lambda(N_1)^\univ}{\lra} (\cA^\univ/C[p^{N_1}])^\vee  \lra \cA^{\univ,\vee} \]
equals $p^{N_1} \lambda^\univ$. 

We have
\[ \beta^{\oplus m} \circ i_{\cA^\univ}=i_{(\cA^\univ)'}. \]

The composite map
\[ \begin{array}{rcl} \eta_{n,U^p(N_1,N_2)}^{(m)}:\Hom_{\cO_{F}}(\cO_{F}^m,\Lambda_n) \otimes_\Z \A^{p,\infty} &\lra& V^p(\cA^\univ)^m \\ & \stackrel{p^{-N_1}}{\lra} &V^p(\cA^\univ/C^\univ[p^{N_1}])^m \\ & \lra& VA_{n,U}^{(m)}, \end{array} \]
where the first maps sends 
\[ f  \longmapsto (\eta^\univ(f(e_1)),...,\eta^\univ(f(e_m))) \]
and the third map sends
\[ x\longmapsto V(i_{A^\univ})^{-1} x, \]
is an isomorphism, which does not depend on the choice of $\cG^\univ$. It satisfies
\[ \eta_{n,U^p(N_1,N_2)}^{(m)} M^p = T^p\cA_{n,U^p(N_1,N_2)}^{(m),\ord}. \]
 \end{itemize}
 (See lemmas 5.2.4.7 and 7.1.2.1, propositions 5.2.4.13, 5.2.4.25 and 7.1.2.5, remarks 7.1.2.38 and 7.1.4.7, and theorem 7.1.4.1 of \cite{kw2}.)
 
 We deduce the following additional properties:
\begin{itemize}
\item If $g \in G_n^{(m)}(\A^\infty)^{\ord,\times}$ then the map $g:  \cA^{(m),\ord}_{n,U_1^p(N_{11},N_{12})} \ra \cA^{(m),\ord}_{n,U_2^p(N_{21},N_{22})}$ is etale. If further $N_{12}=N_{22}$, then it is finite etale.

\item If $g \in G_n^{(m)}(\A^\infty)^\ord$, if $N_{22}>0$, and if $p^{N_{12}-N_{22}}\nu(g_p) \in \Z_p^\times$ then
\[ g:\cA^{(m),\ord}_{n,U_1^p(N_{11},N_{12})} \lra \cA^{(m),\ord}_{n,U_2^p(N_{21},N_{22})}\]
 is finite. If $N_2>0$ then the finite flat map
 \[ \varsigma_p:  \cA^{(m),\ord}_{n,U^p(N_1,N_2+1)} \lra \cA^{(m),\ord}_{n,U^p(N_1,N_2)} \]
 has degree $p^{n(n+2m)[F^+:\Q]}$. 
 
  \item
 \[ i_{\cA^\univ} \circ g = i_{g^* \cA^\univ} \]
 and
 \[ i_{\cA^\univ} \circ \gamma = {}^t\gamma^{-1} \circ i_{\gamma^* \cA^\univ}. \]

 \end{itemize}

Also in this case define
\[ i_\lambda^{(m)} : \cO_{F,(p)}^m \otimes_{\cO_{F,(p)},c}\cO_{F,(p)}^m \lra \Hom(\cA_{n,U^p(N_1,N_2)}^{(m),\ord}, A_{n,U^p(N_1,N_2)}^{(m),\ord,\vee})_{\Z_{(p)}} \]
by 
\[ i_\lambda^{(m)}(x \otimes y) = ||\eta_0^{p,\univ}||^{-1} i^{(m)}_{\cA^\univ}(x)^\vee \circ (\lambda(N_1)^\univ)^{-1} \circ i^{(m)}_{\cA^\univ}(y). \]
This does not depend on the choice of $\cA^\univ$. We have
\[ i_\lambda^{(m)}(x \otimes y)^\vee = i_\lambda^{(m)}(y \otimes x). \]
Moreover
\[ (i_{\cA^\univ}^{-1})^\vee \circ i_\lambda^{(m)}(x \otimes y) \circ i_{\cA^\univ}^{-1} = (\lambda(N_1)^{\univ})^{ \oplus m} \circ i^\univ ({}^{c,t}xy).  \]
If $a \in (\cO_{F,(p)}^m \otimes_{\cO_{F,(p)},c} \cO_{F,(p)}^m)^{\sw =1}$ has image in $S(\cO_{F,(p)}^m)$ lying in $S(\cO_{F,(p)}^m)^{>0}$ then 
\[ (i_{\cA^\univ}^{-1})^\vee \circ i_\lambda^{(m)}(a) \circ i_{\cA^\univ}^{-1} = (\lambda(N_1)^\univ)^{\oplus m} \circ i^\univ (a') \]
for some matrix $a' \in M_{m \times m}(\cO_{F,(p)})^{t=c}$ all whose eigenvalues are positive real numbers. Thus $i_\lambda^{(m)}(a)$ is a quasi-polarization. (See the end of section 21 of \cite{mumav}.)

The completion of $\cA_{U^p(N_1,N_2)}^{(m),\ord}$ along its $\F_p$-fibre does not depend on $N_2$, so we will denote it
\[ \gA^{(m),\ord}_{U^p(N_1)}. \]
(See theorem 7.1.4.1 of \cite{kw2}.)
Then $\{ \gA^{(m),\ord}_{U^p(N)}\}$ is a system of $p$-adic formal schemes with a right $G_n^{(m)}(\A^{\infty})^\ord$-action and a left $GL_m(\cO_{F,(p)})$-action. There is an equivariant map 
\[  \{ \gA^{(m),\ord}_{n,U^p(N)}\} \lra \{ \gX^\ord_{n,(U')^p(N)} \}. \]
We will write $\barA^{(m),\ord}_{n,U^p(N)}$ for the reduced sub-scheme of  $\gA_{n,U^p(N)}^{(m),\ord}$. 

 \newpage \subsection{Some mixed Shimura varieties.}

If $\tU$ (resp. $\tU^p$) is a neat open compact subgroup of $\tG^{(m)}_n(\A^\infty)$ (resp. $\tG^{(m)}_n(\A^{p,\infty})$) we will denote by $S_{n,\tU}^{(m)}$ (resp. $\cS_{n,\tU^p}^{(m),\ord}$) the split torus over $\Spec \Q$ (resp. $\Spec \Z_{(p)}$) with 
\[ X_*(S_{n,\tU}^{(m)}) = Z(N_n^{(m)})(\Q) \cap \tU \subset \Herm^{(m)}(\Q) \]
(resp.
\[ X_*(\cS_{n,\tU^p}^{(m),\ord}) = Z(N_n^{(m)})(\Z_{(p)}) \cap \tU^p \subset \Herm^{(m)}(\Z_{(p)})). \]
If $g \in \tG^{(m)}_n(\A^{\infty})$ (resp. $\tG^{(m)}_n(\A^{\infty})^\ord$) and $\tU_2 \supset g^{-1}\tU_1g$ (resp. $\tU_2^p \supset g^{-1}\tU_1^pg$) we get a map
\[ g: S_{n,\tU_1}^{(m)} \lra S_{n,\tU_2}^{(m)} \]
(resp.
\[ g: \cS_{n,\tU_1^p}^{(m),\ord} \lra \cS_{n,\tU^p_2}^{(m),\ord}) \]
corresponding to
\[ ||\nu(g)||: X_*(S_{n,\tU_1}^{(m)}) \lra X_*(S_{n,\tU_2}^{(m)}) \]
(resp.
\[ ||\nu(g)||: X_*(\cS_{n,\tU_1^p}^{(m),\ord}) \lra X_*(\cS_{n,\tU^p_2}^{(m),\ord})), \]
where we think of the domain and codomain both as subspaces of $\Herm^{(m)}$.
If $\gamma \in GL_m(\Q)$ (resp. $GL_m(\Z_{(p)})$) and $\tU_2 \supset \gamma \tU_1$ (resp. $\tU_2^p \supset \gamma \tU_1^p$) we get a map
\[ \gamma: S_{n,\tU_1}^{(m)} \lra S_{n,\tU_2}^{(m)} \]
(resp.
\[ \gamma: \cS_{n,\tU_1^p}^{(m),\ord} \lra \cS_{n,\tU^p_2}^{(m),\ord}) \]
corresponding to
\[ \gamma: X_*(S_{n,\tU_1}^{(m)}) \lra X_*(S_{n,\tU_2}^{(m)}) \]
(resp.
\[ \gamma: X_*(\cS_{n,\tU_1^p}^{(m),\ord}) \lra X_*(\cS_{n,\tU^p_2}^{(m),\ord})), \]
where again we think of the domain and codomain both as subspaces of $\Herm^{(m)}$. If $m_1 \geq m_2$ and if $\tU_2$ (resp. $\tU_2^p$) is the image of $\tU_1$ (resp. $\tU_1^p$) in $\tG_n^{(m_2)}(\A^\infty)$ (resp. $\tG_n^{(m_2)}(\A^{p,\infty})$), then our chosen map $\Herm^{(m_1)} \ra \Herm^{(m_2)}$ induces a map
\[ S_{n,\tU_1}^{(m_1)} \lra S_{n,\tU_2}^{(m_2)} \]
(resp.
\[ \cS_{n,\tU_1^p}^{(m_1),\ord} \lra \cS_{n,\tU_2^p}^{(m_2),\ord}). \]

As $\tU$ runs over neat open compact subgroups of $\tG^{(m)}_n(\A^\infty)$, there  is a system of $S_{n,\tU}^{(m)}$-torsors 
\[ T^{(m)}_{n,\tU} = \underline{\Spec} \bigoplus_{\chi \in X^*(S^{(m)}_{n,\tU})} \cL^{(m)}_{n,\tU}(\chi) \]
over $A^{(m)}_{n,\tU}$  together with the following extra structures:
 \begin{itemize}
 \item If $g \in \tG_n^{(m)}(\A^\infty)$ and $\tU_1,\tU_2$ are neat open compact subgroups of $\tG_n^{(m)}(\A^\infty)$ with $\tU_2 \supset g^{-1}\tU_1g$ then there is a finite etale map 
 \[ g: T^{(m)}_{n,\tU_1} \lra T^{(m)}_{n,\tU_2} \]
 compatible with the maps $g:A^{(m)}_{n,\tU_1} \lra A^{(m)}_{n,\tU_2}$ and $g:S^{(m)}_{n,\tU_1} \lra S^{(m)}_{n,\tU_2}$. 
  
 \item If $\gamma \in GL_m(F)$ and $\tU_1,\tU_2$ are neat open compact subgroups of $\tG_n^{(m)}(\A^\infty)$ with $\tU_2 \supset \gamma \tU_1$ then there is a finite etale map
 \[ \gamma: T^{(m)}_{n,\tU_1} \lra T^{(m)}_{n,\tU_2},\]
  compatible with the maps $\gamma:A^{(m)}_{n,\tU_1} \lra A^{(m)}_{n,\tU_2}$ and $\gamma:S^{(m)}_{n,\tU_1} \lra S^{(m)}_{n,\tU_2}$.
  
  \item If $m_1 \geq m_2$ and $\tU_2$ is the image of $\tU_1$ in $\tG^{(m_2)}(\A^\infty)$, then there is a map
  \[ T^{(m_1)}_{n,\tU_1} \lra T^{(m_2)}_{n,\tU_2} \]
  compatible with the maps $S^{(m_1)}_{n,\tU_1} \lra S^{(m_2)}_{n,\tU_2}$ and $A^{(m_1)}_{n,\tU_1} \lra A^{(m_2)}_{n,\tU_2}$.
 \end{itemize}
 These enjoy the following properties:
  \begin{itemize}
  
  \item $g_1 \circ g_2=g_2g_1$ (i.e. this is a right action) and $\gamma_1 \circ \gamma_2 = \gamma_1 \gamma_2$ (i.e. this is a left action) and $\gamma \circ g = \gamma(g) \circ \gamma$.
  
 \item If $\tU_1 \subset \tU_2$ is an open normal subgroup of a neat open compact subgroup of $\tG_n^{(m)}(\A^\infty)$, then $T_{n,\tU_1}^{(m)}/T_{n,\tU_2}^{(m)}$ is Galois with group $\tU_2/\tU_1$.
 
 \item The maps $T^{(m_1)}_{n,\tU_1} \lra T^{(m_2)}_{n,\tU_2}$ are compatible with the actions of $\tG_n^{(m_1)}(\A^\infty)$ and $\tG_n^{(m_2)}(\A^\infty)$ and the map
 $\tG_n^{(m_1)}(\A^\infty)\ra\tG_n^{(m_2)}(\A^\infty)$, and also with the action of $Q_{m,(m-m')}(F)$.
 
\item 
Suppose that $\tU=U' \ltimes M$ with $U' \subset G_n(\A^\infty)$ and $M \subset N_n^{(m)}(\A^\infty)$. Also suppose that 
\[ \chi \in X^*(S^{(m)}_{n,\tU}) \subset S(F^m) \]
is sufficiently divisible. Then we can can find $a \in F^m \otimes_{F,c} F^m$ lifting $\chi$ such that
\[ i^{(m)}_\lambda(a): A^{(m)}_{n,\tU} \lra  (A^{(m)}_{n,\tU})^\vee \]
is an isogeny. For any such $a$
\[ \cL^{(m)}_{n,\tU}(\chi) = (1,i^{(m)}_\lambda(a))^* \cP_{A^{(m)}_{n,\tU}}. \]

\item If $\chi \in X^*(S^{(m)}_{n,\tU}) \cap S(F^m)^{>0}$ then $\cL_{n,\tU}^{(m)}(\chi)$ is relatively ample for $A^{(m)}_{n,\tU}/X_{n,\tU}$.

\item There are $\tG^{(m)}_{n}(\A^\infty)$ and $GL_m(F)$ equivariant homeomorphisms
\[ T^{(m)}_{n,\tU}(\C) \cong \tG^{(m)}_{n}(\Q) \backslash \tG^{(m)}_{n}(\A) \Herm^{(m)}(\C)/(\tU \times U_{n,\infty}^0A_n(\R)^0) . \] 
 \end{itemize}
 (See lemmas 1.3.2.25 and 1.3.2.72, and propositions 1.3.2.31, 1.3.2.45 and 1.3.2.90 of \cite{kw2}; section 3.6 of \cite{kwan}; and the second paragraph of section \ref{s32} above.)

Similarly as $\tU^p$ runs over neat open compact subgroups of $\tG^{(m)}_n(\A^{p,\infty})$ and $N_1,N_2$ run over integers with $N_2 \geq N_1 \geq 0$, there  is a system of $\cS_{n,\tU^p}^{(m),\ord}$-torsors 
\[ \cT^{(m),\ord}_{n,\tU^p(N_1,N_2)} = \underline{\Spec} \bigoplus_{\chi \in X^*(\cS^{(m),\ord}_{n,\tU^p(N_1,N_2)})} \cL^{(m),\ord}_{n,\tU^p(N_1,N_2)}(\chi) \]
over $\cA^{(m),\ord}_{n,\tU^p(N_1,N_2)}$  together with the following extra structures:
 \begin{itemize}
 \item If $g \in \tG_n^{(m)}(\A^{\infty})^\ord $ and $\tU_2^p(N_{21},N_{22}) \supset g^{-1}\tU_1^p(N_{11},N_{12})g$ then there is a  quasi-finite, flat map 
 \[ g: \cT^{(m),\ord}_{n,\tU_1^p(N_{11},N_{12})} \lra \cT^{(m),\ord}_{n,\tU_2^p(N_{21},N_{22})} \]
 compatible with the maps $g:\cA^{(m),\ord}_{n,\tU_1^p(N_{11},N_{12})} \lra \cA^{(m),\ord}_{n,\tU_2^p(N_{21},N_{22})}$ and $g:\cS^{(m),\ord}_{n,\tU_1^p} \lra \cS^{(m),\ord}_{n,\tU_2^p}$.

 \item If $\gamma \in GL_m(\cO_{F,(p)})$ and $\tU_2^p \supset \gamma \tU_1^p$ then there is a finite etale map
 \[ \gamma: \cT^{(m),\ord}_{n,\tU_1^p(N_1,N_2)} \lra \cT^{(m),\ord}_{n,\tU_2^p(N_1,N_2)},\]
  compatible with the maps 
  \[ \gamma:\cA^{(m),\ord}_{n,\tU_1^p(N_1,N_2)} \lra \cA^{(m)}_{n,\tU_2^p(N_1,N_2)}\]
   and 
   \[ \gamma:\cS^{(m),\ord}_{n,\tU_1^p} \lra S^{(m)}_{n,\tU_2^p}. \]
   
   \item If $m_1 \geq m_2$ and $\tU_2^p$ is the image of $\tU_1^p$ in $\tG^{(m_2)}(\A^{p,\infty})$, then there is a map
  \[ \cT^{(m_1),\ord}_{n,\tU_1^p(N_1,N_2)} \lra \cT^{(m_2),\ord}_{n,\tU_2^p(N_1,N_2)} \]
  compatible with the maps $\cS^{(m_1),\ord}_{n,\tU_1^p} \lra \cS^{(m_2),\ord}_{n,\tU_2^p}$ and $\cA^{(m_1),\ord}_{n,\tU_1^p(N_1,N_2)} \lra \cA^{(m_2),\ord}_{n,\tU_2^p(N_1,N_2)}$.
 \end{itemize}
 These enjoy the following properties:
  \begin{itemize}
  \item $g_1 \circ g_2=g_2g_1$ (i.e. this is a right action) and $\gamma_1 \circ \gamma_2 = \gamma_1 \gamma_2$ (i.e. this is a left action) and $\gamma \circ g = \gamma(g) \circ \gamma$.
  
  \item  If $g \in \tG_n^{(m)}(\A^\infty)^{\ord,\times}$ then the map $g: \cT^{(m),\ord}_{n,\tU_1^p(N_{11},N_{12})} \ra \cT^{(m),\ord}_{n,\tU_2^p(N_{21},N_{22})}$ is etale. If further $N_{12}=N_{22}$,  then it is finite etale.
  
\item The maps $\cT^{(m_1),\ord}_{n,\tU_1^p(N_1,N_2)} \lra \cT^{(m_2),\ord}_{n,\tU_2^p(N_1,N_2)}$ are compatible with the actions of $\tG_n^{(m_1)}(\A^\infty)^\ord$ and $\tG_n^{(m_2)}(\A^\infty)^\ord$ and the map
 $\tG_n^{(m_1)}(\A^\infty)\ra\tG_n^{(m_2)}(\A^\infty)$, and with the action of $Q_{m,(m-m')}(\cO_{F,(p)})$.

 \item If $\tU_1^p \subset \tU_2^p$ is an open normal subgroup of a neat open compact of $\tG_n^{(m)}(\A^{p,\infty})$, and if $N_{11} \geq N_{21}$ then $\cT_{n,\tU_1^p(N_{11},N_2)}^{(m),\ord}/\cT_{n,\tU_2^p(N_{21},N_2)}^{(m),\ord}$ is Galois with group $\tU_2^p(N_{21})/\tU_1^p(N_{11})$.

\item If $g \in G_n^{(m)}(\A^\infty)^\ord$, if $N_{22}>0$, and if $p^{N_{12}-N_{22}}\nu(g_p) \in \Z_p^\times$, then $g: \cT^{(m),\ord}_{n,\tU_1^p(N_{11},N_{12})} \ra \cT^{(m),\ord}_{n,\tU_2^p(N_{21},N_{22})}$ is
finite. If $N_{2}>0$ then the finite flat map
\[ \varsigma_p: \cT^{(m),\ord}_{n,\tU_1^p(N_{1},N_{2}+1)} \ra \cT^{(m),\ord}_{n,\tU_2^p(N_{1},N_{2})} \]
has degree $p^{(n+m)^2[F^+:\Q]}$.
 
 \item On the $\F_p$-fibre 
 \[ \varsigma_p: \cT^{(m),\ord}_{n,\tU^p(N_1,N_2+1)}\times \Spec \F_p \lra \cT^{(m),\ord}_{n,\tU^p(N_1,N_2)} \times \Spec \F_p \]
 equals the composition of the absolute Frobenius map with the forgetful map (for any $N_2 \geq N_1 \geq 0$). 
    
\item 
Suppose that $\tU^p=(U^p)' \ltimes M^p$ with $(U^p)' \subset G_n(\A^{p,\infty})$ and $M^p \subset N_n^{(m)}(\A^{p,\infty})$. Also suppose that 
\[ \chi \in X^*(\cS^{(m),\ord}_{n,\tU^p}) \subset S(\cO_{F,(p)}^m) \]
is sufficiently divisible. Then we can can find $a \in \cO_{F,(p)}^m \otimes_{\cO_{F,(p)}} \cO_{F,(p)}^m$ lifting $\chi$ such that
\[ i^{(m)}_\lambda(a): \cA^{(m),\ord}_{n,\tU^p(N_1,N_2)} \lra  (\cA^{(m),\ord}_{n,\tU^p(N_1,N_2)})^\vee \]
is an isogeny. For any such $a$
\[ \cL^{(m),\ord}_{n,\tU^p(N_1,N_2)}(\chi) = (1,i^{(m)}_\lambda(a))^* \cP_{\cA^{(m),\ord}_{n,\tU^p(N_1,N_2)}}. \]

\item If $\chi \in X^*(\cS^{(m)}_{n,\tU^p}) \cap S(\cO_{F,(p)}^m)^{>0}$ then $\cL_{n,\tU^p(N_1,N_2)}^{(m),\ord}(\chi)$ is relatively ample for $\cA^{(m),\ord}_{n,\tU^p(N_1,N_2)}/\cX_{n,\tU^p(N_1,N_2)}^\ord$.

\item There are natural identifications 
\[ \cT^{(m),\ord}_{n,\tU^p(N_1,N_2)} \times \Spec \Q \cong T^{(m)}_{n,U^p(N_1,N_2)}. \]
These identifications are compatible with the identifications 
\[ \cA_{n,\tU^p(N_1,N_2)}^{(m),\ord} \times \Spec \Q \cong A^{(m)}_{n,\tU^p(N_1,N_2)}\]
 and the maps 
 \[ \cT^{(m),\ord}_{n,\tU^p(N_1,N_2)} \lra \cA_{n,\tU^p(N_1,N_2)}^{(m),\ord}\]
  and 
  \[ T^{(m)}_{n,\tU^p(N_1,N_2)} \lra A^{(m)}_{n,\tU^p(N_1,N_2)}.\]
   The identifications are also equivariant for the actions actions of the semi-group $\tG_n^{(m)}(\A^{\infty})^\ord$ and the group $GL_m(\cO_{F,(p)})$. 
 \end{itemize}
 (See lemmas 5.2.4.26 and 7.1.2.22, propositions 5.2.4.30, 5.2.4.41 and 7.1.2.36, and remark 7.1.2.38 of \cite{kw2}.)

\newpage \subsection{Vector bundles.}\label{vb1}

Suppose that $U$ is a neat open compact subgroup of $G_n(\A^\infty)$.
We will let $\Omega_{n,U}$ denote the pull back by the identity section of the sheaf of relative differentials $\Omega^1_{A^\univ/X_{n,U}}$. This is a locally free sheaf of rank $n[F:\Q]$. Up to unique isomorphism its definition does not depend on the choice of $A^\univ$. (Because, by the neatness of $U$, there is a unique isogeny between any two universal four-tuples $(A^\univ, i^\univ, \lambda^\univ,[\eta^\univ])$.) The system of sheaves $\{ \Omega_{n,U}\}$ has an action of $G_n(\A^\infty)$. There is a natural isomorphism between $\Omega^1_{A^\univ/X_{n,U}}$ and the pull back of $\Omega_{n,U}$ from $X_{n,U}$ to $A^\univ$.

Similarly, if $\pi: A^\univ \ra X_{n,U}$ is the structural map, then the sheaf
\[ R^i\pi_* \Omega^j_{A^\univ/X_{n,U}} \cong (\wedge^j \Omega_{n,U}) \otimes R^i\pi_* \cO_{A^\univ} \]
 is locally free and canonically independent of the choice of $A^\univ$. These sheaves again have an action of $G_n(\A^\infty)$.

We will also write $\Xi_{n,U}=\cO_{X_{n,U}}(||\nu||)$ for the sheaf $\cO_{X_{n,U}}$ but with the $G_n(\A^\infty)$-action multiplied by $||\nu||$. 

The line bundle $(1,\lambda^\univ)^* \cP_{A^\univ}$ is represented by an element of 
 \[ \begin{array}{rcl} [ (1,\lambda^\univ)^* \cP_{A^{\univ}}] &\in& H^1(A^{\univ},\cO^\times_{A^{\univ}}) \\ &\lra & H^0(X_{n,U},R^1\pi_* \cO^\times_{A^{\univ}}) \\ &\stackrel{d \log}{\lra}& H^0(X_{n,U}, R^1\pi_* \Omega^1_{A^\univ/X_{n,U}}). \end{array} \] 
We obtain an embedding
\[ \Xi_{n,U} \into R^1\pi_* \Omega^1_{A^\univ/X_{n,U}} \]
sending $1$  to $||\eta^\univ|| [ (1,\lambda^\univ)^* \cP_{A^{\univ}}]$. (See section \ref{introshim} for the definition of $||\eta^\univ||$.) These maps are compatible with the isomorphisms 
\[ R^1\pi_* \Omega^1_{A^\univ/X_{n,U}} \liso R^1\pi_* \Omega^1_{A^{\univ,\prime}/X_{n,U}}\]
 induced by the unique isogeny between two universal $4$-tuples. They are also $G_n(\A^\infty)$-equivariant. 

The induced maps
\[ \begin{array}{rcl} \Hom(\Omega_{n,U}, \Xi_{n,U}) &\into& \Hom(\Omega_{n,U}, R^1\pi_* \Omega^1_{A^\univ/X_{n,U}})  \\ & \stackrel{\sim}{\lla} & \Hom(\Omega_{n,U},\Omega_{n,U} \otimes R^1\pi_* \cO_{A^\univ}) \\ & \stackrel{\tr}{\lra} & R^1\pi_* \cO_{A^\univ} \end{array} \]
 are $G_n(\A^\infty)$-equivariant isomorphisms, independent of the choice of $A^\univ$. Moreover the short exact sequence
 \[ (0) \lra \Omega^1_{X_{n,U}} \otimes \cO_{A^\univ} \lra \Omega^1_{A^\univ} \lra \Omega_{n,U} \otimes \cO_{A^\univ} \lra (0) \]
 gives rise to a map
 \[ \begin{array}{rcl} \Omega_{n,U} &\lra& \Omega^1_{X_{n,U}} \otimes R^1\pi_* \cO_{A^\univ}  \\ & \stackrel{\sim}{\lla} & \Omega^1_{X_{n,U}} \otimes \Hom(\Omega_{n,U}, \Xi_{n,U}) \end{array} \]
 and hence to a map
 \[ \Omega_{n,U}^{\otimes 2} \lra \Omega^1_{X_{n,U}}\otimes \Xi_{n,U}. \]
 These maps do not depend on the choice of $A^\univ$ and are $G_n(\A^\infty)$-equivariant. They further induce $G_n(\A^\infty)$-equivariant isomorphisms
 \[ S(\Omega_{n,U}) \liso \Omega^1_{X_{n,U}}\otimes \Xi_{n,U}, \]
 which again do not depend on the choice of $A^\univ$. (See for instance propositions 2.1.7.3 and 2.3.5.2 of \cite{kw1}. This is referred to as the `Kodaira-Spencer isomorphism'.)
 
 Let $\cE_{U}$ denote the principal
$L_{n,(n)}$-bundle on $X_{n,U}$ in the
Zariski topology defined by setting, for $W \subset X_{n,U}$ a Zariski open, 
$\cE_{U}(W)$  to
be the set of pairs $(\xi_0,\xi_1)$, where
\[ \xi_0: \Xi_{n,U}|_W \liso \cO_W \]
and
\[ \xi_1: \Omega_{n,U} \liso \Hom_\Q( V_n/V_{n,(n)}, \cO_W). \]
We define the $L_{n,(n)}$-action on $\cE_{U}$ by
\[ h(\xi_0,\xi_1)=(\nu(h)^{-1}\xi_0, (\circ h^{-1}) \circ \xi_1). \] 
The inverse system $\{ \cE_{U} \}$ has an action of $G_n(\A^\infty)$. 

Suppose that $R_0$ is an irreducible noetherian $\Q$-algebra and that $\rho$ is a
representation of $L_{n,(n)}$ on a finite, locally free $R$-module $W_\rho$. We define a locally
free sheaf $\cE_{U,\rho}$ over
$X_{n,U} \times \Spec R_0$ by setting $\cE_{U,\rho}(W)$ to be the set of $L_{n,(n)}(\cO_W)$-equivariant
maps of Zariski sheaves of sets
\[ \cE_{U}|_W \ra W_\rho \otimes_{R_0} \cO_W. \]
Then $\{ \cE_{U,\rho} \}$ is a system of locally free sheaves with $G_n(\A^\infty)$-action over the system of schemes 
$\{ X_{n,U} \times \Spec R_0\}$. If $g \in G_n(\A^\infty)$, then the natural map
\[ g^*\cE_{U,\rho} \lra \cE_{U',\rho} \]
is an isomorphism.

In the case $R_0=\C$, the holomorphic vector bundle on $X_{n,U}(\C)$ associated to $\cE_{U,\rho}$ is 
\[ \gE_{U,\rho} = G_n(\Q) \backslash \left( G_n(\A^\infty)/U \times \gE_\rho \right) \]
over
\[ X_{n,U}(\C) = G_n(\Q) \backslash \left( G_n(\A^\infty)/U \times \gH_n^\pm \right). \]
(See section \ref{s1.1} for the definition of the holomorphic vector bundle $\gE_\rho/\gH_n^{\pm}$.)

Note that
\[ \cE_{U,\Std^\vee} \cong \Omega_{n,U} \]
and
\[ \cE_{U,\nu^{-1}} \cong \Xi_{n,U} \]
and
\[ \cE_{U,\wedge^{n[F:\Q]} \Std^\vee} \cong \omega_U \]
and
\[ \cE_{U,\KS} \cong \Omega^1_{X_{n,U}}. \]
(See section \ref{mps} for the definition of the representation $\KS$.)
 
 Suppose now that $U$ is a neat open compact subgroup of $G_n^{(m)}(\A^\infty)$ with image $U'$ in $G_n(\A^\infty)$. We will let $\Omega_{n,U}^{(m)}$ denote the pull back by the identity section of the sheaf of relative differentials $\Omega^1_{G^\univ/A_{n,U}^{(m)}}$. This is a locally free sheaf of rank $(n+m)[F:\Q]$. Up to unique isomorphism its definition does not depend on the choice of $G^\univ$. The system of sheaves $\{ \Omega^{(m)}_{n,U}\}$ has actions of $G_n^{(m)}(\A^\infty)$ and of $GL_m(F)$. Moreover there is an exact sequence
 \[ (0) \lra \pi_{A_n^{(m)}/X_n}^* \Omega_{n,U'} \lra \Omega_{n,U}^{(m)} \lra F^m \otimes_\Q \cO_{A_{n,U}^{(m)}} \lra (0) \] 
 which is equivariant for the actions of $G_n^{(m)}(\A^\infty)$ and $GL_m(F)$. 

 Let $\cE_{U}^{(m)}$ denote the principal
$R^{(m)}_{n,(n)}$-bundle on $A^{(m)}_{n,U}$ in the
Zariski topology defined by setting, for $W \subset A^{(m)}_{n,U}$ a Zariski open, 
$\cE_{U}^{(m)}(W)$  to
be the set of pairs $(\xi_0,\xi_1)$, where
\[ \xi_0: \Xi_{n,U}|_W \liso \cO_W \]
and
\[ \xi_1: \Omega_{n,U}^{(m)} \liso \Hom_\Q( V_n/V_{n,(n)}\oplus \Hom_\Q(F^m,\Q) , \cO_W)  \]
satisfies
\[ \xi_1: \Omega_{n,U} \liso \Hom_\Q( V_n/V_{n,(n)}, \cO_W) \]
and induces the canonical isomorphism
\[ F^m \otimes_\Q \cO_W \lra \Hom_\Q(\Hom_\Q(F^m,\Q) , \cO_W). \]
We define the $R^{(m)}_{n,(n)}$-action on $\cE_{U}^{(m)}$ by
\[ h(\xi_0,\xi_1)=(\nu(h)^{-1}\xi_0, (\circ h^{-1}) \circ \xi_1). \] 
The inverse system $\{ \cE_{U}^{(m)} \}$ has an action of $G_n^{(m)}(\A^\infty)$ and of $GL_m(F)$. 

Suppose that $R_0$ is an irreducible noetherian $\Q$-algebra and that $\rho$ is a
representation of $R^{(m)}_{n,(n)}$ on a finite, locally free $R$-module $W_\rho$. We define a locally
free sheaf $\cE^{(m)}_{U,\rho}$ over
$A^{(m)}_{n,U} \times \Spec R_0$ by setting $\cE^{(m)}_{U,\rho}(W)$ to be the set of $R^{(m)}_{n,(n)}(\cO_W)$-equivariant
maps of Zariski sheaves of sets
\[ \cE^{(m)}_{U}|_W \lra W_\rho \otimes_{R_0} \cO_W. \]
Then $\{ \cE^{(m)}_{U,\rho} \}$ is a system of locally free sheaves with both $G^{(m)}_n(\A^\infty)$-action and $GL_m(F)$-action over the system of schemes 
$\{ A^{(m)}_{n,U} \times \Spec R_0\}$. 
If $g \in G_n^{(m)}(\A^\infty)$ and $\gamma \in GL_m(F)$, then the natural maps
\[ g^*\cE_{U,\rho}^{(m)} \lra \cE_{U',\rho}^{(m)} \]
and
\[ \gamma^*\cE_{U,\rho}^{(m)} \lra \cE_{U',\rho}^{(m)} \]
are isomorphisms.
If $\rho$ factors through $R^{(m)}_{n,(n)} \onto L_{n,(n)}$ then $\cE^{(m)}_{U,\rho}$ is canonically isomorphic to the pull-back of $\cE_{U,\rho}$ from $X_{n,U}$. In general $W_\rho$ has a filtration by $R_{n,(n)}^{(m)}$-invariant local direct-summands  such that the action of $R^{(m)}_{n,(n)}$ on each graded piece factors through $L_{n,(n)}$. Thus $\cE^{(m)}_{U,\rho}$ has a $G_n^{(m)}(\A^\infty)$ and $GL_m(F)$ invariant filtration by local direct summands such that each graded piece is the pull back of some $\cE_{U,\rho'}$ from $X_{n,U}$.

 Similarly suppose that $U^p$ is a neat open compact subgroup of $G_n(\A^{p,\infty})$, and that $N_2 \geq N_1 \geq 0$ are integers.
We will let $\Omega_{n,U^p(N_1,N_2)}^\ord$ denote the pull back by the identity section of $\Omega^1_{\cA^\univ/\cX_{n,U^p(N_1,N_2)}^\ord}$. This is a locally free sheaf of rank $n[F:\Q]$. Up to unique isomorphism its definition does not depend on the choice of $\cA^\univ$.  (Because, by the neatness of $U^p$, there is a unique prime-to-$p$ isogeny between any two universal four-tuples $(\cA^\univ, i^\univ, \lambda^\univ,[\eta^\univ])$.) The system of sheaves $\{ \Omega_{n,U^p(N_1,N_2)}^\ord \}$ has an action of $G_n(\A^\infty)^\ord$. There is a natural isomorphism between $\Omega^1_{\cA^\univ/\cX^\ord_{n,U^p(N_1,N_2)}}$ and the pull back of $\Omega_{n,U^p(N_1,N_2)}^\ord$.

We will also write $\Xi_{n,U^p(N_1,N_2)}=\cO_{\cX^\ord_{n,U^p(N_1,N_2)}}(||\nu||)$ for the sheaf $\cO_{\cX_{n,U^p(N_1,N_2)}^\ord}$ but with the $G_n(\A^\infty)^\ord$-action multiplied by $||\nu||$. 

The line bundle $(1,\lambda^\univ)^* \cP_{\cA^\univ}$ is represented by an element of 
 \[ \begin{array}{rcl} [ (1,\lambda^\univ)^* \cP_{\cA^{\univ}}] &\in& H^1(\cA^{\univ},\cO^\times_{\cA^{\univ}}) \\ &\lra & H^0(\cX_{n,U^p(N_1,N_2)},R^1\pi_* \cO^\times_{\cA^{\univ}}) \\ &\stackrel{d \log}{\lra}& H^0(\cX_{n,U^p(N_1,N_2)}^\ord, R^1\pi_* \Omega^1_{\cA^\univ/\cX^\ord_{n,U^p(N_1,N_2)}}). \end{array} \]
 We obtain an embedding
\[ \Xi_{n,U^p(N_1,N_2)}^\ord \into R^1\pi_* \Omega^1_{\cA^\univ/\cX^\ord_{n,U^p(N_1,N_2)}} \]
sending $1$  to $||\eta^\univ|| [ (1,\lambda^\univ)^* \cP_{\cA^{\univ}}]$. These maps are compatible with the isomorphisms 
\[ R^1\pi_* \Omega^1_{\cA^\univ/\cX^\ord_{n,U^p(N_1,N_2)}} \iso R^1\pi_* \Omega^1_{\cA^{\univ,\prime}/\cX^\ord_{n,U^p(N_1,N_2)}} \]
 induced by the unique prime-to-$p$ isogeny between two universal $4$-tuples. They are also $G_n(\A^\infty)^\ord$-equivariant.

The induced maps
\[ \begin{array}{rl} & \Hom(\Omega_{n,U^p(N_1,N_2)}^\ord, \Xi_{n,U^p(N_1,N_2)}^\ord) \\ 
\into& \Hom(\Omega^\ord_{n,U^p(N_1,N_2)}, R^1\pi_* \Omega^1_{\cA^\univ/\cX^\ord_{n,U^p(N_1,N_2)}})  \\ 
 \stackrel{\sim}{\lla} & \Hom(\Omega^\ord_{n,U^p(N_1,N_2)},\Omega_{n,U^p(N_1,N_2)}^\ord \otimes R^1\pi_* \cO_{\cA^\univ}) \\ 
  \stackrel{\tr}{\lra} & R^1\pi_* \cO_{\cA^\univ} \end{array} \]
 are $G_n(\A^\infty)^\ord$-equivariant 
 isomorphisms, independent of the choice of $\cA^\univ$. Moreover the short exact sequence
 \[ (0) \lra \Omega^1_{\cX^\ord_{n,U^p(N_1,N_2)}} \otimes \cO_{\cA^\univ} \lra \Omega^1_{\cA^\univ} \lra \Omega^\ord_{n,U^p(N_1,N_2)} \otimes \cO_{\cA^\univ} \lra (0) \]
 gives rise to a map
 \[ \begin{array}{rcl} \Omega_{n,U^p(N_1,N_2)}^\ord &\lra& \Omega^1_{\cX^\ord_{n,U^p(N_1,N_2)}} \otimes R^1\pi_* \cO_{\cA^\univ}  \\ & \stackrel{\sim}{\lla} & \Omega^1_{\cX^\ord_{n,U^p(N_1,N_2)}} \otimes \Hom(\Omega^\ord_{n,U^p(N_1,N_2)}, \Xi^\ord_{n,U^p(N_1,N_2)}) \end{array} \]
 and hence to a map
 \[ (\Omega^\ord_{n,U^p(N_1,N_2)})^{\otimes 2} \lra \Omega^1_{\cX^\ord_{n,U^p(N_1,N_2)}}\otimes \Xi^\ord_{n,U^p(N_1,N_2)}. \]
 These maps do not depend on the choice of $\cA^\univ$ and are $G_n(\A^\infty)^\ord$-equivariant.
  They further induce $G_n(\A^\infty)^\ord$ 
 isomorphisms
 \[ S(\Omega^\ord_{n,U^p(N_1,N_2)}) \liso \Omega^1_{\cX^\ord_{n,U^p(N_1,N_2)}}\otimes \Xi^\ord_{n,U^p(N_1,N_2)}, \]
 which again do not depend on the choice of $\cA^\univ$. (See for instance  proposition 3.4.3.3 of \cite{kw2}.)
 
 Let $\cE_{U^p(N_1,N_2)}^{\ord}$ denote the principal
$L_{n,(n)}$-bundle on $\cX^\ord_{n,U^p(N_1,N_2)}$ in the
Zariski topology defined by setting, for $W \subset \cX^\ord_{n,U^p(N_1,N_2)}$) a Zariski open, 
$\cE^{\ord}_{U^p(N_1,N_2)}(W)$ to
be the set of pairs $(\xi_0,\xi_1)$, where
\[ \xi_0: \Xi^\ord_{\cA^\ord/\cX^\ord,U^p(N_1,N_2)}|_W \liso \cO_W) \]
and
\[ \xi_1: \Omega^\ord_{\cA^\ord/\cX^\ord,U^p(N_1,N_2)} \liso \Hom_{\Z} 
(\Lambda_n/\Lambda_{n,(n)}, \cO_W)). \]
We define the $L_{n,(n)}$-action on $\cE_{U^p(N_1,N_2)}^{\ord}$) by
\[ h(\xi_0,\xi_1)=(\nu(h)^{-1}\xi_0, (\circ h^{-1}) \circ \xi_1). \] 
The inverse system $\{ \cE_{U^p(N_1,N_2)}^{\ord} \}$ has an action of $G_n(\A^\infty)^\ord$. 

Suppose that $R_0$ is an irreducible noetherian $\Z_{(p)}$-algebra and that $\rho$ is a
representation of $L_{n,(n)}$ on a finite, locally free $R$-module $W_\rho$. We define a locally
free sheaf $\cE_{U^p(N_1,N_2),\rho}^{\ord}$ over
$\cX^\ord_{n,U^p(N_1,N_2)}\times \Spec R_0$ by setting $\cE_{U^p(N_1,N_2),\rho}^{\ord}(W)$ to be the set of $L_{n,(n)}(\cO_W)$-equivariant
maps of Zariski sheaves of sets
\[ \cE^{\ord}_{U^p(N_1,N_2)}|_W \ra W_\rho \otimes_{R_0} \cO_W.\]
Then $\{ \cE^{\ord}_{U^p(N_1,N_2),\rho} \}$ is a system of locally free sheaves with $G_n(\A^\infty)^\ord$-action over the system of schemes 
$\{ \cX^\ord_{n,U^p(N_1,N_2),\Delta} \times \Spec R_0 \} $. The pull-back of $\cE^{\ord}_{U^p(N_1,N_2),\rho}$ to 
\[ \cX^\ord_{n,U^p(N_1,N_2),\Delta} \times \Spec R_0[1/p] \]
 is canonically identified with the sheaf $\cE_{U^p(N_1,N_2),\rho \otimes_{R_0} R_0[1/p]}$. This identification is $G_n(\A^\infty)^\ord$-equivariant. If $g \in G_n(\A^\infty)^{\ord,\times}$, then the natural map
\[ g^*\cE_{U^p(N_1,N_2),\rho} \lra \cE_{(U^p)'(N_1',N_2'),\rho} \]
is an isomorphism.

Note that
\[ \cE^{\ord}_{U^p(N_1,N_2),\Std^\vee} \cong \Omega^\ord_{\cA^\ord/\cX^\ord,U^p(N_1,N_2)} \]
and
\[ \cE^{\ord}_{U^p(N_1,N_2),\nu^{-1}} \cong \Xi^\ord_{\cA^\ord/\cX^\ord,U^p(N_1,N_2)} \]
and
\[ \cE^{\ord}_{U^p(N_1,N_2),\wedge^{n[F:\Q]} \Std^\vee} \cong \omega_{U^p(N_1,N_2)} \]
and
\[ \cE^{\ord}_{U^p(N_1,N_2),\KS} \cong \Omega^1_{\cX^\ord_{n,U^p(N_1,N_2)}}. \]

  Suppose now that $U^p$ is a neat open compact subgroup of $G_n^{(m)}(\A^{p,\infty})$ with image $(U^p)'$ in $G_n(\A^{p,\infty})$. We will let $\Omega_{n,U^p(N_1,N_2)}^{(m),\ord}$ denote the pull back by the identity section of the sheaf of relative differentials $\Omega^1_{\cG^\univ/\cA_{n,U^p(N_1,N_2)}^{(m),\ord}}$. This is a locally free sheaf of rank $(n+m)[F:\Q]$. Up to unique isomorphism its definition does not depend on the choice of $\cG^\univ$. The system of sheaves $\{ \Omega^{(m),\ord}_{n,U^p(N_1,N_2)}\}$ has actions of $G_n^{(m)}(\A^\infty)^\ord$ and of $GL_m(\cO_{F,(p)})$. Moreover there is an exact sequence
 \[ (0) \ra \pi_{\cA_n^{(m),\ord}/\cX^\ord_n}^* \Omega^\ord_{n,(U^p)'(N_1,N_2)} \ra \Omega_{n,U^p(N_1,N_2)}^{(m)} \ra \cO_{F,(p)}^m \otimes_\Q \cO_{\cA_{n,U^p(N_1,N_2)}^{(m)}} \ra (0) \] 
 which is equivariant for the actions of $G_n^{(m)}(\A^\infty)^\ord$ and $GL_m(\cO_{F,(p)})$. 

 Let $\cE_{U^p(N_1,N_2)}^{(m),\ord}$ denote the principal
$R^{(m)}_{n,(n)}$-bundle on $\cA^{(m),\ord}_{n,U^p(N_1,N_2)}$ in the
Zariski topology defined by setting, for $W \subset \cA^{(m),\ord}_{n,U^p(N_1,N_2)}$ a Zariski open, 
$\cE_{U^p(N_1,N_2)}^{(m),\ord}(W)$  to
be the set of pairs $(\xi_0,\xi_1)$, where
\[ \xi_0: \Xi^\ord_{n,U^p(N_1,N_2)}|_W \liso \cO_W \]
and
\[ \xi_1: \Omega_{n,U^p(N_1,N_2)}^{(m),\ord} \liso \Hom( \Lambda_n/\Lambda_{n,(n)}\oplus \Hom(\cO_F^m,\Z) , \cO_W)  \]
satisfies
\[ \xi_1: \Omega_{n,U^p(N_1,N_2)}^\ord \liso \Hom( \Lambda_n/\Lambda_{n,(n)}, \cO_W) \]
and induces the canonical isomorphism
\[ \cO_{F,(p)}^m \otimes_{\Z_{(p)}} \cO_W \lra \Hom(\Hom(\cO_F^m,\Z) , \cO_W). \]
We define the $R^{(m)}_{n,(n)}$-action on $\cE_{U^p(N_1,N_2)}^{(m),\ord}$ by
\[ h(\xi_0,\xi_1)=(\nu(h)^{-1}\xi_0, (\circ h^{-1}) \circ \xi_1). \] 
The inverse system $\{ \cE_{U^p(N_1,N_2)}^{(m),\ord} \}$ has an action of $G_n^{(m)}(\A^\infty)^\ord$ and of $GL_m(\cO_{F,(p)})$. 

Suppose that $R_0$ is an irreducible noetherian $\Z_{(p)}$-algebra and that $\rho$ is a
representation of $R^{(m)}_{n,(n)}$ on a finite, locally free $R_0$-module $W_\rho$. We define a locally
free sheaf $\cE^{(m),\ord}_{U^p(N_1,N_2),\rho}$ over
$\cA^{(m),\ord}_{n,U^p(N_1,N_2)} \times \Spec R_0$ by setting $\cE^{(m)}_{U,\rho}(W)$ to be the set of $R^{(m)}_{n,(n)}(\cO_W)$-equivariant
maps of Zariski sheaves of sets
\[ \cE^{(m),\ord}_{U^p(N_1,N_2)}|_W \lra W_\rho \otimes_{R_0} \cO_W. \]
Then $\{ \cE^{(m),\ord}_{U^p(N_1,N_2),\rho} \}$ is a system of locally free sheaves with $G^{(m)}_n(\A^\infty)^\ord$-action and $GL_m(\cO_{F,(p)})$-action over the system of schemes 
$\{ \cA^{(m),\ord}_{n,U^p(N_1,N_2)} \times \Spec R_0\}$. 
If $g \in G_n^{(m)}(\A^\infty)^{\ord,\times}$ and $\gamma \in GL_m(\cO_{F,(p)})$, then the natural maps
\[ g^*\cE_{U^p(N_1,N_2),\rho}^{(m),\ord} \lra \cE_{(U^p)'(N_1',N_2'),\rho}^{(m),\ord} \]
and
\[ \gamma^*\cE_{U^p(N_1,N_2),\rho}^{(m),\ord} \lra \cE_{(U^p)'(N_1',N_2'),\rho}^{(m),\ord} \]
are isomorphisms.
If $\rho$ factors through $R^{(m)}_{n,(n)} \onto L_{n,(n)}$ then $\cE^{(m),\ord}_{U^p(N_1,N_2),\rho}$ is canonically isomorphic to the pull-back of $\cE^\ord_{U^p(N_1,N_2),\rho}$ from $\cX^\ord_{n,U^p(N_1,N_2)}$. In general $W_\rho$ has a filtration by $R_{n,(n)}^{(m)}$-invariant local direct-summands  such that the action of $R^{(m)}_{n,(n)}$ on each graded piece factors through $L_{n,(n)}$. Thus $\cE^{(m),\ord}_{U^p(N_1,N_2),\rho}$ has a $G_n^{(m)}(\A^\infty)$ and $GL_m(\cO_{F,(p)})$ invariant filtration by local direct summands such that each graded piece is the pull back of some $\cE^\ord_{U^p(N_1,N_2),\rho'}$ from $\cX^\ord_{n,U^p(N_1,N_2)}$.

If $m \geq m'$ and if $U$ is a neat open compact subgroup of $G^{(m)}_n(\A^\infty)$ with image $U'$ in $G_n^{(m')}(\A^\infty)$ then the sheaf 
\[ R^j\pi_{A^{(m)}_{n}/A^{(m')}_n,*} \Omega^i_{A^{(m)}_{n,U}/A^{(m')}_{n,U'}} \]
depends only on $U'$ and not on $U$. We will denote it
\[ (R^j\pi_* \Omega^i_{A^{(m)}_n/A^{(m')}_n})_{U'}. \]
If $g \in G_n^{(m)}(\A^\infty)$ and $g^{-1}U_1g \subset U_2$ then there is a natural isomorphism
\[ g:  (g')^*(R^j\pi_* \Omega^i_{A^{(m)}_n/A^{(m')}_n})_{U_2'} \liso (R^j\pi_* \Omega^i_{A^{(m)}_n/A^{(m')}_n})_{U_1'}, \]
where $g'$ (resp. $U_1'$, resp. $U_2'$) denotes the image of $g$ (resp. $U_1$, resp. $U_2$) in $G_n^{(m')}(\A^\infty)$. This isomorphism only depends on $g'$, $U_1'$ and $U_2'$ and not on $g$, $U_1$ and $U_2$. This gives the system of sheaves $\{  (R^j\pi_* \Omega^i_{A^{(m)}_n/A^{(m')}_n})_{U'} \}$ a left action of $G_n^{(m')}(\A^\infty)$. Also if $\gamma \in Q_{m,(m-m')}(F)$ then $\gamma: A^{(m)}_{n,U} \ra A^{(m)}_{n,\gamma U}$ gives a natural isomorphism
\[ \gamma: (R^j\pi_* \Omega^i_{A^{(m)}_n/A^{(m')}_n})_{U'} \liso (R^j\pi_* \Omega^i_{A^{(m)}_n/A^{(m')}_n})_{U'}, \]
which depends only on $U'$ and not on $U$.
 This gives the system of sheaves $\{  (R^j\pi_* \Omega^i_{A^{(m)}_n/A^{(m')}_n})_{U'} \}$ a right action of $Q_{m,(m-m')}(F)$. We have $\gamma \circ g = \gamma (g) \circ \gamma$. 

If $U_1' \supset U_2'$ and $g' \in U_2'$ normalizes $U_1'$ then on 
\[ (R^j\pi_* \Omega^i_{A^{(m)}_n/A^{(m')}_n})_{U_2'} \cong (R^j\pi_* \Omega^i_{A^{(m)}_n/A^{(m')}_n})_{U_1'} \otimes_{\cO_{A^{(m')}_{n,U_1'}}} \cO_{A^{(m')}_{n,U_2'}} \]
the actions of $g$ and $1 \otimes g$ agree. Moreover if $U$ is a neat open compact subgroup of $G^{(m)}_n(\A^\infty)$ with image $U'$ in $G_n^{(m')}(\A^\infty)$ then the natural map 
 \[ \pi_{A^{(m)}_n/A^{(m')}_n}^* (\pi_* \Omega^1_{A^{(m)}_n/A^{(m')}_n})_{U'} \lra \Omega^1_{A^{(m)}_{n,U}/A^{(m')}_{n,U'}} \] 
 is an isomorphism. These isomorphisms are equivariant for the actions of the groups  
 $G_n^{(m)}(\A^\infty)$ and $Q_{m,(m-m')}(F)$.

The natural maps 
 \[ \wedge^i ( \pi_* \Omega^1_{A_n^{(m)}/A^{(m')}_n})_{U'} \otimes  \wedge^j (R^1\pi_* \cO_{A_n^{(m)}})_{U'} \lra (R^j\pi_* \Omega^i_{A_n^{(m)}/A^{(m')}_n})_{U'} \]  
  are $G_n^{(m')}(\A^\infty)$ and $Q_{m,(m-m')}(F)$ equivariant isomorphisms.

Suppose that $U$ is a neat open compact subgroup of $G_n^{(m)}(\A^\infty)$ with image $U'$ in $G_n^{(m')}(\A^\infty)$ and $U''$ in $G_n(\A^\infty)$. 
If $U$ is of the form $U' \ltimes M$, then the quasi-isogeny $i_{A^\univ}: A^{(m)}_{n,U} \ra (A^\univ)^{m-m'}$ over $A^{(m')}_{n,U'}$ gives rise to an isomorphism 
\[ \Hom_F( F^{m-m'},\Omega_{n,U''}) \otimes \cO_{A^{(m)}_{n,U}} \cong \Omega^1_{A^{(m)}_{n,U}/A^{(m')}_{n,U'}} \]
and a canonical embeddding
\[ \Xi_{n,U''} \otimes \cO_{A^{(m')}_{n,U'}}\into \Xi_{n,U''}^{\oplus (m-m')} \otimes \cO_{A^{(m')}_{n,U'}} \into (R^1\pi_* \Omega^1_{A^{(m)}/A^{(m')}})_{U'}, \]
where the first map denotes the diagonal embedding. These maps do not depend on the choice of $A^\univ$. They are $G_n^{(m)}(\A^\infty)$-equivariant. The first map is also $Q_{m,(m-m')}(F)$-equivariant, where 
an element $\gamma\in Q_{m,(m-m')}(F)$ acts on the left hand sides by composition with the inverse of the projection of $\gamma$ to $GL_{m-m'}(F)$. This remains true if we do not assume that $U$ has the form $U' \ltimes M$.

This gives rise to canonical $G_n^{(m')}(\A^\infty)$-equivariant isomorphisms
\[ \Hom_F(F^{m-m'},\Omega_{n,U''} ) \otimes \cO_{A^{(m')}_{n,U'}} \cong (\pi_* \Omega^1_{A^{(m)}/A^{(m')}})_{n,U'}. \]
Moreover the composite maps
 \[ \begin{array}{rl} &\Hom((\pi_* \Omega^1_{A^{(m)}_n/A^{(m')}_n})_{U'}, \Xi_{n,U''} \otimes \cO_{A^{(m')}_{n,U'}}) \\
  \into & \Hom((\pi_* \Omega^1_{A_n^{(m)}/A_n^{(m')}})_{U'},(R^1\pi_* \Omega^1_{A_n^{(m)}/A_n^{(m')}})_{U'}) \\ 
   \stackrel{\sim}{\lla} & \Hom((\pi_* \Omega^1_{A_n^{(m)}/A_n^{(m')}})_{U'},(\pi_* \Omega^1_{A_n^{(m)}/A_n^{(m')}})_{U'} \otimes (R^1\pi_* \cO_{A_n^{(m)}})_{U'}) \\ 
    \stackrel{\tr}{\lra} & (R^1\pi_* \cO_{A_n^{(m)}})_{U'} \end{array} \]
 are $G_n^{(m')}(\A^\infty)$-equivariant isomorphisms. 
 
Next we turn to the mixed characteristic case. If $m \geq m'$ and if $U^p$ is a neat open compact subgroup of $G^{(m)}_n(\A^{p,\infty})$ with image $(U^p)'$ in $G_n^{(m')}(\A^{p,\infty})$, and if  $0 \leq N_1 \leq N_2$ are integers, then the sheaf 
\[ R^j\pi_{\cA^{(m),\ord}_{n}/\cA^{(m'),\ord}_n,*} \Omega^i_{\cA^{(m),\ord}_{n,U^p(N_1,N_2)}/\cA^{(m'),\ord}_{n,(U^p)'(N_1,N_2)}} \]
depends only on $(U^p)'$ and not on $U^p$. We will denote it
\[ (R^j\pi_* \Omega^i_{\cA^{(m),\ord}_n/\cA^{(m'),\ord}_n})_{(U^p)'(N_1,N_2)}. \]
If $g \in G_n^{(m)}(\A^{\infty})^\ord$ and $g^{-1}U^p_1(N_{11},N_{12})g \subset U^p_2(N_{21},N_{22})$, then there is a natural map
\[ g: (g')^*(R^j\pi_* \Omega^i_{\cA^{(m),\ord}_n/\cA^{(m'),\ord}_n})_{(U^p_2)'(N_{21},N_{22})} \ra (R^j\pi_* \Omega^i_{\cA^{(m),\ord}_n/\cA^{(m'),\ord}_n})_{(U^p_1)'(N_{11},N_{12})}, \]
where $(U_i^p)'$ denotes the image of $U^p_i$ in $G_n^{(m')}(\A^{p,\infty})$ and $g'$ denotes the image of $g$ in $G_n^{(m')}(\A^{\infty})^\ord$. If $g \in G_n^{(m)}(\A^{\infty})^{\ord,\times}$ then it is an isomorphism.
Moreover this map only depends on $g'$, $(U_1^p)'(N_{11},N_{12})$ and $(U^p_2)'(N_{21},N_{22})$ and not on $g$, $U_1^p(N_{11},N_{12})$ and $U_2^p(N_{21},N_{22})$. This gives the system of sheaves 
\[ \{  (R^j\pi_* \Omega^i_{\cA^{(m),\ord}_n/\cA^{(m'),\ord}_n})_{(U^p)'(N_1,N_2)} \}\]
 a left action of $G_n^{(m')}(\A^{\infty})^\ord$. 
 
 If $\gamma \in Q_{m,(m-m')}(\cO_{F,(p)})$ then $\gamma: A^{(m)}_{n,U^p(N_1,N_2)} \ra A^{(m)}_{n,\gamma U^p(N_1,N_2)}$ gives a natural isomorphism
\[ \gamma: (R^j\pi_* \Omega^i_{\cA^{(m),\ord}_n/\cA^{(m'),\ord}_n})_{(U^p)'(N_1,N_2)} \liso (R^j\pi_* \Omega^i_{\cA^{(m),\ord}_n/\cA^{(m'),\ord}_n})_{(U^p)'(N_1,N_2)}, \]
which depends only on $(U^p)'(N_1,N_2)$ and not on $U^p(N_1,N_2)$.
 This gives the system of sheaves 
 \[ \{  (R^j\pi_* \Omega^i_{\cA^{(m),\ord}_n/\cA^{(m'),\ord}_n})_{(U^p)'(N_1,N_2)} \}\]
  a right action of $Q_{m,(m-m')}(\cO_{F,(p)})$. We have $\gamma \circ g = \gamma (g) \circ \gamma$.

If $(U^p_1)'(N_{11},N_{12}) \supset (U_2^p)'(N_{21},N_{22})$ and $g \in (U^p_1)' (N_{11},N_{12})$ normalizes the subgroup $(U^p_2)'(N_{21},N_{22})$, then on 
\[ \begin{array}{l} (R^j\pi_* \Omega^i_{\cA^{(m),\ord}_n/\cA^{(m'),\ord}_n})_{(U^p_2)'(N_{21},N_{22})} \cong \\ (R^j\pi_* \Omega^i_{\cA^{(m),\ord}_n/\cA^{(m'),\ord}_n})_{(U^p_1)'(N_{11},N_{12})} \otimes_{\cO_{\cA^{(m'),\ord}_{n,(U^p_1)'(N_{11},N_{12})}}} \cO_{\cA^{(m'),\ord}_{n,(U^p_2)'(N_{21},N_{22})}} \end{array} \]
the actions of $g$ and $1 \otimes g$ agree. Moreover if $U^p$ is a neat open compact subgroup of $G^{(m)}_n(\A^{p,\infty})$ with image $(U^p)'$ in $G_n^{(m')}(\A^{p,\infty})$, and if $0 \leq N_1 \leq N_2$ then the natural map 
 \[ \pi_{\cA^{(m),\ord}_n/\cA^{(m'),\ord}_n}^* (\pi_* \Omega^1_{\cA^{(m),\ord}_n/\cA^{(m'),\ord}_n})_{(U^p)'(N_1,N_2)} \lra \Omega^1_{\cA^{(m),\ord}_{n,U^p(N_1,N_2)}/\cA^{(m'),\ord}_{n,(U^p)'(N_1,N_2)}} \] 
 is an isomorphism. These isomorphisms are $G_n^{(m)}(\A^{\infty})^\ord$ and $Q_{m,(m-m')}(\cO_{F,(p)})$ equivariant.

 The natural maps 
 \[ \begin{array}{r} \wedge^i ( \pi_* \Omega^1_{\cA_n^{(m),\ord}/\cA^{(m'),\ord}_n})_{(U^p)'(N_1,N_2)} \otimes  \wedge^j (R^1\pi_* \cO_{\cA_n^{(m),\ord}})_{(U^p)'(N_1,N_2)} \lra \\ (R^j\pi_* \Omega^i_{\cA_n^{(m),\ord}/\cA^{(m'),\ord}_n})_{(U^p)'(N_1,N_2)} \end{array} \]  
  are $G_n^{(m')}(\A^{\infty})^\ord$ and $Q_{m,(m-m')}(\cO_{F,(p)})$ 
  equivariant  isomorphisms.

Under the identification 
\[ \cX^\ord_{n, (U^p)'(N_1,N_2)} \times \Spec \Q \cong X_{n,(U^p)'(N_1,N_2)}\]
the sheaves $\Omega^\ord_{n, (U^p)'(N_1,N_2)}$ (resp. $\Xi^\ord_{n,(U^p)'(N_1,N_2)}$) are naturally identified with the sheaves $\Omega_{n, (U^p)'(N_1,N_2)}$ (resp. $\Xi_{n,(U^p)'(N_1,N_2)}$).
 Moreover, under the identification
 \[ \cA^{(m'),\ord}_{n,U^p(N_1,N_2)} \times \Spec \Q \cong A^{(m')}_{n,U^p(N_1,N_2)}\]
  the sheaf $(R^j\pi_* \Omega^i_{\cA^{(m),\ord}/\cA_n^{(m'),\ord}})_{(U^p)'(N_1,N_2)}$ is naturally identified with the sheaf $(R^j\pi_* \Omega^i_{A^{(m)}_n/A^{(m')}_n})_{(U^p)'(N_1,N_2)}$. These identifications are equivariant for the actions of $G_n(\A^{\infty})^\ord$ and $Q_{m,(m-m')}(\cO_{F,(p)})$.

\newpage

\section{Generalized Shimura Varieties}

We will introduce certain disjoint unions of mixed Shimura varieties, which are associated to $L_{n,(i),\lin}$ and $L_{n,(i)}$ and $P_{n,(i)}^+/Z(N_{n,(i)})$ and $P_{n,(i)}^+$; to $L^{(m)}_{n,(i),\lin}$ and $L^{(m)}_{n,(i)}$ and $P^{(m),+}_{n,(i)}/Z(N^{(m)}_{n,(i)})$ and $P^{(m),+}_{n,(i)}$; and to $\tP^{(m),+}_{n,(i)}$. The differences with the last section are purely book keeping. We then describe certain torus embeddings for these generalized Shimura varieties and discuss their completion along the boundary. These completions will serve as local models near the boundary of the toroidal compactifications of the $X_{n,U}$ and the $A^{(m)}_{n,\tU}$ to be discussed in the next section.

\subsection{Generalized Shimura varieties.}

If $U$ is a neat open compact subgroup of $L^{(m)}_{n,(i),\lin}(\A^\infty)$ we set 
\[ Y_{n,(i),U}^{(m),+}=\coprod_{L^{(m)}_{n,(i),\lin}(\A^\infty)/U} \Spec \Q. \]
In the case $m=0$ we will write simply $Y_{n,(i),U}^+$. Then $\{ Y^{(m),+}_{n,(i),U} \}$ is a system of schemes (locally of finite type over $\Spec \Q$) with right $L^{(m)}_{n,(i),\lin}(\A^\infty)$-action. Each $Y_{n,(i),U}^{(m),+}$ also has a left action of $L^{(m)}_{n,(i),\lin}(\Q)$, which commutes with the right  $L^{(m)}_{n,(i),\lin}(\A^\infty)$-action. If $\delta \in GL_m(F)$ we get a map
\[ \delta: Y_{n,(i),U}^{(m)} \lra Y_{n,(i),\delta (U)}^{(m)} \]
which sends $(\Spec \Q)_{hU} \ra (\Spec \Q)_{\delta (h) \delta (U)}$ via the identity. This gives a left action of $GL_m(F)$ on the inverse system of the $Y_{n,(i),U}^{(m)}$. If $\delta \in GL_m(F)$ and $\gamma \in L^{(m)}_{n,(i),\lin}(\Q)$ and $g \in L^{(m)}_{n,(i),\lin}(\A^\infty)$ then $\delta \circ \gamma = \delta(\gamma) \circ \delta$ and $\delta \circ g = \delta(g) \circ \delta$. 
 If $U'$ denotes the image of $U$ in $L_{n,(i),\lin}(\A^\infty)$ then there is a natural map 
\[ Y_{n,(i),U}^{(m),+} \onto Y_{n,(i),U'}^+. \]
These maps are equivariant for 
\[ L^{(m)}_{n,(i),\lin}(\Q) \times L^{(m)}_{n,(i),\lin}(\A^\infty) \lra L_{n,(i),\lin}(\Q) \times L_{n,(i),\lin}(\A^\infty). \]
The naive quotient
\[  L^{(m)}_{n,(i),\lin}(\Q) \backslash Y^{(m),+}_{n,(i),U} \]
makes sense. We will denote this space
\[ Y^{(m),\natural}_{n,(i),U} \]
and drop the $(m)$ if $m=0$. The inverse system of these spaces has a right action of $L^{(m)}_{n,(i),\lin}(\A^\infty)$, and an action of $GL_m(F)$.
The induced map
\[  Y^{(m),\natural}_{n,(i),U} \lra Y^{\natural}_{n,(i),U} \]
is an isomorphism, and $GL_m(F)$ acts trivially on this space.
(Use the fact that 
\[ (U \cap (\Hom_F(F^m,F^i) \otimes_\Q \A^\infty)) + \Hom_F(F^m,F^i)=\Hom_F(F^m,F^i) \otimes_\Q \A^\infty. ) \]

Similarly if $U^p$ is a neat open compact subgroup of $L^{(m)}_{n,(i),\lin}(\A^{p,\infty})$ and if $N \in \Z_{\geq 0}$ we set
\[ \cY_{n,(i),U^p(N)}^{(m),\ord,+}=\coprod_{L^{(m)}_{n,(i),\lin}(\A^{\infty})^{\ord, \times}/U^p(N)} \Spec \Z_{(p)}. \]
In the case $m=0$ we drop it from the notation.
Each $\cY_{n,(i),U^p\times U_p(N)_{(i)}^{(m)}}^{(m),\ord,+}$ has a left action of $L^{(m)}_{n,(i),\lin}(\Z_{(p)})$ and the inverse system of the $\cY_{n,(i),U^p(N)}^{(m),\ord,+}$ has a commuting right action of $L^{(m)}_{n,(i),\lin}(\A^{\infty})^\ord$. It also has a left action of $GL_m(\cO_{F,(p)})$. If $\delta \in GL_m(\cO_{F,(p)})$ and $\gamma \in L^{(m)}_{n,(i),\lin}(\Z_{(p)})$ and $g \in L^{(m)}_{n,(i),\lin}(\A^\infty)^\ord$ then $\delta \circ \gamma = \delta(\gamma) \circ \delta$ and $\delta \circ g = \delta(g) \circ \delta$. There are equivariant maps 
\[ \cY_{n,(i),U^p(N)}^{(m),\ord,+} \lra \cY_{n,(i),U^p(N)}^{\ord,+}. \]
We set
\[ \cY_{n,(i),U^p(N)}^{\ord,\natural}= L_{n,(i),\lin}(\Z_{(p)}) \backslash \cY_{n,(i),U^p(N)}^{\ord,+}= L^{(m)}_{n,(i),\lin}(\Z_{(p)}) \backslash \cY_{n,(i),U^p(N)}^{(m),\ord,+}. \]
There are maps 
\[ \cY_{n,(i),U^p(N)}^{(m),\ord,+} \times \Spec \Q \into Y_{n,(i),U^p(N)}^{(m),+} \]
which are equivariant for the actions of $L^{(m)}_{n,(i),\lin}(\Z_{(p)})$ and $L^{(m)}_{n,(i),\lin}(\A^{\infty})^\ord $ and $GL_m(\cO_{F,(p)})$. Moreover the  maps $\cY_{n,(i),U^p(N)}^{(m),\ord,+} \ra \cY_{n,(i),U^p(N)}^{\ord,+}$ and $Y_{n,(i),U^p(N)}^{(m),+}\ra Y_{n,(i),U^p(N)}^{+}$ are compatible. The induced maps
\[ \cY_{n,(i),U^p(N)}^{(m),\ord,\natural} \times \Spec \Q \liso Y_{n,(i),U^p(N)}^{(m),\natural} \]
are isomorphisms.

Suppose now that $U$ is a neat open compact subgroup of $L^{(m)}_{n,(i)}(\A^\infty)$ we set 
\[ X_{n,(i),U}^{(m),+}= \left. \left(X_{n-i,U \cap G_{n-i}(\A^\infty)} \times Y^{(m),+}_{n,(i),U \cap L^{(m)}_{n,(i),\lin}(\A^\infty)}\right)\right/U  \]
In the case $m=0$ we will write simply $X^+_{n,(i),U}$. Then $\{ X^{(m),+}_{n,(i),U} \}$ is a system of schemes (locally of finite type over $\Spec \Q$) with right $L^{(m)}_{n,(i)}(\A^\infty)$-action via finite etale maps. Each $X_{n,(i),U}^{(m),+}$ has a left action of $L^{(m)}_{n,(i),\lin}(\Q)$, which commutes with the right  $L^{(m)}_{n,(i)}(\A^\infty)$-action. The system also has a left action of $GL_m(F)$. If $\delta \in GL_m(F)$ and $\gamma \in L^{(m)}_{n,(i),\lin}(\Q)$ and $g \in L^{(m)}_{n,(i)}(\A^\infty)$ then $\delta \circ \gamma = \delta(\gamma) \circ \delta$ and $\delta \circ g = \delta(g) \circ \delta$. If $U'$ is an open normal subgroup of $U$ then $X_{n,(i),U}^{(m),+}$ is identified with $X_{n,(i),U'}^{(m),+}/U$. Projection to the second factor gives $L^{(m)}_{n,(i),\lin}(\Q) \times L^{(m)}_{n,(i)}(\A^\infty)$ and $GL_m(F)$ equivariant maps
\[ X_{n,(i),U}^{(m),+} \lra Y_{n,(i),U}^{(m),+}. \]
The fibre over $g \in L_{n,(i),\lin}(\A^\infty)$ is simply $X_{n-i,U \cap G_{n-i}(\A^\infty)}$.
If $U'$ denotes the image of $U$ in $L_{n,(i)}(\A^\infty)$ then there is a natural, $L^{(m)}_{n,(i),\lin}(\Q) \times L^{(m)}_{n,(i)}(\A^\infty)$-equivariant, commutative diagram 
\[ \begin{array}{ccc}  X_{n,(i),U}^{(m),+} &\onto & X_{n,(i),U'}^+ \\ \da && \da \\ Y_{n,(i),U}^{(m),+} & \onto & Y_{n,(i),U'}^+. \end{array} \]

We have
\[ X^{(m),+}_{n,(i),U}(\C)= L_{n,(i),\herm}(\Q) \backslash (L_{n,(i)}(\A^\infty)/U \times \gH_{n-i}^\pm) \]
and
\[  \pi_0(X^{(m),+}_{n,(i),U} \times \Spec \barQQ)  \cong  \left( L^{(m)}_{n,(i),\lin}(\A^\infty) \times (C_{n-i}(\Q)\backslash C_{n-i}(\A)/C_{n-i}(\R)^0 ) \right)/U. \]

The naive quotient
\[  X^{(m),\natural}_{n,(i),U}=L^{(m)}_{n,(i),\lin}(\Q) \backslash X^{(m),+}_{n,(i),U} \]
makes sense and fibres over $Y^{(m),\natural}_{n,(i),U}$, the fibre over $g$ being $X_{n-i,U_1}$, where $U_1$ denotes the projection to $G_{n-i}(\A^\infty)$ of the subgroup $U_2 \subset U$ consisting of elements whose projection to $L^{(m)}_{n,(i),\lin}(\A^\infty)$ lies in 
$g^{-1}L^{(m)}_{n,(i),\lin}(\Q)g$. If $U'$ denotes the projection of $U$ to $L_{n,(i),\lin}(\A^\infty)$, then the induced map
\[  X^{(m),\natural}_{n,(i),U} \liso  X_{n,(i),U'}^{\natural} \]
is an isomorphism. 
The action of $L_{n,(i)}^{(m)}(\A^\infty)$ is by finite etale maps and if $U'$ is an open normal subgroup of $U$ then $X_{n,(i),U}^{(m),\natural}$ is identified with $X_{n,(i),U'}^{(m),\natural}/U$.
We have
\[  \pi_0(X^{(m),\natural}_{n,(i),U} \times \Spec \barQQ)  \cong  (F^\times \times C_{n-i}(\Q)) \backslash (\A_F^\times \times C_{n-i}(\A)) / U(F_\infty^\times \times C_{n-i}(\R)^0). \]

We define sheaves $\Omega^{+}_{n,(i),U}$ and $\Xi^{+}_{n,(i),U}$ over $X^+_{n,(i),U}$ as the quotients of 
\[ \Omega_{n-i,U \cap G_{n-i}(\A^\infty)} / X_{n-i,U \cap G_{n-i}(\A^\infty)} \times Y^{+}_{n,(i),U \cap L_{n,(i),\lin}(\A^\infty)} \]
and
\[ \Xi_{n-i,U \cap G_{n-i}(\A^\infty)} / X_{n-i,U \cap G_{n-i}(\A^\infty)} \times Y^{+}_{n,(i),U \cap L_{n,(i),\lin}(\A^\infty)} \]
by $U$. Then $\{ \Omega^{+}_{n,(i),U} \}$ and $\{ \Xi^{+}_{n,(i),U} \}$ are systems of locally free sheaves on $X^{+}_{n,(i),U}$ with left $L_{n,(i)}(\A^\infty)$-action and commuting right $L_{n,(i),\lin}(\Q)$-action. 

 Let $\cE_{(i),U}^{+}$ denote the principal
$R_{n,(n),(i)}/N(R_{n,(n),(i)})$-bundle on $X^+_{n,(i),U}$ in the
Zariski topology defined by setting, for $W \subset X^+_{n,(i),U}$ a Zariski open, 
$\cE_{(i),U}^{+}(W)$  to
be the set of pairs $(\xi_0,\xi_1)$, where
\[ \xi_0: \Xi^+_{n,(i),U}|_W \liso \cO_W \]
and
\[ \xi_1: \Omega^+_{n,(i),U} \liso \Hom_\Q( V_n/V_{n,(n)} , \cO_W).  \]
We define the $R_{n,(n),(i)}/N(R_{n,(n),(i)})$-action on $\cE_{(i),U}^{+}$ by
\[ h(\xi_0,\xi_1)=(\nu(h)^{-1}\xi_0, (\circ h^{-1}) \circ \xi_1). \] 
The inverse system $\{ \cE_{(i),U}^{+} \}$ has an action of $L_{n,(i)}(\A^\infty)$ and of $L_{n,(i),\lin}(\Q)$. 

Suppose that $R_0$ is an irreducible noetherian $\Q$-algebra and that $\rho$ is a
representation of $R_{n,(n),(i)}$ on a finite, locally free $R_0$-module $W_\rho$. We define a locally
free sheaf $\cE^{+}_{(i),U,\rho}$ over
$X^{+}_{n,(i),U} \times \Spec R_0$ by setting $\cE^{+}_{(i),U,\rho}(W)$ to be the set of $(R^{(m)}_{n,(n)}/N(R_{n,(n),(i)}))(\cO_W)$-equivariant
maps of Zariski sheaves of sets
\[ \cE^{+}_{(i),U}|_W \lra W_\rho \otimes_{R_0} \cO_W. \]
Then $\{ \cE^{+}_{(i),U,\rho} \}$ is a system of locally free sheaves with $L_{n,(i)}(\A^\infty)$-action and $L_{n,(i),\lin}(\Q)$-action over the system of schemes 
$\{ X^{+}_{n,(i),U} \times \Spec R_0\}$. The restriction of $\cE^+_{(i),U,\rho}$ to $X^{(i)}_{n-i,hUh^{-1} \cap G_{n-i}(\A^\infty)}$ can be identified with $\cE_{hUh^{-1} \cap G_{n-i}(\A^\infty),\rho|_{L_{n-i,(n-i)}}}$. However the description of the actions of $L_{n,(i)}(\A^\infty)$ and $L_{n,(i),\lin}(\Q)$ involve $\rho$ and not just $\rho|_{L_{n-i,(n-i)}}$. 
If $g \in L_{n,(i)}(\A^\infty)$ and $\gamma \in L_{n,(i),\lin}(\Q)$, then the natural maps
\[ g^*\cE_{(i),U,\rho}^+ \lra \cE_{(i),U',\rho}^+ \]
and
\[ \gamma^*\cE_{(i),U,\rho}^+ \lra \cE_{(i),U',\rho}^+ \]
are isomorphisms.

We will also write
\[ \Omega^{\natural}_{n,(i),U} = L_{n,(i),\lin}(\Q) \backslash \Omega^{+}_{n,(i),U} \]
and
\[ \Xi^{\natural}_{n,(i),U} = L_{n,(i),\lin}(\Q) \backslash \Xi^{+}_{n,(i),U}, \]
locally free sheaves on $X^{\natural}_{n,(i),U}$. (If $\rho$ is trivial on $L_{n,(i),\lin}$ then one can also form, the quotient of $\cE^+_{(i),U,\rho}$ by $L_{n,(i),\lin}(\Q)$, but in general this quotient does not make sense.)

If $U^p$ is a neat open compact subgroup of $L^{(m)}_{n,(i)}(\A^{p,\infty})$ and $N_2 \geq N_1 \geq 0$ we 
set
\[ \cX_{n,(i),U^p(N_1,N_2)}^{(m),\ord,+}= \left. \left(\cX^\ord_{n-i,(U^p \cap G_{n-i}(\A^{p,\infty}))(N_1,N_2)} \times \cY^{(m),\ord,+}_{n,(i),(U^p \cap L^{(m)}_{n,(i),\lin}(\A^{p,\infty}))(N_1)}\right)\right/U^p.  \]
In the case $m=0$ we drop it from the notation. Each $\cX_{n,(i),U^p(N_1,N_2)}^{(m),\ord,+}$ has a left action of $L_{n,(i),\lin}^{(m)}(\Z_{(p)})$ and the inverse system has a commuting right action of $L_{n,(i)}^{(m)}(\A^{\infty})^\ord$. There is also a left action of $GL_m(\cO_{F,(p)})$. If $\delta \in GL_m(\cO_{F,(p)})$ and $\gamma \in L^{(m)}_{n,(i),\lin}(\Z_{(p)})$ and $g \in L^{(m)}_{n,(i),\lin}(\A^{\infty})^\ord$ then $\delta \circ \gamma = \delta(\gamma) \circ \delta$ and $\delta \circ g = \delta(g) \circ \delta$. 
If $g \in L_{n,(i)}^{(m)}(\A^{\infty})^\ord$ and if
\[ g: \cX_{n,(i),U^p(N_1,N_2)}^{(m),\ord,+} \lra \cX_{n,(i),(U^p)'(N_1',N_2')}^{(m),\ord,+}, \]
then this map is quasi-finite and flat. If $g \in L_{n,(i)}^{(m)}(\A^{\infty})^{\ord,\times}$ then it is etale, and, if further $N_2=N_2'$, then it is finite etale. If $N_2'>0$ and $p^{N_2-N_2'}\nu(g_p) \in \Z_p^\times$ then the map is finite. On $\F_p$-fibres the map $\varsigma_p$ is absolute Frobenius composed with the forgetful map.
If $(U^p)'$ is an open normal subgroup of $U^p$ and if $N_1 \leq N_1'\leq N_2$ then
\[ \cX_{n,(i),(U^p)'(N_1',N_2)}^{(m),\ord,+}/U^p(N_1,N_2) \liso \cX_{n,(i),U^p(N_1,N_2)}^{(m),\ord,+}. \]
There are commutative diagrams
\[ \begin{array}{ccc}  \cX_{n,(i),U^p(N_1,N_2)}^{(m),\ord,+} &\onto & \cX_{n,(i),U^(N_1,N_2)}^{\ord,+} \\ \da && \da \\ \cY_{n,(i),U^p(N_1,N_2)}^{(m),\ord,+} & \onto & \cY_{n,(i),U^p(N_1,N_2)}^{\ord,+}. \end{array} \]

We set 
\[ \cX_{n,(i),U^p(N_1,N_2)}^{\ord,\natural}=L_{n,(i),\lin}(\Z_{(p)}) \backslash \cX_{n,(i),U^p(N_1,N_2)}^{\ord,+}=L_{n,(i),\lin}^{(m)}(\Z_{(p)}) \backslash \cX_{n,(i),U^p(N_1,N_2)}^{(m),\ord,+}. \]
The system of these spaces has a right action of $L_{n,(i)}^{(m)}(\A^{\infty})^\ord$ and a left action of $GL_m(\cO_{F,(p)})$. If $\delta \in GL_m(\cO_{F,(p)})$ and $g \in L^{(m)}_{n,(i),\lin}(\A^{\infty})^\ord$ then $\delta \circ g = \delta(g) \circ \delta$. 
If $g \in L_{n,(i)}^{(m)}(\A^{\infty})^\ord$ and if
\[ g: \cX_{n,(i),U^p(N_1,N_2)}^{(m),\ord,\natural} \lra \cX_{n,(i),(U^p)'(N_1',N_2')}^{(m),\ord,\natural}, \]
then this map is quasi-finite and flat. If $g \in L_{n,(i)}^{(m)}(\A^{\infty})^{\ord,\times}$ then it is etale, and, if further $N_2=N_2'$, then it is finite etale. If $N_2'>0$ and $p^{N_2-N_2'}\nu(g_p) \in \Z_p^\times$ then the map is finite. On $\F_p$-fibres the map $\varsigma_p$ is absolute Frobenius composed with the forgetful map.
If $(U^p)'$ is an open normal subgroup of $U^p$ and if $N_1 \leq N_1'\leq N_2$ then
\[ \cX_{n,(i),(U^p)'(N_1',N_2)}^{(m),\ord,\natural}/U^p(N_1,N_2) \liso \cX_{n,(i),U^p(N_1,N_2)}^{(m),\ord,\natural}. \]

We define sheaves $\Omega^{\ord,+}_{n,(i),U^p(N_1,N_2)}$ and $\Xi^{\ord,+}_{n,(i),U^p(N_1,N_2)}$ over $\cX^{\ord,+}_{n,(i),U^p(N_1,N_2)}$ as the quotients of 
\[ \Omega_{n-i,(U^p \cap G_{n-i}(\A^{p,\infty}))(N_1,N_2)}^{\ord} / \cX_{n-i,(U^p \cap G_{n-i}(\A^{p,\infty}))(N_1,N_2)} \times \cY^{\ord,+}_{n,(i),(U^p \cap L_{n,(i),\lin}(\A^{p,\infty}))(N_1)} \]
and
\[ \Xi^{\ord}_{n-i,(U^p \cap G_{n-i}(\A^{p,\infty}))(N_1,N_2)} / \cX_{n-i,(U^p \cap G_{n-i}(\A^{p,\infty}))(N_1,N_2)} \times \cY^{\ord,+}_{n,(i),(U^p \cap L_{n,(i),\lin}(\A^{p,\infty}))(N_1)} \]
by $U^p$. Then the systems of sheaves $\Omega^{\ord,+}_{n,(i),U^p(N_1,N_2)}$ and $\Xi^{\ord,+}_{n,(i),U^p(N_1,N_2)}$ have commuting actions of  $L_{n,(i),\lin}(\Z_{(p)})$ and $L_{n,(i)}(\A^{\infty})^\ord$. 

 Let $\cE_{(i),U^p(N_1,N_2)}^{\ord,+}$ denote the principal
$R_{n,(n),(i)}/\!N(R_{n,(n),(i)})$-bundle for the Zariski topology on $\cX^{\ord,+}_{n,(i),U^p(N_1,N_2)}$  defined by setting, for $W \subset \cX^{\ord,+}_{n,(i),U^p(N_1,N_2)}$ a Zariski open, 
$\cE_{(i),U^p(N_1,N_2)}^{\ord,+}(W)$  to
be the set of pairs $(\xi_0,\xi_1)$, where
\[ \xi_0: \Xi^{\ord,+}_{n,(i),U^p(N_1,N_2)}|_W \liso \cO_W \]
and
\[ \xi_1: \Omega^{\ord,+}_{n,(i),U^p(N_1,N_2)} \liso \Hom( \Lambda_n/\Lambda_{n,(n)} , \cO_W).  \]
We define the $R_{n,(n),(i)}/N(R_{n,(n),(i)})$-action on $\cE_{(i),U^p(N_1,N_2)}^{\ord,+}$ by
\[ h(\xi_0,\xi_1)=(\nu(h)^{-1}\xi_0, (\circ h^{-1}) \circ \xi_1). \] 
The inverse system $\{ \cE_{(i),U^p(N_1,N_2)}^{+} \}$ has an action of $L_{n,(i)}(\A^\infty)^\ord$ and an action of $L_{n,(i),\lin}(\Z_{(p)})$. 

Suppose that $R_0$ is an irreducible noetherian $\Z_{(p)}$-algebra and that $\rho$ is a
representation of $R_{n,(n),(i)}$ on a finite, locally free $R_0$-module $W_\rho$. We define a locally
free sheaf $\cE^{\ord,+}_{(i),U^p(N_1,N_2),\rho}$ over
$\cX^{\ord,+}_{n,(i),U^p(N_1,N_2)} \times \Spec R_0$ by setting $\cE^{\ord,+}_{(i),U^p(N_1,N_2),\rho}(W)$ to be the set of $(R^{(m)}_{n,(n)}/N(R_{n,(n),(i)}))(\cO_W)$-equivariant
maps of Zariski sheaves of sets
\[ \cE^{\ord,+}_{(i),U^p(N_1,N_2)}|_W \lra W_\rho \otimes_{R_0} \cO_W. \]
Then $\{ \cE^{\ord,+}_{(i),U^p(N_1,N_2),\rho} \}$ is a system of locally free sheaves with $L_{n,(i)}(\A^\infty)^\ord$-action and $L_{n,(i),\lin}(\Z_{(p)})$-action over the system of schemes 
$\{ \cX^{\ord,+}_{n,(i),U^p(N_1,N_2)} \times \Spec R_0\}$. The restriction of $\cE^{\ord,+}_{(i),U^p(N_1,N_2),\rho}$ to $\cX^{(i),\ord}_{n-i,(hU^ph^{-1} \cap G_{n-i}(\A^{p,\infty})(N_1,N_2))}$ can be identified with $\cE^{\ord}_{(hU^ph^{-1} \cap G_{n-i}(\A^{p,\infty}))(N_1,N_2),\rho|_{L_{n-i,(n-i)}}}$. However the description of the actions of $L_{n,(i)}(\A^\infty)^\ord$ and $L_{n,(i),\lin}(\Z_{(p)})$ involve $\rho$ and not just $\rho|_{L_{n-i,(n-i)}}$. 
If $g \in L_{n,(i)}(\A^\infty)^{\ord,\times}$ and $\gamma \in L_{n,(i),\lin}(\Z_{(p)})$, then the natural maps
\[ g^*\cE_{(i),U^p(N_1,N_2),\rho}^{\ord,+}l \lra \cE_{(i),(U^p)'(N_1',N_2'),\rho}^{\ord,+} \]
and
\[ \gamma^*\cE_{(i),U^p(N_1,N_2),\rho}^{\ord,+}l \lra \cE_{(i),(U^p)'(N_1',N_2'),\rho}^{\ord,+} \]
are isomorphisms.

We will also write
\[ \Omega^{\ord,\natural}_{n,(i),U^p(N_1,N_2)} = L_{n,(i),\lin}(\Z_{(p)}) \backslash \Omega^{\ord,+}_{n,(i),U^p(N_1,N_2)} \]
and
\[ \Xi^{\ord,\natural}_{n,(i),U^p(N_1,N_2)} = L_{n,(i),\lin}(\Z_{(p)}) \backslash \Xi^{\ord,+}_{n,(i),U^p(N_1,N_2)}, \]
locally free sheaves on $\cX^{\ord,\natural}_{n,(i),U^p(N_1,N_2)}$. 

There are maps
\[ \cX^{(m),\ord,+}_{n,(i),U^p(N_1,N_2)} \times \Spec \Q \into X^{(m),+}_{n,(i),U^p(N_1,N_2)} \]
which are equivariant for the actions of the groups $L_{n,(i)}^{(m)}(\A^\infty)^\ord$ and $L_{n,(i),\lin}^{(m)}(\Z_{(p)})$ and $GL_m(\cO_{F,(p)})$. Under these maps $\Omega^{\ord,+}_{n,(i),U^p(N_1,N_2)}$ (resp. $\Xi^{\ord,+}_{n,(i),U^p(N_1,N_2)}$, resp. $\cE^{\ord,+}_{(i),U,\rho}$) corresponds to $\Omega^+_{n,(i),U^p(N_1,N_2)}$ (resp. $\Xi^+_{n,(i),U^p(N_1,N_2)}$, resp. $\cE^{\ord,+}_{(i),U,\rho \otimes \Q}$). The induced maps
\[ \cX^{(m),\ord,\natural}_{n,(i),U^p(N_1,N_2)} \times \Spec \Q \into X^{(m),\natural}_{n,(i),U^p(N_1,N_2)} \]
are isomorphisms.

\newpage \subsection{Generalized Kuga-Sato varieties.}\label{gksv}

Now suppose that $U$ is a neat open compact subgroup of $(P^{(m),+}_{n,(i)}/Z(N^{(m)}_{n,(i)}))(\A^\infty)=(\tP^{(m),+}_{n,(i)}/Z(\tN^{(m)}_{n,(i)}))(\A^\infty)$. We set 
\[ A_{n,(i),U}^{(m),+}=  \coprod_{h \in L^{(m)}_{n,(i),\lin}(\A^\infty)/U}  A^{(i+m)}_{n-i, h U h^{-1} \cap G_{n-i}^{(i+m)}(\A^\infty) }.  \]
In the case $m=0$ we will write simply $A^+_{n,(i),U}$. 

If $g \in (P^{(m),+}_{n,(i)}/Z(N^{(m)}_{n,(i)}))(\A^\infty)$ and $g^{-1} U g \subset U'$, then we define a finite etale map
\[ g: A_{n,(i),U}^{(m),+} \lra A_{n,(i),U'}^{(m),+} \]
to be the coproduct of the maps
\[ g': A^{(i+m)}_{n-i, h U h^{-1} \cap G_{n-i}^{(i+m)}(\A^\infty) } \lra A^{(i+m)}_{n-i, h' U' (h')^{-1} \cap G_{n-i}^{(i+m)}(\A^\infty) }, \]
where $h,h' \in L^{(m)}_{n,(i),\lin}(\A^\infty)$ and $g' \in G^{(i+m)}_{n-i}(\A^\infty)$ satisfy $hg=g' h'$. This makes $\{ A^{(m),+}_{n,(i),U} \}$ a system of schemes (locally of finite type over $\Spec \Q$) with right action of the group $(P^{(m),+}_{n,(i)}/Z(N^{(m)}_{n,(i)}))(\A^\infty)$. If $U'$ is an open normal subgroup of $U$ then $A_{n,(i),U}^{(m),+}$ is identified with $A_{n,(i),U'}^{(m),+}/U$.

If $\gamma \in L^{(m)}_{n,(i),\lin}(\Q)$, then we define
\[ \gamma: A_{n,(i),U}^{(m),+} \lra A_{n,(i),U}^{(m),+} \]
to be the coproduct of the maps
\[ \gamma: A^{(i+m)}_{n-i, h U h^{-1} \cap G_{n-i}^{(i+m)}(\A^\infty) } \lra A^{(i+m)}_{n-i, (\gamma h)U (\gamma h)^{-1} \cap G_{n-i}^{(i+m)}(\A^\infty) }. \]
This gives a left action of $L^{(m)}_{n,(i),\lin}(\Q)$ on each $A_{n,(i),U}^{(m),+}$, which commutes with the action of $(P^{(m),+}_{n,(i)}/Z(N^{(m)}_{n,(i)}))(\A^\infty)$. 

If $\delta \in GL_m(F)$ define a map 
\[ \delta : A_{n,(i),U}^{(m),+} \lra A_{n,(i),\delta (U)}^{(m),+} \]
as the coproduct of the maps
\[ \delta: A^{(i+m)}_{n-i, h U h^{-1} \cap G_{n-i}^{(i+m)}(\A^\infty) } \lra A^{(i+m)}_{n-i, \delta (hUh^{-1}) \cap G_{n-i}^{(i+m)}(\A^\infty) }. \]
This gives a left $GL_m(F)$-action on the system of the $A_{n,(i),U}^{(m),+}$. 
If $\delta \in GL_m(F)$ and $\gamma \in L^{(m)}_{n,(i),\lin}(\Q)$ and $g \in (P^{(m),+}_{n,(i)}/Z(N^{(m)}_{n,(i)}))(\A^\infty)$ then $\delta \circ \gamma = \delta(\gamma) \circ \delta$ and $\delta \circ g = \delta(g) \circ \delta$.

There are natural maps
\[ A_{n,(i),U}^{(m),+} \lra X_{n,(i),U}^{(m),+}, \]
which are equivariant for the actions of 
$(P^{(m),+}_{n,(i)}/Z(N^{(m)}_{n,(i)}))(\A^\infty)$ and $L^{(m)}_{n,(i),\lin}(\Q)$ and $GL_m(F)$.

If $U'$ denotes the image of $U$ in $(P_{n,(i)}^+/Z(N_{n,(i)}))(\A^\infty)$ then there is a natural, $L^{(m)}_{n,(i),\lin}(\Q)$-equivariant and $(P^{(m),+}_{n,(i)}/Z(N^{(m)}_{n,(i)}))(\A^\infty)$-equivariant, commutative diagram:
\[ \begin{array}{ccc}  A_{n,(i),U}^{(m),+} &\onto & A_{n,(i),U'}^+ \\ \da && \da \\ X_{n,(i),U}^{(m),+} & \onto & X_{n,(i),U'}^+ \\ \da && \da \\ Y_{n,(i),U}^{(m),+} & \onto & Y_{n,(i),U'}^+. \end{array} \]

We have
\[ A_{n,(i),U}^{(m),+}(\C) = (P_{n,(i)}^{(m)}/Z(N_{n,(i)}^{(m)}))(\Q) \backslash (P_{n,(i)}^{(m),+}/Z(N_{n,(i)}^{(m)}))(\A)/(UU_{n-i,\infty}^0A_{n-i}(\R)^0). \] 

Note that it does not make sense to divide $A_{n,(i),U}^{(m),+}$ by $L_{n,(i),\lin}^{(m)}(\Q)$, so we don't do so.

We define a semi-abelian scheme $\tG^\univ/A^+_{n,(i),U}$ by requiring that over the open and closed subscheme $A^{(i)}_{n-i,hUh^{-1} \cap G_{n-i}^{(i)}(\A^\infty)}$ it restricts to $G^\univ$. It is unique up to unique quasi-isomorphism. 
We also define a sheaf $\tOmega^{+}_{n,(i),U}$ (resp. $\tXi^+_{n,(i),U}$) over $A^+_{n,(i),U}$ to be the unique sheaf which, for each $h$, restricts to $\Omega^{(i)}_{n-i,hUh^{-1} \cap G_{n-i}^{(i)}(\A^\infty)}$ (resp. $\Xi_{n-i,hUh^{-1} \cap G_{n-i}^{(i)}(\A^\infty)}$) on $A^{(i)}_{n-i,hUh^{-1} \cap G_{n-i}^{(i)}(\A^\infty)}$. Thus $\tOmega^+_{n,(i),U}$ is the pull back by the identity section of $\Omega^1_{\tG^\univ/A^+_{n,(i),U}}$. Then $\{ \tOmega^{+}_{n,(i),U} \}$ (resp. $\{ \tXi^{+}_{n,(i),U} \}$)  is a  system of locally free sheaves on $A^{+}_{n,(i),U}$ with a left $(P_{n,(i)}^+/Z(N_{n,(i)}))(\A^\infty)$-action and a commuting right $L_{n,(i),\lin}(\Q)$-action. There are equivariant exact sequences
\[ (0) \lra \pi^* \Omega^+_{n,(i),U} \lra \tOmega^+_{n,(i),U} \lra F^i \otimes_\Q \cO_{A^+_{n,(i),U}} \lra (0), \]
where $\pi$ denotes the map $A^+_{n,(i),U} \ra X^+_{n,(i),U}$.

 Let $\tcE_{(i),U}^{+}$ denote the principal
$R_{n,(n),(i)}$-bundle on $A^+_{n,(i),U}$ in the
Zariski topology defined by setting, for $W \subset A^+_{n,(i),U}$ a Zariski open, 
$\tcE_{(i),U}^{+}(W)$  to
be the set of pairs $(\xi_0,\xi_1)$, where
\[ \xi_0: \Xi^+_{n,(i),U}|_W \liso \cO_W \]
and
\[ \xi_1: \tOmega^+_{n,(i),U} \liso \Hom_\Q( V_n/V_{n,(n)}\oplus \Hom_\Q(F^m,\Q) , \cO_W)  \]
satisfies
\[ \xi_1: \Omega^+_{n,(i),U} \liso \Hom_\Q( V_n/V_{n,(n)}, \cO_W). \]
We define the $R_{n,(n),(i)}$-action on $\cE_{(i),U}^{+}$ by
\[ h(\xi_0,\xi_1)=(\nu(h)^{-1}\xi_0, (\circ h^{-1}) \circ \xi_1). \] 
The inverse system $\{ \tcE_{(i),U}^{+} \}$ has an action of $P_{n,(i)}^+(\A^\infty)$ and of $L_{n,(i),\lin}(\Q)$. 

Suppose that $R_0$ is an irreducible noetherian $\Q$-algebra and that $\rho$ is a
representation of $R_{n,(n),(i)}$ on a finite, locally free $R_0$-module $W_\rho$. We define a locally
free sheaf $\cE^{+}_{(i),U,\rho}$ over
$A^{+}_{n,(i),U} \times \Spec R_0$ by setting $\cE^{+}_{(i),U,\rho}(W)$ to be the set of $R^{(m)}_{n,(n)}(\cO_W)$-equivariant
maps of Zariski sheaves of sets
\[ \tcE^{+}_{(i),U}|_W \lra W_\rho \otimes_{R_0} \cO_W. \]
Then $\{ \cE^{+}_{(i),U,\rho} \}$ is a system of locally free sheaves with $P_{n,(i)}^+(\A^\infty)$-action and $L_{n,(i),\lin}(\Q)$-action over the system of schemes 
$\{ A^{+}_{n,(i),U} \times \Spec R_0\}$. The restriction of $\cE^+_{(i),U,\rho}$ to $A^{(i)}_{n-i,hUh^{-1} \cap G_{n-i}^{(i)}(\A^\infty)}$ can be identified with $\cE^{(i)}_{hUh^{-1} \cap G_{n-i}^{(i)}(\A^\infty),\rho|_{R_{n-i,(n-i)}^{(i)}}}$. However the description of the actions of $P_{n,(i)}^+(\A^\infty)$ and $L_{n,(i),\lin}(\Q)$ involve $\rho$ and not just $\rho|_{R^{(i)}_{n-i,(n-i)}}$. 
If $g \in P_{n,(i)}^+(\A^\infty)$ and $\gamma \in L_{n,(i),\lin}(\Q)$, then the natural maps
\[ g^*\cE_{(i),U,\rho}^+ \lra \cE_{(i),U',\rho}^+ \]
and
\[ \gamma^*\cE_{(i),U,\rho}^+ \lra \cE_{(i),U',\rho}^+ \]
are isomorphisms.
If $\rho$ factors through $R_{n,(n),(i)}/N(R_{n,(n),(i)})$ then $\cE^{+}_{(i),U,\rho}$ is canonically isomorphic to the pull-back of $\cE^+_{(i),U,\rho}$ from $X^+_{n,(i),U}$. In general $W_\rho$ has a filtration by $R_{n,(n),(i)}$-invariant local direct-summands  such that the action of $R_{n,(n),(i)}$ on each graded piece factors through $R_{n,(n),(i)}/N(R_{n,(n),(i)})$. Thus $\cE^{+}_{(i),U,\rho}$ has a $P_{n,(i)}^+(\A^\infty)$ and $L_{n,(i),\lin}(\Q)$ invariant filtration by local direct summands such that each graded piece is the pull back of some $\cE^+_{(i),U,\rho'}$ from $X^+_{n,(i),U}$.
 
Similarly if $U^p$ is a neat open compact subgroup of $(P^{(m),+}_{n,(i)}/Z(N^{(m)}_{n,(i)}))(\A^{p,\infty})=(\tP^{(m),+}_{n,(i)}/Z(\tN^{(m)}_{n,(i)}))(\A^{p,\infty})$ we set
\[ \cA_{n,(i),U^p(N_1,N_2)}^{(m),\ord,+}=  \coprod_{h \in L^{(m)}_{n,(i),\lin}(\A^{\infty})^{\ord,\times}/U^p(N_1,N_2)}  \cA^{(i+m),\ord}_{n-i, (h U^p h^{-1} \cap G_{n-i}^{(i+m)}(\A^{p,\infty}))(N_1,N_2) }.  \]
In the case $m=0$ we will write simply $\cA^{\ord,+}_{n,(i),U^p(N_1,N_2)}$. 
The inverse system of the $\cA_{n,(i),U^p(N_1,N_2)}^{(m),\ord,+}$ has a right action of $(P^{(m),+}_{n,(i)}/Z(N^{(m)}_{n,(i)}))(\A^{\infty})^\ord$ and a commuting left action of $L_{n,(i),\lin}^{(m)}(\Z_{(p)})$. 
If $g \in (P^{(m),+}_{n,(i)}/Z(N^{(m)}_{n,(i)}))(\A^{\infty})^\ord$ then the map
\[ g: \cA_{n,(i),U^p(N_1,N_2)}^{(m),\ord,+} \lra \cA_{n,(i),(U^p)'(N_1',N_2')}^{(m),\ord,+}, \]
is quasi-finite and flat. If $g \in (P^{(m),+}_{n,(i)}/Z(N^{(m)}_{n,(i)}))(\A^{\infty})^{\ord,\times}$ then it is etale, and, if further $N_2=N_2'$, then it is finite etale. If $N_2'>0$ and $p^{N_2-N_2'}\nu(g_p) \in \Z_p^\times$ then the map is finite. On $\F_p$-fibres the map $\varsigma_p$ is absolute Frobenius composed with the forgetful map.
If $(U^p)'$ is an open normal subgroup of $U^p$ and if $N_1 \leq N_1' \leq N_2$ then $\cA_{n,(i),(U^p)'(N_1',N_2)}^{(m),\ord,+}/U^p(N_1,N_2)$ is identified with $\cA_{n,(i),U^p(N_1,N_2)}^{(m),\ord,+}$.
Further there is a left action of $GL_m(\cO_{F,(p)})$ such that if $\delta \in GL_m(\cO_{F,(p)})$ and $\gamma \in L^{(m)}_{n,(i),\lin}(\Z_{(p)})$ and $g \in (P^{(m),+}_{n,(i)}/Z(N^{(m)}_{n,(i)}))(\A^{\infty})^\ord$, $\delta \circ \gamma = \delta(\gamma) \circ \delta$ and $\delta \circ g = \delta(g) \circ \delta$.
There are natural equivariant maps
\[ \cA_{n,(i),U^p(N_1,N_2)}^{(m),\ord,+} \lra \cX_{n,(i),U^p(N_1,N_2)}^{(m),\ord,+}. \]
If $(U^p)'$ denotes the image of $U^p$ in $(P_{n,(i)}^+/Z(N_{n,(i)}))(\A^{p,\infty})$ then there is a natural equivariant, commutative diagram:
\[ \begin{array}{ccc}  \cA_{n,(i),U^p(N_1,N_2)}^{(m),\ord,+} &\onto & \cA_{n,(i),(U^p)'(N_1,N_2)}^{\ord,+} \\ \da && \da \\ \cX_{n,(i),U^p(N_1,N_2)}^{(m),\ord,+} & \onto & \cX_{n,(i),(U^p)'(N_1,N_2)}^{\ord,+} \\ \da && \da \\ \cY_{n,(i),U^p(N_1,N_2)}^{(m),\ord,+} & \onto & \cY_{n,(i),(U^p)'(N_1,N_2)}^{\ord,+}. \end{array} \]
There are equivariant embeddings
\[ \cA_{n,(i),U^p(N_1,N_2)}^{(m),\ord,+} \times \Spec \Q \into A_{n,(i),U^p(N_1,N_2)}^{(m),+}. \]

We define a semi-abelian scheme $\tcG^\univ/\cA^{\ord,+}_{n,(i),U^p(N_1,N_2)}$ over $\cA^{\ord,+}_{n,(i),U^p(N_1,N_2)}$ by requiring that over $\cA^{(i),\ord}_{n-i,(hU^ph^{-1} \cap G_{n-i}^{(i)}(\A^{p,\infty}))(N_1,N_2)}$ it restricts to $\cG^\univ$. It is unique up to unique prime-to-$p$ quasi-isomorphism. 
We define a sheaf $\tOmega^{\ord,+}_{n,(i),U^p(N_1,N_2)}$ (resp. $\tXi^{\ord,+}_{n,(i),U^p(N_1,N_2)}$) over $\cA^{\ord,+}_{n,(i),U^p(N_1,N_2)}$ to be the sheaf which, for each $h$, restricts to $\Omega^{(i),\ord}_{n-i,(hU^ph^{-1} \cap G_{n-i}^{(i)}(\A^{p,\infty}))(N_1,N_2)}$ (resp. $\Xi^{(i),\ord}_{n-i,(hU^ph^{-1} \cap G_{n-i}^{(i)}(\A^{p,\infty}))(N_1,N_2)}$) on $\cA^{(i),\ord}_{n-i,(hU^ph^{-1} \cap G_{n-i}^{(i)}(\A^{p,\infty}))(N_1,N_2)}$. 
Then $\tOmega^{\ord,+}_{n,(i),U^p(N_1,N_2)}$ is the pull back by the identity section of $\Omega^1_{\tcG^\univ/\cA^{\ord,+}_{n,(i),U^p(N_1,N_2)}/\cA^{\ord,+}_{n,(i),U^p(N_1,N_2)}}$. 
Then the collection $\{ \tOmega^{\ord,+}_{n,(i),U^p(N_1,N_2)} \}$ (resp. $\{ \tOmega^{\ord,+}_{n,(i),U^p(N_1,N_2)} \}$) is a  system of locally free sheaves on $\cA^{\ord,+}_{n,(i),U^p(N_1,N_2)}$ with a left $(P_{n,(i)}^+/Z(N_{n,(i)}))(\A^\infty)^\ord$-action and a commuting right $L_{n,(i),\lin}(\Z_{(p)})$-action. Also there are equivariant exact sequences
\[ (0) \lra \pi^* \Omega^{\ord,+}_{n,(i),U^p(N_1,N_2)} \lra \tOmega^{\ord,+}_{n,(i),U^p(N_1,N_2)} \lra F^i \otimes_\Q \cO_{\cA^{\ord,+}_{n,(i),U^p(N_1,N_2)}} \lra (0), \]
where $\pi$ denotes the map $\cA^{\ord,+}_{n,(i),U} \ra \cX^{\ord,+}_{n,(i),U}$.

 Let $\tcE_{(i),U^p(N_1,N_2)}^{\ord,+}$ denote the principal
$R_{n,(n),(i)}$-bundle on $\cA^{\ord,+}_{n,(i),U^p(N_1,N_2)}$ in the
Zariski topology defined by setting, for $W \subset \cA^{\ord,+}_{n,(i),U^p(N_1,N_2)}$ a Zariski open, 
$\tcE_{(i),U^p(N_1,N_2)}^{\ord,+}(W)$  to
be the set of pairs $(\xi_0,\xi_1)$, where
\[ \xi_0: \Xi^{\ord,+}_{n,(i),U^p(N_1,N_2)}|_W \liso \cO_W \]
and
\[ \xi_1: \tOmega^{\ord,+}_{n,(i),U^p(N_1,N_2)} \liso \Hom_\Q( V_n/V_{n,(n)}\oplus \Hom_\Q(F^m,\Q) , \cO_W)  \]
satisfies
\[ \xi_1: \Omega^{\ord,+}_{n,(i),U^p(N_1,N_2)} \liso \Hom_\Q( V_n/V_{n,(n)}, \cO_W). \]
We define the $R_{n,(n),(i)}$-action on $\tcE_{(i),U^p(N_1,N_2)}^{\ord,+}$ by
\[ h(\xi_0,\xi_1)=(\nu(h)^{-1}\xi_0, (\circ h^{-1}) \circ \xi_1). \] 
The inverse system $\{ \tcE_{(i),U^p(N_1,N_2)}^{\ord,+} \}$ has an action both of $P_{n,(i)}^+(\A^\infty)^\ord$ and of $L_{n,(i),\lin}(\Z_{(p)})$. 

Suppose that $R_0$ is an irreducible noetherian $\Q$-algebra and that $\rho$ is a
representation of $R_{n,(n),(i)}$ on a finite, locally free $R_0$-module $W_\rho$. We define a locally
free sheaf $\cE^{\ord,+}_{(i),U^p(N_1,N_2),\rho}$ over
$\cA^{\ord,+}_{n,(i),U^p(N_1,N_2)} \times \Spec R_0$ by setting $\cE^{\ord,+}_{(i),U^p(N_1,N_2),\rho}(W)$ to be the set of $R^{(m)}_{n,(n)}(\cO_W)$-equivariant
maps of Zariski sheaves of sets
\[ \tcE^{\ord,+}_{(i),U^p(N_1,N_2)}|_W \lra W_\rho \otimes_{R_0} \cO_W. \]
Then $\{ \cE^{\ord,+}_{(i),U^p(N_1,N_2),\rho} \}$ is a system of locally free sheaves with $P_{n,(i)}^+(\A^\infty)^\ord$-action and $L_{n,(i),\lin}(\Z_{(p)})$-action over the system of schemes 
$\{ \cA^{\ord,+}_{n,(i),U^p(N_1,N_2)} \times \Spec R_0\}$. The restriction of $\cE^{\ord,+}_{(i),U^p(N_1,N_2),\rho}$ to $\cA^{(i),\ord}_{n-i,(hU^ph^{-1} \cap G_{n-i}^{(i)}(\A^{p,\infty}))(N_1,N_2)}$ can be identified with $\cE^{(i),\ord}_{(hU^ph^{-1} \cap G_{n-i}^{(i)}(\A^{p,\infty}))(N_1,N_2),\rho|_{R_{n-i,(n-i)}^{(i)}}}$. However the description of the actions of $P_{n,(i)}^+(\A^\infty)^\ord$ and $L_{n,(i),\lin}(\Z_{(p)})$ involve $\rho$ and not just $\rho|_{R^{(i)}_{n-i,(n-i)}}$. 
If $g \in P_{n,(i)}^+(\A^\infty)^{\ord,\times}$ and $\gamma \in L_{n,(i),\lin}(\Z_{(p)})$, then the natural maps
\[ g^*\cE_{(i),U^p(N_1,N_2),\rho}^{\ord,+} \lra \cE_{(i),(U^p)'(N_1',N_2'),\rho}^{\ord,+} \]
and
\[ \gamma^*\cE_{(i),U^p(N_1,N_2),\rho}^{\ord,+} \lra \cE_{(i),(U^p)'(N_1',N_2'),\rho}^{\ord,+} \]
are isomorphisms.
If $\rho$ factors through $R_{n,(n),(i)}/N(R_{n,(n),(i)})$ then $\cE^{\ord,+}_{(i),U^p(N_1,N_2),\rho}$ is canonically isomorphic to the pull-back of $\cE^{\ord,+}_{(i),U^p(N_1,N_2),\rho}$ from $\cX^{\ord,+}_{n,(i),U^p(N_1,N_2)}$. In general $W_\rho$ has a filtration by $R_{n,(n),(i)}$-invariant local direct-summands  such that the action of $R_{n,(n),(i)}$ on each graded piece factors through $R_{n,(n),(i)}/N(R_{n,(n),(i)})$. Thus $\cE^{\ord,+}_{(i),U^p(N_1,N_2),\rho}$ has a $P_{n,(i)}^+(\A^\infty)^\ord$ and $L_{n,(i),\lin}(\Z_{(p)})$ invariant filtration by local direct summands such that each graded piece is the pull back of some $\cE^{\ord,+}_{(i),U^p(N_1,N_2),\rho'}$ from $\cX^{\ord,+}_{n,(i),U^p(N_1,N_2)}$.

The next lemma follows from the discussion in section \ref{vb1}.
\begin{lem} If $U'$ 
is the image of $U$ (resp. $U^p$) and if $\pi$ denotes the map $A_{n,(i),U}^{(m),+} \ra A_{n,(i),U'}^+$ 
then there are $L_{n,(i),\lin}(\Q)$-equivariant, 
 $(P^{+}_{n,(i)}/Z(N_{n,(i)}))(\A^\infty)$-equivariant 
 and $GL_m(F)$-equivariant 
 isomorphisms
\[ \begin{array}{l} R^j\pi_* \Omega^k_{A_{n,(i),U}^{(m),+}/ A_{n,(i),U'}^+} \cong \\\left( \wedge^k (F^m \otimes_F \Omega^{+}_{n,(i),U'})\right) \otimes \left( \wedge^j (F^m\otimes_F \Hom(\Omega^{+}_{n,(i),U'},\Xi^{+}_{n,(i),U'}) )\right). \end{array} \]
\end{lem}

\newpage \subsection{Generalized mixed Shimura varieties.}

Next suppose $\tU$ is an open compact subgroup of $\tP^{(m),+}_{n,(i)}(\A^\infty)$. We define a split torus 
$\tS^{(m),+}_{n,(i),\tU}/Y_{n,(i),\tU}^{(m),+}$ as
\[ \coprod_{h \in L^{(m)}_{n,(i),\lin}(\A^\infty)/\tU}  S^{(m+i)}_{n-i,h \tU h^{-1} \cap \tG^{(m+i)}_{n-i}(\A^\infty)} .  \]
If $\tg \in \tP^{(m),+}_{n,(i)}(\A^\infty)$ and $\tg^{-1} \tU \tg \subset \tU'$, then we define
\[ \tg: \tS_{n,(i),\tU}^{(m),+} \lra \tS_{n,(i),\tU'}^{(m),+} \]
to be the coproduct of the maps
\[ \tg':  S^{(m+i)}_{n-i,h \tU h^{-1} \cap \tG^{(m+i)}_{n-i}(\A^\infty)} \lra  S^{(m+i)}_{n-i,h' \tU' (h')^{-1} \cap \tG^{(m+i)}_{n-i}(\A^\infty)}, \]
where $h,h' \in L^{(m)}_{n,(i),\lin}(\A^\infty)$ and $\tg' \in \tG^{(m+i)}_{n-i}(\A^\infty)$ satisfy $h\tg=\tg' h'$. This makes $\{ \tS^{(m),+}_{n,(i),\tU} \}$ a system of relative tori with right $\tP^{(m),+}_{n,(i)}(\A^\infty)$-action. If $\gamma \in L^{(m)}_{n,(i),\lin}(\Q)$, then we define
\[ \gamma: \tS_{n,(i),\tU}^{(m),+} \lra \tS_{n,(i),\tU}^{(m),+} \]
to be the coproduct of the maps
\[ \gamma: S^{(m+i)}_{n-i,h \tU h^{-1} \cap \tG^{(m+i)}_{n-i}(\A^\infty)} \lra S^{(m+i)}_{n-i,(\gamma h) \tU (\gamma h)^{-1} \cap \tP^{(m),\prime}_{n,(i)}(\A^\infty)}. \]
This gives a left action of $L^{(m)}_{n,(i),\lin}(\Q)$ on each $\tS_{n,(i),\tU}^{(m),+}$, which commutes with the action of $\tP^{(m),+}_{n,(i)}(\A^\infty)$. 

Similarly suppose $\tU^p$ is an open compact subgroup of $\tP^{(m),+}_{n,(i)}(\A^{p,\infty})$ and that $N$ is a non-negative integer. We define a split torus 
$\tcS^{(m),\ord,+}_{n,(i),\tU^p(N)}/Y_{n,(i),\tU^p(N)}^{(m),+}$ as
\[ \coprod_{h \in L^{(m)}_{n,(i),\lin}(\A^{\infty})^{\ord,\times}/\tU^p(N)}  \cS^{(m+i)}_{n-i,(h \tU^p h^{-1} \cap \tG^{(m+i)}_{n-i}(\A^{p,\infty}))(N,N')}  \]
for any $N'\geq N$.
If $\tg \in \tP^{(m),+}_{n,(i)}(\A^{\infty})^\ord$ and $\tg^{-1} \tU^p(N) \tg \subset (\tU^p)'(N')$, then we define
\[ \tg: \tcS_{n,(i),\tU^p(N)}^{(m),\ord,+} \lra \tcS_{n,(i),(\tU^p)'(N')}^{(m),\ord,+} \]
to be the coproduct of the maps
\[ \tg':  \cS^{(m+i)}_{n-i,(h \tU^p h^{-1} \cap \tG^{(m+i)}_{n-i}(\A^{p,\infty}))(N,N')} \lra  \cS^{(m+i)}_{n-i,(h' (\tU^p)' (h')^{-1} \cap \tG^{(m+i)}_{n-i}(\A^{p,\infty}))(N',N')}, \]
where $h,h' \in L^{(m)}_{n,(i),\lin}(\A^{\infty})^\ord$ and $\tg' \in \tG^{(m+i)}_{n-i}(\A^\infty)$ satisfy $h\tg=\tg' h'$. This makes $\{ \tcS^{(m),+}_{n,(i),\tU^p(N)} \}$ a system of relative tori with right $\tP^{(m),+}_{n,(i)}(\A^{\infty})^\ord$-action. If $\gamma \in L^{(m)}_{n,(i),\lin}(\Z_{(p)})$, then we define
\[ \gamma: \tcS_{n,(i),\tU^p(N)}^{(m),+} \lra \tcS_{n,(i),\tU^p(N)}^{(m),+} \]
to be the coproduct of the maps
\[ \gamma: \cS^{(m+i)}_{n-i,(h \tU^p h^{-1} \cap \tG^{(m+i)}_{n-i}(\A^{p,\infty}))(N,N)} \lra \cS^{(m+i)}_{n-i,((\gamma h) \tU^p (\gamma h)^{-1} \cap \tG^{(m),\prime}_{n,(i)}(\A^{p,\infty}))(N,N)}. \]
This gives a left action of $L^{(m)}_{n,(i),\lin}(\Z_{(p)})$ on each $\tcS_{n,(i),\tU^p(N)}^{(m),+}$, which commutes with the action of $\tP^{(m),+}_{n,(i)}(\A^{\infty})^\ord$. 

The sheaves 
$X^*(\tS_{n,(i),\tU}^{(m),+})$ and $X_*(\tS_{n,(i),\tU}^{(m),+})$ have actions of $L^{(m)}_{n,(i),\lin}(\Q)$. The sheaves 
 $X^*(\tcS_{n,(i),\tU^p(N)}^{(m),\ord,+})$ and $X_*(\tcS_{n,(i),\tU^p(N)}^{(m),\ord,+})$ have actions of $L^{(m)}_{n,(i),\lin}(\Z_{(p)})$.
The systems of sheaves $\{ X^*(\tS_{n,(i),\tU}^{(m),+}) \}$ and $\{ X_*(\tS_{n,(i),\tU}^{(m),+}) \}$ (resp. $\{ X^*(\tcS_{n,(i),\tU^p(N)}^{(m),\ord,+})\}$ and $\{X_*(\tcS_{n,(i),\tU^p(N)}^{(m),\ord,+})\}$) have actions of $\tP^{(m),+}_{n,(i)}(\A^\infty)$ (resp. $\tP^{(m),+}_{n,(i)}(\A^\infty)^\ord$).

The sheaf
\[ (X_*(\tS_{n,(i),\tU}^{(m),+}) \cap \Herm_{F^m})=\coprod_{h \in L^{(m)}_{n,(i),\lin}(\A^\infty)/\tU}  (X_*(\tS^{(m+i)}_{n-i,h \tU h^{-1} \cap \tG^{(m+i)}_{n-i}(\A^\infty)})\cap \Herm_{F^m}) \]
is a subsheaf of $X_*(\tS_{n,(i),\tU}^{(m),+})$. It is invariant by the actions of $L^{(m)}_{n,(i),\lin}(\Q)$ and $\tP^{(m),+}_{n,(i)}(\A^\infty)$. We define
a split torus
\[ \hS_{n,(i),\tU}^{(m),+}/Y_{n,(i),\tU}^{(m),+} \]
by
\[ X_*(\hS_{n,(i),\tU}^{(m),+}) = X_*(\tS_{n,(i),\tU}^{(m),+}) \cap \Herm_{F^m}, \]
and set
\[ S_{n,(i),\tU}^{(m),+}=\tS_{n,(i),\tU}^{(m),+}/\hS_{n,(i),\tU}^{(m),+}.\]
In the case $m=0$ we will write simply $S^+_{n,(i),\tU}$. The tori  $\tS_{n,(i),\tU}^{(m),+}$ and $S_{n,(i),\tU}^{(m),+}$ inherit a left action of $L^{(m)}_{n,(i),\lin}(\Q)$ and a right action of $\tP^{(m),+}_{n,(i)}(\A^\infty)$. In the case of $S_{n,(i),\tU}^{(m),+}$ the latter factors through $P^{(m),+}_{n,(i)}(\A^\infty)$.
If $\tU$ is a neat open compact subgroup of $\tP^{(m),+}_{n,(i)}(\A^\infty)$ with image $U$ in $P^{(m),+}_{n,(i)}(\A^\infty)$ and image $U'$ in $P^{+}_{n,(i)}(\A^\infty)$, then there is a natural, $L^{(m)}_{n,(i),\lin}(\Q)$-equivariant and $\tP^{(m),+}_{n,(i)}(\A^\infty)$-equivariant, commutative diagram:
\[ \begin{array}{ccccc}  \tS_{n,(i),\tU}^{(m),+} &\onto & S_{n,(i),U}^{(m),+} &\onto & S_{n,(i),U'}^+ \\ \da && \da && \da \\  Y_{n,(i),\tU}^{(m),+} & = & Y_{n,(i),U}^{(m),+} & \onto & Y_{n,(i),U'}^+. \end{array} \]

Similarly the sheaf
\[ \begin{array}{l} (X_*(\tcS_{n,(i),\tU^p(N_1,N_2)}^{(m),\ord,+}) \cap \Herm_{F^m})=\\\coprod_{h \in L^{(m)}_{n,(i),\lin}(\A^\infty)^{\ord,\times}/\tU^p(N)}  (X^*(\tcS^{(m+i),\ord}_{n-i,(h \tU^p h^{-1} \cap \tG^{(m+i)}_{n-i}(\A^{p,\infty}))(N,N')})\cap \Herm_{F^m}) \end{array}\]
is a sub-sheaf of $X_*(\tcS_{n,(i),\tU}^{(m),\ord,+})$. It is invariant by the actions of $L^{(m)}_{n,(i),\lin}(\Z_{(p)})$ and $\tP^{(m),+}_{n,(i)}(\A^\infty)^\ord$. We define
a split torus
\[ \hcS_{n,(i),\tU^p(N)}^{(m),\ord,+}/Y_{n,(i),\tU^p(N)}^{(m),\ord,+} \]
by
\[ X_*(\hcS_{n,(i),\tU^p(N)}^{(m),\ord,+}) = X_*(\tcS_{n,(i),\tU^p(N)}^{(m),\ord+}) \cap \Herm_{F^m}, \]
and set
\[ \cS_{n,(i),\tU^p(N)}^{(m),\ord,+}=\tcS_{n,(i),\tU^p(N)}^{(m),\ord,+}/\hcS_{n,(i),\tU^p(N)}^{(m),\ord,+}.\]
In the case $m=0$ we will write simply $\cS^+_{n,(i),\tU^p(N)}$. The tori  $\tcS_{n,(i),\tU^p(N)}^{(m),\ord,+}$ and $\cS_{n,(i),\tU^p(N)}^{(m),\ord,+}$ inherit a left action of $L^{(m)}_{n,(i),\lin}(\Z_{(p)})$ and a right action of $\tP^{(m),+}_{n,(i)}(\A^\infty)^\ord$. In the case of $\cS_{n,(i),\tU^p(N)}^{(m),\ord,+}$ the latter factors through $P^{(m),+}_{n,(i)}(\A^\infty)^\ord$.
If $\tU^p$ is a neat open compact subgroup of $\tP^{(m),+}_{n,(i)}(\A^{p,\infty})$ with image $U^p$ in $P^{(m),+}_{n,(i)}(\A^{p,\infty})$ and image $(U^p)'$ in $P^{+}_{n,(i)}(\A^{p,\infty})$, then there is a natural, $L^{(m)}_{n,(i),\lin}(\Z_{(p)})$-equivariant and $\tP^{(m),+}_{n,(i)}(\A^\infty)^\ord$-equivariant, commutative diagram:
\[ \begin{array}{ccccc}  \tcS_{n,(i),\tU^p(N)}^{(m),\ord,+} &\onto & \cS_{n,(i),U^p(N)}^{(m),\ord,+} &\onto & \cS_{n,(i),(U^p)'(N)}^{\ord,+} \\ \da && \da && \da \\  \cY_{n,(i),\tU^p(N)}^{(m),\ord,+} & = & \cY_{n,(i),U^p(N)}^{(m),\ord,+} & \onto & \cY_{n,(i),(U^p)'(N)}^{\ord,+}. \end{array} \]
There are natural equivariant embeddings
\[ \tcS^{(m),\ord,+}_{n,(i),\tU^p(N)} \times \Spec \Q \into \tS^{(m),+}_{n,(i),\tU^p(N)} \]
and
\[ \cS^{(m),\ord,+}_{n,(i),U^p(N)} \times \Spec \Q \into S^{(m),+}_{n,(i),U^p(N)} \]
and
\[ \hcS^{(m),\ord,+}_{n,(i),(U^p)'(N)} \times \Spec \Q \into \hS^{(m),+}_{n,(i),(U^p)'(N)}. \]

We write $X_*(S_{n,(i),\tU}^{(m),+})_\R^{\succ 0}$ (resp. $X_*(S_{n,(i),\tU}^{(m),+})_\R^{>0}$) for the sub-sheaves (of monoids) of $X_*(S_{n,(i),\tU}^{(m),+})_\R$ 
corresponding to the subset $\gC^{(m),\succ 0}_{(i)} \subset Z(N^{(m)}_{n,(i)})(\R)$ (resp. to the subset $\gC^{(m),> 0}_{(i)} \subset Z(N^{(m)}_{n,(i)})(\R)$).

We will also write $X^*(S_{n,(i),\tU}^{(m),+})_\R^{\geq 0}$ (resp. $X^*(S_{n,(i),\tU}^{(m),+})_\R^{>0}$, resp. $X^*(S_{n,(i),\tU}^{(m),+})^{\geq 0}$, resp. $X^*(S_{n,(i),\tU}^{(m),+})^{>0}$) for the sub-sheaves (of monoids) of $X^*(S_{n,(i),\tU}^{(m),+})_\R$ (resp. $X^*(S_{n,(i),\tU}^{(m),+})_\R$, resp. $X^*(S_{n,(i),\tU}^{(m),+})$, resp. $X^*(S_{n,(i),\tU}^{(m),+})$) consisting of sections that have non-negative (resp. strictly positive, resp. non-negative, resp. strictly positive) pairing with each section of $X_*(S_{n,(i),\tU}^{(m),+})_\R^{>0}$.
All these sheaves have (compatible) actions of $L^{(m)}_{n,(i),\lin}(\Q)$.
The system of sheaves $\{ X_*(S_{n,(i),\tU}^{(m),+}) \}$ has an action of $P^{(m),+}_{n,(i)}(\A^\infty)$, and the same is true for all the other systems of sheaves we are considering in this paragraph.

We may take the quotients of the sheaves $X^*(S_{n,(i),\tU}^{(m),+})$ (resp. $X^*(S_{n,(i),\tU}^{(m),+})^{>0}$, resp. $X^*(S_{n,(i),\tU}^{(m),+})^{\geq 0}$) by $L^{(m)}_{n,(i),\lin}(\Q)$ to give sheaves of sets on $Y_{n,(i),\tU}^{(m),\natural}$, which we will denote $X^*(S_{n,(i),\tU}^{(m),+})^\natural$ (resp. $X^*(S_{n,(i),\tU}^{(m),+})^{>0,\natural}$, resp. $X^*(S_{n,(i),\tU}^{(m),+})^{\geq 0,\natural}$).
If $y = h \tU$ lies in $Y^{(m),+}_{n,(i),\tU}$ above $y^\natural\in Y^{(m),\natural}_{n,(i),\tU}$ then the fibre of $X^*(S_{n,(i),\tU}^{(m),+})^\natural$ at $y^\natural$ equals
\[ \{ \gamma \in L^{(m)}_{n,(i),\lin}(\Q):\,\, \gamma y = y\} \backslash X^*(\tS^{(m+i)}_{n-i,h \tU h^{-1} \cap \tG^{(m+i)}_{n-i}(\A^\infty)}). \]

Similarly we will write $X_*(\cS_{n,(i),\tU^p(N)}^{(m),\ord,+})_\R^{\geq 0}$ (resp. $X_*(\cS_{n,(i),\tU^p(N)}^{(m),\ord,+})_\R^{>0}$) for the sub-sheaves (of monoids) of $X_*(\cS_{n,(i),\tU^p(N)}^{(m),\ord,+})_\R$  
corresponding to $\gC^{(m),\succ 0}_{(i)} \subset Z(N^{(m)}_{n,(i)})(\R)$ (resp. $\gC^{(m),> 0}_{(i)} \subset Z(N^{(m)}_{n,(i)})(\R)$).

Again we will write $X^*(\cS_{n,(i),\tU^p(N)}^{(m),\ord,+})_\R^{\geq 0}$ (resp. $X^*(\cS_{n,(i),\tU^p(N)}^{(m),\ord,+})^{\geq 0}$) for the sub-sheaves (of monoids) of $X^*(\cS_{n,(i),\tU^p(N)}^{(m),\ord,+})_\R$ (resp. $X^*(\cS_{n,(i),\tU^p(N)}^{(m),\ord,+})$) consisting of sections that have non-negative pairing with each section of $X_*(\cS_{n,(i),\tU^p(N)}^{(m),\ord,+})_\R^{>0}$. We will also write
$X^*(\cS_{n,(i),\tU^p(N)}^{(m),\ord,+})_\R^{>0}$ (resp. $X^*(\cS_{n,(i),\tU^p(N)}^{(m),\ord+})^{>0}$) for the sub-sheaves (of monoids) of $X^*(\cS_{n,(i),\tU^p(N)}^{(m),\ord,+})_\R$ (resp. $X^*(\cS_{n,(i),\tU^p(N)}^{(m),\ord,+})$) consisting of sections that have strictly positive pairing with each section of $X_*(\cS_{n,(i),\tU^p(N)}^{(m),\ord,+})_\R^{>0}$.
All these sheaves have (compatible) actions of $L^{(m)}_{n,(i),\lin}(\Z_{(p)})$.
The system of sheaves $\{ X_*(\cS_{n,(i),\tU^p(N)}^{(m),\ord,+}) \}$ has an action of $P^{(m),+}_{n,(i)}(\A^\infty)^\ord$, and the same is true for all the other systems of sheaves we are considering in this paragraph.

We may take the quotients of the sheaves $X^*(\cS_{n,(i),\tU^p(N)}^{(m),\ord,+})$ and $X^*(\cS_{n,(i),\tU^p(N)}^{(m),\ord,+})^{>0}$ and $X^*(\cS_{n,(i),\tU^p(N)}^{(m),\ord,+})^{\geq 0}$ by $L^{(m)}_{n,(i),\lin}(\Z_{(p)})$ to give sheaves of sets on $\cY_{n,(i),\tU^p(N)}^{(m),\ord,\natural}$, which we will denote $X^*(\cS_{n,(i),\tU^p(N)}^{(m),\ord,+})^\natural$ and $X^*(\cS_{n,(i),\tU^p(N)}^{(m),\ord,+})^{>0,\natural}$ and $X^*(\cS_{n,(i),\tU^p(N)}^{(m),\ord,+})^{\geq 0,\natural}$.

Suppose again that $\tU$ is a neat open compact subgroup of $\tP^{(m),+}_{n,(i)}(\A^\infty)$ and set 
\[ \tT_{n,(i),\tU}^{(m),+}=  \coprod_{h \in L^{(m)}_{n,(i),\lin}(\A^\infty)/\tU}  \tT^{(m+i)}_{n-i,h \tU h^{-1} \cap \tG^{(m+i)}_{n-i}(\A^\infty)} .  \]
It is an $\tS_{n,(i),\tU}^{(m),+}$-torsor over $A_{n,(i),\tU}^{(m),+}$. If $U$ denotes the image of $\tU$ in $P^{(m),+}_{n,(i)}(\A^\infty)$ then the push-out of $\tT_{n,(i),\tU}^{(m),+}$ under $\tS_{n,(i),\tU}^{(m),+} \onto S_{n,(i),U}^{(m),+}$ is an $S_{n,(i),U}^{(m),+}$-torsor over $A_{n,(i),U}^{(m),+}=A_{n,(i),\tU}^{(m),+}$, which only depends on $U$ (and not $\tU$), and which we will denote $T_{n,(i),U}^{(m),+}$. 
In the case $m=0$ we will write simply $T^+_{n,(i),U}$. Note that $\tT_{n,(i),\tU}^{(m),+}$ is an $\hS^{(m),+}_{n,(i),\tU}$-torsor over $T_{n,(i),\tU}^{(m),+}$.

If $\tg \in \tP^{(m),+}_{n,(i)}(\A^\infty)$ and $\tg^{-1} \tU \tg \subset \tU'$, then we define
\[ \tg: \tT_{n,(i),\tU}^{(m),+} \lra \tT_{n,(i),\tU'}^{(m),+} \]
to be the coproduct of the maps
\[ \tg': \tT^{(m+i)}_{n-i,h \tU h^{-1} \cap \tG^{(m+i)}_{n-i}(\A^\infty)}  \lra \tT^{(m+i)}_{n-i,h' \tU' (h')^{-1} \cap \tG^{(m+i)}_{n-i}(\A^\infty)} , \]
where $h,h' \in L^{(m)}_{n,(i),\lin}(\A^\infty)$ and $\tg' \in \tG^{(m+i)}_{n-i}(\A^\infty)$ satisfy $h\tg=\tg' h'$. This makes $\{ \tT^{(m),+}_{n,(i),\tU} \}$ a system of $\{\tS_{n,(i),\tU}^{(m),+} \} $-torsors over $\{ A_{n,(i),\tU}^{(m),+}\} $ with right $\tP^{(m),+}_{n,(i)}(\A^\infty)$-action. It also induces an action of $P^{(m),+}_{n,(i)}(\A^\infty)$ on $\{ T^{(m),+}_{n,(i),U} \}$, which makes $\{ T^{(m),+}_{n,(i),U} \}$ a system of $\{S_{n,(i),U}^{(m),+} \} $-torsors over $\{ A_{n,(i),U}^{(m),+}\} $ with right $P^{(m),+}_{n,(i)}(\A^\infty)$-action.
If $\gamma \in L^{(m)}_{n,(i),\lin}(\Q)$, then we define
\[ \gamma: \tT_{n,(i),\tU}^{(m),+} \lra \tT_{n,(i),\tU}^{(m),+} \]
to be the coproduct of the maps
\[ \gamma:  \tT^{(m+i)}_{n-i,h \tU h^{-1} \cap \tG^{(m+i)}_{n-i}(\A^\infty)}  \lra  \tT^{(m+i)}_{n-i,(\gamma h) \tU (\gamma h)^{-1} \cap \tG^{(m+i),\prime}_{n-i}(\A^\infty)} . \]
This gives a left action of $L^{(m)}_{n,(i),\lin}(\Q)$ on each $\tT_{n,(i),\tU}^{(m),+}$, which commutes with the action of $\tP^{(m),+}_{n,(i)}(\A^\infty)$. It induces a left action of $L^{(m)}_{n,(i),\lin}(\Q)$ on each $T_{n,(i),U}^{(m),+}$, which commutes with the action of $P^{(m),+}_{n,(i)}(\A^\infty)$.
Suppose that $\tU$ is a neat open compact subgroup of $\tP^{(m),+}_{n,(i)}(\A^\infty)$ with image $U$ in $P^{(m),+}_{n,(i)}(\A^\infty)$ and image $U'$ in $P^{+}_{n,(i)}(\A^\infty)$. Then there is a commutative diagram
\[ \begin{array}{ccc}  T_{n,(i),U}^{(m),+} &\onto &T_{n,(i),U'}^+ \\  \da&& \da \\  A_{n,(i),U}^{(m),+} &\onto & A_{n,(i),U'}^+ \\ \da&& \da \\  X_{n,(i),U}^{(m),+} &\onto & X_{n,(i),U'}^+ \\  \da&& \da \\  Y_{n,(i),U}^{(m),+} & \onto & Y_{n,(i),U'}^+. \end{array} \]
\[ \begin{array}{ccccc}  \tT_{n,(i),\tU}^{(m),+} &\onto & T_{n,(i),U}^{(m),+} &\onto &T_{n,(i),U'}^+ \\ \da && \da&& \da \\ A_{n,(i),U}^{(m),+} &=& A_{n,(i),U}^{(m),+} &\onto & A_{n,(i),U'}^+ \\ \da && \da&& \da \\ X_{n,(i),U}^{(m),+} &=& X_{n,(i),U}^{(m),+} &\onto & X_{n,(i),U'}^+ \\ \da && \da&& \da \\ Y_{n,(i),U}^{(m),+} &=& Y_{n,(i),U}^{(m),+} & \onto & Y_{n,(i),U'}^+. \end{array} \]
These diagrams are $L^{(m)}_{n,(i),\lin}(\Q)$-equivariant and $\tP^{(m),+}_{n,(i)}(\A^\infty)$-equivariant. We have
\[ T_{n,(i),U}^{(m),+}(\C)= P^{(m)}_{n,(i)}(\Q) \backslash (P^{(m),+}_{n,(i)}(\A) Z(N_{n,(i)}^{(m)})(\C))/ (U U_{n-i,\infty}^0 A_{n-i}(\R)^0).\]

Similarly if $\tU^p$ is a neat open compact subgroup of $\tP^{(m),+}_{n,(i)}(\A^{p,\infty})$ and $0 \leq N_1 \leq N_2$ we set 
\[ \tcT_{n,(i),\tU^p(N_1,N_2)}^{(m),\ord,+}=  \coprod_{h \in L^{(m)}_{n,(i),\lin}(\A^{\infty})^\ord/\tU^p(N_1,N_2)}  \cT^{(m+i),\ord}_{n-i,(h \tU^p h^{-1} \cap \tG^{(m+i)}_{n-i}(\A^{p,\infty}))(N_1,N_2)} .  \]
It is a $\tcS_{n,(i),\tU^p(N_1)}^{(m),\ord,+}$-torsor over $\cA_{n,(i),\tU^p(N_1,N_2)}^{(m),\ord,+}$. If $U^p$ denotes the image of $\tU^p$ in $P^{(m),+}_{n,(i)}(\A^{p,\infty})$ then the push-out of $\tcT_{n,(i),\tU^p(N_1,N_2)}^{(m),\ord,+}$ under $\tcS_{n,(i),\tU^p(N_1)}^{(m),\ord,+} \onto \cS_{n,(i),U^p(N_1)}^{(m),\ord,+}$ is a $\cS_{n,(i),U^p(N_1)}^{(m),\ord,+}$-torsor over $\cA_{n,(i),U^p(N_1,N_2)}^{(m),\ord,+}$, which only depends on $U^p$ (and not $\tU^p$) and $N_1,N_2$, and which we will denote $\cT_{n,(i),U^p(N_1,N_2)}^{(m),\ord,+}$. 
In the case $m=0$ we will write simply $\cT^{\ord,+}_{n,(i),U^p(N_1,N_2)}$. Note that $\tcT_{n,(i),\tU^p(N_1,N_2)}^{(m),\ord,+}$ is a $\hcS^{(m),\ord,+}_{n,(i),\tU^p(N_1)}$-torsor over $\cT_{n,(i),\tU^p(N_1,N_2)}^{(m),\ord,+}$.

As above the system $\{ \tcT_{n,(i),\tU^p(N_1,N_2)}^{(m),\ord,+} \}$ has a right action of $\tP_{n,(i)}^{(m),+}(\A^{\infty})^\ord$ and a commuting left action of $L_{n,(i),\lin}^{(m)}(\Z_{(p)})$. If $g \in \tP_{n,(i)}^{(m),+}(\A^{\infty})^{\ord,\times}$ then the map $g$ is finite etale. The map
\[ \varsigma_p: \tcT_{n,(i),\tU^p(N_1,N_2)}^{(m),\ord,+}\times \Spec \F_p \lra \tcT_{n,(i),\tU^p(N_1,N_2-1)}^{(m),\ord,+} \times \Spec \F_p\]
equals absolute Frobenius composed with the forgetful map. If $N_2>1$ then the map
\[ \varsigma_p: \tcT_{n,(i),\tU^p(N_1,N_2)}^{(m),\ord,+} \lra \tcT_{n,(i),\tU^p(N_1,N_2-1)}^{(m),\ord,+} \]
 is finite flat. Further there is a left action of $GL_m(\cO_{F,(p)})$ such that if $\delta \in GL_m(\cO_{F,(p)})$ and $\gamma \in L^{(m)}_{n,(i),\lin}(\Z_{(p)})$ and $g \in \tP_{n,(i)}^{(m),+}(\A^{\infty})^\ord$, then $\gamma$ followed by $\delta$ equals $\delta$ followed by $\delta \gamma \delta^{-1}$, and $g$ followed by $\delta$ equals $\delta$ followed by $\delta g \delta^{-1}$. These actions are also all compatible with the actions on $\{ \tcS_{n,(i),\tU^p(N_1)}^{(m),\ord,+} \}$. There are induced actions of the groups $GL_m(\cO_{F,(p)})$ and $L^{(m)}_{n,(i),\lin}(\Z_{(p)})$ and $P_{n,(i)}^{(m),+}(\A^{\infty})^\ord$ on $\{  \cT_{n,(i),\tU^p(N_1,N_2)}^{(m),\ord,+} \}$. There are equivariant commutative diagrams
 \[ \begin{array}{ccc} \cT_{n,(i),\tU^p(N_1,N_2)}^{(m),\ord,+} &\onto &\cT_{n,(i),\tU^p(N_!,N_2)}^{\ord,+} \\  \da&& \da \\  \cA_{n,(i),\tU^p(N_1,N_2)}^{(m),\ord,+} &\onto & \cA_{n,(i),\tU^p(N_1,N_2)}^{\ord,+} \\  \da&& \da \\  \cX_{n,(i),\tU^p(N_1,N_2)}^{(m),\ord,+} &\onto & \cX_{n,(i),\tU^p(N_1,N_2)}^{\ord,+} \\  \da&& \da \\  \cY_{n,(i),\tU^p(N_1)}^{(m),\ord,+} & \onto & \cY_{n,(i),\tU(N_1)}^{\ord,+}. \end{array} \]
\[ \begin{array}{ccccc}  \tcT_{n,(i),\tU^p(N_1,N_2)}^{(m),\ord,+} &\onto & \cT_{n,(i),\tU^p(N_1,N_2)}^{(m),\ord,+} &\onto &\cT_{n,(i),\tU^p(N_!,N_2)}^{\ord,+} \\ \da && \da&& \da \\ \cA_{n,(i),\tU^p(N_1,N_2)}^{(m),\ord,+} &=& \cA_{n,(i),\tU^p(N_1,N_2)}^{(m),\ord,+} &\onto & \cA_{n,(i),\tU^p(N_1,N_2)}^{\ord,+} \\ \da && \da&& \da \\ \cX_{n,(i),\tU^p(N_1,N_2)}^{(m),\ord,+} &=& \cX_{n,(i),\tU^p(N_1,N_2)}^{(m),\ord,+} &\onto & \cX_{n,(i),\tU^p(N_1,N_2)}^{\ord,+} \\ \da && \da&& \da \\ \cY_{n,(i),\tU^p(N_1)}^{(m),\ord,+} &=& \cY_{n,(i),\tU^p(N_1)}^{(m),\ord,+} & \onto & \cY_{n,(i),\tU(N_1)}^{\ord,+}. \end{array} \]
There are natural equivariant embeddings
\[ \tcT_{n,(i),\tU^p(N_1,N_2)}^{(m),\ord,+} \times \Spec \Q \into \tT_{n,(i),\tU^p(N_1,N_2)}^{(m),+} \]
and
\[ \cT_{n,(i),U^p(N_1,N_2)}^{(m),\ord,+} \times \Spec \Q \into T_{n,(i),U^p(N_1,N_2)}^{(m),\ord,+}. \]

If $a$ is a global section of 
$X^*(\tS^{(m),+}_{n,(i),\tU})$ (resp. $X^*(\hS^{(m),+}_{n,(i),\tU})$, resp. 
$X^*(S^{(m),+}_{n,(i),\tU})$ over $Y^{(m),+}_{n,(i),\tU}$ then we can associate to it a line bundle
\[ \cL_{\tU}^+(a) \]
over $A^{(m),+}_{n,(i),\tU}$. 
(resp. $T^{(m),+}_{n,(i),\tU}$, resp. $A^{(m),+}_{n,(i),\tU}$). 
There are natural isomorphisms 
\[ \cL_{\tU}^+(a)\otimes \cL_{\tU}^+(a') \cong \cL_{\tU}^+(a+a'). \]
If $\ta \in X^*(\tS^{(m),+}_{n,(i),\tU})$ has image $\hata \in X^*(\hS^{(m),+}_{n,(i),\tU})$ then $\cL_{\tU}^+(\ta)/A^{(m),+}_{n,(i),\tU}$ pulls back to $\cL_{\tU}^+(\hata)$ over
$T^{(m),+}_{n,(i),\tU}$. 

Suppose that $R_0$ is an irreducible, noetherian $\Q$-algebra. Suppose also that $U$ is a neat open compact subgroup of $P_{n,(i)}^+(\A^\infty)$.
If $a$ is a global section of $X^*(S^{+}_{n,(i),U})^{>0}$ then $\cL_{U}^+(a)$ is relatively ample for $A^{+}_{n,(i),U}/X^{+}_{n,(i),U}$. 
If $\pi^+$ denotes the map
\[ A^+_{n,(i),U} \times \Spec R_0 \lra X^+_{n,(i),U} \times \Spec R_0, \]
then we see that
\[ R^i\pi^+_* \cL_{U}^+(a) = (0) \]
for $i>0$. (Because $A_{n,(i),U}^+/X_{n,(i),U}^+$ is a torsor for an abelian scheme and $\cL_U^+(a)$ is relatively ample for this morphism.) We will denote by $(\pi_{A^+/X^+,*} \cL)^+_U(a)$ the image
$\pi^+_* \cL_U^+(a)$. Suppose further that $\cF$ is a locally free sheaf on $X^+_{n,(i),U} \times \Spec R_0$ with $L_{n,(i),\lin}(\Q)$-action. If $a^\natural$ is a section of $X^*(S^{+}_{n,(i),U})^{>0,\natural}$ we will define 
\[ (\pi_{A^+/X^\natural,*} \cL \otimes \cF)^+_U(a^\natural) \]
as follows: Over a point $y^\natural$ of $Y_{n,(i),U}^\natural$ we take the sheaf
\[ \prod_{y,a} (\pi_{A^+/X^+,*} \cL)^+_{U}(a)_y \otimes \cF_y \]
over $X^\natural_{n,(i),U,y^\natural} \times \Spec R_0$, where $y$ runs over points of $Y_{n,(i),U}^+ $ above $y^\natural$ and $a$ runs over sections of $X^*(S_{n,(i),U}^+)_y$ above $a^\natural$. It is a sheaf with an action of $L_{n,(i),\lin}(\Q)$.

\begin{lem} \label{subvan1} Keep the notation and assumptions of the previous paragraph. \begin{enumerate}
\item 
\[ (\pi_{A^+/X^\natural,*}\cL \otimes \cF)_{U}^+(a^\natural) \cong \Ind_{\{1\}}^{L_{n,(i),\lin}(\Q)} (\pi_{A^+/X^\natural,*}\cL \otimes \cF)_{U}^+(a^\natural)^{L_{n,(i),\lin}(\Q)} \]
as a sheaf on $X^{\natural}_{n,(i),U} \times \Spec R_0$ with $L_{n,(i),\lin}(\Q)$-action.

\item If 
\[ \pi: A^{+}_{n,(i),U} \times \Spec R_0 \lra X^{\natural}_{n,(i),U} \times \Spec R_0 \]
then
\[ \begin{array}{l} R^i\pi_* \prod_{a \in X^*(S^{+}_{n,(i),U})^{>0}} (\cL^+_{U}(a)\otimes \cF) \\ \cong \left\{ \begin{array}{ll} \prod_{a^\natural \in X^*(S^{+}_{n,(i),U})^{>0,\natural}} (\pi_{A^+/X^\natural,*}\cL \otimes \cF)_{U}^+(a^\natural)& {\rm if}\,\,\, i=0 \\ (0) & {\rm otherwise}. \end{array} \right. \end{array} \]
\end{enumerate} \end{lem}

\pfbegin For the first part note that if $y$ in $Y^+_{n,(i),U}$ and if $a \in X^*(S^+_{n,(i),U})_y$ then the stabilizer of $a$ in $\{ \gamma \in L_{n,(i)}(\Q): \,\, \gamma y=y\}$  is finite, and that if $U$ is neat then it is trivial. The second part follows from the observations of the previous paragraph together with proposition 0.13.3.1 of \cite{ega3}. 
\pfend

Similarly if $a$ is a global section of 
$X^*(\tcS^{(m),\ord,+}_{n,(i),\tU^p(N_1)})$ (resp. $X^*(\hcS^{(m),\ord,+}_{n,(i),\tU^p(N_1)})$, resp. 
$X^*(\cS^{(m),\ord,+}_{n,(i),\tU^p(N_1)})$ over $\cY^{(m),\ord,+}_{n,(i),\tU^p(N_1)}$ then we can associate to it a line bundle
\[ \cL_{\tU^p(N_1,N_2)}^+(a) \]
over $\cA^{(m),\ord,+}_{n,(i),\tU^p(N_1,N_2)}$. 
(resp. $\cT^{(m),\ord,+}_{n,(i),\tU^p(N_1,N_2)}$, resp. $\cA^{(m),\ord,+}_{n,(i),\tU^p(N_1,N_2)}$). 
There are natural isomorphisms 
\[ \cL_{\tU^p(N_1,N_2)}^+(a)\otimes \cL_{\tU^p(N_1,N_2)}^+(a') \cong \cL_{\tU^p(N_1,N_2)}^+(a+a'). \]
Suppose a section $\ta$ of $X^*(\tcS^{(m),\ord,+}_{n,(i),\tU(N_1)})$ has image $\hata \in X^*(\hcS^{(m),\ord,+}_{n,(i),\tU^p(N_1,N_2)})$ then the sheaf $\cL_{\tU^p(N_1,N_2)}^+(\ta)$ over $\cA^{(m),\ord,+}_{n,(i),\tU^p(N_1,N_2)}$ pulls back to the sheaf $\cL_{\tU^p(N_1,N_2)}^+(\hata)$ over
$\cT^{(m),\ord,+}_{n,(i),\tU^p(N_1,N_2)}$. 

Suppose that $R_0$ is an irreducible, noetherian $\Z_{(p)}$-algebra. Suppose that $U^p$ is a neat open compact subgroup of $P^+_{n,(i)}(\A^{p,\infty})$ and that $0 \leq N_1 \leq N_2$. If $a$ is a section of $X^*(\cS^{+}_{n,(i),U^p(N_1)})^{>0}$ then $\cL_{U^p(N_1,N_2)}^+(a)$ is relatively ample for $\cA^{\ord,+}_{n,(i),U^p(N_1,N_2)}$ over $\cX^{\ord,+}_{n,(i),U^p(N_1,N_2)}$.
If $\pi^+$ denotes the map
\[ \cA^{\ord,+}_{n,(i),U^p(N_1,N_2)} \times \Spec R_0 \lra \cX_{n,(i),U^p(N_1,N_2)}^{\ord,+}  \times \Spec R_0\]
then we see that
\[ R^i\pi_*^+ \cL_{U^p(N_1,N_2)}^+(a)=(0) \]
for $i>0$. (Again because $\cA_{n,(i),U^p(N_1,N_2)}^{\ord,+}/\cX_{n,(i),U^p(N_1,N_2)}^{\ord,+}$ is a torsor for an abelian scheme and $\cL_{U^p(N_1,N_2)}^+(a)$ is relatively ample for this morphism.) We will denote by $(\pi_{\cA^{\ord,+}/\cX^{\ord,+},*} \cL)^+_{U^p(N_1,N_2)}(a)$ the image $\pi^+_* \cL^+(a)$. Suppose further that $\cF$ is a locally free sheaf on $\cX^{\ord,+}_{n,(i),U^p(N_1,N_2)} \times \Spec R_0$ with $L_{n,(i),\lin}(\Z_{(p)})$-action. If $a^\natural$ is a section of  $X^*(\cS^{\ord,+}_{n,(i),U^p(N_1)})^{>0,\natural}$ we define a sheaf
\[ (\pi_{\cA^{\ord,+}/\cX^{\ord,\natural},*}\cL \otimes \cF)_{U^p(N_1,N_2)}^+(a^\natural) \]
as follows: Over a point $y^\natural$ of $\cY_{n,(i),U^p(N_1,N_2)}^{\ord,\natural}$ we take the sheaf
\[ \prod_{y,a} (\pi_{\cA^{\ord,+}/\cX^{\ord,+},*} \cL)^+_{U^p(N_1,N_2)}(a)_y \otimes \cF_y \]
over $\cX^{\ord,\natural}_{n,(i),U^p(N_1,N_2),y^\natural} \times \Spec R_0$, where $y$ runs over points of $\cY^{\ord,\natural}_{n,(i),U^p(N_1,N_2)}$ above $y^\natural$ and $a$ runs over sections of $X^*(\cS^{\ord,+}_{n,(i),U^p(N_1)})_y$ above $a^\natural$. It is a sheaf with an action of $L_{n,(i),\lin}(\Z_{(p)})$. As above we have the following lemma.

\begin{lem} \label{subvan1ord} Keep the notation and assumptions of the previous paragraph.
\begin{enumerate}
\item 
\[ \begin{array}{l} (\pi_{\cA^{\ord,+}/\cX^{\ord,\natural},*}\cL \otimes \cF)_{U^p(N_1,N_2)}^+(a^\natural) \cong \\ \Ind_{\{1\}}^{L_{n,(i),\lin}(\Z_{(p)})} (\pi_{\cA^{\ord,+}/\cX^{\ord,+},*}\cL \otimes\cF)_{U^p(N_1,N_2)}^+(a^\natural)^{L_{n,(i),\lin}(\Z_{(p)})} \end{array} \]
as a sheaf on $\cX^{\ord,\natural}_{n,(i),U^p(N_1,N_2)} \times \Spec R_0$ with $L_{n,(i),\lin}(\Z_{(p)})$-action.

\item If 
\[ \pi: \cA^{\ord,+}_{n,(i),U^p(N_1,N_2)}  \times \Spec R_0\lra \cX^{\ord,\natural}_{n,(i),U^p(N_1,N_2)}  \times \Spec R_0\]
then
\[ R^i\pi_* \prod_{a \in X^*(\cS^{\ord,+}_{n,(i),U^p(N_1)})^{>0}} (\cL^+_{U^p(N_1,N_2)}(a) \otimes \cF)\]
is isomorphic to
\[ \prod_{a^\natural \in X^*(\cS^{\ord,+}_{n,(i),U^p(N_1)})^{>0,\natural}} (\pi_{\cA^{\ord,+}/\cX^{\ord,\natural},*}\cL\otimes \cF)_{U^p(N_1,N_2)}^+(a^\natural) \]
if $i=0$, and otherwise is $(0)$.
\end{enumerate} \end{lem}

\newpage \subsection{Partial compactifications.}\label{partcomp}

We will now turn to the partial compactification of the generalized Shimura varieties, $T^{(m)}_{n,(i),U}$, we discussed in the last section. These will serve as models for the full compactification of the $A^{(m)}_{n,U}$, which near the boundary can be formally modeled on the partial compactifications of the $T^{(m)}_{n,(i),U}$. 

Suppose that $U$ (resp. $U^p$) is a neat open compact subgroup of $L_{n,(i),\lin}^{(m)}(\A^\infty)$ (resp. $L_{n,(i),\lin}^{(m)}(\A^\infty)$) and that $N$ is a non-negative integer.
By an {\em admissible cone decomposition} $\Sigma_0$ for $X_*(S^{(m),+}_{n,(i),U})_\R^{\succ 0}$ (resp. $X_*(\cS^{(m),\ord,+}_{n,(i),U^p(N)})_\R^{\succ 0}$)  we shall mean a partial fan $\Sigma_0$ in $X_*(S^{(m),+}_{n,(i),U})_\R$ (resp. $X_*(\cS^{(m),\ord,+}_{n,(i),U^p(N)})_\R$) such that
\begin{itemize}
\item $|\Sigma_0|=X_*(S^{(m),+}_{n,(i),U})_\R^{\succ 0}$ (resp. $X_*(\cS^{(m),\ord,+}_{n,(i),U^p(N)})_\R^{\succ 0}$);
\item $|\Sigma_0|^0=X_*(S^{(m),+}_{n,(i),U})_\R^{> 0}$ (resp. $X_*(\cS^{(m),\ord,+}_{n,(i),U^p(N)})_\R^{> 0}$);
\item $\Sigma_0$ is invariant under the left action of $L_{n,(i),\lin}^{(m)}(\Q)$ (resp. $L_{n,(i),\lin}^{(m)}(\Z_{(p)})$);
\item $L_{n,(i),\lin}^{(m)}(\Q) \backslash \Sigma_0$ (resp. $L_{n,(i),\lin}^{(m)}(\Z_{(p)}) \backslash \Sigma_0$) is a finite set;
\item if $\sigma \in \Sigma_0$ and $1 \neq \gamma \in L_{n,(i),\lin}^{(m)}(\Q)$ (resp. $L_{n,(i),\lin}^{(m)}(\Z_{(p)})$) then
\[ \sigma \cap \gamma \sigma \not\in  \Sigma_0. \]
\end{itemize}
(Many authors would not include the last condition in the definition of an `admissible cone decomposition'.)
In concrete terms $\Sigma_0$ consists of a partial fan $\Sigma_{g,0}$ in $Z(N_{n,(i)}^{(m)})(\R)$ for each $g \in L_{n,(i),\lin}^{(m)}(\A^\infty)$ (resp. $L_{n,(i),\lin}^{(m)}(\A^\infty)^{\ord,\times}$), such that
\begin{itemize}
\item $\Sigma_{\gamma g u} = \gamma \Sigma_g$
for all $\gamma \in L_{n,(i),\lin}^{(m)}(\Q)$ (resp. $L_{n,(i),\lin}^{(m)}(\Z_{(p)})$) and $u \in U$ (resp. $U^p(N)$);
\item $|\Sigma_{g,0}|=\gC^{(m),\succ 0}(V_{n,(i)})$ and $|\Sigma_{g,0}|^0=\gC^{(m),>0}(V_{n,(i)})$ for each $g$;
\item $(L_{n,(i),\lin}^{(m)}(\Q) \cap gUg^{-1})\backslash \Sigma_{g,0}$ (resp. $(L_{n,(i),\lin}^{(m)}(\Z_{(p)}) \cap gU^p(N)g^{-1})\backslash \Sigma_g$) is finite for all $g$;
\item for each $g$ and each $\sigma \in \Sigma_{g,0}$, if $1 \neq \gamma \in (L_{n,(i),\lin}^{(m)}(\Q) \cap gUg^{-1})$ (resp. $(L_{n,(i),\lin}^{(m)}(\Z_{(p)}) \cap gU^p(N)g^{-1})$) and 
\[ \sigma \cap \gamma \sigma \not\in \Sigma_{g,0}. \]
\end{itemize}

Note that an admissible cone decomposition for  $X_*(S^{(m),+}_{n,(i),U^p(N)})_\R^{\succ 0}$ induces (by restriction) one for $X_*(\cS^{(m),\ord,+}_{n,(i),U^p(N)})_\R^{\succ 0}$. This sets up a bijection between admissible cone decompositions for $X_*(S^{(m),+}_{n,(i),U^p(N)})_\R^{\succ 0}$ and for $X_*(\cS^{(m),\ord,+}_{n,(i),U^p(N)})_\R^{\succ 0}$. 

\begin{lem} Suppose that $U$ (resp. $U^p$) is a neat open compact subgroup of $L_{n,(i),\lin}^{(m)}(\A^\infty)$ (resp. $L_{n,(i),\lin}^{(m)}(\A^\infty)$) and that $N$ is a non-negative integer. Suppose also that $\Sigma_0$ is
an admissible cone decomposition for $X_*(S^{(m),+}_{n,(i),U})_\R^{\succ 0}$ (resp. $X_*(\cS^{(m),\ord,+}_{n,(i),U^p(N)})_\R^{\succ 0}$). Also suppose that $\tau \subset |\Sigma_0|$ is a rational polyhedral cone. Then
the set 
\[ \{ \gamma \in L_{n,(i),\lin}^{(m)}(\Q): \,\, \gamma \tau \cap \tau \cap |\Sigma_0|^0 \neq \emptyset \}  \]
(resp.
\[ \{ \gamma \in L_{n,(i),\lin}^{(m)}(\Z_{(p)}): \,\, \gamma \tau \cap \tau \cap |\Sigma_0|^0 \neq \emptyset \} ) \] 
is finite.
\end{lem}

\pfbegin
We treat the case of $X_*(S^{(m),+}_{n,(i),U})_\R^{\succ 0}$, the other being exactly similar. 
Suppose that $\tau$ has support $y=hU$ and set $\Gamma = L_{n,(i),\lin}^{(m)}(\Q) \cap hUh^{-1}$ a discrete subgroup of $L_{n,(i),\lin}^{(m)}(\Q)$. We certainly have
\[ \{ \gamma \in L_{n,(i),\lin}^{(m)}(\Q): \,\, \gamma \tau \cap \tau \cap |\Sigma_0|^0 \neq \emptyset \}=\{ \gamma \in \Gamma: \,\, \gamma \tau \cap \tau \cap |\Sigma_0|^0(y) \neq \emptyset \}. \]
That this set is finite follows from theorem II.4.6 and the remark (ii) at the end of section II.4.1
of \cite{amrt}. 
\pfend

\begin{cor} If $\Sigma_0$ is an admissible cone decomposition for $X_*(S^{(m),+}_{n,(i),U})_\R^{\succ 0}$ or $X_*(\cS^{(m),\ord,+}_{n,(i),U^p(N)})_\R^{\succ 0}$, then $\Sigma_0$ is locally finite.\end{cor}

\pfbegin Let $\tau \subset |\Sigma_0|$ be a rational polyhedral cone. Let $\sigma_1,...,\sigma_r$ be representatives for $L_{n,(i),\lin}^{(m)}(\A^\infty)\backslash \Sigma_0$ (resp. $L_{n,(i),\lin}^{(m)}(\Z_{(p)}) \backslash \Sigma_0$); and suppose they are chosen with the same support as $\tau$ whenever possible. Let $\tau'$ be the rational polyhedral cone spanned by $\tau$ and those $\sigma_i$ with the same support as $\tau$. If $\gamma \in L_{n,(i),\lin}^{(m)}(\A^\infty)$ (resp. $L_{n,(i),\lin}^{(m)}(\Z_{(p)})$) and 
\[  \gamma \sigma_i \cap \tau \cap |\Sigma_0|^0 \neq \emptyset, \]
then 
\[ \gamma \tau' \cap \tau' \cap |\Sigma_0|^0 \neq \emptyset \]
and so by the previous lemma $\gamma$ lies in a finite set. The corollary follows.
\pfend

If $g \in P^{(m),+}_{n,(i)}(\A^\infty)$, if $U' \supset g^{-1} U g$ are neat open compact subgroups of the group $P^{(m),+}_{n,(i)}(\A^\infty)$, and if $\Sigma'_0$ is a $U'$-admissible cone decomposition for $X_*(S_{n,(i),U'}^{(m),+})_\R^{\succ 0}$, then $\Sigma'_0 g^{-1}$ is a $U$-admissible cone decomposition for $X_*(S_{n,(i),U}^{(m),+})_\R^{\succ 0}$. 
We will call a $U$-admissible cone decomposition $\Sigma_0$ for $X_*(S_{n,(i),U}^{(m),+})_\R^{\succ 0}$ {\em compatible} with $\Sigma_0'$ with respect to $g$ if $\Sigma_0$ refines $\Sigma'_0 g^{-1}$.
Similarly if $g \in P^{(m),+}_{n,(i)}(\A^{\infty})^\ord$, if $(U^p)' (N')\supset (g^{-1} U^p g)(N)$, and if $\Sigma'_0$ is a $(U^p)'(N')$-admissible cone decomposition for $X_*(\cS_{n,(i),(U^p)'(N')}^{(m),\ord,+})_\R^{\succ 0}$, then $(\Sigma' g^{-1},\Sigma'_0 g^{-1})$ is an admissible cone decomposition for $X_*(\cS_{n,(i),U^p(N)}^{(m),\ord,+})_\R^{\succ 0}$. 
We will call a $U^p(N)$-admissible cone decomposition $\Sigma_0$ for $X_*(\cS_{n,(i),U^p(N)}^{(m),\ord,+})_\R^{\succ 0}$ {\em compatible} with $\Sigma_0'$ with respect to $g$ if $\Sigma_0$ refines $\Sigma'_0 g^{-1}$.

If $U'$ is a neat open compact subgroup of $P_{n,(i)}^+(\A^\infty)$ which contains the image of $U$, we will call an admissible cone decomposition $\Sigma_0$ of $X_*(S_{n,(i),U}^{(m),+})_\R^{\succ 0}$ and an admissible cone decomposition $\Delta_0$ of $X_*(S_{n,(i),U'}^{+})_\R^{\succ 0}$ {\em compatible} if, under the natural map
\[ X_*(S_{n,(i),U}^{(m),+})_\R^{\succ 0} \onto X_*(S_{n,(i),U'}^{+})_\R^{\succ 0}, \]
the image of each $\sigma \in \Sigma_0$ is contained in some element of $\Delta_0$.
Similarly if $(U^p)'$ is a neat open compact subgroup of $P_{n,(i)}^+(\A^{p,\infty})$ which contains the image of $U^p$ and if $N' \geq N$, we will call an admissible cone decomposition $\Sigma_0$ of $X_*(\cS_{n,(i),U^p(N)}^{(m),\ord,+})_\R^{\succ 0}$ and an admissible cone decomposition $\Delta_0$ of $X_*(\cS_{n,(i),(U^p)'(N')}^{\ord,+})_\R^{\succ 0}$ {\em compatible} if, under the natural map
\[ X_*(\cS_{n,(i),U^p(N)}^{(m),\ord,+})_\R^{\succ 0} \onto X_*(\cS_{n,(i),(U^p)'(N')}^{\ord,+})_\R^{\succ 0}, \]
the image of each $\sigma \in \Sigma_0$ is contained in some element of $\Delta_0$.

If $\Sigma_0$ is a smooth admissible cone decomposition of $X_*(S_{n,(i),U}^{(m),+})_\R^{\succ 0}$ (resp. of $X_*(\cS_{n,(i),U^p(N_1,N_2))}^{(m),\ord,+})_\R^{\succ 0}$), then the log smooth, log scheme $(T^{(m),+}_{n,(i),U,\tSigma_0},\cM_{\tSigma_0})$ (resp. $(\cT^{(m),\ord,+}_{n,(i),U^p(N_1,N_2),\tSigma_0},\cM_{\tSigma_0})$) has a left action of $L_{n,(i),\lin}^{(m)}(\Q)$ (resp. $L_{n,(i),\lin}^{(m)}(\Z_{(p)})$) extending that on $T^{(m),+}_{n,(i),U}$ (resp. $\cT^{(m),\ord,+}_{n,(i),U^p(N_1,N_2)}$). If $g \in P^{(m),+}_{n,(i)}(\A^\infty)$ (resp. $g \in P^{(m),+}_{n,(i)}(\A^{\infty})^\ord$) and if $\Sigma_0$ is compatible with $\Sigma'_0$ with respect to $g$ then the map
\[ g: T^{(m),+}_{n,(i),U} \lra T^{(m),+}_{n,(i),U'} \]
(resp.
\[ g: \cT^{(m),\ord,+}_{n,(i),U^p(N_1,N_2)} \lra \cT^{(m),\ord,+}_{n,(i),(U^p)'(N_1',N_2')}) \]
uniquely extends to an $L^{(m)}_{n,(i),\lin}(\Q)$-equivariant (resp. $L^{(m)}_{n,(i),\lin}(\Z_{(p)})$-equivariant) log smooth map
\[ g: (T^{(m),+}_{n,(i),U,\Sigma_0},\cM_{\Sigma_0}) \lra (T^{(m),+}_{n,(i),U',\Sigma_0'},\cM_{\Sigma_0'}) \]
(resp.
\[ g: (\cT^{(m),\ord,+}_{n,(i),U^p(N_1,N_2),\tSigma_0},\cM_{\tSigma_0}) \lra (\cT^{(m),\ord,+}_{n,(i),(U^p)'(N'_1,N_2'),\tSigma_0'},\cM_{\tSigma_0'})). \]
This makes $\{ (T^{(m),+}_{n,(i),U,\tSigma_0},\cM_{\tSigma_0}) \}$ (resp. 
$ \{ ( \cT^{(m),\ord,+}_{n,(i),U^p(N_1,N_2),\tSigma_0},\cM_{\tSigma_0})\}$) a system of log schemes with 
$P^{(m),+}_{n,(i)}(\A^\infty)$-action (resp. $P^{(m),+}_{n,(i)}(\A^{\infty})^\ord$-action).
There are equivariant embeddings
\[ ( \cT^{(m),\ord,+}_{n,(i),U^p(N_1,N_2),\tSigma_0} \times \Spec \Q,\cM_{\tSigma_0}) \into ( T^{(m),+}_{n,(i),U^p(N_1,N_2),\tSigma_0},\cM_{\tSigma_0}). \]

We have
\[ \begin{array}{l} |\cS(\partial  T^{(m),+}_{n,(i),U,\tSigma_0})|-|\cS(\partial  T^{(m),+}_{n,(i),U,\tSigma_0-\Sigma_0})|= \\ (L_{n,(i),\lin}^{(m)}(\A^\infty) \times (C_{n-i}(\Q)  \backslash C_{n-i}(\A)/ C_{n-i}(\R)^0)) /U \times (\gC^{(m),> 0}_{(i)}/\R^\times_{>0}). \end{array} \]

If $U'$ (resp. $(U')^p$) is a neat open compact subgroup of $P_{n,(i)}^+(\A^\infty)$ (resp. $P_{n,(i)}^+(\A^{p,\infty})$) which contains the image of $U$ (resp. $U^p$), if $\Delta_0$ is a smooth admissible cone decomposition of $X_*(S_{n,(i),U'}^{+})_\R^{\succ 0}$ (resp. $X_*(\cS_{n,(i),(U')^p(N_1)}^{\ord,+})_\R^{\succ 0}$), and if $\Sigma_0$ is a compatible smooth admissible cone decomposition of $X_*(S_{n,(i),U}^{(m),+})_\R^{\succ 0}$ (resp. $X_*(\cS_{n,(i),U^p(N_1)}^{(m),\ord,+})_\R^{\succ 0}$),  then map
\[ T^{(m),+}_{n,(i),U} \lra T^{+}_{n,(i),U'} \]
(resp.
\[ \cT^{(m),\ord,+}_{n,(i),U^p(N_1,N_2)} \lra \cT^{\ord,+}_{n,(i),(U')^p(N_1,N_2)}) \]
extends to a $L^{(m)}_{n,(i),\lin}(\Q)$-equivariant (resp. $L^{(m)}_{n,(i),\lin}(\Z_{(p)})$-equivariant) log smooth map
\[ (T^{(m),+}_{n,(i),U,\tSigma_0},\cM_{\tSigma_0}) \lra (T^{+}_{n,(i),U',\tDelta_0},\cM_{\tDelta_0}) \]
(resp.
\[ (\cT^{(m),\ord,+}_{n,(i),U^p(N_1,N_2),\tSigma_0},\cM_{\tSigma_0}) \lra (\cT^{\ord,+}_{n,(i),(U')^p(N_1,N_2),\tDelta_0},\cM_{\tDelta_0})). \]
This gives rise to a $P^{(m)}_{n,(i)}(\A^\infty)$-equivariant (resp. $P^{(m),+}_{n,(i)}(\A^{\infty})^\ord$-equivariant) map of systems of log schemes
\[ \{ (T^{(m),+}_{n,(i),U,\tSigma_0},\cM_{\tSigma_0}) \} \lra \{ (T^{+}_{n,(i),U',\tDelta_0},\cM_{\tDelta_0}) \} \]
(resp.
\[ \{ (\cT^{(m),\ord,+}_{n,(i),U^p(N_1,N_2),\tSigma_0},\cM_{\tSigma_0}) \} \lra \{ (\cT^{\ord,+}_{n,(i),(U')^p(N'_1,N_2'),\tDelta_0},\cM_{\tDelta_0}) \} ).\]
These maps are compatible with the embeddings
\[ ( \cT^{(m),\ord,+}_{n,(i),U^p(N_1,N_2),\tSigma_0} \times \Spec \Q,\cM_{\tSigma_0}) \into ( T^{(m),+}_{n,(i),U^p(N_1,N_2),\tSigma_0},\cM_{\tSigma_0}) \]
and
\[ ( \cT^{\ord,+}_{n,(i),U^p(N_1,N_2),\tDelta_0} \times \Spec \Q,\cM_{\tDelta_0}) \into ( T^{+}_{n,(i),U^p(N_1,N_2),\tDelta_0},\cM_{\tDelta_0}). \]

\newpage \subsection{Completions.}

If $\Sigma_0$ denotes a smooth admissible cone decomposition of $X_*(S_{n,(i),U}^{(m),+})_\R^{\succ 0}$ (resp. $X_*(\cS_{n,(i),U^p(N_1,N_2))}^{(m),\ord,+})_\R^{\succ 0}$), then the associated log formal scheme $(T^{(m),+,\wedge}_{n,(i),U,\Sigma_0},\cM_{\Sigma_0}^\wedge)$ (resp. $(\cT^{(m),\ord,+,\wedge}_{n,(i),U^p(N_1,N_2),\Sigma_0},\cM_{\Sigma_0}^\wedge)$) inherits a left action of the group $L_{n,(i),\lin}^{(m)}(\Q)$ (resp. $L_{n,(i),\lin}^{(m)}(\Z_{(p)})$). If $g \in P^{(m),+}_{n,(i)}(\A^\infty)$ (resp. $P^{(m),+}_{n,(i)}(\A^{\infty})^\ord$) and if $\Sigma_0$ is compatible with $\Sigma'_0$ with respect to $g$, then there is an induced
 $L^{(m)}_{n,(i),\lin}(\Q)$-equivariant (resp. $L^{(m)}_{n,(i),\lin}(\Z_{(p)})$-equivariant) map
\[ g: (T^{(m),+,\wedge}_{n,(i),U,\Sigma_0},\cM_{\Sigma_0}^\wedge) \lra (T^{(m),+,\wedge}_{n,(i),U',\Sigma_0'},\cM_{\Sigma_0'}^\wedge) \]
(resp.
\[ g: (\cT^{(m),\ord,+,\wedge}_{n,(i),U^p(N_1,N_2),\Sigma_0}, \cM_{\Sigma_0}^\wedge) \lra (\cT^{(m),\ord,+,\wedge}_{n,(i),(U^p)'(N'_1,N_2'),\Sigma_0'},\cM_{\Sigma_0'}^\wedge)). \]
This makes $\{ (T^{(m),+,\wedge}_{n,(i),U,\Sigma_0},\cM_{\Sigma_0}^\wedge) \}$ (resp.$ \{ (\cT^{(m),\ord,+,\wedge}_{n,(i),U^p(N_1,N_2),\Sigma_0},\cM_{\Sigma_0}^\wedge) \}$) a system of log formal schemes with $P^{(m),+}_{n,(i)}(\A^\infty)$-action (resp. $P^{(m),+}_{n,(i)}(\A^{\infty})^\ord$-action). 

Similarly the schemes $\partial_{\Sigma_0} T^{(m),+}_{n,(i),U}$ (resp. $\partial_{\Sigma_0} \cT^{(m),\ord,+}_{n,(i),U^p(N_1,N_2)}$) inherit a left action of the group $L_{n,(i),\lin}^{(m)}(\Q)$ (resp. $L_{n,(i),\lin}^{(m)}(\Z_{(p)})$). If $g \in P^{(m),+}_{n,(i)}(\A^\infty)$ (resp. $P^{(m),+}_{n,(i)}(\A^{\infty})^\ord$) and if $\Sigma_0$ is compatible with $\Sigma'_0$ with respect to $g$, then there is an induced
 $L^{(m)}_{n,(i),\lin}(\Q)$-equivariant (resp. $L^{(m)}_{n,(i),\lin}(\Z_{(p)})$-equivariant) map
\[ g: \partial_{\Sigma_0} T^{(m),+}_{n,(i),U} \lra \partial_{\Sigma_0'} T^{(m),+}_{n,(i),U'} \]
(resp.
\[ g: \partial_{\Sigma_0} \cT^{(m),\ord,+}_{n,(i),U^p(N_1,N_2)} \lra \partial_{\Sigma_0'} \cT^{(m),\ord,+}_{n,(i),(U^p)'(N_1',N_2')}). \]
This makes $\{ \partial_{\Sigma_0} T^{(m),+}_{n,(i),U} \}$ (resp.$ \{ \partial_{\Sigma_0} \cT^{(m),\ord,+}_{n,(i),U^p(N_1,N_2)} \}$) a system of log formal schemes with $P^{(m),+}_{n,(i)}(\A^\infty)$-action (resp. $P^{(m),+}_{n,(i)}(\A^{\infty})^\ord$-action). 

If $U'$ (resp. $(U^p)'$) is a neat open compact subgroup of $P_{n,(i)}^+(\A^\infty)$ (resp. $P_{n,(i)}^+(\A^{p,\infty})$) which contains the image of $U$ (resp. $U^p$), if $\Delta_0$ is a smooth admissible cone decomposition of $X_*(S_{n,(i),U'}^{+})_\R^{\succ 0}$ (resp. $X_*(\cS_{n,(i),(U^p)'(N_1)}^{\ord,+})_\R^{\succ 0}$), and if $\Sigma_0$ is a compatible smooth admissible cone decomposition of $X_*(S_{n,(i),U}^{(m),+})_\R^{\succ 0}$ (resp. $X_*(\cS_{n,(i),U^p(N_1)}^{(m),\ord,+})_\R^{\succ 0}$),  then there are induced
$L^{(m)}_{n,(i),\lin}(\Q)$-equivariant (resp. $L^{(m)}_{n,(i),\lin}(\Z_{(p)})$-equivariant) maps
\[ (T^{(m),+,\wedge}_{n,(i),U,\Sigma_0},\cM_{\Sigma_0}^\wedge) \lra (T^{+,\wedge}_{n,(i),U',\Delta_0},\cM_{\Delta_0}^\wedge) \]
and
\[ \partial_{\Sigma_0} T^{(m),+}_{n,(i),U} \lra \partial_{\Delta_0} T^{+}_{n,(i),U'} \]
(resp.
\[ (\cT^{(m),\ord,+,\wedge}_{n,(i),U^p(N_1,N_2),\Sigma_0},\cM_{\Sigma_0}^\wedge) \lra (\cT^{\ord,+,\wedge}_{n,(i),(U^p)'(N_1,N_2),\Delta_0},\cM_{\Delta_0}^\wedge) \]
and
\[ \partial_{\Sigma_0} \cT^{(m),\ord,+}_{n,(i),U^p(N_1,N_2)} \lra \partial_{\Delta_0} \cT^{\ord,+}_{n,(i),(U^p)'(N_1,N_2)} ). \]
This gives rise to $P^{(m)}_{n,(i)}(\A^\infty)$-equivariant (resp. $P^{(m),+}_{n,(i)}(\A^{\infty})^\ord$-equivariant) maps of systems of log formal schemes
\[ \{ (T^{(m),+,\wedge}_{n,(i),U,\Sigma_0},\cM_{\Sigma_0}^\wedge) \} \lra \{ (T^{+,\wedge}_{n,(i),U',\Delta_0},\cM_{\Delta_0}^\wedge) \} \]
and of systems of schemes
\[ \{ \partial_{\Sigma_0} T^{(m),+}_{n,(i),U} \} \lra \{ \partial_{\Delta_0} T^{+}_{n,(i),U'} \} \]
(resp.
\[ \{ (\cT^{(m),\ord,+,\wedge}_{n,(i),U^p(N_1,N_2),\Sigma_0},\cM_{\Sigma_0}^\wedge) \} \lra \{ (\cT^{\ord,+,\wedge}_{n,(i),(U')^p(N'_1,N_2'),\Delta_0},\cM_{\Delta_0}^\wedge) \} \]
and
\[ \{ \partial_{\Sigma_0} \cT^{(m),\ord,+}_{n,(i),U^p(N_1,N_2)}\} \lra \{ \partial_{\Delta_0} \cT^{\ord,+}_{n,(i),(U^p)'(N_1,N_2)} \} ).\]

If $\sigma \in \Sigma_0$ and if $1 \neq \gamma \in L_{n,(i),\lin}^{(m)}(\Q)$ (resp. $L_{n,(i),\lin}^{(m)}(\Z_{(p)})$) then $\sigma \cap \gamma \sigma \not\in \Sigma_0$.
Thus
\[  (T^{(m),+,\wedge}_{n,(i),U,\Sigma_0})_\sigma \cap (T^{(m),+,\wedge}_{n,(i),U,\Sigma_0})_{\gamma \sigma} = \emptyset \]
and
\[  (\partial_{\Sigma_0} T^{(m),+}_{n,(i),U})_\sigma \cap (\partial_{\Sigma_0} T^{(m),+}_{n,(i),U})_{\gamma \sigma} = \emptyset \]
(resp.
\[ (\cT^{(m),\ord,+,\wedge}_{n,(i),U^p(N_1,N_2),\Sigma_0})_\sigma \cap (\cT^{(m),\ord,+,\wedge}_{n,(i),U^p(N_1,N_2),\Sigma_0})_{\gamma \sigma} = \emptyset \]
and
\[  (\partial_{\Sigma_0} \cT^{(m),\ord,+}_{n,(i),U^p(N_1,N_2)})_\sigma \cap (\partial_{\Sigma_0} \cT^{(m),\ord,+}_{n,(i),U^p(N_1,N_2)})_{\gamma \sigma} = \emptyset    ). \]
It follows we can form log formal schemes
\[ (T^{(m),\natural,\wedge}_{n,(i),U,\Sigma_0},\cM_{\Sigma_0}^\wedge) = L_{n,(i),\lin}^{(m)}(\Q) \backslash (T^{(m),+,\wedge}_{n,(i),U,\Sigma_0},\cM_{\Sigma_0}^\wedge) \]
(resp.
\[ (\cT^{(m),\ord,\natural,\wedge}_{n,(i),U^p(N_1,N_2),\Sigma_0},\cM_{\Sigma_0}^\wedge) = L_{n,(i),\lin}^{(m)}(\Z_{(p)}) \backslash (\cT^{(m),\ord,+,\wedge}_{n,(i),U^p(N_1,N_2),\Sigma_0},\cM_{\Sigma_0}^\wedge) )\]
and
\[ (T^{(m),\natural+,\wedge}_{n,(i),U,\Sigma_0},\cM_{\Sigma_0}^\wedge) = \Hom_F(F^m,F^i) \backslash (T^{(m),+,\wedge}_{n,(i),U,\Sigma_0},\cM_{\Sigma_0}^\wedge) \]
(resp.
\[ (\cT^{(m),\ord,\natural+,\wedge}_{n,(i),U^p(N_1,N_2),\Sigma_0},\cM_{\Sigma_0}^\wedge) = \Hom_{\cO_{F,(p)}}(\cO_{F,(p)}^m,\cO_{F,(p)}^i) \backslash (\cT^{(m),\ord,+,\wedge}_{n,(i),U^p(N_1,N_2),\Sigma_0},\cM_{\Sigma_0}^\wedge) ).\]
We can also form schemes
\[ \partial_{\Sigma_0} T^{(m),\natural}_{n,(i),U} = L_{n,(i),\lin}^{(m)}(\Q) \backslash \partial_{\Sigma_0} T^{(m),+}_{n,(i),U} \]
(resp.
\[ \partial_{\Sigma_0} \cT^{(m),\ord,\natural}_{n,(i),U^p(N_1,N_2)} = L_{n,(i),\lin}^{(m)}(\Z_{(p)}) \backslash \partial_{\Sigma_0} \cT^{(m),\ord,+}_{n,(i),U^p(N_1,N_2)} ).\]
The quotient maps
\[ (T^{(m),+,\wedge}_{n,(i),U,\Sigma_0},\cM_{\Sigma_0}^\wedge) \onto (T^{(m),\natural+,\wedge}_{n,(i),U,\Sigma_0},\cM_{\Sigma_0}^\wedge) \onto (T^{(m),\natural,\wedge}_{n,(i),U,\Sigma_0},\cM_{\Sigma_0}^\wedge) \]
and
\[ \partial_{\Sigma_0} T^{(m),+}_{n,(i),U} \onto \partial_{\Sigma_0} T^{(m),\natural}_{n,(i),U} \]
(resp. 
\[ (\cT^{(m),\ord,+,\wedge}_{n,(i),U^p(N_1,N_2),\Sigma_0},\cM_{\Sigma_0}^\wedge) \onto (\cT^{(m),\ord,\natural+,\wedge}_{n,(i),U^p(N_1,N_2),\Sigma_0},\cM_{\Sigma_0}^\wedge) \onto (\cT^{(m),\ord,\natural,\wedge}_{n,(i),U^p(N_1,N_2),\Sigma_0},\cM_{\Sigma_0}^\wedge) \]
and
\[    \partial_{\Sigma_0} \cT^{(m),\ord,+}_{n,(i),U^p(N_1,N_2)} \onto \partial_{\Sigma_0} \cT^{(m),\ord,\natural}_{n,(i),U^p(N_1,N_2)}     ) \]
are Zariski locally isomorphisms. The log formal scheme $(T^{(m),\natural+,\wedge}_{n,(i),U,\Sigma_0},\cM_{\Sigma_0}^\wedge)$ (resp. $(\cT^{(m),\ord,\natural+,\wedge}_{n,(i),U^p(N_1,N_2),\Sigma_0},\cM_{\Sigma_0}^\wedge)$) inherits an action of $L_{n,(i),\lin}(\Q)$ (resp. $L_{n,(i),\lin}(\Z_{(p)})$).

If $g \in P^{(m),+}_{n,(i)}(\A^\infty)$ (resp. $P^{(m),+}_{n,(i)}(\A^{\infty})^\ord$) and if $\Sigma_0$ is compatible with $\Sigma'_0$ with respect to $g$ then there are induced maps
\[ g: (T^{(m),\natural,\wedge}_{n,(i),U,\Sigma_0},\cM_{\Sigma_0}^\wedge) \lra (T^{(m),\natural,\wedge}_{n,(i),U',\Sigma_0'},\cM_{\Sigma_0'}^\wedge) \]
(resp.
\[ g: (\cT^{(m),\ord,\natural,\wedge}_{n,(i),U^p(N_1,N_2),\Sigma_0}, \cM_{\Sigma_0}^\wedge) \lra (\cT^{(m),\ord,\natural,\wedge}_{n,(i),(U^p)'(N'_1,N_2'),\Sigma_0'},\cM_{\Sigma_0'}^\wedge)) \]
and
\[ g: (T^{(m),\natural+,\wedge}_{n,(i),U,\Sigma_0},\cM_{\Sigma_0}^\wedge) \lra (T^{(m),\natural+,\wedge}_{n,(i),U',\Sigma_0'},\cM_{\Sigma_0'}^\wedge) \]
(resp.
\[ g: (\cT^{(m),\ord,\natural+,\wedge}_{n,(i),U^p(N_1,N_2),\Sigma_0}, \cM_{\Sigma_0}^\wedge) \lra (\cT^{(m),\ord,\natural+,\wedge}_{n,(i),(U^p)'(N'_1,N_2'),\Sigma_0'},\cM_{\Sigma_0'}^\wedge)) \]
and
\[ g: \partial_{\Sigma_0} T^{(m),\natural}_{n,(i),U} \lra \partial_{\Sigma_0'} T^{(m),\natural}_{n,(i),U'} \]
(resp.
\[ g: \partial_{\Sigma_0} \cT^{(m),\ord,\natural}_{n,(i),U^p(N_1,N_2)} \lra \partial_{\Sigma_0'} \cT^{(m),\ord,\natural}_{n,(i),(U^p)'(N_1',N_2')}  ). \]
This makes the collections $\{ (T^{(m),\natural,\wedge}_{n,(i),U,\Sigma_0},\cM_{\Sigma_0}^\wedge) \}$ (resp. $\{ (\cT^{(m),\ord,\natural,\wedge}_{n,(i),U^p(N_1,N_2),\Sigma_0},\cM_{\Sigma_0}^\wedge)\}$) and $\{ (T^{(m),\natural+,\wedge}_{n,(i),U,\Sigma_0},\cM_{\Sigma_0}^\wedge) \}$ (resp. $ \{ (\cT^{(m),\ord,\natural+,\wedge}_{n,(i),U^p(N_1,N_2),\Sigma_0},\cM_{\Sigma_0}^\wedge) \}$)
systems of log formal schemes with $P^{(m),+}_{n,(i)}(\A^\infty)$-action (resp. $P^{(m),+}_{n,(i)}(\A^{\infty})^\ord$-action). It also makes the collections $\{ \partial_{\Sigma_0} T^{(m),\natural}_{n,(i),U} \}$ (resp. $\{ \partial_{\Sigma_0} \cT^{(m),\ord,\natural}_{n,(i),U^p(N_1,N_2)} \}$)
systems of schemes with $P^{(m),+}_{n,(i)}(\A^\infty)$-action (resp. $P^{(m),+}_{n,(i)}(\A^{\infty})^\ord$-action).

If $U'$ (resp. $(U^p)'$) is a neat open compact subgroup of $P_{n,(i)}^+(\A^\infty)$ (resp. $P_{n,(i)}^+(\A^{p,\infty})$) which contains the image of $U$ (resp. $U^p$), if $\Delta_0$ is a smooth admissible cone decomposition of $X_*(S_{n,(i),U'}^{+})_\R^{\succ 0}$ (resp. $X_*(\cS_{n,(i),(U^p)'(N_1)}^{\ord,+})_\R^{\succ 0}$), and if $\Sigma_0$ is a compatible smooth admissible cone decomposition of $X_*(S_{n,(i),U}^{(m),+})_\R^{\succ 0}$ (resp. $X_*(\cS_{n,(i),U^p(N_1)}^{(m),\ord,+})_\R^{\succ 0}$),  then there are induced
 maps
\[ (T^{(m),\natural,\wedge}_{n,(i),U,\Sigma_0},\cM_{\Sigma_0}^\wedge) \lra (T^{\natural,\wedge}_{n,(i),U',\Delta_0},\cM_{\Delta_0}^\wedge) \]
(resp.
\[ (\cT^{(m),\ord,\natural,\wedge}_{n,(i),U^p(N_1,N_2),\Sigma_0},\cM_{\Sigma_0}^\wedge) \lra (\cT^{\ord,\natural,\wedge}_{n,(i),(U^p)'(N_1,N_2),\Delta_0},\cM_{\Delta_0}^\wedge)) \]
and
\[ (T^{(m),\natural+,\wedge}_{n,(i),U,\Sigma_0},\cM_{\Sigma_0}^\wedge) \lra (T^{+,\wedge}_{n,(i),U',\Delta_0},\cM_{\Delta_0}^\wedge) \]
(resp.
\[ (\cT^{(m),\ord,\natural+,\wedge}_{n,(i),U^p(N_1,N_2),\Sigma_0},\cM_{\Sigma_0}^\wedge) \lra (\cT^{\ord,+,\wedge}_{n,(i),(U^p)'(N_1,N_2),\Delta_0},\cM_{\Delta_0}^\wedge)) \]
and
\[ \partial_{\Sigma_0} T^{(m),\natural}_{n,(i),U} \lra \partial_{\Delta_0} T^{\natural}_{n,(i),U'} \]
(resp.
\[ \partial_{\Sigma_0} \cT^{(m),\ord,\natural}_{n,(i),U^p(N_1,N_2)} \lra \partial_{\Delta_0} \cT^{\ord,\natural}_{n,(i),(U^p)'(N_1,N_2)}  ). \]
These give rise to $P^{(m)}_{n,(i)}(\A^\infty)$-equivariant (resp. $P^{(m),+}_{n,(i)}(\A^{\infty})^\ord$-equivariant) maps of systems of log formal schemes
\[ \{ (T^{(m),\natural,\wedge}_{n,(i),U,\Sigma_0},\cM_{\Sigma_0}^\wedge) \} \lra \{ (T^{\natural,\wedge}_{n,(i),U',\Delta_0},\cM_{\Delta_0}^\wedge) \} \]
(resp.
\[ \{ (\cT^{(m),\ord,\natural,\wedge}_{n,(i),U^p(N_1,N_2),\Sigma_0},\cM_{\Sigma_0}^\wedge) \} \lra \{ (\cT^{\ord,\natural,\wedge}_{n,(i),(U^p)'(N'_1,N_2'),\Delta_0},\cM_{\Delta_0}^\wedge) \} )\]
and
\[ \{ (T^{(m),\natural+,\wedge}_{n,(i),U,\Sigma_0},\cM_{\Sigma_0}^\wedge) \} \lra \{ (T^{+,\wedge}_{n,(i),U',\Delta_0},\cM_{\Delta_0}^\wedge) \} \]
(resp.
\[ \{ (\cT^{(m),\ord,\natural+,\wedge}_{n,(i),U^p(N_1,N_2),\Sigma_0},\cM_{\Sigma_0}^\wedge) \} \lra \{ (\cT^{\ord,+,\wedge}_{n,(i),(U^p)'(N'_1,N_2'),\Delta_0},\cM_{\Delta_0}^\wedge) \} ).\]
They also give rise to a $P^{(m)}_{n,(i)}(\A^\infty)$-equivariant (resp. $P^{(m),+}_{n,(i)}(\A^{\infty})^\ord$-equivariant) map of systems of schemes
\[ \{ \partial_{\Sigma_0} T^{(m),\natural}_{n,(i),U} \} \lra \{ \partial_{\Delta_0} T^{\natural}_{n,(i),U'} \} \]
(resp.
\[ \{  \partial_{\Sigma_0} \cT^{(m),\ord,\natural}_{n,(i),U^p(N_1,N_2)} \} \lra \{  \partial_{\Delta_0} \cT^{\ord,\natural}_{n,(i),(U^p)'(N_1,N_2)} \} ).\]

We will write
\[ \partial_{\Sigma_0} \barT^{(m),\natural}_{n,(i),U^p(N)} = \partial_{\Sigma_0} \cT^{(m),\ord,\natural}_{n,(i),U^p(N_1,N_2)} \times \Spec \F_p. \]
It is independent of $N_2$. 

We also get a commutative diagram
\[ \begin{array}{ccc} T^{(m),+,\wedge}_{n,(i),U,\Sigma_0} && \\ \da && \\
T^{(m),\natural+,\wedge}_{n,(i),U,\Sigma_0} & \lra & T^{+,\wedge}_{n,(i),U',\Delta_0} \\ \da && \da \\
T^{(m),\natural,\wedge}_{n,(i),U,\Sigma_0} & \lra & T^{\natural,\wedge}_{n,(i),U',\Delta_0} \\ \da && \da \\
X^{(m),\natural}_{n,(i),U} &=& X^{\natural}_{n,(i),U'}  \\ \da && \da \\ Y^{(m),\natural}_{n,(i),U} &=& Y^{\natural}_{n,(i),U'} \end{array} \]
(resp.
\[ \begin{array}{ccc} \cT^{(m),\ord,+,\wedge}_{n,(i),U^p(N_1,N_2),\Sigma_0} && \\ \da && \\
\cT^{(m),\ord,\natural+,\wedge}_{n,(i),U^p(N_1,N_2),\Sigma_0} & \lra & \cT^{\ord,+,\wedge}_{n,(i),(U^p)'(N_1,N_2),\Delta_0} \\ \da && \da \\
\cT^{(m),\ord,\natural,\wedge}_{n,(i),U^p(N_1,N_2),\Sigma_0} & \lra & \cT^{\natural,\wedge}_{n,(i),(U^p)'(N_1,N_2),\Delta_0}
\\ \da && \da \\  \cX^{(m),\ord,\natural}_{n,(i),U^p(N_1,N_2)} &=& \cX^{\ord,\natural}_{n,(i),(U')^p(N_1,N_2)}  \\ \da && \da \\ \cY^{(m),\ord,\natural}_{n,(i),U^p(N_1,N_2)} &=& \cY^{\ord,\natural}_{n,(i),(U')^p(N_1,N_2)}). \end{array} \]

We will let $\cI^{(m),+,\wedge}_{\partial, n,(i),U,\Sigma_0}$ denote the formal completion of the ideal sheaf defining
\[ \partial T^{(m),+}_{n,(i),U,\tSigma_0} \subset T^{(m),+}_{n,(i),U,\tSigma_0}. \]
We will let $\cI^{(m),\natural+,\wedge}_{\partial, n,(i),U,\Sigma_0}$ denote its quotient by $\Hom_F(F^m,F^i)$ and $\cI^{(m),\natural,\wedge}_{\partial, n,(i),U,\Sigma_0}$ denote its quotient by $L_{n,(i),\lin}^{(m)}(\Q)$. Similarly we will let $\cI^{(m),\ord,+,\wedge}_{\partial, n,(i),U^p(N_1,N_2),\Sigma_0}$ denote the formal completion of the  ideal sheaf defining
\[ \partial \cT^{(m),\ord,+}_{n,(i),U^p(N_1,N_2),\tSigma_0} \subset \cT^{(m),\ord,+}_{n,(i),U^p(N_1,N_2),\tSigma_0}. \]
We will let $\cI^{(m),\ord,\natural+,\wedge}_{\partial, n,(i),U^p(N_1,N_2),\Sigma_0}$ denote its quotient by $\Hom_{\cO_{F,(p)}}(\cO_{F,(p)}^m,\cO_{F,(p)}^i)$ and $\cI^{(m),\ord,\natural,\wedge}_{\partial, n,(i),U^p(N_1,N_2),\Sigma_0}$ denote its quotient by $L_{n,(i),\lin}^{(m)}(\Z_{(p)})$. 

There are $P^{(m),+}_{n,(i)}(\A^\infty)^\ord$ and $L_{n,(i),\lin}(\Z_{(p)})$ equivariant maps
\[ \cT^{(m),\ord,+,\wedge}_{n,(i),U^p(N_1,N_2),\Sigma_0^\ord} \times \Spf \Q \into T^{(m),+,\wedge}_{n,(i),U^p(N_1,N_2),\Sigma_0} , \]
if $\Sigma_0^\ord$ and $\Sigma_0$ correspond under the bijection of section
\ref{partcomp}. These embeddings are compatible with the maps
\[ \cT^{(m),\ord,+,\wedge}_{n,(i),U^p(N_1,N_2),\Sigma_0^\ord}  \lra \cT^{\ord,+,\wedge}_{n,(i),U^p(N_1,N_2),\Delta_0^\ord} \]
and
\[ T^{(m),+,\wedge}_{n,(i),U^p(N_1,N_2),\Sigma_0} \lra T^{(m),+,\wedge}_{n,(i),U^p(N_1,N_2),\Delta_0}. \]
Moreover they are also compatible with the log structures and with the sheaves $\cI^{(m),\ord,+,\wedge}_{\partial, n,(i),U^p(N_1,N_2),\Sigma_0^\ord}$
and $\cI^{(m),+,\wedge}_{\partial, n,(i),U^p(N_1,N_2),\Sigma_0}$. They induce isomorphisms
\[ \cT^{(m),\ord,\natural,\wedge}_{n,(i),U^p(N_1,N_2),\Sigma_0^\ord} \times \Spf \Q \liso T^{(m),\natural,\wedge}_{n,(i),U^p(N_1,N_2),\Sigma_0} . \]

\begin{lem}\label{ptctp} Suppose that $R_0$ is an irreducible noetherian $\Q$-algebra (resp. $\Z_{(p)}$-algebra) with the discrete topology. Suppose also that $U \supset U'$ (resp. $U^p \supset (U^p)'$) are neat open compact subgroups of $P^{(m),+}_{n,(i)}(\A^\infty)$ (resp. $P^{(m),+}_{n,(i)}(\A^{p,\infty})$), that $N_2' \geq N_1' \geq 0$ and $N_2 \geq N_1 \geq 0$ are integers with $N_2' \geq N_2$ and $N_1' \geq N_1$, and that $\Sigma_0$ and $\Sigma'_0$ are compatible smooth admissible cone decompositions for $X_*(S_{n,(i),U}^{(m),+})_\R^{\succ 0}$ and $X_*(S_{n,(i),U'}^{(m),+})_\R^{\succ 0}$ (resp. $X_*(\cS_{n,(i),U^p(N_1,.N_2)}^{(m),\ord,+})_\R^{\succ 0}$ and $X_*(\cS_{n,(i),(U^p)'(N_1',N_2')}^{(m),\ord,+})_\R^{\succ 0}$). Let $\pi_{(U',\Sigma'),(U,\Sigma)}$ (resp. $\pi_{((U^p)'(N_1',N_2'),\Sigma'),(U^p(N_1,N_2),\Sigma)}$) denote the map 
\[ 1_*: T^{(m),\natural,\wedge}_{n,(i),U',\Sigma_0'} \ra T^{(m),\natural,\wedge}_{n,(i),U,\Sigma_0} \]
(resp.
\[ 1_*: \cT^{(m),\ord,\natural,\wedge}_{n,(i),(U^p)'(N_1',N_2'),\Sigma_0'} \ra \cT^{(m),\ord,\natural,\wedge}_{n,(i),U^p(N_1,N_2),\Sigma_0}). \]
 
\begin{enumerate}
\item If $i>0$ then
\[ R^i\pi_{(U',\Sigma'),(U,\Sigma),*} (\cI^{(m),\natural,\wedge}_{\partial, n,(i),U',\Sigma_0'} \hatotimes R_0)= R^i\pi_{(U',\Sigma'),(U,\Sigma),*} \cO_{ T^{(m),\natural,\wedge}_{n,(i),U',\Sigma_0'} \times \Spf R_0}=(0) \]
(resp.
\[ R^i\pi_{((U^p)'(N_1',N_2'),\Sigma'),(U^p(N_1,N_2),\Sigma),*} (\cI^{(m),\ord,\natural,\wedge}_{\partial, n,(i),(U^p)'(N_1',N_2') ,\Sigma_0'} \hatotimes R_0) = (0) \]
and
\[ R^i\pi_{((U^p)'(N_1',N_2'),\Sigma'),(U^p(N_1,N_2),\Sigma),*} \cO_{ \cT^{(m),\ord,\natural,\wedge}_{n,(i),(U^p)'(N_1',N_2'),\Sigma_0'} \times \Spf R_0}=(0) ).\]

\item Suppose further that $U'$ (resp. $(U^p)'$) is a normal subgroup of $U$ (resp. $U^p$) and that $\Sigma_0'$ is $U$-invariant (resp. $U^p(N_1,N_2)$-invariant). Then the natural maps
\[ \cO_{T^{(m),\natural,\wedge}_{n,(i),U,\Sigma_0} \times \Spf R_0} \lra ( \pi_{(U',\Sigma'),(U,\Sigma),*} \cO_{T^{(m),\natural,\wedge}_{n,(i),U',\Sigma_0'} \times \Spf R_0})^{U} \]
(resp.
\[ \begin{array}{l} \cO_{\cT^{(m),\ord,\natural,\wedge}_{n,(i),U^p(N_1,N_2),\Sigma_0} \times \Spf R_0} \lra \\ ( \pi_{(((U^p)'(N_2,N_2),\Sigma'),(U^p(N_1,N_2),\Sigma),*} \cO_{\cT^{(m),\ord,\natural,\wedge}_{n,(i),(U^p)'(N_2,N_2),\Sigma_0'} \times \Spf R_0})^{U^p(N_1,N_2)} )\end{array} \]
and
\[ \cI^{(m),\natural,\wedge}_{\partial ,n,(i),U,\Sigma_0}\hatotimes R_0 \lra ( \pi_{(U',\Sigma'),(U,\Sigma),*} (\cI^{(m),\natural,\wedge}_{\partial ,n,(i),U',\Sigma_0'} \hatotimes R_0) )^{U} \]
(resp. 
\[ \begin{array}{l} \cI^{(m),\ord,\natural,\wedge}_{\partial,n,(i),U^p(N_1,N_2), \Sigma_0} \hatotimes R_0 \lra \\ ( \pi_{((U^p)'(N_2,N_2),\Sigma'),(U^p(N_1,N_2),\Sigma),*} (\cI^{(m),\ord,\natural,\wedge}_{\partial ,n,(i),(U^p)'(N_2,N_2),\Sigma_0'}\hatotimes R_0))^{U^p(N_1,N_2)} )\end{array} \]
are isomorphisms.
\end{enumerate} 

The same statements are true with $\natural$ replaced by $+$ or by $\natural+$.
\end{lem}

\pfbegin
It suffices to treat the case of $+$. 
We treat the case of $T^{(m),+,\wedge}_{n,(i),U',\Sigma_0'} \times \Spec R_0$, the case of $\cT^{(m),\ord,+,\wedge}_{n,(i),(U^p)'(N_1,N_2),\Sigma_0'} \times \Spec R_0$ being exactly similar. 

Let $U''$ denote the open compact subgroup of $P^{(m),+}_{n,(i)}(\A^\infty)$ generated by $U'$ and $U \cap Z(N_{n,(i)}^{(m)})(\A^\infty)$. Then $\Sigma_0$ is a $U''$ admissible smooth cone decomposition of $X_*(S^{(m),+}_{n,(i),U''})_\R^{\succ 0}$. Moreover
\[ T^{(m),+,\wedge}_{n,(i),U'',\tSigma_0}\times \Spf R_0 \lra T^{(m),+,\wedge}_{n,(i),U,\tSigma_0}  \times \Spf R_0\]
is finite etale, and if $U'$ is normal in $U$ then it is Galois with group $U/U''$. Thus we may replace $U$ by $U''$ and reduce to the case that $U$ and $U'$ have the same projection to $(P^{(m),+}_{n,(i)}/Z (N_{n,(i)}^{(m)}))(\A^\infty)$. In this case the result follows from lemma \ref{refine}.
\pfend

Define $\tOmega_{n,(i),U,\Delta_0}^\natural$ on $T^{\natural,\wedge}_{n,(i),U,\Delta_0}$ as the quotient by $L_{n,(i),\lin}(\Q)$ of the pull-back of $\tOmega_{n,(i),U,\Delta_0}^+$ from $A^+_{n,(i),U}$ to $T^{+,\wedge}_{n,(i),U,\Delta_0}$. Similarly define a sheaf $\tOmega_{n,(i),U^p(N_1,N_2),\Delta_0}^{\ord,\natural}$ on $\cT^{\ord,\natural,\wedge}_{n,(i),U^p(N_1,N_2),\Delta_0}$ as the quotient by $L_{n,(i),\lin}(\Z_{(p)})$ of the pull-back of the sheaf $\tOmega_{n,(i),U^p(N_1,N_2),\Delta_0}^{\ord,+}$ from $\cA^{\ord,+}_{n,(i),U^p(N_1,N_2)}$ to $\cT^{\ord,+,\wedge}_{n,(i),U^p(N_1,N_2),\Delta_0}$.

Suppose that $R_0$ is an irreducible noetherian $\Q$-algebra and that $\rho$ is a representation of $R_{n,(n),(i)}$ on a finite, locally free $R_0$-module $W_\rho$. Then we define a locally free sheaf $\cE_{(i),U,\Delta_0,\rho}^\natural$ on $T^{\natural,\wedge}_{n,(i),U,\Delta_0}$ as the quotient by $L_{n,(i),\lin}(\Q)$ of the pull-back of $\cE^+_{(i),U,\rho}$ from $A^+_{n,(i),U}$ to $T^{+,\wedge}_{n,(i),U,\Delta_0}$. Then the system of sheaves $\{ \cE_{(i),U,\Delta_0,\rho}^\natural\}$ over $\{ T^{\natural,\wedge}_{n,(i),U,\Delta_0} \}$ has an action of $P_{n,(i)}^+(\A^\infty)$. If $g \in P_{n,(i)}^+(\A^\infty)$, then the natural map
\[ g^*\cE_{(i),U,\Delta_0,\rho}^\natural \lra \cE_{(i),U',\Delta'_0,\rho}^\natural \]
is an isomorphism. The sheaves $\cE_{(i),U,\Delta_0,\rho}^\natural$ have $P_{n,(i)}^+(\A^\infty)$-invariant filtrations by local direct summands whose graded pieces pull-backed to $T^{+,\wedge}_{n,(i),U,\Delta_0}$ are equivariantly isomorphic to the pull-backs of sheaves of the form $\cE_{(i),U,\rho'}^+$ on $X^+_{n,(i),U}$. 

Similarly in the case of mixed characteristic
suppose that $R_0$ is an irreducible noetherian $\Z_{(p)}$-algebra and that $\rho$ is a representation of $R_{n,(n),(i)}$ on a finite, locally free $R_0$-module $W_\rho$. Then we define a locally free sheaf $\cE_{(i),U^p(N_1,N_2),\Delta_0,\rho}^{\ord,\natural}$ on $\cT^{\ord,\natural,\wedge}_{n,(i),U^p(N_1,N_2),\Delta_0}$ as the quotient by $L_{n,(i),\lin}(\Z_{(p)})$ of the pull-back of $\cE^{\ord,+}_{(i),U^p(N_1,N_2),\rho}$ from $\cA^{\ord,+}_{n,(i),U^p(N_1,N_2)}$ to $\cT^{\ord,+,\wedge}_{n,(i),U^p(N_1,N_2),\Delta_0}$. Then the collection $\{ \cE_{(i),U^p(N_1,N_2),\Delta_0,\rho}^{\ord,\natural}\}$ is a system of sheaves over $\{ \cT^{\ord,\natural,\wedge}_{n,(i),U^p(N_1,N_2),\Delta_0} \}$ has an action of $P_{n,(i)}^+(\A^\infty)^\ord$. If $g \in P_{n,(i)}^+(\A^\infty)^{\ord,\times}$, then the natural map
\[ g^*\cE_{(i),U^p(N_1,N_2),\Delta_0,\rho}^{\ord,\natural} \lra \cE_{(i),(U^p)'(N_1,N_2),\Delta'_0,\rho}^{\ord,\natural} \]
is an isomorphism. The sheaves $\cE_{(i),U^p(N_1,N_2),\Delta_0,\rho}^{\ord,\natural}$ have $P_{n,(i)}^+(\A^\infty)^\ord$-invariant filtrations by local direct summands whose graded pieces pull-backed to the formal scheme $\cT^{\ord,+,\wedge}_{n,(i),U^p(N_1,N_2),\Delta_0}$ are equivariantly isomorphic to the pull-backs of sheaves of the form $\cE^{\ord,+}_{(i),U^p(N_1,N_2),\rho'}$ on $\cX^{\ord,+}_{n,(i),U^p(N_1,N_2)}$.

\begin{cor} Suppose that $R_0$ is an irreducible noetherian $\Q$-algebra (resp. $\Z_{(p)}$-algebra) with the discrete topology. Let $\rho$ be a representation of $R_{n,(n),(i)}$ on a finite, locally free $R_0$-module $W_\rho$.
Suppose also that $U \supset U'$ (resp. $U^p \supset (U^p)'$) are neat open compact subgroups of $P^{(m),+}_{n,(i)}(\A^\infty)$ (resp. $P^{(m),+}_{n,(i)}(\A^{p,\infty})$), that $N_2' \geq N_1' \geq 0$ and $N_2 \geq N_1 \geq 0$ are integers with $N_2' \geq N_2$ and $N_1' \geq N_1$, 
and that $\Sigma_0$ and $\Sigma'_0$ are compatible smooth admissible cone decompositions for $X_*(S_{n,(i),U}^{(m),+})_\R^{\succ 0}$ and $X_*(S_{n,(i),U'}^{(m),+})_\R^{\succ 0}$ (resp. $X_*(\cS_{n,(i),U^p(N_1,.N_2)}^{(m),\ord,+})_\R^{\succ 0}$ and $X_*(\cS_{n,(i),(U^p)'(N_1',N_2')}^{(m),\ord,+})_\R^{\succ 0}$). Let $\pi_{(U',\Sigma_0'),(U,\Sigma_0)}$ (resp. $\pi_{((U^p)'(N_1',N_2'),\Sigma_0'),(U^p(N_1,N_2),\Sigma_0)}$) denote the map 
\[ 1_*: T^{(m),\natural,\wedge}_{n,(i),U',\Sigma_0'} \ra T^{(m),\natural,\wedge}_{n,(i),U,\Sigma_0} \]
(resp.
\[ 1_*: \cT^{(m),\ord,\natural,\wedge}_{n,(i),(U^p)'(N_1',N_2'),\Sigma_0'} \ra \cT^{(m),\ord,\natural,\wedge}_{n,(i),U^p(N_1,N_2),\Sigma_0}). \]
 
\begin{enumerate}
\item If $i>0$ then
\[ R^i\pi_{(U',\Sigma_0'),(U,\Sigma_0),*} (\cI^{(m),\natural,\wedge}_{\partial, n,(i),U',\Sigma_0'} \hatotimes \cE_{(i),U',\Sigma_0',\rho}^\natural)= R^i\pi_{(U',\Sigma_0'),(U,\Sigma_0),*} \cE_{(i),U',\Sigma_0',\rho}^\natural=(0) \]
(resp.
\[ R^i\pi_{((U^p)'(N_1',N_2'),\Sigma_0'),(U^p(N_1,N_2),\Sigma_0),*} (\cI^{(m),\ord,\natural,\wedge}_{\partial, n,(i),(U^p)'(N_1',N_2') ,\Sigma_0'} \hatotimes \cE_{(i),(U^p)'(N_1',N_2'),\Sigma_0',\rho}^{\ord,\natural}) = (0) \]
and
\[ R^i\pi_{((U^p)'(N_1',N_2'),\Sigma_0'),(U^p(N_1,N_2),\Sigma_0),*} \cE_{(i),(U^p)'(N_1',N_2'),\Sigma_0',\rho}^{\ord,\natural}=(0) ).\]

\item Suppose further that $U'$ (resp. $(U^p)'$) is a normal subgroup of $U$ (resp. $U^p$) and that $\Sigma_0'$ is $U$-invariant (resp. $U^p(N_1,N_2)$-invariant). Then the natural maps
\[ \cE_{(i),U,\Sigma_0,\rho}^\natural \lra ( \pi_{(U',\Sigma_0'),(U,\Sigma_0),*} \cE_{(i),U',\Sigma_0',\rho}^\natural)^{U} \]
(resp.
\[ \begin{array}{l} \cE_{(i),U^p(N_1,N_2),\Sigma_0,\rho}^{\ord,\natural} \lra \\ ( \pi_{(((U^p)'(N_2,N_2),\Sigma_0'),(U^p(N_1,N_2),\Sigma_0),*} \cE_{(i),(U^p)'(N_2,N_2),\Sigma_0',\rho}^{\ord,\natural})^{U^p(N_1,N_2)} )\end{array} \]
and
\[ \cI^{(m),\natural,\wedge}_{\partial ,n,(i),U,\Sigma_0}\hatotimes \cE_{(i),U,\Sigma_0,\rho}^\natural \lra ( \pi_{(U',\Sigma_0'),(U,\Sigma_0),*} (\cI^{(m),\natural,\wedge}_{\partial ,n,(i),U',\Sigma_0'} \hatotimes \cE_{(i),U',\Sigma'_0,\rho}^\natural) )^{U} \]
(resp. 
\[ \begin{array}{l} \cI^{(m),\ord,\natural,\wedge}_{\partial,n,(i),U^p(N_1,N_2), \Sigma_0} \hatotimes \cE_{(i),U^p(N_1,N_2),\Sigma_0,\rho}^{\ord,\natural} \lra \\ ( \pi_{((U^p)'(N_2,N_2),\Sigma'),(U^p(N_1,N_2),\Sigma),*} (\cI^{(m),\ord,\natural,\wedge}_{\partial ,n,(i),(U^p)'(N_2,N_2),\Sigma_0'}\hatotimes \cE_{(i),(U^p)'(N_2,N_2),\Sigma_0',\rho}^{\ord,\natural}))^{U^p(N_1,N_2)} )\end{array} \]
are isomorphisms.
\end{enumerate} 
\end{cor}

\begin{lem}\label{rc} Suppose that $U$ is a neat open compact subgroup of $P_{n,(i)}^{(m),+}(\A^\infty)$ and let $U'$ denote the image of $U$ in $P_{n,(i)}^+(\A^\infty)$. Let $\Delta_0$ be a smooth admissible cone decomposition for $X_*(S^+_{n,(i),U'})$ and let $\Sigma_0$ be a compatible smooth admissible cone decomposition for $X_*(S^{(m),+}_{n,(i),U})$. Let $\pi^+=\pi_{(U,\Sigma_0),(U',\Delta_0)}^+$ denote the map
\[ T^{(m),\natural+,\wedge}_{n,(i),U,\Sigma_0} \lra T^{+,\wedge}_{n,(i),U',\Delta_0} \]
and let $\pi^\natural=\pi_{(U,\Sigma_0),(U',\Delta_0)}^\natural$ denote the map
\[ T^{(m),\natural,\wedge}_{n,(i),U,\Sigma_0} \lra T^{\natural,\wedge}_{n,(i),U',\Delta_0}. \]
\begin{enumerate}
\item The maps $\pi_{(U,\Sigma_0),(U',\Delta_0)}^+$ and $\pi_{(U,\Sigma_0),(U',\Delta_0)}^\natural$ are proper.

\item The natural maps
\[ \cO_{T^{+,\wedge}_{n,(i),U',\Delta_0}} \lra \pi_{(U,\Sigma_0),(U',\Delta_0),*}^+ \cO_{T^{(m),\natural+,\wedge}_{n,(i),U,\Sigma_0}} \]
and
\[ \cI^{+,\wedge}_{\partial,n,(i),U',\Delta_0} \lra \pi_{(U,\Sigma_0),(U',\Delta_0),*}^+ \cI^{(m),\natural+,\wedge}_{\partial,n,(i),U,\Sigma_0} \]
and
\[ \cO_{T^{\natural,\wedge}_{n,(i),U',\Delta_0}} \lra \pi_{(U,\Sigma_0),(U',\Delta_0),*}^\natural \cO_{T^{(m),\natural,\wedge}_{n,(i),U,\Sigma_0}} \]
and
\[ \cI^{\natural,\wedge}_{\partial,n,(i),U',\Delta_0} \lra \pi_{(U,\Sigma_0),(U',\Delta_0),*}^\natural \cI^{(m),\natural,\wedge}_{\partial,n,(i),U,\Sigma_0} \]
are isomorphisms.

\item The natural maps
\[ \cI^{+,\wedge}_{\partial,n,(i),U',\Delta_0} \otimes R^j\pi_{(U,\Sigma_0),(U',\Delta_0),*}^+ \cO_{T^{(m),\natural+,\wedge}_{n,(i),U,\Sigma_0}} \lra R^j\pi_{(U,\Sigma_0),(U',\Delta_0),*}^+ \cI^{(m),\natural+,\wedge}_{\partial,n,(i),U,\Sigma_0} \]
and
\[ \cI^{\natural,\wedge}_{\partial,n,(i),U',\Delta_0} \otimes R^j\pi_{(U,\Sigma_0),(U',\Delta_0),*}^\natural \cO_{T^{(m),\natural,\wedge}_{n,(i),U,\Sigma_0}} \lra R^j\pi_{(U,\Sigma_0),(U',\Delta_0),*}^\natural \cI^{(m),\natural,\wedge}_{\partial,n,(i),U,\Sigma_0} \]
are isomorphisms.
\end{enumerate}\end{lem}

\pfbegin
It suffices to treat the $+$ case.

The first part follows from lemma \ref{dproper}. We deduce that all the sheaves mentioned in the remaining parts are coherent. 

Thus, by theorem 4.1.5 of \cite{ega3} (`the theorem on formal functions'), it suffices to prove the remaining assertions after completing at a point of $T^{+,\wedge}_{n,(i),U',\Delta_0}$. The set points where the assertions are true after completing at that point is open. (Again because the sheaves involved are all coherent.) This open set is $S^{+}_{n,(i),U'}$-invariant. (The sheaves in question do not all have $S^{+}_{n,(i),U'}$-actions. However locally on $T^{+,\wedge}_{n,(i),U',\Delta_0}$ they do.) Thus it will do to prove the lemma after completion at $\partial_\sigma T^{+}_{n,(i),U',\tDelta_0}$, for $\sigma \in \Delta_0$ maximal. We will add a subscript $\sigma$ to denote completion along $\partial_\sigma T^{+}_{n,(i),U',\tDelta_0}$.

We write $\tpi$ for the map
\[ T^{(m),+,\wedge}_{n,(i),U,\Sigma_0} \lra T^{+,\wedge}_{n,(i),U',\Delta_0} \]
and factor $\tpi= \pi_2 \circ \pi_1$, where 
\[ \pi_1: T^{(m),+,\wedge}_{n,(i),U,\Sigma_0} \lra T^{+,\wedge}_{n,(i),U',\Delta_0} \times_{A^+_{n,(i),U'}} A^{(m),+}_{n,(i),U} \]
and
\[ \pi_2: T^{+,\wedge}_{n,(i),U',\Delta_0} \times_{A^+_{n,(i),U'}} A^{(m),+}_{n,(i),U} \lra T^{+,\wedge}_{n,(i),U',\Delta_0}. \]
Also write $\pi_3$ for the other projection
\[ \pi_3: T^{+,\wedge}_{n,(i),U',\Delta_0} \times_{A^+_{n,(i),U'}} A^{(m),+}_{n,(i),U} \lra A^{(m),+}_{n,(i),U}. \]

We will first show that
\[ R^j\pi_{1,\sigma,*} \cO_{T^{(m),+,\wedge}_{n,(i),U,\Sigma_0,\sigma}} = \left\{ \begin{array}{ll} (0) & {\rm if}\,\, i>0 \\ \cO_{T^{+,\wedge}_{n,(i),U',\Delta_0,\sigma} \times_{A^+_{n,(i),U'}} A^{(m),+}_{n,(i),U}} & {\rm if}\,\, i=0 \end{array} \right. \]
and
\[ R^j\pi_{1,\sigma,*} \cI^{(m),+,\wedge}_{\partial,n,(i),U,\Sigma_0,\sigma} = \left\{ \begin{array}{ll} (0) & {\rm if}\,\, i>0 \\ \pi_{2,\sigma}^*\cI^{+,\wedge}_{\partial,n,(i),U',\Delta_0,\sigma}  & {\rm if}\,\, i=0. \end{array} \right. \]
As $T^{+,\wedge}_{n,(i),U',\Delta_0,\sigma} \times_{A^+_{n,(i),U'}} A^{(m),+}_{n,(i),U}$ has the same underlying topological space as $A^{(m),+}_{n,(i),U}$, i.e. $\pi_{3,\sigma}$ is a homeomorphism on the underlying topological space, it suffices to show that
\[ R^j(\pi_3 \circ \pi_1)_{\sigma,*} \cO_{T^{(m),+,\wedge}_{n,(i),U,\Sigma_0,\sigma}} = \left\{ \begin{array}{ll} (0) & {\rm if}\,\, i>0 \\ \pi_{3,\sigma,*} \cO_{T^{+,\wedge}_{n,(i),U',\Delta_0,\sigma} \times_{A^+_{n,(i),U'}} A^{(m),+}_{n,(i),U}} & {\rm if}\,\, i=0 \end{array} \right. \]
and
\[ R^j\pi_{1,\sigma,*} \cI^{(m),+,\wedge}_{\partial,n,(i),U,\Sigma_0,\sigma} = \left\{ \begin{array}{ll} (0) & {\rm if}\,\, i>0 \\ \pi_{3,\sigma,*}\pi_{2,\sigma}^*\cI^{+,\wedge}_{\partial,n,(i),U',\Delta_0,\sigma}  & {\rm if}\,\, i=0. \end{array} \right. \]
This would follow from lemma \ref{dctorem3} as long as we can show that, for all $y \in Y^{(m),+}_{n,(i),U}$ with image $y'$ in $Y^{+}_{n,(i),U'}$, we have $|\Sigma_0|^\vee(y)= |\Delta_0|^\vee (y')$
and  $|\Sigma_0|^{\vee,0}(y)  = |\Delta_0|^{\vee, 0} (y')$. Concretely these required equalities 
are
\[ \gC^{(m),\succ 0}(V_{n,(i))})^\vee = \gC^{\succ 0}(V_{n,(i))})^\vee  \]
and
\[  \gC^{(m),\succ 0}(V_{n,(i))})^{\vee,0} = \gC^{\succ 0}(V_{n,(i))})^{\vee,0}. \]
However, we know that
\[  \gC^{>0}(V_{n,(i))})^\vee = \gC^{\succ 0}(V_{n,(i))})^\vee = \gC^{\geq 0}(V_{n,(i))})^\vee. \]
(They all equal the positive semi-definite cone.)
Also
\[ \gC^{>0}(V_{n,(i))})^{\vee,0} = \gC^{\succ 0}(V_{n,(i))})^{\vee,0} = \gC^{\geq 0}(V_{n,(i))})^{\vee,0}. \]
(They all equal the positive definite cone.) Moreover
\[ \begin{array}{ccccccc} \gC^{>0}(V_{n,(i))})^\vee &=& \gC^{(m),>0}(V_{n,(i))})^\vee &\supset & \gC^{(m),\succ 0}(V_{n,(i))})^\vee &\supset &  \gC^{(m),\geq 0}(V_{n,(i))})^\vee \\ &&&&& = &\gC^{\geq 0}(V_{n,(i))})^\vee \end{array} \]
and
\[ \begin{array}{ccccc} \gC^{>0}(V_{n,(i))})^{\vee,0} &=& \gC^{(m),>0}(V_{n,(i))})^{\vee,0} &\supset & \gC^{(m),\succ 0}(V_{n,(i))})^{\vee,0} \\ &\supset & \gC^{(m),\geq 0}(V_{n,(i))})^{\vee,0} &= &\gC^{\geq 0}(V_{n,(i))})^{\vee,0}. \end{array} \]
Thus
\[ \gC^{(m),\succ 0}(V_{n,(i))})^\vee = \gC^{\succ 0}(V_{n,(i))})^\vee \]
and
\[ \gC^{(m),\succ 0}(V_{n,(i))})^{\vee,0} = \gC^{\succ 0}(V_{n,(i))})^{\vee,0}, \]
as desired.

We deduce that
\[ R^j\tpi_{\sigma,*} \cO_{T^{(m),+,\wedge}_{n,(i),U,\Sigma_0,\sigma}} =  (\wedge^j \Hom_F( \Omega^+_{n,(i),U'}, F^m\otimes_\Q \Xi^+_{n,(i),U'})) \otimes_{\cO_{X^{+}_{n,(i),U'}}} \cO_{T^{+,\wedge}_{n,(i),U',\Delta_0,\sigma}}  \]
and
\[ R^j\tpi_{\sigma,*} \cI^{(m),+,\wedge}_{\partial,n,(i),U,\Sigma_0,\sigma} =  (\wedge^j \Hom_F(\Omega^+_{n,(i),U'}, F^m\otimes_\Q \Xi^+_{n,(i),U'})) \otimes_{\cO_{X^{+}_{n,(i),U'}}} \cI^{+,\wedge}_{\partial,n,(i),U',\Delta_0,\sigma} . \]
As $T^{(m),\natural+,\wedge}_{n,(i),U,\Sigma_0,\sigma}$ is the quotient of $T^{(m),+,\wedge}_{n,(i),U,\Sigma_0,\sigma}$ by $\Hom_F(F^m,F^i)$, we obtain spectral sequences
\[ \begin{array}{l} \!\!\!\! H^{j_1}(\Hom_F(F^m,F^i), (\wedge^{j_2} \Hom_F( \Omega^+_{n,(i),U'}, F^m\otimes_\Q \Xi^+_{n,(i),U'}))) \otimes_{\cO_{X^{+}_{n,(i),U'}}} \cO_{T^{+,\wedge}_{n,(i),U',\Delta_0,\sigma}} \\ \Rightarrow R^{j_1+j_2}\pi_*^+ \cO_{T^{(m),\natural+,\wedge}_{n,(i),U,\Sigma_0,\sigma}} \end{array} \]
and
\[ \begin{array}{l}  \!\!\!\! H^{j_1}(\Hom_F(F^m,F^i), (\wedge^{j_2} \Hom_F( \Omega^+_{n,(i),U'}, F^m\otimes_\Q \Xi^+_{n,(i),U'}))) \otimes_{\cO_{X^{+}_{n,(i),U'}}} \!\!\!\cI^{+,\wedge}_{\partial,n,(i),U',\Delta_0,\sigma} \\ \Rightarrow R^{j_1+j_2}\pi_{*}^+ \cI^{(m),\natural+,\wedge}_{\partial,n,(i),U,\Sigma_0,\sigma}. \end{array}\]
These can also be written
\[ \begin{array}{l}  \!\! \Hom(\wedge^{j_1}\Hom_F(F^m,F^i), \wedge^{j_2} \Hom_F(\Omega^+_{n,(i),U'}, F^m\otimes_\Q \Xi^+_{n,(i),U'})) \!\! \otimes_{\cO_{X^{+}_{n,(i),U'}}} \!\!\!\!\!\!\!\! \cO_{T^{+,\wedge}_{n,(i),U',\Delta_0,\sigma}} \\ \Rightarrow R^{j_1+j_2}\pi_{*}^+ \cO_{T^{(m),\natural+,\wedge}_{n,(i),U,\Sigma_0,\sigma}} \end{array} \]
and
\[ \begin{array}{l}  \!\!\! \Hom(\wedge^{j_1}\Hom_F(F^m,F^i), \wedge^{j_2} \Hom_F(\Omega^+_{n,(i),U'}, \! F^m\otimes_\Q \Xi^+_{n,(i),U'})) \!\!\otimes_{\cO_{X^{+}_{n,(i),U'}}} \!\!\!\!\!\!\!\! \cI^{+,\wedge}_{\partial,n,(i),U',\Delta_0,\sigma} \\ \Rightarrow R^{j_1+j_2}\pi_{*}^+ \cI^{(m),\natural+,\wedge}_{\partial,n,(i),U,\Sigma_0,\sigma}. \end{array} \]
The lemma follows (as $\cI^{+,\wedge}_{\partial,n,(i),U',\Delta_0,\sigma}$ is flat over $\cO_{T^{+,\wedge}_{n,(i),U',\Delta_0,\sigma}}$).
\pfend

The following lemma is lemma 1.3.2.79 of \cite{kw2}.
\begin{lem} Suppose that $U$ is a neat open compact subgroup of $P_{n,(i)}^{(m),+}(\A^\infty)$ and let $U'$ denote the image of $U$ in $P_{n,(i)}^+(\A^\infty)$. Let $\Delta_0$ be a smooth admissible cone decomposition for $X_*(S^+_{n,(i),U'})$ and let $\Sigma_0$ be a compatible smooth admissible cone decomposition for $X_*(S^{(m),+}_{n,(i),U})$. There are canonical equivariant isomorphisms
\[ \Hom_F(F^m,\tOmega_{n,(i),U'}^+) \otimes_{\cO_{A^+_{n,(i),U'}}} \cO_{T^{(m),+,\wedge}_{n,(i),U,\Sigma_0}} \liso  \Omega^1_{T^{(m),+,\wedge}_{n,(i),U,\Sigma_0}/T^{+,\wedge}_{n,(i),U',\Delta_0}} (\log \infty). \] 
\end{lem}

We deduce the following lemmas.
\begin{lem}\label{rc2} Suppose that $U$ is a neat open compact subgroup of $P_{n,(i)}^{(m),+}(\A^\infty)$ and let $U'$ denote the image of $U$ in $P_{n,(i)}^+(\A^\infty)$. Let $\Delta_0$ be a smooth admissible cone decomposition for $X_*(S^+_{n,(i),U'})$ and let $\Sigma_0$ be a compatible smooth admissible cone decomposition for $X_*(S^{(m),+}_{n,(i),U})$. Let $\pi^+=\pi_{(U,\Sigma_0),(U',\Delta_0)}^+$ denote the map
\[ T^{(m),\natural+,\wedge}_{n,(i),U,\Sigma_0} \lra T^{+,\wedge}_{n,(i),U',\Delta_0} \]
and let $\pi^\natural=\pi_{(U,\Sigma_0),(U',\Delta_0)}^\natural$ denote the map
\[ T^{(m),\natural,\wedge}_{n,(i),U,\Sigma_0} \lra T^{\natural,\wedge}_{n,(i),U',\Delta_0}. \]
\begin{enumerate}

\item $\pi^\natural_*\Omega^1_{T^{(m),\natural,\wedge}_{n,(i),U,\Sigma_0}/T^{\natural,\wedge}_{n,(i),U',\Delta_0}} 
(\log \infty) \cong \Hom_F(F^m,\tOmega_{n,(i),U'}^\natural) $ is locally free of finite rank.

\item The natural map
\[ \begin{array}{l} \pi_{(U,\Sigma_0),(U',\Delta_0)}^{\natural,*}\pi_{(U,\Sigma_0),(U',\Delta_0),*}^\natural \Omega^1_{T^{(m),\natural,\wedge}_{n,(i),U,\Sigma_0}/T^{\natural,\wedge}_{n,(i),U',\Delta_0}} 
(\log \infty) \lra \\ \Omega^1_{T^{(m),\natural,\wedge}_{n,(i),U,\Sigma_0}/T^{\natural,\wedge}_{n,(i),U',\Delta_0}} (\log \infty) \end{array} \]
is an isomorphism. 

\item The natural maps 
\[ \begin{array}{l} (R^{j_1}\pi^\natural_* \cO_{T^{(m),\natural,\wedge}_{n,(i),U,\Sigma_0}}) \otimes (\wedge^{j_2} \pi^\natural_* \Omega^1_{T^{(m),\natural,\wedge}_{n,(i),U,\Sigma_0}/T^{\natural,\wedge}_{n,(i),U',\Delta_0}} 
(\log \infty)) \lra \\ R^{j_1}\pi^\natural_* \Omega^{j_2}_{T^{(m),\natural,\wedge}_{n,(i),U,\Sigma_0}/T^{\natural,\wedge}_{n,(i),U',\Delta_0}} 
(\log \infty) \end{array} \]
and
\[ \begin{array}{l} (R^{j_1}\pi^\natural_* \cO_{T^{(m),\natural,\wedge}_{n,(i),U,\Sigma_0}}) \otimes (\wedge^{j_2} \pi^\natural_* \Omega^1_{T^{(m),\natural,\wedge}_{n,(i),U,\Sigma_0}/T^{\natural,\wedge}_{n,(i),U',\Delta_0}} 
(\log \infty)) \otimes \cI^{\natural,\wedge}_{\partial,n,(i),U',\Delta_0} \lra \\ R^{j_1}\pi^\natural_* (\Omega^{j_2}_{T^{(m),\natural,\wedge}_{n,(i),U,\Sigma_0}/T^{\natural,\wedge}_{n,(i),U',\Delta_0}} (\log \infty) \otimes \cI^{(m),\natural,\wedge}_{\partial,n,(i),U,\sigma_0})
 \end{array} \]
are isomorphisms.

\end{enumerate}\end{lem}

\begin{lem}\label{rc3} Suppose that $U \supset U'$ are neat open compact subgroups of the group $P_{n,(i)}^{(m),+}(\A^\infty)$ and let $V$ and $V'$ denote the images of $U$ and $U'$ in $P_{n,(i)}^+(\A^\infty)$. Let $\Delta_0$ (resp. $\Delta_0'$) be a smooth admissible cone decomposition for $X_*(S^+_{n,(i),V})$ (resp. $X_*(S^+_{n,(i),V'})$) and let $\Sigma_0$ (resp. $\Sigma_0'$) be a compatible smooth admissible cone decomposition for $X_*(S^{(m),+}_{n,(i),U})$ (resp. $X_*(S^{(m),+}_{n,(i),U'})$). Further suppose that $\Sigma_0$ and $\Sigma_0'$ are compatible and that $\Delta_0$ and $\Delta_0'$ are compatible. 
\begin{enumerate}
\item The natural map
\[ \pi_{(U',\Sigma_0'),(U,\Sigma_0)}^* \Omega^1_{T^{(m),\natural,\wedge}_{n,(i),U,\Sigma_0}/T^{\natural,\wedge}_{n,(i),V,\Delta_0}}(\log \infty) \lra \Omega^1_{T^{(m),\natural,\wedge}_{n,(i),U',\Sigma_0'}/T^{\natural,\wedge}_{n,(i),V',\Delta_0'}} (\log \infty)\]
is an isomorphism. 

\item The natural map
\[ \begin{array}{r} \pi_{(V',\Delta_0'),(V,\Delta_0)}^* \pi_{(U,\Sigma_0),(V,\Delta_0),*}     \Omega^1_{T^{(m),\natural,\wedge}_{n,(i),U,\Sigma_0}/T^{\natural,\wedge}_{n,(i),V,\Delta_0}}(\log \infty) \lra \\ \pi_{(U',\Sigma_0'),(V',\Delta_0'),*} \Omega^1_{T^{(m),\natural,\wedge}_{n,(i),U',\Sigma_0'}/T^{\natural,\wedge}_{n,(i),V',\Delta_0'}}(\log \infty) \end{array}\]
is an isomorphism.
\end{enumerate} \end{lem}

Similarly we have the following lemma.
\begin{lem}\label{rcord} Suppose that $U^p$ is a neat open compact subgroup of $P_{n,(i)}^{(m),+}(\A^{p,\infty})$ and let $(U^p)'$ denote the image of $U^p$ in $P_{n,(i)}^+(\A^{p,\infty})$. Also suppose that $N_2 \geq N_1 \geq 0$ are integers. Let $\Delta_0$ be a smooth admissible cone decomposition for $X_*(\cS^{\ord,+}_{n,(i),(U^p)'(N_1,N_2)})$ and let $\Sigma_0$ be a compatible smooth admissible cone decomposition for $X_*(\cS^{(m),\ord,+}_{n,(i),U^p(N_1,N_2)})$. Let 
$\pi^\natural =\pi_{(U^p(N_1,N_2),\Sigma_0),((U^p)'(N_1,N_2),\Delta_0)}^\natural$ denote the map
\[ \cT^{(m),\ord,\natural,\wedge}_{n,(i),U^p(N_1,N_2),\Sigma_0} \lra \cT^{\ord,\natural,\wedge}_{n,(i),(U^p)'(N_1,N_2),\Delta_0}. \]
\begin{enumerate}
\item The map 
$\pi_{(U^p(N_1,N_2),\Sigma_0),((U^p)'(N_1,N_2),\Delta_0)}^\natural$ is proper.

\item The natural maps
\[ \cO_{\cT^{\ord,\natural,\wedge}_{n,(i),(U^p)'(N_1,N_2),\Delta_0}} \lra \pi_{(U^p(N_1,N_2),\Sigma_0),((U^p)'(N_1,N_2),\Delta_0),*}^\natural \cO_{\cT^{(m),\ord,\natural,\wedge}_{n,(i),U^p(N_1,N_2),\Sigma_0}} \]
and
\[ \cI^{\ord,\natural,\wedge}_{\partial,n,(i),(U^p)'(N_1,N_2),\Delta_0} \lra \pi_{(U^p(N_1,N_2),\Sigma_0),((U^p)'(N_1,N_2),\Delta_0),*}^\natural \cI^{(m),\ord,\natural,\wedge}_{\partial,n,(i),U^p(N_1,N_2),\Sigma_0} \]
are isomorphisms.

\item The natural map
\[ \cI^{\ord,\natural,\wedge}_{\partial,n,(i),(U^p)'(N_1,N_2),\Delta_0} \otimes R^j\pi_{*}^\natural \cO_{\cT^{(m),\ord,\natural,\wedge}_{n,(i),U^p(N_1,N_2),\Sigma_0}} \lra R^j\pi_{*}^\natural \cI^{(m),\ord,\natural,\wedge}_{\partial,n,(i),U^p(N_1,N_2),\Sigma_0} \]
is an isomorphism.

\end{enumerate}\end{lem}

We finish this section with an important vanishing result.
\begin{lem}\label{subvan5} Suppose that $R_0$ is an irreducible, noetherian $\Q$-algebra (resp. $\Z_{(p)}$-algebra) with the discrete topology. Suppose also that $U$ (resp. $U^p$) is a neat open compact subgroup of $P^{+}_{n,(i)}(\A^\infty)$ (resp. $P^{+}_{n,(i)}(\A^{p,\infty})$), that $N_2 \geq N_1 \geq 0$ are integers, and that $\Delta_0$ is a smooth admissible cone decomposition for $X_*(S_{n,(i),U}^{+})_\R^{\succ 0}$ (resp. $X_*(\cS_{n,(i),U}^{\ord,+})_\R^{\succ 0}$). Let $\pi$ denote the map
\[ \pi: T^{\natural,\wedge}_{n,(i),U,\Delta_0} \lra X^\natural_{n,(i),U} \]
(resp.
\[ \pi: \cT^{\ord,\natural,\wedge}_{n,(i),U^p(N_1,N_2),\Delta_0} \lra \cX^{\ord,\natural}_{n,(i),U^p(N_1,N_2)}).\]

Further suppose that $\cE$ is a coherent sheaf on the formal scheme $T^{\natural,\wedge}_{n,(i),U,\Delta_0} \times \Spf R_0$ (resp. $\cT^{\ord,\natural,\wedge}_{n,(i),U^p(N_1,N_2),\Delta_0} \times \Spf R_0$) with an exhaustive separated filtration, such that the pull back to $T^{+,\wedge}_{n,(i),U,\Delta_0} \times \Spf R_0$ (resp. $\cT^{\ord,+,\wedge}_{n,(i),U^p(N_1,N_2),\Delta_0} \times \Spf R_0$) of each
\[ \gr^i \cE \]
 is $L_{n,(i),\lin}(\Q)$-equivariantly (resp. $L_{n,(i),\lin}(\Z_{(p)})$-equivariantly) isomorphic to the pull back to $T^{+,\wedge}_{n,(i),U,\Delta_0} \times \Spf R_0$ (resp. $\cT^{\ord,+,\wedge}_{n,(i),U^p(N_1,N_2),\Delta_0} \times \Spf R_0$) of a locally free sheaf $\cF_i$ with $L_{n,(i),\lin}(\Q)$-action (resp. $L_{n,(i),\lin}(\Z_{(p)})$-action) over $X_{n,(i),U}^+ \times \Spec R_0$ (resp. $\cX_{n,(i),U}^{\ord,+} \times \Spec R_0$). 

Then for $i>0$
\[ R^i\pi_* (\cE \otimes \cI^{\natural,\wedge}_{\partial, n,(i),U, \Delta_0}) =  (0) \]
(resp.
\[ R^i\pi_* (\cE \otimes \cI^{\ord,\natural,\wedge}_{\partial, n,(i),U^p(N_1,N_2),\Delta_0}) = (0)). \]
\end{lem}

\pfbegin
We will treat the case of $T^{\natural,\wedge}_{n,(i),U,\Delta_0} \times \Spf R_0$, the other case 
being exactly similar. We can immediately reduce to the case that the pull back to $T^{+,\wedge}_{n,(i),U,\Delta_0} \times \Spf R_0$ of $\cE$ is $L_{n,(i),\lin}(\Q)$-equivariantly isomorphic to the pull back to $T^{+,\wedge}_{n,(i),U,\Delta_0} \Spf R_0$ of a locally free sheaf $\cF$ with $L_{n,(i),\lin}(\Q)$-action over $X_{n,(i),U}^+ \times \Spec R_0$.

Let $\pi^+$ denote the map
\[ \pi^+: T^{+,\wedge}_{n,(i),U,\Delta_0} \times \Spf R_0\lra X^\natural_{n,(i),U}. \times \Spec R_0\]
Also write $\pi^+=\pi_1^+ \circ \pi_2^+$, where
\[ \pi_1^+: A^{+}_{n,(i),U} \times \Spec R_0\lra X^\natural_{n,(i),U} \times \Spec R_0\]
and
\[ \pi_2^+: T^{+,\wedge}_{n,(i),U,\Delta_0} \times \Spec R_0\lra A^+_{n,(i),U} \times \Spec R_0. \]
By lemma \ref{dctorem3} we have that
\[ R^i\pi_{2,*}^+ (\cF \otimes \cI^\wedge_{\partial, n,(i),U, \Delta_0}) = \left\{ \begin{array}{ll} \cF \otimes  \prod_{a \in X^*(S_{n,(i),U}^+)^{>0}} \cL_U^+(a) 
& {\rm if}\,\, i=0 \\ (0) & {\rm otherwise.} \end{array} \right. \]
Then by lemma \ref{subvan1} (or in the case of $\cT^{\ord,\natural,\wedge}_{n,(i),U^p(N_1,N_2),\Delta_0} \times \Spf R_0$ lemma \ref{subvan1ord}) we deduce that
\[ \begin{array}{l} R^i\pi_{*}^+ (\cF \otimes \cI^\wedge_{\partial, n,(i),U, \Delta_0}) \\ = \left\{ \begin{array}{ll} \Ind_{\{1\}}^{L_{n,(i),\lin}(\Q)} \prod_{a^\natural \in X^*(S_{n,(i),U}^+)^{>0,\natural}}  (\pi_{A^+/X^\natural,*}\cL \otimes \cF)_U^+ (a^\natural)^{L_{n,(i),\lin}(\Q)} 
& {\rm if}\,\, i=0 \\ (0) & {\rm otherwise} \end{array} \right. \end{array} \]
Finally there is a spectral sequence
\[ H^i(L_{n,(i),\lin}(\Q), R^j\pi_*^+(\cF \otimes \cI^\wedge_{\partial, n,(i),U, \Delta_0})) \Rightarrow R^{i+j}\pi_* (\cF \otimes \cI^\wedge_{\partial, n,(i),U, \Delta_0}), \]
and so the present lemma follows on applying Shapiro's lemma.
\pfend

\begin{cor}\label{subvan5.5} Suppose that $U$ (resp. $U^p$) is a neat open compact subgroup of $P^{+}_{n,(i)}(\A^\infty)$ (resp. $P^{+}_{n,(i)}(\A^{p,\infty})$), that $N_2 \geq N_1 \geq 0$ are integers, and that $\Delta_0$ is a smooth admissible cone decomposition for $X_*(S_{n,(i),U}^{+})_\R^{\succ 0}$ (resp. $X_*(\cS_{n,(i),U}^{\ord,+})_\R^{\succ 0}$). Let $\pi$ denote the map
\[ \pi: T^{\natural,\wedge}_{n,(i),U,\Delta_0} \lra X^\natural_{n,(i),U} \]
(resp.
\[ \pi: \cT^{\ord,\natural,\wedge}_{n,(i),U^p(N_1,N_2),\Delta_0} \lra \cX^{\ord,\natural}_{n,(i),U^p(N_1,N_2)}).\]

Also suppose that $R_0$ is an irreducible noetherian $\Q$-algebra (resp. $\Z_{(p)}$-algebra) with the discrete
topology and that $\rho$ is a representation of $R_{n,(n),(i)}$ on a finite locally free $R_0$-module. 

Then for $i>0$
\[ R^i\pi_* (\cE_{n,(i),U,\Delta_0,\rho}^\natural \otimes \cI^{\natural,\wedge}_{\partial, n,(i),U, \Delta_0}) =  (0) \]
(resp.
\[ R^i\pi_* (\cE_{n,(i),U^p(N_1,N_2),\Delta_0,\rho}^{\ord,\natural} \otimes \cI^{\ord,\natural,\wedge}_{\partial, n,(i),U^p(N_1,N_2),\Delta_0}) = (0)). \]
\end{cor}

\newpage

\section{Compactification of Shimura Varieties.}

We now turn to the compactification of the $X_{n,U}$ and the $A^{(m)}_{n,U}$. 

\subsection{The minimal compactification.}\label{mincomp}

There is a canonically defined system of normal projective schemes with
$G_n(\A^\infty)$-action, $\{ X_{n,U}^\mini/\Spec
\Q\}$ (for $U \subset G_n(\A^\infty)$ a neat open compact subgroup), together with compatible,
$G_n(\A^\infty)$-equivariant, dense open
embeddings 
\[ j_U^\mini:X_{n,U} \into X_{n,U}^\mini. \]
These schemes are referred to as the {\em minimal} (or sometimes `Baily-Borel')
compactifications. (The introduction to \cite{pink} asserts that the scheme $X_{n,U}^\mini$ is the minimal normal compactification of $X_{n,U}$, although we won't need this fact.) For $g \in G_n(\A^\infty)$ and $g^{-1}Ug \subset U'$ the maps
\[ g: X_{n,U}^\mini \lra X_{n,U'}^\mini \]
are finite.

Write
\[ \partial X_{n,U}^\mini = X_{n,U}^\mini-X_{n,U}. \]
There is a family of closed sub-schemes
\[ \partial_0 X_{n,U}^\mini=X_{n,U}^\mini \supset \partial_1X_{n,U}^\mini=\partial X_{n,U}^\mini \supset \partial_2 X_{n,U}^\mini \supset ... \supset \partial_n X_{n,U}^\mini \supset \partial_{n+1} X_{n,U}^\mini=\emptyset \]
such that each
\[ \partial_i^0X_{n,U}^\mini=\partial_iX_{n,U}^\mini - \partial_{i+1}X_{n,U}^\mini \]
is smooth of dimension $(n-i)^2[F^+:\Q]$. The families $\{ \partial_i X_{n,U}^\mini \}$ and $\{ \partial_i^0X_{n,U}^\mini \}$ are families of schemes with $G_n(\A^\infty)$-action. Moreover we have a decomposition 
\[ \partial_i^0 X_{n,U}^\mini = \coprod_{h \in P^+_{n,(i)}(\A^\infty) \backslash G_n(\A^\infty) /U} X_{n,(i),hUh^{-1} \cap P_{n,(i)}^+(\A^\infty)}^\natural . \]
If $g \in G_n(\A^\infty)$ and if $g^{-1} U g \subset U'$ then the map 
\[ g: \partial_i^0 X_{n,U}^\mini \lra \partial_i^0 X_{n,U'}^\mini \]
is the coproduct of the maps
\[ g': X_{n,(i),hUh^{-1} \cap P^+_{n,(i)}(\A^\infty)}^\natural \lra X_{n,(i),h'U'(h')^{-1} \cap P^+_{n,(i)}(\A^\infty)}^\natural \]
where $hg=g'h'$ with $g' \in P_{n,(i)}^+(\A^\infty)$.
We will write $X_{n,U,i}^{\mini,\wedge}$ for the completion of $X_{n,U}^\mini$ along $\partial_i^0X_{n,U}^\mini$. (See theorem 7.2.4.1 and proposition 7.2.5.1 of \cite{kw1}.)

There is also a canonically defined system of normal quasi-projective schemes with
$G_n(\A^\infty)^\ord$-action, $\{ \cX_{n,U^p(N_1,N_2)}^{\ord,\mini}/\Spec
\Z_{(p)} \}$, together with compatible, dense open
embeddings 
\[ j_{U^p(N_1,N_2)}^\mini:\cX_{n,U^p(N_1,N_2)}^\ord \into \cX_{n,U^p(N_1,N_2)}^{\ord,\mini}, \]
which are $G_n(\A^\infty)^\ord$-equivariant. Suppose that $g \in G_n(\A^\infty)^\ord$ and that 
\[ g^{-1}U^p(N_1,N_2) g \subset (U^p)'(N_1',N_2'),\]
 then
\[ g: \cX^{\ord,\mini}_{n,U^p(N_1,N_2)} \lra \cX^{\ord,\mini}_{n,(U^p)'(N_1',N_2')} \]
is quasi-finite. If $p^{N_2-N_2'}\nu(g) \in \Z_p^\times$ and either $N_2'=N_2$ or $N_2'>0$, then it is also finite.  
On $\F_p$ fibres $\varsigma_p$ acts as absolute Frobenius composed with the forgetful map. (See theorem 6.2.1.1, proposition 6.2.2.1  and corollary 6.2.2.9 of \cite{kw2}.)

Write
\[ \partial \cX_{n,U^p(N_1,N_2)}^{\ord,\mini} = \cX_{n,U^p(N_1,N_2)}^{\ord,\mini}-\cX_{n,U^p(N_1,N_2)}^\ord. \]
There is a family of closed sub-schemes
\[ \begin{array}{r} \partial_0 \cX_{n,U^p(N_1,N_2)}^{\ord,\mini} \!\!\! =\cX_{n,U^p(N_1,N_2)}^{\ord,\mini} \supset \partial_1\cX_{n,U^p(N_1,N_2)}^{\ord,\mini}=\partial \cX_{n,U^p(N_1,N_2)}^{\ord,\mini} \supset \partial_2 \cX_{n,U^p(N_1,N_2)}^{\ord,\mini} \supset \\ ... \supset \partial_n \cX_{n,U^p(N_1,N_2)}^{\ord,\mini} \supset \partial_{n+1} \cX_{n,U^p(N_1,N_2)}^{\ord,\mini}=\emptyset \end{array} \]
such that each
\[ \partial_i^0\cX_{n,U^p(N_1,N_2)}^{\ord,\mini}=\partial_i\cX_{n,U^p(N_1,N_2)}^{\ord,\mini} - \partial_{i+1}\cX_{n,U^p(N_1,N_2)}^{\ord,\mini} \]
is smooth over $\Z_{(p)}$ of relative dimension $(n-i)^2[F^+:\Q]$. Then
\[ \{ \partial_i\cX_{n,U^p(N_1,N_2)}^{\ord,\mini} \} \]
 and 
 \[ \{ \partial_i^0\cX_{n,U^p(N_1,N_2)}^{\ord,\mini} \}\]
  are families of schemes with $G_n(\A^\infty)^\ord$-action. We will write $\cX_{n,U^p(N_1,N_2),i}^{\ord,\mini,\wedge}$ for the completion of $\cX_{n,U^p(N_1,N_2)}^{\ord,\mini}$ along $\partial_i^0\cX_{n,U^p(N_1,N_2)}^{\ord,\mini}$. We have a decomposition 
\[ \begin{array}{l} \partial_i^0 \cX_{n,U^p(N_1,N_2)}^{\ord,\mini} \!\! = \coprod_{h \in P^+_{n,(i)}(\A^\infty)^{\ord,\times} \backslash G_n(\A^\infty)^{\ord,\times} /U^p(N_1)} \cX_{n,(i),(hU^ph^{-1} \cap P_{n,(i)}^+(\A^{p,\infty}))(N_1,N_2)}^{\ord,\natural} \\  \amalg  \coprod_{h}  X_{n,(i),hU^p(N_1,N_2)h^{-1} \cap P_{n,(i)}^+(\A^\infty)}^\natural , \end{array} \]
where the second coproduct runs over
\[ h \in (P^+_{n,(i)}(\A^\infty) \backslash G_n(\A^\infty) /U^p(N_1,N_2))-(P^+_{n,(i)}(\A^\infty)^{\ord,\times} \backslash G_n(\A^\infty)^{\ord,\times} /U^p(N_1)).\]
(Again see theorems 6.2.1.1 and proposition 6.2.2.1 of \cite{kw2}.)
  
[We explain why the map
\[ P^+_{n,(i)}(\A^\infty)^\ord \backslash G_n(\A^\infty)^\ord /U^p(N_1) \lra P^+_{n,(i)}(\A^\infty) \backslash G_n(\A^\infty) /U^p(N_1,N_2) \]
is injective. It suffices to check that
\[ \begin{array}{rl} & (P^+_{n,(i)} \cap P^+_{n,(n)})(\Z_p) \backslash P^+_{n,(n)}(\Z_p) /  U_p(N_1,N_1)_{n,(n)}^+ \\ \into &P^+_{n,(i)}(\Q_p) \backslash G_n(\Q_p) /U_p(N_1,N_2)_n \\ =& P^+_{n,(i)}(\Z_p) \backslash G_n(\Z_p) /U_p(N_1,N_2)_n, \end{array} \]
or even that
\[ \begin{array}{l}  (P^+_{n,(i)} \cap P^+_{n,(n)})(\Z/p^{N_2}\Z) \backslash P^+_{n,(n)}(\Z/p^{N_2}\Z) / V \\ \into  P^+_{n,(i)}(\Z/p^{N_2}\Z) \backslash G_n(\Z/p^{N_2}\Z) /V, \end{array} \]
where
\[ V = \ker(P^+_{n,(n)}(\Z/p^{N_2}\Z) \ra L_{n,(n),\lin}(\Z/p^{N_1}\Z) ). \]
This is clear.]

If $g \in G_n(\A^\infty)^\ord$ and if $g^{-1} U^p(N_1,N_2) g \subset (U^p)'(N_1',N_2')$ then the map 
\[ g: \partial_i^0 \cX_{n,U^p(N_1,N_2)}^{\ord,\mini} \lra \partial_i^0 \cX_{n,(U^p)'(N_1',N_2')}^{\ord,\mini} \]
is the coproduct of the maps
\[ g': \cX_{n,(i),(hU^ph^{-1} \cap P^+_{n,(i)}(\A^{p,\infty}))(N_1,N_2)}^{\ord,\natural} \lra \cX_{n,(i),(h'(U^p)'(h')^{-1} \cap P^+_{n,(i)}(\A^{p,\infty}))(N_1',N_2')}^{\ord,\natural }\]
where $hg=g'h'$ with $g' \in P_{n,(i)}^+(\A^\infty)^\ord$, and of the maps
\[ g': X_{n,(i),hU^p(N_1,N_2)h^{-1} \cap P^+_{n,(i)}(\A^\infty)}^\natural \lra X_{n,(i),h'(U^p)'(N_1',N_2')(h')^{-1} \cap P^+_{n,(i)}(\A^\infty)}^\natural \]
where $hg=g'h'$ with $g' \in P_{n,(i)}^+(\A^\infty)$. (Again see theorems 6.2.1.1 and proposition 6.2.2.1 of \cite{kw2}.)

If $N_2' \geq N_2\geq N_1$ then the natural map
\[ \cX^{\ord,\mini}_{n,U^p(N_1,N_2')} \lra \cX^{\ord,\mini}_{n,U^p(N_1,N_2)} \]
is etale in a Zariski neighborhood of the $\F_p$-fibre, and the natural map
\[ \gX^{\ord,\mini}_{n,U^p(N_1,N_2')} \lra \gX^{\ord,\mini}_{n,U^p(N_1,N_2)} \]
between formal completions along the $\F_p$-fibres is an isomorphism. (See corollary 6.2.2.8 and example 3.4.5.5 of \cite{kw2}.) We will denote this $p$-adic formal scheme
\[ \gX^{\ord,\mini}_{n,U^p(N_1)} \]
and will denote its reduced subscheme
\[ \barX^{\ord,\mini}_{n,U^p(N_1)}. \]
We will also write 
\[ \partial \barX_{n,U^p(N_1)}^{\ord,\mini} = \barX_{n,U^p(N_1)}^{\ord,\mini}-\barX^\ord_{n,U^p(N_1)}. \]
The families $\{ \gX_{n,U^p(N)}^{\ord,\mini} \}$ and $\{ \barX_{n,U^p(N)}^{\ord,\mini} \}$ and $\{ \partial \barX_{n,U^p(N)}^{\ord,\mini} \}$ have $G_n(\A^\infty)^\ord$-actions.
There is a family of closed sub-schemes
\[ \begin{array}{l} \partial_0 \barX_{n,U^p(N)}^{\ord,\mini}=\barX_{n,U^p(N)}^{\ord,\mini} \supset \partial_1\barX_{n,U^p(N)}^{\ord,\mini}=\partial \barX_{n,U^p(N)}^{\ord,\mini} \supset \partial_2 \barX_{n,U^p(N)}^{\ord,\mini} \supset ... \\ \multicolumn{1}{r}{... \supset \partial_n \barX_{n,U^p(N)}^{\ord,\mini} \supset \partial_{n+1} \barX_{n,U^p(N)}^{\ord,\mini}=\emptyset} \end{array}  \]
such that each
\[ \partial_i^0\barX_{n,U^p(N)}^{\ord,\mini}=\partial_i\barX_{n,U^p(N)}^{\ord,\mini}- \partial_{i+1}\barX_{n,U^p(N)}^{\ord,\mini} \]
is smooth of dimension $(n-i)^2[F^+:\Q]$. Then $\{ \partial_i \barX_{n,U^p(N)}^\mini \}$ and $\{ \partial_i^0\barX_{n,U^p(N)}^\mini \}$ are families of schemes with $G_n(\A^\infty)^\ord$-action. Moreover we have a decomposition 
\[ \partial_i^0 \barX_{n,U^p(N)}^{\ord,\mini} = \coprod_{h \in P^+_{n,(i)}(\A^\infty)^{\ord,\times} \backslash G_n(\A^\infty)^{\ord,\times} /U^p(N)} \barX^{\ord,\natural}_{n,(i),(hU^ph^{-1} \cap P_{n,(i)}^+(\A^{p,\infty}))(N)} . \]
If $g \in G_n(\A^\infty)^\ord$ and if $g^{-1} U^p(N) g \subset (U^p)'(N')$ then the map 
\[ g: \partial_i^0 \barX_{n,U^p(N)}^{\ord,\mini} \lra \partial_i^0 \barX_{n,(U^p)'(N')}^{\ord,\mini} \]
is the coproduct of the maps
\[ g': \barX_{n,(i),(hU^ph^{-1} \cap P^+_{n,(i)}(\A^{p,\infty}))(N)}^{\ord,\natural} \lra \barX_{n,(i),(h'(U^p)'(h')^{-1} \cap P^+_{n,(i)}(\A^{p,\infty}))(N')}^{\ord,\natural} \]
where $hg=g'h'$ with $g' \in P_{n,(i)}^+(\A^\infty)^\ord$. 
In particular $\varsigma_p$ acts as absolute Frobenius.

The schemes $\cX^{\ord,\min}_{n,U^p(N_1,N_2)}$ are not proper. There are proper integral models of the schemes $X_{n,U}^\mini$, but we have less control over them.

More specifically suppose that $U \subset G_n(\A^{p,\infty} \times \Z_p)$ is an open compact subgroup whose projection to $G_n(\A^{p,\infty})$ is neat.  Then there is a normal, projective, flat $\Z_{(p)}$-scheme $\cX_{n,U}^\mini$ with generic fibre $X_{n,U}^\mini$. If $g \in G_n(\A^{p,\infty} \times \Z_p)$ and if 
\[ g^{-1}Ug\subset U'\]
 then there is a map
\[g:\cX_{n,U}^\mini \lra \cX_{n,U'}^\mini \]
extending the map $g: X_{n,U}^\mini \ra X_{n,U'}^\mini$. This gives the system $\{ \cX^\mini_{n,U}\}$ an action of $G_n(\A^{p,\infty} \times \Z_p)$. 
We set 
\[ \barX_{n,U}^\mini=\cX_{n,U}^\mini \times_{\Z_{(p)}} \F_p. \] 
On $\cX_{n,U}^\mini$ there is an ample line bundle $\omega_{U}$, and the system of line bundles $\{ \omega_{U}\}$ over $\{ \cX_{n,U}^\mini\}$ has an action of $G_n(\A^{p,\infty}\times \Z_p )$. The pull back of $\omega_U$ to $X_{n,U}$ is $G_n(\A^{p,\infty}\times \Z_p)$-equivariantly identified with $\wedge^{n[F:\Q]} \Omega_{n,U}$.(See propositions 2.2.1.2 and 2.2.3.1 of \cite{kw2}.)

Moreover there are canonical sections 
\[ \hasse_{U} \in H^0(\barX_{U}^\mini,\omega_{U}^{\otimes(p-1)}) \]
such that 
\[ g^* \hasse_{U'} = \hasse_{U} \]
whenever $g \in  G_n(\A^{p,\infty}\times \Z_p)$ and $U' \supset g^{-1}Ug$. 
We will write $\barX_{n,U}^{\mini,\nord}$ for the zero locus in $\barX_{n,U}^\mini$ of $\hasse_{U}$. (See corollaries 6.3.1.7 and 6.3.1.8 of \cite{kw2}.)
Then $\barX_{n,U}^\mini - \barX_{n,U}^{\min,\nord}$ is relatively affine over $\barX_{n,U}^\mini$ associated to the sheaf of algebras
\[ \left( \bigoplus_{i=0}^\infty \omega^{\otimes (p-1)ai} \right)/(\hasse_U^a-1) \]
for any $a \in \Z_{>0}$. It is also affine over $\F_p$ associated to the algebra
\[ \left( \bigoplus_{i=0}^\infty H^0(\barX_{n,U}^\mini, \omega^{\otimes (p-1)ai}) \right)/(\hasse_U^a-1) \]
for any $a \in \Z_{>0}$. 

There are $G_n(\A^{\infty})^{\ord,\times}$-equivariant open embeddings
\[ \cX_{n,U^p(N_1,N_2)}^{\ord,\mini} \into \cX_{n,U^p(N_1,N_2)}^\mini . \]
The induced map on $\F_p$-fibres is an open and closed embedding
\[ \barX_{n,U^p(N_1,N_2)}^{\ord,\mini} \into \barX_{n,U^p(N_1,N_2)}^\mini-\barX_{n,U^p(N_1,N_2)}^{\mini,\nord}. \]
(See proposition 6.3.2.2 of \cite{kw2}.) In the case $N_1=N_2=0$ this is in fact an isomorphism. (See lemmas 6.3.2.7 and 6.3.2.9 of \cite{kw2}.) We remark that for $N_2>0$ this map is not an isomorphism. the definition of $\barX_{n,U^p(N_1,N_2)}^{\ord}$ requires not only that the universal abelian variety is ordinary, the condition that defines $\barX_{n,U^p(N_1,N_2)}-\barX_{n,U^p(N_1,N_2)}^{\mini,\nord}$, but also that the universal subgroup of $\cA^{\univ}[p^{N_2}]$ is connected. 

Also the pull back of $\omega_{U^p(N_1,N_2)}$ to $\cX^\ord_{n,U^p(N_1,N_2)}$ is $G_n(\A^{\infty})^{\ord,\times}$-equivariantly identified with the sheaf $\wedge^{n[F:\Q]} \Omega^\ord_{n,U^p(N_1,N_2)}$. If $g \in G_n(\A^{\infty})^{\ord,\times}$ and
\[ g^{-1}(U^p)'(N_1',N_2')g \subset U^p(N_1,N_2), \]
 then the commutative square
\[ \begin{array}{ccc} \cX^{\ord,\mini}_{n,(U^p)'(N_1,N_2)} & \stackrel{g}{\lra} & \cX^{\ord,\mini}_{n,U^p(N_1,N_2)} \\
\da && \da \\ \cX_{n,(U^p)'(N_1,N_2)}^{\mini} & \stackrel{g}{\lra} & \cX_{n,U^p(N_1,N_2)}^{\mini}
\end{array} \]
is a pull-back square.
(See theorem 6.2.1.1 and proposition 6.2.2.1 of \cite{kw2}.)

\newpage \subsection{Cone decompositions.}

Let $U \subset G^{(m)}(\A^\infty)$ be an open compact subgroup.
By a {\em $U$-admissible cone decomposition} $\Sigma$ of $G_n^{(m)}(\A^\infty) \times \pi_0(G_n(\R)) \times \gC^{(m)}$ we shall mean a set of closed subsets $\sigma \subset G_n^{(m)}(\A^\infty) \times \pi_0(G_n(\R)) \times \gC^{(m)}$ such that
\begin{enumerate}
\item each $\sigma$ is contained in $\{(g,\delta) \} \times \gC_W^{(m),\succ 0}$ for some isotropic subspace $W \subset V_n$ and some $(g,\delta)  \in G_n^{(m)}(\A^\infty) \times \pi_0(G_n(\R))$ and is the set of $\R_{\geq 0}$-linear combinations of a finite set of elements of $\Herm_{V/W^\perp} \times W^m$;

\item if $\sigma \in \Sigma$ then any face of $\sigma$ also lies in $\Sigma$;

\item if $\sigma, \sigma' \in \Sigma$ then either $\sigma \cap \sigma'=\emptyset$ or $\sigma \cap \sigma'$ is a face of $\sigma$ and $\sigma'$;

\item $G_n^{(m)}(\A^\infty) \times \pi_0(G_n(\R)) \times \gC^{(m)} = \bigcup_{\sigma \in \Sigma} \sigma$;

\item $\Sigma$ is left invariant by the diagonal action of $G_n^{(m)}(\Q)$ on $G_n^{(m)}(\A^\infty) \times \pi_0(G_n(\R)) \times \gC^{(m)}$;

\item $\Sigma$ is invariant by the right action of $U$ on $G_n^{(m)}(\A^\infty) \times \pi_0(G_n(\R)) \times \gC^{(m)}$ (acting only on the first factor);

\item $G_n^{(m)}(\Q) \backslash \Sigma / U$ is a finite set.

\end{enumerate}
Note that if $U' \subset U$ and if $\Sigma$ is a $U$-admissible cone decomposition of $G_n^{(m)}(\A^\infty) \times \pi_0(G_n(\R)) \times \gC^{(m)}$ then $\Sigma$ is also $U'$-admissible. We will call a set $\Sigma$ of closed subsets of $G_n^{(m)}(\A^\infty) \times \pi_0(G_n(\R)) \times \gC^{(m)}$ an {\em admissible cone decomposition} of $G_n^{(m)}(\A^\infty) \times \pi_0(G_n(\R)) \times \gC^{(m)}$ if it is $U$-admissible for some open compact subgroup $U$. 

We remark that different authors use the term `$U$-admissible cone decomposition' in somewhat different ways. 

We call $\Sigma'$ a {\em refinement} of $\Sigma$ if every element of $\Sigma$ is a union of elements of $\Sigma'$. 
We define a partial order on the set of pairs $(U,\Sigma)$, where $U \subset G_n^{(m)}(\A^\infty)$ is an open compact subgroup and $\Sigma$ is a $U$-admissible cone decomposition of $G_n^{(m)}(\A^\infty) \times \pi_0(G_n(\R)) \times \gC^{(m)}$, as follows: we set
\[ (U' ,\Sigma') \geq (U,\Sigma) \]
if and only if $U' \subset U$ and $\Sigma'$ is a refinement of $\Sigma$. If $g \in G_n^{(m)}(\A^\infty)$ and $\Sigma$ is a $U$-admissible cone decomposition of $G_n^{(m)}(\A^\infty) \times \pi_0(G_n(\R)) \times \gC^{(m)}$, then 
\[ \Sigma g = \{ \sigma (g\times 1): \,\, \sigma \in \Sigma\} \]
is a $g^{-1}U g$-admissible cone decomposition of $G_n^{(m)}(\A^\infty) \times \pi_0(G_n(\R)) \times \gC^{(m)}$. The action of $G_n^{(m)}(\A^\infty)$ preserves $\geq$. 

There is a natural projection
\[ G_n^{(m)}(\A^\infty) \times  \pi_0(G_n(\R)) \times \gC^{(m)} \onto G_n(\A^\infty) \times \pi_0(G_n(\R)) \times \gC. \]
We will call admissible cone decompositions $\Sigma$ of $G_n^{(m)}(\A^\infty) \times \pi_0(G_n(\R)) \times \gC^{(m)}$ and $\Delta$ of $G_n(\A^\infty) \times \pi_0(G_n(\R)) \times \gC$ {\em compatible} if the image of every $\sigma \in \Sigma$ is contained in an element of $\Delta$. If in addition $\Sigma$ is $U$-admissible, $\Delta$ is $U'$-admissible and $U'$ contains the image of $U$ in $G_n(\A^\infty)$ we will say that $(U,\Sigma)$ and $(U',\Delta)$ are {\em compatible} and write
\[ (U,\Sigma) \geq (U',\Delta'). \]

Now let $U^p \subset G^{(m)}(\A^{p,\infty})$ be an open compact subgroup and let $N \geq 0$ be an integer and consider $U^p(N) \subset G_n^{(m)}(\A^\infty)^{\ord,\times}$.
By a {\em $U^p(N)$-admissible cone decomposition}) $\Sigma$ of $(G_n^{(m)}(\A^\infty) \times \pi_0(G_n(\R)) \times \gC^{(m)})^\ord$ we shall mean a set of closed subsets $\sigma \subset (G_n^{(m)}(\A^\infty) \times \pi_0(G_n(\R)) \times \gC^{(m)})^\ord$ such that
\begin{enumerate}
\item each $\sigma$ is contained in $\{(g,\delta) \} \times \gC_W^{(m),\succ 0}$ for some isotropic subspace $W \subset V_n$ and some $(g,\delta)  \in G_n^{(m)}(\A^\infty) \times \pi_0(G_n(\R))$ and is the set of $\R_{\geq 0}$-linear combinations of a finite set of elements of $\Herm_{V/W^\perp} \times W^m$;

\item if $\sigma \in \Sigma$ then any face of $\sigma$ also lies in $\Sigma$;

\item if $\sigma, \sigma' \in \Sigma$ then either $\sigma \cap \sigma'=\emptyset$ or $\sigma \cap \sigma'$ is a face of $\sigma$ and $\sigma'$;

\item $(G_n^{(m)}(\A^\infty) \times \pi_0(G_n(\R)) \times \gC^{(m)})^\ord = \bigcup_{\sigma \in \Sigma} \sigma$;

\item if $\sigma \in \Sigma$ , if $\gamma \in G_n^{(m)}(\Q)$ and if $u \in U^p(N,N)$ are such that $\gamma \sigma u \subset (G_n^{(m)}(\A^\infty) \times \pi_0(G_n(\R)) \times \gC^{(m)})^\ord$, then $\gamma \sigma u \in \Sigma$;

\item there is a finite subset of $\Sigma$ such that any element of $\Sigma$ has the form $\gamma \sigma u$ with $\gamma \in G_n^{(m)}(\Q)$ and $u \in U^p(N,N)$ and $\sigma$ in the given finite subset.

\end{enumerate}
Note that if $(U^p)'(N') \subset U^p(N)$  and if $\Sigma$ is a $U^p(N)$-admissible cone decomposition of $(G_n^{(m)}(\A^\infty) \times \pi_0(G_n(\R)) \times \gC^{(m)})^\ord$ then $\Sigma$ is also $(U^p)'(N')$-admissible. We will call a set $\Sigma$ of closed subsets of $(G_n^{(m)}(\A^\infty) \times \pi_0(G_n(\R)) \times \gC^{(m)})^\ord$ an {\em admissible cone decomposition} of $(G_n^{(m)}(\A^\infty) \times \pi_0(G_n(\R)) \times \gC^{(m)})^\ord$ if it is $U^p(N)$-admissible for some open compact subgroup $U^p$ and for some $N$. 

If $\Sigma$ is a $U^p(N_1,N_2)$-admissible cone decomposition of $G_n^{(m)}(\A^\infty) \times \pi_0(G_n(\R)) \times \gC^{(m)}$ then
\[ \Sigma^\ord =\{ \sigma \in \Sigma: \,\,  \sigma \subset (G_n^{(m)}(\A^\infty) \times \pi_0(G_n(\R)) \times \gC^{(m)})^\ord \} \]
is a $U^p(N_1)$-admissible cone decomposition for $(G_n^{(m)}(\A^\infty) \times \pi_0(G_n(\R)) \times \gC^{(m)})^\ord$.

We call $\Sigma'$ a {\em refinement} of $\Sigma$ if every element of $\Sigma$ is a union of elements of $\Sigma'$. 
We define a partial order on the set of pairs $(U^p(N),\Sigma)$, where $U^p \subset G_n^{(m)}(\A^{p,\infty})$ is an open compact subgroup, $N \in \Z_{\geq 0}$ and $\Sigma$ is a $U^p(N)$-admissible cone decomposition of $(G_n^{(m)}(\A^\infty) \times \pi_0(G_n(\R)) \times \gC^{(m)})^\ord$, as follows: we set
\[ ((U^p)'(N') ,\Sigma') \geq (U^p(N),\Sigma) \]
if and only if $(U^p)'(N') \subset U^p(N)$ and $\Sigma'$ is a refinement of $\Sigma$. If $g \in G_n^{(m)}(\A^\infty)^\ord$ and $\Sigma$ is a $U^p(N)$-admissible cone decomposition of $(G_n^{(m)}(\A^\infty) \times \pi_0(G_n(\R)) \times \gC^{(m)})^\ord$, then 
\[ \Sigma g = \{ \sigma (g\times 1): \,\, \sigma \in \Sigma\} \]
is a $g^{-1}U^p(N) g$-admissible cone decomposition of 
\[ (G_n^{(m)}(\A^\infty) \times \pi_0(G_n(\R)) \times \gC^{(m)})^\ord.\]
 The action of $G_n^{(m)}(\A^\infty)^\ord$ preserves $\geq$. 

There is a natural projection
\[ (G_n^{(m)}(\A^\infty) \times  \pi_0(G_n(\R)) \times \gC^{(m)})^\ord \onto (G_n(\A^\infty) \times \pi_0(G_n(\R)) \times \gC)^\ord. \]
We will call admissible cone decompositions $\Sigma$ of $(G_n^{(m)}(\A^\infty) \times \pi_0(G_n(\R)) \times \gC^{(m)})^\ord$ and $\Delta$ of $(G_n(\A^\infty) \times \pi_0(G_n(\R)) \times \gC)^\ord$ {\em compatible} if the image of every $\sigma \in \Sigma$ is contained in an element of $\Delta$. If in addition $\Sigma$ is $U^p(N)$-admissible, $\Delta$ is $(U^p)'(N')$-admissible and $(U^p)'(N')$ contains the image of $U^p(N)$ in $G_n(\A^\infty)^\ord$ we will say that $(U^p(N),\Sigma)$ and $((U^p)'(N'),\Delta)$ are {\em compatible} and write
\[ (U^p(N),\Sigma) \geq ((U^p)'(N'),\Delta'). \]

If $\Sigma$ is a $U$-admissible cone decomposition of $G_n^{(m)}(\A^\infty) \times  \pi_0(G_n(\R)) \times \gC^{(m)}$ and if $h \in G_n^{(m)}(\A^\infty)$ then we define an admissible cone decomposition $\Sigma(h)_0$ for 
\[ X_*(S^{(m),+}_{n,(i),hU h^{-1}\cap P^{(m),+}_{n,(i)}(\A^\infty)})^{\succ 0}_\R \]
as follows: The cones in $\Sigma(h)_0$ over an element 
\[ y=[h'((hU h^{-1}\cap P^{(m),+}_{n,(i)}(\A^\infty))/ (hU h^{-1}\cap P^{(m)}_{n,(i)}(\A^\infty)))] \in Y_{n,(i),hU h^{-1}\cap P^{(m),+}_{n,(i)}(\A^\infty)}^{(m),+}\]
 are the cones 
 \[ \sigma \subset \gC^{(m),\succ 0}(V_{n,(i)}) \cong X_*(S^{(m),+}_{n,(i),hU h^{-1}\cap P^{(m),+}_{n,(i)}(\A^\infty)})_{\R,y}^{\succ 0} \]
  with
\[ \{(h'h,1)\} \times \sigma \in \Sigma. \]
This does not depend on the representative $h'$ we choose for $y$. It also only depends on 
\[ h \in P_{n,(i)}^{(m)}(\A^\infty) \backslash G_n^{(m)}(\A^\infty)/U. \]
If $h_1 \in L_{n,(i),\lin}^{(m)}(\A^\infty)$ then under the natural isomorphism
\[ h_1:Y_{n,(i),hU h^{-1}\cap P^{(m),+}_{n,(i)}(\A^\infty)}^{(m),+} \liso Y_{n,(i),h_1hU (h_1h)^{-1}\cap P^{(m),+}_{n,(i)}(\A^\infty)}^{(m),+} \]
we see that $\Sigma(h)_0$ and $\Sigma(h_1h)_0$ correspond. 

Similarly if $\Sigma$ is a $U^p(N)$-admissible cone decomposition of 
\[ (G_n^{(m)}(\A^\infty) \times  \pi_0(G_n(\R)) \times \gC^{(m)})^\ord\]
 and if $h \in G_n^{(m)}(\A^\infty)^\ord$ then we define an admissible cone decomposition $\Sigma(h)_0$ for 
\[ X_*\left(\cS^{(m),\ord,+}_{n,(i),(h^pU^p(N) (h^p)^{-1}\cap P^{(m),+}_{n,(i)}(\A^{p,\infty}))(N)}\right)^{\succ 0}_\R \]
as follows: The cones in $\Sigma(h)_0$ over an element $y$ given as
\[ \begin{array}{l}  [h'(h^pU^p(N) (h^p)^{-1}\cap P^{(m),+}_{n,(i)}(\A^{p,\infty}))(N)/ (h^pU^p(N) (h^p)^{-1}\cap P^{(m),+}_{n,(i)}(\A^{p,\infty}))(N)] \\ \multicolumn{1}{r}{\in \cY_{n,(i),(h^pU^p(N) (h^p)^{-1}\cap P^{(m),+}_{n,(i)}(\A^{p,\infty}))(N)}^{(m),\ord,+} }\end{array}\]
 are the cones 
 \[ \sigma \subset \gC^{(m),\succ 0}(V_{n,(i)}) \cong X_*(\cS^{(m),\ord,+}_{n,(i),(h^pU^p(N) (h^p)^{-1}\cap P^{(m),+}_{n,(i)}(\A^{p,\infty}))(N)})_{\R,y}^{\succ 0} \]
  with
\[ \{(h'h,1)\} \times \sigma \in \Sigma. \]
This does not depend on the representative $h'$ we choose for $y$. It also only depends on 
\[ h \in P_{n,(i)}^{(m)}(\A^\infty)^\ord \backslash G_n^{(m)}(\A^\infty)^\ord/U^p(N). \]
If $h_1 \in L_{n,(i),\lin}^{(m)}(\A^\infty)^\ord$ then under the natural isomorphism
\[ h_1:\cY_{n,(i),(h^pU^p(N) (h^p)^{-1}\cap P^{(m),+}_{n,(i)}(\A^{p,\infty}))(N)}^{(m),\ord,+} \liso \cY_{n,(i),(h_1^ph^pU^p(N) (h_1^ph^p)^{-1}\cap P^{(m),+}_{n,(i)}(\A^{p,\infty}))(N)}^{(m),\ord,+} \]
we see that $\Sigma(h)_0$ and $\Sigma(h_1h)_0$ correspond. 

There are sets $\cJ^{(m),\tor}_n$ (resp. $\cJ^{(m),\tor,\ord}_n$) of pairs $(U,\Sigma)$ (resp. $(U^p(N),\Sigma)$) where $U \subset G_n^{(m)}(\A^\infty)$ is a neat open compact subgroup (resp. $U^p \subset G_n^{(m)}(\A^{p,\infty})$ is a neat open compact subgroup and $N\in \Z_{\geq 0}$) and $\Sigma$ is a $U$-admissible (resp. $U^p(N)$-admissible) cone decomposition of $G_n^{(m)}(\A^\infty) \times \pi_0(G_n(\R)) \times \gC^{(m)}$ (resp. $(G_n^{(m)}(\A^\infty) \times \pi_0(G_n(\R)) \times \gC^{(m)})^\ord$), with a number of properties which will be listed in this section and the next section. (See \cite{kw2}.)

Firstly we have the following properties:
\begin{enumerate}
\item The sets $\cJ^{(m),\tor}_n$ (resp. $\cJ^{(m),\tor,\ord}_n$) are invariant under the action of $G_n^{(m)}(\A^\infty)$ (resp. $G_n^{(m)}(\A^\infty)^{\ord,\times}$). 

\item If $U$ is any neat open compact subgroup of $G_n^{(m)}(\A^\infty)$, then there is some $\Sigma$ with $(U,\Sigma) \in \cJ^{(m),\tor}_n$.

\item If $U^p$ is any neat open compact subgroup of $G_n^{(m)}(\A^{p,\infty})$ and if $N \in \Z_{\geq 0}$, then there is some $\Sigma$ with $(U^p(N),\Sigma) \in \cJ^{(m),\tor,\ord}_n$.

\item If $(U,\Sigma) \in \cJ^{(m),\tor}_n$ and if $U' \subset U$ then there exists $(U',\Sigma') \in \cJ^{(m),\tor}_n$ with $(U',\Sigma') \geq (U,\Sigma)$.

\item If $(U^p(N),\Sigma') \in \cJ^{(m),\tor,\ord}_n$, if $N' \geq N$ and if $(U^p)'(N') \subset U^p(N)$ then there exists an element $((U^p)'(N'),\Sigma') \in \cJ^{(m),\tor,\ord}_n$ with $((U^p)'(N'),\Sigma') \geq (U^p(N),\Sigma)$.

\item  If $(U',\Sigma') \geq (U,\Sigma)$ are elements of $\cJ^{(m),\tor}_n$ and if moreover $U'$ is a normal subgroup of $U$, then we may choose $(U',\Sigma'')\in \cJ^{(m),\tor}_n$ such that $\Sigma''$ is $U$-invariant and such that $(U',\Sigma'') \geq (U',\Sigma')$.

\item  If $((U^p)'(N'),\Sigma') \geq (U^p(N),\Sigma)$ are elements of $\cJ^{(m),\tor,\ord}_n$ and if moreover $(U^p)'$ is a normal subgroup of $U^p$, then we may choose an element $((U^p)'(N'),\Sigma'')\in \cJ^{(m),\tor}_n$ such that $\Sigma''$ is $U^p(N)$-invariant and such that $((U^p)'(N'),\Sigma'') \geq ((U^p)'(N'),\Sigma')$.

\item If $(U,\Sigma)$ and $(U,\Sigma') \in \cJ^{(m),\tor}_n$ (resp.  $(U^p(N),\Sigma)$ and $(U^p(N),\Sigma') \in \cJ^{(m),\tor.\ord}_n$) then there exists $(U,\Sigma'') \in \cJ^{(m),\tor}_n$ (resp. $(U^p(N),\Sigma'') \in \cJ^{(m),\tor,\ord}_n$) with $(U,\Sigma'')\geq (U,\Sigma)$ and $(U,\Sigma'')\geq (U,\Sigma')$ (resp. with $(U^p(N),\Sigma'')\geq (U^p(N),\Sigma)$ and $(U^p(N),\Sigma'')\geq (U^p(N),\Sigma')$).

\item If $(U',\Delta) \in \cJ^\tor_n$ (resp. $((U^p)'(N'),\Delta) \in \cJ^{\tor,\ord}_n$) and if $U$ is a neat open compact subgroup of $G_n^{(m)}(\A^\infty)$ mapping into $U'$ (resp. $U^p$ is a neat open compact subgroup of $G_n^{(m)}(\A^{p,\infty})$ mapping into $(U^p)'$ and $N \geq N'$), then there exists $(U,\Sigma) \in \cJ^{(m),\tor}_n$ (resp. $(U^p(N),\Sigma) \in \cJ^{(m),\tor,\ord}_n$) compatible with $(U',\Delta)$ (resp. $((U^p)'(N'),\Delta)$).

\item If $(U^p(N_1,N_2),\Sigma) \in \cJ^{(m),\tor}_n$ then $(U^p(N_1),\Sigma^\ord) \in \cJ_n^{(m),\tor,\ord}$.

\item If $(U^p(N),\Sigma') \in \cJ^{(m),\tor,\ord}_n$ and if $N' \geq N$, then there is an element $(U^p(N,N'),\Sigma) \in \cJ_n^{(m),\tor}$ with $\Sigma^\ord=\Sigma'$.

\item If $(U^p(N_1,N_2),\Sigma)$ and $(U^p(N_1,N_2),\Sigma') \in \cJ^{(m),\tor}_n$ with $\Sigma^\ord = (\Sigma')^\ord$, then there is an element $(U^p(N_1,N_2),\Sigma'') \in \cJ^{(m),\tor}_n$ with $(\Sigma'')^\ord=\Sigma^\ord=(\Sigma')^\ord$ and
$(U^p(N_1,N_2),\Sigma'') \geq (U^p(N_1,N_2),\Sigma)$ and $(U^p(N_1,N_2),\Sigma'') \geq (U^p(N_1,N_2),\Sigma')$.

\item If $(U^p(N_1,N_2),\Sigma)$ and $((U^p)'(N_1',N_2'),\Sigma') \in \cJ^{(m),\tor}_n$ with $(U^p)'(N_1',N_2') \subset U^p(N_1,N_2)$ and $(\Sigma')^\ord$ refining $\Sigma^\ord$, then there also exists another pair $((U^p)'(N_1',N_2'),\Sigma'') \in \cJ^{(m),\tor}_n$ with $\Sigma''$ refining both $\Sigma$ and $\Sigma'$ and with $(\Sigma'')^\ord=(\Sigma')^\ord$.

\item If $(U^p(N_1,N_2),\Delta) \in \cJ^\tor_n$ and $((U^p)'(N_1',N_2'),\Sigma') \in \cJ^{(m),\tor}_n$ are such that $(U^p)'(N_1',N_2') \subset U^p(N_1,N_2)$ and $(\Sigma')^\ord$ is compatible with $\Delta^\ord$, then there exists $((U^p)'(N_1',N_2'),\Sigma'') \in \cJ^{(m),\tor}_n$ with $\Sigma''$ refining $\Sigma'$ and compatible with $\Delta$ and with $(\Sigma'')^\ord=(\Sigma')^\ord$.

\item If $(U^p(N_1,N_2),\Sigma) \in \cJ_n^{(m),\tor}$ and if $N_2'\geq N_2$ then there exists a pair $(U^p(N_1,N_2'),\Sigma') \in \cJ_n^{(m),\tor}$ with $(\Sigma')^\ord =\Sigma^\ord$.
\end{enumerate}
(See proposition 7.1.1.21 of \cite{kw2}.) 

Secondly if $(U,\Sigma)\in \cJ^{(m),\tor}_n$ (resp. $(U^p(N),\Sigma)\in \cJ^{(m),\tor,\ord}_n$) and if $h \in G_n^{(m)}(\A^\infty)$ (resp. $h \in G_n^{(m)}(\A^\infty)^\ord$) then $\Sigma(h)_0$ is smooth.

Thirdly if $(U,\Sigma) \in \cJ^{(m),\tor}_n$, then there is a simplicial complex $\cS(U,\Sigma)$ whose simplices are in bijection with the cones in  
\[ G_n^{(m)}(\Q) \backslash \Sigma /U \]
which have dimension bigger than $0$, and have the same face relations. We will write $\cS(U,\Sigma)_{\leq i}$ for the subcomplex of $\cS(U,\Sigma)$ consisting of simplices associated to the orbits of cones $(g,\delta) \times \sigma \in \Sigma$ with $\sigma \subset \gC^{(m), >0}(W)$ for some $W$ with $\dim_F W \leq i$. We will also set
\[ |\cS(U,\Sigma)|_{=i}= | \cS(U,\Sigma)_{\leq i} |- |  \cS(U,\Sigma)_{\leq i-1} |, \]
an open subset of $| \cS(U,\Sigma)_{\leq i} |$.
Then one sees that
\[ | \cS(U,\Sigma) | \cong
G_n^{(m)}(\Q) \backslash \left( (G_n^{(m)}(\A^\infty)/U) \times \pi_0(G_n(\R)) \times ((\gC^{(m)}-\gC^{(m)}_{=0})/\R^\times_{>0}) \right) \]
and
\[  \begin{array}{ll} &  |\cS(U,\Sigma)|_{=i} \\ \cong & 
G_n^{(m)}(\Q) \backslash \left( (G_n^{(m)}(\A^\infty)/U) \times \pi_0(G_n(\R)) \times ((\gC^{(m)}_{=i})/\R^\times_{>0}) \right)
\\ \cong & \coprod_{h \in P_{n,(i)}^{(m),+}(\A^\infty) \backslash G_n^{(m)}(\A^\infty)/U} 
L_{n,(i)}^{(m)}(\Q) \backslash L_{n,(i)}^{(m)}(\A) / \\ &  
(hUh^{-1} \cap P_{n,(i)}^{(m),+}(\A^\infty))L_{n,(i),\herm}(\R)^+ (L_{n,(i),\lin}^{(m)}(\R) \cap U_{n,\infty}^0)A_{n,(i)}(\R)^0 . \end{array}\]
(See section \ref{herm}.)

If $(U^p(N),\Sigma) \in \cJ^{(m),\tor,\ord}_n$ then there is a simplicial complex $\cS(U^p(N),\Sigma)^\ord$ whose simplices are in bijection with equivalence classes of cones of dimension greater than $0$ in  $\Sigma$, where $\sigma$ and $\sigma'$ are considered equivalent if $\sigma'=\gamma \sigma u$ for some $\gamma \in G_n^{(m)}(\Q)$ and some $u\in U^p(N,N)$.
We will write $\cS(U^p(N),\Sigma)^\ord_{\leq i}$ for the subcomplex of $\cS(U^p(N),\Sigma)^\ord$ consisting of simplices associated to the orbits of cones $(g,\delta) \times \sigma \in \Sigma$ with $\sigma \subset \gC^{(m), >0}(W)$ for some $W$ with $\dim_F W \leq i$. We will also set
\[ |\cS(U^p(N),\Sigma)^\ord|_{=i}= | \cS(U^p(N),\Sigma)^\ord_{\leq i} |- |  \cS(U^p(N),\Sigma)^\ord_{\leq i-1} |, \]
an open subset of $| \cS(U^p,\Sigma)^\ord_{\leq i} |$.
Then we see that
\[\begin{array}{rl} & | \cS(U^p(N),\Sigma)^\ord |\\  \cong &
G_n^{(m)}(\Q) \backslash \left( (G_n^{(m)}(\A^\infty)/U^p(N)) \times \pi_0(G_n(\R)) \times (\gC^{(m)}-\gC^{(m)}_{=0}) /\R^\times_{>0}) \right)^\ord,\end{array} \]
where 
\[ \left( (G_n^{(m)}(\A^\infty)/U^p(N)) \times \pi_0(G_n(\R)) \times (\gC^{(m)}-\gC^{(m)}_{=0}) /\R^\times_{>0}) \right)^\ord \]
denotes the image of 
\[ \left( G_n^{(m)}(\A^\infty) \times \pi_0(G_n(\R)) \times \gC^{(m)} \right)^\ord - \left( G_n^{(m)}(\A^\infty) \times \pi_0(G_n(\R)) \times \gC^{(m)}_{=0} \right)^\ord \]
in
\[ G_n^{(m)}(\Q) \backslash \left( (G_n^{(m)}(\A^\infty)/U^p(N,N)) \times \pi_0(G_n(\R)) \times (\gC^{(m)}-\gC^{(m)}_{=0}) /\R^\times_{>0}) \right) .\]
Moreover
\[  \begin{array}{ll} & |\cS(U^p(N),\Sigma)^\ord|_{=i} \\ \cong & 
G_n^{(m)}(\Q) \backslash \left( (G_n^{(m)}(\A^\infty)/U^p(N,N)) \times \pi_0(G_n(\R)) \times ((\gC^{(m)}_{=i})/\R^\times_{>0}) \right)^\ord
\\ \cong & \coprod_{h \in P_{n,(i)}^{(m),+}(\A^\infty)^{\ord,\times} \backslash G_n^{(m)}(\A^\infty)^{\ord,\times}/U^p(N)} L_{n,(i)}^{(m)}(\Q) \backslash L_{n,(i)}^{(m)}(\A) /  \\ \multicolumn{2}{r}{ (hU^p(N) h^{-1} \cap P_{n,(i)}^{(m),+}(\A^\infty)^\ord)L^-_{n,(i),\herm}(\Z_p)L_{n,(i),\herm}(\R)^+ (L_{n,(i),\lin}^{(m)}(\R) \cap U_{n,\infty}^0) } \\ \multicolumn{2}{r}{A_{n,(i)}(\R)^0 .} \end{array}\]
(Use the same argument as in the proof of lemma \ref{orddeco}.) In particular
\[ |\cS(U^p(N),\Sigma)^\ord|_{=n} \cong \gT^{(m),\ord}_{U^p(N),=n}. \]

\newpage \subsection{Toroidal compactifications.}\label{torcomp}

If $(U,\Sigma) \in \cJ^{(m),\tor}_n$,
then there is a smooth projective scheme $A^{(m)}_{n,U,\Sigma}$ and a dense open embedding
\[ j_{U,\Sigma}^{(m)}: A^{(m)}_{n,U} \into A^{(m)}_{n,U,\Sigma} \]
and a projection
\[ \pi_{A^{(m),\tor}/X^\mini}: A^{(m)}_{n,U,\Sigma} \lra X_{n,U}^\mini \]
such that
\[ \begin{array}{ccc} A_{n,U}^{(m)} & \into & A^{(m)}_{n,U,\Sigma}  \\ \da && \da \\  X_{n,U}  & \into &  X_{n,U}^\mini \end{array} \]
is a commutative pull-back square and such that
\[ \partial A^{(m)}_{n,U,\Sigma} = A^{(m)}_{n,U,\Sigma}-j_{U,\Sigma}^{(m)}A^{(m)}_{n,U} \]
is a divisor with simple normal crossings. This induces a log structure $\cM_{\Sigma}$ on $A^{(m)}_{n,U,\Sigma}$.

If $(U,\Sigma) \in \cJ^{(m),\tor}_n$ and $(U',\Delta) \in \cJ^\tor_n$ with $(U,\Sigma) \geq (U',\Delta)$ then there is a log smooth map
\[ \pi_{A^{(m),\tor}/X^\tor}: (A^{(m)}_{n,U,\Sigma},\cM_\Sigma) \lra (X_{n,U',\Delta},\cM_\Delta) \]
over $X_{n,u'}^\mini$ extending the map 
\[ \pi_{A^{(m)}/X}: A^{(m)}_{n,U} \lra X_{n,U'}. \]

If $(U',\Sigma')$ and $(U,\Sigma) \in \cJ^{(m),\tor}_n$; if $g \in G_n^{(m)}(\A^\infty)$; if $U' \supset g^{-1}Ug$; and if $\Sigma g$ is a refinement of $\Sigma'$ then the map $g:A^{(m)}_{n,U} \ra A^{(m)}_{n,U'}$ extends to a log smooth morphism
\[ g: (A^{(m)}_{n,U,\Sigma},\cM_\Sigma) \lra (A^{(m)}_{n,U',\Sigma'},\cM_{\Sigma'}). \]
The collection $\{ A^{(m)}_{n,U,\Sigma} \}$ becomes a system of schemes with right $G_n^{(m)}(\A^\infty)$-action, indexed by $\cJ^{(m),\tor}_n$.  The maps $j_{U,\Sigma}^{(m)}$ and $\pi_{A^{(m),\tor}/X^\mini}$ and $\pi_{A^{(m),\tor}/X^\tor}$ are all $G_n^{(m)}(\A^\infty)$-equivariant. 
If $(U,\Sigma) \geq (U',\Sigma')$ we will write $\pi_{(U,\Sigma), (U',\Sigma')}$ for the map $1:A^{(m)}_{n,U,\Sigma}\ra A^{(m)}_{n,U',\Sigma'}$. (See theorem 1.3.3.15 of \cite{kw2} for the assertions of the last three paragraphs.)

Any of the (canonically quasi-isogenous) universal abelian varieties $A^\univ/X_{n,U}$ extend uniquely to semi-abelian varieties $A^\univ_\Delta/X_{n,U,\Delta}$. The quasi-isogenies between the $A^\univ$ extend uniquely to quasi-isogenies between the $A^\univ_\Delta$. If $g \in G_n(\A^\infty)$ and $(U,\Delta) \geq (U',\Delta')g$ then $g^* A_{\Delta'}^\univ$ is one of the $A_{\Delta}^\univ$. (See remarks 1.1.2.1 and 1.3.1.4 of \cite{kw2}.)

We will write $\partial_i A_{n,U,\Sigma}^{(m)}$ for the pre-image under $\pi_{A^{(m),\tor}/X^\mini}$ of $\partial_i X_U^\mini$. We also set
\[ \partial_i^0 A_{n,U,\Sigma}^{(m)}=\partial_i A_{n,U,\Sigma}^{(m)}-\partial_{i+1} A_{n,U,\Sigma}^{(m)}. \]
We will also write $A_{n,U,\Sigma,i}^{(m),\wedge}$ for the formal completion of $A_{n,U,\Sigma}^{(m)}$ along $\partial_{i}^0 A_{n,U,\Sigma}^{(m)}$ and $\cM_{\Sigma,i}^\wedge$ for the log structure induced on $A_{n,U,\Sigma,i}^{(m),\wedge}$ by $\cM_\Sigma$. There are isomorphisms
\[ (A_{n,U,\Sigma,i}^{(m),\wedge},\cM_{\Sigma,i}^\wedge) \cong \coprod_{h \in P_{n,(i)}^{(m),+}(\A^\infty)\backslash G_n^{(m)}(\A^\infty) / U} (T^{(m),\natural,\wedge}_{n,(i),hU h^{-1}\cap P^{(m),+}_{n,(i)}(\A^\infty), \Sigma(h)_0}, \cM_{\Sigma(h)_0}^\wedge). \]
Suppose that $g^{-1} U g \subset U'$ and that $\Sigma g$ is a refinement of $\Sigma'$. Suppose also that
$h , h' \in G^{(m)}_n(\A^\infty)$ with
\[ h g (h')^{-1} \in P^{(m),+}_{n,(i)}(\A^\infty). \]
Then the diagram 
\[ \begin{array}{ccc} T^{(m),\natural,\wedge}_{n,(i),h U h^{-1} \cap P^{(m),+}_{n,(i)}(\A^\infty), \Sigma(h)_0} & \stackrel{hg (h')^{-1}}{\lra} & T^{(m),\natural,\wedge}_{n,(i),h' U' (h')^{-1} \cap P^{(m),+}_{n,(i)}(\A^\infty), \Sigma'(h')_0} \\ \da && \da \\ A^{(m),\wedge}_{n,U,\Sigma,i} &\stackrel{g}{\lra} & A^{(m),\wedge}_{n,U',\Sigma',i} \end{array} \]
commutes, and is compatible with the log structures on each of these formal schemes. (See theorem 1.3.3.15 of \cite{kw2}.)

If $U'$ is a neat subgroup of $G_n(\A^\infty)$ containing the image of $U$, if $(U',\Delta) \in \cJ^{\tor}_n$ and if $\Sigma$ and $\Delta$ are compatible; then for all $h \in P^{(m),+}_{n,(i)}(\A^\infty)$ with image $h' \in P_{n,(i)}^+(\A^\infty)$ the cone decompositions $\Sigma(h)_0$ and $\Delta(h')_0$ are compatible we have a diagram 
\[ \begin{array}{ccc} T^{(m),\natural,\wedge}_{n,(i),h U h^{-1} \cap P^{(m),+}_{n,(i)}(\A^\infty), \Sigma(h)_0} & \into & A_{n,U,\Sigma,i}^{(m),\wedge} \\ \da && \da \\ T^{\natural,\wedge}_{n,(i),h' U (h')^{-1} \cap P_{n,(i)}^+(\A^\infty), \Delta(h')_0} & \into & X_{n,U',\Delta,i}^{\wedge} \\ \da && \da \\ X^\natural_{n,(i),h' U (h')^{-1} \cap P_{n,(i)}^+(\A^\infty)} & \into & X^{\mini,\wedge}_{n,U',i}, \end{array} \]
which is commutative as a diagram of topological spaces (but not as a diagram of locally ringed spaces). The top square is commutative as a diagram of formal schemes and is compatible with the log structures. (Again see theorem 1.3.3.15 of \cite{kw2}.)

The pull back of $A_\Delta^\univ$ from $X^\wedge_{n,U',\Delta,i}$ to $T^{+,\wedge}_{n,(i),h' U (h')^{-1} \cap P_{n,(i)}^+(\A^\infty), \Delta(h')_0}$ is canonically isogenous to the pull back of $\tG^\univ$ from $A^{+,\wedge}_{n,(i),h' U (h')^{-1} \cap P_{n,(i)}^+(\A^\infty), \Delta(h')_0}$.

We will write
\[ |\cS(\partial A_{n,U,\Sigma}^{(m)})|_{=i} = |\cS(\partial A_{n,U,\Sigma}^{(m)} -\partial_{i+1} A_{n,U,\Sigma}^{(m)})|-  |\cS(\partial A_{n,U,\Sigma}^{(m)} -\partial_{i} A_{n,U,\Sigma}^{(m)})|. \]
Then there are compatible identifications
\[ \cS(\partial A_{n,U,\Sigma}^{(m)}) \cong \cS(U,\Sigma) \]
and
\[ \cS(\partial A_{n,U,\Sigma}^{(m)} - \partial_{i+1} A_{n,U,\Sigma}^{(m)}) \cong \cS(U,\Sigma)_{\leq i} \]
and
\[  |\cS(\partial A_{n,U,\Sigma}^{(m)})|_{=i} \cong |\cS(U,\Sigma)|_{=i}; \]
and the latter is compatible with the identifications
\[ \begin{array}{rl} & |\cS(\partial A_{n,U,\Sigma}^{(m)})|_{=i} \\

\cong & \coprod_{h \in P_{n,(i)}^{(m),+}(\A^\infty)\backslash G_n^{(m)}(\A^\infty)/U}  L_{n,(i),\lin}^{(m)}(\Q) \backslash \\ 

 \multicolumn{2}{l}{  \left( |\cS(\partial T^{(m),+}_{n,(i), hUh^{-1}\cap P_{n,(i)}^{(m),+}(\A^\infty), \widetilde{\Sigma(h)}_0})|-|\cS(\partial T^{(m),+}_{n,(i), hUh^{-1}\cap P_{n,(i)}^{(m),+}(\A^\infty), \widetilde{\Sigma(h)}_0 - \Sigma(h)_0})|\right)} \\

\cong  & \coprod_{h \in P_{n,(i)}^{(m),+}(\A^\infty) \backslash G_n^{(m)}(\A^\infty)/U} L_{n,(i)}^{(m)}(\Q) \backslash L_{n,(i)}^{(m)}(\A) / \\ & (hUh^{-1} \cap P_{n,(i)}^{(m),+}(\A^\infty))L_{n,(i),\herm}(\R)^+ (L_{n,(i),\lin}^{(m)}(\R) \cap U_{n,\infty}^0)A_{n,(i)}(\R)^0 \\

\cong & | \cS(U,\Sigma)|_{=i}. \end{array} \]
(See theorem 1.3.3.15 of \cite{kw2}.)
If $[\sigma] \in \cS(U,\Sigma)$ we will write
\[ \partial_{[\sigma]} A_{n,U,\Sigma}^{(m)} \]
for the corresponding closed boundary stratum of $A^{(m)}_{n,U,\Sigma}$.

Similarly if $(U^p(N_1,N_2),\Sigma) \in \cJ^{(m),\tor}_n$, then there is a smooth quasi-projective scheme $\cA^{(m),\ord}_{n,U^p(N_1,N_2),\Sigma}$ and a dense open embedding
\[ j_{U^p(N_1,N_2),\Sigma}^{(m),\ord}: \cA^{(m),\ord}_{n,U^p(N_1,N_2)} \into \cA^{(m),\ord}_{n,U^p(N_1,N_2),\Sigma} \]
and a projection
\[ \pi_{\cA^{(m),\ord,\tor}/\cX^{\ord,\mini}}: \cA^{(m),\ord}_{n,U^p(N_1,N_2),\Sigma} \lra \cX_{n,(U^p)'(N_1,N_2)}^{\ord,\mini} \]
such that
\[ \begin{array}{ccc} \cA_{n,U^p(N_1,N_2)}^{(m),\ord} & \into & \cA^{(m),\ord}_{n,U^p(N_1,N_2),\Sigma}  \\ \da && \da \\  \cX_{n,(U^p)'(N_1,N_2)}^\ord  & \into &  \cX_{n,(U^p)'(N_1,N_2)}^{\ord,\mini} \end{array} \]
is a commutative pull-back square and such that
\[ \partial \cA^{(m),\ord}_{n,U^p(N_1,N_2),\Sigma} = \cA^{(m),\ord}_{n,U^p(N_1,N_2),\Sigma}-j_{U^p(N_1,N_2),\Sigma}^{(m),\ord}\cA^{(m),\ord}_{n,U^p(N_1,N_2)} \]
is a divisor with simple normal crossings. This induces a log structure $\cM_\Sigma$ on $\cA^{(m),\ord}_{n,U^p(N_1,N_2),\Sigma}$.

If $(U^p(N_1,N_2),\Sigma) \in \cJ^{(m),\tor}_n$ and $((U^p)'(N_1,N_2),\Delta) \in \cJ^\tor_n$ satisfy 
\[ (U^p(N_1,N_2),\Sigma) \geq ((U^p)'(N_1,N_2),\Delta) \]
 then there is a log smooth map
\[ \pi_{\cA^{(m),\ord,\tor}/\cX^{\ord,\tor}}: (\cA^{(m),\ord}_{n,U^p(N_1,N_2),\Sigma} ,\cM_\Sigma) \lra (\cX_{n,(U^p)'(N_1,N_2),\Delta}^\ord,\cM_\Delta) \]
over $\cX_{n,(U^p)'(N_1,N_2)}^{\ord,\mini}$ extending the map 
\[ \pi_{\cA^{(m),\ord}/X^\ord}: \cA^{(m),\ord}_{n,U^p(N_1,N_2)} \lra \cX_{n,(U^p)'(N_1,N_2)}^\ord. \]

If $((U^p)'(N_1,N_2),\Sigma')$ and $(U^p)'(N_1,N_2),\Sigma) \in \cJ^{(m),\tor}_n$; if $g \in G_n^{(m)}(\A^\infty)^\ord$; if $(U^p)'(N_1',N_2') \supset g^{-1}U^p(N_1,N_2)g$; and if $\Sigma g$ is a refinement of $\Sigma'$ then the map $g:\cA^{(m),\ord}_{n,U^p(N_1,N_2)} \ra \cA^{(m),\ord}_{n,(U^p)'(N_1,N_2)}$ extends to a log smooth morphism
\[ g: (\cA^{(m),\ord}_{n,U^p(N_1,N_2),\Sigma},\cM_\Sigma) \lra (\cA^{(m),\ord}_{n,(U^p)'(N_1,N_2),\Sigma'},\cM_{\Sigma'}). \]
Then $\{ \cA^{(m),\ord}_{n,U^p(N_1,N_2),\Sigma} \}$ is a system of schemes with right $G_n^{(m)}(\A^\infty)^\ord$-action, indexed by the subset of $\cJ^{(m),\tor}_n$ consisting of elements of the form $(U^p(N_1,N_2),\Sigma)$.  
The maps $j_{U,\Sigma}^{(m),\ord}$ and $\pi_{\cA^{(m),\ord,\tor}/\cX^{\ord,\mini}}$ and $\pi_{\cA^{(m),\ord,\tor}/\cX^{\ord,\tor}}$ are $G_n^{(m)}(\A^\infty)^\ord$ equivariant.
If $(U^p(N_1,N_2),\Sigma) \geq ((U^p)'(N_1',N_2'),\Sigma')$, then we will denote the map  $1:\cA^{(m),\ord}_{n,U^p(N_1,N_2),\Sigma}\ra \cA^{(m),\ord}_{n,(U^p)'(N_1,N_2),\Sigma'}$ by $\pi_{(U^p(N_1,N_2),\Sigma), ((U^p)'(N_1,N_2),\Sigma')}$. (See theorem 7.1.4.1 of \cite{kw2} for the assertions of the last three paragraphs.)

Any of the (canonically prime-to-$p$ quasi-isogenous) universal abelian varieties $\cA^\univ/\cX^\ord_{n,U^p(N_1,N_2)}$ extend uniquely to semi-abelian varieties $\cA^\univ_\Delta/\cX_{n,U^p(N_1,N_2),\Delta}$. The prime-to-$p$ quasi-isogenies between the $\cA^\univ$ extend uniquely to prime-to-$p$ quasi-isogenies between the $\cA^\univ_\Delta$. If $g \in G_n(\A^\infty)^{\ord,\times}$ and $(U^p(N_1,N_2),\Delta) \geq ((U^p)'(N_1,N_2),\Delta')g$ then $g^* \cA_{\Delta'}^\univ$ is one of the $\cA_{\Delta}^\univ$. (See remarks 3.4.2.8 and 5.2.1.5 of \cite{kw2}.)

We will write $\partial_i \cA_{n,U^p(N_1,N_2),\Sigma}^{(m),\ord}$ for the pre-image under $\pi_{\cA^{(m),\ord,\tor}/\cX^{\ord,\mini}}$ of $\partial_i \cX_{n,U^p(N_1,N_2)}^{\ord,\mini}$ and set
\[ \partial_i^0 \cA_{n,U^p(N_1,N_2),\Sigma}^{(m),\ord} = \partial_i^0\cX_{n,U^p(N_1,N_2)}^{\ord,\mini} - \partial_{i-1}^0\cX_{n,U^p(N_1,N_2)}^{\ord,\mini}.\] 
We will also write $\cA_{n,U^p(N_1,N_2),\Sigma,i}^{(m),\ord,\wedge}$ for the formal completion of $\cA_{n,U^p(N_1,N_2),\Sigma}^{(m),\ord}$ along $\partial_{i}^0 \cA_{n,U^p(N_1,N_2),\Sigma}^{(m),\ord}$, and $\cM_{\Sigma,i}^\wedge$ for the log structure induced on $\cA_{n,U^p(N_1,N_2),\Sigma,i}^{(m),\ord,\wedge}$ by $\cM_\Sigma$. There are isomorphisms
\[ \begin{array}{l} (\cA_{n,U^p(N_1,N_2),\Sigma,i}^{(m),\ord,\wedge},\cM_{\Sigma,i}^\wedge) \cong 
\coprod_{h \in P_{n,(i)}^{(m),+}(\A^\infty)^{\ord,\times} \backslash G_n^{(m)}(\A^\infty)^{\ord,\times} / U^p(N_1)} \\ 
\multicolumn{1}{r}{(\cT^{(m),\ord,\natural,\wedge}_{n,(i),(hU^p h^{-1}\cap P^{(m),+}_{n,(i)}(\A^{p,\infty}))(N_1,N_2), \Sigma^\ord(h)_0},\cM^\wedge_{\Sigma^\ord(h)_0}) } \\
\amalg \coprod_{h \in (P_{n,(i)}^{(m),+}(\A^\infty)\backslash G_n^{(m)}(\A^\infty) / U^p(N_1,N_2))-(P_{n,(i)}^{(m),+}(\A^\infty)^{\ord,\times} \backslash G_n^{(m)}(\A^\infty)^{\ord,\times} / U^p(N_1))} \\
\multicolumn{1}{r}{(T^{(m),\natural,\wedge}_{n,(i),hU^p(N_1,N_2) h^{-1}\cap P^{(m),+}_{n,(i)}(\A^\infty), \Sigma(h)_0},\cM_{\Sigma(h)_0}^\wedge). } \end{array} \]
Suppose that $g \in G_n^{(m)}(\A^\infty)^\ord$ and $g^{-1} U^p(N_1,N_2) g \subset (U^p)'(N_1',N_2')$ and that $\Sigma g$ is a refinement of $\Sigma'$. Suppose also that
$h , h' \in G^{(m)}_n(\A^\infty)^\ord$ with
\[ h g (h')^{-1} \in P^{(m),+}_{n,(i)}(\A^\infty)^\ord. \]
Then the diagram 
\[ \begin{array}{ccc} \cT^{(m),\ord,\natural,\wedge}_{n,(i),V, \Sigma(h)^\ord} & \stackrel{hg (h')^{-1}}{\lra} & \cT^{(m),\ord,\natural,\wedge}_{n,(i),V', \Sigma'(h')^\ord} \\ \da && \da \\ \cA^{(m),\ord,\wedge}_{n,U^p(N_1,N_2),\Sigma,i} &\stackrel{g}{\lra} & \cA^{(m),\ord,\wedge}_{n,(U^p)'(N_1',N_2'),\Sigma',i} \end{array} \]
commutes, where 
\[ V= (h U^p h^{-1} \cap P^{(m),+}_{n,(i)}(\A^{p,\infty}))(N_1,N_2) \]
and
\[ V'= (h' (U^p)' (h')^{-1} \cap P^{(m),+}_{n,(i)}(\A^{p,\infty}))(N_1',N_2'). \]
Moreover this is compatible with the log structures defined on each of the four formal schemes. (See theorem 7.1.4.1 of \cite{kw2}.)

If $[\sigma]\in \cS(U^p(N_1,N_2),\Sigma)$ we will write
\[ \partial_{[\sigma]} \cA^{(m),\ord}_{n,U^p(N_1,N_2),\Sigma} \]
for the closure of $\partial_{[\sigma]} A^{(m)}_{n,U^p(N_1,N_2),\Sigma}$ in $\cA^{(m),\ord}_{n,U^p(N_1,N_2),\Sigma}$. The special fibre
\[ ( \partial_{[\sigma]} \cA^{(m),\ord}_{n,U^p(N_1,N_2),\Sigma} ) \times \Spec \F_p \]
is non-empty if and only if $[\sigma] \in \cS(U^p(N_1),\Sigma^\ord)^\ord$. (We remind the reader that the first superscript ${}^\ord$ associates the `ordinary' cone decomposition $\Sigma^\ord$ to the cone decomposition $\Sigma$, while the second superscript ${}^\ord$ is the notation we are using for the simplicial complex associated to an `ordinary' cone decomposition.)
We will write 
\[ (\cA^{(m),\ord}_{n,U^p(N_1,N_2),\Sigma})^0=\cA^{(m),\ord}_{n,U^p(N_1,N_2),\Sigma} - \!\!\!\! \bigcup_{[\sigma] \in \cS(U^p(N_1,N_2),\Sigma) - \cS(U^p(N_1),\Sigma^\ord)^\ord} \partial_{[\sigma]} \cA^{(m),\ord}_{n,U^p(N_1,N_2),\Sigma}. \]
This only depends on $\Sigma^\ord$.

If $(U^p)'$ is a neat subgroup of $G_n(\A^{p,\infty})$ containing the image of $U^p$, and if $((U^p)'(N_1,N_2),\Delta) \in \cJ_n^{\tor}$, and if $\Sigma$ and $\Delta$ are compatible; then for all $h \in P^{(m),+}_{n,(i)}(\A^\infty)^\ord$ with image $h' \in P_{n,(i)}^+(\A^\infty)^\ord$ the cone decompositions $\Sigma^\ord(h)_0$ and $\Delta^\ord(h')_0$ are compatible and we have a diagram
\[ \begin{array}{ccc} \cT^{(m),\ord,\natural,\wedge}_{n,(i),(h U^p h^{-1} \cap P^{(m),+}_{n,(i)}(\A^{p,\infty}))(N_1,N_2), \Sigma^\ord(h)_0} & \into & \cA_{n,U^p(N_1,N_2),\Sigma,i}^{(m),\ord,\wedge} \\ \da && \da \\ \cT^{\ord,\natural,\wedge}_{n,(i),(h' (U^p)' (h')^{-1} \cap P_{n,(i)}^+(\A^{p,\infty}))(N_1,N_2), \Delta^\ord(h')_0} & \into & \cX_{n,(U^p)'(N_1,N_2),\Delta,i}^{\ord,\wedge} \\ \da && \da \\ \cX^{\ord,\natural}_{n,(i),(h' (U^p)' (h')^{-1} \cap P_{n,(i)}^+(\A^{p,\infty}))(N_1,N_2)} & \into & \cX^{\ord,\mini,\wedge}_{n,(U^p)'(N_1,N_2),i}, \end{array} \]
which is commutative as a diagram of topological spaces (but not as a diagram of locally ringed spaces). The top square is commutative as a diagram of formal schemes and is compatible with the log structures. (See theorem 7.1.4.1 of \cite{kw2}.)

The pull back of $\cA^\univ_\Delta$ to $\cT^{\ord,+,\wedge}_{n,(i),(h' (U^p)' (h')^{-1} \cap P_{n,(i)}^+(\A^{p,\infty}))(N_1,N_2), \Delta^\ord(h')_0}$ is canonically quasi-isomorphic to the pull back of $\tcG^\univ$ from 
\[ \cA^{\ord,+,\wedge}_{n,(i),(h' (U^p)' (h')^{-1} \cap P_{n,(i)}^+(\A^{p,\infty}))(N_1,N_2), \Delta^\ord(h')_0}. \]

All this is compatible with passage to the generic fibre and our previous discussion. (Again see theorem 7.1.4.1 of \cite{kw2}.)

If $N_2' \geq N_2\geq N_1$, if $\Sigma'$ is a refinement of $\Sigma$ and if $\Sigma^\ord =(\Sigma')^\ord$ then the natural map
\[ \cA^{(m),\ord}_{n,U^p(N_1,N_2'),\Sigma'} \lra \cA^{(m),\ord}_{n,U^p(N_1,N_2),\Sigma} \]
is etale in a neighbourhood of the $\F_p$-fibre of $\cA^{(m),\ord}_{n,U^p(N_1,N_2'),\Sigma'}$ and induces an isomorphism between the formal completions of these schemes along their $\F_p$-fibres. (See theorem 7.1.4.1(4) of \cite{kw2}.)
We will denote this $p$-adic formal scheme
\[ \gA^{(m),\ord}_{n,U^p(N_1),\Sigma^\ord} \]
and will denote its reduced subscheme
\[ \barA^{(m),\ord}_{n,U^p(N_1),\Sigma^\ord}. \]
We will also write 
\[ \partial \barA_{n,U^p(N_1),\Sigma^\ord}^{(m),\ord} = \barA_{n,U^p(N_1),\Sigma^\ord}^{(m),\ord}-\barA^{(m),\ord}_{n,U^p(N_1)}. \]
The family $\{ \gA_{n,U^p(N),\Sigma^\ord}^{(m),\ord} \}$ (resp. $\{ \barA_{n,U^p(N),\Sigma^\ord}^{(m),\ord} \}$, resp. $\{ \partial \barA_{n,U^p(N),\Sigma^\ord}^{(m),\ord} \}$) is a family of formal schemes (resp. schemes, resp. schemes) indexed by $\cJ^{(m),\tor,\ord}$ with $G_n(\A^\infty)^\ord$ action. Let
\[ \partial_i \barA_{n,U^p(N),\Sigma^\ord}^{(m),\ord} \]
denote the pre-image of $\partial_i\barX_{n,U^p(N)}^{\ord,\mini}$ in $\partial \barA_{n,U^p(N),\Sigma^\ord}^{(m),\ord}$, and set
\[ \partial_i^0 \barA_{n,U^p(N),\Sigma^\ord}^{(m),\ord}=\partial_i \barA_{n,U^p(N),\Sigma^\ord}^{(m),\ord}-\partial_{i+1} \barA_{n,U^p(N),\Sigma^\ord}^{(m),\ord}. \]
The families $\{ \partial_i \barA_{n,U^p(N),\Sigma^\ord}^{(m),\ord} \}$ and $\{ \partial_i^0\barA_{n,U^p(N),\Sigma^\ord}^{(m),\ord} \}$ are families of schemes with $G_n(\A^\infty)^\ord$ action. 
Moreover we have a decomposition 
\[ \begin{array}{l} \partial_i^0 \barA_{n,U^p(N),\Sigma^\ord}^{(m),\ord} = \\ \coprod_{h \in P^+_{n,(i)}(\A^\infty)^\ord \backslash G_n(\A^\infty)^\ord /U^p(N)} 
\partial_{\Sigma^\ord(h)_0} \barT^{(m),\ord,\natural}_{n,(i),(hU^ph^{-1} \cap P_{n,(i)}^+(\A^{p,\infty}))(N)} . \end{array} \]
If $g \in G_n^{(m)}(\A^\infty)^\ord$, if $g^{-1} U^p(N) g \subset (U^p)'(N')$ and if $\Sigma^\ord g$ is a refinement of $(\Sigma')^\ord$, then the map 
\[ g: \partial_i^0 \barA_{n,U^p(N),\Sigma^\ord}^{(m),\ord} \lra \partial_i^0 \barA_{n,(U^p)'(N'),(\Sigma')^\ord}^{(m),\ord} \]
is the coproduct of the maps
\[ \begin{array}{l} g': \partial_{\Sigma^\ord(h)_0} \barT^{(m),\ord,\natural}_{n,(i),(hU^ph^{-1} \cap P_{n,(i)}^+(\A^{p,\infty}))(N)} \lra \\\partial_{\Sigma^\ord(h)_0} \barT^{(m),\ord,\natural}_{n,(i),(h'(U')^p(h')^{-1} \cap P_{n,(i)}^+(\A^{p,\infty}))(N')} \end{array} \]
where $hg=g'h'$ with $g' \in P_{n,(i)}^+(\A^\infty)^\ord$. 

The map
\[ \varsigma_p: \gA^{(m),\ord}_{n,U^p(N),\Sigma^\ord} \lra  \gA^{(m),\ord}_{n,U^p(N),\Sigma^\ord} \]
is finite flat of degree $p^{(2m+n)n[F^+:\Q]}$ and on $\F_p$-fibres it is identified with absolute Frobenius.

If $N_2' \geq N_2\geq N_1$, if $\Sigma'$ is a refinement of $\Sigma$ and if $\sigma \in \Sigma^\ord =(\Sigma')^\ord$ then the natural map
\[ \partial_{[\sigma]} \cA^{(m),\ord}_{n,U^p(N_1,N_2'),\Sigma'} \lra \partial_{[\sigma]} \cA^{(m),\ord}_{n,U^p(N_1,N_2),\Sigma} \]
is etale in a neighbourhood of the $\F_p$-fibre of $\partial_{[\sigma]} \cA^{(m),\ord}_{n,U^p(N_1,N_2'),\Sigma'}$ and so induces an isomorphism of the formal completions of these schemes along their $\F_p$-fibres.
We will denote this $p$-adic formal scheme
\[ \partial_{[\sigma]} \gA^{(m),\ord}_{n,U^p(N_1),\Sigma^\ord} \]
and will denote its reduced subscheme
\[ \partial_{[\sigma]} \barA^{(m),\ord}_{n,U^p(N_1),\Sigma^\ord}. \]
We will write
\[ \partial^{(s)} \gA^{(m),\ord}_{n,U^p(N_1),\Sigma^\ord} = \coprod_{\substack{[\sigma] \in\cS(U^p(N_1),\Sigma^\ord)^\ord \\ \dim \sigma = s-1}} \partial_{[\sigma]} \gA^{(m),\ord}_{n,U^p(N_1),\Sigma^\ord} \]
and
\[ \partial^{(s)} \barA^{(m),\ord}_{n,U^p(N_1),\Sigma^\ord} = \coprod_{\substack{[\sigma] \in\cS(U^p(N_1),\Sigma^\ord)^\ord \\ \dim \sigma = s-1}} \partial_{[\sigma]} \barA^{(m),\ord}_{n,U^p(N_1),\Sigma^\ord}. \]
The maps 
\[ \varsigma_p: \partial^{(s)} \gA^{(m),\ord}_{n,U^p(N_1),\Sigma^\ord} \lra \partial^{(s)} \gA^{(m),\ord}_{n,U^p(N_1),\Sigma^\ord} \]
are finite flat of degree $p^{(2m+n)n[F^+:\Q]-s}$. 

Then $\partial \barA_{n,U^p(N_1),\Sigma^\ord}^{(m),\ord}$ is stratified by the $\partial_{[\sigma]} \barA^{(m),\ord}_{n,U^p(N_1),\Sigma^\ord}$ with $[\sigma]$ running over $\cS(U^p(N_1),\Sigma^\ord)^\ord$. If $\sigma \in\Sigma^\ord$ but $\sigma$ is not contained in 
\[ \bigcup_{i<n} (G_n^{(m)}(\A^\infty) \times \pi_0(G_n^{(m)}(\R))  \times \gC^{(m)}_{=i})^\ord\]
then $\partial_{[\sigma]} \barA^{(m),\ord}_{n,U^p(N_1),\Sigma^\ord}$ is irreducible. (Because $\partial_{[\sigma]} \barA^{(m),\ord}_{n,U^p(N_1),\Sigma^\ord}$ is a toric variety over $\F_p$. It is presumably true that $\partial_{[\sigma]} \barA^{(m),\ord}_{n,U^p(N_1),\Sigma^\ord}$ is irreducible for any $\sigma$, but to prove it one would need an irreducibility statement about the special fibre of a Shimura variety. In many cases such a theorem has been proved by Hida in \cite{hidabook}, but not in the full generality in which we are working here.)

We will write
\[ \begin{array}{l} |\cS(\partial \barA_{n,U^p(N),\Sigma^\ord}^{(m),\ord})|_{=i}= \\ |\cS(\partial \barA_{n,U^p(N),\Sigma^\ord}^{(m),\ord}-\partial_{ i+1} \barA_{n,U^p(N),\Sigma^\ord}^{(m),\ord})|- |\cS(\partial \barA_{n,U^p(N),\Sigma^\ord}^{(m),\ord}-\partial_{ i} \barA_{n,U^p(N),\Sigma^\ord}^{(m),\ord})| \end{array} \]
an open subset of $|\cS(\partial \barA_{n,U^p(N),\Sigma^\ord}^{(m),\ord}-\partial_{ i+1} \barA_{n,U^p(N),\Sigma^\ord}^{(m),\ord})|$.
Then there are natural surjections
\[ \cS(\partial \barA_{n,U^p(N),\Sigma^\ord}^{(m),\ord}) \onto \cS(U^p(N),\Sigma^\ord)^\ord \]
which restrict to surjections 
\[ \cS(\partial \barA_{n,U^p(N),\Sigma^\ord}^{(m),\ord}-\partial_{\leq i+1} \barA_{n,U^p(N),\Sigma^\ord}^{(m),\ord}) \onto \cS(U^p(N),\Sigma^\ord)^\ord_{ i}. \]
This gives rise to surjections
\[ |\cS(\partial \barA_{n,U^p(N),\Sigma^\ord}^{(m),\ord})|_{=i} \onto |\cS(U^p(N),\Sigma^\ord)^\ord|_{=i}. \]
In the case $n=i$ this is actually a homeomorphism
\[ |\cS(\partial \barA_{n,U^p(N),\Sigma^\ord}^{(m),\ord})|_{=n} 
\cong  |\cS(U^p(N),\Sigma^\ord)^\ord|_{=n}  \cong  \gT^{(m),\ord}_{U^p(N),=n}.  \]
This is compatible with the identifications
\[ \begin{array}{rl} &  |\cS(\partial \barA_{n,U^p(N),\Sigma^\ord}^{(m),\ord})|_{=n} \\

\cong & \coprod_{h \in P_{n,(n)}^{(m),+}(\A^\infty)^{\ord,\times} \backslash G_n^{(m)}(\A^\infty)^{\ord,\times}/U^p(N)} L_{n,(n),\lin}^{(m)}(\Z_{(p)}) \backslash 
\\ & \left( |\cS(\partial \barT^{(m),\ord,+}_{n,(n), hU^p(N)h^{-1}\cap P_{n,(n)}^{(m),+}(\A^\infty)^\ord, \widetilde{\Sigma^\ord(h)}_0})|- \right. \\ & \multicolumn{1}{r}{\left. |\cS(\partial \barT^{(m),\ord,+}_{n,(n), hU^p(N)h^{-1}\cap P_{n,(n)}^{(m),+}(\A^\infty)^\ord, \widetilde{\Sigma^\ord(h)}_0 - \Sigma^\ord(h)_0})|\right) }\\

\cong  & \coprod_{h \in P_{n,(i)}^{(m),+}(\A^\infty)^{\ord,\times} \backslash G_n^{(m)}(\A^\infty)^{\ord,\times}/U^p(N)} L_{n,(n)}^{(m)}(\Q) \backslash L_{n,(n)}^{(m)}(\A) / 
\\ & (hU^p(N)h^{-1} \cap P_{n,(n)}^{(m),+}(\A^\infty))L_{n,(n),\herm}(\R)^+ (L_{n,(n),\lin}^{(m)}(\R) \cap U_{n,\infty}^0)A_{n,(n)}(\R)^0 \\
\cong & | \cS(U^p(N),\Sigma^\ord)^\ord |_{=n}. \end{array} \]

\newpage \subsection{Vector bundles}

We will write $\cI_{\partial X_{n,U}^\mini}$ (resp. $\cI_{\partial X_{n,U,\Delta}}$, resp. $\cI_{\partial A_{U,\Sigma}^{(m)}}$) for the ideal sheaf in $\cO_{X_{n,U}^\mini}$ (resp. $\cO_{X_{n,U,\Delta}}$, resp. $\cO_{A_{U,\Sigma}^{(m)}}$) defining the boundary $\partial X_{n,U}^\mini$ (resp. ${\partial X_{n,U,\Delta}}$, resp. ${\partial A_{U,\Sigma}^{(m)}}$).

\begin{lem} \label{l101} Suppose that $R_0$ is an irreducible, noetherian $\Q$-algebra.
\begin{enumerate}

\item \label{l1012}If $i>0$ then
\[ R^i\pi_{(U,\Sigma),(U',\Sigma'),*}  \cO_{A^{(m)}_{n,U,\Sigma} \times \Spec R_0} = (0) \]
and
\[ R^i\pi_{(U,\Sigma),(U',\Sigma'),*} \cI_{\partial A^{(m)}_{n,U,\Sigma} \times \Spec R_0} = (0). \]

\item \label{l1013} If $(U,\Sigma) \geq (U',\Sigma')$ and $U$ is a normal subgroup of $U'$, then the natural maps
\[ \cO_{A^{(m)}_{U',\Sigma'} \times \Spec R_0} \lra (\pi_{(U,\Sigma),(U',\Sigma'),*}\cO_{A^{(m)}_{U,\Sigma}\times \Spec R_0})^{U'} \]
and
\[ \cI_{\partial A^{(m)}_{U',\Sigma'}\times \Spec R_0} \lra (\pi_{(U,\Sigma),(U',\Sigma'),*}\cI_{\partial A^{(m)}_{U,\Sigma}\times \Spec R_0})^{U'} \]
are isomorphisms.

\item If $U'$ is the image in $G_n(\A^\infty)$ of $U \subset G_n^{(m)}(\A^\infty)$ and if $\Sigma$ and $\Delta$ are compatible, then
\[ \pi_{A^{(m),\tor}/X^\tor,*} \cO_{A^{(m)}_{n,U,\Sigma}} = \cO_{X_{n,U',\Delta}}. \]

\end{enumerate} \end{lem}

\pfbegin If $\Sigma$ is $U'$ invariant the first two parts follow from lemma \ref{ptctp}. In the general case we choose $(U,\Sigma'') \geq (U,\Sigma)$ with $\Sigma''$ being $U'$-invariant, and apply the cases of the lemma already proved to the pairs $((U,\Sigma''),(U',\Sigma'))$ and $((U,\Sigma''),(U,\Sigma))$.

The third part follows from lemma \ref{rc}. 
\pfend

Similarly we will write $\cI_{\partial \cX_{n,U^p(N_1,N_2)}^{\ord,\mini}}$ (resp. $\cI_{\partial \cX_{n,U^p(N_1,N_2),\Delta}^\ord}$, resp. $\cI_{\partial \cA_{U^p(N_1,N_2),\Sigma}^{(m),\ord}}$) for the ideal sheaf in $\cO_{\cX_{n,U^p(N_1,N_2)}^{\ord,\mini}}$ (resp. $\cO_{\cX_{n,U^p(N_1,N_2),\Delta}^\ord}$, resp. $\cO_{\cA_{U^p(N_1,N_2),\Sigma}^{(m),\ord}}$) defining the boundary $\partial \cX_{n,U^p(N_1,N_2)}^{\ord,\mini}$ (resp. ${\partial \cX_{n,U^p(N_1,N_2),\Delta}^\ord}$, resp. ${\partial \cA_{U^p(N_1,N_2),\Sigma}^{(m),\ord}}$). The next lemma follows from lemmas \ref{ptctp} and \ref{rcord}.

\begin{lem} \label{l101ord} Suppose that $R_0$ is an irreducible, noetherian $\Z_{(p)}$-algebra.\begin{enumerate}

\item If $i>0$ then
\[ R^i\pi_{(U^p(N_1,N_2),\Sigma),((U^p)'(N_1',N_2'),\Sigma'),*}  \cO_{\cA^{(m),\ord}_{n,U^p(N_1,N_2),\Sigma}\times \Spec R_0} = (0) \]
and
\[ R^i\pi_{(U^p(N_1,N_2),\Sigma),((U^p)'(N_1',N_2'),\Sigma'),*} \cI_{\partial \cA^{(m),\ord}_{n,U^p(N_1,N_2),\Sigma} \times \Spec R_0} = (0). \]

\item If $(U^p(N_1,,N_2),\Sigma) \geq ((U^p)'(N_1',N_2),\Sigma')$ and $U^p$ is a normal subgroup of $(U^p)'$, then the natural maps
\[ \cO_{A^{(m)}_{(U^p)'(N_1',N_2'),\Sigma'}\times \Spec R_0} \!\!\! \ra (\pi_{(U^p(N_1,N_2),\Sigma),((U^p)'(N_1',N_2'),\Sigma'),*}\cO_{A^{(m)}_{U^p(N_1,N_2),\Sigma} \times \Spec R_0})^{(U^p)'(N_1')} \]
and
\[ \cI_{\partial A^{(m)}_{(U^p)'(N_1',N_2'),\Sigma'}\times \Spec R_0} \!\!\! \ra ( \pi_{(U^p(N_1,N_2),\Sigma),((U^p)'(N_1',N_2'),\Sigma'),*}\cI_{\partial A^{(m)}_{U^p(N_1,N_2),\Sigma}\times \Spec R_0})^{(U^p)'(N_1')} \]
are isomorphisms.

\item If $(U^p)'$ is the image in $G_n(\A^{p,\infty})$ of $U^p \subset G_n^{(m)}(\A^{p,\infty})$ and if $\Sigma$ and $\Delta$ are compatible, then
\[ \pi_{\cA^{(m),\ord,\tor}/X^{\ord,\tor},*} \cO_{\cA^{(m),\ord}_{n,U^p(N_1,N_2),\Sigma}} = \cO_{\cX_{n,(U^p)'(N_1,N_2),\Delta}}. \]

\end{enumerate} \end{lem}

The pull back by the identity section of $\Omega^1_{A^\univ_\Delta/X_{n,U,\Delta}}$ (resp. $\Omega^1_{\cA^\univ_\Delta/\cX^\ord_{n,U^p(N_1,N_2),\Delta}}$) is a locally free sheaf, which is canonically independent of the choice of $A^\univ$ (resp. $\cA^\univ$). We will denote it $\Omega_{n,U,\Delta}$ (resp. $\Omega^\ord_{n,U^p(N_1,N_2),\Delta}$). If $g \in G_n(\A^\infty)$ (resp. $g \in G_n(\A^\infty)^{\ord,\times}$) and $(U,\Delta)g \geq (U',\Delta')$ (resp. $(U^p(N_1,N_2),\Delta)g \geq ((U^p)'(N_1',N_2'),\Delta')$) then there is a natural isomorphism
\[ g^* \Omega_{n,U',\Delta'} \lra \Omega_{n,U,\Delta} \]
(resp. 
\[ g^* \Omega^\ord_{n,(U^p)'(N_1,N_2),\Delta'} \lra \Omega^\ord_{n,U^p(N_1,N_2),\Delta} ).\]
This gives the inverse system $\{ \Omega_{n,U,\Delta} \}$ (resp. $\! \{ \Omega^\ord_{n,U^p(N_1,N_2),\Delta} \}$) an action of $G_n(\A^\infty)$ (resp. $G_n(\A^\infty)^{\ord,\times}$). 
There is also a natural map
\[ \varsigma_p: \varsigma_p^{*}\Omega^\ord_{n,U^p(N_1,N_2-1),\Delta} \lra
\Omega^\ord_{n,U^p(N_1,N_2),\Delta}. \]
There is a canonical identification 
\[ \Omega_{n,U,\Delta}|_{X_{n,U}} \cong \Omega_{n,U} \]
(resp. 
\[ \Omega^\ord_{n,U^p(N_1,N_2),\Delta}|_{X_{n,U}} \cong \Omega^\ord_{n,U^p(N_1,N_2)}). \]

The pull-back of $\Omega_{n,U,\Delta}$ to $T^{\natural,\wedge}_{n,(i),h U h^{-1} \cap P_{n,(i)}^+(\A^{\infty}), \Delta(h)_0}$ is canonically and equivariantly identified with the sheaf $\tOmega^{\natural}_{n,(i), h U h^{-1} \cap P_{n,(i)}^+(\A^{\infty}),\Delta(h)_0}$.
Similarly the pull-back of $\Omega^\ord_{n,U^p(N_1,N_2),\Delta}$ to $\cT^{\ord,\natural,\wedge}_{n,(i),(h U^p h^{-1} \cap P_{n,(i)}^+(\A^{p,\infty}))(N_1,N_2), \Delta^\ord(h)_0}$ is canonically and equivariantly identified with the sheaf $\tOmega^{\ord,\natural}_{n,(i), (h U^p h^{-1} \cap P_{n,(i)}^+(\A^{p,\infty}))(N_1,N_2),\Delta^\ord(h)_0}$.
(See lemmas 1.3.2.41 and 5.2.4.38 of \cite{kw2}.)

We will write
\[ \Xi_{n,U,\Delta} = \cO_{X_{n,U,\Delta}}(||\nu||) \]
(resp.
\[ \Xi^\ord_{n,U^p(N_1,N_2),\Delta} = \cO_{\cX^\ord_{n,U^p(N_1,N_2),\Delta}}(||\nu||) )\]
for the structure sheaf of $X_{n,U,\Delta}$ (resp. $\cX^\ord_{n,U^p(N_1,N_2),\Delta}$) with the $G_n(\A^\infty)$ (resp. $G_n(\A^\infty)^\ord$) action twisted by $||\nu||$. 
If $g \in G_n(\A^\infty)$ (resp. $g \in G_n(\A^\infty)^{\ord,\times}$) then the maps
\[ g^* \Xi_{n,U,\Delta} \lra \Xi_{n,U',\Delta'} \]
(resp.
\[ g^* \Xi^\ord_{n,U^p(N_1,N_2),\Delta} \lra \Xi^\ord_{n,(U^p)'(N_1',N_2'),\Delta'}) \]
are isomorphisms. 

The pull-back of $\Xi_{n,U,\Delta}$ to $T^{\natural,\wedge}_{n,(i),h U h^{-1} \cap P_{n,(i)}^+(\A^{\infty}), \Delta(h)_0}$ equals the pull back of the sheaf $\Xi^{\natural}_{n,(i), h U h^{-1} \cap P_{n,(i)}^+(\A^{\infty})}$ from $X^{\natural}_{n,(i),h U h^{-1} \cap P_{n,(i)}^+(\A^{\infty})}$. Similarly the pull-back of $\Xi^\ord_{n,U^p(N_1,N_2),\Delta}$ to $\cT^{\ord,\natural,\wedge}_{n,(i),(h U^p h^{-1} \cap P_{n,(i)}^+(\A^{p,\infty}))(N_1,N_2), \Delta^\ord(h)_0}$ is naturally isomorphic to the pull back of $\Xi^{\ord,\natural}_{n,(i), (h U^p h^{-1} \cap P_{n,(i)}^+(\A^{p,\infty}))(N_1,N_2)}$ from $\cX^{\ord,\natural}_{n,(i),(h U^p h^{-1} \cap P_{n,(i)}^+(\A^{p,\infty}))(N_1,N_2)}$ to $\cT^{\ord,\natural,\wedge}_{n,(i),(h U^p h^{-1} \cap P_{n,(i)}^+(\A^{p,\infty}))(N_1,N_2), \Delta^\ord(h)_0}$.

Let $\cE^\can_{U,\Delta}$ (resp. $\cE_{U,\Delta}^{\can,\ord}$) denote the principal
$L_{n,(n)}$-bundle on $X_{n,U,\Delta}$ (resp. $\cX^\ord_{n,U^p(N_1,N_2),\Delta}$) in the
Zariski topology defined by setting, for $W \subset X_{n,U,\Delta}$ (resp. $\cX^\ord_{n,U^p(N_1,N_2),\Delta}$) a Zariski open, 
$\cE^\can_{U,\Delta}(W)$ (resp. $\cE^{\can,\ord}_{U^p(N_1,N_2),\Delta}(W)$) to
be the set of pairs $(\xi_0,\xi_1)$, where
\[ \xi_0: \Xi_{n,U,\Delta}|_W \liso \cO_W \]
(resp.
\[ \xi_0: \Xi^\ord_{n,U^p(N_1,N_2),\Delta}|_W \liso \cO_W) \]
and
\[ \xi_1: \Omega_{n,U,\Delta} \liso \Hom_\Q( V_n/V_{n,(n)}, \cO_W) \]
(resp.
\[ \xi_1: \Omega^\ord_{n,U^p(N_1,N_2),\Delta} \liso \Hom_{\Z} 
(\Lambda_n/\Lambda_{n,(n)}, \cO_W)). \]
We define the $L_{n,(n)}$-action on $\cE^\can_{U,\Delta}$ (resp. $\cE_{U^p(N_1,N_2),\Delta}^{\can,\ord}$) by
\[ h(\xi_0,\xi_1)=(\nu(h)^{-1}\xi_0, (\circ h^{-1}) \circ \xi_1). \] 
The inverse system $\{ \cE^\can_{U,\Delta} \}$ (resp. $\{ \cE_{U^p(N_1,N_2),\Delta}^{\can,\ord} \}$) has an action of $G_n(\A^\infty)$ (resp. $G_n(\A^\infty)^\ord$). 

Suppose that $R_0$ is an irreducible noetherian $\Q$-algebra (resp. $\Z_{(p)}$-algebra) and that $\rho$ is a
representation of $L_{n,(n)}$ on a finite, locally free $R$-module $W_\rho$. We define a locally
free sheaf $\cE^\can_{U,\Delta,\rho}$ (resp. $\cE_{U^p(N_1,N_2),\Delta,\rho}^{\can,\ord}$) over
$X_{n,U,\Delta} \times \Spec R_0$ (resp. $\cX^\ord_{n,U^p(N_1,N_2),\Delta}\times \Spec R_0$) by setting $\cE^\can_{U,\Delta,\rho}(W)$ (resp. $\cE_{U^p(N_1,N_2),\Delta,\rho}^{\can,\ord}(W)$) to be the set of $L_{n,(n)}(\cO_W)$-equivariant
maps of Zariski sheaves of sets
\[ \cE^\can_{U,\Delta}|_W \ra W_\rho \otimes_{R_0} \cO_W\]
(resp. 
\[ \cE^{\can,\ord}_{U^p(N_1,N_2),\Delta}|_W \ra W_\rho \otimes_{R_0} \cO_W).\]
Then $\{ \cE^\can_{U,\Delta,\rho} \}$ (resp. $\{ \cE^{\can,\ord}_{U^p(N_1,N_2),\Delta,\rho} \}$) is a system of locally free sheaves with $G_n(\A^\infty)$-action (resp. $G_n(\A^\infty)^\ord$-action) over the system of schemes 
$\{ X_{n,U,\Delta} \times \Spec R_0\}$ (resp. $\{ \cX^\ord_{n,U^p(N_1,N_2),\Delta} \times \Spec R_0 \} $).

Note that
\[ \cE^\can_{U,\Delta,\Std^\vee} \cong \Omega_{A/X,(U,\Delta)} \]
and
\[ \cE^\can_{U,\Delta,\nu^{-1}} \cong \Xi_{A/X,(U,\Delta)} \]
and
\[ \cE^\can_{U,\Delta,\wedge^{n[F:\Q]} \Std^\vee} \cong \omega_U. \]
Similarly
\[ \cE^{\can,\ord}_{U^p(N_1,N_2),\Delta,\Std^\vee} \cong \Omega^\ord_{\cA^\ord/\cX^\ord,(U^p(N_1,N_2),\Delta)} \]
and
\[ \cE^{\can,\ord}_{U^p(N_1,N_2),\Delta,\nu^{-1}} \cong \Xi^\ord_{\cA^\ord/\cX^\ord,(U^p(N_1,N_2),\Delta)} \]
and
\[ \cE^{\can,\ord}_{U^p(N_1,N_2),\Delta,\wedge^{n[F:\Q]} \Std^\vee} \cong \omega_{U^p(N_1,N_2)}. \]

Also note that the pull back of  $\cE^\can_{U,\Delta,\rho}$ (resp. $\cE_{U^p(N_1,N_2),\Delta,\rho}^{\can,\ord}$) to $X_{n,U} \times \Spec R_0$ (resp. $\cX^\ord_{n,U^p(N_1,N_2)} \times \Spec R_0$) is canonically identified with $\cE_{U,\rho}$ (resp. $\cE_{U^p(N_1,N_2),\rho}^{\ord}$). These identifications are $G_n(\A^\infty)$ (resp. $G_n(\A^\infty)^\ord$) equivariant.

Moreover note that the pull back of $\cE^\can_{U,\Delta,\rho}$ to $T^{\natural,\wedge}_{n,(i),h U h^{-1} \cap P_{n,(i)}^+(\A^{\infty}), \Delta(h)_0}$ is canonically and equivariantly identified with the sheaf $\cE^{\natural}_{n,(i), h U h^{-1} \cap P_{n,(i)}^+(\A^{\infty}),\Delta(h)_0,\rho|_{R_{n,(n),(i)}}}$. Similarly the pull-back of $\cE^{\can,\ord}_{U^p(N_1,N_2),\Delta,\rho}$ to $\cT^{\ord,\natural,\wedge}_{n,(i),(h U^p h^{-1} \cap P_{n,(i)}^+(\A^{p,\infty}))(N_1,N_2), \Delta^\ord(h)_0}$ is canonically and equivariantly identified with 
\[ \cE^{\ord,\natural}_{n,(i), (h U^p h^{-1} \cap P_{n,(i)}^+(\A^{p,\infty}))(N_1,N_2),\Delta^\ord(h)_0,\rho|_{R_{n,(n),(i)}}}. \]

Set 
\[ \cE_{U,\Delta,\rho}^\sub = \cI_{\partial X_{n,U,\Delta}}
\cE_{U,\Delta,\rho} \cong \cI_{\partial X_{n,U,\Delta}} \otimes
\cE_{U,\Delta,\rho} \]
and
\[ \cE_{U^p(N_1,N_2),\Delta,\rho}^{\ord,\sub} = \cI_{\partial \cX^\ord_{n,U^p(N_1,N_2),\Delta}}
\cE_{U^p(N_1,N_2),\Delta,\rho}^{\ord} \cong \cI_{\partial \cX^\ord_{n,U^p(N_1,N_2),\Delta}} \otimes
\cE_{U^p(N_1,N_2),\Delta,\rho}^{\ord} \]
Then $\{ \cE^\sub_{U,\Delta,\rho} \}$ (resp. $\{ \cE^{\ord,\sub}_{U^p(N_1,N_2),\Delta,\rho} \}$) is also a system of locally free sheaves with $G_n(\A^\infty)$-action (resp. $G_n(\A^\infty)^\ord$-action) over
$\{ X_{n,U,\Delta} \times \Spec R_0\}$ (resp. $\{ \cX^\ord_{n,U^p(N_1,N_2),\Delta} \times \Spec R_0\} $). 

\begin{lem}\label{l59} \begin{enumerate}
\item If $g \in G_n(\A^\infty)$ (resp. $G_n(\A^\infty)^{\ord,\times}$) and $g:X_{n,U,\Delta} \ra X_{n,U',\Delta'}$ (resp. $g:\cX_{n,U^p(N_1,N_2),\Delta}^\ord \ra \cX_{n,(U^p)'(N_1',N_2'),\Delta'}^\ord$) then
\[ g^* \cE^\can_{U',\Delta',\rho} \liso \cE^\can_{U,\Delta,\rho} \]
(resp.
\[ g^* \cE^{\can,\ord}_{(U^p)'(N_1',N_2'),\Delta',\rho} \liso \cE^{\can,\ord}_{U^p(N_1,N_2),\Delta,\rho} ).\]

\item If $i>0$ then
\[ R^i\pi_{(U,\Delta),(U',\Delta'),*}  \cE^\can_{U,\Delta,\rho} = (0) \]
and
\[ R^i\pi_{(U,\Delta),(U',\Delta'),*}  \cE_{U,\Delta,\rho}^\sub = (0). \]
Similarly, for $i>0$ we have
\[ R^i\pi_{(U^p(N_1,N_2),\Delta),((U^p)'(N_1',N_2'),\Delta'),*}  \cE_{U^p(N_1,N_2),\Delta,\rho}^{\can,\ord} = (0) \]
and
\[ R^i\pi_{(U^p(N_1,N_2),\Delta),((U^p)'(N_1',N_2'),\Delta'),*}  \cE_{U^p(N_1,N_2),\Delta,\rho}^{\ord,\sub} = (0). \]

\item \[ (\lim_{\ra (U,\Delta)} \pi_{(U,\Delta),(U',\Delta'),*}  \cE^\can_{U,\Delta,\rho})^{U'} = \cE_{U',\Delta',\rho} \]
and
\[ (\lim_{\ra (U,\Delta)} \pi_{(U,\Delta),(U',\Delta'),*}  \cE_{U,\Delta,\rho}^\sub)^{U'} = \cE_{U',\Delta',\rho}^\sub  \]
and
\[ \begin{array}{l} \cE_{(U^p)'(N_1',N_2),\Delta',\rho}^{\can,\ord}= \\ (\lim_{\ra (U^p(N_1,N_2),\Delta)} \pi_{(U^p(N_1,N_2),\Delta),((U^p)'(N_1',N_2),\Delta'),*}  \cE_{U^p(N_1,N_2),\Delta,\rho}^{\ord})^{(U^p)'(N_1')} \end{array}  \]
and
\[ \begin{array}{l} \cE_{(U^p)'(N_1',N_2),\Delta',\rho}^{\ord,\sub} = \\ (\lim_{\ra (U^p(N_1,N_2),\Delta)} \pi_{(U^p(N_1,N_2),\Delta),((U^p)'(N_1',N_2),\Delta'),*}  \cE_{U^p(N_1,N_2),\Delta,\rho}^{\ord,\sub})^{(U^p)'(N_1')}  . \end{array} \]

\end{enumerate} \end{lem}

\pfbegin
the first part follows easily from the corresponding facts for $\Omega_{n,U,\Delta}$ and $\Xi_{n,U,\Delta}$ (resp. $\Omega^\ord_{n,U^p(N_1,N_2),\Delta}$ and $\Xi_{n,U^p(N_1,N_2),\Delta}^\ord$). The second and third parts follow from the first part and parts \ref{l1012} and \ref{l1013} of lemma \ref{l101} (resp. lemma \ref{l101ord}). 
\pfend

We next deduce our first main observation. 
\begin{thm}  \label{ttbb} 
 If $i>0$ and $U$ is neat then $R^i\pi_{X^\tor/X^\mini,*}
\cE_{U,\Delta,\rho}^\sub =(0)$.

Similarly if $i>0$ and $U^p$ is neat then $R^i\pi_{\cX^{\ord,\tor}/\cX^{\ord,\mini},*} \cE_{U^p(N_1,N_2),\Delta,\rho}^{\ord,\sub} =(0)$.
\end{thm}

\pfbegin
The argument is the same in both cases, so we explain the argument only in the first case. Write $X^\wedge_{n,U,\Delta,i,h}$ (resp. $X^{\mini,\wedge}_{n,U,h,\partial_{i}^0 X^\mini_{n,U}}$) for the open and closed subset of $X^\wedge_{n,U,\Delta,i}$ (resp. $X^{\mini,\wedge}_{n,U,\partial_{i}^0 X^\mini_{n,U}}$) corresponding to $T^{\natural,\wedge}_{n,(i),h U h^{-1} \cap P^+_{n,(i)}(\A^\infty), \Delta(h)}$ (resp. $X^\natural_{n,(i),h U h^{-1} \cap P^+_{n,(i)}(\A^\infty)}$). (Recall that $X^\wedge_{n,U,\Delta,i}$ is the completion of a smooth toroidal compactification of the  Shimura variety $X_{n,U}$ along the locally closed subspace of the boundary corresponding to the parabolic subgroup $P^+_{n,(i)} \subset G_n$. The formal scheme $X^{\mini,\wedge}_{n,U,\partial_{i}^0 X^\mini_{n,U}}$ is the completion of the minimal (Baily-Borel) compactification of the same Shimura variety along the locally closed subspace of the boundary corresponding to the same parabolic. Each of these formal schemes is a disjoint union of sub-formal schemes indexed by certain elements $h \in G_n(\A^\infty)$.)

We have maps of locally ringed spaces
\[ \begin{array}{ccc} T^{\natural,\wedge}_{n,(i),h U h^{-1} \cap P^+_{n,(i)}(\A^\infty), \Delta(h)} & \liso & X^\wedge_{n,U,\Delta,i,h} \\ \da && \da \\ X^\natural_{n,(i),h U h^{-1} \cap P^+_{n,(i)}(\A^\infty)} & \into & X^{\mini,\wedge}_{n,U,h,\partial_{i}^0 X^\mini_{n,U}}. \end{array} \]
(Recall that $T^{\natural,\wedge}_{n,(i),h U h^{-1} \cap P^+_{n,(i)}(\A^\infty), \Delta(h)}$ is a formal local model for the boundary of the toroidal compactification. It is the quotient by a discrete group of the formal completion of a toroidal embedding over a principal homogeneous space for an abelian scheme over a disjoint union of smaller Shimura varieties. The scheme $X^\natural_{n,(i),h U h^{-1} \cap P^+_{n,(i)}(\A^\infty)}$ is a disjoint union of smaller Shimura varieties, and also a locally closed subscheme of the boundary of the minimal compactification of $X_{n,U}$.)

This diagram is commutative as a diagram of topological spaces (but not of locally ringed spaces)
and the lower horizontal map is an isomorphism on the underlying topological spaces. 
It suffices to show that the higher direct images from the topological space $T^{\natural,\wedge}_{n,(i),h U h^{-1} \cap P^+_{n,(i)}(\A^\infty), \Delta(h)}$ to the topological space $X^\natural_{n,(i),h U h^{-1} \cap P^+_{n,(i)}(\A^\infty)}$ of the pull-back of 
$\cE_{U,\Delta,\rho}^\sub$ vanishes. The theorem follows on combining the last part of the previous lemma with corollary \ref{subvan5.5}.
\pfend

We set
\[ \cE_{U,\rho}^\sub = \pi_{X^\tor/X^\mini,*} \cE_{U,\Delta,\rho}^\sub \]
(resp.
\[ \cE_{U^p(N_1,N_2),\rho}^{\ord,\sub} = \pi_{\cX^{\ord,\tor}/\cX^{\ord,\mini},*} \cE_{U^p(N_1,N_2),\Delta,\rho}^{\ord,\sub}) \]
a coherent sheaf on $X_{n,U}^\mini \times \Spec R_0$ (resp. $\cX_{n,U^p(N_1,N_2)}^{\ord,\mini} \times \Spec R_0$). These definitions are
independent of $\Delta$. Note that
\[ \cE_{U,\rho}^\sub \otimes \omega_{U}^{\otimes N} \cong \cE_{U,\rho \otimes (\wedge^{n[F:\Q]} \Std^\vee)^{\otimes N}}^\sub \]
and
\[ \cE_{U^p(N_1,N_2),\rho}^{\ord,\sub} \otimes (\omega_{U^p(N_1,N_2)})^{ \otimes N} \cong \cE_{U^p(N_1,N_2),\rho \otimes (\wedge^{n[F:\Q]} \Std^\vee)^{\otimes N}}^{\ord,\sub}. \]

We will let $\cE^{\ord,\can}_{U^p(N),\Delta^\ord,\rho}$ (resp. $\cE^{\ord,\sub}_{U^p(N),\Delta^\ord,\rho}$, resp. $\cE^{\ord,\sub}_{U^p(N),\rho}$) denote the pull-back of $\cE^{\can,\ord}_{U^p(N,N'),\Delta,\rho}$ (resp. $\cE^{\ord,\sub}_{U^p(N,N'),\Delta,\rho}$, resp. $\cE^{\ord,\sub}_{U^p(N,N'),\rho}$) to $\gX^{\ord}_{U^p(N),\Delta^\ord}$ (resp. $\gX^{\ord}_{U^p(N),\Delta^\ord}$, resp. $\gX^{\ord,\mini}_{U^p(N)}$). It is independent of the choice of $N'$ and $\Delta$.

If $\rho$ is a representation of $L_{n,(n)}$ on a finite $\Q$-vector space, we will set
\[ \begin{array}{rcl} H^i(X_n^\mini, \cE_\rho^\sub) &=& \lim_{\substack{\lra \\U'}} H^i(X_{n,U'}^\mini, \cE_{U',\rho}^\sub) \\ &=& \lim_{\substack{\lra \\ U',\Delta} }H^i(X_{n,U',\Delta}, \cE_{U',\Delta,\rho}^\sub). \end{array} \]
It is an admissible $G_n(\A^\infty)$-module with 
\[ H^i(X_n^\mini, \cE_\rho^\sub)^{U'}=H^i(X_{n,U'}^\mini, \cE_{U',\rho}^\sub). \]

Similarly, if $\rho$ is a representation of $L_{n,(n)}$ on a finite free $\Z_{(p)}$-module, we will set
\[ \begin{array}{rl}& H^0(\cX^{\ord,\mini}_n, \cE_\rho^{\ord,\sub} \otimes \Z/p^r\Z) \\
=& \lim_{\substack{\lra \\ U^p,N_1,N_2}} H^0(\cX_{n,U^p(N_1,N_2)}^{\ord,\mini}, \cE_{U^p(N_1,N_2),\rho}^{\ord,\sub} \otimes \Z/p^r\Z) \\ 
=& \lim_{\substack{\lra \\ U^p,N_1,N_2,\Delta}} H^0(\cX^\ord_{n,U^p(N_1,N_2),\Delta}, \cE_{U^p(N_1,N_2),\Delta,\rho}^{\ord,\sub}\otimes \Z/p^r\Z) \end{array} \]
and
\[ \begin{array}{rl}& H^0(\gX^{\ord,\mini}_n, \cE_\rho^{\ord,\sub}) \\
=& \lim_{\substack{\lra \\ U^p,N}} H^0(\gX_{n,U^p(N)}^{\ord,\mini}, \cE_{U^p(N),\rho}^{\ord,\sub}) \\ 
=& \lim_{\substack{\lra \\ U^p,N,\Delta}} H^0(\gX^\ord_{n,U^p(N),\Delta}, \cE_{U^p(N),\Delta,\rho}^{\ord,\sub}) \end{array} \]
They are smooth $G_n(\A^\infty)^\ord$-modules with
\[ H^0(\cX^{\ord,\mini}_n, \cE_\rho^{\ord,\sub} \otimes \Z/p^r\Z)^{U^p(N_1)} = H^0(\cX_{n,U^p(N_1,N_2)}^{\ord,\mini}, \cE_{U^p(N_1,N_2),\rho}^{\ord,\sub} \otimes \Z/p^r\Z)  \]
and
\[ H^0(\gX^{\ord,\mini}_n, \cE_\rho^{\ord,\sub})^{U^p(N)} = H^0(\gX_{n,U^p(N)}^{\ord,\mini}, \cE_{U^p(N),\rho}^{\ord,\sub}).  \]
Note that there is a $G_n(\A^\infty)^\ord$-equivariant embedding
\[ H^0(\gX^{\ord,\mini}_n, \cE_\rho^{\ord,\sub}) \otimes_{\Z_p} \Z/p^r\Z \into H^0(\cX^{\ord,\mini}_n, \cE_\rho^{\ord,\sub} \otimes \Z/p^r\Z). \]
Finally set
\[ H^0(\gX^{\ord,\mini}, \cE_\rho^{\ord,\sub})_{\barQQ_p} = H^0(\gX^{\ord,\mini}, \cE_\rho^{\ord,\sub}) \otimes_{\Z_p} \Q_p, \]
a smooth representation of $G_n(\A^\infty)^\ord$.

We record the following result from \cite{kw2}.
\begin{lem}  If  $\rho$ is a representation of $L_{n,(n)}$ on a finite locally free $\Z_{(p)}$-module then there is a unique system $\{ \cE^\sub_{U^p(N_!,N_2),\rho} \}$ of $\cO_{\cX^\mini_{n,U^p(N_1,N_2)}}$-torsion free coherent sheaves with $G_n(\A^{\infty})^{\ord,\times}$-action over $\{ \cX^\mini_{n,U^p(N_1,N_2)} \}$ with the following properties.
\begin{enumerate}
\item $\{ \cE_{U^p(N_1,N_2),\rho}^\sub \}$ pulls back to $\{ \cE_{U^p(N_1,N_2),\rho\otimes_{\Z_{(p)}} \Q}^\sub \}$ on $\{ X_{n,U^p(N_1,N_2)}^\mini \}$;
\item $\{ \cE_{U^p(N_1,N_2),\rho}^\sub \}$ pulls back to $\{ \cE_{U^p(N_1,N_2),\rho}^{\ord,\sub} \}$ on $\{ \cX_{n,U^p(N_1,N_2)}^{\ord,\mini}\}$;
\item if $U^p$ is a normal subgroup of $(U^p)'$ and if $g \in (U^p)'(N_1',N_2)$ then
\[ g:g^*\cE_{U^p(N_1,N_2),\rho}^\sub \iso \cE_{U^p(N_1,N_2),\rho}^\sub; \]
\item if $U^p$ is a normal subgroup of $(U^p)'$ then 
\[ \cE_{U^p(N_1,N_2),\rho}^\sub \iso (\pi_{(U^p)'(N_1',N_2),U^p(N_1,N_2),*} \cE_{(U^p)'(N_1',N_2),\rho}^\sub)^{U^p(N_1,N_2)};\]

\item $\{ \cE_{U^p(N_1,N_2),\rho \otimes \wedge^{n[F:\Q]} \Std^\vee}^\sub \} \cong \{ \omega_{U^p(N_1,N_2),\Sigma}\otimes \cE_{U^p(N_1,N_2),\rho}^\sub \}$.
\end{enumerate} \end{lem}

\pfbegin
For the definition of $\cE^\sub_{U^p(N_!,N_2),\rho}$ see definition 8.3.5.1 of \cite{kw2}. For the $\cO_{\cX^\mini_{n,U^p(N_1,N_2)}}$-torsion freeness see corollary 8.3.5.8 of \cite{kw2}. For the $G_n(\A^{\infty})^{\ord,\times}$-action see corollary 8.3.6.5 of \cite{kw2}. For part one of the lemma see lemma 8.3.5.2 and corollary 8.3.6.5 of \cite{kw2}. For the second part see corollary 8.3.5.4 of \cite{kw2}. The third part is clear. For the fourth part see proposition 8.3.6.9 of \cite{kw2}, and for the final part see lemma 8.3.5.10 of \cite{kw2}. \pfend

We will write 
$\Omega^i_{A_{n,U,\Sigma}^{(m)}}(\log \infty)$ (resp. $\Omega^i_{A_{n,U,\Sigma}^{(m)}/X_{n,U',\Delta}}(\log \infty)$) as shorthand for 
$\Omega^i_{A_{n,U,\Sigma}^{(m)}/\Spec \Q}(\log \cM_\Sigma)$ (resp. $\Omega^i_{A_{n,U,\Sigma}^{(m)}/X_{n,U',\Delta}}(\log \cM_\Sigma/\cM_\Delta)$).
Then the collection 
$\{ \Omega^1_{A_{n,U,\Sigma}^{(m)}}(\log \infty)\}$ (resp. $\{\Omega^1_{A_{n,U,\Sigma}^{(m)}/X_{n,U',\Delta}}(\log \infty)\}$) is a system of locally free sheaves (for the Zariski topology) with 
$G_n^{(m)}(\A^\infty)$-action. 

There are natural differentials
\[ d:\Omega^i_{A_{n,U,\Sigma}^{(m)}}(\log \infty) \lra \Omega^{i+1}_{A_{n,U,\Sigma}^{(m)}}(\log \infty), \]
(resp.
\[ d:\Omega^i_{A_{n,U,\Sigma}^{(m)}/X_{n,U',\Delta}}(\log \infty) \lra \Omega^{i+1}_{A_{n,U,\Sigma}^{(m)}/X_{n,U',\Delta}}(\log \infty) )\]
making 
$\Omega^\bullet_{A_{n,U,\Sigma}^{(m)}}(\log \infty)$ (resp. $\Omega^\bullet_{A_{n,U,\Sigma}^{(m)}/X_{n,U',\Delta}}(\log \infty)$) a complex. The tensor products 
$\Omega^\bullet_{A_{n,U,\Sigma}^{(m)}}(\log \infty) \otimes \cI_{\partial A_{n,U,\Sigma}^{(m)}}$ (resp. $\Omega^\bullet_{A_{n,U,\Sigma}^{(m)}/X_{n,U',\Delta}}(\log \infty) \otimes \cI_{\partial A_{n,U,\Sigma}^{(m)}}$) is a sub-complex.

\begin{lem} \label{cscor}\begin{enumerate}

\item If $(U,\Sigma) \geq (U',\Delta) \geq (U'',\Delta')$ then the natural morphism $\Omega^1_{A_{n,U,\Sigma}^{(m)}/X_{n,U'',\Delta'}}(\log \infty) \iso \Omega^1_{A_{n,U,\Sigma}^{(m)}/X_{n,U',\Delta}}(\log \infty)$ is an isomorphism, so we will simply write $\Omega^1_{A_{n,U,\Sigma}^{(m)}/X}(\log \infty)$ for this sheaf.

\item \label{cscor3} If $(U',\Sigma')\geq (U,\Sigma)$ then  $\pi_{(U',\Sigma'),(U,\Sigma)}^* \Omega^1_{A^{(m)}_{n,U,\Sigma}}(\log \infty) \iso \Omega^1_{A^{(m)}_{n,U',\Sigma'}}(\log \infty)$ and $\pi_{(U',\Sigma'),(U,\Sigma)}^* \Omega^1_{A^{(m)}_{n,U,\Sigma}/X}(\log \infty) \iso \Omega^1_{A^{(m)}_{n,U',\Sigma'}/X}(\log \infty)$. 

\item\label{cscor4} If $(U,\Sigma) \geq (U',\Delta)$ then there is an exact sequence
\[ (0) \ra \pi_{(U,\Sigma),(U',\Delta)}^* \Omega^1_{X_{n,U',\Delta}}(\log \infty) \ra \Omega^1_{A^{(m)}_{n,U,\Sigma}}(\log \infty) \ra \Omega^1_{A^{(m)}_{n,U,\Sigma}/X}(\log \infty) \ra (0). \]

\item\label{cscor5} Suppose that $(U_1,\Sigma_1) \geq (U_2,\Sigma_2) \geq (U',\Delta)$, and that $U'$ is the image of both $U_1$ and $U_2$ in $G_n(\A^\infty)$. Then the natural maps
\[ R^i\pi_{A^{(m),\tor}/X^\tor,*} \Omega^j_{A_{n,U_2,\Sigma_2}^{(m)}/X}(\log \infty) \lra R^i\pi_{A^{(m),\tor}/X^\tor,*} \Omega^j_{A_{n,U_1,\Sigma_1}^{(m)}/X}(\log \infty) \]
and
\[ \begin{array}{l} R^i\pi_{A^{(m),\tor}/X^\tor,*} (\Omega^j_{A_{n,U_2,\Sigma_2}^{(m)}/X}(\log \infty) \otimes \cI_{\partial A^{(m)}_{n,U_2,\Sigma_2}}) \\ \lra R^i\pi_{A^{(m),\tor}/X^\tor,*} (\Omega^j_{A_{n,U_1,\Sigma_1}^{(m)}/X}(\log \infty) \otimes \cI_{\partial A^{(m)}_{n,U_1,\Sigma_1}}) \end{array} \]
on $X_{n,U',\Delta}$ are isomorphisms. 
We will write simply
\[ (R^i\pi_*\Omega^j_{A^{(m)}/X}(\log \infty))_{(U',\Delta)} \]
and
\[ (R^i\pi_*(\Omega^j_{A^{(m)}/X}(\log \infty) \otimes \cI_{\partial A^{(m)}}))_{(U',\Delta)} \]
for these sheaves.

\item \label{cscor6}  $\{ (R^i\pi_*\Omega^j_{A^{(m)}/X}(\log \infty))_{(U',\Delta)} \}$ and $\{ (R^i\pi_*(\Omega^j_{A^{(m)}/X}(\log \infty) \otimes \cI_{\partial A^{(m)}}))_{(U',\Delta)} \}$ are systems of coherent sheaves with $G_n^{(m)}(\A^\infty)$-action over $\{ X_{n,U',\Delta} \}$. Moreover the maps
\[ g: g^* (R^i\pi_*\Omega^j_{A^{(m)}/X}(\log \infty))_{(U',\Delta)} \lra (R^i\pi_*\Omega^j_{A^{(m)}/X}(\log \infty))_{(U'',\Delta')} \]
are isomorphisms.

\item \label{cscor7}  The $G_n^{(m)}(\A^\infty)$-actions on the systems $\{ (R^i\pi_*\Omega^j_{A^{(m)}/X}(\log \infty))_{(U',\Delta)} \}$ and  $\{ (R^i\pi_*\Omega^j_{A^{(m)}/X}(\log \infty) \otimes \cI_{\partial A^{(m)}})_{(U',\Delta)} \}$ factor through $G_n(\A^\infty)$.

\item \label{locmod6} The pull-back of $(\pi_* \Omega^1_{A^{(m)}/X}(\log \infty))_{(U,\Delta)}$ to $T^{\natural,\wedge}_{n,(i),h U h^{-1} \cap P_{n,(i)}^+(\A^\infty), \Delta(h)_0}$ is isomorphic to
 \[ \pi_{(U',\Sigma_0),(h U h^{-1} \cap P_{n,(i)}^+(\A^\infty), \Delta(h)_0),*} \Omega^1_{T^{(m),\natural,\wedge}_{n,(i),U',\Sigma_0}/T^{\natural,\wedge}_{n,(i),h U h^{-1} \cap P_{n,(i)}^+(\A^\infty), \Delta(h)_0}} (\log \infty) \]
 for some $U'$ and $\Sigma_0$.

\end{enumerate} \end{lem}

\pfbegin
This follows from the properties of log differentials for log smooth maps (see section \ref{logstruct}). For part \ref{cscor5} we also use lemma \ref{l101}. For part \ref{cscor7} we also use the discussion of section \ref{vb1} and a density argument.
\pfend

The next lemma follows from lemma \ref{rc2}. 

\begin{lem}\label{locmod} \begin{enumerate}

\item The natural maps
\[ (\pi_* \Omega^1_{A^{(m)}/X}(\log \infty))_{(U',\Delta)} \otimes_{\cO_{X_{n,U',\Delta}}} \cO_{A^{(m)}_{n,U,\Sigma}} \lra \Omega^1_{A^{(m)}_{n,U,\Sigma} /X}(\log \infty) \]
are $G^{(m)}_n(\A^\infty)$-equivariant isomorphisms.

\item The natural maps
\[  (\wedge^j (\pi_*\Omega^1_{A^{(m)}/X}(\log \infty))_{(U',\Delta)}) \otimes (R^i\pi_*\cO_{A^{(m)}})_{(U',\Delta)} \lra (R^i\pi_{*} \Omega^j_{A^{(m)}/X}(\log \infty))_{(U',\Delta)} \]
and
\[  \begin{array}{l} (\wedge^j (\pi_*\Omega^1_{A^{(m)}/X}(\log \infty))_{(U',\Delta)}) \otimes (R^i\pi_*\cO_{A^{(m)}})_{(U',\Delta)}\otimes \cI_{\partial X_{n,U',\Delta}}\\ \lra (R^i\pi_{*} \Omega^j_{A^{(m)}/X}(\log \infty) \otimes \cI_{\partial A^{(m)}})_{(U',\Delta)} \end{array} \]
are $G_n(\A^\infty)$-equivariant isomorphisms.

\item $(\pi_* \Omega^1_{A^{(m)}/X}(\log \infty))_{(U,\Delta)}$ is a flat coherent $\cO_{X_{n,U,\Delta}}$-module, and hence locally free of finite rank.

\end{enumerate} \end{lem}

Next we record some results of one of us (K.-W.L). 

\begin{lem}\label{kods} \begin{enumerate} 

\item\label{cscor8} There are natural $G_n(\A^\infty)$-equivariant isomorphisms
\[ \Hom_F(F^m , \Omega_{n,U',\Delta}) \liso (\pi_* \Omega^1_{A_n^{(m)}/X_n}(\log \infty))_{(U',\Delta)}. \]

\item\label{wedge} The cup product maps
\[ \wedge^i (R^1\pi_* \cO_{A^{(m)}})_{(U',\Delta)} \lra (R^i\pi_* \cO_{A^{(m)}})_{(U',\Delta)} \]
are $G_n(\A^\infty)$-equivariant isomorphisms.

\item\label{kods2} There is a unique embedding 
\[ \Xi_{(U',\Delta)} \into (R^1\pi_* \Omega^1_{A^{(m)}/X}(\log \infty))_{(U',\Delta)} \]
extending
\[ \Xi_{U'} \into (R^1\pi_* \Omega^1_{A^{(m)}/X})_{U'}. \]
It is $G_n(\A^\infty)$-equivariant.

\item\label{kods5} The composite maps
\[ \begin{array}{rl} & \Hom((\pi_* \Omega^1_{A_n^{(m)}/X_n}(\log \infty))_{(U',\Delta)}, \Xi_{(U',\Delta)}) \\
\lra &\Hom \left( ((\pi_* \Omega^1_{A_n^{(m)}/X_n}(\log \infty))_{(U',\Delta)}, \right.
\\ &\left. (\pi_* \Omega^1_{A_n^{(m)}/X_n}(\log \infty))_{(U',\Delta)} \otimes (R^1\pi_* \cO_{A^{(m)}})_{(U',\Delta)}\right)\\ 
 \stackrel{\tr}{\lra} &(R^1\pi_* \cO_{A^{(m)}})_{(U',\Delta)}\end{array} \]
are $G_n(\A^\infty)$-equivariant isomorphisms.

\item\label{kods6} The boundary maps 
\[ \begin{array}{rcl} \Omega_{(U',\Delta)} &\lra &R^1\pi_{A^{(1)}/X,*} (\pi_{A^{(1)}/X}^* \Omega^1_{X_{n,U',\Delta}}(\log \infty)) \\ & \cong &\Omega^1_{X_{n,U',\Delta}}(\log \infty) \otimes \Hom(\Omega_{(U',\Delta)}, \Xi_{(U',\Delta)}) \end{array} \]
associated to the short exact sequence of part \ref{cscor4} of lemma \ref{cscor}, give rise to isomorphisms
\[ S(\Omega_{(U',\Delta)}) \liso \Omega^1_{X_{n,U',\Delta}}(\log \infty) \otimes \Xi_{(U',\Delta)}. \]

\item There are $G_n(\A^{p,\infty} \times \Z_p)$-equivariant identifications between the pull back of $\omega_U$ from $\cX^\mini_{n,U}$ to $X_{n,U,\Delta}$ and $\wedge^{n[F:\Q]} \Omega_{(U,\Delta)}$.

\end{enumerate}\end{lem} 

\pfbegin
For the first four parts see theorem 2.5 and proposition 6.9 of \cite{kw1.5} and theorem 1.3.3.15 of \cite{kw2}. 
For the fifth part see theorem 6.4.1.1 (4) of \cite{kw1}. For the sixth part see theorem 7.2.4.1 and proposition 7.2.5.1 of \cite{kw1}. 
\pfend

\begin{cor} There are equivariant isomorphisms $\cE_{U,\Delta,\KS}^\can \cong \Omega^1_{X_{n,U,\Delta}}(\log \infty)$.
(See section \ref{mps} for the definition of the representation $\KS$.)
\end{cor}

\begin{lem}\label{spect} Suppose that $U$ is a neat open compact subgroup of $G^{(m)}_n(\A^\infty)$ with image $U'$ in $G_n(\A^\infty)$. Then there are representations $\rho_{m,r}^{i,j}$ of $L_{n,(n)}$ and a spectral sequence of sheaves on $X_{n,U'}^\mini$ with first page
\[ E_1^{i,j}=\cE_{U',\rho_{m,r}^{i,j}}^\sub \Rightarrow R^{i+j}\pi_{A^{(m),\tor}/X^\mini,*} (\Omega^r_{A^{(m)}_{n,U,\Sigma}}(\log \infty)\otimes \cI_{\partial A^{(m)}_{n,U,\Sigma}}). \]
This spectral sequence is $G_n(\A^\infty)$-equivariant. 
\end{lem}

\pfbegin Using part \ref{cscor3} of corollary \ref{cscor} and parts \ref{l1012} and \ref{l1013} of lemma \ref{l101}, we may reduce to the case that there is a cone decomposition 
$\Delta$ compatible with $\Sigma$. By the preceding theorem it suffices to find $\rho^{i,j}_{m,r}$ such that there is a spectral sequence of sheaves on $X_{n,U',\Delta}$ with first page
\[ E_1^{i,j}=\cE_{U',\Delta,\rho_{m,r}^{i,j}}^\sub \Rightarrow R^{i+j}\pi_{A^{(m),\tor}/X^\tor,*} (\Omega^r_{A^{(m)}_{n,U,\Sigma}}(\log \infty)\otimes \cI_{\partial A^{(m)}_{n,U,\Sigma}}). \]
However we may filter $\Omega^r_{A^{(m)}_{n,U,\Sigma}}(\log \infty)\otimes \cI_{\partial A^{(m)}_{n,U,\Sigma}}$ with graded pieces 
\[ (\pi_{A^{(m),\tor}/X^\tor}^* \Omega^j_{X_{n,U',\Delta}}(\log \infty))\otimes 
\Omega^{r-j}_{A_{n,U,\Sigma}^{(m)}/X}(\log \infty) \otimes \cI_{\partial A^{(m)}_{n,U,\Sigma}}. \]
Moreover by lemma \ref{locmod} we have that
\[ \begin{array}{l} (\wedge^j\Omega^1_{X_{n,U',\Delta}}(\log \infty)) \otimes (\wedge^{r-j} (\pi_* \Omega^1_{A^{(m)}/X})_{(U',\Delta)}) \otimes (R^i\pi_{A^{(m),\tor}/X^\tor,*} \cO_{A^{(m)}})_{(U',\Delta)} \\ \multicolumn{1}{r}{\otimes \cI_{\partial X_{n,U',\Delta}}} \\ \liso R^i\pi_{A^{(m),\tor}/X^\tor,*} (\pi_{A^{(m),\tor}/X^\tor}^* \Omega^j_{X_{n,U',\Delta}}(\log \infty))\otimes 
\Omega^{r-j}_{A_{n,U,\Sigma}^{(m)}/X}(\log \infty) \otimes \cI_{\partial A^{(m)}_{n,U,\Sigma}} . \end{array} \]
The result follows on combining this with parts \ref{cscor8}, \ref{wedge}, \ref{kods5} and \ref{kods6} of lemma \ref{kods}.
\pfend

\newpage \subsection{Connection to the complex theory}

\begin{lem}\label{ssauto} Suppose that 
\[ b=(b_0,(b_{\tau,i})_{\tau \in \Hom(F,\C)}) \in X^*(T_n/\C)_{(n)}^+ \]
satisfies 
\[ -2n \geq b_{\tau,1}+b_{\tau c,1} \]
for all $\tau \in \Hom(F,\C)$.
Then $H^0(X^\mini,\cE_{\rho_{(n),b}}^\sub)$ is a semi-simple $G_n(\A^\infty)$-module. If $\pi$ is an irreducible sub-quotient of $H^0(X^\mini,\cE_{\rho_{(n),b}}^\sub)$, then $\pi$ is the finite part of a cohomological, cuspidal automorphic representation of $G_n(\A)$. \end{lem}

\pfbegin
According to proposition 5.4.2 and lemma 5.2.3 of \cite{dbar} and theorems 4.1.1, 5.1.1 and 5.2.12 of \cite{kwan} we have an isomorphism
\[ H^0(X_n^\mini,\cE_{\rho_{(n),b}}^\sub) \cong \bigoplus_\Pi \Pi^\infty \otimes H^0( \gq_n ,U_{n,\infty}^0A_n(\R)^0, \Pi_\infty \otimes \rho_{(n),b})\]
where $\Pi$ runs over cuspidal automorphic representations of $G_n(\A)$ taken with their multiplicity in the space of cuspidal automorphic forms. 

Thus $\pi \cong \Pi^\infty$ for some cuspidal automorphic representation $\Pi$ of $G_n(\A)$ with
\[ H^0( \gq_n ,U_{n,\infty}^0A_n(\R)^0, \Pi_\infty \otimes \rho_{(n),b}) \neq (0). \]
It follows from theorem 2.6 of \cite{co} that the Harish-Chandra parameter of the infinitesimal character of $\Pi_\infty$ equals
\[ \varrho_n - 2 \varrho_{n,(n)} -b. \]
As we have assumed that
\[ b - 2(\varrho_n-\varrho_{n,(n)}) \in X^*(T_n/\C)^+, \]
we see that $\Pi_\infty$ has the same infinitesimal character as $\rho_{b - 2(\varrho_n-\varrho_{n,(n)})}^\vee$. Moreover proposition 4.5 of \cite{dbar} tells us that
\[ \Hom_{U_{n,\infty}^0A_n(\R)^0}(\rho_{(n),b}^\vee, \Pi_\infty) \neq (0). \]
We deduce that
\[ \Hom_{U_{n,\infty}^0A_n(\R)^0}(\rho_{(n),-2(\varrho_n-\varrho_{n,(n)})} , \Pi_\infty \otimes \rho_{b - 2(\varrho_n-\varrho_{n,(n)})}) \neq (0). \]
However $\rho_{(n),-2(\varrho_n-\varrho_{n,(n)})}$ is the representation of $U_{n,\infty}^0A_n(\R)^0$ on $\wedge^{[F^+:\Q]n^2} \gp^+$. Thus
\[ \Hom_{U_{n,\infty}^0A_n(\R)^0}( \wedge^{[F^+:\Q]n^2} \gp \otimes_\R \C, \Pi_\infty \otimes \rho_{b - 2(\varrho_n-\varrho_{n,(n)})}) \neq (0). \]
Proposition II.3.1 of \cite{bw} then tells us that 
\[ H^{[F^+:\Q]n^2}((\Lie G_n(\R)) \otimes_\R \C, U_{n,\infty}^0A_n(\R)^0, \Pi_\infty \otimes \rho_{b - 2(\varrho_n-\varrho_{n,(n)})}) \neq (0), \]
and the lemma follows.
\pfend

\begin{cor}\label{unitgr} Suppose that 
\[ b=(b_0,(b_{\tau,i})_{\tau \in \Hom(F,\barQQ_p)}) \in X^*(T_n/\barQQ_p)_{(n)}^+ \]
satisfies 
\[ -2n \geq b_{\tau,1}+b_{\tau c,1} \]
for all $\tau \in \Hom(F,\barQQ_p)$. If $\Pi$ is an irreducible sub-quotient of $H^0(X_n^\mini,\cE_{\rho_{(n),b}}^\sub)$, then there is a continuous representation
\[ R_{p}(\Pi): G_F \lra GL_{2n}(\barQQ_p) \]
which is de Rham above $p$ and has the following property: 
Suppose that $v\ndiv p$ is a prime of $F$ which is
\begin{itemize}
\item either split over $F^+$,
\item or inert but unramified over $F^+$ and $\Pi$ is unramified at $v$;
\end{itemize}
 then
\[ \WD(R_{p}(\Pi)|_{G_{F_v}})^{\Fsemis} \cong \rec_{F_v}(\BC(\Pi_q)_v |\det|_v^{(1-2n)/2}), \]
where $q$ is the rational prime below $v$.
\end{cor}

\pfbegin 
By the lemma $\imath \Pi$ is the finite part of a cohomological, square integrable, automorphic representation of $G_n(\A)$. The result now follows from corollary \ref{galois1}.
\pfend

\newpage

\section{The Ordinary Locus.}

\subsection{P-adic automorphic forms.}

Zariski locally on $\cX_{U}^\mini$ we may lift $\hasse_{U}$ to a (non-canonical) section $\widetilde{\hasse}_{U}$ of $\omega^{\otimes(p-1)}$ over (an open subset of) $\cX_{U}^\mini$. Although
$\widetilde{\hasse}_{U}$ is non-canonical $\widetilde{\hasse}_{U}^{p^{M-1}} \bmod p^M$ is canonical, and so these glue to give a canonical element
\[ \hasse_{M,U} \in H^0(\cX_{U}^\mini \times \Spec \Z/p^M\Z, \omega_{U}^{\otimes (p-1)p^{M-1}}). \]
Again if $g \in G_n(\A^{\infty,p}\times \Z_p)$ and $U' \supset g^{-1}Ug$ then 
\[ g  \hasse_{M,U'} =\hasse_{M,U}. \]
We will denote by $\hasse_{M,U^p(N)}$ the restriction of $\hasse_{M,U^p(N,N')}$ to
\[ H^0(\gX_{U^p(N)}^{\ord,\mini} \times \Spec \Z/p^M\Z, (\omega_{\cA^\ord/\cX^\ord,U^p(N)}^\ord)^{\otimes (p-1)p^{M-1}})\]
This is independent of $N'$.

If $\rho$ is a representation of $L_{n,(n)}$ on a finite free $\Z_{(p)}$-module then, for any integer $i$, there is a natural map
\[ \begin{array}{r} H^0(\cX^\mini_{U^p(N_1,N_2)}, \cE^\sub_{\rho \otimes (\wedge^{n[F:\Q]}\Std^\vee)^{ip^{M-1}(p-1)}} ) \cong  
H^0(\cX^\mini_{U^p(N_1,N_2)}, \cE^\sub_\rho \otimes \omega_{U}^{\otimes i(p-1)p^{M-1}}) \\
\lra  H^0(\cX^{\ord,\mini}_{U^p(N_1,N_2)}, \cE^{\ord,\sub}_\rho \otimes \Z/p^M\Z), \end{array}\]
which sends $f$ to
\[ (f|_{\cX^{\ord,\mini}_{U^p(N_1,N_2)}})/\hasse_{M,U^p(N_1,N_2)}^i. \]
These maps are $G_n(\A^\infty)^{\ord,\times}$-equivariant. 

\begin{lem}\label{katzl} For any $r$ the induced map
\[ \begin{array}{rl} & \bigoplus_{j=r}^\infty H^0(\cX^\mini_{U^p(N_1,N_2)}, \cE^\sub_{U^p(N_1,N_2),\rho \otimes (\wedge^{n[F:\Q]}\Std^\vee)^{jp^{M-1}(p-1)}} )\\ \lra &H^0(\cX^{\ord,\mini}_{U^p(N_1,N_2)}, \cE^{\ord,\sub}_{U^p(N_1,N_2),\rho} \otimes \Z/p^M\Z) \end{array}\]
is surjective. \end{lem}

\pfbegin
To simplify the formulae in this proof, for the duration of the proof  we will write $U$ for $U^p(N_1,N_2)$.

Multiplying by a power of $\hasse_{M,U}$ we may replace $\rho$ by
\[ \rho \otimes (\wedge^{n[F:\Q]}\Std^\vee)^{tp^{M-1}(p-1)} \]
and $r$ by $r-t$ for any $t$. Thus, using the ampleness of $\omega_{U}$ over $\cX^\mini_{U}$, we may suppose that
\[ H^i(\cX^\mini_{U}, \cE_{U,\rho}^\sub \otimes \omega_{U}^{\otimes j}) =(0) \]
for all $i > 0$ and $j\geq 0$. We may also suppose that $r\leq 0$. Then we may replace $r$ by $0$.

Because $\cX^{\ord,\mini}_{U} \times \Spec \Z/p^M\Z$ is a union of connected components of 
\[ \cY=\cX^{\mini}_{U}\times \Spec \Z/p^M\Z -\barX^{\mini , \nord}_{U}\]
 it suffices to replace $\cX^{\ord,\mini}_{U} \times \Spec \Z/p^M\Z$ by $\cY$. 

Now we need to show that
\[ \bigoplus_{j=0}^\infty H^0(\cX^\mini_U, \cE^\sub_{U,\rho} \otimes \omega_{U}^{\otimes j(p-1)p^{M-1}}) \onto 
H^0(\cY, \cE^{\ord,\sub}_{U,\rho} \otimes \Z/p^M\Z), \]
under the assumption that
\[ H^i(\cX^\mini_{U}, \cE_{U,\rho}^\sub \otimes \omega_{U}^{\otimes j}) =(0) \]
for all $i > 0$ and $j\geq 0$.

The scheme $\cY$ is relatively affine over $\cX^{\mini}_{U}$ corresponding to the sheaf of algebras
\[ (\bigoplus_{j=0}^\infty \omega_{U}^{\otimes jp^{M-1}(p-1)})/(\hasse_{M,U}-1, p^M). \]
Hence
\[ H^0(\cY, \cE^\sub_{U,\rho}) \cong H^0\left(\cX^\mini_{U}, \left( \bigoplus_{j=0}^\infty \cE^\sub_{U,\rho} \otimes  \omega_{U}^{\otimes p^{M-1}(p-1)j}   \right)/(\hasse_{M,U}-1, p^M) \right)  \]
and the map
\[ \bigoplus_{j=0}^\infty H^0(\cX^\mini_U, \cE^\sub_{U,\rho} \otimes \omega_{U}^{\otimes j(p-1)p^{M-1}}) \lra 
H^0(\cY, \cE^{\ord,\sub}_{U,\rho} \otimes \Z/p^M\Z) \]
is induced by the map 
\[ \bigoplus_{j=0}^\infty  \cE^\sub_{U,\rho}  \otimes  \omega_{U}^{\otimes j(p-1)p^{M-1}}\onto \left. \left( \bigoplus_{j=0}^\infty \cE^\sub_{U,\rho} \otimes  \omega_{U}^{\otimes p^{M-1}(p-1)j}   \right)\right/ (\hasse_{M,U}-1, p^M)   \]
of sheaves over $\cX_U^\mini$.

Because 
\[ H^i(\cX^\mini_{U}, \cE_{U,\rho}^\sub \otimes \omega_{U}^{\otimes j}) =(0) \]
for all $i > 0$ and $j\geq 0$; we see that
\[ H^0(\cX^\mini_{U}, \cE_{U^p,\rho}^\sub \otimes \omega_{U}^{\otimes j}) \otimes \Z/p^M\Z  \liso  H^0(\cX^\mini_{U}, \cE_{U,\rho}^\sub \otimes \omega_{U}^{\otimes j} \otimes \Z/p^M\Z)  \]
for all $j\geq 0$, and
\[ H^i(\cX^\mini_{U}, \cE_{U,\rho}^\sub \otimes \omega_{U}^{\otimes j} \otimes \Z/p^M\Z) =(0) \]
for all $i>0$ and $j \geq 0$. Thus it suffices to check that 
\[ \begin{array}{c} \left. H^0\left(\cX^\mini_{U},  \bigoplus_{j=0}^\infty \cE^\sub_{U,\rho} \otimes  \omega_{U}^{\otimes p^{M-1}(p-1)j}   \otimes \Z/p^M\Z \right) \right/(\hasse_{M,U}-1) \\ \da \\
H^0\left(\cX^\mini_{U},  ( \bigoplus_{j=0}^\infty \cE^\sub_{U,\rho} \otimes  \omega_{U}^{\otimes p^{M-1}(p-1)j}   \otimes \Z/p^M\Z)/ (\hasse_{M,U}-1) \right) \end{array} \]
is surjective. This follows using the long exact sequence in cohomology associated to the short exact sequence
\[ \begin{array}{l} (0) \lra  \bigoplus_{j=0}^\infty \cE^\sub_{U,\rho} \otimes  \omega_{U}^{\otimes p^{M-1}(p-1)j}   \otimes \Z/p^M\Z 
\stackrel{\hasse_{M,U}-1}{\lra} \\ \bigoplus_{j=0}^\infty \cE^\sub_{U,\rho} \otimes  \omega_{U}^{\otimes p^{M-1}(p-1)j}   \otimes \Z/p^M\Z 
\lra \\ \left(\bigoplus_{j=0}^\infty \cE^\sub_{U,\rho} \otimes  \omega_{U}^{\otimes p^{M-1}(p-1)j}   \otimes \Z/p^M\Z \right)/(\hasse_{M,U}-1) \lra (0)\end{array}  \]
and the vanishing
\[ H^1\left(\cX^\mini_{U},  \bigoplus_{j=0}^\infty \cE^\sub_{U,\rho} \otimes  \omega_{U}^{\otimes p^{M-1}(p-1)j}   \otimes \Z/p^M\Z \right)=(0).\]
\pfend

Let $S$ denote the set of rational primes consisting of $p$ and the primes where $F$ ramifies. Also choose a neat open compact subgroup 
\[ U^p=G_n(\hatZ^S) \times U^p_S  \subset G_n(\A^{p,\infty}). \]

Suppose that $v$ is a place of $F$ above a rational prime $q \not\in S$ and let $i \in \Z$. There is a unique element $\gt_v^{(i)}$ in the Bernstein centre of $G_n(\Q_q)$ such that 
\begin{itemize}
\item $\gt_v^{(i)}$ acts as $0$ on any irreducible smooth representation of $G_n(\Q_q)$ over $\C$ which is not a subquotient of an unramified principal series;
\item on an unramified representation $\Pi_q$ of $G_n(\Q_q)$ the eigenvalue of $\gt_v^{(i)}$ on $\Pi_q$ equals $\tr \rec_{F_v}(\BC(\Pi_q)_v|\det|_v^{(1-2n)/2})(\Frob_v^i)$. 
\end{itemize}
Multiplying $\gt_v^{(i)}$ by the characteristic function of $G_n(\Z_q)$ we obtain a unique element $T_v^{(i)} \in \C[G_n(\Z_q)\backslash G_n(\Q_q) /G_n(\Z_q)]$ such that if $\Pi_q$ is an unramified representation of $G_n(\Q_q)$ and if $T_v^{(i)}$ has eigenvalue $t_v^{(i)}(\Pi_q)$ on $\Pi_q^{G_n(\Z_q)}$ then 
\[ \tr \rec_{F_v}(\BC(\Pi_q)_v|\det|_v^{(1-2n)/2})(\Frob_v^i) = t_v^{(i)}(\Pi_q). \]
If $\sigma \in \Aut(\C)$ we see that ${}^\sigma T_v^{(i)}=T_v^{(i)}$. (Use the fact that
\[ {}^\sigma \rec_{F_v}(\BC(\Pi_q)_v|\det|_v^{(1-2n)/2}) \cong \rec_{F_v}(\BC({}^\sigma\Pi_q)_v|\det|_v^{(1-2n)/2}).) \]
Thus
\[ T_v^{(i)} \in \Q[G_n(\Z_q)\backslash G_n(\Q_q) /G_n(\Z_q)]. \]
Choose $d_v^{(i)} \in \Q^\times$ such that
\[ d_v^{(i)}T_v^{(i)} \in \Z[G_n(\Z_q)\backslash G_n(\Q_q) /G_n(\Z_q)]. \]

Suppose that $q\not\in S$ is a rational prime. Let $u_1,...,u_r$ denote the primes of $F^+$ above $\Q$ which split $u_i=w_i{}^cw_i$ in $F$, and let $v_1,...,v_s$ denote the primes of $F^+$ above $q$ which do not split in $F$. 
Then under the identification
\[ G_n(\Q_q) \cong \prod_{i=1}^r GL_{2n}(F_{w_i}) \times H \]
of section \ref{secbc}, the Hecke operator $T_{w_i}^{(1)}$ is identified with the double coset 
\[ G_n(\Z_q) a_i G_n(\Z_q), \]
 where $a_i \in GL_n(F_{w_i})$ is the diagonal matrix $\diag(1,...,1,\varpi_{w_i})$, and we may take $d_{w_i}^{(1)}=1$.

 We will call a topological $\Z_p[G_n(\Zhat^S) \backslash G_n(\A^S) /G_n(\Zhat^S)]$-algebra $\T$ {\em of Galois type} if for every there is a continuous pseudo-representation (see \cite{pseudo}) 
\[ T: G_F^S \lra \T \]
such that
\[ d_v^{(i)}T(\Frob_v^i) = d_v^{(i)}T_v^{(i)} \]
for all $v|q \not\in S$ and all $i \in \Z$. 

Let $\T_{U^p(N_1,N_2),\rho}^{S}$ denote the image of $\Z_p[G_n(\Zhat^S)\backslash G_n(\A^S) /G_n(\Zhat^S)]$ in the endomorphism algebra $\End(H^0(\cX^{\mini}_{U^p(N_1,N_2)}, \cE_\rho^{\sub}))$, which is also the image in the endomorphism algebra $\End(H^0(X^{\mini}_{U^p(N_1,N_2)}, \cE_\rho^{\sub}))$.

\begin{lem}\label{l62}  For $t$ sufficiently large $\T_{U^p(N_1,N_2),\rho\otimes (\wedge^{n[F:\Q]} \std)^{\otimes t}}^{S}$ is of Galois type.
\end{lem}

\pfbegin 
Write 
\[ \rho_t =  \rho\otimes (\wedge^{n[F:\Q]} \std)^{\otimes t}. \]
It  suffices to show that there is a continuous pseudo-representation
\[ T: G_F^S \lra \T_{U^p(N_1,N_2),\rho_t}^{S} \otimes \barQQ_p \]
which is unramified outside $S$ and satisfies 
\[ T(\Frob_v^i) = T_v^{(i)} \]
for all $v |q \not\in S$ and all $i \in \Z$. (Because $T$ will then automatically be valued in $\T_{U^p(N_1,N_2),\rho_t}^{S}$, by the Cebotarev density theorem. Note that if $v$ is a prime of $F$ split over $F^+$ and lying above a rational prime $q \not\in S$, then 
\[ T(\Frob_v) = T_v^{(1)} \in \T_{U^p(N_1,N_2),\rho_t}^{S}.) \]
We may then reduce to the case that $\rho \otimes \barQQ_p$ is irreducible. Let 
\[ (b_0,(b_{\tau,i})) \in X^*(T_n/\barQQ_p)_{(n)}^+ \]
denote the highest weight of $\rho \otimes \barQQ_p$. 

Suppose that $t$ satisfies the inequality
\[ -2n \geq (b_{\tau,1}-t(p-1))+(b_{\tau c,1}-t(p-1)). \]
By lemma \ref{ssauto},
\[ \T_{U^p(N_1,N_2),\rho_t}^{S} \otimes \barQQ_p \cong \bigoplus_{\Pi} \barQQ_p \]
where the sum runs over irreducible admissible representations of  $G_n(\A^{\infty})$ with $\Pi^{U^p(N_1,N_2)} \neq (0)$ which occur in $H^0(X^\mini \times \Spec \barQQ_p,\cE_{\rho_t  }^\sub)$. Further, from corollary \ref{unitgr}, we deduce that 
there is a continuous representation 
\[ r: G_F^S \lra GL_{2n}(\T_{U^p(N_1,N_2),\rho_t}^{S}) \]
such that if $v|q \not\in S$ then $r$ is unramified at $v$ and
\[ \tr r(\Frob_v^i)=T_v^{(i)} \]
for all $i \in \Z$. Taking $T=\tr r$ completes the proof of the lemma.
\pfend

If
\[ W \subset H^0(\gX^{\ord,\mini}_{U^p(N)}, \cE^{\ord,\sub}_{\rho} ) \]
(resp.
\[ W \subset H^0(\cX^{\ord,\mini}_{U^p(N_1,N_2)}, \cE^{\ord,\sub}_{\rho} \otimes \Z/p^M\Z )) \]
is a finitely generated $\Z_p$-submodule invariant under the action of the algebra $\Z_p[G_n(\Zhat^S)\backslash G_n(\A^S) /G_n(\Zhat^S)]$, then let $\T_{U^p(N),\rho}^{\ord,S}(W)$ (resp. $\T_{U^p(N_1,N_2),\rho}^{\ord,S}(W)$) denote the image of $\Z_p[G_n(\Zhat^S)\backslash G_n(\A^S) /G_n(\Zhat^S)]$ in $\End_{\Z_p}(W)$. The next corollary follows from lemmas \ref{katzl} and \ref{l62}.

\begin{cor}  If 
\[ W \subset H^0(\cX^{\ord,\mini}_{U^p(N_1,N_2)}, \cE^{\ord,\sub}_{\rho} \otimes \Z/p^M\Z )\]
is a finitely generated $\Z_p$-submodule invariant under the action of the algebra $\Z_p[G_n(\Zhat^S)\backslash G_n(\A^S) /G_n(\Zhat^S)]$, then $\T_{U^p(N_1,N_2),\rho}^{\ord,S}(W)$ is of Galois type.\end{cor}

We deduce from this the next corollary.

\begin{cor}
If 
\[ W \subset H^0(\gX^{\ord,\mini}_{U^p(N)}, \cE^{\ord,\sub}_{\rho} ) \]
is a finitely generated $\Z_p$-submodule invariant under the action of the algebra $\Z_p[G_n(\Zhat^S)\backslash G_n(\A^S) /G_n(\Zhat^S)]$, then $\T_{U^p(N),\rho}^{\ord,S}(W)$ is of Galois type.
\end{cor}

Finally we deduce the following proposition.
\begin{prop}\label{padicgalois} Suppose that $\rho$ is a representation of $L_{n,(n)}$ over $\Z_{(p)}$. Suppose also that $\Pi$ is an irreducible quotient of an {\em admissible} $G_n(\A^{\infty})^{\ord,\times}$-sub-module $\Pi'$ of $H^0(\gX^{\ord,\mini}, \cE_\rho^{\ord,\sub})_{\barQQ_p}$. Then there is a continuous semi-simple representation
\[ R_p({\Pi}):G_F \lra GL_{2n}(\barQQ_p) \]
with the following property: Suppose that $q\neq p$ is a rational prime above which $F$ and $\Pi$ are unramified, and suppose that $v|q$ is a prime of $F$. Then \[ \WD(R_{p}(\Pi)|_{G_{F_v}})^{\Fsemis} \cong \rec_{F_v}(\BC(\Pi_q)_v |\det|_v^{(1-2n)/2}), \]
where $q$ is the rational prime below $v$.
\end{prop}

\pfbegin
Let $S$ denote the set of rational primes consisting of $p$ and the primes where $F$ or $\Pi$ ramifies. Also choose a neat open compact subgroup 
\[ U^p=G_n(\hatZ^S) \times U^p_S \]
and integer $N$ such that 
\[ \Pi^{U^p(N)} \neq (0). \]
As $(\Pi')^{U^p(N)}$ is a finite dimensional, and hence closed, subspace of  the topological vector space
$H^0(\gX^{\ord,\mini}, \cE_\rho^{\ord,\sub})_{\barQQ_p}$ preserved by 
$\Z_p[G_n(\Zhat^S)\backslash G_n(\A^S) /G_n(\Zhat^S)]$ and, as there is a
$\Z_p[G_n(\Zhat^S)\backslash G_n(\A^S) /G_n(\Zhat^S)]$-equivariant map $(\Pi')^{U^p(N)} \onto \Pi^{U^p(N)}$, there is a continuous homomorphism
\[ \theta: \T_{U^p(N),\rho}^{\ord,S}((\Pi')^{U^p(N)}) \lra \barQQ_p \]
which for $v|q \not\in S$ sends $T_v^{(i)}$ to its eigenvalue on $\Pi^{G_n(\Z_q)}$. Proposition \ref{padicgalois} now follows from the above corollary and the main theorem on pseudo-representations (see \cite{pseudo}).
\pfend

We remark that we don't know how to prove this proposition for a general irreducible subquotient of $H^0(\gX^{\ord,\mini}, \cE_\rho^{\ord,\sub})_{\barQQ_p}$ (or indeed whether the corresponding statement remains true).

\newpage

\subsection{Interlude concerning linear algebra.}

Suppose that $K$ is an algebraic extension of $\Q_p$. For $a \in \Q$, we say that a polynomial $P(X) \in K(X)$ has {\em slopes $\leq a$} if $P(X) \neq 0$ and every root of $P(X)$ in $\barK$ has $p$-adic valuation $\leq a$. (We normalize the $p$-adic valuation so that $p$ has valuation $1$.) If $V$ is a $K$-vector space and $T$ is an endomorphism of $V$, then we say that $V$ {\em admits slope decompositions} for $T$, if for each $a \in \Q$ there is a decomposition 
\[ V=V_{\leq a} \oplus V_{>a} \]
with the following properties:
\begin{itemize}
\item $T$ preserves $V_{\leq a}$ and $V_{>a}$;
\item $V_{\leq a}$ is finite dimensional;
\item if $P(X) \in K[X]$ has slopes $\leq a$ then the endomorphism $P(T)$ restricts to an automorphism of $V_{>a}$;
\item there is a non-zero polynomial $P(X) \in K[X]$ with slopes $\leq a$ such that the endomorphism $P(T)$ restricts to $0$ on $V_{\leq a}$.
\end{itemize}
In this case $V_{\leq a}$ and $V_{>a}$ are unique, and we refer to them as the {\em slope $a$ decomposition} of $V$ with respect to $T$.

\begin{lem}\label{slope} \begin{enumerate}
\item If $V$ is finite dimensional then it always admits slope decompositions. 

\item If $K$ is a finite extension of $\Q_p$, if $V$ is a $K$-Banach space, and if $T$ is a completely continuous (see \cite{serreecc}) endomorphism of $V$ then $V$ admits slope decompositions for $T$. 

\item Suppose that $L/K$ is an algebraic extension and that $V$ is a $K$ vector space which admits slope decompositions with respect to an endomorphism $T$. Then $V \otimes_KL$ also admits slope decompositions with respect to $T$.

\item Suppose that $V_1$ admits slope decompositions with respect to $T_1$; that $V_2$ admits a slope decomposition with respect to $T_2$ and that $d:V_1 \ra V_2$ is a linear map such that
\[ d \circ T_1 = T_2 \circ d. \]
Then for all $a \in \Q$ we have 
\[ d V_{1,\leq a} \subset V_{2,\leq a} \]
and
\[ d V_{1,> a} \subset V_{2,>a}. \]
Moreover $\ker d$ admits slope decompositions for $T_1$, while $\Im d$ and $\coker d$ admit slope decompositions for $T_2$. More specifically 
\[ (\ker d)_{\leq a} = (\ker d) \cap V_{1,\leq a} \]
and
\[ (\ker d)_{> a} = (\ker d) \cap V_{1,> a} \]
and
\[ (\Im d)_{\leq a} = V_{1,\leq a}/(\ker d)_{\leq a} \]
and
\[ (\Im d)_{> a} = V_{1,> a}/(\ker d)_{> a} \]
and
\[ (\coker d)_{\leq a} = V_{2,\leq a}/(\Im d)_{\leq a} \]
and
\[ (\coker d)_{> a} = V_{2,> a}/(\Im d)_{> a}. \]

\item Suppose that
\[ V_1 \subset V_2 \subset V_3 \subset .... \subset V_\infty \]
are vector spaces with 
\[ V_\infty = \bigcup_{i=1}^\infty V_i. \]
Suppose also that $T$ is an endomorphism of $V_\infty$ such that for all $i>1$ 
\[ T V_i \subset V_{i-1}. \]
If for each $i$ the space $V_i$ admits slope decompositions for $i$, then $V_\infty$ admits slope decompositions for $T$. 

\item Suppose that 
\[ (0) \lra V_1 \lra V \lra V_2 \lra (0) \]
is an exact sequence of $K$ vector spaces and that $T$ is an endomorphism of $V$ that preserves $V_1$. If $V_1$ and $V_2$ both admit slope decompositions with respect to $T$, then so does $V$. Moreover we have short exact sequences
\[ (0) \lra V_{1,\leq a} \lra V_{\leq a} \lra V_{2,\leq a} \lra (0) \]
and
\[ (0) \lra V_{1,> a} \lra V_{> a} \lra V_{2,> a} \lra (0) \]
\end{enumerate} \end{lem}

\pfbegin
The first and third and fourth parts are straightforward. The second part follows from \cite{serreecc}.

For the fifth part
one checks that $V_{i, \leq a}$ is independent of $i$. If we set
\[ V_{\infty,\leq a}=V_{i,\leq a} \]
for any $i$, and
\[ V_{\infty,>a} = \bigcup_{i=1}^\infty V_{i,>a}, \]
then these provide the slope $a$ decomposition of $V_\infty$ with respect to $T$.

Finally we turn to the sixth part. Choose non-zero polynomials $P_i(X) \in K[X]$ with slopes $\leq a$ such that $P_i(T)V_{i,\leq a}=(0)$, for $i=1,2$. Set $P(X)=P_1(X)P_2(X)$. Also set $V_{\leq a} = \ker P(T)$ and $V_{>a}= \Im P(T)$. 
We have complexes
\[ (0) \lra V_{1,>a} \lra V_{>a} \lra V_{2,>a} \lra (0) \]
and
\[ (0) \lra  V_{1,\leq a} \lra V_{\leq a} \lra V_{2,\leq a} \lra (0). \]
It suffices to show that these complexes are both short exact sequences. For then we see that, if $Q(X) \in K[X]$ has slopes $\leq a$, then the restriction of $Q(T)$ to $V_{>a}$ is an automorphism of $V_{>a}$. Applying this to $P(T)$, we see that $V_{\leq a} \cap V_{>a}=(0)$. Moreover $V_{\leq a}+V_{>a}$ contains $V_1$ and maps onto $V_2$, so that $V=V_{\leq a}+V_{>a}$.

To show the first complex is short exact we need only check that $V_{1,>a}=V_{>a} \cap V_1$, i.e. that $V_{1,\leq a} \cap V_{>a} =(0)$. So suppose that $v \in V_{1,\leq a} \cap V_{>a}$ then $v=P(T)v'$ and $P_1(T)v=0$. Thus $P_1(T)^2P_2(T)v'=0$ so the image of $v'$ in $V_2$ lies in $V_{2,\leq a}$ and so $P_2(T)v' \in V_1$, and in fact $P_2(T)v' \in V_{1,\leq a}$. Finally we see that $v=P_1(T)P_2(T)v'=0$, as desired.

To show the second complex is short exact we have only to show that $V_{\leq a} \ra V_{2,\leq a}$ is surjective. So suppose that $\barv \in V_{2,\leq a}$ and suppose that $v \in V$ lifts $\barv$. Then $P(T) v\in V_{1,>a}$. Set
\[ v' = v -  (P(T)|_{V_{1,>a}}^{-1})P(T) v \in v+V_{1,>a} \]
Then $v'$ maps to $\barv \in V_2$, while 
\[ P(T)v' = P(T) v - P(T)v =0, \]
so that $v' \in V_{\leq a}$. 
\pfend

We warn the reader that to the best of our knowledge it is not in general true that if $V_1 \subset V$ is $T$-invariant then either $V_1$ or $V/V_1$ admits slope decompositions for $T$.

\newpage

\subsection{The ordinary locus of a toroidal compactification as a dagger space.}

We first review some general facts about dagger spaces. 
We refer to \cite{gkcrelle} for the basic facts. 

Suppose that $K/\Q_p$ is a finite extension with ring of integers $\cO_K$ and residue field $k$. Suppose also that $\cY/\cO_K$ is quasi-projective. Let $Y$ denote the generic fibre $\cY \times \Spec K$, let $\barY$ denote the special fibre $\cY \times \Spec k$ and let $\cY^\wedge$ denote the formal completion of $\cY$ along $\barY$. Let $Y^\an$ (resp. $Y^\dag$) denote the rigid analytic (resp. dagger) space associated to $Y$. (For the latter see section 3.3 of \cite{gkcrelle}.) Thus $Y^\an$ and $Y^\dag$ share the same underlying G-topological space, and in fact the completion $(Y^\dag)'$ (see theorem 2.19 of \cite{gkcrelle}) of $Y^\dag$ equals $Y^\an$. Let $\cY^\wedge_\eta$ denote the rigid analytic space associated to $\cY^\wedge$, its `generic fibre'. Then $\cY^\wedge_\eta$ is identified with an admissible open subset $]\barY[ \subset Y^\an$. We will denote by $\cY^\dag$ the admissible open dagger subspace of $Y^\dag$ with the same underlying topological space as $]\barY[$. 

\begin{lem}\label{berth} If $\cY$ and $\cY'$ are two quasi-projective $\cO_K$-schemes as described in the previous paragraph and if $f:\cY \ra \cY'$ is a morphism, then there is an induced map $f^\dag: \cY^\dag \ra (\cY')^\dag$. 

If further $f: \barY \iso \barY'$ and $f$ is etale in a neighbourhood of $\barY$ then $f^\dag$ is an isomorphism. \end{lem}

\pfbegin The first part of the lemma is clear. 

For the second part, let $\cY \into \PP_{\cO_K}^M$ and $\cY' \into \PP_{\cO_K}^{M'}$ be closed embeddings. Let $\cP'$ denote the closure of $\cY'$ in $\PP_{\cO_K}^{M'}$. Also let $\cP$ denote the closure of $\cY$ in $\PP_{\cO_K}^{M} \times \PP_{\cO_K}^{M'}$. Then $f$ extends to a map $\cP \ra \cP'$. The second part of the lemma follows from theorem 1.3.5 of \cite{berthelot} applied to $\barY \subset \cP$ and $\barY' \subset \cP'$. \pfend

We will let $H^i_\rig(\barY)$ denote the rigid cohomology of $\barY$ in the sense of Berthelot - see for instance \cite{lestum}.

\begin{lem}\label{dagrig} \begin{enumerate}
\item If  $\cY/\cO_K$ is a smooth and quasi-projective scheme, then there is a canonical isomorphism
\[ H^i_\rig(\barY) \cong \HH^i(\cY^\dag, \Omega^\bullet_{\cY^\dag}). \]
\item If $f:\cY \ra \cZ$ is a morphism of smooth quasi-projective schemes over $\cO_K$ then the following diagram is commutative:
\[ \begin{array}{ccc}  H^i_\rig(\barZ) &\stackrel{f^*}{\lra} & H^i_\rig(\barY) \\ ||\wr && || \wr \\ \HH^i(\cZ^\dag, \Omega^\bullet_{\cZ^\dag})& \stackrel{f^*}{\lra}& \HH^i(\cY^\dag, \Omega^\bullet_{\cY^\dag}).
\end{array} \]
\end{enumerate} \end{lem}

\pfbegin For the first part apply theorem 5.1 of \cite{gkcrelle} to the closure of $\cY$ in some projective space over $\cO_K$. For the second part choose embeddings $i: \cY \into \PP_{\cO_K}^{M}$ and $i': \cZ \into \PP_{\cO_K}^{M'}$. Let $\cP'$ denote the closure of $\cZ$ in $\PP_{\cO_K}^{M'}$ and $\cP$ the closure of $\cY$ in $\PP_{\cO_K}^M \times \cP'$, so that $f$ extends to a map $\cP \ra \cP'$. The desired result again follows from theorem 5.1 of \cite{gkcrelle}, because the isomorphisms of theorem 5.1 of \cite{gkcrelle} are functorial under morphisms of the set up in that theorem. 
\pfend

[It is unclear to us whether this functoriality is supposed to be implied by the word `canonical' in the statement of theorem 5.1 of \cite{gkcrelle}. For safety's sake we sketch the argument for this functoriality.  More precisely if $f: \cX_1\ra \cX_2$ is a morphism of proper admissible formal $\Spf R$-schemes which takes $Y_1 \subset \cX_{1,s}$ to $Y_2 \subset \cX_{2,s}$, then we will show that the isomorphisms of theorem 5.1 of \cite{gkcrelle} are compatible with the maps in cohomology induced by $f$. For part (a) we also suppose that we are given a map $f^*: f^* \cF_2 \ra \cF_1$.

Using the notation of part (a) of theorem 5.1 of \cite{gkcrelle}, it suffices to show that the diagram
\[ \begin{array}{ccc} H^q(X_2,\cF_{2,X_2}) & \stackrel{f^*}{\lra} & H^q(X_1,\cF_{1,X_1}) \\ \da && \da \\ H^q(]\barY_2[_{\cX_2},j_2^\dag \cF_2') &   \stackrel{f^*}{\lra} &  H^q(]\barY_1[_{\cX_2},j_2^\dag \cF_2')     \end{array} \]
commutes. (The functoriality of parts (b) and (c) follow easily from the functoriality of part (a).) The vertical morphisms arise from maps $L_k^\bullet \ra K_k^\bullet$ of resolutions of the sheaves $Ri_* \cF_{k,X_k}$  and $j_k^\dag \cF_k'$ respectively. To define these resolutions one needs to choose
 affine covers $\{ Y_{k,i} \}$ of $Y_k$. We may suppose these are chosen so that $f$ carries $Y_{1,i}$ to $Y_{2,i}$ for all $i$. Then $L_k^\bullet$ and $K_k^\bullet$ are the Cech complexes with
 \[ L^q_k = \bigoplus_{\# J=q}   i_{J*} \cF_{k,]Y_{k,J}[_{\cX_k}} \]
 and
 \[ K^q_k = \bigoplus_{\# J=q}  j_{k,J}^\dag \cF_k'. \]
 The maps $L_k^\bullet \ra K_k^\bullet$ arise from maps
 \[ (i_{J*} \cF_{k,]Y_{k,J}[_{\cX_k}})(U) \cong \lim_{\ra V} \cF_k'(V) \lra \lim_{\ra V'} \cF_k'(V' \cap U) = (j_{k,J}^\dag \cF_k')(U). \]
 Here $V$ runs over strict neighbourhoods of $U \cap ]Y_{k,J}[_{\cX_k}$ in $]\barY_k[_{\cX_k}$ and $V'$ runs over strict neighbourhoods of $]Y_{k,J}[_{\cX_k}$ in $]\barY_k[_{\cX_k}$. The first isomorphism is justified in section 2.23 of \cite{gkcrelle}. The second morphism arises because, for every $V$, we can find a $V'$ so that 
 \[ V' \cap U \subset V. \]
 It suffices to show that if $fU_1 \subset U_2$, then the diagrams
 \[ \begin{array}{ccc} (i_{J*} \cF_{2,]Y_{2,J}[_{\cX_2}})(U)_2 & \stackrel{f^*}{\lra} & (i_{J*} \cF_{1,]Y_{1,J}[_{\cX_1}})(U_1) \\ \da && \da \\
 (j_{2,J}^\dag \cF_2')(U_2) & \stackrel{f^*}{\lra} & (j_{1,J}^\dag \cF_1')(U_1) \end{array} \]
 are commutative. But this is now clear.]

\begin{lem}\label{daggaga} Suppose that $f:X \ra Y$ is a proper morphism between $\Q_p$-schemes of finite type and that $\cF/X$ is a coherent sheaf. Denote by $f^\dag:X^\dag \ra Y^\dag$ the corresponding map of dagger spaces and by $\cF^\dag$ the coherent sheaf on $X^\dag$ corresponding to $\cF/X$. Suppose also that $V$ is an admissible open subset of  $Y^\dag$ and that $U$ is its pre-image in $X^\dag$. Then
\[ R^i(f^\dag|_U)_* \cF^\dag|_U \cong (R^if_* \cF)^\dag|_V, \]
where $(R^if_* \cF)^\dag$ denotes the coherent sheaf on $Y^\dag$ corresponding to $(R^if_* \cF)/Y$.
\end{lem}

\pfbegin It suffices to check this in the case $V=Y^\dag$. There is a chain of isomorphisms
\[ \left[(R^i f_\ast \cF)^\dag \right]^\text{an} \to (R^i f_\ast \cF)^\text{an} \to R^i f^\text{an}_\ast \cF^\text{an} \to (R^i f^\dag_\ast \cF^\dag )^\text{an}. \]
The first arrow is the transitivity of dagger and rigid analytification. The second arrow is Theorem 6.5 of \cite{kopf}. The third arrow is Theorem 3.5 of \cite{gkcrelle}. Since $Y^\dag$ is partially proper, Theorem 2.26 of \cite{gkcrelle} implies that there is a unique isomorphism $(R^i f_\ast \cF)^\dag \cong R^i f^\dag_\ast \cF^\dag$ which recovers the above map after passage to rigid spaces. 
\pfend

Now we return to our Shimura and Kuga-Sato varieties.

If $U^p$ is a neat open compact subgroup of $G_n^{(m)}(\A^{\infty,p})$, if   
$N_2 \geq N_1 \geq 0$ and if $(U^p(N_1,N_2),\Sigma) \in \cJ^{(m),\tor}$, we will write 
\[ \cA^{(m),\ord,\dag}_{U^p(N_1,N_2),\Sigma} \]
(resp.
\[ \partial \cA^{(m),\ord,\dag}_{U^p(N_1,N_2),\Sigma}, \]
resp.
\[ \partial_{[\sigma]} \cA^{(m),\ord,\dag}_{U^p(N_1,N_2),\Sigma} \]
for ${[\sigma]} \in \cS(U^p(N_1,N_2),\Sigma)$)
for the dagger space associated to $\cA^{(m),\ord}_{U^p(N_1,N_2),\Sigma}$ (resp. $\partial \cA^{(m),\ord}_{U^p(N_1,N_2),\Sigma}$, resp. $\partial_{[\sigma]} \cA^{(m),\ord}_{U^p(N_1,N_2),\Sigma}$) as described in the paragraph before lemma \ref{berth}. For $s>0$ also write
\[ \partial^{(s)} \cA^{(m),\ord,\dag}_{U^p(N_1,N_2),\Sigma} = \coprod_{\substack{{[\sigma]} \in \cS(U^p(N_1,N_2),\Sigma) \\ \dim {[\sigma]} = s-1}} \partial_{[\sigma]} \cA^{(m),\ord,\dag}_{U^p(N_1,N_2),\Sigma} \]
and $i^{(s)}$ for the finite map
\[ \partial^{(s)} \cA^{(m),\ord,\dag}_{U^p(N_1,N_2),\Sigma} \lra \partial \cA^{(m),\ord,\dag}_{U^p(N_1,N_2),\Sigma} \into \cA^{(m),\ord,\dag}_{U^p(N_1,N_2),\Sigma}. \]
We set 
\[ \partial^{(0)} \cA^{(m),\ord,\dag}_{U^p(N_1,N_2),\Sigma} =\cA^{(m),\ord,\dag}_{U^p(N_1,N_2),\Sigma} .\]
Then the  various systems of dagger spaces $\{ \cA^{(m),\ord,\dag}_{U^p(N_1,N_2),\Sigma} \}$ and $\{ \partial \cA^{(m),\ord,\dag}_{U^p(N_1,N_2),\Sigma} \}$ and $\{ \partial^{(s)} \cA^{(m),\ord,\dag}_{U^p(N_1,N_2),\Sigma} \}$ have compatible actions of $G_n^{(m)}(\A^\infty)^\ord$. If $(U^p)'$ denotes the image of $U^p$ in $G_n(\A^{p,\infty})$ then there is a natural map 
\[ \cA^{(m),\ord,\dag}_{U^p(N_1,N_2),\Sigma} \lra \cX^{\mini,\ord,\dag}_{(U^p)'(N_1,N_2)}. \]
These maps are $G_n^{(m)}(\A^\infty)^\ord$-equivariant.

If $N_2' \geq N_2$ and if $\Sigma'$ is a refinement of $\Sigma$ with $\Sigma^\ord=(\Sigma')^\ord$ then the natural map
\[ \cA^{(m),\ord}_{U^p(N_1,N_2'),\Sigma'} \lra \cA^{(m),\ord}_{U^p(N_1,N_2),\Sigma} \]
restricts to an isomorphism
\[ \barA^{(m),\ord}_{U^p(N_1,N_2'),\Sigma'} \liso \barA^{(m),\ord}_{U^p(N_1,N_2),\Sigma} \]
and is etale in a neighbourhood of $\barA^{(m),\ord}_{U^p(N_1,N_2'),\Sigma'}$. 
It follows from lemma \ref{berth} that
\[ \cA^{(m),\ord,\dag}_{U^p(N_1,N_2'),\Sigma'} \lra \cA^{(m),\ord,\dag}_{U^p(N_1,N_2),\Sigma} \]
is an isomorphism. We will denote this dagger space simply
\[ \cA^{(m),\ord,\dag}_{U^p(N_1),\Sigma^\ord}. \]
Similarly $\partial \cA^{(m),\ord,\dag}_{U^p(N_1,N_2),\Sigma}$ and $\partial_{[\sigma]} \cA^{(m),\ord,\dag}_{U^p(N_1,N_2),\Sigma}$ and $\partial^{(s)} \cA^{(m),\ord,\dag}_{U^p(N_1,N_2),\Sigma}$ depend only on $U^p(N_1)$ and $\Sigma^\ord$ and we will denote them $\partial \cA^{(m),\ord,\dag}_{U^p(N_1),\Sigma^\ord}$ and $\partial_{[\sigma]} \cA^{(m),\ord,\dag}_{U^p(N_1),\Sigma^\ord}$ and $\partial^{(s)} \cA^{(m),\ord,\dag}_{U^p(N_1),\Sigma^\ord}$ respectively. If ${[\sigma]} \not\in \cS(U^p(N_1),\Sigma^\ord)^\ord$ then 
\[ \partial_{[\sigma]} \cA^{(m),\ord,\dag}_{U^p(N_1),\Sigma^\ord}=\emptyset.\]
 Thus for $s>0$
\[ \partial^{(s)} \cA^{(m),\ord,\dag}_{U^p(N),\Sigma^\ord} = \coprod_{\substack{{[\sigma]} \in \cS(U^p(N),\Sigma^\ord)^\ord \\ \dim {[\sigma]} = s-1}} \partial_{[\sigma]} \cA^{(m),\ord,\dag}_{U^p(N),\Sigma^\ord} \]
The three projective systems of dagger spaces $\{ \cA^{(m),\ord,\dag}_{U^p(N),\Sigma^\ord} \}$ and $\{ \partial \cA^{(m),\ord,\dag}_{U^p(N),\Sigma^\ord} \}$ and $\{ \partial^{(s)} \cA^{(m),\ord,\dag}_{U^p(N),\Sigma^\ord} \}$
have actions of $G^{(m)}_n(\A^\infty)^\ord$.

If $(U^p)'$ contains the projection of $U^p$ and if $\Delta^\ord$ and $\Sigma^\ord$ are compatible, then there are  $G_n^{(m)}(\A^{\infty})^\ord$-equivariant maps
\[ \cA^{(m),\ord,\dag}_{U_p(N),\Sigma} \lra \cX_{(U^p)'(N),\Delta}^{\ord,\dag}. \]

The maps
\begin{itemize}
\item $\varsigma_p:\cA^{(m),\ord,\dag}_{U^p(N),\Sigma^\ord} \ra \cA^{(m),\ord,\dag}_{U^p(N),\Sigma^\ord}$,
\item and $\varsigma_p: \partial^{(s)} \cA_{U^p(N),\Sigma^\ord}^{(m),\ord,\dag} \ra \partial^{(s)} \cA_{U^p(N),\Sigma^\ord}^{(m),\ord,\dag}$
\end{itemize}
are finite, flat of degrees $p^{(2m+n)n[F^+:\Q]}$ and $p^{(2m+n)n[F^+:\Q]-s}$, respectively. (Use the finite flatness of 
 \[ \varsigma_p:\gA^{(m),\ord}_{U^p(N),\Sigma^\ord} \ra \gA^{(m),\ord}_{U^p(N),\Sigma^\ord}\] 
  and 
  \[ \varsigma_p: \partial^{(s)} \gA_{U^p(N),\Sigma^\ord}^{(m),\ord} \ra \partial^{(s)} \gA_{U^p(N),\Sigma^\ord}^{(m),\ord},\]
   together with theorems 1.7(1) and 1.12 of \cite{gkcrelle}.) 

We will write $\Omega^j_{\cA^{(m),\ord,\dag}_{U^p(N),\Sigma}}(\log \infty)$ (resp. $\Omega^j_{\cA^{(m),\ord,\dag}_{U^p(N),\Sigma}}(\log \infty) \otimes \cI_{\partial \cA^{(m),\ord,\dag}_{U^p(N),\Sigma}}$) for the locally free sheaf on $\cA^{(m),\ord,\dag}_{U^p(N),\Sigma}$ induced by $\Omega^j_{\cA^{(m),\ord}_{U^p(N,N'),\Sigma'}}(\log \infty)$ (resp. induced by $\Omega^j_{\cA^{(m),\ord}_{U^p(N,N'),\Sigma'}}(\log \infty)\otimes \cI_{\partial \cA^{(m),\ord}_{U^p(N,N'),\Sigma'}}$) for any $N' \geq N$ and $\Sigma' \in \cJ_n^{(m),\tor}$ with $(\Sigma')^\ord = \Sigma$. This is canonically independent of the choices of $N'$ and $\Sigma'$. The systems of sheaves $\{ \Omega^j_{\cA^{(m),\ord,\dag}_{U^p(N),\Sigma}}(\log \infty)\}$ and $\{ \Omega^j_{\cA^{(m),\ord,\dag}_{U^p(N),\Sigma}}(\log \infty)\otimes \cI_{\partial \cA^{(m),\ord,\dag}_{U^p(N),\Sigma}} \}$ over $\{ \cA^{(m),\ord,\dag}_{U^p(N),\Sigma}\}$ have actions of $G_n^{(m)}(\A^\infty)^\ord$. For $g \in G_n^{(m)}(\A^\infty)^\ord$ the map
\[ g: g^* \Omega^j_{\cA^{(m),\ord,\dag}_{(U^p)'(N'),\Sigma'}}(\log \infty)  \lra \Omega^j_{\cA^{(m),\ord,\dag}_{U^p(N),\Sigma}}(\log \infty) \]
is an isomorphism. The inverse of $\varsigma_p^*$ gives maps
\[ \varsigma_{p,*} \Omega^j_{\cA^{(m),\ord,\dag}_{U^p(N),\Sigma}}(\log \infty) \liso \Omega^j_{\cA^{(m),\ord,\dag}_{U^p(N),\Sigma}}(\log \infty) \otimes_{\cO_{\cA^{(m),\ord,\dag}_{U^p(N),\Sigma}} ,\varsigma_p^{*}} \cO_{\cA^{(m),\ord,\dag}_{U^p(N),\Sigma}} \]
and
\[ \begin{array}{l} \varsigma_{p,*}( \Omega^j_{\cA^{(m),\ord,\dag}_{U^p(N),\Sigma}}(\log \infty) \otimes \cI_{\partial \cA^{(m),\ord,\dag}_{U^p(N),\Sigma}} ) \liso 
\\ \Omega^j_{\cA^{(m),\ord,\dag}_{U^p(N),\Sigma}}(\log \infty) \otimes_{\cO_{\cA^{(m),\ord,\dag}_{U^p(N),\Sigma}} ,\varsigma_p^{*}} \cI_{\partial \cA^{(m),\ord,\dag}_{U^p(N),\Sigma}}. \end{array}\]

As $\varsigma_p:\cA^{(m),\ord,\dag}_{U^p(N),\Sigma^\ord} \ra \cA^{(m),\ord,\dag}_{U^p(N),\Sigma^\ord}$ is finite and flat we get a trace map
\[ \tr_{\varsigma_p} : \varsigma_{p,*} \cO_{\cA^{(m),\ord,\dag}_{U^p(N),\Sigma}} \lra \cO_{\cA^{(m),\ord,\dag}_{U^p(N),\Sigma}}. \]
Because $\partial  \cA^{(m),\ord,\dag}_{U^p(N),\Sigma}$ has the same support as 
\[ \cA^{(m),\ord,\dag}_{U^p(N),\Sigma} \times_{\varsigma_p, \cA^{(m),\ord,\dag}_{U^p(N),\Sigma}} \partial  \cA^{(m),\ord,\dag}_{U^p(N),\Sigma}, \]
this trace map restricts to a map
\[ \tr_{\varsigma_p} : \varsigma_{p,*} \cI_{\partial \cA^{(m),\ord,\dag}_{U^p(N),\Sigma}} \lra \cI_{\partial \cA^{(m),\ord,\dag}_{U^p(N),\Sigma}}. \]
(This is a consequence of the following fact: if $R$ is a noetherian ring, if $S$ is an $R$-algebra, finite and free as an $R$-module, and if $I$ and $J$ are ideals of $R$ and $S$ respectively with
\[ \sqrt{J} = \sqrt{IS}, \]
then the trace map $\tr_{S/R}$ maps $J$ to $I$. To see this we may reduce to the case $I=0$. In this case every element of $J$ is nilpotent and so has trace $0$.)

Composing $(\varsigma_p^*)^{-1}$ with $\tr_{\varsigma_p}$ we get
$G_n^{(m)}(\A^\infty)^{\ord,\times}$-equivariant maps
\[ \tr_F: \varsigma_{p,*} \Omega^j_{\cA^{(m),\ord,\dag}_{U^p(N),\Sigma}}(\log \infty)  \lra \Omega^j_{\cA^{(m),\ord,\dag}_{U^p(N),\Sigma}}(\log \infty) .\]
and
\[ \tr_F: \varsigma_{p,*} (\Omega^j_{\cA^{(m),\ord,\dag}_{U^p(N),\Sigma}}(\log \infty)\otimes \cI_{\partial \cA^{(m),\ord,\dag}_{U^p(N),\Sigma}})  \lra \Omega^j_{\cA^{(m),\ord,\dag}_{U^p(N),\Sigma}}(\log \infty) \otimes \cI_{\partial \cA^{(m),\ord,\dag}_{U^p(N),\Sigma}}.\]
We have
\[ \tr_F \circ \varsigma_p^{*} = p^{(n+2m)n[F^+:\Q]}. \]
This induces endomorphisms
\[ \tr_F \in \End( H^i(\cA^{(m),\ord,\dag}_{U^p(N),\Sigma}, \Omega^j_{\cA^{(m),\ord,\dag}_{U^p(N),\Sigma}}(\log \infty)\otimes \cI_{\partial \cA^{(m),\ord,\dag}_{U^p(N),\Sigma}})) \]
which commute with the action of $G_n(\A^\infty)^{\ord,\times}$ and satisfy
\[ \tr_F \circ \varsigma_p = p^{(n+2m)n[F^+:\Q]}. \]

We will also write $\Omega^j_{ \partial^{(s)} \cA^{(m),\ord,\dag}_{U^p(N),\Sigma} }$ for the sheaf of $j$-forms on $\partial^{(s)} \cA^{(m),\ord,\dag}_{U^p(N),\Sigma}$. The system $\{ \Omega^j_{\partial^{(s)} \cA^{(m),\ord,\dag}_{U^p(N),\Sigma}} \}$ over $\{  \partial^{(s)} \cA^{(m),\ord,\dag}_{U^p(N),\Sigma}\}$ has an action of $G_n^{(m)}(\A^\infty)^\ord$.

Furthermore if $\rho$ is a representation of $L_{n,(n)}$ on a finite dimensional $\Q_p$-vector space, there is a locally free sheaf $\cE^{\can,\dag}_{U^p(N),\Delta,\rho}$ (resp. $\cE^{\sub,\dag}_{U^p(N),\Delta,\rho}$) on $\cX^{\ord,\dag}_{U^p(N),\Delta}$ induced by $\cE^\can_{U^p(N,N'),\Delta',\rho}$ (resp. $\cE^\sub_{U^p(N,N'),\Delta',\rho}$) for any $N' \geq N$ and $\Delta' \in \cJ_n^{\tor}$ with $(\Delta')^\ord = \Delta$. This is canonically independent of the choices of $N'$ and $\Delta'$. The systems of sheaves $\{ \cE^{\can,\dag}_{U^p(N),\Delta,\rho}\}$ and $\{ \cE^{\sub,\dag}_{U^p(N),\Delta,\rho}\}$ over $\{ \cX^{\ord,\dag}_{U^p(N),\Delta}\}$ have  actions of $G_n(\A^\infty)^\ord$. There are equivariant identifications
\[ \cE^{\sub,\dag}_{U^p(N),\Delta,\rho} \cong \cE^{\can,\dag}_{U^p(N),\Delta,\rho} \otimes \cI_{\partial \cX^{\ord,\dag}_{U^p(N),\Delta}}, \]
where $\cI_{\partial \cX^{\ord,\dag}_{U^p(N),\Delta}}$ denotes the sheaf of ideals in $\cO_{\cX^{\ord,\dag}_{U^p(N),\Delta}}$ defining $\partial \cX^{\ord,\dag}_{U^p(N),\Delta}$.
For $g \in G_n(\A^\infty)^\ord$ the map
\[ g: g^* \cE^{\can,\dag}_{(U^p)'(N'),\Delta',\rho}  \lra \cE^{\can,\dag}_{U^p(N),\Delta,\rho} \]
is an isomorphism. (Because the same is true over $X_{U^p(N,N'),\Delta'}$ and hence over $X_{U^p(N,N'),\Delta'}^\dag$.) The inverse of $\varsigma_p^*$ gives maps
\[ \varsigma_{p,*} \cE^{\can,\dag}_{U^p(N),\Delta,\rho} \liso \cE^{\can,\dag}_{U^p(N),\Delta,\rho} \otimes_{\cO_{\cX^{\ord,\dag}_{U^p(N),\Delta}},\varsigma_p^{*}} \cO_{\cX^{\ord,\dag}_{U^p(N),\Delta}} \]
and
\[ \varsigma_{p,*} \cE^{\sub,\dag}_{U^p(N),\Delta,\rho} \liso \cE^{\can,\dag}_{U^p(N),\Delta,\rho} \otimes_{\cO_{\cX^{\ord,\dag}_{U^p(N),\Delta}},\varsigma_p^{*}} \cI_{\partial \cX^{\ord,\dag}_{U^p(N),\Delta}}. \]
Composing $(\varsigma_p^*)^{-1}$ with $\tr_{\varsigma_p}$ we get
$G_n^{(m)}(\A^\infty)^{\ord,\times}$-equivariant maps
\[ \tr_F: \varsigma_{p,*} \cE^{\can,\dag}_{U^p(N),\Delta,\rho}  \lra \cE^{\can,\dag}_{U^p(N),\Delta,\rho} .\]
and
\[ \tr_F: \varsigma_{p,*} \cE^{\sub,\dag}_{U^p(N),\Delta,\rho}  \lra \cE^{\sub,\dag}_{U^p(N),\Delta,\rho}.\]
We have
\[ \tr_F \circ \varsigma_p^{*} = p^{n^2[F^+:\Q]}. \]
This induces compatible endomorphisms
\[ \tr_F \in \End( H^0(\cX^{\ord,\dag}_{U^p(N),\Delta}, \cE^{\can,\dag}_{U^p(N),\Delta,\rho})) \]
and
\[ \tr_F \in \End( H^0(\cX^{\ord,\dag}_{U^p(N),\Delta}, \cE^{\sub,\dag}_{U^p(N),\Delta,\rho})) \]
which commute with the action of $G_n(\A^\infty)^{\ord,\times}$ and satisfy
\[ \tr_F \circ \varsigma_p = p^{n^2[F^+:\Q]}. \]

We define
\[ \begin{array}{l} H^i(\cA^{(m),\ord,\dag} , \Omega^j(\log \infty)\otimes \cI) = \\
\lim_{\substack{\lra \\ U^p,N,\Sigma}} H^i(\cA^{(m),\ord,\dag}_{U^p(N),\Sigma}, \Omega^j_{\cA^{(m),\ord,\dag}_{U^p(N),\Sigma}}(\log \infty)\otimes \cI_{\partial \cA^{(m),\ord,\dag}_{U^p(N),\Sigma}}) \end{array} \]
and
\[ H^i(\partial^{(s)}\cA^{(m),\ord,\dag} , \Omega^j) = \lim_{\substack{\lra \\ U^p,N,\Sigma}} H^i(\partial^{(s)} \cA^{(m),\ord,\dag}_{U^p(N),\Sigma}, \Omega^j_{\partial^{(s)} \cA^{(m),\ord,\dag}_{U^p(N),\Sigma}}) \]
and
\[ H^0(\cX^{\ord,\dag}, \cE^{\sub}_\rho) = \lim_{\substack{\lra \\ U^p,N,\Delta} }H^0(\cX^{\ord,\dag}_{U^p(N),\Delta}, \cE^{\sub,\dag}_{U^p(N),\Delta,\rho}), \]
smooth $G_n(\A^\infty)^\ord$-modules. We obtain an element
\[ \tr_F \in \End(H^i(\cA^{(m),\ord,\dag} , \Omega^j(\log \infty)\otimes \cI)) \]
which commutes with the $G_n(\A^\infty)^{\ord,\times}$-action and satisfies
\[ \tr_F \circ \varsigma_p= p^{(n+2m)n[F^+:\Q]}. \]
We also obtain an element
\[ \tr_F \in \End(H^0(\cX^{\ord,\dag} , \cE^\sub_\rho)) \]
which commutes with the $G_n(\A^\infty)^{\ord,\times}$-action and satisfies
\[ \tr_F \circ \varsigma_p= p^{n^2[F^+:\Q]}. \]

\begin{lem}\label{daginv}
There are natural isomorphisms
\[ \begin{array}{rl} & H^i(\cA^{(m),\ord,\dag}_{U^p(N),\Sigma}, \Omega^j_{\cA^{(m),\ord,\dag}_{U^p(N),\Sigma}}(\log \infty)\otimes \cI_{\partial \cA^{(m),\ord,\dag}_{U^p(N),\Sigma}}) \\ \liso & H^i(\cA^{(m),\ord,\dag} , \Omega^j(\log \infty)\otimes \cI)^{U^p(N)} \end{array} \]
and
\[ H^0(\cX^{\ord,\dag}_{U^p(N),\Delta}, \cE^{\sub,\dag}_{U^p(N),\Delta,\rho}) \liso H^0(\cX^{\ord,\dag}, \cE^{\sub}_\rho)^{U^p(N)}. \]
\end{lem}

\pfbegin Use lemmas \ref{l101}, \ref{cscor}, \ref{locmod}, \ref{l59} and \ref{daggaga}. 
\pfend
\newpage

\subsection{The ordinary locus of the minimal compactification as a dagger space.}

Suppose that $U^p$ is a neat open compact subgroup of $G_n(\A^{\infty,p})$ and that  
$N_2 \geq N_1 \geq 0$. We will write 
\[ X^{\mini,\ord,\dag}_{n,U^p(N_1,N_2)} \]
for the dagger space associated to $\cX^{\mini,\ord}_{n,U^p(N_1,N_2)}$ as described in the paragraph before lemma \ref{berth}. Then the system of dagger spaces $\{ X^{\mini,\ord,\dag}_{n,U^p(N_1,N_2)} \}$ has an action of $G_n(\A^\infty)^\ord$.

Recall from section \ref{mincomp} that, if $N_2' \geq N_2$, then the natural map
\[ \cX^{\ord, \mini}_{n,U^p(N_1,N_2')} \lra \cX^{\ord,\mini}_{n,U^p(N_1,N_2)} \]
restricts to an isomorphism
\[ \barX^{\ord, \mini}_{n,U^p(N_1,N_2')} \liso \barX^{\ord,\mini}_{n,U^p(N_1,N_2)} \]
and is etale in a neighbourhood of $\barX^{\ord,\mini}_{n,U^p(N_1,N_2')}$. 
It follows from lemma \ref{berth} that
\[ X^{\mini,\ord,\dag}_{n,U^p(N_1,N_2')} \lra X^{\mini,\ord,\dag}_{n,U^p(N_1,N_2)} \]
is an isomorphism. We will denote this dagger space simply
\[ X^{\mini,\ord,\dag}_{n,U^p(N_1)}. \]
The system of dagger spaces $\{ X^{\mini,\ord,\dag}_{n,U^p(N)} \}$ has an action of $G_n(\A^\infty)^{\ord}$. 

Let $\bare_{U^p(N_1,N_2)}$ denote the idempotent in 
\[ \left. \left( \bigoplus_{i=0}^\infty H^0(\barX_{n,U^p(N_1,N_2)}^\mini, \omega^{\otimes (p-1)i}) \right)\right/ (\hasse_{U^p(N_1,N_2)}-1) \]
which is $1$ on $\barX_{n,U^p(N_1,N_2)}^{\mini,\ord}$ and $0$ on
\[ \barX_{n,U^p(N_1,N_2)}^\mini - \barX_{n,U^p(N_1,N_2)}^{\mini,\nord} - \barX_{n,U^p(N_1,N_2)}^{\mini,\ord}. \]
Multiplying the terms of $\bare_{U^p(N_1,N_2)}$ by suitable powers of $\hasse_{U^p(N_1,N_2)}$, we may suppose that $\bare_{U^p(N_1,N_2)}$ lies in $H^0(\barX_{n,U^p(N_1,N_2)}^\mini, \omega^{\otimes (p-1)a})$ for any sufficiently large $a$, and that
\[ \bare_{U^p(N_1,N_2)}/\hasse_{U^p(N_1,N_2)} \in H^0(\barX_{n,U^p(N_1,N_2)}^\mini, \omega^{\otimes (p-1)(a-1)}). \]
Then 
\[ \barX^\ord_{n,U^p(N_1,N_2)} = \Spec \left. \left( \bigoplus_{i=0}^\infty H^0(\barX_{n,U^p(N_1,N_2)}^\mini, \omega^{\otimes (p-1)ai}) \right)\right/(\bare_{U^p(N_1,N_2)}-1).  \]
For $a$ sufficiently large we have $H^1(\cX^\mini_{n,U^p(N_1,N_2)},\omega^{\otimes (p-1)a})=(0)$. In that case we can lift $\bare_{U^p(N_1,N_2)}$ to a non-canonical element 
\[ e_{U^p(N_1,N_2)} \in H^0(\cX_{n,U^p(N_1,N_2)}^\mini, \omega^{\otimes (p-1)a}). \]
Let $\cX^\mini_{n,U^p(N_1,N_2)}[1/e_{U^p(N_1,N_2)}]$ denote the open subscheme of $\cX^\mini_{n,U^p(N_1,N_2)}$, where $e_{U^p(N_1,N_2)} \neq 0$. As $\omega^{\otimes (p-1)a}$ is ample, $\cX^\mini_{n,U^p(N_1,N_2)}[1/e_{U^p(N_1,N_2)}]$ is affine and so has the form
\[ \Spec \Z_{(p)}[T_1,...,T_s]/I \]
for some $s$ and $I$. It is normal and flat over $\Z_{(p)}$.

For $r \in p^{\Q_{\geq 0}}$ let $||\,\,\,||_r$ denote the norm on $\Z_{(p)}[T_1,...,T_s]$ defined by
\[ ||\sum_{\vec{i}} a_{\vec{i}}T^{\vec{i}} ||_r = \sup_{\vec{i}} |a_{\vec{i}}|_pr^{|\vec{i}|}, \]
where $\vec{i}$ runs over $\Z_{\geq 0}^s$ and $|(i_1,...,i_s)|=i_1+...+i_s$. We will write
$\Z_p\langle T_1,...,T_s\rangle_r$ for the completion of $\Z_{(p)}[T_1,...,T_s]$ with respect to $||\,\,\,||_r$. Thus $\Z_p\langle T_1,...,T_s \rangle_1$ is the $p$-adic completion of $\Z_{(p)}[T_1,...,T_s]$ and also the $p$-adic completion of $\Z_p\langle T_1,...,T_s\rangle_{r}$ for any $r\geq 1$. 
Set $\Q_p\langle T_1,...,T_s\rangle_r=\Z_p\langle T_1,...,T_s\rangle_r[1/p]$, the completion of $\Q[T_1,...,T_s]$ with respect to $||\,\,\,||_r$. In the case $r=1$ we will drop it from the notation. 
We will write
$\Z_p\langle T_1/r,...,T_s/r\rangle_1$ for the $||\,\,\,||_r$ unit-ball in $\Q_p\langle T_1,...,T_s \rangle_r$, i.e. for the set of power series 
\[ \sum_{\vec{i} \in \Z_{\geq 0}^s} a_{\vec{i}} \vec{T}^{\vec{i}} \]
where $a_{\vec{i}} \in \Q_p$, and $|a_{\vec{i}}|_p \leq r^{-|\vec{i}|}$ for all $\vec{i}$, and $|a_{\vec{i}}|_p r^{|\vec{i}|} \ra 0$ as $|\vec{i}|\ra \infty$. 
We will also write
\[ \Q_p\langle T_1,...,T_s\rangle^\dag = \bigcup_{r>1} \Q_p\langle T_1,...,T_s\rangle_r . \]

Let $\langle I\rangle_r$ denote the ideal of $\Z_p\langle T_1,...,T_s\rangle_r$ generated by $I$ and let
$\langle I \rangle_r'$ denote the intersection of $\langle I \rangle_1$ with $\Z_p\langle T_1,...,T_s\rangle_r$. 
Then $\langle I \rangle_1$ is the $p$-adic completion of $I$. Moreover
\[ \Z_p\langle T_1,...,T_s\rangle_1/\langle I \rangle_1 \]
is normal and flat over $\Z_p$, and
\[ \gX^{\mini,\ord}_{U^p(N_1)} = \Spf \Z_p\langle T_1,...,T_s  \rangle_1 /\langle I \rangle_1. \]
Note that
\[ \Z_{(p)}[T_1,...,T_s]/(I,p) \liso \Z_p\langle T_1,...,T_s\rangle_r/(\langle I\rangle_r, p) \]
for all $r \geq 1$. Thus $(\langle I \rangle_r,p)=(\langle I \rangle_r', p)$. 

We will also write $\langle I \rangle_{r,\Q_p}$ (resp. $\langle I \rangle_{r,\Q_p}'$) for the $\Q_p$ span of $\langle I \rangle_r$ (resp. $\langle I \rangle_r'$) in $\Q_p\langle T_1,...,T_s\rangle_r$. 
Then 
\[ \Sp \Q_p\langle T_1,...,T_s\rangle_1/\langle I \rangle_{1,\Q_p} \subset  \Sp \Q_p\langle T_1,...,T_s\rangle_r/\langle I \rangle_{r,\Q_p}'  \subset \Sp \Q_p\langle T_1,...,T_s\rangle_r/\langle I \rangle_{r,\Q_p} \]
are all affinoid subdomians of $X^{\mini,\an}_{U^p(N,N_2)}$, the rigid analytic space associated to $X_{U^p(N,N_2)}^\mini \times \Spec \Q_p$. Thus they are normal. Moreover $\Sp \Q_p\langle T_1,...,T_s\rangle_r/\langle I \rangle_{r,\Q_p}' $ and $\Sp \Q_p\langle T_1,...,T_s\rangle_r/\langle I \rangle_{r,\Q_p} -\Sp \Q_p\langle T_1,...,T_s\rangle_r/\langle I \rangle_{r,\Q_p}' $ forms an admissible open cover of $\Sp \Q_p\langle T_1,...,T_s\rangle_r/\langle I \rangle_{r,\Q_p}$. ($\Sp \Q_p\langle T_1,...,T_s\rangle_r/\langle I \rangle_{r,\Q_p}' $ is the union of the connected components of $\Sp \Q_p\langle T_1,...,T_s\rangle_r/\langle I \rangle_{r,\Q_p} $, which contain a component of $\Sp \Q_p\langle T_1,...,T_s\rangle_1/\langle I \rangle_{1,\Q_p}$. See Proposition 8 of section 9.1.4 of \cite{bgr}.) Moreover $\Sp \Q_p\langle T_1,...,T_s\rangle_1/\langle I \rangle_{1}$ is Zariski dense in $\Sp \Q_p\langle T_1,...,T_s\rangle_r/\langle I \rangle_{r,\Q_p}$. Indeed 
\[ X_{U^p(N,N_2)}^\an \cap \Sp \Q_p\langle T_1,...,T_s\rangle_1/\langle I \rangle_{1}\]
 is Zariski dense in $\Sp \Q_p\langle T_1,...,T_s\rangle_r/\langle I \rangle_{r,\Q_p}$, where $X^{\an}_{U^p(N,N_2)}$, the rigid analytic space associated to $X_{U^p(N,N_2)} \times \Spec \Q_p$.

If $1 \leq r' < r$ then
\[ \Sp \Q_p\langle T_1,...,T_s\rangle_{r'}/\langle I \rangle_{r',\Q_p}' \subset \Sp \Q_p\langle T_1,...,T_s\rangle_r/\langle I \rangle_{r,\Q_p}' \]
and
\[ \Sp \Q_p\langle T_1,...,T_s\rangle_{r'}/\langle I \rangle_{r',\Q_p} \subset \Sp \Q_p\langle T_1,...,T_s\rangle_r/\langle I \rangle_{r,\Q_p} , \]
and these are strict neighbourhoods.
The natural maps
\[ i_{r,r'}: \Q_p\langle T_1,...,T_s\rangle_r/\langle I\rangle_{r,\Q_p} \lra \Q_p\langle T_1,...,T_s\rangle_{r'}/\langle I \rangle_{r',\Q_p} \]
and
\[ i_{r,r'}': \Q_p\langle T_1,...,T_s\rangle_r/\langle I\rangle_{r,\Q_p}' \into \Q_p\langle T_1,...,T_s\rangle_{r'}/\langle I \rangle_{r',\Q_p}' \]
are completely continuous. The latter is an inclusion.
Moreover
\[ (i_{r,1}')^{-1} \Z_p\langle T_1,...,T_s\rangle_r/\langle I\rangle_1=\Z_p\langle T_1,...,T_s\rangle_r/\langle I\rangle_{r}'. \]

Also write $\langle I \rangle^\dag$ for the ideal of $\Q_p\langle T_1,...,T_s \rangle^\dag$ generated by $I$. Thus
\[ \langle I \rangle^\dag = \bigcup_{r>1} \langle I \rangle_{r,\Q_p} = \bigcup_{r>1} \langle I \rangle_{r,\Q_p}'. \]
Moreover
\[ \Q_p\langle T_1,...,T_s \rangle^\dag /\langle I \rangle^\dag = \lim_{\substack{\ra \\ r>1}} \Q_p\langle T_1,...,T_s \rangle_r/\langle I \rangle_{r,\Q_p} = \lim_{\substack{\ra \\ r>1}} \Q_p\langle T_1,...,T_s \rangle_r/\langle I \rangle_{r,\Q_p}', \]
and
\[ \cX^{\mini,\ord,\dag}_{U^p(N_1)} = \Sp \Q_p\langle T_1,...,T_s\rangle^\dag/\langle I \rangle^\dag. \]
(See for instance proposition 3.3.7 of \cite{lestum}. For the meaning of $\Sp$ in the context of dagger algebras see section 2.11 of \cite{gkcrelle}.) Thus we have the following lemma.

\begin{lem}\label{affinoid} $\cX^{\mini,\ord,\dag}_{U^p(N)}$ is affinoid. \end{lem}

We have a map
\[ \varsigma_p^*: \Z_p\langle T_1,...,T_s \rangle_1 /\langle I \rangle_1 \lra \Z_p \langle T_1,...,T_s\rangle_1 /\langle I \rangle_1 \]
such that
\begin{itemize}
\item $\varsigma_p^*(T_i) \equiv (T_i)^p \bmod p$,
\item and there exists an $r_1 \in p^{\Q_{>0}}$ such that for all $j=1,...,s$
the element $\varsigma_p^*(T_j)$ is in the image of $\Q_p\langle T_1,...,T_s\rangle_{r_1}/\langle I \rangle_{r_1}$.
\end{itemize}
Thus $(\varsigma_p^*(T_j)-T_j^p)/p \in \Z_p\langle T_1,...,T_s\rangle_r/\langle I\rangle_{r_1}'$, and so is the image of some element 
$G_j(\vec{T}) \in \Z_p\langle T_1,...,T_s\rangle_{r_1}$. We have 
\[ \varsigma_p^*(T_j)\equiv (T_j)^p+pG_j(T_1,...,T_s) \bmod \langle I \rangle_1. \]
This formula then to define a map $\varsigma_p^*: \Z_p[T_1,...,T_s] \ra \Z_p\langle T_1,...,T_s\rangle_{r_1}$ such that
\[ \begin{array}{ccc} \Z_p[T_1,...,T_s] & \stackrel{\varsigma_p^*}{\lra} &\Z_p\langle T_1,...,T_s\rangle_{r_1} \\
\da && \da \\ \Z_p\langle T_1,...,T_s\rangle_1 /\langle I \rangle_1 & \stackrel{\varsigma_p^*}{\lra} &\Z_p\langle T_1,...,T_s\rangle_{1}/\langle I \rangle_1 \end{array} \]
commutes.
Write $G_j(\vec{T})=\sum_{\vec{i}} g_{j,\vec{i}}\vec{T}^{\vec{i}}$. 
Choose $I_0 \in \Z^{>0}$ such that
\[ p^{-1}||G_j||_{r_1}  < (\sqrt{r_1})^{I_0} \]
for all $j=1,...,s$ and then choose $r_2 \in (1,\sqrt{r_1}) \cap p^\Q$ with 
\[ r_2^{I_0}< p. \]
If $r \in [1,r_2] \cap p^\Q$ we have 
\[ ||\varsigma_p^*(T_j)-(T_j)^p||_r < 1. \]
(Because if $|\vec{i}| \geq I_0$ then $||pg_{j,\vec{i}}\vec{T}^{\vec{i}}||_r\leq (1/p)||G_j||_{r_1} (r/r_1)^{I_0} < 1$, while
for $|\vec{i}| \leq I_0$ we have $||pg_{j,\vec{i}}\vec{T}^{\vec{i}}||_r\leq (1/p)r^{I_0} < 1$.) 
If $r \in (1,r_2] \cap p^\Q$ and $H \in \Z_p[ T_1,...,T_s]$ we deduce that
\[ ||\varsigma_p^*H - H(\vec{T}^p)||_r \leq r^{-p}||H||_{r^p}. \]
(One only need check this on monomials. Hence we only need check that if it is true for $H_1$ and $H_2$ then it is also true for $H_1H_2$. For this one uses the formula
\[ \begin{array}{r} \varsigma_p^*(H_1H_2)-(H_1H_2)(\vec{T}^p) = (\varsigma_p^*H_1-H_1(\vec{T}^p))(\varsigma_p^*H_2-H_2(\vec{T}^p))+ \\ (\varsigma_p^*H_1-H_1(\vec{T}^p)) H_2(\vec{T}^p)+
(\varsigma_p^*H_2-H_2(\vec{T}^p)) H_1(\vec{T}^p). )\end{array} \]
Hence, if $r \in (1,r_2] \cap p^\Q$ and $H \in \Z_p[ T_1,...,T_s]$ we deduce that
\[ ||\varsigma_p^*H||_r=||H||_{r^p}, \]
and so $\varsigma_p^*$ extends to an isometric homomorphism
\[ \varsigma_p^*: \Z_p\langle T_1/r^p,...,T_s/r^p\rangle_{1} \lra \Z_p\langle T_1/r,...,T_s/r\rangle_1.  \]
Modulo $p$ this map reduces to the Frobenius, which is finite and so
\[ \varsigma_p^*: \Z_p\langle T_1/r^p,...,T_s/r^p\rangle_{1} \lra \Z_p\langle T_1/r,...,T_s/r\rangle_1  \]
is finite. (See section 6.3.2 of \cite{bgr}.)
Thus we get an isometric, finite homomorphism between normal rings
\[ \varsigma_p^*: \Q_p\langle T_1,...,T_s\rangle_{r^p}/\langle I \rangle_{r^p,\Q_p}' \lra \Q_p\langle T_1,...,T_s\rangle_r/\langle I \rangle_{r,\Q_p}', \]
such that the diagram
\[ \begin{array}{ccc}  \Q_p\langle T_1,...,T_s\rangle_{r^p}/\langle I \rangle_{r^p,\Q_p}' &\stackrel{\varsigma_p^*}{\lra} &\Q_p\langle T_1,...,T_s\rangle_r/\langle I \rangle_{r,\Q_p}' \\ \da && \da \\
 \Q_p\langle T_1,...,T_s\rangle^\dag/\langle I \rangle^\dag &\stackrel{\varsigma_p^*}{\lra}& \Q_p\langle T_1,...,T_s\rangle^\dag /\langle I \rangle^\dag \\ \da && \da \\
 \Q_p\langle T_1,...,T_s\rangle_{1}/\langle I \rangle_{1,\Q_p} &\stackrel{\varsigma_p^*}{\lra}& \Q_p\langle T_1,...,T_s\rangle_1/\langle I \rangle_{1,\Q_p} \end{array} \]
 commutes.
 
 The map
 \[ \varsigma_p: \Sp \Q_p\langle T_1,...,T_s\rangle_r/\langle I \rangle_{r,\Q_p}' \lra \Sp \Q_p\langle T_1,...,T_s\rangle_{r^p}/\langle I \rangle_{r^p,\Q_p}' \]
 is compatible with the  map
 \[ \varsigma_p: X_{U^p(N,N_2)}^{\mini,\an} \lra X_{U^p(N,N_2-1)}^{\mini,\an}. \]
This latter map is finite, and away from the boundary is flat of degree $p^{n^2[F^+:\Q]}$. Thus the pre-image of $\Sp \Q_p\langle T_1,...,T_s\rangle_{r^p}/\langle I \rangle_{r^p,\Q_p}'$ has the form $\Sp B$ where $B$ is a normal, finite $\Q_p\langle T_1,...,T_s\rangle_{r^p}/\langle I \rangle_{r^p,\Q_p}'$ algebra, and we have a factorization
\[ \varsigma_p^*: \Q_p\langle T_1,...,T_s\rangle_{r^p}/\langle I \rangle_{r^p,\Q_p}' \lra B \lra \Q_p\langle T_1,...,T_s\rangle_r/\langle I \rangle_{r,\Q_p}'. \]
For $\wp$ a maximal ideal of $\Q_p\langle T_1,...,T_s\rangle_{r^p}/\langle I \rangle_{r^p,\Q_p}'$ corresponding to a point of $X_{U^p(N,N_2)}^{\an} \cap \Sp \Q_p\langle T_1,...,T_s\rangle_{1}/\langle I \rangle_{1}$ we see that 
\[ B/\gm = ( \Q_p\langle T_1,...,T_s\rangle_{1}/\langle I \rangle_{1})/\varsigma_p^* \gm = (\Q_p\langle T_1,...,T_s\rangle_r/\langle I \rangle_{r,\Q_p}')/\varsigma_p^*\gm. \]
Thus for a Zariski dense set of maximal ideals $\gm \in \Sp \Q_p\langle T_1,...,T_s\rangle_{r^p}/\langle I \rangle_{r^p,\Q_p}'$ the map
\[ B \lra \Q_p\langle T_1,...,T_s\rangle_{r}/\langle I \rangle_{r,\Q_p}'\]
becomes an isomorphism modulo $\gm$. Hence for any minimal prime $\wp$ of the ring $\Q_p\langle T_1,...,T_s\rangle_{r^p}/\langle I \rangle_{r^p,\Q_p}'$ we have
\[ B_\wp/\wp = (\Q_p\langle T_1,...,T_s\rangle_{r}/\langle I \rangle_{r,\Q_p}')_\wp/\wp . \]
(Choose bases over $A_\wp/\wp$, then this map being an isomorphism is equivalent to some matrix having full rank. For $\gm$ in a dense Zariski open set these bases reduce to bases modulo $\gm$. So modulo a Zariski dense set of $\gm$ this matrix has full rank, so it has full rank.)
As $B$ is normal and $\Q_p\langle T_1,...,T_s\rangle_{r}/\langle I \rangle_{r,\Q_p}'$ is finite over $B$, we see that
\[ B = \Q_p\langle T_1,...,T_s\rangle_{r}/\langle I \rangle_{r,\Q_p}' , \]
i.e. 
\[ \varsigma_p^{-1} \Sp \Q_p\langle T_1,...,T_s\rangle_{r^p}/\langle I \rangle_{r^p,\Q_p}'  = \Sp \Q_p\langle T_1,...,T_s\rangle_{r}/\langle I \rangle_{r,\Q_p}' . \]

The sheaf $\cE_{U^p(N_1,N_2),\rho}^\sub$ induces a coherent sheaf $\cE_{U^p(N_1),\rho}^{\sub,\dag}$ on $X_{n,U^p(N_1)}^{\mini,\ord,\dag}$, which does not depend on $N_2$.  It equals the push forward from any $X^{\ord,\dag}_{n,U^p(N_1),\Delta}$ of the sheaf $\cE_{U^p(N_1),\rho}^{\sub,\dag}$. 
The inverse system $\{ \cE_{U^p(N),\rho}^{\sub,\dag} \}$ is a system of coherent sheaves with $G_n(\A^{\infty})^\ord$-action on $\{X_{n,U^p(N_1)}^{\mini,\ord,\dag}\}$. The map
\[ \tr_F: \varsigma_{p,*} \cE_{U^p(N),\Delta,\rho}^{\sub,\dag} \lra \cE_{U^p(N),\Delta,\rho}^{\sub,\dag} \]
over $X^{\ord,\dag}_{n,U^p(N_1),\Delta}$ induces a map
\[ \tr_F: \varsigma_{p,*} \cE_{U^p(N),\rho}^{\sub,\dag} \lra \cE_{U^p(N),\rho}^{\sub,\dag} \]
over $X_{n,U^p(N_1)}^{\mini,\ord,\dag}$. This map does not depend on the choice of $\Delta$ and is $G_n(\A^\infty)^{\ord,\times}$-equivariant. It satisfies
\[ \tr_F \circ \varsigma_p = p^{n^2[F^+:\Q]}. \]
It induces a map
\[ \tr_F \in \End(H^0(X^{\mini,\ord,\dag}_{n,U^p(N)}, \cE^{\sub,\dag}_{U^p(N),\rho})) \]
also satisfying
\[ \tr_F \circ \varsigma_p = p^{n^2[F^+:\Q]}. \]

There are $G_n(\A^\infty)^\ord$ and $\tr_F$ equivariant isomorphisms 
\[ H^0(\cX^{\mini,\ord,\dag}_{U^p(N)}, \cE^{\sub,\dag}_{U^p(N),\rho}) \liso H^0(\cX^{\ord,\dag}_{U^p(N),\Delta}, \cE^{\sub,\dag}_{U^p(N),\rho}). \]
There are also natural $G_n(\A^\infty)^\ord$-equivariant embeddings
\[ H^0(\cX^{\mini,\ord,\dag}_{U^p(N)}, \cE^{\sub,\dag}_{U^p(N),\rho}) \into H^0(\gX^{\ord,\mini}_{U^p(N)}, \cE^\sub_{U^p(N),\rho}) \otimes_{\Z_p} \Q_p. \]

 We will set
 \[ H^0(\cX^{\mini,\ord,\dag},\cE_\rho^\sub)_{\barQQ_p}=\left( \lim_{\ra U^p,N} H^0(\cX_{U^p(N)}^{\mini,\ord,\dag}, \cE_{U^p(N),\rho}^{\sub,\dag})\right) \otimes_{\Q_p} \barQQ_p, \]
 a smooth $G_n(\A^\infty)^\ord$-module with an endomorphism $\tr_F$, which commutes with $G_n(\A^\infty)^{\ord,\times}$ and satisfies $\tr_F \circ \varsigma_p = p^{n^2[F^+:\Q]}$. From lemma \ref{daginv} and the first observation of the last paragraph, we see that
 \[ H^0(\cX^{\mini,\ord,\dag},\cE_\rho^\sub)_{\barQQ_p}^{U^p(N)}= H^0(\cX^{\mini,\ord,\dag}_{U^p(N)}, \cE^{\sub,\dag}_{U^p(N),\rho}). \]
There is a $G_n(\A^\infty)^\ord$-equivariant embedding
 \[ H^0(\cX^{\mini,\ord,\dag},\cE_\rho^\sub)_{\barQQ_p} \into H^0(\gX^{\ord,\mini},\cE_\rho^\sub)_{\barQQ_p}. \]

 Similarly the coherent sheaf $\cE_{U^p(N_1,N_2),\rho}^\sub$ gives rise to a coherent sheaf $\cE_{U^p(N_1,N_2),\rho}^{\sub,\an}$ on $X_{U^p(N_1,N_2)}^{\mini,\ord,\an}$. 
The inverse system $\{ \cE_{U^p(N_1,N_2),\rho}^{\sub,\an} \}$ is a system of coherent sheaves with $G_n(\A^{\infty})^\ord$-action on $\{X_{U^p(N_1,N_2)}^{\mini,\ord,\an}\}$. 

Let us make some of this more explicit.
The sheaf $\cE_{U^p(N_1,N_2),\rho}^{\sub,\an}$ restricted to $\Sp \Q_p\langle T_1,...,T_s\rangle_r/\langle I \rangle_{r,\Q_p}'$ corresponds to a finitely generated module $E_r$ over the ring $\Q_p\langle T_1,...,T_s\rangle_r/\langle I \rangle_{p,\Q_p}'$, which is naturally a Banach module. If $r' < r$ then 
\[ E_{r'} = E_r \otimes_{\Q_p\langle T_1,...,T_s\rangle_r/\langle I \rangle_{r,\Q_p}', i_{r,r'}'} \Q_p\langle T_1,...,T_s\rangle_{r'}/\langle I \rangle_{r',\Q_p}'. \]
Then the map $E_r \ra E_{r'}$, which we will also denote $i_{r,r'}'$, is completely continuous. 
The map $\tr_F$ extends to a continuous $\Q_p\langle T_1,...,T_s\rangle_{r^p}/\langle I \rangle_{r^p,\Q_p}'$ linear map
\[ t_r: E_r \lra E_{r^p} \]
for $r \in [1,r_2] \cap p^\Q$. We set
\[ E^\dag = \bigcup_{r>1} E_r, \]
so that
\[ E^\dag = H^0(\cX_{U^p(N_1)}^{\mini,\ord,\dag}, \cE_{U^p(N_1),\rho}^{\sub,\dag}). \]
We have that $\tr_F|_{E_r}=t_r$. As $t_r$ is continuous and $i_{r^p,r}'$ is completely continuous we see that 
\[ \tr_F: E_r \lra E_r \]
and that this map is completely continuous. Thus each $E_r$ admits slope decompositions for $\tr_F$ and hence by lemma \ref{slope} so does $E^\dag$ and $E^\dag \otimes \barQQ_p$.

   If $a \in \Q$ we thus have a well defined, finite dimensional subspace
\[ H^0(\cX_{U^p(N)}^{\mini,\ord,\dag}, \cE_{U^p(N),\rho}^{\sub,\dag})_{\barQQ_p,\leq a}  \subset H^0(\cX_{U^p(N)}^{\mini,\ord,\dag}, \cE_{U^p(N),\rho}^{\sub,\dag}) \otimes_{\Q_p} \barQQ_p. \]
(Defined with respect to $\tr_F$.)
We set
\[ H^0(\cX^{\mini,\ord,\dag},\cE_\rho^\sub)_{\barQQ_p,\leq a}=\lim_{\ra U^p,N} H^0(\cX_{U^p(N)}^{\mini,\ord,\dag}, \cE_{U^p(N),\rho}^{\sub,\dag})_{\barQQ_p, \leq a}, \]
so that there are $G_n(\A^\infty)^{\ord,\times}$-equivariant embeddings
\[ H^0(\cX^{\mini,\ord,\dag},\cE_\rho^\sub)_{\barQQ_p,\leq a} \subset H^0(\cX^{\mini,\ord,\dag},\cE_\rho^\sub)_{\barQQ_p} \into H^0(\gX^{\ord,\mini},\cE_\rho^\sub)_{\barQQ_p}. \]

We have proved the following lemma.

\begin{lem}\label{admis} $H^0(\cX^{\mini,\ord,\dag}, \cE_{\rho}^\sub)_{\barQQ_p,\leq a}$ is an admissible $G_n(\A^\infty)^{\ord,\times}$-module. \end{lem}

Combining this with corollary \ref{padicgalois} we obtain the following result.

\begin{cor}\label{ocgalois} Suppose that $\rho$ is a representation of $L_{n,(n)}$ over $\Z_{(p)}$, that $a \in \Q$ and that $\Pi$ is an irreducible $G_n(\A^{\infty})^{\ord,\times}$-subquotient of 
\[ H^0(\cX^{\mini,\ord,\dag},\cE_\rho^\sub)_{\barQQ_p,\leq a}. \]
Then there is a continuous semi-simple representation
\[ R_p({\Pi}):G_F \lra GL_{2n}(\barQQ_p) \]
with the following property: Suppose that $q\neq p$ is a rational prime above which $F$ and $\Pi$ are unramified, and suppose that $v|q$ is a prime of $F$. Then 
\[ \WD(R_{p}(\Pi)|_{G_{F_v}})^{\Fsemis} \cong \rec_{F_v}(\BC(\Pi_q)_v |\det|_v^{(1-2n)/2}), \]
where $q$ is the rational prime below $v$.
\end{cor}

We will next explain the consequences of these results for sheaves of differentials on $A^{(m),\ord,\dag}_{n,U^p(N_1,N_2),\Sigma}$. But we first need to record a piece of commutative algebra.
\begin{lem}\label{trcomp} Suppose that $A \ra B \ra C$ are reduced noetherian rings, with $B$ a finite flat $A$ module of rank $r_B$ and $C$ a finite flat $A$-module of rank $r_C$. Suppose also that the total ring of fractions of $C$ is finite flat over the total ring of fractions of $B$. Then $r_B|r_C$ and 
\[ (r_C/r_B)\tr_{B/A} = \tr_{C/A}: B \lra A. \]
\end{lem}

\pfbegin It suffices to check this after passing to total rings of fractions (i.e. localizations at the set of non-zero divisors). In this case $B$ is free over $A$ and $C$ is free over $B$ so the lemma is clear. \pfend

\begin{prop}\label{dagspect} There are representations $\rho^{i,j}_{m,s}$ of $L_{n,(n)}$ with the following property. 
If $(U^p(N),\Sigma) \in \cJ^{(m),\tor,\ord}$ and if $(U^p)'$ denotes the image of $U^p $ in $G_n(\A^{p,\infty})$, then there is a spectral sequence with first page
\[ \begin{array}{l} E_1^{i,j} = H^0(\cX^{\mini,\ord,\dag}_{(U^p)'(N)}, \cE^{\sub,\dag}_{(U^p)'(N),\rho_{m,s}^{i,j}}) \Rightarrow \\  H^{i+j}(\cA^{(m),\ord,\dag}_{U^p(N)},\Omega^s_{\cA^{(m),\ord,\dag}_{n, U^p(N)}}(\log \infty) \otimes \cI_{\partial \cA_{n,U^p(N)}^{(m),\ord,\dag}}). \end{array} \]
These spectral sequences are equivariant for the action of $G_n(\A^{\infty})^\ord$.  The map $\tr_F$ on the 
$H^{i+j}(\cA^{(m),\ord,\dag}_{ U^p(N)},\Omega^s_{A^{(m),\ord,\dag}_{n, U^p(N)}}(\log \infty) \otimes \cI_{\partial \cA_{U^p(N)}^{(m),\ord,\dag}})$ is compatible with the map $p^{nm[F:\Q]}\tr_F$ on the $H^0(\cX^{\mini,\ord,\dag}_{(U^p)'(N)}, \cE^{\sub,\dag}_{(U^p)'(N),\rho_{m,s}^{i,j}})$.
\end{prop}

\pfbegin We may replace $\Sigma$ by a refinement and so reduce to the case that there is a $\Delta$ with $((U^p)',\Delta) \in \cJ_n^{\tor,\ord}$ and $((U^p)'(N), \Delta) \leq (U^p(N),\Sigma)$. Let $\pi$ denote the map $\cA^{(m),\ord,\dag}_{U^p(N),\Sigma} \ra \cX^{\mini,\ord,\dag}_{(U^p)'(N)}$. Lemma \ref{spect} and lemma \ref{daggaga} tell us that there is a spectral sequence of coherent sheaves on $\cX^{\mini,\ord,\dag}_{U^p(N)}$ with first page
\[ E_1^{i,j} = \cE^{\sub,\dag}_{U^p(N),\rho_{m,s}^{i,j}} \Rightarrow R^{i+j}\pi_* (\Omega^s_{\cA^{(m),\ord,\dag}_{(U^p)'(N),\Sigma}}(\log \infty) \otimes \cI_{\partial \cA^{(m),\ord,\dag}_{(U^p)'(N),\Sigma}}).\]
The first assertion follows from lemma \ref{affinoid} and proposition 3.1 of \cite{gkcrelle} (which tell us that
\[ H^k(\cX^{\mini,\ord,\dag}, \cE^{\sub,\dag}_{\rho_{m,s}^{i,j}})=(0) \]
for $k>0$). 

Let $(U^p)'$ denote the image of $U^p$ in $G_n(\A^{p,\infty})$.
To avoid confusion we will write $\varsigma_{p,\cA}$ or $\varsigma_{p,\cX}$ depending on whether $\varsigma_p$ is acting on $\cA^{(m),\ord,\dag}_{(U^p)(N),\Sigma^\ord}$ or $\cX^{\ord,\dag}_{(U^p)'(N),\Delta^\ord}$. We will also factorize
$\varsigma_{p,\cA}$ as
\[ \cA^{(m),\ord,\dag}_{(U^p)(N),\Sigma^\ord} \stackrel{\Phi}{\lra} \varsigma_{p,\cX}^*\cA^{(m),\ord,\dag}_{(U^p)(N),\Sigma^\ord} 
\stackrel{\Psi}{\lra} \cA^{(m),\ord,\dag}_{(U^p)(N),\Sigma^\ord}. \]
Write $\pi'$ for the map
\[ \pi': \cA^{(m),\ord,\dag}_{(U^p)(N),\Sigma^\ord} \ra \cX^{\ord,\dag}_{(U^p)'(N),\Delta^\ord} \]
and $\pi''$ for the map
\[ \pi'': \varsigma_{p,\cX}^* \cA^{(m),\ord,\dag}_{(U^p)(N),\Sigma^\ord} \ra \cX^{\ord,\dag}_{(U^p)'(N),\Delta^\ord}. \]
The sheaf $\cE^\can_{\rho_{m,s}^{i,j}}$ on $\cX^{\ord,\dag}_{(U^p)'(N),\Delta^\ord}$ is
$R^i\pi_*' \cF_{j}$, where 
\[ \cF_{j} = \Omega^j_{\cX^{\ord,\dag}_{U^p(N),\Delta^\ord}} (\log \infty) \otimes \Omega^{s-j}_{\cA^{(m),\ord,\dag}_{(U^p)(N),\Sigma^\ord} /\cX^{\ord,\dag}_{(U^p)'(N),\Delta^\ord} }(\log \infty) .\]

To prove the last sentence of the lemma it suffices to show that the diagrams
\[ \begin{array}{ccc} \varsigma_{p,\cA,*}( \cF_j \otimes (\pi')^*  \cI_{\partial \cX^{\ord,\dag}_{(U^p)'(N),\Delta^\ord}}) & \lla & \cF_j \otimes_{\cO_{\cX^{\ord,\dag}_{(U^p)'(N),\Delta^\ord}}, \varsigma_{p,\cX}^*} \cI_{\partial \cX^{\ord,\dag}_{(U^p)'(N),\Delta^\ord}} 
 \\
\da && \da  \\
\varsigma_{p,\cA,*} (\cF_j \otimes \cI_{\partial \cA^{(m),\ord,\dag}_{(U^p)(N),\Sigma^\ord}}) & \lla& \cF_j \otimes_{\cO_{\cA^{(m),\ord,\dag}_{(U^p)(N),\Sigma^\ord}}, \varsigma_{p,\cA}^*} \cI_{\partial \cA^{(m),\ord,\dag}_{(U^p)(N),\Sigma^\ord}} 
 \end{array} \]
and
\[ \begin{array}{ccc}   \cF_j \otimes_{\cO_{\cX^{\ord,\dag}_{(U^p)'(N),\Delta^\ord}}, \varsigma_{p,\cX}^*} \cI_{\partial \cX^{\ord,\dag}_{(U^p)'(N),\Delta^\ord}} 
& \!\!\! \stackrel{1 \otimes p^{nm[F:\Q]} \tr}{\lra} & \!\!\!
\cF_j \otimes_{\cO_{\cX^{\ord,\dag}_{(U^p)'(N),\Delta^\ord}}} \cI_{\partial \cX^{\ord,\dag}_{(U^p)'(N),\Delta^\ord}} \\
\da && \da \\
 \cF_j \otimes_{\cO_{\cA^{(m),\ord,\dag}_{(U^p)(N),\Sigma^\ord}}, \varsigma_{p,\cA}^*} \cI_{\partial \cA^{(m),\ord,\dag}_{(U^p)(N),\Sigma^\ord}} 
& \stackrel{1 \otimes \tr}{\lra} & 
\cF_j \otimes_{\cO_{\cA^{(m),\ord,\dag}_{(U^p)(N),\Sigma^\ord}}} \cI_{\partial \cA^{(m),\ord,\dag}_{(U^p)(N),\Sigma^\ord}} \end{array} \]
commute. In the first diagram the upper horizontal map is the composite
\[ \begin{array}{cl} & \cF_j \otimes_{\cO_{\cX^{\ord,\dag}_{(U^p)'(N),\Delta^\ord}}, \varsigma_{p,\cX}^*} \cI_{\partial \cX^{\ord,\dag}_{(U^p)'(N),\Delta^\ord}} \\ = &\Psi_*((\Psi^* \cF_j) \otimes (\pi'')^*  \cI_{\partial \cX^{\ord,\dag}_{(U^p)'(N),\Delta^\ord}}) \\ \lra & \Psi_* \Phi_* 
(( \Phi^* \Psi^* \cF_j) \otimes (\Phi^* (\pi'')^*  \cI_{\partial \cX^{\ord,\dag}_{(U^p)'(N),\Delta^\ord}})) \\ =& \varsigma_{p,\cA,*}((\varsigma_{p,\cA}^* \cF_j) \otimes (\pi')^*  \cI_{\partial \cX^{\ord,\dag}_{(U^p)'(N),\Delta^\ord}}) \\ \stackrel{\varsigma_{p,\cA}^*}{\lra}& \varsigma_{p,\cA,*}( \cF_j \otimes (\pi')^*  \cI_{\partial \cX^{\ord,\dag}_{(U^p)'(N),\Delta^\ord}}), \end{array} \]
and the lower horizontal map is 
\[ \begin{array}{rcl} \cF_j \otimes_{\cO_{\cA^{(m),\ord,\dag}_{(U^p)(N),\Sigma^\ord}}, \varsigma_{p,\cA}^*} \cI_{\partial \cA^{(m),\ord,\dag}_{(U^p)(N),\Sigma^\ord}}  &\cong &\varsigma_{p,\cA,*} ((\varsigma_{p,\cA}^* \cF_j )\otimes \cI_{\partial \cA^{(m),\ord,\dag}_{(U^p)(N),\Sigma^\ord}}) \\ &\stackrel{\varsigma_{p,\cA}^*}{\lra}& \varsigma_{p,\cA,*} (\cF_j \otimes \cI_{\partial \cA^{(m),\ord,\dag}_{(U^p)(N),\Sigma^\ord}}). \end{array}\]
We see that the first square tautologically commutes. The second square commutes because the two maps
\[ p^{nm[F:\Q]} \tr: \Psi_* \cO_{\varsigma_{p,\cX}^* \cA^{(m),\ord,\dag}_{(U^p)(N),\Sigma^\ord}} \lra \cO_{\cA^{(m),\ord,\dag}_{(U^p)(N),\Sigma^\ord}} \]
and
\[ \Psi_* \cO_{\varsigma_{p,\cX}^* \cA^{(m),\ord,\dag}_{(U^p)(N),\Sigma^\ord}} \stackrel{\Phi^*}{\lra} \varsigma_{p,\cA,*} \cO_{\cA^{(m),\ord,\dag}_{(U^p)(N),\Sigma^\ord}}\stackrel{\tr}{\lra} \cO_{\cA^{(m),\ord,\dag}_{(U^p)(N),\Sigma^\ord}} \]
are equal. This in turn follows from lemma \ref{trcomp}.
\pfend

\begin{cor}\label{sd1} For all $i$ and $s$ the vector space $H^{i}(\cA^{(m),\ord,\dag}_{ U^p(N)},\Omega^s_{\cA^{(m),\ord,\dag}_{U^p(N)}}(\log \infty) \otimes \cI_{\partial \cA_{U^p(N)}^{(m),\ord,\dag}})$ admits slope decompositions for $\tr_F$. \end{cor}

We write 
\[ \begin{array}{rl} & H^i(\cA^{(m),\ord,\dag} , \Omega^j(\log \infty)\otimes \cI_{\partial \cA^{(m),\ord,\dag}})_{\leq a} \\ = &\lim_{\ra U^p,N,\Sigma} H^{i}(\cA^{(m),\ord,\dag}_{ U^p(N)},\Omega^s_{\cA^{(m),\ord,\dag}_{ U^p(N)}}(\log \infty) \otimes \cI_{\partial \cA_{U^p(N)}^{(m),\ord,\dag}})_{\leq a}. \end{array} \]
The next corollary now follows from the proposition and lemma \ref{slope}.

\begin{cor} 
For any $a \in \Q$ there is a $G_n(\A^\infty)^{\ord,\times}$-equivariant spectral sequence with first page
\[ E_1^{i,j} = H^0(\cX^{\mini,\ord,\dag}, \cE^\sub_{\rho_{m,s}^{i,j}})_{\leq a} \Rightarrow H^{i+j}(\cA^{(m),\ord,\dag},\Omega^s_{\cA^{(m),\ord,\dag}}(\log \infty) \otimes \cI_{\partial \cA^{(m),\ord,\dag}})_{\leq a}. \]
\end{cor}

Combining this with corollary \ref{ocgalois} we obtain the following corollary.

\begin{cor}\label{ocgalois2} Suppose that $\Pi$ is an irreducible $G_n(\A^{\infty})^{\ord,\times}$-subquotient of 
\[ H^{i}(\cA^{(m),\ord,\dag},\Omega^s_{\cA^{(m),\ord,\dag}}(\log \infty) \otimes \cI_{\partial \cA^{(m),\ord,\dag}})_{\leq a} \otimes_{\Q_p} \barQQ_p \]
for some $a \in \Q$.
Then there is a continuous semi-simple representation
\[ R_p({\Pi}):G_F \lra GL_{2n}(\barQQ_p) \]
with the following property: Suppose that $q\neq p$ is a rational prime above which $F$ and $\Pi$ are unramified, and suppose that $v|q$ is a prime of $F$. Then 
\[ \WD(R_{p}(\Pi)|_{G_{F_v}})^{\Fsemis} \cong \rec_{F_v}(\BC(\Pi_q)_v |\det|_v^{(1-2n)/2}), \]
where $q$ is the rational prime below $v$.
\end{cor}

\newpage

\subsection{Rigid cohomology.}\label{secrc}

Our main object of study will be the groups
\[ H^i_{c-\partial}(\barA^{(m),\ord}_{U^p(N),\Sigma})=\HH^i(\cA^{(m),\ord,\dag}_{U^p(N),\Sigma}, \Omega^\bullet_{\cA^{(m),\ord,\dag}_{U^p(N),\Sigma}}(\log \infty) \otimes \cI_{\partial \cA^{(m),\ord,\dag}_{U^p(N),\Sigma}}), \]
where $(U^p(N),\Sigma) \in \cJ_n^{(m),\tor,\ord}$.  This can be thought of as a sort of rigid cohomology of $\barA^{(m),\ord}_{U^p(N),\Sigma}$ with compact supports towards
the toroidal boundary, but not towards the non-ordinary locus. It seems plausible to us that this can be intrinsically attached to the pair $\barA^{(m),\ord}_{U^p(N)} \supset \partial \barA^{(m),\ord}_{U^p(N)}$. Hence our notation. However we will not prove this, so the reader is cautioned that our notation is nothing more than a short-hand, and the group $H^i_{c-\partial}(\barA^{(m),\ord}_{U^p(N),\Sigma})$ must be assumed to depend on the pair $\cA^{(m),\ord,\dag}_{U^p(N),\Sigma} \supset \partial \cA^{(m),\ord,\dag}_{U^p(N),\Sigma}$. We will also set
\[ H^i_{c-\partial}(\barA^{(m),\ord}) = \lim_{\substack{\lra \\ U^p,N,\Sigma}} H^i_{c-\partial}(\barA^{(m),\ord}_{U^p(N),\Sigma}). \]
It has a smooth action of $G_n(\A^{\infty})^\ord$. The maps 
\[ \tr_F: \varsigma_{p,*} \Omega^j_{\cA^{(m),\ord,\dag}_{U^p(N),\Sigma}} (\log \infty) \lra \Omega^j_{\cA^{(m),\ord,\dag}_{U^p(N),\Sigma}} (\log \infty) \]
induce endomorphisms 
\[ \tr_F \in \End( H^i_{c-\partial}(\barA^{(m),\ord}_{U^p(N),\Sigma})) \]
which commute with the action of $G_n(\A^\infty)^{\ord,\times}$ and satisfy
\[ \tr_F \circ \varsigma_p = p^{(n+2m)n[F^+:\Q]}. \] 

\begin{lem}\label{rigcoh}
There are natural isomorphisms
\[ H^i_{c-\partial}(\barA^{(m),\ord}_{U^p(N),\Sigma}) \liso H^i_{c-\partial}(\barA^{(m),\ord})^{U^p(N)}. \]
\end{lem}

\pfbegin Use lemmas \ref{l101}, \ref{cscor}, \ref{locmod} and \ref{daggaga}. 
\pfend

We will compute the group $H^i_{c-\partial}(\barA^{(m),\ord}_{U^p(N),\Sigma})$ in two ways. The first way will be in terms of $p$-adic cusp forms and will allow us to attach Galois representations to irreducible $G_n(\A^{\infty})^{\ord,\times}$-sub-quotients of $H^i_{c-\partial}(\barA^{(m),\ord}) \otimes_{\Q_p} \barQQ_p$. The second way will be geometrical, in terms of the stratification of the boundary. In this second approach the cohomology of the locally symmetric spaces associated to $L_{n,(n),\lin}^{(m)}$ will appear. 

Here is our first calculation.
\begin{lem}\label{lhdr} The vector spaces $H^i_{c-\partial}(\barA^{(m),\ord}_{U^p(N),\Sigma})$ admit slope decompositions for $\tr_F$. If moreover we set
\[ H^i_{c-\partial}(\barA^{(m),\ord})_{\leq a} = \lim_{\substack{\lra \\ U^p,N,\Sigma}} H^i_{c-\partial}(\barA^{(m),\ord}_{U^p(N),\Sigma})_{\leq a}, \] 
then there is a $G_n(\A^\infty)^{\ord,\times}$-spectral sequence with first page
\[   E^{i,j}_1 = H^i(\cA^{(m),\ord,\dag}, \Omega^j(\log \infty) \otimes \cI_{\partial \cA^{(m),\ord,\dag}})_{\leq a} \Rightarrow  H^{i+j}_{c-\partial}(\barA^{(m),\ord})_{\leq a}.  \]
\end{lem}

\pfbegin This follows from lemma \ref{slope}, corollary \ref{sd1} and the spectral sequence
\[   E^{i,j}_1 = H^i(\cA^{(m),\ord,\dag}_{U^p(N),\Sigma}, \Omega^j_{\cA^{(m),\ord,\dag}_{U^p(N),\Sigma}}(\log \infty) \otimes \cI_{\partial \cA^{(m),\ord,\dag}_{U^p(N),\Sigma}}) \Rightarrow  H^{i+j}_{c-\partial}(\barA^{(m),\ord}_{U^p(N),\Sigma}).  \]
\pfend

And here is our second calculation.
\begin{lem} There are $G_n(\A^\infty)^{\ord,\times}$-equivariant spectral sequences with first page
\[ E_1^{i,j}=H^i_{\rig}(\partial^{(j)} \barA^{(m),\ord}_{U^p(N),\Sigma}) \Rightarrow H^{i+j}_{c-\partial}(\barA^{(m),\ord}_{U^p(N),\Sigma}). \]
Moreover the action of Frobenius on the left hand side is compatible with the action of $\varsigma_p$ on the right hand side.
\end{lem}

\pfbegin
By lemmas \ref{log} and \ref{daggaga} the group $H^{i}_{c-\partial}(\barA^{(m),\ord}_{U^p(N),\Sigma})$ is isomorphic to the hyper-cohomology of the double complex
\[ \HH^i(\cA^{(m),\ord,\dag}_{U^p(N),\Sigma}, i^{(\bullet)}_* \Omega^\bullet_{\partial^{(\bullet)} \cA^{(m),\ord,\dag}_{U^p(N),\Sigma}}), \]
and so there is a spectral sequence with first page
\[ E_1^{i,j} =\HH^i(\partial^{(j)} \cA^{(m),\ord,\dag}_{U^p(N),\Sigma},  \Omega^\bullet_{\partial^{(j)} \cA^{(m),\ord,\dag}_{U^p(N),\Sigma}}) \Rightarrow H^{i}_{c-\partial}(\barA^{(m),\ord}_{U^p(N),\Sigma}). \]
However, by lemma \ref{dagrig} and the quasi-projectivity of $\partial^{(j)} \cA^{(m),\ord}_{U^p(N),\Sigma}$,
we see that there are $G_n^{(m)}(\A^\infty)^\ord$-equivariant isomorphisms
\[ \HH^i(\partial^{(j)} \cA^{(m),\ord,\dag}_{U^p(N),\Sigma},  \Omega^\bullet_{\partial^{(j)} \cA^{(m),\ord,\dag}_{U^p(N),\Sigma}}) \cong H^i_{\rig}(\partial^{(j)} \barA^{(m),\ord}_{U^p(N),\Sigma}), \]
and that under this identification $\varsigma_p$ corresponds to Frobenius (because $\varsigma_p$ equals Frobenius on the special fibre).
\pfend

\begin{cor} $H^{i}_{c-\partial}(\barA^{(m),\ord}_{U^p(N),\Sigma})$ is finite dimensional. Moreover 
\[ H^{i}_{c-\partial}(\barA^{(m),\ord}_{U^p(N),\Sigma})=H^{i}_{c-\partial}(\barA^{(m),\ord}_{U^p(N),\Sigma})_{\leq a}, \]
for some $a$, and so
\[ H^{i}_{c-\partial}(\barA^{(m),\ord}) = \bigcup_{a \in \Q} H^{i}_{c-\partial}(\barA^{(m),\ord})_{\leq a}. \]
 \end{cor}

\pfbegin The first assertion follows from the lemma and theorem 3.1 of \cite{berthfin}. The second assertion follows because $\tr_F \circ \varsigma_p = p^{n(n+2m)[F^+:\Q]}$ and so by the first part $\tr_F$ must be an automorphism of $H^{i}_{c-\partial}(\barA^{(m),\ord}_{U^p(N),\Sigma})$.\pfend

Combining this with corollary \ref{ocgalois2} and lemma \ref{lhdr} we obtain the following corollary.

\begin{cor}\label{riggalois} Suppose that $\Pi$ is an irreducible $G_n(\A^{\infty})^{\ord,\times}$-subquotient of 
\[ H^{i}_{c-\partial}(\barA^{(m),\ord}) \otimes_{\Q_p} \barQQ_p. \]
Then there is a continuous semi-simple representation
\[ R_p({\Pi}):G_F \lra GL_{2n}(\barQQ_p) \]
with the following property: Suppose that $q\neq p$ is a rational prime above which $F$ and $\Pi$ are unramified, and suppose that $v|q$ is a prime of $F$. Then 
\[ \WD(R_{p}(\Pi)|_{G_{F_v}})^{\Fsemis} \cong \rec_{F_v}(\BC(\Pi_q)_v |\det|_v^{(1-2n)/2}), \]
where $q$ is the rational prime below $v$.
\end{cor}

\begin{cor} \label{toptorig} The eigenvalues of $\varsigma_p$ on $H^{i}_{c-\partial}(\barA^{(m),\ord})_{\barQQ_p}$ are Weil $p^w$-numbers for some $w \in \Z_{\geq 0}$ (depending on the eigenvalue). We will write 
\[ W_0 H^{i}_{c-\partial}(\barA^{(m),\ord})_{\barQQ_p} \]
for the subspace of $H^{i}_{c-\partial}(\barA^{(m),\ord})_{\barQQ_p}$ spanned by generalized eigenspaces of $\varsigma_p$ with eigenvalue a $p^0$-Weil number.

For $i>0$ there is a $G_n(\A^{\infty})^\ord$-equivariant isomorphism
\[ \lim_{\substack{\lra \\ U^p,N,\Sigma}} H^{i}(|\cS(\partial \barA^{(m),\ord}_{U^p(N),\Sigma })|, \barQQ_p) \liso W_0 H^{i+1}_{c-\partial}(\barA^{(m),\ord})_{\barQQ_p}. \]
(For $i=0$ there is still a surjection.) 
\end{cor}

\pfbegin 
By theorem 2.2 of \cite{chiar}, the eigenvalues of the Frobenius endomorphism on
$H^i_{\rig}(\partial^{(j)} \barA^{(m),\ord}_{U^p(N),\Sigma})$ 
are all Weil $p^w$-numbers for some $w \in \Z_{\geq i}$ (depending on the eigenvalue). The first part of the corollary follows.

It follows moreover that $W_0  H^{i}_{c-\partial}(\barA^{(m),\ord}_{U^p(N),\Sigma})_{\barQQ_p}$ is the cohomology of the complex
\[ ... 
\lra H^0_\rig(\partial^{(i)} \barA^{(m),\ord}_{U^p(N),\Sigma},\barQQ_p) \lra H^0_\rig(\partial^{(i+1)} \barA^{(m),\ord}_{U^p(N),\Sigma},\barQQ_p) \lra  ... \] 
However by proposition 8.2.15 of \cite{lestum}
\[ H^0_\rig(\partial^{(i)} \barA^{(m),\ord}_{U^p(N),\Sigma },\barQQ_p) \cong   \barQQ_p^{\pi_0(\partial^{(i)} \barA^{(m),\ord}_{U^p(N),\Sigma } \times \Spec \barFF_p)}, \]
and so the cohomology of the above complex becomes
\[ \begin{array}{ll} \ker (H^0_\rig(\barA^{(m),\ord}_{U^p(N),\Sigma},\barQQ_p) \lra H^0(|\cS(\partial \barA^{(m),\ord}_{U^p(N),\Sigma })|, \barQQ_p)) & {\rm in \,\,\, degree \,\,\,} 0 \\ H^0(|\cS(\partial \barA^{(m),\ord}_{U^p(N),\Sigma })|, \barQQ_p)/\Im H^0_\rig(\barA^{(m),\ord}_{U^p(N),\Sigma},\barQQ_p) & {\rm in \,\,\, degree \,\,\,} 1 \\ H^{i-1}(|\cS(\partial \barA^{(m),\ord}_{U^p(N),\Sigma })|, \barQQ_p) & {\rm in \,\,\, degree \,\,\,} i>1. \end{array} \]
The last part of the corollary follows.
\pfend

According to the discussion at the end of section \ref{torcomp} we see that there are $G_n(\A^\infty)^\ord$-equivariant open embeddings
\[ \gT^{(m),\ord}_{U^p(N),=n} \into |\cS(\partial \barA^{(m),\ord}_{U^p(N),\Sigma })|.\]
Thus the following corollary follows by applying lemma \ref{intcoh} and corollary \ref{indu}.

\begin{cor} For $i>0$,
\[ H^i_\Int(\gT^{(m),\ord}_{=n}, \barQQ_p) \cong \Ind^{G_n^{(m)}(\A^{p,\infty})}_{P_{n,(n)}^{(m),+}(\A^{p,\infty})} H^i_\Int(\gT_{(n)}^{(m)},\barQQ_p)^{\Z_p^\times} \]
 is a $G_n(\A^\infty)^\ord$ subquotient of $W_0 H^{i+1}_{c-\partial}(\barA^{(m),\ord}_{n})_{\barQQ_p}$. \end{cor}

Combining this proposition with corollary \ref{riggalois} (and using lemma \ref{indbc}) we deduce the following consequence.

\begin{cor}\label{compgalois} Suppose that $i>0$ and that $\pi$ is an irreducible $L_{n,(n),\lin}(\A^\infty)$-subquotient of $H^i_{\Int}(\gT_{(n)}^{(m)},\barQQ_p)$. Then there is a continuous semi-simple representation
\[ R_{p}({\pi}):G_F \lra GL_{2n}(\barQQ_p) \]
such that, if $q \neq p$ is a rational prime above which $\pi$ and $F$ are unramified and if $v|q$ is a prime of $F$, then $R_{p}({\pi})$ is unramified at $v$ and
\[ R_p(\pi)|_{W_{F_v}}^\semis \cong \rec_{F_v}(\pi_v |\det|_v^{(1-n)/2}) \oplus \rec_{F_{{}^cv}}(\pi_{{}^cv}|\det|_{{}^cv}^{(1-n)/2})^{\vee,c} \epsilon_p^{1-2n}. \]
\end{cor}

Combining this with corollary \ref{toptocone} we obtain the following result.

\begin{cor}\label{mt1}
Suppose that $n>1$, that $\rho$ is an irreducible algebraic representation of $L_{n,(n),\lin}$ on a finite dimensional $\C$-vector space, and that
$\pi$ is a cuspidal automorphic representation of $L_{n,(n),\lin}(\A)$ so that $\pi_\infty$ has the same infinitesimal character as $\rho^\vee$. Then, for all sufficiently large integers $N$, 
there is a continuous, semi-simple representation
\[ R_{p,\imath}(\pi ,N): G_F \lra GL_{2n}(\barQQ_p) \]
such that if $q \neq p$ is a prime above which $\pi$ and $F$ are unramified and if $v|q$ is a prime of $F$ then $R_{p,\imath}(\pi,N)$ is unramified at $v$ and
\[ R_{p,\imath}(\pi,N)|_{W_{F_v}}^\semis = \imath^{-1}\rec_{F_v}(\pi_v |\det|_v^{(1-n)/2})  \oplus (\imath^{-1}\rec_{F_{{}^cv}}(\pi_{{}^cv} |\det|_{{}^cv}^{(1-n)/2}))^{\vee,c} \epsilon_p^{1-2n-2N}. \]
\end{cor}

\pfbegin Take
\[ R_{p,\imath}(\pi,N)= R_p(\imath^{-1} (\pi^\infty ||\det||^N)) \otimes \epsilon_p^{-N}. \]
\pfend

\newpage

\section{Galois representations.}

In order to improve upon corollary \ref{mt1} it is necessary to apply some simple group theory. To this end, let $\Gamma$ be a topological group and let $\gF$ be a dense set of elements of $\Gamma$. Let $k$ be an algebraically closed, topological field of characteristic $0$ and let $d \in \Z_{>0}$. 

Let 
\[ \mu: \Gamma \lra k^\times \]
be a continuous homomorphism such that $\mu(f)$ has infinite order for all $f \in \gF$. For $f \in \gF$ let $\cE_f^1$ and $\cE_f^2$ be two $d$-element multisets of elements of $k^\times$. Let $\cM$ be an infinite subset of $\Z$. For $m \in \cM$ suppose that
\[ \rho_m: \Gamma \lra GL_{2d}(k) \]
be a continuous semi-simple representation such that for every $f \in \gF$ the multi-set of roots of the characteristic polynomial of $\rho_m(f)$ equals
\[ \cE_f^1 \amalg \cE_f^2 \mu(f)^m. \]

Suppose that $\cM'$ is a finite subset of $\cM$. Let $G_{\cM'}$ denote the Zariski closure in $\G_m \times GL_{2d}^{\cM'}$ of the image of 
\[ \mu \oplus \bigoplus_{m \in \cM'} \rho_m. \]
It is a, possibly disconnected, reductive group. There is a natural continuous homomorphism
\[ \rho_{\cM'}=\mu \times \prod_{m \in \cM'} \rho_m : \Gamma \lra G_{\cM'}(k). \]
Note that $\rho_{\cM'} \gF$ is Zariski dense in $G_{\cM'}$. 
We will use $\mu$ for the character of $G_{\cM'}$ which is projection to $\G_m$. For $m \in \cM'$ we will let
\[ R_m: G_{\cM'} \lra GL_{2d} \]
denote the projection to the factor indexed by $m$. 

\begin{lem} For every $g \in G_{\cM'}(k)$ there are two $d$-element multisets $\Sigma_g^1$ and $\Sigma_g^2$ of elements of $k^\times$ such that for every $m \in \cM'$ the multiset of roots of the characteristic polynomial of $R_m(g)$ equals
\[ \Sigma^1_g \amalg \Sigma^2_g\mu(g)^m. \]
\end{lem}

\pfbegin
It suffices to show that the subset of $k^\times \times GL_{2d}^{\cM'}(k)$ consisting of elements $(t,(g_m)_{m \in \cM'})$ such that there are $d$-element multisets $\Sigma^1$ and $\Sigma^2$ of elements of $k^\times$ so that for all $m \in \cM'$ the multiset of roots of the characteristic polynomial of $g_m$ equals $\Sigma^1 \amalg \Sigma^2t^m$, is Zariski closed. Let $\Pol_{2d}$ denote the space of monic polynomials of degree $2d$. It even suffices to show that the subset $X$ of $k^\times \times \Pol_{2d}^{\cM'}(k)$ consisting of elements $(t,(P_m)_{m \in \cM'})$ such that there are $d$-element multisets $\Sigma^1$ and $\Sigma^2$ of elements of $k$ so that for all $m \in \cM'$ the multiset of roots of $P_m$ equals $\Sigma^1 \amalg \Sigma^2t^m$, is Zariski closed.

There is a natural finite map
\[ \begin{array}{rcl} \pi: \Aff^{2d} &\lra &\Pol_{2d} \\ (\alpha_i) & \longmapsto & \prod_i (T-\alpha_i).
\end{array} \]
If 
\[ (\sigma_m)\in S_{2d}^{\cM'}, \]
where $S_{2d}$ denotes the symmetric group on $2d$ letters,
define $V_{(\sigma_m)}$ to be the set of 
\[ (t,(a_{m,i}) ) \in \G_m \times (\Aff^{2d})^{\cM'} \]
such that, for all $m,m' \in \cM'$ we have 
\[ a_{m,\sigma_m i}=a_{m',\sigma_{m'} i} \]
if $i=1,..,d$ and
\[ a_{m,\sigma_m i}=a_{m',\sigma_{m'} i} t^{m'-m}\]
if $i=d+1,...,2d$. Then $V_{(\sigma_m)}$ is closed in $\G_m \times (\Aff^{2d})^{\cM'}$. Moreover 
\[ X = \bigcup_{(\sigma_m) \in S_{2d}^{\cM'}}(1 \times \pi^{\cM'}) V_{(\sigma_m)}. \]
The lemma now follows from the finiteness of $1 \times \pi^{\cM'}$.
\pfend

\begin{cor} If $\emptyset \neq \cM' \subset \cM''$ are finite subsets of $\cM$ then $G_{\cM''} \iso G_{\cM'}$. \end{cor}

\pfbegin Suppose that $g$ is in the kernel of the natural map
\[ G_{\cM''} \onto G_{\cM'}. \]
Then for all $m \in \cM''$ the only eigenvalue of $R_m(g)$ is $1$. Thus $g$ must be unipotent. However $\ker (G_{\cM''} \onto G_{\cM'})$ is reductive and so must be trivial.
\pfend

Thus we can write $G$ for $G_{\cM'}$ without danger of confusion.

\begin{cor}\label{113} For every $g \in G(k)$ there are two $d$-element multisets $\Sigma_g^1$ and $\Sigma_g^2$ of elements of $k^\times$ such that for every $m \in \cM$ the multiset of roots of the characteristic polynomial of $R_m(g)$ equals
\[ \Sigma^1_g \amalg \Sigma^2_g\mu(g)^m. \]
Moreover if $\mu(g)$ has infinite order then the multisets $\Sigma^1_g$ and $\Sigma^2_g$ are unique. 
\end{cor}

\pfbegin Choose non-empty finite subsets 
\[ \cM_1' \subset \cM_2' \subset ... \subset \cM \]
with
\[ \cM = \bigcup_{i=1}^\infty \cM_i'. \]
For each $i$ we can find two $d$-element multisets $\Sigma_{g,i}^1$ and $\Sigma_{g,i}^2$ of elements of $k^\times$ such that for every $m \in \cM_i'$ the multiset of roots of the characteristic polynomial of $R_m(g)$ equals
\[ \Sigma^1_{g,i} \amalg \Sigma^2_{g,i}\mu(g)^m. \]

Let $m_1 \in \cM_1'$ and let $\Sigma$ denote the set of eigenvalues of $R_{m_1}(g)$. Then,  for every $i$, the multiset $\Sigma^1_{g,i}$ consists of elements of $\Sigma$ and the multiset $\Sigma^2_{g,i}$ consists of elements of $\Sigma \mu(g)^{-m_1}$. Thus there are only finitely many possibilities for the pair of multisets $(\Sigma_{g,i}^1,\Sigma_{g,i}^2)$ as $i$ varies. Hence some such pair $(\Sigma_{g}^1,\Sigma_{g}^2)$ occurs infinitely often. This pair satisfies the requirements of the lemma. 

For uniqueness suppose that $\Sigma_g^{1,\prime}$ and $\Sigma_g^{2,\prime}$ is another such pair of multisets. Choose $m \in \cM$ with $\mu(g)^m \neq \alpha/\beta$ for any $\alpha,\beta \in \Sigma_g^1 \amalg \Sigma_g^2 \amalg \Sigma_g^{1,\prime} \amalg \Sigma_g^{2,\prime}$. Then the equality
\[  \Sigma^1_{g} \amalg \Sigma^2_{g}\mu(g)^m= \Sigma^{1,\prime}_{g} \amalg \Sigma^{2,\prime}_{g}\mu(g)^m \]
implies that $\Sigma^{1,\prime}_g=\Sigma_g^1$ and $\Sigma^{2,\prime}_g=\Sigma_g^2$.
\pfend

The connected component $Z(G)^0$ of the centre of $G$ is a torus.

\begin{lem} The character $\mu$ is non-trivial on $Z(G)^0$. \end{lem}

\pfbegin If $\mu$ were trivial on $Z(G^0)^0$ then it would be trivial on $G^0$ (because $G^0/Z(G^0)^0$ is semi-simple), and so $\mu$ would have finite order, a contradiction. Thus $\mu|_{Z(G^0)^0}$ is non-trivial. 

The space
\[ X^*(Z(G^0)^0) \otimes_\Z \Q \]
is a representation of the finite group $G/G^0$ and we can decompose
\[ X^*(Z(G^0)^0) \otimes_\Z \Q = (X^*(Z(G)^0) \otimes_\Z \Q) \oplus Y \]
where $Y$ is a $\Q[G/G^0]$-module with 
\[ Y^{G/G^0}=(0). \]
But 
\[ \mu|_{Z(G^0)^0} \in X^*(Z(G^0)^0)^{G/G^0} \subset X^*(Z(G)^0) \otimes_\Z \Q\]
is non-trivial, and so $\mu|_{Z(G)^0}$ is non-trivial. \pfend

For $m \in \cM$ let $\gX_m$ denote the $2d$-element multiset of characters of $Z(G)^0$ which occur in $R_m$ (taken with their multiplicity). If $g \in G$ then we will write $\gY(g)_m$ for the $2d$-element multiset of pairs $(\chi,a)$, where $\chi$ is a character of $Z(G)^0$ and $a$ is a root of the characteristic polynomial of $g$ acting on the $\chi$ eigenspace of $Z(G)^0$ in $R_m$. (The pair $(\chi,a)$ occurs with the same multiplicity as $a$ has as a root of the characteristic polynomial of $g$ acting on the $\chi$-eigenspace of $R_m$.) 

If $\gY \subset \gY(g)_m$ and if $\psi \in X^*(G)$ then we will set
\[ \gY\psi=\{(\chi\psi,a\psi(g)): \,\, (\chi,a) \in \gY\}. \]
We warn the reader that this depends on $g$ and not just on the set $\gY$.

\begin{lem} Suppose that $T/k$ is a torus and that $\gX$ is a finite set of non-trivial characters of $T$.
Let $A$ be a finite subset of $k^\times$. Then we can find $t \in T(k)$ such that $\chi(t) \neq a$ for all $\chi \in \gX$ and $a \in A$. 
\end{lem}

\pfbegin Let $(\,\, ,\,\,)$ denote the usual perfect pairing
\[ X^*(T) \times X_*(T) \lra \Z. \]
We can find $\nu \in X_*(T)$ such that $(\chi,\nu) \neq 0$ for all $\chi \in\gX$. Thus we are reduced to the case $T=\G_m$, in which case we may take $t$ to be any element of $k^\times$ that does not lie in the divisible hull of the subgroup $H$ of $k^\times$ generated by $A$. (For example we can take $t$ to be a rational prime such that all elements of a finite set of generators of $H \cap \Q^\times$ are units at $t$.) \pfend

\begin{cor} Suppose that $T/k$ is a torus and that $\gX$ is a finite set of characters of $T$. Then we can find $t \in T(k)$ such that if $\chi \neq \chi'$ lie in $\gX$ then
\[ \chi(t) \neq \chi'(t). \]
\end{cor}

\begin{lem} 
If $m,m',m'' \in \cM$, then we can decompose
\[ \gY(g)_m = \gY(g)_{m,m',m''}^1 \amalg \gY(g)_{m,m',m''}^2 \]
into two $d$-element multisets, such that
\[ \gY(g)_{m'} = \gY(g)_{m,m',m''}^1 \amalg \gY(g)_{m,m',m''}^2 \mu^{m'-m} \]
and
\[ \gY(g)_{m''} = \gY(g)_{m,m',m''}^1 \amalg \gY(g)_{m,m',m''}^2 \mu^{m''-m}. \]
If $\mu^{m-m'} \neq \chi/\chi'$ for all $\chi, \chi' \in \gX_m$ then the equation
\[ \gY(g)_{m'} = \gY(g)_{m,m',m''}^1 \amalg \gY(g)_{m,m',m''}^2 \mu^{m'-m} \]
uniquely determines this decomposition. 
 \end{lem}

\pfbegin
Choose $t \in Z(G)^0(k)$ such that $a\chi(t) \neq a'\chi'(t)$ for $(\chi,a) \neq (\chi',a')$ with
\[ (\chi,a), (\chi',a') \in \gY(g)_m \cup \gY(g)_m \mu^{m'-m} \cup \gY(g)_m \mu^{m''-m} \cup \gY(g)_{m'} \cup \gY(g)_{m''}. \]
(Note that it suffices to choose $t \in Z(G)^0(k)$ such that for
\[ (\chi,a), (\chi',a') \in \gY(g)_m \cup \gY(g)_m \mu^{m'-m} \cup \gY(g)_m \mu^{m''-m} \cup \gY(g)_{m'} \cup \gY(g)_{m''}, \]
with $\chi \neq \chi'$ we have $(\chi/\chi')(t) \neq a'/a$.)
We can decompose 
\[ \gY(g)_m = \gY(g)_{m,m',m''}^1 \amalg \gY(g)_{m,m',m''}^2 \]
into two $d$-element multisets, such that
\[ \{ a\chi(t):\,\, (\chi,a) \in \gY(g)_{m,m',m''}^1\} = \Sigma_{gt}^1 \]
and
\[ \{ a\chi(t): \,\, (\chi,a) \in  \gY(g)_{m,m',m''}^2 \mu^{-m}\} = \Sigma_{gt}^2. \]
Then
\[ \begin{array}{l} \{ a\chi(t):\,\, (\chi,a) \in \gY(g)_{m'}\} =\\ \{ a\chi(t):\,\, (\chi,a) \in \gY(g)_{m,m',m''}^1\} \amalg \{ a\chi(t): \,\, (\chi,a) \in  \gY(g)_{m,m',m''}^2 \mu ^{m'-m} \} \end{array} \]
and
\[ \begin{array}{l} \{ a\chi(t):\,\, (\chi,a) \in \gY(g)_{m''}\} = \\ \{ a\chi(t):\,\, (\chi,a) \in \gY(g)_{m,m',m''}^1\} \amalg \{ a\chi(t): \,\, (\chi,a) \in  \gY(g)_{m,m',m''}^2 \mu ^{m''-m} \}. \end{array} \]
It follows that
\[ \gY(g)_{m'} = \gY(g)_{m,m',m''}^1 \amalg \gY(g)_{m,m',m''}^2 \mu^{m'-m} \]
and
\[ \gY(g)_{m'} = \gY(g)_{m,m',m''}^1 \amalg \gY(g)_{m,m',m''}^2 \mu^{m''-m}. \]

If $\mu^{m-m'} \neq \chi/\chi'$ for all $\chi, \chi' \in \gX_m$ then
\[ \gY(g)_{m,m',m''}^1= \gY(g)_{m} \cap \gY(g)_{m'}, \]
so the uniqueness assertion is clear.
\pfend

\begin{cor} If $m \in \cM$, then we can uniquely decompose
\[ \gY(g)_m = \gY(g)_{m}^1 \amalg \gY(g)_{m}^2 \]
into two $d$-element multisets, such that for all $m'\in \cM$ we have
\[ \gY(g)_{m'} = \gY(g)_{m}^1 \amalg \gY(g)_{m}^2 \mu^{m'-m}. \]
\end{cor}

\pfbegin Choose $m'$ such that $\mu^{m-m'} \neq \chi/\chi'$ for all $\chi, \chi' \in \gX_m$. Then we see that for all $m'', m''' \in \cM$ we have
\[  \gY(g)_{m,m',m''}^1= \gY(g)_{m,m',m'''}^1 \]
and
\[  \gY(g)_{m,m',m''}^2= \gY(g)_{m,m',m'''}^2. \]
Then we can simply take $\gY(g)_m^i=\gY(g)_{m,m',m''}^i$.
\pfend

\begin{cor} For all $m,m' \in \cM$ we have
\[ \gY(g)_{m'}^1=\gY(g)_m^1 \]
and
\[ \gY(g)_{m'}^2=\gY(g)_m^2 \mu^{m'-m}. \]
\end{cor}

\pfbegin It is immediate from the previous corollary that $\gY(g)_m^1$ and $\gY(g)_m^2 \mu^{m'-m}$ have the properties that uniquely characterize $\gY(g)_{m'}^1$ and $\gY(g)_{m'}^2$. \pfend

\begin{cor} For all $g \in G$ and $m \in \cM$ and for $i=1,2$ we have
\[ \gY(1)_m^i = \{ (\chi,1):\,\, (\chi,a) \in \gY(g)_m^i\} . \]
\end{cor}

\pfbegin It is again immediate that $\{ (\chi,1):\,\, (\chi,a) \in \gY(g)_m^1\}$ and $\{ (\chi,1):\,\, (\chi,a) \in \gY(g)_m^2\}$ have the properties that uniquely characterize $\gY(1)_{m}^1$ and $\gY(1)_{m}^2$. 
\pfend

We set
\[ \gX_m^i=\{ \chi: \,\, (\chi,1) \in \gY(1)_m^i \}. \]
Note that
\[ \gX_{m'}^1=\gX_m^1 \]
and that
\[ \gX_{m'}^2=\gX_m^2 \mu^{m'-m}. \]

\begin{cor} For all but finitely many $m \in \cM$ the multisets $\gX_m^1$ and $\gX_m^2$ are disjoint. 
 \end{cor}
 
 Let $\cM'$ denote the set of $m \in \cM$ such that $\gX_m^1$ and $\gX_m^2$ are disjoint. Then we see that for $m \in \cM'$ we have
\[ \gY(g)_m^i = \{ (\chi,a)\in \gY(g)_m:\,\, \chi \in \gX_m^i\}. \] 
Moreover for $m \in \cM'$ we may decompose
\[ R_m = R_m^1 \oplus R_m^2 \]
where $R_m^i$ is the sum of the $\chi$-eigenspaces of $Z(G)^0$ for $\chi \in \gX_m^i$. We see that
the multi-set of roots of the characteristic polynomial of $R_m^i(g)$ equals
\[ \{ a: \,\, (\chi,a) \in \gY(g)_m^i\}. \]
Thus $R_m^1$ is independent of $m \in \cM'$, as is $R_m^2 \mu^{-m}$. Denote these representations of $G$ by $r_1$ and $r_2$, so that
\[ R_m \cong r_1 \oplus r_2 \mu^m \]
for all $m \in \cM'$. From corollary \ref{113} (applied to $\cM'$) we see that if $g \in G$ and $\mu(g)$ has infinite order then $\Sigma^i_g$ is the multiset of roots of the characteristic polynomial of $r^i(g)$. Thus we have proved the following result.

\begin{prop}\label{grpth} Keep the notation and assumptions of the first two paragraphs of this section. Then there are continuous semi-simple representations
\[ \rho^i: \Gamma \lra GL_d(k) \]
for $i=1,2$ such that for all $f \in \gF$ the multiset of roots of the characteristic polynomial of $\rho^i(f)$ equals $\cE_f^i$. \end{prop}

This proposition allows us to deduce our main theorem from corollary \ref{mt1}.

\begin{thm} \label{mt2} Suppose that $\pi$ is a cuspidal automorphic representation of $GL_n(\A_F)$ such that $\pi_{\infty}$ has the same infinitesimal character as an algebraic representation of $\RS^F_\Q GL_n$. Then there is a continuous semi-simple representation 
\[ r_{p,\imath}(\pi):G_F \lra GL_{n}(\barQQ_p) \]
such that, if $q \neq p$ is a prime above which $\pi$ and $F$ are unramified and if $v|q$ is a prime of $F$, then $r_{p,\imath}(\pi)$ is unramified at $v$ and
\[ r_{p,\imath}(\pi)|_{W_{F_v}}^\semis = \imath^{-1} \rec_{F_v}(\pi_{v}|\det|_v^{(1-n)/2}). \]
\end{thm}

\pfbegin 
We may suppose that $n>1$, as in the case $n=1$ the result is well known.
Let $S$ denote the set of rational primes above which $F$ or $\pi$ ramifies together with $p$; and let $G_{F,S}$ denote the Galois group over $F$ of the maximal extension of $F$ unramified outside $S$.
Apply proposition \ref{grpth} to $\Gamma=G_{F,S}$, and $k=\barQQ_p$, and $\mu=\epsilon_p^{-2}$, and $\cM$ consisting of all sufficiently large integers, and $\rho_m=R_{p,\imath}(\pi,m)$ (as in theorem \ref{mt1}), and $\cF$ the set of Frobenius elements at primes not above $S$, and $\cE_{\Frob_v}^1$ equal to the multiset of roots of the characteristic polynomial of $
\imath^{-1}\rec_{F_v}(\pi_{v}|\det|_v^{(1-n)/2})(\Frob_v)$, and $\cE_{\Frob_v}^2$ equal to the multiset of roots of the characteristic polynomial of $\imath^{-1}\rec_{F_{{}^cv}}(\pi_{{}^cv}|\det|_{{}^cv}^{(1+3n)/2})(\Frob_{{}^cv}^{-1})$. 
\pfend

\begin{cor}\label{mt3}  Suppose that $E$ is a totally real or CM field and that $\pi$ is a cuspidal automorphic representation such that $\pi_{\infty}$ has the same infinitesimal as an algebraic representation of $\RS^E_\Q GL_n$. Then there is a continuous semi-simple representation 
\[ r_{p,\imath}:G_E \lra GL_{n}(\barQQ_p) \]
such that, if $q \neq p$ is a prime above which $\pi$ and $E$ are unramified and if $v|q$ is a prime of $E$, then $r_{p,\imath}(\pi)$ is unramified at $v$ and
\[ r_{p,\imath}(\pi)|_{W_{E_v}}^\semis = \imath^{-1} \rec_{E_v}(\pi_{v}|\det|_v^{(1-n)/2}). \]
\end{cor}

\pfbegin This can be deduced from theorem \ref{mt2} by using lemma 1 of \cite{sorensen}. (This is the same argument used in the proof of theorem VII.1.9 of \cite{ht}.)
\pfend


\newpage



\begin{thebibliography}{BLGGT} 

\bibitem[AIP1]{iovita1} F. Andreatta, A. Iovita and V. Pilloni, {\em $p$-adic families of Siegel modular cuspforms}, preprint available at \verb+http://www.mathstat.concordia.ca/faculty/iovita/research.html+

\bibitem[AIP2]{iovita2} F. Andreatta, A. Iovita and V. Pilloni, {\em $p$-adic families of Hilbert modular cuspforms}, in preparation.

\bibitem[AMRT]{amrt} A. Ash, D. Mumford, M. Rapoport and Y.-S. Tai, {\em Smooth compactifications of locally symmetric varieties}, 2nd edition, CUP 2010. 

\bibitem[Be1]{berthelot} P. Berthelot, {\em Cohomologie rigide et cohomologie rigide \`{a} supports propres}, pr\'{e}publication IRMAR 96-03 (1996), available at \verb+http://perso.univ-rennes1.fr/pierre.berthelot/+

\bibitem[Be2]{berthfin} P. Berthelot, {\em Finitude et puret\'{e} cohomologique en cohomologie rigide}, Inv. Math. 128 (1997), 329--377.

\bibitem[BGR]{bgr} S. Bosch, U. G\"{u}nzer and R. Remmert, {\em Non-archimedean analysis}, Springer 1984.

\bibitem[BLGGT]{blggt} T.Barnet-Lamb, T.Gee, D.Geraghty and R.Taylor, {\em Potential automorphy and change of weight}, preprint available at \verb+http://www.math.ias.edu/~rtaylor/+

\bibitem[BLGGT2]{blggtcomp2} T.Barnet-Lamb, T.Gee, D.Geraghty and R.Taylor, {\em Local-global compatibility for $l=p$ II}, preprint available at \verb+http://www.math.ias.edu/~rtaylor/+

\bibitem[BLGHT]{blght} T.Barnet-Lamb, D.Geraghty, M.Harris and R.Taylor, {\em A family of Calabi-Yau varieties and potential automorphy II}, P.R.I.M.S. 47 (2011), 29--98. 

\bibitem[BLR]{blr} S. Bosch, W. L\"{u}tkebohmert and M. Raynaud, {\em N\'{e}ron models}, Springer 1980.

\bibitem[Bo]{borel} A. Borel, {\em Stable real cohomology of arithmetic groups II}, in ``Manifolds and Lie groups, papers  in honor of Y. Matsushima'', Progress in Math. 14, Birkha\"{u}ser, 1981.

\bibitem[Br]{bredon} G. Bredon, {\em Sheaf theory}, $2^{nd}$ edition, Springer 1997.

\bibitem[BW]{bw} A. Borel and N. Wallach, {\em Continuous cohomology, discrete subgroups, and representations of reductive}, Annals of Math. Studies 94, PUP 1980.

\bibitem[Ca]{ana} A. Caraiani, {\em Local-global compatibility and the action of monodromy on nearby cycles}, to appear Duke Math. J.

\bibitem[CF]{cf} C.-L. Chai and G. Faltings, {\em Degeneration of abelian varieties}, Springer 1990.

\bibitem[CH]{ch} G. Chenevier and M. Harris, {\em Construction of automorphic Galois representations II}, to appear Cambridge Math. J. 1 (2013).

\bibitem[Ch]{chiar} B. Chiarellotto, {\em Weights in rigid cohomology applications to unipotent F-isocrystals}, Ann. Sci. E.N.S. 31 (1998), 683--715.

\bibitem[Cl]{clozelaa} L. Clozel, {\em Motifs et formes automorphes: applications du principe de fonctorialit\'{e}}, in `Automorphic forms, Shimura varieties, and L-functions I' Perspect. Math., 10, Academic Press 1990.

\bibitem[CO]{co} W. Casselman and M.S. Osborne, {\em The $n$-cohomology of representations with an infinitesimal character}, Comp. Math. 31 (1975), 219--227.

\bibitem[EGA3]{ega3} J. Dieudonn\'{e} and A. Grothendieck, {\em  El\'{e}ments de
g\'{e}om\'{e}trie alg\'{e}brique III}, Pub. Math. IHES 11 (1961), 5--167 and 17 (1963), 5--91.

\bibitem[EGA4]{ega4} J. Dieudonn\'{e} and A. Grothendieck, {\em  El\'{e}ments de
g\'{e}om\'{e}trie alg\'{e}brique IV}, Pub. Math. IHES 20 (1964) 5--259 and 24 (1965) 5--231 and 28 (1966) 5--255 and 32 (1967) 5--361.

\bibitem[Fr]{rs} G. Friedman, {\em An elementary illustrated introduction to simplicial sets}, preprint available at arXiv:0809.4221

\bibitem[Fu]{fulton} W. Fulton, {\em Introduction to toric varieties}, Annals of Math. Studies 131, PUP 1993.

\bibitem[GK]{gkcrelle} E. Grosse-Kl\"{o}nne, {\em Rigid analytic spaces with overconvergent structure sheaf}, J. Reine Angew. Math. 519 (2000), 73--95.

\bibitem[Ha]{dbar} M. Harris, {\em Automorphic forms of $\overline{\partial}$-cohomology type as coherent cohomology classes}, J. Diff. Geom. 32 (1990), 1--63.

\bibitem[Hi]{hidabook} H. Hida, {\em $p$-adic automorphic forms on Shimura varieties}, Springer 2004. 

\bibitem[HT]{ht} M. Harris and R. Taylor, {\em The geometry and cohomology of some simple Shimura varieties}, Annals of Math. Studies, 151, PUP,  2001.

\bibitem[HZ1]{hz1} M. Harris and S. Zucker, {\em Boundary cohomology of Shimura varieties, I:   coherent cohomology on the toroidal boundary},  Annales Scient. de l'Ec. Norm. Sup. 27 (1994), 249--344. 

\bibitem[HZ2]{hz2} M. Harris and S. Zucker, {\em Boundary cohomology of Shimura varieties, II:  mixed Hodge structures},  Inventiones Math.116 (1994), 243-307; erratum, Inventiones Math.  123 (1995), 437. 

\bibitem[HZ3]{hz3} M. Harris and S. Zucker, {\em Boundary cohomology of Shimura varieties, III:  Coherent cohomology on higher-rank boundary strata and applications to Hodge theory}, M\'{e}moires de la SMF 85 (2001). 

\bibitem[Ka1]{katz} N. Katz, {\em P-adic properties of modular schemes and modular forms}, in ``Modular functions of one variable, III'', LNM 350, Springer 1973.

\bibitem[Ka2]{katzst} N. Katz, {\em Serre-Tate local moduli}, in `Algebraic surfaces (Orsay, 1976--78)' LNM 868, Springer 1981.

\bibitem[KM]{km} N. Katz and B. Mazur, {\em Arithmetic moduli of elliptic curves}, Annals of Mathematics Studies 108, PUP 1985.

\bibitem[K\"o]{kopf} U. K\"opf, {\em \"{U}ber eigentliche Familien algebraischer Variet\"{a}ten \"{u}ber affinoiden R\"{a}umen}, Schr. Math. Inst. Univ. M\"{u}nster (2) Heft 7 (1974)

\bibitem[Ku]{kunz} E. Kunz, {\em K\"{a}hler differentials}, Viehweg 1986.

\bibitem[La1]{kw1} K.-W. Lan, {\em Arithmetic compactifications of PEL-type
Shimura varieties}, London Mathematical Society Monographs, vol. 36, PUP 2013. 

\bibitem[La2]{kw1.5} K.-W. Lan, {\em Toroidal compactifications of PEL-type Kuga families}, Algebra and Number Theory 6 (2012), 885--966.

\bibitem[La3]{kwan} K.-W. Lan, {\em Comparison between analytic and algebraic constructions of toroidal compactifications of PEL-type Shimura varieties}, J. Reine Angew. Math. 664 (2012), 163--228.

\bibitem[La4]{kw2} K.-W. Lan, {\em Compactifications of PEL-type Shimura varieties and Kuga families with ordinary loci}, preprint available at \verb+http://www.math.umn.edu/~kwlan/academic.html+

\bibitem[LeS]{lestum} B. Le Stum, {\em Rigid cohomology}, Cambridge Tracts in Math. 172, CUP, 2007

\bibitem[M]{mumav} D. Mumford, {\em Abelian varieties}, Tata Institute of Fundamental Research Studies in Mathematics 5, OUP, 1970.

\bibitem[MW]{mw} C. Moeglin and J.-L. Waldspurger, {\em Le spectre r\'{e}siduel de $GL(n)$}, Ann. Sci. Ecole Norm. Sup. (4) 22 (1989), 605--674.

\bibitem[Pi]{pink} R. Pink, {\em Arithmetical compactification of mixed Shimura varieties}, Bonner Mathematische Schriften 209, Universit\"{a}t Bonn, Mathematisches Institut, Bonn, 1990.

\bibitem[Se]{serreecc} J.-P. Serre, {\em Endomorphismes compl\`{e}tements continues des espaces de
Banach p-adiques}, Publ. Math. I.H.E.S., 12 (1962) 69--85.

\bibitem[Sh1]{shin} S.-W. Shin, {\em Galois representations arising from some compact Shimura varieties}, Ann. Math. 173 (2011), 1645--1741.

\bibitem[Sh2]{sws} S.-W. Shin, {\em On the cohomological base change for unitary similitude groups}, (an appendix to the forthcoming paper by W. Goldring), preprint available at \verb+http://math.mit.edu/~swshin/+

\bibitem[So]{sorensen} C. Sorensen, {\em A patching lemma}, preprint available at \verb+http://fa.institut.math.jussieu.fr/node/45+

\bibitem[T]{pseudo} R. Taylor, {\em Galois representations associated to Siegel modular forms of low weight}, Duke Math. J. 63 (1991),  281--332.

\bibitem[TY]{ty} R. Taylor and T. Yoshida, {\em Compatibility of local and global Langlands correspondences}, 
J.A.M.S. 20 (2007), 467--493. 

\bibitem[W]{weibel} C. Weibel, {\em An introduction to homological algebra}, Cambridge studies in advanced mathematics 38, CUP 1994.
 
\end{thebibliography}
\end{document}